%% file: thesis_arXiv.tex
\documentclass[12pt,a4paper,ngerman,english,intoc,bibliography=totoc,index=totoc,BCOR10mm,titlepage,captions=tableheading,fceqn]{scrbook}
\usepackage{lmodern}

\usepackage[T1]{fontenc}
\usepackage[latin9]{inputenc}
\usepackage{fancyhdr}
\usepackage{color}
\usepackage{babel}
\usepackage{amsmath}
\usepackage{amssymb}
\usepackage{graphicx}
\usepackage{esint}
\usepackage{nomencl}
\providecommand{\printnomenclature}{\printglossary}
\providecommand{\makenomenclature}{\makeglossary}
\makenomenclature
\usepackage[unicode=true,
 bookmarks=true,bookmarksnumbered=true,bookmarksopen=true,bookmarksopenlevel=0,
 breaklinks=false,pdfborder={0 0 0},backref=false,colorlinks=true]
 {hyperref}
\hypersetup{pdftitle={Optimal trajectory tracking},
 pdfauthor={Jakob Löber},
 pdfsubject={Dissertation Theoretical Physics},
 pdfkeywords={trajectory tracking, optimal control, nonlinear control, exact linearization},
 pdfpagelayout=OneColumn, pdfnewwindow=true, pdfstartview=XYZ, plainpages=false}
\usepackage{breakurl}

\makeatletter

\@ifundefined{date}{}{\date{}}
\usepackage{custom}
\usepackage{doi}
\usepackage[round, semicolon, authoryear]{natbib}
\renewcommand\cite{\citep}
\usepackage{caption}

\usepackage[pass]{geometry}
\usepackage{pdfpages}

\@ifundefined{showcaptionsetup}{}{%
 \PassOptionsToPackage{caption=false}{subfig}}
\usepackage{subfig}
\makeatother

\begin{document}
\pagestyle{plain}

\lhead{\rightmark}

\rhead[\leftmark]{}

\lfoot[\thepage]{}

\cfoot{}

\rfoot[]{\thepage}

\title{Optimal trajectory tracking}

\author{Jakob Löber}

\lowertitleback{2nd corrected version\\TU Berlin, 16 December 2015}

\maketitle
\cleardoublepage{}

\include{Summary}\cleardoublepage{}\include{Zusammenfassung}\cleardoublepage{}\pdfbookmark{\contentsname}{Contents}
\pagestyle{plain}

\lhead{\rightmark}

\rhead[\leftmark]{}

\lfoot[\thepage]{}

\cfoot{}

\rfoot[]{\thepage}

\tableofcontents{}\cleardoublepage{}\pagestyle{plain}

\lhead{\rightmark}

\rhead[\leftmark]{}

\lfoot[\thepage]{}

\cfoot{}

\rfoot[]{\thepage}

\listoffigures
\cleardoublepage{}

\pagestyle{fancy}
\mainmatter

\lhead[\chaptername~\thechapter]{\rightmark}

\include{chapter-1}

\include{chapter-2}

\include{chapter-3}

\include{chapter-4}

\include{chapter-5}

\appendix
\include{Appendix}

\cleardoublepage{}

\lhead[]{\rightmark}

\rhead[\leftmark]{}

\bibliographystyle{myabbrvnat}
\bibliography{literature}

\cleardoublepage{}

\lhead[]{Nomenclature}

\rhead[Nomenclature]{}

\printnomenclature[2.5cm]{}
\end{document}

%% file: Summary.tex
\chapter*{Summary}

\addcontentsline{toc}{chapter}{Summary} 

\input{SummaryText.tex}

%% file: SummaryText.tex
This thesis investigates optimal trajectory tracking of nonlinear
dynamical systems with affine controls. The control task is to enforce
the system state to follow a prescribed desired trajectory as closely
as possible. The concept of so-called exactly realizable trajectories
is proposed. For exactly realizable desired trajectories exists a
control signal which enforces the state to exactly follow the desired
trajectory.\\
This approach does not only yield an explicit expression for the control
signal in terms of the desired trajectory, but also identifies a particularly
simple class of nonlinear control systems. Systems in this class satisfy
the so-called linearizing assumption and share many properties with
linear control systems. For example, conditions for controllability
can be formulated in terms of a rank condition for a controllability
matrix analogously to the Kalman rank condition for linear time invariant
systems.\\
Furthermore, exactly realizable trajectories arise as solutions to
unregularized optimal control problems. Based on that insight, the
regularization parameter is used as the small parameter for a perturbation
expansion. This results in a reinterpretation of affine optimal control
problems with small regularization term as singularly perturbed differential
equations. The small parameter originates from the formulation of
the control problem and does not involve simplifying assumptions about
the system dynamics. Combining this approach with the linearizing
assumption, approximate and partly linear equations for the optimal
trajectory tracking of arbitrary desired trajectories are derived.\\
For vanishing regularization parameter, the state trajectory becomes
discontinuous and the control signal diverges. On the other hand,
the analytical treatment becomes exact and the solutions are exclusively
governed by linear differential equations. Thus, the possibility of
linear structures underlying nonlinear optimal control is revealed.
This fact enables the derivation of exact analytical solutions to
an entire class of nonlinear trajectory tracking problems with affine
controls. This class comprises, among others, mechanical control systems
in one spatial dimension and the FitzHugh-Nagumo model with a control
acting on the activator.

%% file: chapter-1.tex
\lhead[\chaptername~\thechapter\leftmark]{}

\rhead[]{\rightmark}

\lfoot[\thepage]{}

\cfoot{}

\rfoot[]{\thepage}

\chapter{\label{chap:Introduction}Introduction}

\nomenclature[100]{$n$}{dimensionality of state space}\nomenclature[101]{$p$}{dimension of vector of control signals}\nomenclature[1015]{$m$}{dimension of vector of output}\nomenclature[102]{$t$}{time}\nomenclature[103]{$t_{0}$}{initial time}\nomenclature[104]{$t_{1}$}{terminal time}\nomenclature[105]{$\boldsymbol{x}\left(t\right)$}{$n$-dimensional time-dependent state vector}\nomenclature[106]{$\boldsymbol{x}_{d}\left(t\right)$}{$n$-dimensional desired trajectory}\nomenclature[1060]{$\Delta\boldsymbol{x}\left(t\right)$}{difference between state and desired trajectory}\nomenclature[107]{$\boldsymbol{x}_{0}$}{$n$-dimensional vector of initial states}\nomenclature[108]{$\boldsymbol{x}_{1}$}{$n$-dimensional vector of terminal states}\nomenclature[109]{$\boldsymbol{u}\left(t\right)$}{$p$-dimensional vector of control or input signals}\textit{\nomenclature[110]{$\boldsymbol{\mathcal{A}}$}{$n \times n$ state matrix}\nomenclature[111]{$\boldsymbol{R}\left(\boldsymbol{x}\right)$}{nonlinear kinetics of an $n$-dimensional dynamical system}\nomenclature[112]{$\boldsymbol{\mathcal{B}}\left(\boldsymbol{x}\right)$}{$n \times p$ input or coupling matrix}}\nomenclature[113]{$\boldsymbol{\mathcal{P}}$}{$n\times n$ projector}\nomenclature[114]{$\boldsymbol{\mathcal{Q}}$}{$n\times n$ projector}\nomenclature[1140]{$\boldsymbol{\mathcal{A}}^{T}$}{transpose of matrix $\boldsymbol{\mathcal{A}}$}\nomenclature[115]{$\boldsymbol{\mathcal{A}}^{+}$}{Moore-Penrose pseudo inverse of matrix $\boldsymbol{\mathcal{A}}$}\nomenclature[115]{$\boldsymbol{\mathcal{A}}^{g}$}{generalized inverse of matrix $\boldsymbol{\mathcal{A}}$}\nomenclature[116]{$\nabla\boldsymbol{R}\left(\boldsymbol{x}\right)$}{$n\times n$ Jacobi matrix of $\boldsymbol{R}\left(\boldsymbol{x}\right)$}\nomenclature[117]{$\boldsymbol{\Phi}\left(t,t_{0}\right)$}{$n\times n$ state transition matrix for a linear dynamical system}\nomenclature[118]{$\boldsymbol{\mathcal{K}}$}{Kalman's $n\times np$ controllability matrix}\nomenclature[119]{$\boldsymbol{\mathcal{\tilde{K}}}$}{$n\times n^2$ controllability matrix for exactly realizable trajectories}\textit{\nomenclature[120]{$\boldsymbol{\mathcal{\tilde{K}}}_{\boldsymbol{\mathcal{N}}}$}{$n\times (n^2 + n)$ output trajectory realizability matrix}}\nomenclature[121]{$\boldsymbol{\mathcal{C}}$}{$p\times n$ output matrix}\nomenclature[122]{$\mathcal{J}$}{target functional of optimal control}\nomenclature[122]{$\epsilon$}{small regularization parameter}\nomenclature[123]{$\boldsymbol{\lambda}\left(t\right)$}{$n$-dimensional vector of time-dependent co-states for optimal control}\nomenclature[123]{$\boldsymbol{\mathcal{S}}$}{$n\times n$ symmetric matrix of weighting coefficients}\nomenclature[124]{$\boldsymbol{\mathcal{P}}_{\boldsymbol{\mathcal{S}}}$}{$n\times n$ projector}\nomenclature[125]{$\boldsymbol{\mathcal{Q}}_{\boldsymbol{\mathcal{S}}}$}{$n\times n$ projector}\nomenclature[1260]{$\boldsymbol{x}_{O}\left(t\right)$}{state vector for outer equations}\nomenclature[1261]{$\tau_{L}$}{rescaled time for left inner equations}\nomenclature[1262]{$\tau_{R}$}{rescaled time for right inner equations}\nomenclature[1263]{$\boldsymbol{X}_{L}\left(\tau_{L}\right)$}{rescaled state vector for left inner equations}\nomenclature[1264]{$\boldsymbol{X}_{R}\left(\tau_{R}\right)$}{rescaled state vector for right inner equations}

Science often begins with the discovery of physical phenomena. The
second step is to describe, understand, and predict them, often in
terms of mathematical theories. The final step is to take advantage
of the discovered phenomena. This last step is the topic of control
theory. 

Section \ref{sec:AffineControlSystems} introduces the notation for
control systems. Some examples of affine control systems, which are
used repeatedly throughout the thesis to demonstrate the developed
concepts, are presented in Section \ref{sec:ExamplesOfAffineControlSystems}.
Section \ref{sec:OptimalTrajectoryTracking} illustrates the main
result of this thesis by means of an example.

\section{\label{sec:AffineControlSystems}Affine control systems}

The subject of this thesis are controlled dynamical systems of the
form
\begin{align}
\boldsymbol{\dot{x}}\left(t\right) & =\boldsymbol{R}\left(\boldsymbol{x}\left(t\right)\right)+\boldsymbol{\mathcal{B}}\left(\boldsymbol{x}\left(t\right)\right)\boldsymbol{u}\left(t\right),\label{eq:AffineNonlinearDynamicalSystem}\\
\boldsymbol{x}\left(t_{0}\right) & =\boldsymbol{x}_{0}.\label{eq:AffineNonlinearDynamicalSystemInitCond}
\end{align}
Here, $t$ is the time, $\boldsymbol{x}\left(t\right)=\left(x_{1}\left(t\right),\dots,x_{n}\left(t\right)\right)^{T}\in\mathbb{R}^{n}$
is called the \textit{state vector} with $n$ components and $\boldsymbol{x}^{T}$
denotes the transposed of vector $\boldsymbol{x}$. The dot 
\begin{align}
\boldsymbol{\dot{x}}\left(t\right) & =\dfrac{d}{dt}\boldsymbol{x}\left(t\right)
\end{align}
denotes the time derivative of $\boldsymbol{x}\left(t\right)$. The
vector $\boldsymbol{u}\left(t\right)=\left(u_{1}\left(t\right),\dots,u_{p}\left(t\right)\right)^{T}\in\mathbb{R}^{p}$
with $p\leq n$ components is the vector of \textit{control} or \textit{input}
\textit{signals}. The \textit{nonlinearity} $\boldsymbol{R}$ is a
sufficiently well behaved function mapping $\mathbb{R}^{n}$ to $\mathbb{R}^{n}$,
and $\boldsymbol{\mathcal{B}}$ is a sufficiently well behaved $n\times p$
matrix function called the \textit{coupling matrix }or\textit{ input
matrix}. As a function of the state vector $\boldsymbol{x}$, $\boldsymbol{\mathcal{B}}$
maps from $\mathbb{R}^{n}$ to $\mathbb{R}^{n}$. A \textit{single
input} system has a scalar control signal $u\left(t\right)$, i.e.,
$p=1$, and the coupling matrix $\boldsymbol{\mathcal{B}}\left(\boldsymbol{x}\right)$
is a coupling vector written as $\boldsymbol{B}\left(\boldsymbol{x}\right)$.
The initial condition $\boldsymbol{x}_{0}$ prescribes the value of
the state vector $\boldsymbol{x}$ at the initial time $t_{0}\leq t$.

Regarding the state vector $\boldsymbol{x}$, the system Eq. \eqref{eq:AffineNonlinearDynamicalSystem}
has two possible sources of nonlinearity. First, the nonlinearity
$\boldsymbol{R}\left(\boldsymbol{x}\right)$ typically is a nonlinear
function of the state. This is the nonlinearity encountered in uncontrolled
systems. Second, the coupling matrix $\boldsymbol{\mathcal{B}}\left(\boldsymbol{x}\right)$
may depend nonlinearly on the state $\boldsymbol{x}$. This nonlinearity
is exclusive for control systems. Equation \eqref{eq:AffineNonlinearDynamicalSystem}
is called an \textit{affine control system} because the control signal
$\boldsymbol{u}\left(t\right)$ enters only linearly. Throughout the
thesis, it is assumed that the coupling matrix $\boldsymbol{\mathcal{B}}\left(\boldsymbol{x}\right)$
has full rank for all values of $\boldsymbol{x}$. Because $p\leq n$,
this condition is 
\begin{align}
\text{rank}\left(\boldsymbol{\mathcal{B}}\left(\boldsymbol{x}\right)\right) & =p.\label{eq:BFullRank}
\end{align}
Assumption \eqref{eq:BFullRank} ensures that the maximum number of
$p$ independent control signal acts on the system regardless of the
value of the state vector $\boldsymbol{x}$.

\section{\label{sec:ExamplesOfAffineControlSystems}Examples of affine control
systems}

Some examples of affine control systems are discussed. These examples
are encountered repeatedly to illustrate the developed concepts.

\begin{example}[Mechanical control system in one spatial dimension]\label{ex:OneDimMechSys1}

Newton's equation of motion for a single point mass in one spatial
dimension $x$ is \cite{goldstein2001classical}, 
\begin{align}
\ddot{x}\left(t\right) & =R\left(x\left(t\right),\dot{x}\left(t\right)\right)+B\left(x\left(t\right),\dot{x}\left(t\right)\right)u\left(t\right).\label{eq:OneDimMechanicalSystem}
\end{align}
The point mass is moving in the external force field $R$ which may
depend on position $x$ and velocity $\dot{x}$ of the particle. The
control signal $u\left(t\right)$ couples to the point mass via the
control force $B\left(x\left(t\right),\dot{x}\left(t\right)\right)u\left(t\right)$,
with $u\left(t\right)$ being the control signal and $B$ the coupling
function. Introducing the velocity $y\left(t\right)=\dot{x}\left(t\right)$,
Eq. \eqref{eq:OneDimMechanicalSystem} can be rearranged as an affine
control system,
\begin{align}
\dot{x}\left(t\right) & =y\left(t\right),\label{eq:OneDimMechanicalSystemPosition}\\
\dot{y}\left(t\right) & =R\left(x\left(t\right),y\left(t\right)\right)+B\left(x\left(t\right),y\left(t\right)\right)u\left(t\right).\label{eq:OneDimMechanicalSystemVelocity}
\end{align}
In vector notation, Eqs. \eqref{eq:OneDimMechanicalSystemPosition}
and \eqref{eq:OneDimMechanicalSystemVelocity} become 
\begin{align}
\boldsymbol{\dot{x}}\left(t\right)= & \boldsymbol{R}\left(\boldsymbol{x}\left(t\right)\right)+\boldsymbol{B}\left(\boldsymbol{x}\left(t\right)\right)u\left(t\right),
\end{align}
with
\begin{align}
\boldsymbol{x}\left(t\right) & =\left(x\left(t\right),y\left(t\right)\right)^{T}, & \boldsymbol{R}\left(\boldsymbol{x}\right) & =\left(y,R\left(x,y\right)\right)^{T}, & \boldsymbol{B}\left(\boldsymbol{x}\right) & =\left(0,B\left(x,y\right)\right)^{T}.
\end{align}
The condition of full rank for the coupling vector $\boldsymbol{B}$
is
\begin{align}
\mbox{rank}\left(\boldsymbol{B}\left(\boldsymbol{x}\right)\right) & =1,
\end{align}
which in turn implies
\begin{align}
B\left(x,y\right) & \neq0
\end{align}
for all values of $x$ and $y$. Note that the control force acts
on the nonlinear equation \eqref{eq:OneDimMechanicalSystemVelocity}
for the velocity $y\left(t\right)$, while the remaining equation
\eqref{eq:OneDimMechanicalSystemPosition} for the position $x\left(t\right)$
is linear.

\end{example}

\begin{example}[FitzHugh-Nagumo model]\label{ex:FHN1}

The FitzHugh-Nagumo (FHN) model \cite{fitzhugh1961impulses,nagumo1962active}
is a simple nonlinear model describing a prototype excitable system
\cite{izhikevich2007dynamical}. It arose as a simplified version
of the Hodgkin-Huxley model which describes action-potential dynamics
in neurons \cite{hodgkin1952quantitative,keener1998mathematical1}
and contains the Van der Pol oscillator as a special case \cite{vanderPol1926lxxxviii}.
The model consists of two variables called the inhibitor $x$ and
activator $y$,
\begin{align}
\left(\begin{array}{c}
\dot{x}\left(t\right)\\
\dot{y}\left(t\right)
\end{array}\right) & =\left(\begin{array}{c}
a_{0}+a_{1}x\left(t\right)+a_{2}y\left(t\right)\\
R\left(x\left(t\right),y\left(t\right)\right)
\end{array}\right)+\boldsymbol{\mathcal{B}}\left(\boldsymbol{x}\left(t\right)\right)\boldsymbol{u}\left(t\right).\label{eq:FHNEquation}
\end{align}
The function $R\left(x,y\right)$ is given by 
\begin{align}
R\left(x,y\right) & =R\left(y\right)-x,\label{eq:Rxy}
\end{align}
with $R\left(y\right)$ being a cubic polynomial of the form 
\begin{align}
R\left(y\right) & =y-\frac{1}{3}y^{3}.\label{eq:Ry}
\end{align}
The nonlinearity is linear in the inhibitor $x$ but nonlinear in
the activator $y$. Unless otherwise announced, a set of standard
parameter values
\begin{align}
a_{0} & =0.056, & a_{1} & =-0.064, & a_{2} & =0.08
\end{align}
is used for numerical simulations. Because this model is not a mechanical
system, it is not predefined in which way a control acts on the system.
Several simple choices with a constant coupling matrix $\boldsymbol{\mathcal{B}}\left(\boldsymbol{x}\right)=\boldsymbol{\mathcal{B}}$
are possible. A control acting on the activator equation leads to
a coupling vector $\boldsymbol{B}=\left(\begin{array}{cc}
0, & 1\end{array}\right)^{T}$, while a control acting on the inhibitor equation gives $\boldsymbol{B}=\left(\begin{array}{cc}
1, & 0\end{array}\right)^{T}$. The former is called the \textit{activator-controlled FitzHugh-Nagumo
model}, while the latter is named \textit{inhibitor-controlled FitzHugh-Nagumo
model}. The simplest case occurs if the number of independent control
signals equals the number of state components, $p=n$, and the coupling
matrix $\boldsymbol{\mathcal{B}}$ attains the form
\begin{align}
\boldsymbol{\mathcal{B}} & =\left(\begin{array}{cc}
1 & 0\\
0 & 1
\end{array}\right).
\end{align}
Note that Eq. \eqref{eq:FHNEquation} reduces to a mechanical control
system in one spatial dimension with external force $R\left(x,y\right)$
for the parameter values $a_{0}=a_{1}=0$, and $a_{2}=1$, and a coupling
vector $\boldsymbol{B}\left(\boldsymbol{x}\right)=\left(\begin{array}{cc}
0, & B\left(x,y\right)\end{array}\right)^{T}$.

Strictly speaking, only a model with $R\left(x,y\right)$ given by
a cubic polynomial in $y$ and linear in $x$ is called the FitzHugh-Nagumo.
The approach to control developed here identifies Eq. \eqref{eq:FHNEquation}
with arbitrary nonlinearity $R\left(x,y\right)$ and coupling vector
$\boldsymbol{B}=\left(\begin{array}{cc}
0, & 1\end{array}\right)^{T}$ as a particularly simple form of two-dimensional controlled dynamical
systems. In absence of a better name, this model is occasionally called
the activator-controlled FHN model as well.

\end{example}

\begin{example}[SIR model]\label{ex:SIRModel1}

The SIR-model is a nonlinear dynamical system to describe the transmission
of a disease among a population \cite{bailey1975mathematical,murray1993mathematicalbiology1,murray1993mathematicalbiology2}.
The original model was created by Kermack and McKendrick in 1927 \cite{kermack1927contribution}
and consists of three components 
\begin{align}
\dot{S}\left(t\right) & =-\beta\frac{S\left(t\right)I\left(t\right)}{N},\label{eq:SIRS}\\
\dot{I}\left(t\right) & =\beta\frac{S\left(t\right)I\left(t\right)}{N}-\gamma I\left(t\right),\label{eq:SIRI}\\
\dot{R}\left(t\right) & =\gamma I\left(t\right).\label{eq:SIRR}
\end{align}
The variable $S$ denotes the number of susceptible individuals. If
susceptible individuals $S$ come in contact with infected individuals
$I$, they become infected with a transmission rate $\beta=0.36$.
The average period of infectiousness is set to $1/\gamma=5$ days,
after which infected individuals either recover or die. Both possibilities
are collected in the variable $R$. Recovered or dead individuals
$R$ are immune and do not become susceptible again. The total population
number $N=S\left(t\right)+I\left(t\right)+R\left(t\right)$ is constant
in time because 
\begin{align}
\dot{S}\left(t\right)+\dot{I}\left(t\right)+\dot{R}\left(t\right) & =0.
\end{align}
No exact analytical solution to Eqs. \eqref{eq:SIRS}-\eqref{eq:SIRR}
is known. Figure \ref{fig:SIR} shows a typical time evolution of
an epidemic obtained by numerical simulations.

\begin{minipage}{1.0\linewidth}
\begin{center}
\includegraphics[scale=0.85]{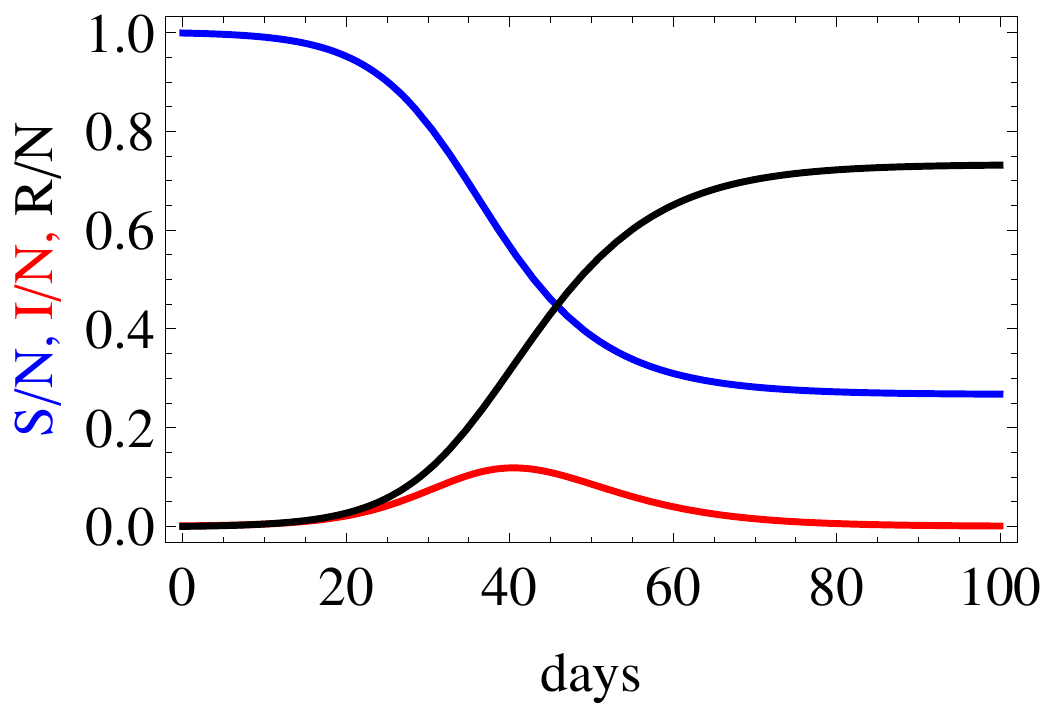}
\captionof{figure}[Time evolution of an epidemic according to the SIR model]{\label{fig:SIR}Time evolution of an epidemic according to the SIR model. Time is measured in days. Initially, almost all individuals are susceptible (blue) and only very few are infected (red). The number of infected individuals reaches a maximum and subsequently decays to zero, with only susceptible and recovered individuals (black) remaining.}
%"/home/jakob/svn/Control/SIR/SIRWithVaccination.nb"
\end{center}
\end{minipage}

The reproductive number $R_{0}$ defined as 
\begin{align}
R_{0} & =\frac{\beta}{\gamma}
\end{align}
is the average number of susceptible individuals an infectious individual
is infecting. To prevent further spreading of the epidemics, this
number must be $R_{0}<1$.

There are two parameters $\beta$ and $\gamma$ in the system which
can be affected by control measures. Culling of infected domestic
animals to increase the rate $\gamma$ is a common procedure but it
is out of question for humans. A common measure during an epidemic
among humans is to separate the infected persons from the susceptible
persons. In the framework of the SIR model, such measures decrease
the transmission rate $\beta\left(t\right)$. Thus, the transmission
rate becomes time dependent and is of the form
\begin{align}
\beta\left(t\right) & =\beta+u\left(t\right).
\end{align}
Here, $\beta$ is the constant transmission rate of the uncontrolled
system and $u\left(t\right)$ is the control signal. The controlled
SIR model investigated in this thesis is
\begin{align}
\boldsymbol{\dot{x}}\left(t\right)= & \boldsymbol{R}\left(\boldsymbol{x}\left(t\right)\right)+\boldsymbol{B}\left(\boldsymbol{x}\left(t\right)\right)u\left(t\right),
\end{align}
with
\begin{align}
\boldsymbol{x}\left(t\right) & =\left(S\left(t\right),I\left(t\right),R\left(t\right)\right)^{T},\\
\boldsymbol{R}\left(\boldsymbol{x}\left(t\right)\right) & =\left(-\beta\frac{S\left(t\right)I\left(t\right)}{N},\beta\frac{S\left(t\right)I\left(t\right)}{N}-\gamma I\left(t\right),\gamma I\left(t\right)\right)^{T},\\
\boldsymbol{B}\left(\boldsymbol{x}\left(t\right)\right) & =\left(-\frac{S\left(t\right)I\left(t\right)}{N},\frac{S\left(t\right)I\left(t\right)}{N},0\right)^{T}.
\end{align}
\end{example}

\section{\label{sec:OptimalTrajectoryTracking}Optimal trajectory tracking}

An important control objective is the guidance of state trajectories
of dynamical systems along a desired reference trajectory. The reference
trajectory is called the \textit{desired trajectory} and denoted as
$\boldsymbol{x}_{d}\left(t\right)\in\mathbb{R}^{n}$. It has the same
number $n$ of components as the system's state $\boldsymbol{x}\left(t\right)$
and is defined for a time interval $t_{0}\leq t\leq t_{1}$. The closer
the controlled state trajectory $\boldsymbol{x}\left(t\right)$ follows
the reference trajectory $\boldsymbol{x}_{d}\left(t\right)$, the
better the control target is achieved. In the ideal case, the desired
trajectory $\boldsymbol{x}_{d}\left(t\right)$ is exactly equal to
the actual controlled state trajectory $\boldsymbol{x}\left(t\right)$.
A convenient measure for the distance between a desired trajectory
$\boldsymbol{x}_{d}\left(t\right)$ and the actual trajectory $\boldsymbol{x}\left(t\right)$
of the controlled system is the squared difference integrated over
the time interval, 
\begin{align}
\mathcal{J}\left[\boldsymbol{x}\left(t\right)\right] & =\frac{1}{2}\intop_{t_{0}}^{t_{1}}dt\left(\boldsymbol{x}\left(t\right)-\boldsymbol{x}_{d}\left(t\right)\right)^{2}.\label{eq:JDef}
\end{align}
The quantity $\mathcal{J}\left[\boldsymbol{x}\left(t\right)\right]$
defined in Eq. \eqref{eq:JDef} is a functional of the state vector
$\boldsymbol{x}\left(t\right)$ over the time interval $t\in\left[t_{0},t_{1}\right]$.
It defines a distance between trajectories, i.e., a distance in function
space. Optimal trajectory tracking aims to find the control signal
$\boldsymbol{u}\left(t\right)$ such that $\mathcal{J}\left[\boldsymbol{x}\left(t\right)\right]$
is minimal. If this control solution exists and is unique, no ``better''
control signal exists. Any other control would result in a controlled
state trajectory $\boldsymbol{x}\left(t\right)$ with a larger distance
to the desired trajectory $\boldsymbol{x}_{d}\left(t\right)$. However,
as will be discussed later on, the functional as defined in Eq. \eqref{eq:JDef}
might lead to an ill-defined control signal $\boldsymbol{u}\left(t\right)$
and state trajectory $\boldsymbol{x}\left(t\right)$. A possible remedy
is to introduce a regularization term
\begin{align}
\mathcal{J}\left[\boldsymbol{x}\left(t\right),\boldsymbol{u}\left(t\right)\right] & =\frac{1}{2}\intop_{t_{0}}^{t_{1}}dt\left(\boldsymbol{x}\left(t\right)-\boldsymbol{x}_{d}\left(t\right)\right)^{2}+\dfrac{\epsilon^{2}}{2}\intop_{t_{0}}^{t_{1}}dt\left|\boldsymbol{u}\left(t\right)\right|^{2}.\label{eq:JDef_1}
\end{align}
The coefficient $\epsilon$ is a regularization parameter. The regularization
term penalizes large controls and guarantees a well-defined solution
to the optimization problem.

Mathematically speaking, the problem of finding the optimal control
signal $\boldsymbol{u}\left(t\right)$ by minimizing Eq. \eqref{eq:JDef_1}
is a constrained minimization problem, with $\boldsymbol{x}\left(t\right)$
 constrained to be the solution to the controlled dynamical system
\eqref{eq:AffineNonlinearDynamicalSystem} with initial condition
\eqref{eq:AffineNonlinearDynamicalSystemInitCond}. The standard approach
to solving constrained optimization problems is to introduce Lagrange
multipliers $\boldsymbol{\lambda}\left(t\right)$. The \textit{co-state}
$\boldsymbol{\lambda}\left(t\right)$ has the same number of components
as the state $\boldsymbol{x}\left(t\right)$. Its time evolution is
governed by the \textit{adjoint equation}. In contrast to solving
an uncontrolled problem, which only involves finding a solution for
the state $\boldsymbol{x}\left(t\right)$ in state space with dimension
$n$, solving an optimal control problem involves finding a solution
to the coupled state and adjoint equations in the \textit{extended
state space} with dimension $2n$. This renders optimal control problems
much more difficult than uncontrolled problems. Numerical solutions
of optimal control problems suffer from inconvenient terminal conditions
for the co-state and necessitate a computationally expensive iterative
algorithm.

The next example shows optimal trajectory tracking in the FHN model.
The solution is obtained numerically with the help of the open source
package \href{http://www.acadotoolkit.org}{ACADO} \cite{acadoManual,Houska2011a,Houska2011b}.

\begin{example}[Optimal trajectory tracking in the FHN model]\label{ex:ControlledFHN1}

Optimal trajectory tracking is discussed for the activator-controlled
FHN model of Example \ref{ex:FHN1}. The desired trajectory is chosen
to be an ellipse,
\begin{align}
x_{d}\left(t\right) & =A_{x}\cos\left(2\pi t/T\right)-\dfrac{1}{2}, & y_{d}\left(t\right) & =A_{y}\sin\left(2\pi t/T\right)+\dfrac{1}{2},\label{eq:DesiredTrajectories}
\end{align}
with $A_{x}=1,\,A_{y}=15$, and $T=1$. The regularization parameter
$\epsilon$ is set to the value $\epsilon=10^{-3}$, such that the
coefficient of the regularization term is $\sim10^{-6}$. Within the
time interval $0=t_{0}\leq t<t_{1}=1$, the controlled state shall
follow the ellipse as closely as possible. The initial and terminal
states lie exactly on the desired trajectory. Figure \ref{fig:ControlledFHN_11}
shows a numerical solution of the optimal trajectory tracking problem.
The controlled state trajectory $\boldsymbol{x}\left(t\right)$ (red
dashed line) looks quite different from the desired trajectory $\boldsymbol{x}_{d}\left(t\right)$
(blue solid line). However, the control is optimal, and not other
control would yield a controlled state trajectory closer to the desired
trajectory as measured by the functional Eq. \eqref{eq:JDef_1}.

\begin{minipage}{1.0\linewidth}
\begin{center}
\includegraphics[scale=0.55]{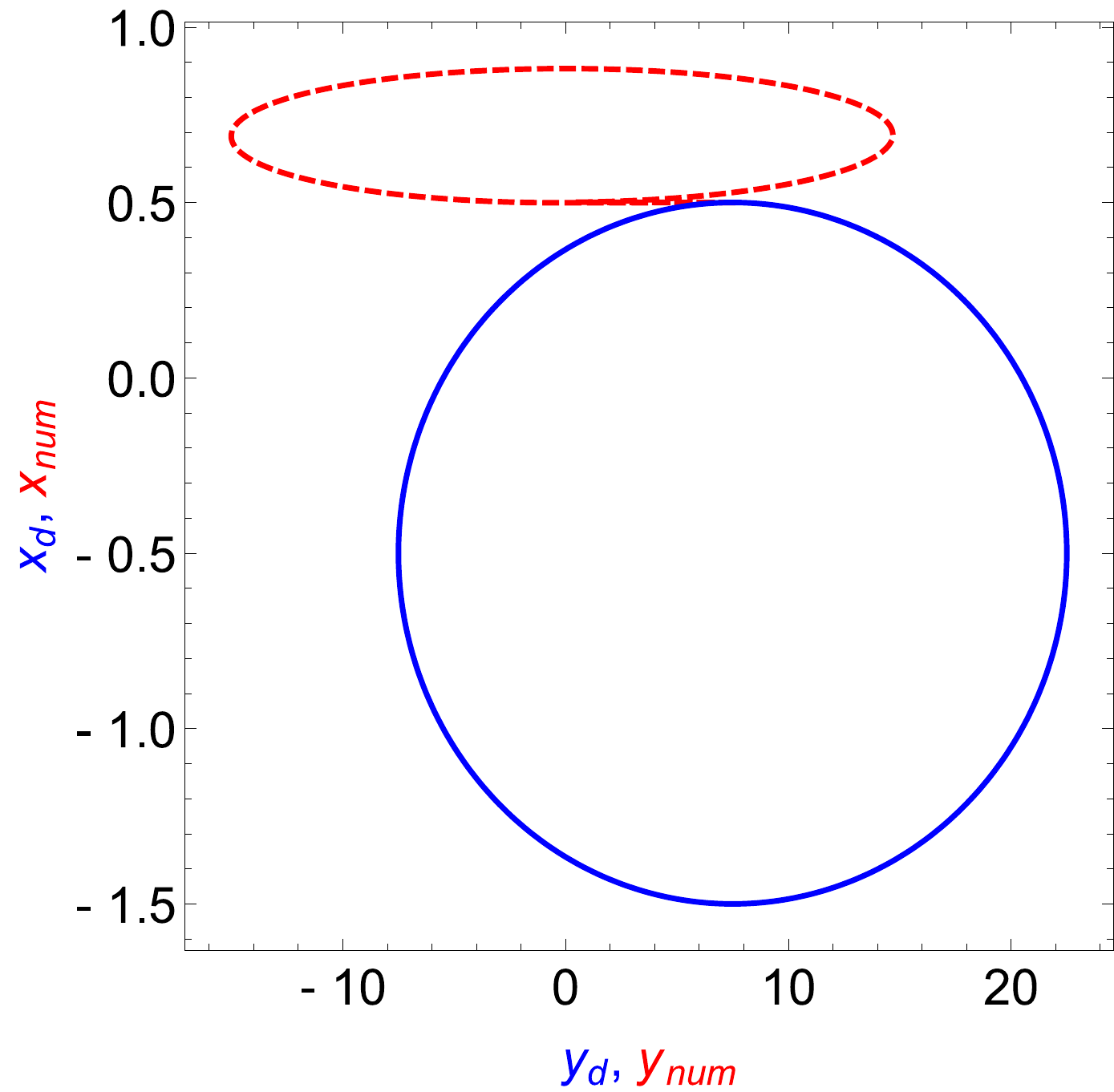}
\captionof{figure}[Optimal trajectory tracking in the  FHN model]{\label{fig:ControlledFHN_11}Optimal trajectory tracking in the activator-controlled FHN model. The state space plot compares the desired trajectory (blue solid line) with the optimally controlled state trajectory (red dashed line). The agreement is not particularly impressive. However, the control is optimal, and no better control exists. Any other control yields a state trajectory $\boldsymbol{x}\left(t\right)$ with a larger distance to the desired trajectory $\boldsymbol{x}_{d}\left(t\right)$ as measured by the functional $\mathcal{J}$ defined in Eq. \eqref{eq:JDef_1}.}
%"/home/jakob/ACADOtoolkit/examples/my_examples/FHN/CompareResults_FHN.nb"
\end{center}
\end{minipage}

Comparing the individual state components with is desired counterparts
in Fig. \ref{fig:ControlledFHN_11} somewhat clarifies the picture.
The controlled activator (red dashed line in Fig. \ref{fig:ControlledFHN_11}
right) is at least similar in shape to the desired activator (blue
solid line) but seems to be shifted by a constant. A very steep initial
transition leads from the initial condition onto the shifted trajectory,
while a similarly steep transition occurs at the terminal time. The
controlled inhibitor (red dashed line in Fig. \ref{fig:ControlledFHN_11}
left) does not show any similarity to the desired inhibitor (blue
solid line). In contrast to the activator component, it does not exhibit
steep initial and terminal transitions.

\begin{minipage}{1.0\linewidth}
\begin{center}
\includegraphics[scale=0.425]{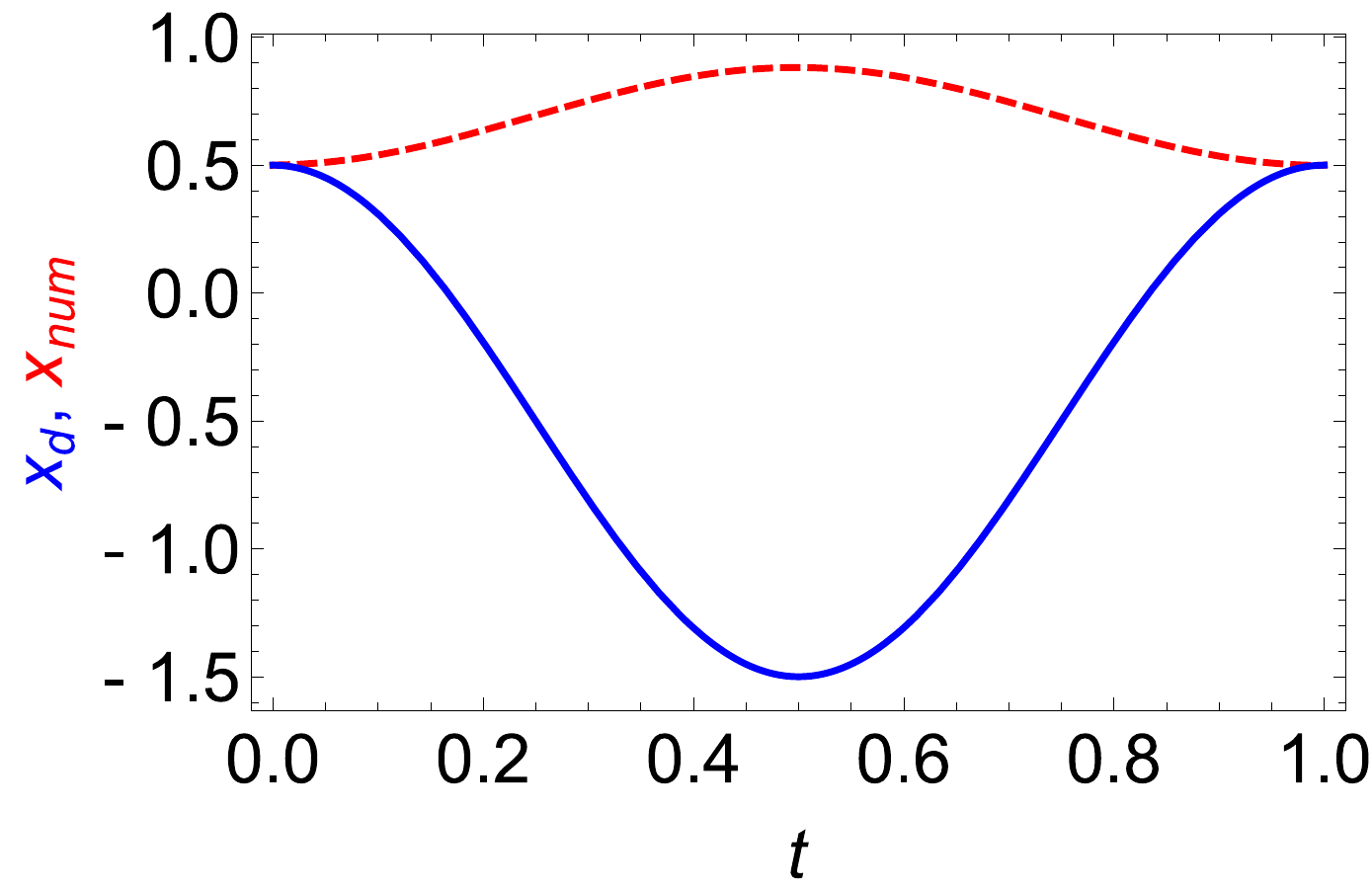}\hspace{0.2cm}\includegraphics[scale=0.425]{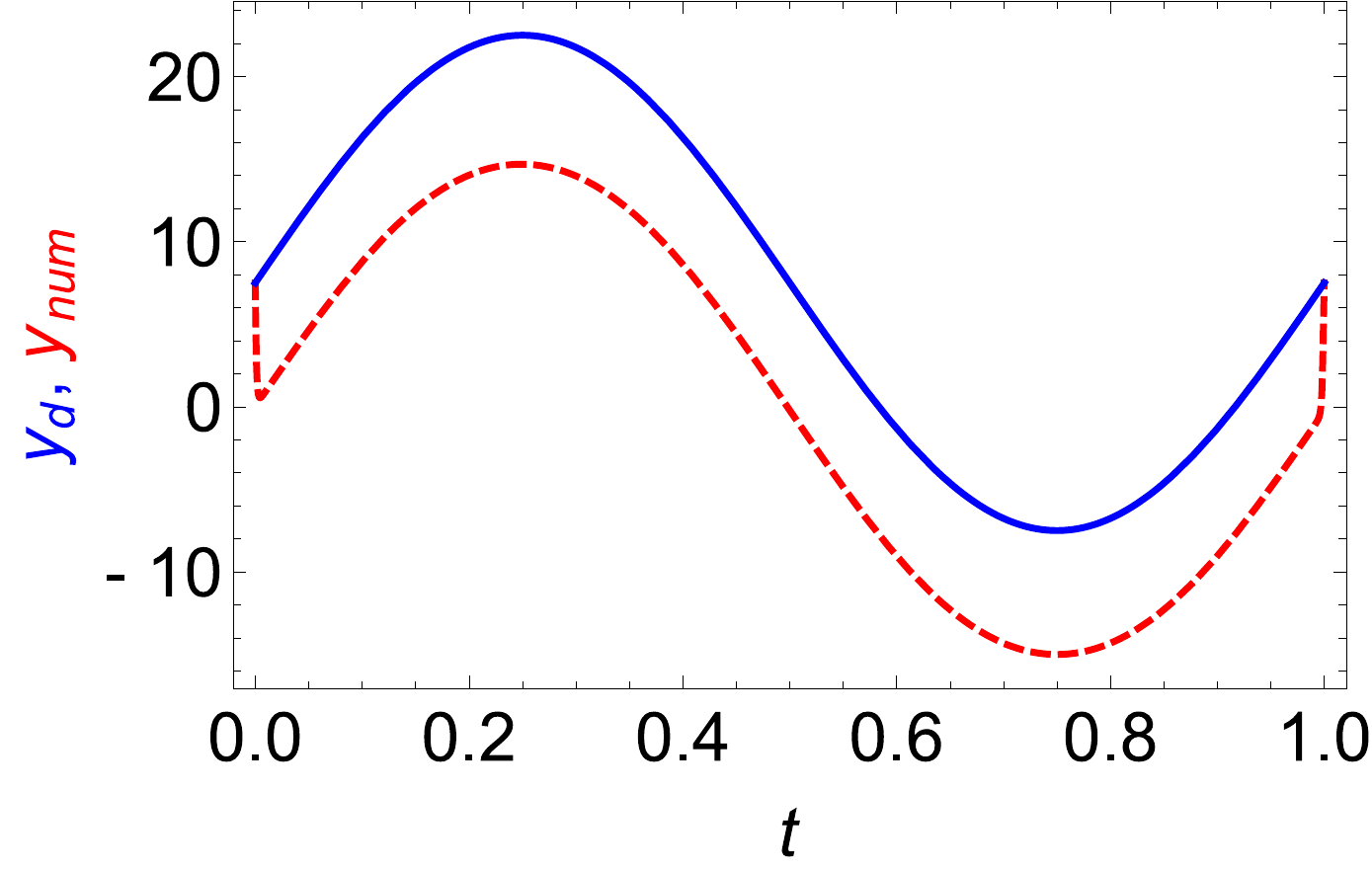}
\captionof{figure}[Optimal activator and inhibtor over time]{\label{fig:ControlledFHN_12}Plotting the individual state components over time reveals a controlled activator component $y$ similar in shape to $y_d$ but shifted by a constant. Note the steep transition regions at the beginning and end of the time interval. However, the controlled inhibitor component $x$ is very different from its desired counterpart $x_d$.}
%"/home/jakob/ACADOtoolkit/examples/my_examples/FHN/CompareResults_FHN.nb"
\end{center}
\end{minipage}

\end{example}

An important lesson is to be learned from Example \ref{ex:ControlledFHN1}.
In general, it is impossible for the controlled state trajectory $\boldsymbol{x}\left(t\right)$
to follow exactly the desired trajectory $\boldsymbol{x}_{d}\left(t\right)$.
The desired trajectory is that what you want, but it is usually not
that what you get. If the numerical solution corresponds to the global
minimum of the functional $\mathcal{J}$, Eq. \eqref{eq:JDef_1},
there is no other control which can enforce controlled state trajectory
$\boldsymbol{x}\left(t\right)$ closer to $\boldsymbol{x}_{d}\left(t\right)$.
Although the value of the functional $\mathcal{J}$ attains its minimally
possible value, this value might still be very large, indicating a
large distance between controlled and desired state trajectory. Naturally,
the following question arises. Under which conditions is the controlled
state trajectory $\boldsymbol{x}\left(t\right)$ identical to the
desired trajectory $\boldsymbol{x}_{d}\left(t\right)$? 

To answer that question, the concept of \textit{exactly realizable
trajectories} is proposed in Chapter \ref{chap:ExactlyRealizableTrajectories}.
A desired trajectory is exactly realizable if it satisfies a condition
called the \textit{constraint equation}. An open loop control signal
$\boldsymbol{\ensuremath{u}}\left(t\right)$ can be determined which
enforces the state to follow the desired trajectory exactly, $\boldsymbol{x}\left(t\right)=\boldsymbol{x}_{d}\left(t\right)$.
This approach does not only yield an explicit expression for the control
signal in terms of the desired trajectory, but also identifies a particularly
simple class of nonlinear affine control systems. Systems in this
class share many properties with linear control systems and satisfy
the so-called \textit{linearizing assumption}. Chapter \ref{chap:OptimalControl}
relates exactly realizable trajectories to optimal control. In particular,
an exactly realizable trajectories, together with the corresponding
control signal, is the solution to an unregularized optimal control
problem. Based on that insight, the regularization parameter $\epsilon$
is used as the small parameter for a singular perturbation expansion
in Chapter \ref{chap:AnalyticalApproximationsForOptimalTrajectoryTracking}.
This results in a reinterpretation of affine optimal trajectory tracking
problems with small regularization term as a system of singularly
perturbed differential equations. Combining this approach with the
linearizing assumption, approximate solutions for optimal trajectory
tracking in terms of mostly linear equations can be derived. The analytical
solutions are valid for arbitrary desired trajectories. This approach
applies, among other systems, to the mechanical control systems from
Example \ref{ex:FHN1} and the activator-controlled FHN model from
Example \ref{ex:ControlledFHN1}. Note that the small parameter $\epsilon$
originates from the formulation of the control problem Eq. \eqref{eq:JDef_1}.
Assuming this parameter to be small does \textit{not} involve any
approximations of the system dynamics. The system dynamics is exactly
taken into account by the perturbative approach. While the analytical
results are obtained for open loop control, they are modified in Section
\ref{sec:OptimalFeedback} to yield solutions for optimal feedback
control. Chapter \ref{chap:ControlOfReactionDiffusion} extends the
notion of exactly realizable trajectories to reaction-diffusion systems.

As a teaser and to demonstrate the accuracy of the analytical approximation,
we compare the numerical solution for the optimal trajectory tracking
in the FHN model from Example \ref{ex:ControlledFHN1} with the analytical
approximation in Example \ref{ex:ControlledFHN2}. The exact analytical
expression and its derivation is quite involved. All details can be
found in Chapter \ref{chap:AnalyticalApproximationsForOptimalTrajectoryTracking}.

\begin{example}[Analytical approximation for the optimally controlled FHN model]\label{ex:ControlledFHN2}

Optimal trajectory tracking in the activator-controlled FHN model
(see Example \ref{ex:FHN1}) can be approximately solved with the
analytical techniques developed in this thesis. The regularization
parameter $\epsilon$ in Eq. \eqref{eq:JDef_1} is used as the small
parameter for a singular perturbation expansion. For the same desired
trajectories as in Example \ref{ex:ControlledFHN1}, and the same
value of the regularization parameter $\epsilon=10^{-3}$, Fig. \ref{fig:ControlledFHN_2}
compares the analytical approximation for the optimally controlled
state trajectory $\boldsymbol{x}\left(t\right)$ with the corresponding
numerical solution. For such a small value of the regularization parameter
$\epsilon=10^{-3}$, the agreement is almost perfect.

\begin{minipage}{1.0\linewidth}
\begin{center}
\includegraphics[scale=0.5]{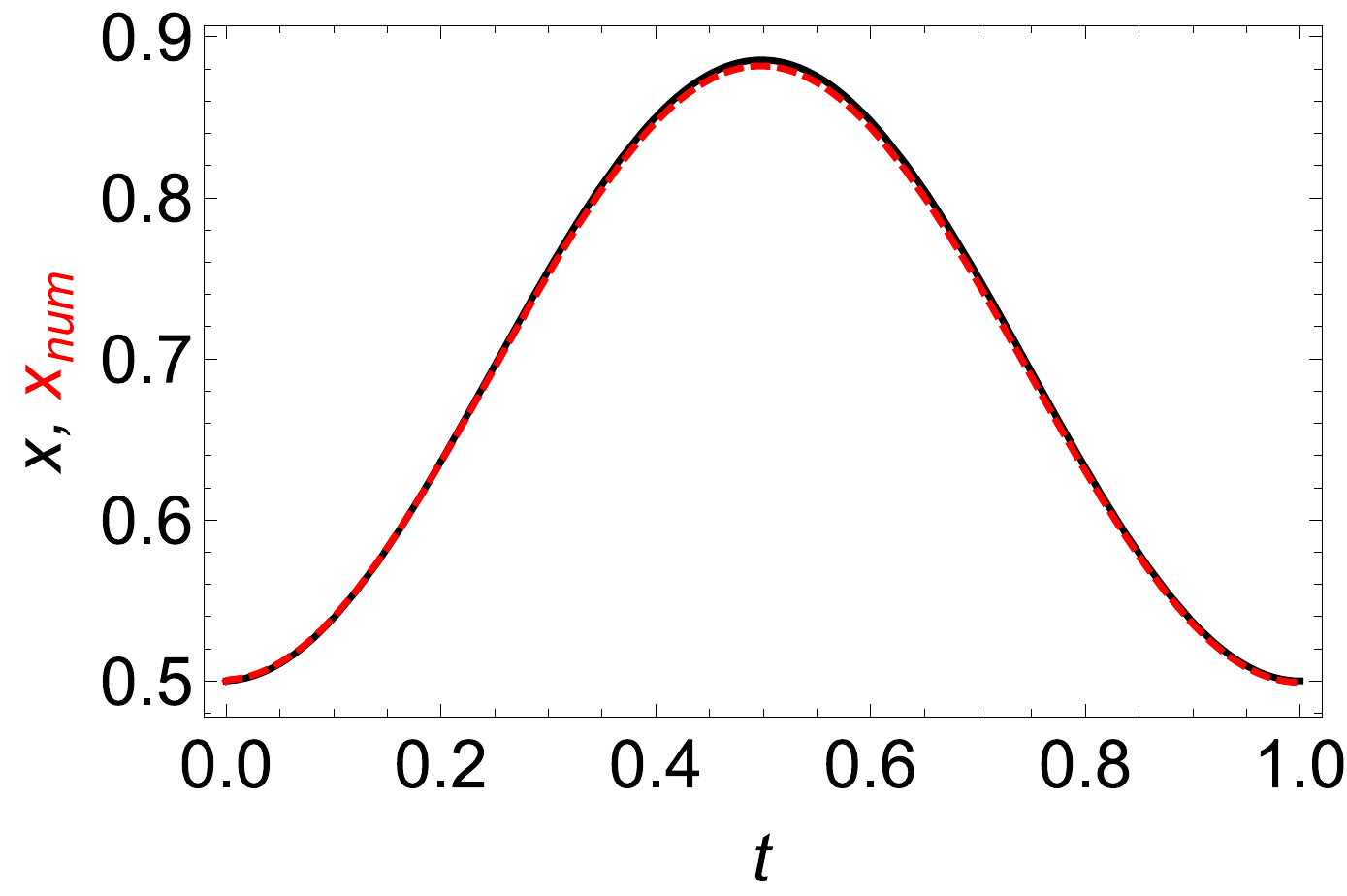}\includegraphics[scale=0.5]{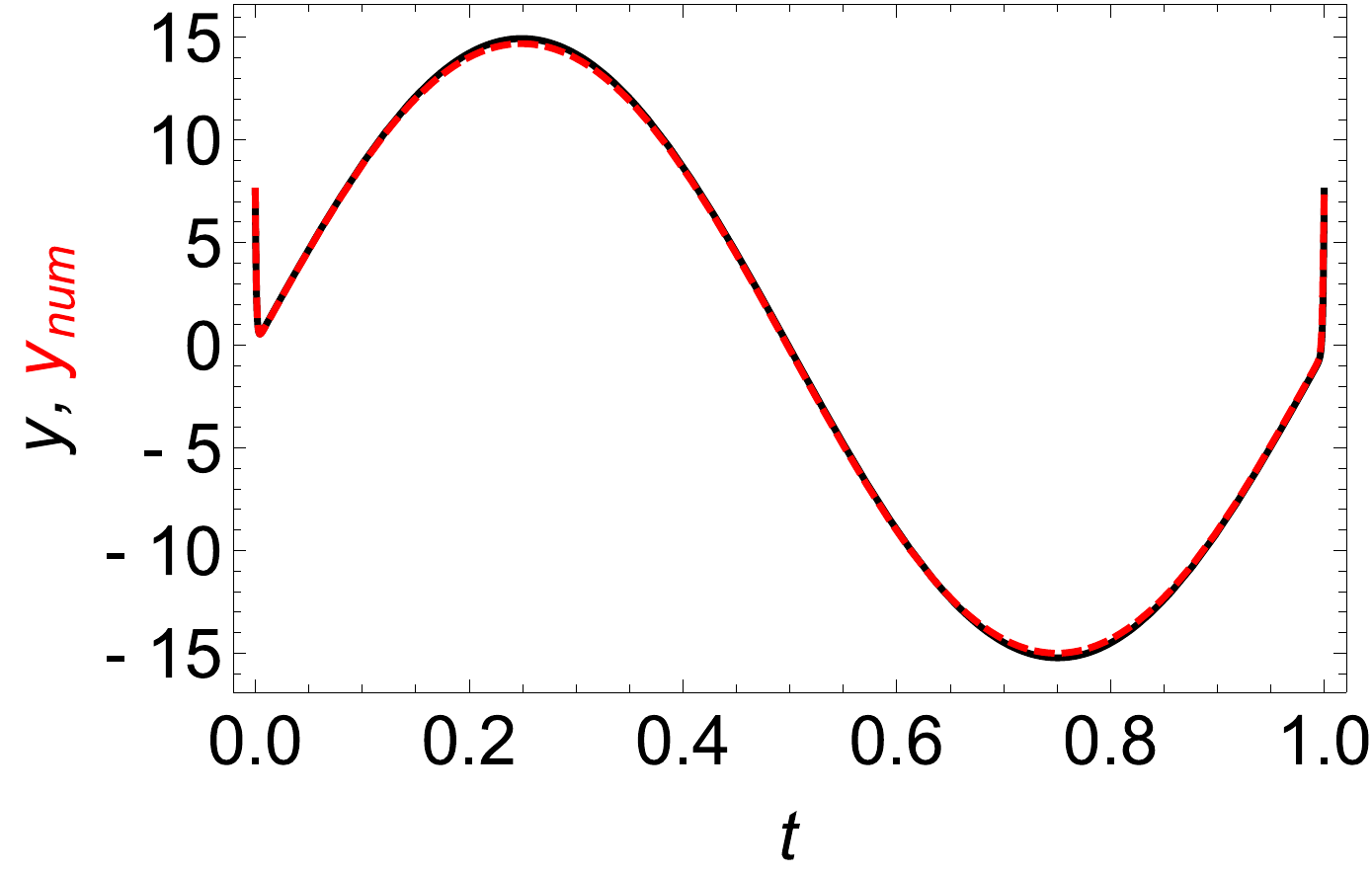}
\captionof{figure}[Analytical and numerical result for the controlled state trajectory]{\label{fig:ControlledFHN_2}Comparison of analytical approximation (black solid line) and numerically obtained (red dashed line) optimally controlled state trajectory of Example \ref{ex:ControlledFHN1} for the activator $y$ (right) and inhibitor $x$ (left) over time.}
%"/home/jakob/ACADOtoolkit/examples/my_examples/FHN/CompareResults_FHN.nb"
\end{center}
\end{minipage}

\end{example}

The analytical result of Chapter \ref{chap:AnalyticalApproximationsForOptimalTrajectoryTracking}
reveals a surprising result. The analytical approximation for the
controlled state trajectory $\boldsymbol{x}\left(t\right)$ of Example
\ref{ex:ControlledFHN1} does not depend on the nonlinearity $R\left(x,y\right)$
(see Example \ref{ex:FHN1} for the model equations)! More precisely,
changing the parameter values of the nonlinearity $R\left(x,y\right)$,
or changing $R\left(x,y\right)$ altogether, has no effect on the
controlled state trajectory. Although the system dynamics is governed
by nonlinear differential equations, the optimally controlled system
can be approximated by solving only linear equations. However, the
analytical solution for the control signal depends strongly on the
nonlinearity $R$. This prediction is verified with an additional
numerical computation for a FHN like model with vanishing nonlinearity
$R\equiv0$ in Example \ref{ex:ControlledFHN3}.

\begin{example}[Optimal trajectory tracking for a linear system]\label{ex:ControlledFHN3}

The affine control system
\begin{align}
\left(\begin{array}{c}
\dot{x}\left(t\right)\\
\dot{y}\left(t\right)
\end{array}\right) & =\left(\begin{array}{c}
a_{0}+a_{1}x\left(t\right)+a_{2}y\left(t\right)\\
0
\end{array}\right)+\left(\begin{array}{c}
0\\
1
\end{array}\right)u\left(t\right),\label{eq:FHNLinear}
\end{align}
has the same form as the FHN model of Example \ref{ex:FHN1} except
for the vanishing nonlinearity,
\begin{align}
R\left(x,y\right) & =0.
\end{align}
The parameter values for $a_{0},\,a_{1},$ and $a_{2}$ are the same
as in Example \ref{ex:FHN1}. The analytical solution predicts that,
in the limit of small regularization parameter $\epsilon\rightarrow0$,
the optimally controlled state trajectory is independent of the actual
form of the nonlinearity $R$. Therefore, choosing the same desired
trajectories Eq. \eqref{eq:DesiredTrajectories} for the optimal control
of Eq. \eqref{eq:FHNLinear} should yield the same optimally controlled
state trajectories $\boldsymbol{x}\left(t\right)$ as in Example \ref{ex:ControlledFHN1}.
Indeed, no discernible difference is visible in the numerical solutions
for both problems, as Fig. \ref{fig:ControlledFHN_31} shows. However,
the corresponding control signals depend on the nonlinearity, as is
shown in Fig. \ref{fig:ControlledFHN_32}.

\begin{minipage}{1.0\linewidth}
\begin{center}
\includegraphics[scale=0.49]{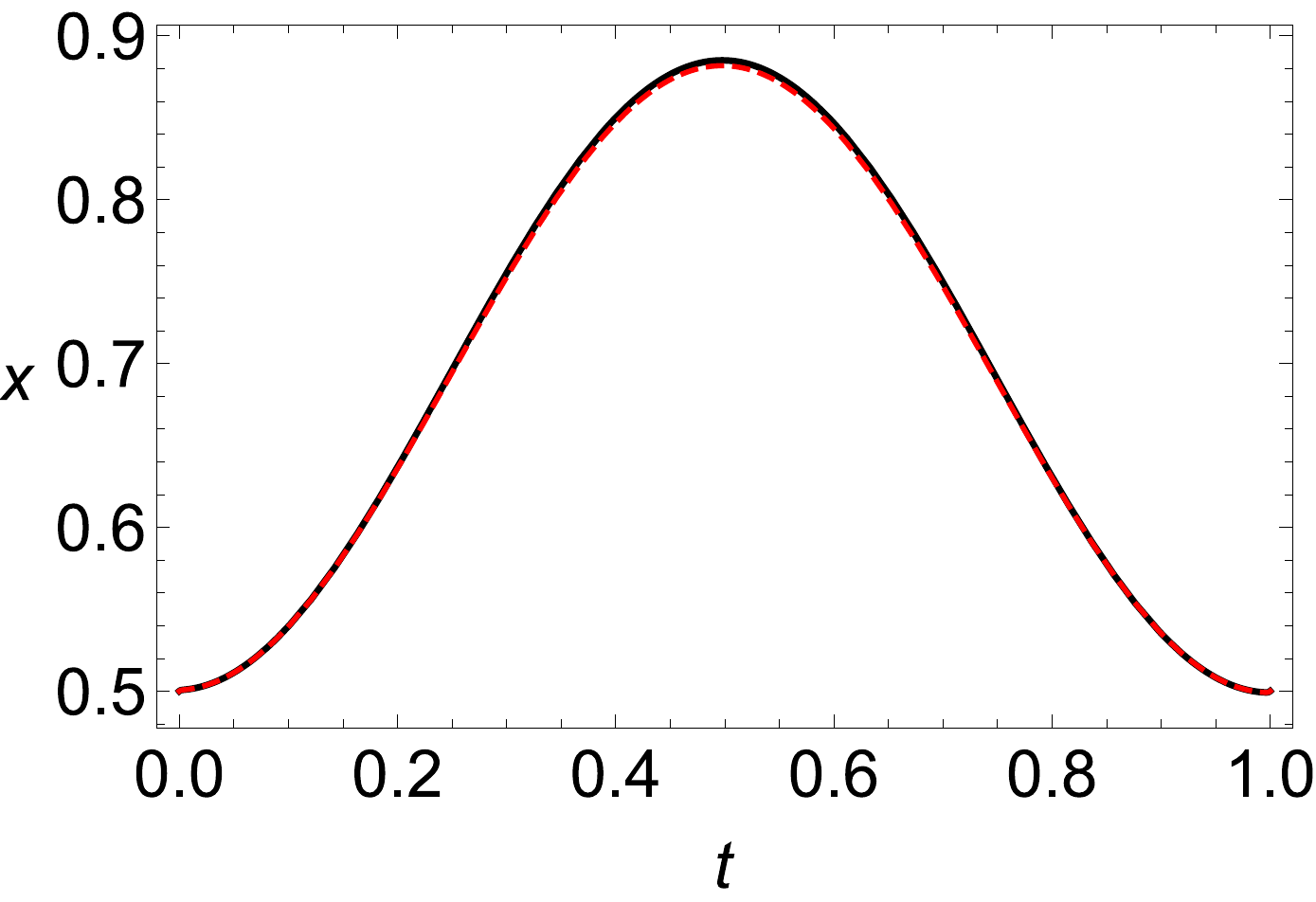}\includegraphics[scale=0.5]{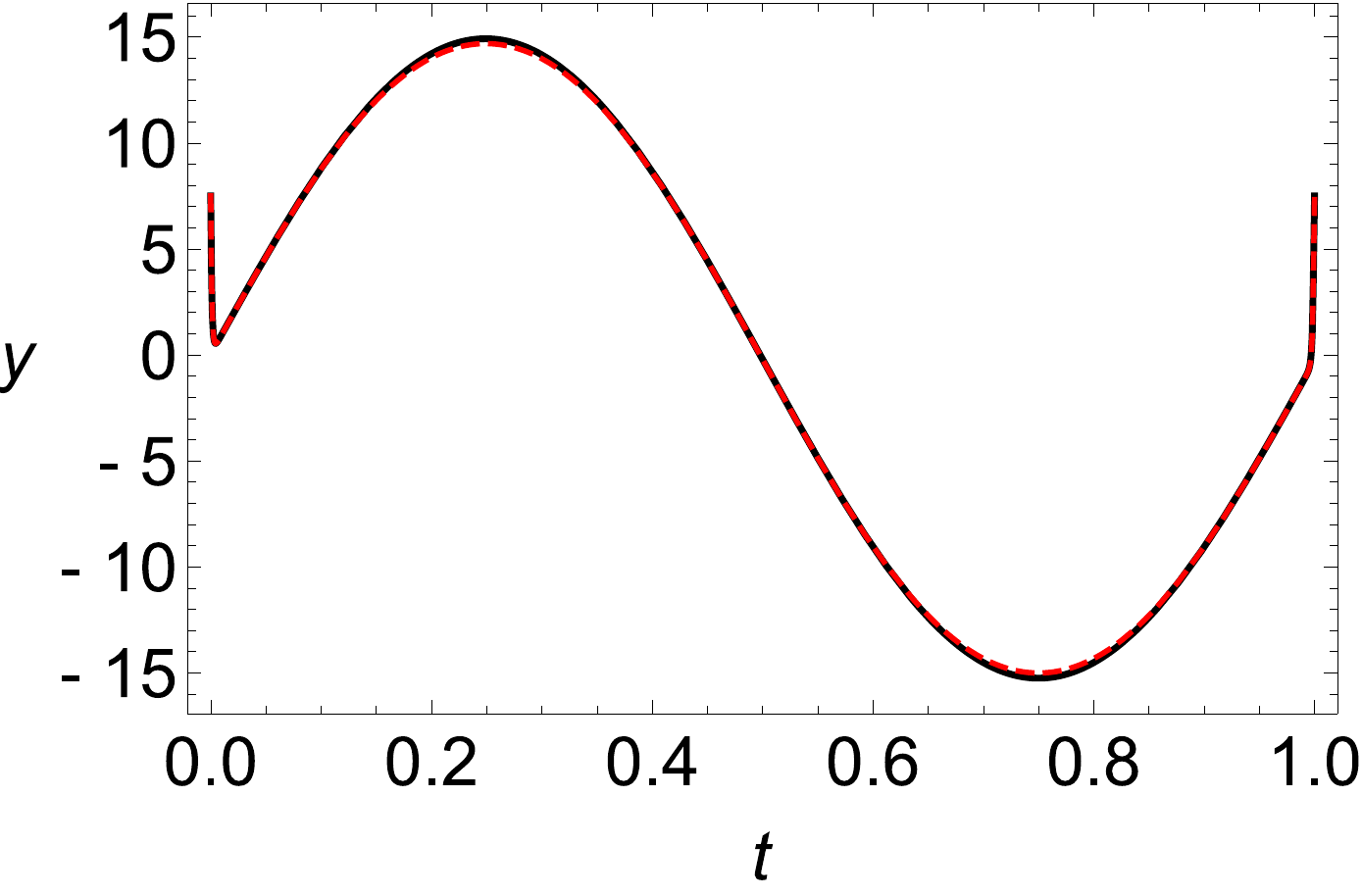}
\captionof{figure}[The controlled state does not depend on the nonlinearity]{\label{fig:ControlledFHN_31} The optimally controlled state trajectory $\boldsymbol{x}=\left( x, y \right)^T$ does not depend on the nonlinearity $R$ in the limit of vanishing regularization parameter, $\epsilon\rightarrow 0$. The black solid line is the numerical result for a vanishing nonlinearity $R\left(x,y\right)=0$, while the red dashed line shows the result for the FHN nonlinearity $R\left(x,y\right)=y-\frac{1}{3}y^{3}-x$.}
%"/home/jakob/ACADOtoolkit/examples/my_examples/FHN/CompareResults_FHN.nb"
\end{center}
\end{minipage}

\begin{minipage}{1.0\linewidth}
\begin{center}
\includegraphics[scale=0.5]{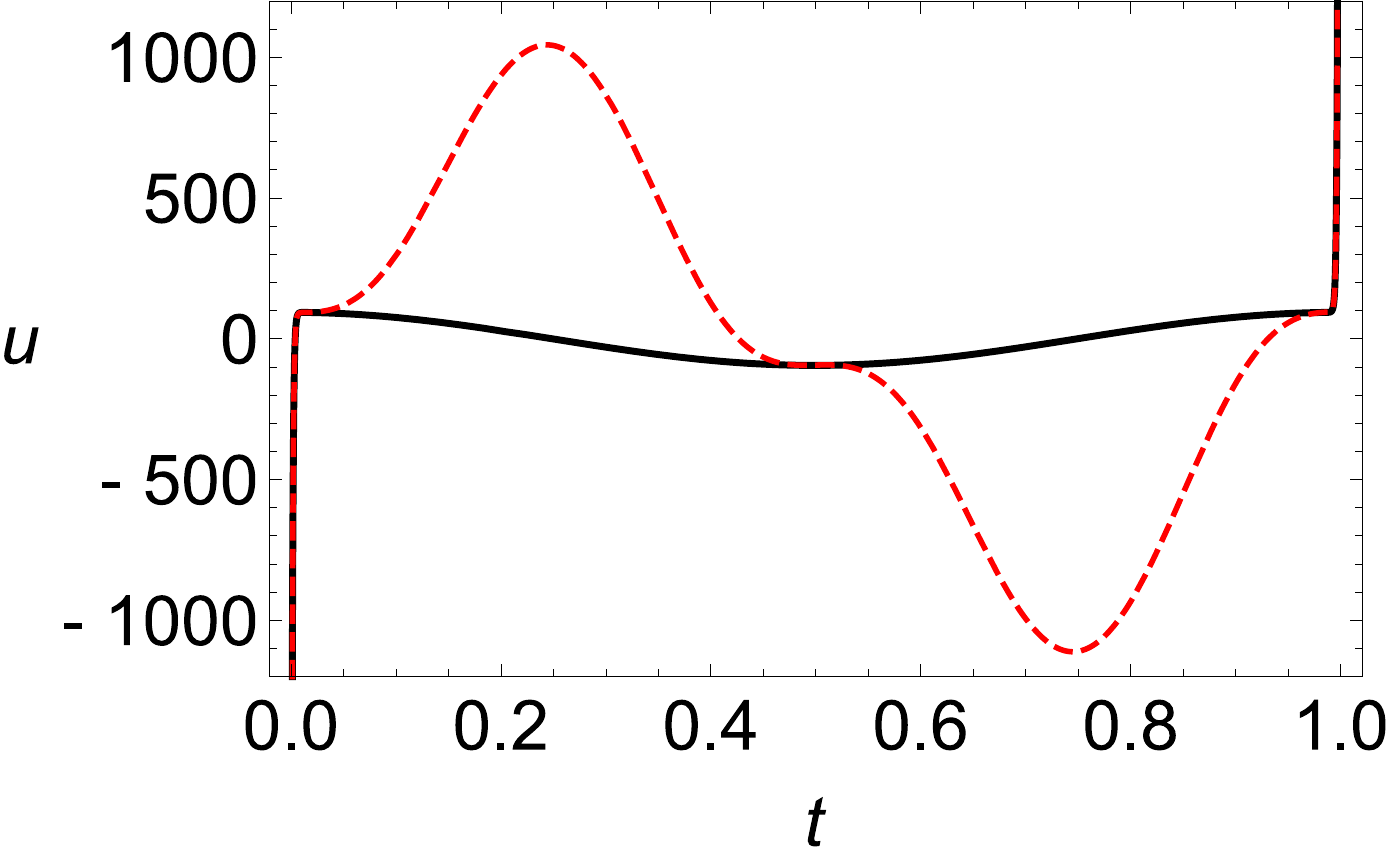}
\captionof{figure}[The control signal depends heavily on the nonlinearity]{\label{fig:ControlledFHN_32} The control signal $u$ depends on the nonlinearity $R$. Black solid line is the numerical result for vanishing nonlinearity $R\left(x,y\right)=0$, while the red dashed line shows the result for the standard FHN nonlinearity $R\left(x,y\right)=y-\frac{1}{3}y^{3}-x$. Note that the control exhibits very large values and steep slopes at the beginning and end of the time interval in both cases.}
%"/home/jakob/ACADOtoolkit/examples/my_examples/FHN/CompareResults_FHN.nb"
\end{center}
\end{minipage}

\end{example}

The last example reveals that, under certain conditions, the nonlinearity
$R\left(x,y\right)$ plays only a minor role for the controlled state
trajectory. In fact, if the regularization parameter $\epsilon$ is
zero, the approximate solution for the state trajectory is exact and
governed solely by linear equations. The analytical treatment uncovers
an underlying linear structure of certain nonlinear optimal trajectory
tracking problems. However, this exact solution to nonlinear optimal
control systems does not come without a price. For $\epsilon=0$,
the optimally controlled state trajectory cannot be expressed in terms
of continuous functions, but involves jumps in at least one component.
These jumps are located at the beginning and the end of the time interval.
Even worse, the corresponding control signal does diverge at the same
points. This behavior can already be anticipated from Fig. \ref{fig:ControlledFHN_31}
right: the $y$-component of the state exhibits a steep transition
region close to the beginning and the end of the time interval for
a small but finite value of $\epsilon=10^{-3}$. These transition
regions degenerate to jumps in the limit $\epsilon\rightarrow0$.
In the language of singular perturbation theory, such transition regions
are known as boundary layers. Similarly, Fig. \ref{fig:ControlledFHN_32}
shows that the corresponding control signal $u\left(t\right)$ exhibits
very large amplitudes located equally at the beginning and the end
of the time interval.

An underlying linear structure of nonlinear control systems might
sound surprising. However, instances of exact linearizations are well
known from mathematical control theory. A prominent example is feedback
linearization, see e.g. \cite{khalil2002nonlinear,slotine1991applied,isidori1995nonlinear}.
As a very simple example system, consider the activator-controlled
FHN model,
\begin{align}
\left(\begin{array}{c}
\dot{x}\left(t\right)\\
\dot{y}\left(t\right)
\end{array}\right) & =\left(\begin{array}{c}
a_{0}+a_{1}x\left(t\right)+a_{2}y\left(t\right)\\
R\left(x\left(t\right),y\left(t\right)\right)
\end{array}\right)+\left(\begin{array}{c}
0\\
1
\end{array}\right)u\left(t\right),\label{eq:FHNfeedbackLinear}
\end{align}
with nonlinearity
\begin{align}
R\left(x,y\right) & =y-\frac{1}{3}y^{3}-x.
\end{align}
As the name implies, feedback linearization assumes a feedback control
from the very beginning, i.e., $u$ may depend on state $\boldsymbol{x}$
as
\begin{align}
u\left(t\right) & =u\left(\boldsymbol{x}\left(t\right),t\right).
\end{align}
With the help of a very simple transform of the control signal, it
is possible to obtain a controlled system linear in state and control.
Introducing a new control signal $v\left(t\right)$ as
\begin{align}
u\left(x\left(t\right),y\left(t\right),t\right) & =-R\left(x\left(t\right),y\left(t\right)\right)+v\left(t\right),\label{eq:FeedbackLinearizeControl}
\end{align}
the controlled system transforms to
\begin{align}
\dot{x}\left(t\right) & =y\left(t\right), & \dot{y}\left(t\right) & =v\left(t\right).\label{eq:FHNfeedbackLinearize}
\end{align}
While the original control signal $u\left(t\right)$ depends on the
nonlinearity $R$, the controlled state trajectory $\boldsymbol{x}\left(t\right)$,
obtained as the solution to Eq. \eqref{eq:FHNfeedbackLinearize},
does not depend on $R$. This is in fact similar to our approach. 

In contrast to approximate linearizations performed to study the linear
stability of solutions, feedback linearization is an exact transformation
of a nonlinear to a linear dynamical system. Exact linearizations
of uncontrolled dynamical systems exist as well. For example, a nonlinear
transformation of the state converts Riccati equations to linear differential
equations \cite{Zaitsev2002Handbook}. However, because the class
of exactly linearizable uncontrolled systems is small, this method
is rarely applied in practice. In general, feedback linearization
applies a combined transformation of state and control to obtain a
linear system. The class of feedback linearizable nonlinear control
systems is huge, and the simple model Eq. \eqref{eq:FHNfeedbackLinear}
is only the trivial case requiring no state transformation \cite{khalil2002nonlinear}.

A disadvantage of feedback linearization is that it assumes a feedback
control from the very beginning and does not easily apply to open
loop control methods. Furthermore, feedback control might lead to
results which are not optimal. In principle, a feedback-controlled
nonlinear system can be much simpler than its corresponding uncontrolled
counterpart. This is in contrast to optimally controlled systems,
which are usually much more difficult than its uncontrolled counterpart
due to the coupling of the state and co-state equations. As demonstrated
by the example  above, in this thesis we develop analytical techniques
which reveal an underlying linear structure for a certain class of
nonlinear optimal control systems. Similar to feedback linearization,
the nonlinearity is absorbed by the control signal, and the time evolution
of the controlled state trajectory is entirely determined by linear
equations. These techniques apply only to a limited class of nonlinear
control systems, and are not as powerful as feedback linearization.
Nevertheless, this class includes some simple but important models
motivated by physics and nonlinear dynamics, as the activator-controlled
FHN model and mechanical control systems in one spatial dimension.

The chapter is concluded with a more philosophical remark. Trajectory
tracking is actually ill-defined because it is a circular task. To
achieve the aim of trajectory tracking, an appropriate control signal
must be applied to the dynamical system. In a universe which consists
exclusively of dynamical systems, this control signal must be the
output of a dynamical system. The only way to obtain an output which
behaves in exactly the way necessary for trajectory tracking is to
control the dynamical system which generates the output. To perform
the task of trajectory tracking, it is necessary to have a second
system for which the task of trajectory tracking is already performed
with sufficient accuracy. Trajectory tracking is a circular task.

%% file: chapter-2.tex
\lhead[\chaptername~\thechapter\leftmark]{}

\rhead[]{\rightmark}

\lfoot[\thepage]{}

\cfoot{}

\rfoot[]{\thepage}

\chapter{\label{chap:ExactlyRealizableTrajectories}Exactly realizable trajectories}

This chapter introduces the notion of exactly realizable trajectories.
The necessary formalism is established in Section \ref{sec:Formalism}.
After the definition of exactly realizable trajectories in Section
\ref{sec:ExactlyRealizableTrajectories}, the linearizing assumption
is introduced in Section \ref{sec:LinearizingAssumption}. This assumption
defines a class of nonlinear control systems which, to a large extent,
behave like linear control systems. Combining the notion of an exactly
realizable trajectory with the linearizing assumption allows one to
extend some well known results about the controllability of linear
systems to nonlinear control systems in Sections \ref{sec:Controllability}
and \ref{sec:OutputControllability}. Output Realizability is discussed
in Section \ref{sec:OutputRealizability}, and Section \ref{sec:2Conclusions}
concludes with a discussion and outlook.

\section{\label{sec:Formalism}Formalism}

This section introduces the formalism which is repeatedly used throughout
the thesis. The main elements are two complementary projection matrices
$\boldsymbol{\mathcal{P}}$ and $\boldsymbol{\mathcal{Q}}$. Projectors
are a useful ingredient for a number of physical theories. Take, for
example, quantum mechanics, which describes measurements as projections
of the state (an element from a Hilbert space) onto a ray or unions
of rays of the Hilbert space \cite{fick1988einfuhrung,cohen2007quantenmechanik}.
Also in non-equilibrium statistical mechanics, projectors have found
widespread application to separate a subsystem of interest from its
bath \cite{balescu1975equilibrium,grabert1982projection}.

To the best of our knowledge, projectors have not been utilized in
the context of control systems. Section \ref{sub:TwoUsefulProjectors}
defines the projectors and Section \ref{sub:SeparationOfTheStateEq}
separates the controlled state equation in two equations. The first
equation involves the control signal, while the second equation is
independent of the control signal. This formalism provides a useful
approach for analyzing general affine control systems. Appendix \ref{sec:OverAndUnderdetSysOfEqs}
demonstrates how projectors arise in the context of overdetermined
and underdetermined systems of linear equations.

\subsection{\label{sub:TwoUsefulProjectors}The projectors \texorpdfstring{$\boldsymbol{\mathcal{P}}$}{P}
and \texorpdfstring{$\boldsymbol{\mathcal{Q}}$}{Q}}

Consider the affine control system with state dependent coupling matrix
$\boldsymbol{\mathcal{B}}\left(\boldsymbol{x}\left(t\right)\right)$
\begin{align}
\boldsymbol{\dot{x}}\left(t\right) & =\boldsymbol{R}\left(\boldsymbol{x}\left(t\right)\right)+\boldsymbol{\mathcal{B}}\left(\boldsymbol{x}\left(t\right)\right)\boldsymbol{u}\left(t\right).
\end{align}
Define two complementary projectors $\boldsymbol{\mathcal{P}}$ and
$\boldsymbol{\mathcal{Q}}$ in terms of the coupling matrix $\boldsymbol{\mathcal{B}}\left(\boldsymbol{x}\left(t\right)\right)$
as 
\begin{align}
\boldsymbol{\mathcal{P}}\left(\boldsymbol{x}\right) & =\boldsymbol{\mathcal{B}}\left(\boldsymbol{x}\right)\left(\boldsymbol{\mathcal{B}}^{T}\left(\boldsymbol{x}\right)\boldsymbol{\mathcal{B}}\left(\boldsymbol{x}\right)\right)^{-1}\boldsymbol{\mathcal{B}}^{T}\left(\boldsymbol{x}\right),\label{eq:TimeDependentP}\\
\boldsymbol{\mathcal{Q}}\left(\boldsymbol{x}\right) & =\boldsymbol{1}-\boldsymbol{\mathcal{P}}\left(\boldsymbol{x}\right).\label{eq:TimeDependentQ}
\end{align}
$\boldsymbol{\mathcal{P}}$ and $\boldsymbol{\mathcal{Q}}$ are $n\times n$
matrices which, in general, do depend on the state $\boldsymbol{x}$.
Note that the $p\times p$ matrix $\boldsymbol{\mathcal{B}}^{T}\left(\boldsymbol{x}\right)\boldsymbol{\mathcal{B}}\left(\boldsymbol{x}\right)$
has full rank $p$ because of assumption Eq. \eqref{eq:BFullRank}
that $\boldsymbol{\mathcal{B}}\left(\boldsymbol{x}\right)$ has full
rank. Therefore, $\boldsymbol{\mathcal{B}}^{T}\left(\boldsymbol{x}\right)\boldsymbol{\mathcal{B}}\left(\boldsymbol{x}\right)$
is a quadratic and non-singular matrix and its inverse exists. The
projectors $\boldsymbol{\mathcal{P}}$ and $\boldsymbol{\mathcal{Q}}$
are also known as \textit{Moore-Penrose projectors}. The rank of $\boldsymbol{\mathcal{P}}\left(\boldsymbol{x}\right)$
and $\boldsymbol{\mathcal{Q}}\left(\boldsymbol{x}\right)$ is
\begin{align}
\text{rank}\left(\boldsymbol{\mathcal{P}}\left(\boldsymbol{x}\right)\right) & =p, & \text{rank}\left(\boldsymbol{\mathcal{Q}}\left(\boldsymbol{x}\right)\right) & =n-p.
\end{align}
Multiplying the $n$-component state vector $\boldsymbol{x}$ by the
$n\times n$ matrix $\boldsymbol{\mathcal{P}}\left(\boldsymbol{x}\right)$
yields an $n$-component vector $\boldsymbol{z}=\boldsymbol{\mathcal{P}}\left(\boldsymbol{x}\right)\boldsymbol{x}$.
However, because $\boldsymbol{\mathcal{P}}\left(\boldsymbol{x}\right)$
has rank $p$, only $p$ components of $\boldsymbol{z}$ are independent.
Similar, only $n-p$ components of $\boldsymbol{y}=\boldsymbol{\mathcal{Q}}\left(\boldsymbol{x}\right)\boldsymbol{x}$
are independent.

From the definitions Eqs. \eqref{eq:TimeDependentP} and \eqref{eq:TimeDependentQ}
follow the projector properties idempotence
\begin{align}
\boldsymbol{\mathcal{Q}}\left(\boldsymbol{x}\right)\boldsymbol{\mathcal{Q}}\left(\boldsymbol{x}\right) & =\boldsymbol{\mathcal{Q}}\left(\boldsymbol{x}\right), & \boldsymbol{\mathcal{P}}\left(\boldsymbol{x}\right)\boldsymbol{\mathcal{P}}\left(\boldsymbol{x}\right) & =\boldsymbol{\mathcal{P}}\left(\boldsymbol{x}\right),
\end{align}
and complementarity
\begin{align}
\boldsymbol{\mathcal{Q}}\left(\boldsymbol{x}\right)\boldsymbol{\mathcal{P}}\left(\boldsymbol{x}\right) & =\boldsymbol{\mathcal{P}}\left(\boldsymbol{x}\right)\boldsymbol{\mathcal{Q}}\left(\boldsymbol{x}\right)=\boldsymbol{0}.\label{eq:PTimesQ}
\end{align}
The projectors are symmetric,
\begin{align}
\boldsymbol{\mathcal{P}}^{T}\left(\boldsymbol{x}\right) & =\boldsymbol{\mathcal{P}}\left(\boldsymbol{x}\right), & \boldsymbol{\mathcal{Q}}^{T}\left(\boldsymbol{x}\right) & =\boldsymbol{\mathcal{Q}}\left(\boldsymbol{x}\right),
\end{align}
because the inverse of the symmetric matrix $\boldsymbol{\mathcal{B}}^{T}\left(\boldsymbol{x}\right)\boldsymbol{\mathcal{B}}\left(\boldsymbol{x}\right)$
is symmetric. Furthermore, matrix multiplication from the right with
the input matrix $\boldsymbol{\mathcal{B}}\left(\boldsymbol{x}\right)$
yields the important relations 
\begin{align}
\boldsymbol{\mathcal{P}}\left(\boldsymbol{x}\right)\boldsymbol{\mathcal{B}}\left(\boldsymbol{x}\right) & =\boldsymbol{\mathcal{B}}\left(\boldsymbol{x}\right), & \boldsymbol{\mathcal{Q}}\left(\boldsymbol{x}\right)\boldsymbol{\mathcal{B}}\left(\boldsymbol{x}\right) & =\boldsymbol{0}.\label{eq:PTimesB}
\end{align}
Similarly, matrix multiplication from the left with the transposed
input matrix $\boldsymbol{\mathcal{B}}^{T}\left(\boldsymbol{x}\right)$
yields
\begin{align}
\boldsymbol{\mathcal{B}}^{T}\left(\boldsymbol{x}\right)\boldsymbol{\mathcal{P}}\left(\boldsymbol{x}\right) & =\boldsymbol{\mathcal{B}}^{T}\left(\boldsymbol{x}\right), & \boldsymbol{\mathcal{B}}^{T}\left(\boldsymbol{x}\right)\boldsymbol{\mathcal{Q}}\left(\boldsymbol{x}\right) & =\boldsymbol{0}.
\end{align}
Some more properties of $\boldsymbol{\mathcal{P}}$ and $\boldsymbol{\mathcal{Q}}$
necessary for later chapters are compiled in Appendix \ref{sec:PropertiesOfTimeDependentProjectors}.

\subsection{\label{sub:SeparationOfTheStateEq}Separation of the state equation}

The projectors defined in Eqs. \eqref{eq:TimeDependentP} and \eqref{eq:TimeDependentQ}
are used to split up the controlled state equation
\begin{align}
\boldsymbol{\dot{x}}\left(t\right)= & \boldsymbol{R}\left(\boldsymbol{x}\left(t\right)\right)+\boldsymbol{\mathcal{B}}\left(\boldsymbol{x}\left(t\right)\right)\boldsymbol{u}\left(t\right).\label{eq:StateEquation}
\end{align}
Multiplying every term by $\boldsymbol{1}=\boldsymbol{\mathcal{P}}\left(\boldsymbol{x}\left(t\right)\right)+\boldsymbol{\mathcal{Q}}\left(\boldsymbol{x}\left(t\right)\right)$,
Eq. \eqref{eq:StateEquation} can be written as
\begin{align}
\frac{d}{dt}\left(\boldsymbol{\mathcal{P}}\left(\boldsymbol{x}\left(t\right)\right)\boldsymbol{x}\left(t\right)+\boldsymbol{\mathcal{Q}}\left(\boldsymbol{x}\left(t\right)\right)\boldsymbol{x}\left(t\right)\right) & =\left(\boldsymbol{\mathcal{P}}\left(\boldsymbol{x}\left(t\right)\right)+\boldsymbol{\mathcal{Q}}\left(\boldsymbol{x}\left(t\right)\right)\right)\boldsymbol{R}\left(\boldsymbol{x}\left(t\right)\right)\nonumber \\
 & +\left(\boldsymbol{\mathcal{P}}\left(\boldsymbol{x}\left(t\right)\right)+\boldsymbol{\mathcal{Q}}\left(\boldsymbol{x}\left(t\right)\right)\right)\boldsymbol{\mathcal{B}}\left(\boldsymbol{x}\left(t\right)\right)\boldsymbol{u}\left(t\right).
\end{align}
Multiplying with $\boldsymbol{\mathcal{Q}}\left(\boldsymbol{x}\left(t\right)\right)$
from the left and using Eq. \eqref{eq:PTimesB} yields an equation
independent of the control signal $\boldsymbol{u}$, 
\begin{align}
\boldsymbol{\mathcal{Q}}\left(\boldsymbol{x}\left(t\right)\right)\left(\boldsymbol{\dot{x}}\left(t\right)-\boldsymbol{R}\left(\boldsymbol{x}\left(t\right)\right)\right) & =\mathbf{0}.\label{eq:ConstraintEquation}
\end{align}
Equation \eqref{eq:ConstraintEq} is called the \textit{constraint
equation}. Multiplying the controlled state equation \eqref{eq:StateEquation}
by $\boldsymbol{\mathcal{B}}^{T}\left(\boldsymbol{x}\left(t\right)\right)$
from the left yields
\begin{align}
\boldsymbol{\mathcal{B}}^{T}\left(\boldsymbol{x}\left(t\right)\right)\boldsymbol{\dot{x}}\left(t\right)= & \boldsymbol{\mathcal{B}}^{T}\left(\boldsymbol{x}\left(t\right)\right)\boldsymbol{R}\left(\boldsymbol{x}\left(t\right)\right)+\boldsymbol{\mathcal{B}}^{T}\left(\boldsymbol{x}\left(t\right)\right)\boldsymbol{\mathcal{B}}\left(\boldsymbol{x}\left(t\right)\right)\boldsymbol{u}\left(t\right).
\end{align}
Multiplying with $\left(\boldsymbol{\mathcal{B}}^{T}\left(\boldsymbol{x}\left(t\right)\right)\boldsymbol{\mathcal{B}}\left(\boldsymbol{x}\left(t\right)\right)\right)^{-1}$,
which exists as long as $\boldsymbol{\mathcal{B}}\left(\boldsymbol{x}\left(t\right)\right)$
has full rank, from the left results in an expression for the vector
of control signals $\boldsymbol{u}\left(t\right)$ in terms of the
controlled state trajectory $\boldsymbol{x}\left(t\right)$,
\begin{align}
\boldsymbol{u}\left(t\right) & =\boldsymbol{\mathcal{B}}^{+}\left(\boldsymbol{x}\left(t\right)\right)\left(\boldsymbol{\dot{x}}\left(t\right)-\boldsymbol{R}\left(\boldsymbol{x}\left(t\right)\right)\right).\label{eq:NonLinearSystemControlSolution}
\end{align}
The abbreviation
\begin{align}
\boldsymbol{\mathcal{B}}^{+}\left(\boldsymbol{x}\right) & =\left(\boldsymbol{\mathcal{B}}^{T}\left(\boldsymbol{x}\right)\boldsymbol{\mathcal{B}}\left(\boldsymbol{x}\right)\right)^{-1}\boldsymbol{\mathcal{B}}^{T}\left(\boldsymbol{x}\right)
\end{align}
is known as the \textit{Moore-Penrose pseudo inverse} of the matrix
$\boldsymbol{\mathcal{B}}\left(\boldsymbol{x}\right)$ \cite{CampbellJr.1991Generalized}.
See also Appendix \ref{sec:OverAndUnderdetSysOfEqs} how to express
a solution to an overdetermined system of linear equations in terms
of the Moore-Penrose pseudo inverse. With the help of $\boldsymbol{\mathcal{B}}^{+}$,
the projector $\boldsymbol{\mathcal{P}}$ can be expressed as
\begin{align}
\boldsymbol{\mathcal{P}}\left(\boldsymbol{x}\right) & =\boldsymbol{\mathcal{B}}\left(\boldsymbol{x}\right)\left(\boldsymbol{\mathcal{B}}^{T}\left(\boldsymbol{x}\right)\boldsymbol{\mathcal{B}}\left(\boldsymbol{x}\right)\right)^{-1}\boldsymbol{\mathcal{B}}^{T}\left(\boldsymbol{x}\right)=\boldsymbol{\mathcal{B}}\left(\boldsymbol{x}\right)\boldsymbol{\mathcal{B}}^{+}\left(\boldsymbol{x}\right).
\end{align}
Note that
\begin{align}
\boldsymbol{\mathcal{B}}^{+}\left(\boldsymbol{x}\left(t\right)\right)\boldsymbol{\dot{x}}\left(t\right) & =\left(\boldsymbol{\mathcal{B}}^{T}\left(\boldsymbol{x}\left(t\right)\right)\boldsymbol{\mathcal{B}}\left(\boldsymbol{x}\left(t\right)\right)\right)^{-1}\boldsymbol{\mathcal{B}}^{T}\left(\boldsymbol{x}\left(t\right)\right)\boldsymbol{\dot{x}}\left(t\right)\nonumber \\
 & =\left(\boldsymbol{\mathcal{B}}^{T}\left(\boldsymbol{x}\left(t\right)\right)\boldsymbol{\mathcal{B}}\left(\boldsymbol{x}\left(t\right)\right)\right)^{-1}\boldsymbol{\mathcal{B}}^{T}\left(\boldsymbol{x}\left(t\right)\right)\boldsymbol{\mathcal{P}}\left(\boldsymbol{x}\left(t\right)\right)\boldsymbol{\dot{x}}\left(t\right),
\end{align}
such that expression \eqref{eq:NonLinearSystemControlSolution} for
the control involves only the time derivative $\boldsymbol{\mathcal{P}}\left(\boldsymbol{x}\left(t\right)\right)\boldsymbol{\dot{x}}\left(t\right)$
and does not depend on $\boldsymbol{\mathcal{Q}}\left(\boldsymbol{x}\left(t\right)\right)\boldsymbol{\dot{x}}\left(t\right)$.

In conclusion, every affine controlled state equation \eqref{eq:StateEquation}
can be split in two equations. The equation \eqref{eq:NonLinearSystemControlSolution}
involving $\boldsymbol{\mathcal{P}}\boldsymbol{\dot{x}}$ determines
the control signal $\boldsymbol{u}$ in terms of the controlled state
trajectory $\boldsymbol{x}$ and its derivative. The constraint equation
\eqref{eq:ConstraintEquation} involves only $\boldsymbol{\mathcal{Q}}\boldsymbol{\dot{x}}$
and does not depend on the control signal. These relations are valid
for any kind of control, be it an open or a closed loop control. The
proposed separation of the state equation plays a central role in
this thesis.

To illustrate the approach, the separation of the state equation is
discussed with the help of two simple examples.

\begin{example}[Mechanical control system in one spatial dimension]\label{ex:OneDimMechSys2}

The controlled state equation for mechanical control systems is (see
Example \ref{ex:OneDimMechSys1}),

\begin{align}
\left(\begin{array}{c}
\dot{x}\left(t\right)\\
\dot{y}\left(t\right)
\end{array}\right) & =\left(\begin{array}{c}
y\left(t\right)\\
R\left(x\left(t\right),y\left(t\right)\right)
\end{array}\right)+\left(\begin{array}{c}
0\\
B\left(x\left(t\right),y\left(t\right)\right)
\end{array}\right)u\left(t\right).
\end{align}
The $2\times1$ coupling matrix is a vector which depends on the state
vector $\boldsymbol{x}\left(t\right)=\left(\begin{array}{cc}
x\left(t\right), & y\left(t\right)\end{array}\right)$, 
\begin{align}
\boldsymbol{B}\left(\boldsymbol{x}\right) & =\left(\begin{array}{c}
0\\
B\left(x,y\right)
\end{array}\right),
\end{align}
while its transpose is a row vector
\begin{align}
\boldsymbol{B}^{T}\left(\boldsymbol{x}\right) & =\left(\begin{array}{cc}
0, & B\left(x,y\right)\end{array}\right).
\end{align}
The computation of the Moore-Penrose pseudo inverse $\boldsymbol{B}^{+}$
involves the inner product 
\begin{align}
\boldsymbol{B}^{T}\left(\boldsymbol{x}\right)\boldsymbol{B}\left(\boldsymbol{x}\right) & =\left(\begin{array}{cc}
0, & B\left(x,y\right)\end{array}\right)\left(\begin{array}{c}
0\\
B\left(x,y\right)
\end{array}\right)=B\left(x,y\right)^{2}.
\end{align}
The pseudo inverse $\boldsymbol{B}^{+}$ of $\boldsymbol{B}$ is given
by
\begin{align}
\boldsymbol{B}^{+}\left(\boldsymbol{x}\right) & =\left(\boldsymbol{B}^{T}\left(\boldsymbol{x}\right)\boldsymbol{B}\left(\boldsymbol{x}\right)\right)^{-1}\boldsymbol{B}^{T}\left(\boldsymbol{x}\right)=B\left(x,y\right)^{-2}\left(\begin{array}{cc}
0, & B\left(x,y\right)\end{array}\right),
\end{align}
while the projectors $\boldsymbol{\mathcal{P}}\left(\boldsymbol{x}\right)$
and $\boldsymbol{\mathcal{Q}}\left(\boldsymbol{x}\right)$ are given
by 
\begin{align}
\boldsymbol{\mathcal{P}}\left(\boldsymbol{x}\right) & =\boldsymbol{\mathcal{P}}=\boldsymbol{B}\left(\boldsymbol{x}\right)\left(\boldsymbol{B}^{T}\left(\boldsymbol{x}\right)\boldsymbol{B}\left(\boldsymbol{x}\right)\right)^{-1}\boldsymbol{B}^{T}\left(\boldsymbol{x}\right)\nonumber \\
 & =\left(\begin{array}{c}
0\\
B\left(x,y\right)
\end{array}\right)B\left(x,y\right)^{-2}\left(\begin{array}{cc}
0, & B\left(x,y\right)\end{array}\right)=\left(\begin{array}{cc}
0 & 0\\
0 & 1
\end{array}\right),\\
\boldsymbol{\mathcal{Q}}\left(\boldsymbol{x}\right) & =\boldsymbol{\mathcal{Q}}=\mathbf{1}-\boldsymbol{\mathcal{P}}\left(\boldsymbol{x}\right)=\left(\begin{array}{cc}
1 & 0\\
0 & 0
\end{array}\right).
\end{align}
Although the coupling vector $\boldsymbol{B}\left(\boldsymbol{x}\right)$
depends on the state $\boldsymbol{x}$, the projectors $\boldsymbol{\mathcal{P}}$
and $\boldsymbol{\mathcal{Q}}$ are actually independent of the state.
With $\boldsymbol{\mathcal{P}}$ and $\boldsymbol{\mathcal{Q}}$,
the state $\boldsymbol{x}$ can be split up in two parts,
\begin{align}
\boldsymbol{\mathcal{P}}\boldsymbol{x}\left(t\right) & =\left(\begin{array}{c}
0\\
y\left(t\right)
\end{array}\right), & \boldsymbol{\mathcal{Q}}\boldsymbol{x}\left(t\right) & =\left(\begin{array}{c}
x\left(t\right)\\
0
\end{array}\right).
\end{align}
Both parts are vectors with two components, but have only one non-vanishing
component. The control signal can be expressed in terms of the controlled
state trajectory $\boldsymbol{x}\left(t\right)$ as 
\begin{align}
u\left(t\right) & =\boldsymbol{B}^{+}\left(\boldsymbol{x}\left(t\right)\right)\left(\boldsymbol{\dot{x}}\left(t\right)-\boldsymbol{R}\left(\boldsymbol{x}\left(t\right)\right)\right)\nonumber \\
 & =B\left(x\left(t\right),y\left(t\right)\right)^{-2}\left(\begin{array}{cc}
0, & B\left(x,y\right)\end{array}\right)\left(\left(\begin{array}{c}
\dot{x}\left(t\right)\\
\dot{y}\left(t\right)
\end{array}\right)-\left(\begin{array}{c}
y\left(t\right)\\
R\left(x\left(t\right),y\left(t\right)\right)
\end{array}\right)\right)\nonumber \\
 & =\frac{1}{B\left(x\left(t\right),y\left(t\right)\right)}\left(\dot{y}\left(t\right)-R\left(x\left(t\right),y\left(t\right)\right)\right).
\end{align}
Note that the assumption of full rank for the coupling vector $\boldsymbol{B}$
implies that the function $B\left(x,y\right)$ does not vanish, i.e.,
$B\left(x,y\right)\neq0$ for all values of $x$ and $y$. Consequently,
$u\left(t\right)$ is well defined for all times.

\end{example}

\begin{example}[Single input diagonal LTI system]\label{ex:DiagonalLTI}

Consider a diagonal $2\times2$ linear time-invariant (LTI) system
for the state vector $\boldsymbol{x}\left(t\right)=\left(\begin{array}{cc}
x_{1}\left(t\right), & x_{2}\left(t\right)\end{array}\right)^{T}$. Let both components be controlled by the same control signal $u\left(t\right)$,
\begin{align}
\dot{x}_{1}\left(t\right) & =\lambda_{1}x_{1}\left(t\right)+u\left(t\right), & \dot{x}_{2}\left(t\right) & =\lambda_{2}x_{2}\left(t\right)+u\left(t\right).
\end{align}
The state matrix $\boldsymbol{\mathcal{A}}$ and input matrix $\boldsymbol{\mathcal{B}}$
are
\begin{align}
\boldsymbol{\mathcal{A}} & =\left(\begin{array}{cc}
\lambda_{1} & 0\\
0 & \lambda_{2}
\end{array}\right), & \boldsymbol{\mathcal{B}} & =\left(\begin{array}{c}
1\\
1
\end{array}\right).
\end{align}
The constant projectors $\boldsymbol{\mathcal{P}}$ and $\boldsymbol{\mathcal{Q}}$
can be computed as
\begin{align}
\boldsymbol{\mathcal{P}} & =\frac{1}{2}\left(\begin{array}{cc}
1 & 1\\
1 & 1
\end{array}\right), & \boldsymbol{\mathcal{Q}} & =\frac{1}{2}\left(\begin{array}{cc}
1 & -1\\
-1 & 1
\end{array}\right).
\end{align}
The two projections of the state $\boldsymbol{x}\left(t\right)$ are
\begin{align}
\boldsymbol{z}\left(t\right) & =\left(\begin{array}{c}
z_{1}\left(t\right)\\
z_{2}\left(t\right)
\end{array}\right)=\boldsymbol{\mathcal{P}}\boldsymbol{x}\left(t\right)=\frac{1}{2}\left(\begin{array}{c}
x_{1}\left(t\right)+x_{2}\left(t\right)\\
x_{1}\left(t\right)+x_{2}\left(t\right)
\end{array}\right)
\end{align}
and
\begin{align}
\boldsymbol{y}\left(t\right) & =\left(\begin{array}{c}
y_{1}\left(t\right)\\
y_{2}\left(t\right)
\end{array}\right)=\boldsymbol{\mathcal{Q}}\boldsymbol{x}\left(t\right)=\frac{1}{2}\left(\begin{array}{c}
x_{1}\left(t\right)-x_{2}\left(t\right)\\
x_{2}\left(t\right)-x_{1}\left(t\right)
\end{array}\right).
\end{align}
While both components of $\boldsymbol{y}\left(t\right)$ are non-zero,
they are not linearly independent. The component $y_{2}\left(t\right)$
is redundant and is simply given by
\begin{align}
y_{2}\left(t\right) & =-y_{1}\left(t\right).
\end{align}
Similarly, the component $z_{2}\left(t\right)$ of vector $\boldsymbol{z}\left(t\right)$
is redundant because of
\begin{align}
z_{2}\left(t\right) & =z_{1}\left(t\right).
\end{align}
\end{example}

Example \ref{ex:DiagonalLTI} shows that the projectors $\boldsymbol{\mathcal{P}}$
and $\boldsymbol{\mathcal{Q}}$ do not necessarily project onto single
components of the state vector. If the projectors $\boldsymbol{\mathcal{P}}\left(\boldsymbol{x}\right)=\boldsymbol{\mathcal{P}}$
and $\boldsymbol{\mathcal{Q}}\left(\boldsymbol{x}\right)=\boldsymbol{\mathcal{Q}}$
are independent of the state $\boldsymbol{x}$, the parts $\boldsymbol{z}=\boldsymbol{\mathcal{P}}\boldsymbol{x}$
and $\boldsymbol{y}=\boldsymbol{\mathcal{Q}}\boldsymbol{x}$ are linear
combinations of the original state components $\boldsymbol{x}$. If
$\boldsymbol{\mathcal{P}}=\boldsymbol{\mathcal{P}}\left(\boldsymbol{x}\right)$
and therefore also $\boldsymbol{\mathcal{Q}}\left(\boldsymbol{x}\right)=\boldsymbol{1}-\boldsymbol{\mathcal{P}}\left(\boldsymbol{x}\right)$
depend on the state $\boldsymbol{x}$ itself, both parts $\boldsymbol{y}$
and $\boldsymbol{z}$ are nonlinear functions of the state $\boldsymbol{x}$.
Only if the projectors are diagonal, constant, and appropriately ordered,
$\boldsymbol{\mathcal{P}}\left(\boldsymbol{x}\right)=\boldsymbol{\mathcal{P}}_{D}$
and $\boldsymbol{\mathcal{Q}}\left(\boldsymbol{x}\right)=\boldsymbol{\mathcal{Q}}_{D}$,
then the two parts $\boldsymbol{y}$ and $\boldsymbol{z}$ attain
the particularly simple form
\begin{align}
\boldsymbol{y} & =\boldsymbol{\mathcal{Q}}_{D}\boldsymbol{x}=\left(\begin{array}{cccccc}
0, & \dots, & 0, & x_{p+1}, & \dots, & x_{n}\end{array}\right)^{T},\label{eq:Eq233}\\
\boldsymbol{z} & =\boldsymbol{\mathcal{P}}_{D}\boldsymbol{x}=\left(\begin{array}{cccccc}
x_{1}, & \dots, & x_{p}, & 0, & \dots, & 0\end{array}\right)^{T}.\label{eq:Eq234}
\end{align}
Only this form allows a clear interpretation which component of $\boldsymbol{x}$
belongs to which part. However, in any case, $\boldsymbol{y}=\boldsymbol{\mathcal{Q}}\left(\boldsymbol{x}\right)\boldsymbol{x}$
has exactly $n-p$ independent components because $\boldsymbol{\mathcal{Q}}\left(\boldsymbol{x}\right)$
has rank $n-p$, while $\boldsymbol{z}=\boldsymbol{\mathcal{P}}\left(\boldsymbol{x}\right)\boldsymbol{x}$
has $p$ independent components because $\boldsymbol{\mathcal{P}}\left(\boldsymbol{x}\right)$
has rank $p$.

Projectors have only zeros and ones as possible eigenvalues. The diagonalization
of an $n\times n$ projector $\boldsymbol{\mathcal{P}}\left(\boldsymbol{x}\right)$
with $\text{rank}\left(\boldsymbol{\mathcal{P}}\left(\boldsymbol{x}\right)\right)=p$
is always possible \cite{fischer2008lineare} and results in a diagonal
$n\times n$ matrix with $p$ entries of value one and $n-p$ entries
of value zero on the diagonal. The transformation of the projectors
$\boldsymbol{\mathcal{P}}\left(\boldsymbol{x}\right)$ and $\boldsymbol{\mathcal{Q}}\left(\boldsymbol{x}\right)$
to their diagonal counterparts defines a transformation of the state
$\boldsymbol{x}$. See Appendix \ref{sec:DiagonalizingTheConstraintEquation}
how to construct this transformation. The transformation is nonlinear
if $\boldsymbol{\mathcal{P}}\left(\boldsymbol{x}\right)$ and $\boldsymbol{\mathcal{Q}}\left(\boldsymbol{x}\right)$
are state dependent. Expressed in terms of the transformed state,
the projectors $\boldsymbol{\mathcal{P}}$ and $\boldsymbol{\mathcal{Q}}$
are constant, diagonal, and appropriately ordered. Consequently, they
yield a state separation of the form Eqs. \eqref{eq:Eq233}-\eqref{eq:Eq234}.
Such a representation defines a normal form of an affine control system.
For a specified affine control system, computations will usually be
simpler after the system is transformed to its normal form. However,
for computations with general affine control systems, it is dispensable
to perform the transformation if $\boldsymbol{z}\left(\boldsymbol{x}\right)=\boldsymbol{\mathcal{P}}\left(\boldsymbol{x}\right)\boldsymbol{x}$
and $\boldsymbol{y}\left(\boldsymbol{x}\right)=\boldsymbol{\mathcal{Q}}\left(\boldsymbol{x}\right)\boldsymbol{x}$
are simply viewed as separate parts. This allows a coordinate-free
treatment of affine control systems.

\section{\label{sec:ExactlyRealizableTrajectories}Exactly realizable trajectories}

As demonstrated in Example \ref{ex:ControlledFHN1}, not every desired
state trajectory $\boldsymbol{x}_{d}\left(t\right)$ can be realized
by control. Here, we answer the question under which conditions a
desired trajectory $\boldsymbol{x}_{d}\left(t\right)$ is exactly
realizable.

Consider the controlled state equation
\begin{align}
\boldsymbol{\dot{x}}\left(t\right) & =\boldsymbol{R}\left(\boldsymbol{x}\left(t\right)\right)+\boldsymbol{\mathcal{B}}\left(\boldsymbol{x}\left(t\right)\right)\boldsymbol{u}\left(t\right),\label{eq:ControlledStateEquation}\\
\boldsymbol{x}\left(t_{0}\right) & =\boldsymbol{x}_{0}.
\end{align}
The notion of \textit{exactly realizable trajectories} is introduced.
A realizable trajectory is a desired trajectory $\boldsymbol{x}_{d}\left(t\right)$
which satisfies two conditions.
\begin{enumerate}
\item \label{enu:RealizableTrajectory1}The desired trajectory $\boldsymbol{x}_{d}\left(t\right)$
satisfies the constraint equation
\begin{align}
\boldsymbol{\mathcal{Q}}\left(\boldsymbol{x}_{d}\left(t\right)\right)\left(\boldsymbol{\dot{x}}_{d}\left(t\right)-\boldsymbol{R}\left(\boldsymbol{x}_{d}\left(t\right)\right)\right) & =\mathbf{0}.\label{eq:ConstraintEquationForRealizableTrajectories}
\end{align}

\item \label{enu:RealizableTrajectory2}The initial value $\boldsymbol{x}_{d}\left(t_{0}\right)$
must equal the initial value $\boldsymbol{x}_{0}$ of the controlled
state equation,
\begin{align}
\boldsymbol{x}_{d}\left(t_{0}\right) & =\boldsymbol{x}_{0}.
\end{align}

\end{enumerate}
The control solution for an exactly realizable trajectory is given
by
\begin{align}
\boldsymbol{u}\left(t\right) & =\boldsymbol{\mathcal{B}}^{+}\left(\boldsymbol{x}_{d}\left(t\right)\right)\left(\boldsymbol{\dot{x}}_{d}\left(t\right)-\boldsymbol{R}\left(\boldsymbol{x}_{d}\left(t\right)\right)\right),\label{eq:UInTermsOfx_r}
\end{align}
with the Moore-Penrose pseudo inverse $p\times n$ matrix $\boldsymbol{\mathcal{B}}^{+}$
defined as 
\begin{align}
\boldsymbol{\mathcal{B}}^{+}\left(\boldsymbol{x}\right) & =\left(\boldsymbol{\mathcal{B}}^{T}\left(\boldsymbol{x}\right)\boldsymbol{\mathcal{B}}\left(\boldsymbol{x}\right)\right)^{-1}\boldsymbol{\mathcal{B}}^{T}\left(\boldsymbol{x}\right).
\end{align}
The notion of an exactly realizable trajectory allows the proof of
the following statement. 

If $\boldsymbol{x}_{d}\left(t\right)$ is an exactly realizable trajectory,
i.e., if $\boldsymbol{x}_{d}\left(t\right)$ satisfies both conditions
\ref{enu:RealizableTrajectory1} and \ref{enu:RealizableTrajectory2},
then the state trajectory $\boldsymbol{x}\left(t\right)$ follows
the desired trajectory $\boldsymbol{x}_{d}\left(t\right)$ exactly,
\begin{align}
\boldsymbol{x}\left(t\right) & =\boldsymbol{x}_{d}\left(t\right).
\end{align}
Using the control solution Eq. \eqref{eq:UInTermsOfx_r} in the controlled
state equation \eqref{eq:ControlledStateEquation} yields the following
equation for the controlled state
\begin{align}
\boldsymbol{\dot{x}}\left(t\right) & =\boldsymbol{R}\left(\boldsymbol{x}\left(t\right)\right)+\boldsymbol{\mathcal{B}}\left(\boldsymbol{x}\left(t\right)\right)\boldsymbol{\mathcal{B}}^{+}\left(\boldsymbol{x}_{d}\left(t\right)\right)\left(\boldsymbol{\dot{x}}_{d}\left(t\right)-\boldsymbol{R}\left(\boldsymbol{x}_{d}\left(t\right)\right)\right).\label{eq:E16}
\end{align}
Note that $\boldsymbol{\mathcal{B}}$ depends on the actual system
state $\boldsymbol{x}\left(t\right)$ while $\boldsymbol{\mathcal{B}}^{+}$
depends on the desired trajectory $\boldsymbol{x}_{d}\left(t\right)$.
The difference $\Delta\boldsymbol{x}\left(t\right)$ between the true
state $\boldsymbol{x}\left(t\right)$ and the desired trajectory $\boldsymbol{x}_{d}\left(t\right)$
is defined as 
\begin{align}
\Delta\boldsymbol{x}\left(t\right) & =\boldsymbol{x}\left(t\right)-\boldsymbol{x}_{d}\left(t\right).
\end{align}
Using the definition for $\Delta\boldsymbol{x}\left(t\right)$ and
Eq. \eqref{eq:E16} results in an ordinary differential equation (ODE)
for $\Delta\boldsymbol{x}\left(t\right)$,
\begin{align}
\Delta\boldsymbol{\dot{x}}\left(t\right) & =\boldsymbol{R}\left(\Delta\boldsymbol{x}\left(t\right)+\boldsymbol{x}_{d}\left(t\right)\right)-\boldsymbol{\dot{x}}_{d}\left(t\right)\nonumber \\
 & +\boldsymbol{\mathcal{B}}\left(\Delta\boldsymbol{x}\left(t\right)+\boldsymbol{x}_{d}\left(t\right)\right)\boldsymbol{\mathcal{B}}^{+}\left(\boldsymbol{x}_{d}\left(t\right)\right)\left(\boldsymbol{\dot{x}}_{d}\left(t\right)-\boldsymbol{R}\left(\boldsymbol{x}_{d}\left(t\right)\right)\right),\label{eq:ControlledStateEquationForDeltax}\\
\Delta\boldsymbol{x}\left(t_{0}\right) & =\boldsymbol{x}\left(t_{0}\right)-\boldsymbol{x}_{d}\left(t_{0}\right).
\end{align}
Assuming $\left|\Delta\boldsymbol{x}\left(t\right)\right|\ll1$ and
expanding Eq. \eqref{eq:ControlledStateEquationForDeltax} in $\Delta\boldsymbol{x}\left(t\right)$
yields 
\begin{align}
\Delta\boldsymbol{\dot{x}}\left(t\right) & =\boldsymbol{R}\left(\boldsymbol{x}_{d}\left(t\right)\right)-\boldsymbol{\dot{x}}_{d}\left(t\right)+\boldsymbol{\mathcal{B}}\left(\boldsymbol{x}_{d}\left(t\right)\right)\boldsymbol{\mathcal{B}}^{+}\left(\boldsymbol{x}_{d}\left(t\right)\right)\left(\boldsymbol{\dot{x}}_{d}\left(t\right)-\boldsymbol{R}\left(\boldsymbol{x}_{d}\left(t\right)\right)\right)\nonumber \\
 & +\nabla\boldsymbol{R}\left(\boldsymbol{x}_{d}\left(t\right)\right)\Delta\boldsymbol{x}\left(t\right)+\left(\nabla\boldsymbol{\mathcal{B}}\left(\boldsymbol{x}_{d}\left(t\right)\right)\Delta\boldsymbol{x}\left(t\right)\right)\boldsymbol{\mathcal{B}}^{+}\left(\boldsymbol{x}_{d}\left(t\right)\right)\left(\boldsymbol{\dot{x}}_{d}\left(t\right)-\boldsymbol{R}\left(\boldsymbol{x}_{d}\left(t\right)\right)\right)\nonumber \\
 & +\mathcal{O}\left(\Delta\boldsymbol{x}\left(t\right)^{2}\right).\label{eq:ControlledStateEquationForDeltaxExpanded}
\end{align}
Note that assuming $\left|\Delta\boldsymbol{x}\left(t\right)\right|\ll1$
and subsequently expanding in $\Delta\boldsymbol{x}\left(t\right)$
does not result in a loss of generality of the final outcome. The
expression $\nabla\boldsymbol{R}\left(\boldsymbol{x}\right)$ denotes
the Jacobian matrix of the nonlinearity $\boldsymbol{R}\left(\boldsymbol{x}\right)$
with components 
\begin{align}
\left(\nabla\boldsymbol{R}\left(\boldsymbol{x}\right)\right)_{ij} & =\dfrac{\partial}{\partial x_{j}}R_{i}\left(\boldsymbol{x}\right),\,i,j\in\left\{ 1,\dots,n\right\} .
\end{align}
The Jacobian $\nabla\boldsymbol{\mathcal{B}}\left(\boldsymbol{x}\right)$
of $\boldsymbol{\mathcal{B}}$ is a third order tensor with components
\begin{align}
\left(\nabla\boldsymbol{\mathcal{B}}\left(\boldsymbol{x}\right)\right)_{ijk} & =\dfrac{\partial}{\partial x_{k}}\mathcal{B}_{ij}\left(\boldsymbol{x}\right),\,i,k\in\left\{ 1,\dots,n\right\} ,\,j\in\left\{ 1,\dots p\right\} .
\end{align}
In the first line of Eq. \eqref{eq:ControlledStateEquationForDeltaxExpanded},
one can recognize the projector $\boldsymbol{\mathcal{B}}\left(\boldsymbol{x}_{d}\left(t\right)\right)\boldsymbol{\mathcal{B}}^{+}\left(\boldsymbol{x}_{d}\left(t\right)\right)=\boldsymbol{\mathcal{P}}\left(\boldsymbol{x}_{d}\left(t\right)\right)=\mathbf{1}-\boldsymbol{\mathcal{Q}}\left(\boldsymbol{x}_{d}\left(t\right)\right)$.
Introducing the $n\times n$ matrix $\boldsymbol{\mathcal{T}}\left(\boldsymbol{x}\right)$
with components
\begin{align}
\left(\boldsymbol{\mathcal{T}}\left(\boldsymbol{x}\right)\right)_{il} & =\sum_{j=1}^{p}\sum_{k=1}^{n}\dfrac{\partial}{\partial x_{l}}\mathcal{B}_{ij}\left(\boldsymbol{x}\right)\mathcal{B}_{jk}^{+}\left(\boldsymbol{x}\right)\left(\dot{x}_{k}\left(t\right)-R_{k}\left(\boldsymbol{x}\right)\right),\,i,l\in\left\{ 1,\dots,n\right\} ,
\end{align}
allows a rearrangement of Eq. \eqref{eq:ControlledStateEquationForDeltaxExpanded}
in the form
\begin{align}
\Delta\boldsymbol{\dot{x}}\left(t\right) & =\boldsymbol{\mathcal{Q}}\left(\boldsymbol{x}_{d}\left(t\right)\right)\left(\boldsymbol{R}\left(\boldsymbol{x}_{d}\left(t\right)\right)-\boldsymbol{\dot{x}}_{d}\left(t\right)\right)+\left(\nabla\boldsymbol{R}\left(\boldsymbol{x}_{d}\left(t\right)\right)+\boldsymbol{\mathcal{T}}\left(\boldsymbol{x}_{d}\left(t\right)\right)\right)\Delta\boldsymbol{x}\left(t\right),\label{eq:Eq122}\\
\Delta\boldsymbol{x}\left(t_{0}\right) & =\boldsymbol{x}\left(t_{0}\right)-\boldsymbol{x}_{d}\left(t_{0}\right).
\end{align}
If $\boldsymbol{x}_{d}\left(t\right)$ is an exactly realizable trajectory,
it satisfies the constraint equation \eqref{eq:ConstraintEquationForRealizableTrajectories}
and the initial condition $\boldsymbol{x}_{d}\left(t_{0}\right)=\boldsymbol{x}\left(t_{0}\right)$,
and Eq. \eqref{eq:Eq122} simplifies to the linear homogeneous equation
for $\Delta\boldsymbol{x}\left(t\right)$, 
\begin{align}
\Delta\boldsymbol{\dot{x}}\left(t\right) & =\left(\nabla\boldsymbol{R}\left(\boldsymbol{x}_{d}\left(t\right)\right)+\boldsymbol{\mathcal{T}}\left(\boldsymbol{x}_{d}\left(t\right)\right)\right)\Delta\boldsymbol{x}\left(t\right)\label{eq:EqForLinearStability}\\
\Delta\boldsymbol{x}\left(t_{0}\right) & =\mathbf{0}.
\end{align}
Clearly, Eq. \eqref{eq:EqForLinearStability} has a vanishing solution
\begin{align}
\Delta\boldsymbol{x}\left(t\right) & \equiv\mathbf{0}.
\end{align}
In summary, it was proven that if the desired trajectory $\boldsymbol{x}_{d}\left(t\right)$
is an exactly realizable trajectory, then the state trajectory $\boldsymbol{x}\left(t\right)$
follows the desired trajectory exactly, i.e., $\boldsymbol{x}\left(t\right)=\boldsymbol{x}_{d}\left(t\right)$.
Equation \eqref{eq:ConstraintEquation} in Section \ref{sub:SeparationOfTheStateEq}
proved already the converse: any controlled state trajectory $\boldsymbol{x}\left(t\right)$
satisfies the constraint equation. This allows the following conclusion:

\textit{The controlled state trajectory $\boldsymbol{x}\left(t\right)$
follows the desired trajectory $\boldsymbol{x}_{d}\left(t\right)$
exactly if and only if $\boldsymbol{x}_{d}\left(t\right)$ is an exactly
realizable trajectory.}

The notion of an exactly realizable trajectory leads to the following
interpretation. Not every desired trajectory $\boldsymbol{x}_{d}\left(t\right)$
can be enforced in a specified controlled dynamical system. In general,
the desired trajectory $\boldsymbol{x}_{d}\left(t\right)$ is that
what you want, but is not what you get. What you get is an exactly
realizable trajectory. Because the control signal $\boldsymbol{u}\left(t\right)$
consists of only $p$ independent components, it is possible to find
at most $p$ one-to-one relations between state components and components
of the control signal. Only $p$ components of a state trajectory
can be prescribed, while the remaining $n-p$ components are free.
The time evolution of these $n-p$ components is given by the constraint
equation \eqref{eq:ConstraintEquationForRealizableTrajectories}.
This motivates the name constraint equation. For an arbitrary desired
trajectory $\boldsymbol{x}_{d}\left(t\right)$ to be exactly realizable,
it has to be constrained by Eq. \eqref{eq:ConstraintEquationForRealizableTrajectories}.
There is still some freedom to choose which state components are actually
prescribed, and which have to be determined by the constraint equation.
Until further notice, we adopt the canonical view that the part $\boldsymbol{\mathcal{P}}\boldsymbol{x}_{d}\left(t\right)$
is prescribed by the experimenter, while the part $\boldsymbol{\mathcal{Q}}\boldsymbol{x}_{d}\left(t\right)$
of the state vector is fixed by the constraint equation \eqref{eq:ConstraintEquationForRealizableTrajectories}.
This, however, is not the only possibility, and many more choices
are possible. The part $\boldsymbol{\mathcal{P}}\boldsymbol{x}_{d}\left(t\right)$
can be seen as an output for the control system which can be enforced
exactly if $\boldsymbol{x}_{d}\left(t\right)$ is exactly realizable.
Section \ref{sec:OutputRealizability} discusses the possibility to
realize general desired outputs not necessarily given by $\boldsymbol{\mathcal{P}}\boldsymbol{x}_{d}\left(t\right)$.

Chapter \ref{chap:OptimalControl} investigates the relation of exactly
realizable trajectories with optimal trajectory tracking. The control
solution Eq. \eqref{eq:UInTermsOfx_r} is the solution to a certain
optimal trajectory tracking problem. This insight is the starting
point in Chapter \ref{chap:AnalyticalApproximationsForOptimalTrajectoryTracking}
to obtain analytical approximations to optimal trajectory tracking
of desired trajectories which are not exactly realizable.

The necessity to satisfy condition \ref{enu:RealizableTrajectory2}
of equal initial condition leaves two possibilities. Either the system
is prepared in the initial state $\boldsymbol{x}\left(t_{0}\right)=\boldsymbol{x}_{0}=\boldsymbol{x}_{d}\left(t_{0}\right)$,
or the desired trajectory $\boldsymbol{x}_{d}\left(t\right)$ is designed
such that it starts from the observed initial system state $\boldsymbol{x}_{0}$.
In any case, the constraint equation \eqref{eq:ConstraintEquationForRealizableTrajectories},
seen as an ODE for $\boldsymbol{\mathcal{Q}}\boldsymbol{x}_{d}\left(t\right)$,
has to be solved with the initial condition $\boldsymbol{\mathcal{Q}}\boldsymbol{x}_{d}\left(t_{0}\right)=\boldsymbol{\mathcal{Q}}\boldsymbol{x}_{0}$.

The control solution as given by Eq. \eqref{eq:UInTermsOfx_r} is
an open loop control. As such, it does not guarantee a stable time
evolution, and the controlled system does not necessarily follow the
realizable trajectory in the presence of perturbations. The linear
equation \eqref{eq:EqForLinearStability} encountered during the proof
implies statements about the linear stability of realizable trajectories.
A non-vanishing initial value $\Delta\boldsymbol{x}\left(t_{0}\right)=\Delta\boldsymbol{x}_{0}\neq\mathbf{0}$
constitutes a perturbation of the initial conditions of an exactly
realizable trajectory. The control approach as proposed here is only
a first step. For a specified exactly realizable trajectory, Eq. \eqref{eq:EqForLinearStability}
has to be investigated to determine its linear stability properties.
If the desired trajectory is linearly unstable, countermeasures in
form of an additional feedback control, for example, have to be applied
to guarantee a successful control. Stability of exactly realizable
trajectories is not discussed in this thesis.

The concept of a realizable trajectory is elucidated with the help
of some examples in the following.

\begin{example}[Controlled FHN model with invertible coupling matrix]\label{ex:FHN2}

Consider the controlled FHN model in the form
\begin{align}
\left(\begin{array}{c}
\dot{x}\left(t\right)\\
\dot{y}\left(t\right)
\end{array}\right) & =\left(\begin{array}{c}
a_{0}+a_{1}x\left(t\right)+a_{2}y\left(t\right)\\
R\left(x\left(t\right),y\left(t\right)\right)
\end{array}\right)+\left(\begin{array}{c}
u_{1}\left(t\right)\\
u_{2}\left(t\right)
\end{array}\right).
\end{align}
The constant coupling matrix $\boldsymbol{\mathcal{B}}$ is identical
to the identity,
\begin{align}
\boldsymbol{\mathcal{B}} & =\left(\begin{array}{cc}
1 & 0\\
0 & 1
\end{array}\right).
\end{align}
This system has two state components $x,\,y$, and two independent
control signals $u_{1},\,u_{2}$. The projector $\boldsymbol{\mathcal{P}}$
is simply the identity, $\boldsymbol{\mathcal{P}}=\boldsymbol{1}$,
and $\boldsymbol{\mathcal{Q}}=\boldsymbol{0}$ the zero matrix. The
constraint equation \eqref{eq:ConstraintEquationForRealizableTrajectories}
is trivially satisfied. Any desired trajectory $\boldsymbol{x}_{d}\left(t\right)$
is a realizable trajectory as long as initially, the desired trajectory
equals the state trajectory,
\begin{align}
\boldsymbol{x}_{d}\left(t_{0}\right) & =\boldsymbol{x}\left(t_{0}\right).
\end{align}
\end{example}

\begin{example}[Mechanical control system in one spatial dimension]\label{ex:OneDimMechSys3}

The control signal realizing a desired trajectory \textit{$\boldsymbol{x}_{d}\left(t\right)$}
of a mechanical control system (see Examples \ref{ex:OneDimMechSys1}
and \ref{ex:OneDimMechSys2} for more details) is
\begin{align}
u\left(t\right) & =\boldsymbol{B}^{+}\left(\boldsymbol{x}_{d}\left(t\right)\right)\left(\boldsymbol{\dot{x}}_{d}\left(t\right)-\boldsymbol{R}\left(\boldsymbol{x}_{d}\left(t\right)\right)\right)\nonumber \\
 & =\frac{1}{B\left(x_{d}\left(t\right),y_{d}\left(t\right)\right)}\left(\dot{y}_{d}\left(t\right)-R\left(x_{d}\left(t\right),y_{d}\left(t\right)\right)\right).\label{eq:ComputedTorque}
\end{align}
The non-vanishing component of the constraint equation \eqref{eq:ConstraintEquationForRealizableTrajectories}
for realizable desired trajectories simply becomes
\begin{align}
\dot{x}_{d}\left(t\right) & =y_{d}\left(t\right).\label{eq:ConstEqOneDimMechSys}
\end{align}
With a scalar control signal $u\left(t\right)$ only one state component
can be controlled. According to our convention, this state component
is
\begin{align}
\boldsymbol{\mathcal{P}}\boldsymbol{x}_{d}\left(t\right) & =\left(\begin{array}{c}
0\\
y_{d}\left(t\right)
\end{array}\right).\label{eq:OutPutConvention}
\end{align}
The desired velocity over time $y_{d}\left(t\right)$ can be arbitrarily
chosen apart from its initial value, which must be identical to the
initial state velocity, $y_{d}\left(t_{0}\right)=y\left(t_{0}\right)$.
The corresponding position over time $x_{d}\left(t\right)$ is given
by the constraint equation \eqref{eq:ConstEqOneDimMechSys}. Because
Eq. \eqref{eq:ConstEqOneDimMechSys} is a linear differential equation
for $x_{d}\left(t\right)$, its solution in terms of the arbitrary
velocity $y_{d}\left(t\right)$ is easily obtained as
\begin{align}
x_{d}\left(t\right) & =x_{d}\left(t_{0}\right)+\intop_{t_{0}}^{t}d\tau y_{d}\left(\tau\right).\label{eq:ConstEqSolOneDimMechSys}
\end{align}
The initial desired position $x_{d}\left(t_{0}\right)$ has to agree
with the initial state position, $x_{d}\left(t_{0}\right)=x\left(t_{0}\right)=x_{0}$.
With the help of solution \eqref{eq:ConstEqSolOneDimMechSys}, the
control \eqref{eq:ComputedTorque} can be entirely expressed in terms
of the prescribed velocity over time $y_{d}\left(t\right)$ as
\begin{align}
u\left(t\right) & =\frac{1}{B\left(x_{0}+\intop_{t_{0}}^{t}d\tau y_{d}\left(\tau\right),y_{d}\left(t\right)\right)}\nonumber \\
 & \times\left(\dot{y}_{d}\left(t\right)-R\left(x_{0}+\intop_{t_{0}}^{t}d\tau y_{d}\left(\tau\right),y_{d}\left(t\right)\right)\right).
\end{align}
Note that an exact solution to the nonlinear controlled state equation
as well as to the control signal $\boldsymbol{u}\left(t\right)$ is
obtained without actually solving any nonlinear equation. The context
of a mechanical control system allows the following interpretation
of our approach. The constraint equation \eqref{eq:ConstEqOneDimMechSys}
is the definition of the velocity of a point particle, and no external
force $R$ or control force $Bu$ can change that definition. With
only a single control signal $u$, position $x$ and velocity $y$
over time cannot be controlled independently from each other.

One might ask if it is possible to control position and velocity independently
of each other by introducing an additional control signal. If both
control signal act as forces, the controlled mechanical system with
state space dimension $n=2$ and control space dimension $p=2$ is
\begin{align}
\left(\begin{array}{c}
\dot{x}\left(t\right)\\
\dot{y}\left(t\right)
\end{array}\right) & =\left(\begin{array}{c}
y\left(t\right)\\
R\left(x\left(t\right),y\left(t\right)\right)
\end{array}\right)\nonumber \\
 & +\left(\begin{array}{cc}
0 & 0\\
B_{1}\left(x\left(t\right),y\left(t\right)\right) & B_{2}\left(x\left(t\right),y\left(t\right)\right)
\end{array}\right)\left(\begin{array}{c}
u_{1}\left(t\right)\\
u_{2}\left(t\right)
\end{array}\right),\label{eq:ConstEqOneDimMechSys2}
\end{align}
such that the $2\times2$ coupling matrix $\boldsymbol{\mathcal{\tilde{B}}}$
becomes
\begin{align}
\boldsymbol{\mathcal{\tilde{B}}}\left(\boldsymbol{x}\right) & =\left(\begin{array}{cc}
0 & 0\\
B_{1}\left(x,y\right) & B_{2}\left(x,y\right)
\end{array}\right).
\end{align}
However, the structure of $\boldsymbol{\mathcal{\tilde{B}}}$ reveals
that it violates the condition of full rank. Indeed, for arbitrary
functions $B_{1}\neq0$ and $B_{2}\neq0$, the rank of $\boldsymbol{\mathcal{\tilde{B}}}$
is $\mbox{rank}\left(\boldsymbol{\mathcal{\tilde{B}}}\left(\boldsymbol{x}\right)\right)=1$
and therefore smaller than the control space dimension $p=2$. The
computation of the projectors $\boldsymbol{\mathcal{P}}$ and $\boldsymbol{\mathcal{Q}}$
as well as the computation of the control signal $u$ requires the
existence of the inverse of $\boldsymbol{\mathcal{\tilde{B}}}^{T}\boldsymbol{\mathcal{\tilde{B}}}$,
which in turn requires $\boldsymbol{\mathcal{\tilde{B}}}$ to have
full rank. Our approach cannot be applied to system \eqref{eq:ConstEqOneDimMechSys2}
because both control signals $u_{1}$ and $u_{2}$ act on the same
state component. The corresponding control forces are not independent
of each other, but can be combined to a single control force $B_{1}u_{1}+B_{2}u_{2}$.

The constraint equation \eqref{eq:ConstEqOneDimMechSys} can also
be regarded as an algebraic equation for the desired position over
time $x_{d}\left(t\right)$. This is an example for a desired output
different from the conventional choice Eq. \eqref{eq:OutPutConvention}.
Eliminating the position from the control solution Eq. \eqref{eq:ComputedTorque}
yields
\begin{align}
u\left(t\right) & =\frac{1}{B\left(x_{d}\left(t\right),\dot{x}_{d}\left(t\right)\right)}\left(\ddot{x}_{d}\left(t\right)-R\left(x_{d}\left(t\right),\dot{x}_{d}\left(t\right)\right)\right).\label{eq:ComputedTorque2}
\end{align}
Equation \eqref{eq:ComputedTorque2} is a special case of the so-called
\textit{computed torque formula}. This approach, also known as inverse
dynamics, is regularly applied in robotics. For further information,
the reader is referred to the literature about robot control \cite{lewis1993control,deWit1996theory,angeles2002fundamentals}.

\end{example}

\begin{example}[Activator-controlled FHN model]\label{ex:FHN3}

Consider the FHN model from Example \ref{ex:FHN1},
\begin{align}
\left(\begin{array}{c}
\dot{x}\left(t\right)\\
\dot{y}\left(t\right)
\end{array}\right) & =\left(\begin{array}{c}
a_{0}+a_{1}x\left(t\right)+a_{2}y\left(t\right)\\
R\left(x\left(t\right),y\left(t\right)\right)
\end{array}\right)+\left(\begin{array}{c}
0\\
1
\end{array}\right)u\left(t\right),\label{eq:Eq267}
\end{align}
with coupling vector $\boldsymbol{B}=\left(\begin{array}{cc}
0, & 1\end{array}\right)^{T}$ and standard FHN nonlinearity $R\left(x,y\right)=y-\frac{1}{3}y^{3}-x$.
The projectors $\boldsymbol{\mathcal{P}}$ and $\boldsymbol{\mathcal{Q}}$
are readily computed as 
\begin{align}
\boldsymbol{\mathcal{P}} & =\left(\begin{array}{cc}
0 & 0\\
0 & 1
\end{array}\right), & \boldsymbol{\mathcal{Q}} & =\left(\begin{array}{cc}
1 & 0\\
0 & 0
\end{array}\right).
\end{align}
Being given by $\boldsymbol{\mathcal{P}}\boldsymbol{x}$, the desired
activator component over time $y_{d}\left(t\right)$ can be prescribed.
For the desired trajectory to be exactly realizable, the desired inhibitor
$x_{d}\left(t\right)$ must be determined from the constraint equation
\begin{align}
\boldsymbol{\mathcal{Q}}\boldsymbol{\dot{x}}_{d}\left(t\right) & =\boldsymbol{\mathcal{Q}}\boldsymbol{R}\left(\boldsymbol{x}_{d}\left(t\right)\right).\label{eq:Eq140_1}
\end{align}
Writing down only the non-vanishing component of Eq. \eqref{eq:Eq140_1}
yields 
\begin{align}
\dot{x}_{d}\left(t\right) & =a_{1}x_{d}\left(t\right)+a_{2}y_{d}\left(t\right)+a_{0}.\label{eq:Eq270}
\end{align}
This linear differential equation for $x_{d}\left(t\right)$ with
an inhomogeneity is readily solved in terms of the desired activator
over time $y_{d}\left(t\right)$,
\begin{align}
x_{d}\left(t\right) & =\frac{a_{0}}{a_{1}}\left(e^{a_{1}\left(t-t_{0}\right)}-1\right)+e^{a_{1}\left(t-t_{0}\right)}x_{d}\left(t_{0}\right)+a_{2}\intop_{t_{0}}^{t}d\tau e^{a_{1}\left(t-\tau\right)}y_{d}\left(\tau\right).\label{eq:Eq140}
\end{align}
For the desired trajectory $\boldsymbol{x}_{d}\left(t\right)$ to
be exactly realizable, it must agree with the initial state $\boldsymbol{x}\left(t_{0}\right)$
of the controlled system. The control is given by
\begin{align}
u\left(t\right) & =\dot{y}_{d}\left(t\right)-R\left(x_{d}\left(t\right),y_{d}\left(t\right)\right)=\dot{y}_{d}\left(t\right)-y_{d}\left(t\right)+\frac{1}{3}y_{d}\left(t\right)^{3}+x_{d}\left(t\right).
\end{align}
Using the solution Eq. \eqref{eq:Eq140}, the inhibitor variable $x_{d}\left(t\right)$
can be eliminated from the control signal. Consequently, the control
can be expressed as a functional of the desired activator variable
$y_{d}\left(t\right)$ and the initial desired inhibitor value $x_{d}\left(t_{0}\right)$
as
\begin{align}
u\left(t\right) & =\dot{y}_{d}\left(t\right)-y_{d}\left(t\right)+\frac{1}{3}y_{d}\left(t\right)^{3}+\frac{a_{0}}{a_{1}}\left(e^{a_{1}\left(t-t_{0}\right)}-1\right)\nonumber \\
 & +e^{a_{1}\left(t-t_{0}\right)}x_{d}\left(t_{0}\right)+a_{2}\intop_{t_{0}}^{t}d\tau e^{a_{1}\left(t-\tau\right)}y_{d}\left(\tau\right).\label{eq:UOverTimeFHN}
\end{align}
To evaluate the performance of the control, the control signal Eq.
\eqref{eq:UOverTimeFHN} is used in Eq. \eqref{eq:Eq267}, and the
resulting controlled dynamical system is solved numerically. The numerically
obtained state trajectory is compared with the desired reference trajectory.
The desired trajectory is chosen as
\begin{align}
y_{d}\left(t\right) & =\sin\left(20t\right)\cos\left(2t\right),\label{eq:ActivatorTimeEvolution}
\end{align}
and the initial conditions are set to $x\left(t_{0}\right)=x_{d}\left(t_{0}\right)=y\left(t_{0}\right)=y_{d}\left(t_{0}\right)=0$.
As expected from Eq. \eqref{eq:ActivatorTimeEvolution}, the controlled
activator $y\left(t\right)$ oscillates wildly, see blue solid line
in Fig. \ref{fig:ControlledFHNa1} left. The numerically obtained
controlled inhibitor $x\left(t\right)$ (red dashed line) increases
almost linearly. This behavior can easily be understood from the smallness
of $a_{1}$ in Eq. \eqref{eq:Eq140} (see Example \ref{ex:FHN1} for
parameter values). Indeed, in the limit of vanishing $a_{1}$,
\begin{align}
\lim_{a_{1}\rightarrow0}x_{d}\left(t\right) & =a_{0}\left(t-t_{0}\right)+x_{d}\left(t_{0}\right)+a_{2}\intop_{t_{0}}^{t}d\tau y_{d}\left(\tau\right),
\end{align}
$x_{d}$ increases linearly in time with coefficient $a_{0}$, while
the integral term over a periodic function $y_{d}\left(t\right)$
with zero mean vanishes on average. The control signal $u\left(t\right)$,
being proportional to $\dot{y}\left(t\right)$, oscillates as well,
see Fig. \ref{fig:ControlledFHNa1} right. Comparing the differences
between the controlled state components and its desired counterparts
reveals agreement within numerical precision, see Fig. \ref{fig:ControlledFHNb1}
left for the activator and Fig. \ref{fig:ControlledFHNb1} right for
the inhibitor component, respectively. However, note that the error
increases in time, which could indicate a developing instability.
The control, being an open loop control, is potentially unstable.
It often must be stabilized to guarantee a successful control. Stabilization
of exactly realizable trajectories is not discussed in this thesis.

\begin{minipage}{1.0\linewidth}
\begin{center}
\includegraphics[scale=0.59]{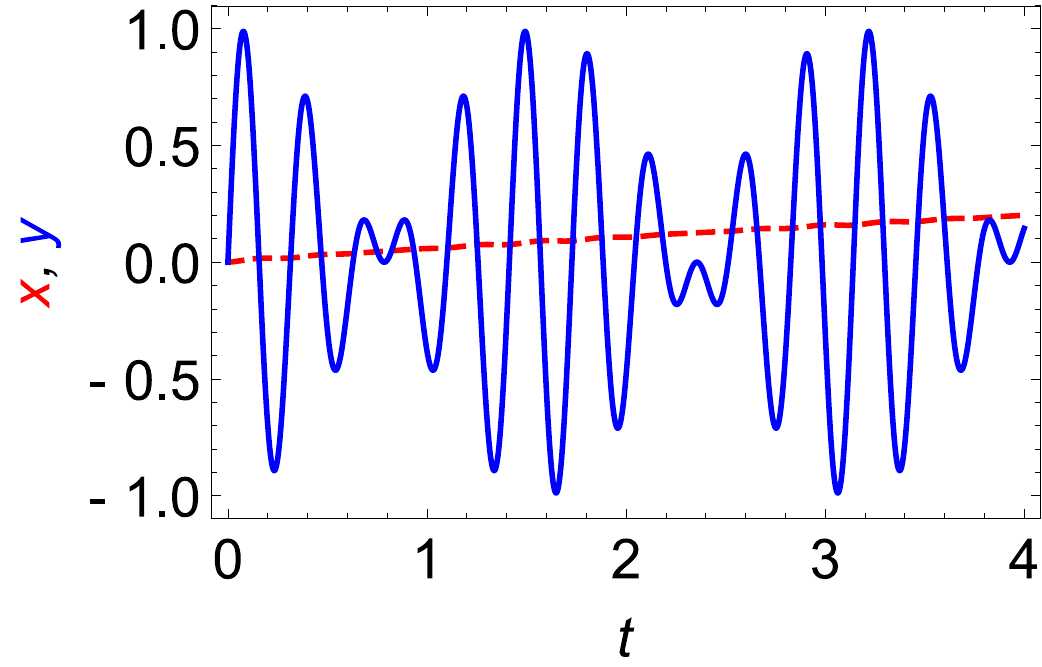}\hspace{1.0cm}\includegraphics[scale=0.59]{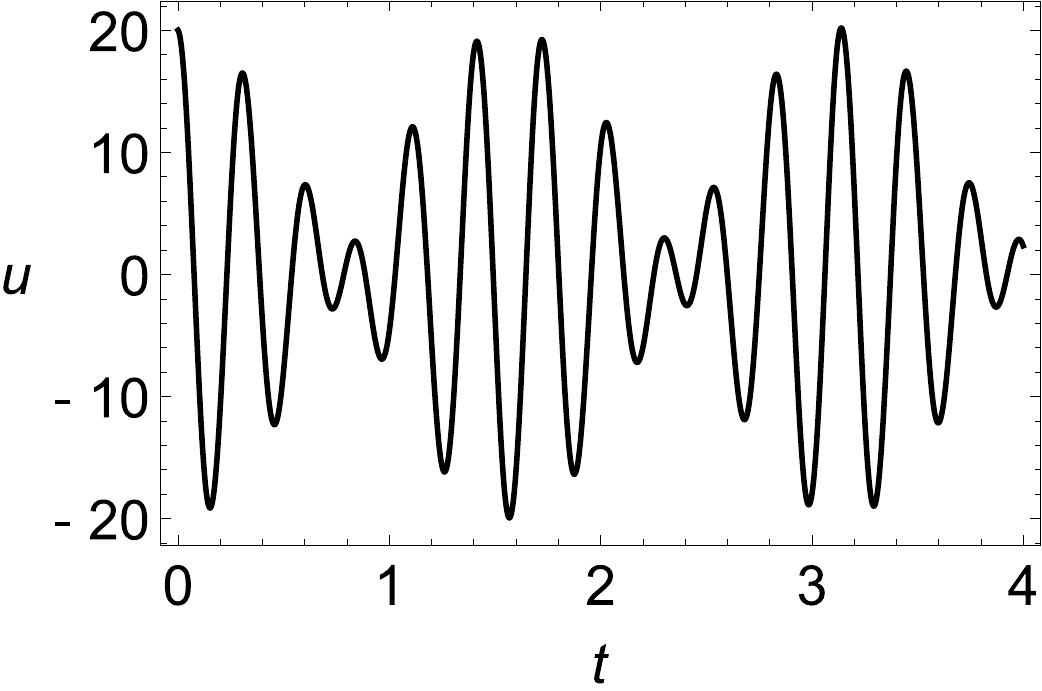}
\captionof{figure}[Exactly realizable trajectory for the activator-controlled FHN  model]{\label{fig:ControlledFHNa1}Activator-controlled FHN model driven along an exactly realizable trajectory. The numerically obtained activator $y$ (blue solid line) and inhibitor $x$ (red dashed line) of the controlled system is shown left. The oscillating activator is prescribed according to Eq. \eqref{eq:ActivatorTimeEvolution}, while the inhibitor cannot be prescribed and is given as the solution to the constraint equation \eqref{eq:Eq270}. The control signal (right) oscillates as well because it is proportional to $\dot{y}_d$.}
%"/home/jakob/svn/Control/FitzHughNagumo/FHNLocalDynamics.nb"
\end{center}
\end{minipage}

\begin{minipage}{1.0\linewidth}
\begin{center}
\includegraphics[scale=0.575]{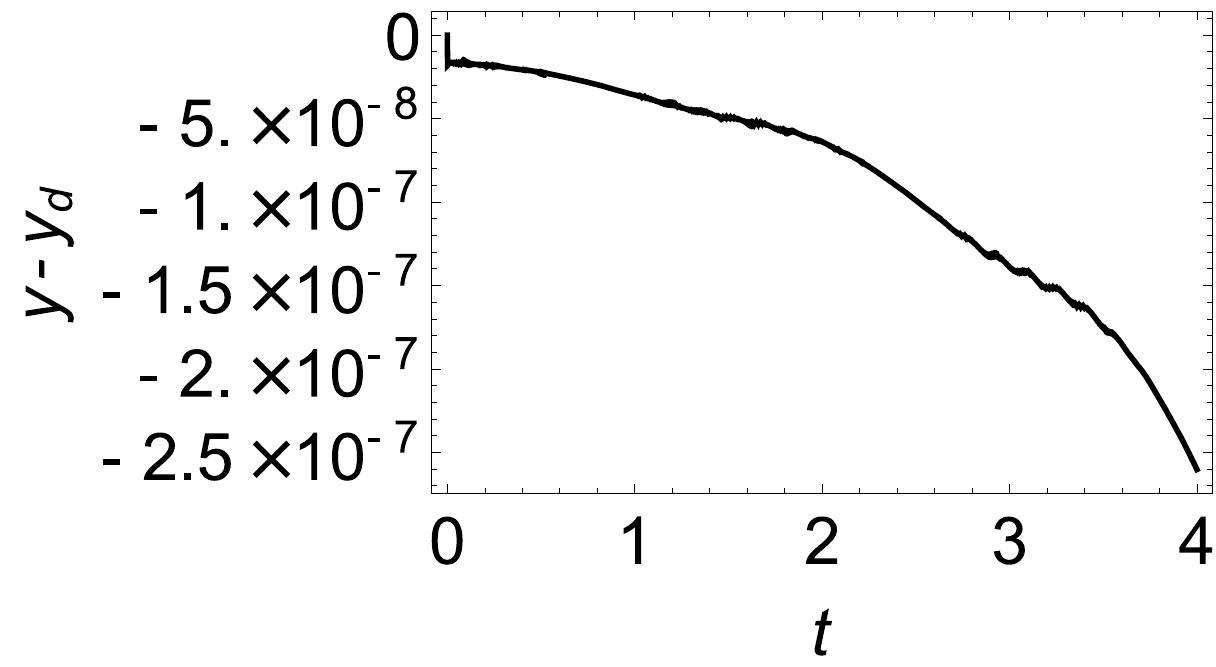}\hspace{0.2cm}\includegraphics[scale=0.55]{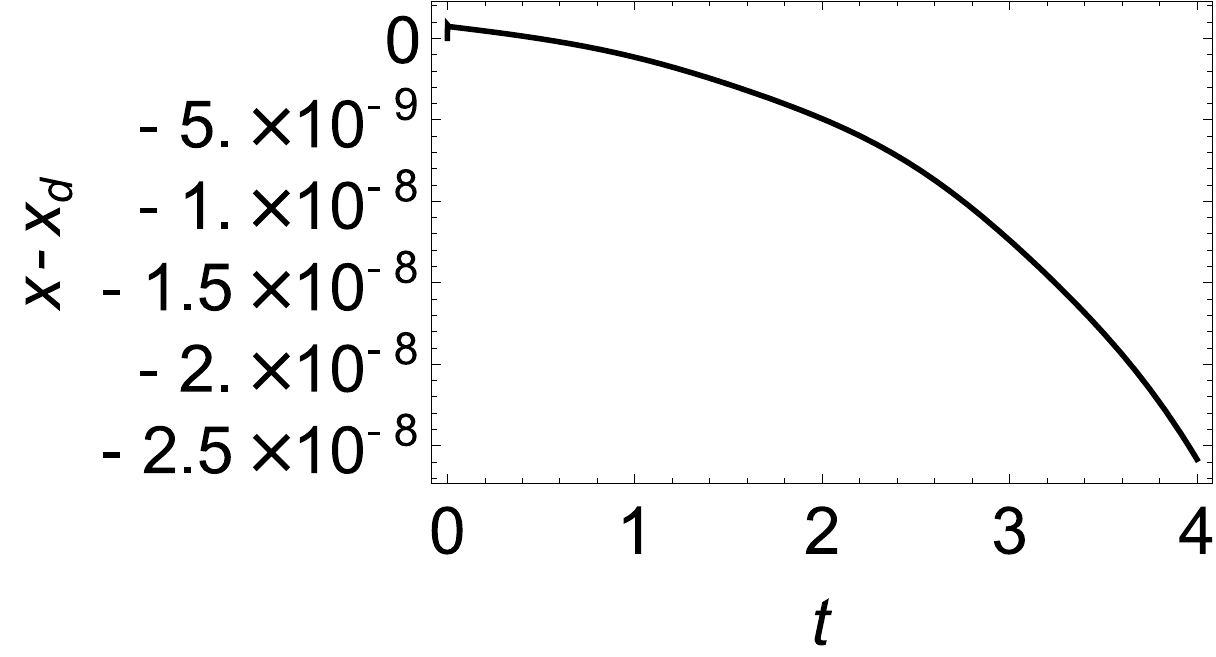}
\captionof{figure}[Difference between desired and controlled trajectory]{\label{fig:ControlledFHNb1}Difference between desired and controlled state components in the activator-controlled FHN model. Plotting the difference between controlled and desired activator $y-y_d$ (left) and controlled and desired inhibitor $x-x_d$ (right) reveals agreement within numerical precision.}
%"/home/jakob/svn/Control/FitzHughNagumo/FHNLocalDynamics.nb"
\end{center}
\end{minipage}

\end{example}

\begin{example}[Inhibitor-controlled FHN model]\label{ex:FHN3_1}

Consider the same model as in Example \ref{ex:FHN3} but with a coupling
vector $\boldsymbol{B}=\left(\begin{array}{cc}
1, & 0\end{array}\right)^{T}$ corresponding to a control acting on the inhibitor equation (see
also Example \ref{ex:FHN1}).

The projectors are $\boldsymbol{\mathcal{P}}=\footnotesize{\left(\begin{array}{cc}
1 & 0\\
0 & 0
\end{array}\right)}$ and $\boldsymbol{\mathcal{Q}}=\footnotesize{\left(\begin{array}{cc}
0 & 0\\
0 & 1
\end{array}\right)}$. The desired inhibitor over time $x_{d}\left(t\right)$ is prescribed
while the activator component $y_{d}\left(t\right)$ must be determined
from (the non-vanishing component of) the constraint equation,
\begin{align}
\dot{y}_{d}\left(t\right) & =R\left(x_{d}\left(t\right),y_{d}\left(t\right)\right)=y_{d}\left(t\right)-\frac{1}{3}y_{d}\left(t\right)^{3}-x_{d}\left(t\right).\label{eq:Eq148}
\end{align}
The constraint equation is a nonlinear non-autonomous differential
equation for $y_{d}\left(t\right)$. An analytical expression for
the solution $y_{d}\left(t\right)$ in terms of the prescribed inhibitor
trajectory $x_{d}\left(t\right)$ is not available. Equation \eqref{eq:Eq148}
must be solved numerically.

Figure \ref{fig:ControlledFHN2} shows the result of a numerical simulation
of the controlled system with a desired inhibitor trajectory
\begin{align}
x_{d}\left(t\right) & =4\sin\left(2t\right),
\end{align}
and initial conditions $x_{d}\left(t_{0}\right)=x\left(t_{0}\right)=y_{d}\left(t_{0}\right)=y\left(t_{0}\right)=0$.
Comparing the desired activator and inhibitor trajectories with the
corresponding controlled state trajectories in the bottom panels demonstrates
a difference in the range of numerical precision over the whole time
interval. Both state components (top left) as well as the control
(top right) are oscillating.

\begin{minipage}{1.0\linewidth}
\begin{center}
\includegraphics[scale=0.62]{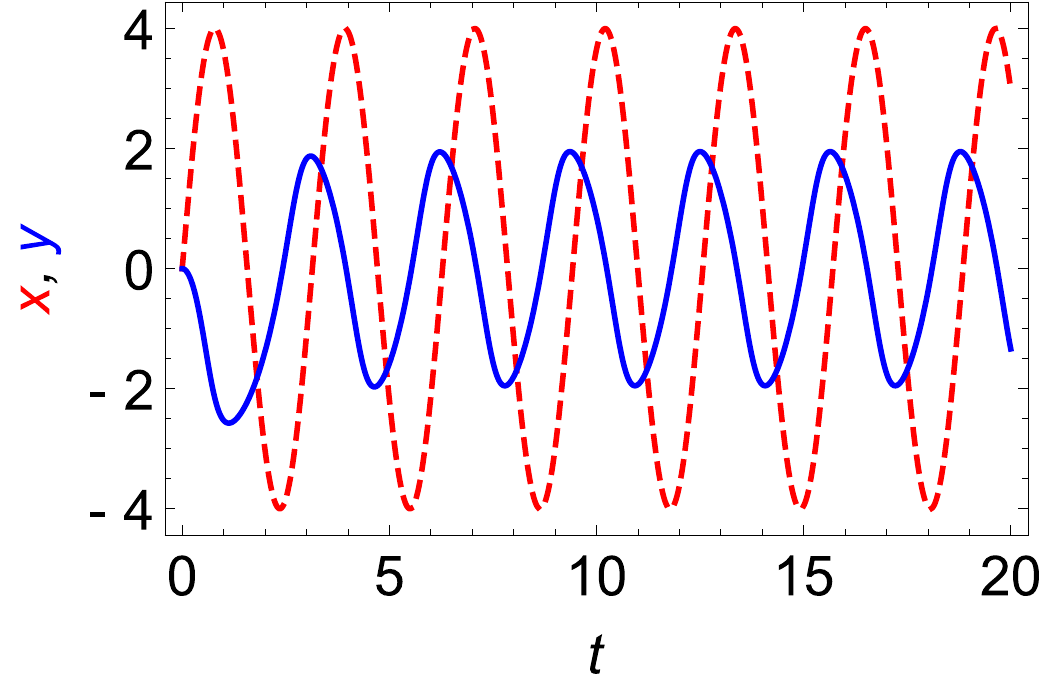}
\includegraphics[scale=0.59]{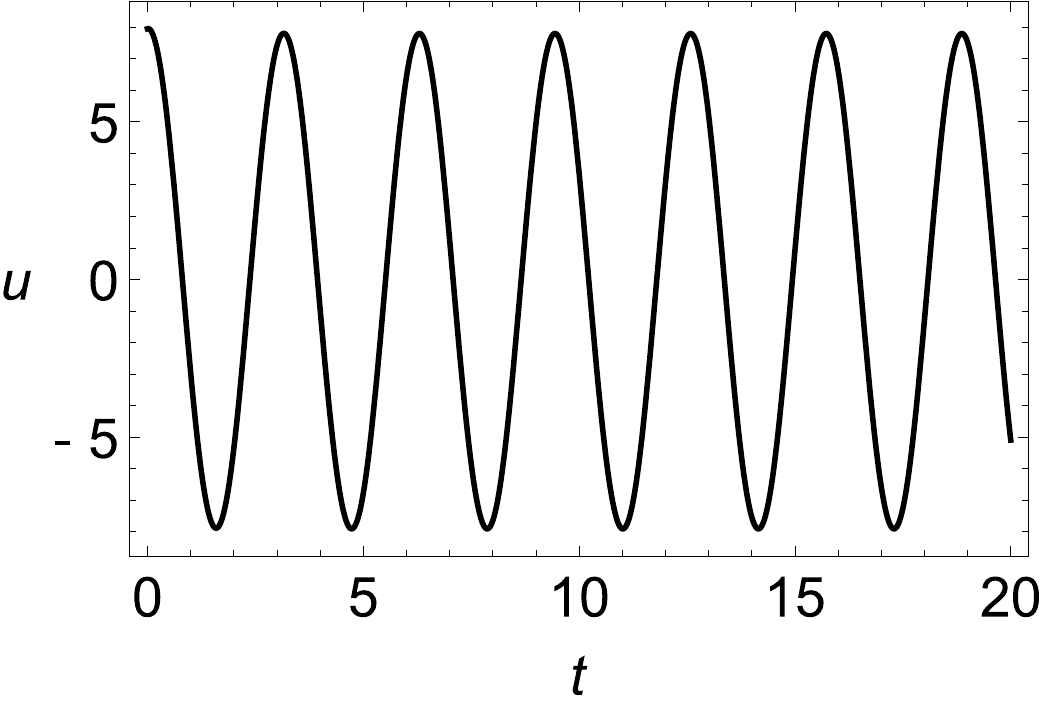}
\includegraphics[scale=0.55]{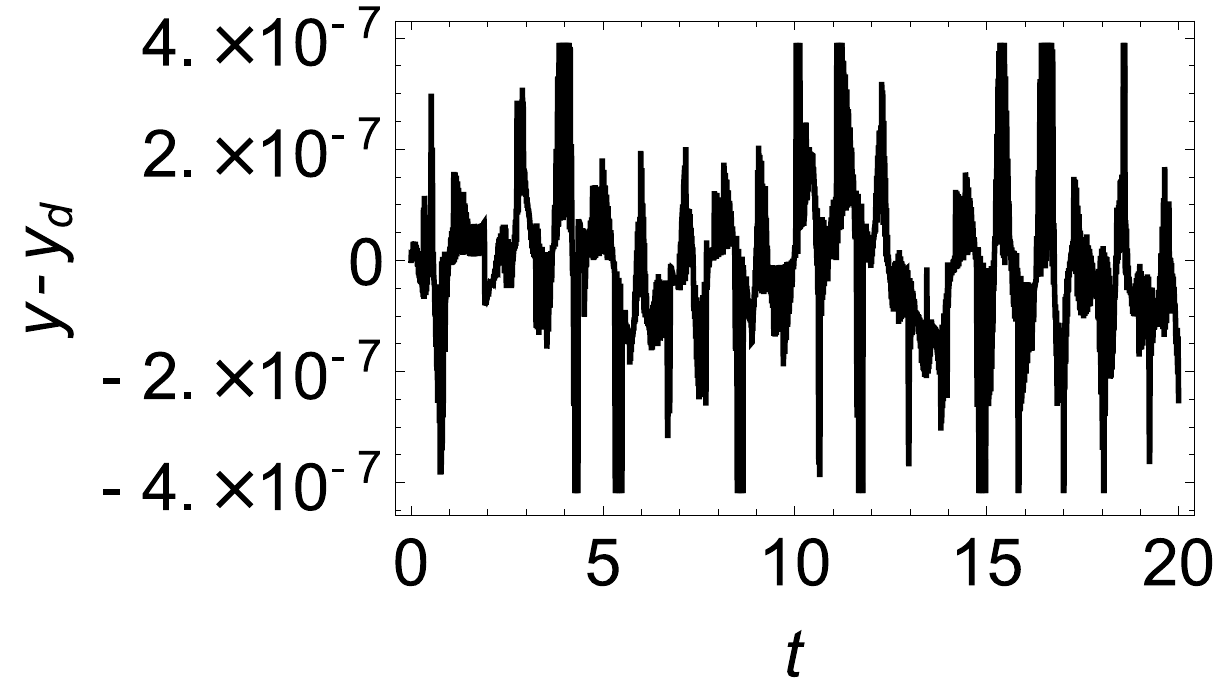}\includegraphics[scale=0.55]{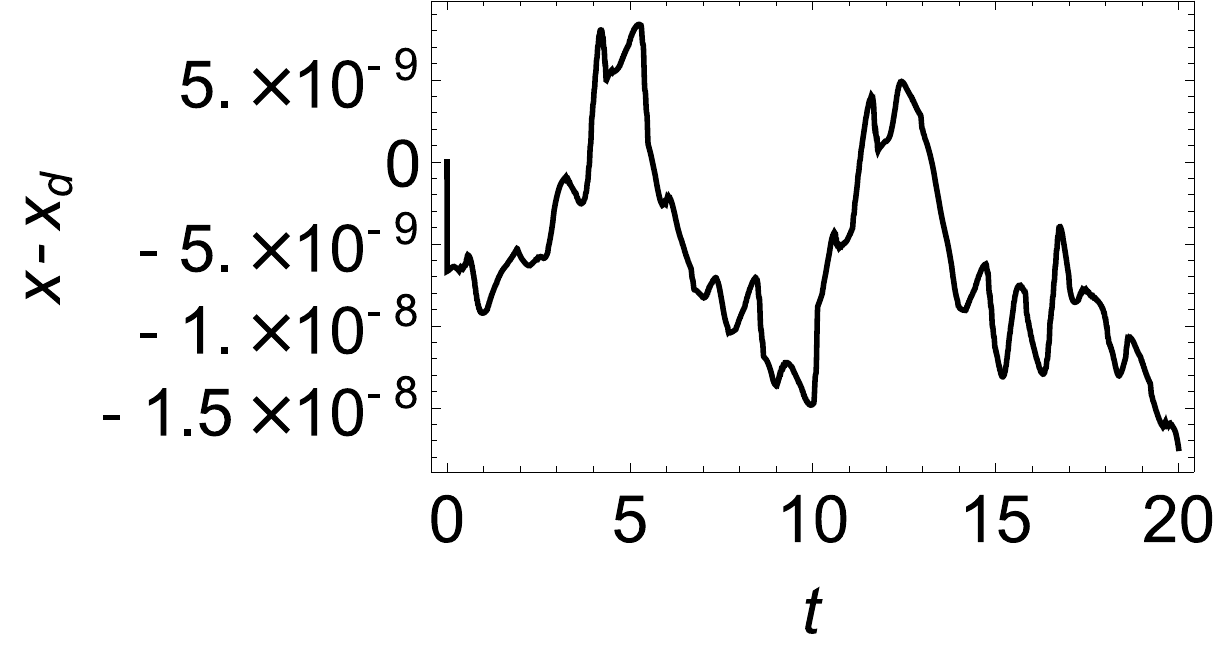}
\captionof{figure}[Exactly realizable trajectory for the inhibitor-controlled FHN model]{\label{fig:ControlledFHN2}Inhibitor-controlled FHN model driven along an exactly realizable trajectory. The numerically obtained solution of the controlled state is shown top left, and the control signal is shown top right. Comparing desired and controlled activator $y_d$ and $y$ (bottom left) as well as desired and controlled inhibitor $x_d$ and $x$ (bottom left) reveals a difference within numerical precision.}
%"/home/jakob/svn/Control/FitzHughNagumo/FHNLocalDynamics.nb"
\end{center}
\end{minipage}

\end{example}

\section{\label{sec:LinearizingAssumption}Linearizing assumption}

An uncontrolled dynamical system requires solving
\begin{align}
\boldsymbol{\dot{x}}\left(t\right) & =\boldsymbol{R}\left(\boldsymbol{x}\left(t\right)\right).
\end{align}
In contrast, control of exactly realizable trajectories requires only
the solution of the constraint equation
\begin{align}
\boldsymbol{\mathcal{Q}}\left(\boldsymbol{x}\left(t\right)\right)\left(\boldsymbol{\dot{x}}\left(t\right)-\boldsymbol{R}\left(\boldsymbol{x}\left(t\right)\right)\right) & =\mathbf{0}.\label{eq:ConstraintEq}
\end{align}
This opens up the possibility to solve a nonlinear control problem
without actually solving any nonlinear equations. If the constraint
equation is linear, the entire controlled system can be regarded,
in some sense and to some extent, as being linear. Two conditions
must be met for Eq. \eqref{eq:ConstraintEq} to be linear. First of
all, the projection matrices $\boldsymbol{\mathcal{P}}\left(\boldsymbol{x}\right)$
and $\boldsymbol{\mathcal{Q}}\left(\boldsymbol{x}\right)$ should
be independent of the state $\boldsymbol{x}$. This condition can
be expressed as
\begin{align}
\boldsymbol{\mathcal{P}}\left(\boldsymbol{x}\right) & =\boldsymbol{1}-\boldsymbol{\mathcal{Q}}\left(\boldsymbol{x}\right)=\boldsymbol{\mathcal{B}}\left(\boldsymbol{x}\right)\left(\boldsymbol{\mathcal{B}}^{T}\left(\boldsymbol{x}\right)\boldsymbol{\mathcal{B}}\left(\boldsymbol{x}\right)\right)^{-1}\boldsymbol{\mathcal{B}}^{T}\left(\boldsymbol{x}\right)=\text{const.}\label{eq:LinearizingAssumption1}
\end{align}
or
\begin{align}
\nabla\left(\boldsymbol{\mathcal{B}}\left(\boldsymbol{x}\right)\left(\boldsymbol{\mathcal{B}}^{T}\left(\boldsymbol{x}\right)\boldsymbol{\mathcal{B}}\left(\boldsymbol{x}\right)\right)^{-1}\boldsymbol{\mathcal{B}}^{T}\left(\boldsymbol{x}\right)\right) & =\boldsymbol{0}.
\end{align}
Note that this condition does not imply that the coupling matrix is
independent of $\boldsymbol{x}$. Second, the nonlinearity $\boldsymbol{R}\left(\boldsymbol{x}\right)$
must satisfy
\begin{align}
\boldsymbol{\mathcal{Q}}\boldsymbol{R}\left(\boldsymbol{x}\right) & =\boldsymbol{\mathcal{Q}}\boldsymbol{\mathcal{A}}\boldsymbol{x}+\boldsymbol{\mathcal{Q}}\boldsymbol{b},\label{eq:LinearizingAssumption2}
\end{align}
with $n\times n$ matrix $\boldsymbol{\mathcal{A}}$ and $n$-component
vector $\boldsymbol{b}$ independent of the state $\boldsymbol{x}$.
Strictly speaking, the projector $\boldsymbol{\mathcal{Q}}$ in front
of $\boldsymbol{\mathcal{A}}$ and $\boldsymbol{b}$ is not really
necessary. It is placed there to make it clear that $\boldsymbol{\mathcal{A}}$
and $\boldsymbol{b}$ do not contain any parts in the direction of
$\boldsymbol{\mathcal{P}}$. Condition Eq. \eqref{eq:LinearizingAssumption1}
combined with condition Eq. \eqref{eq:LinearizingAssumption2} constitute
the\textit{ linearizing assumption}. Control systems satisfying the
linearizing assumption behave, to a large extent, similar to truly
linear control systems. A nonlinear control systems with scalar input
$u\left(t\right)$ satisfying the linearizing assumption is sometimes
said to be in companion form. A system in companion form is trivially
feedback linearizable, see the discussion of feedback linearization
in Chapter \ref{chap:Introduction} and e.g. \cite{khalil2002nonlinear}.

Condition Eq. \eqref{eq:LinearizingAssumption2} is a strong assumption.
It enforces $n-p$ components of $\boldsymbol{R}\left(\boldsymbol{x}\right)$
to depend only linearly on the state. However, some important models
of nonlinear dynamics satisfy the linearizing assumption. Among these
are the mechanical control systems in one spatial dimension, see Examples
\ref{ex:OneDimMechSys1} and \ref{ex:OneDimMechSys3}, as well as
the activator-controlled FHN model discussed in Example \ref{ex:FHN3}.
In both cases, the control signal $\boldsymbol{u}\left(t\right)$
acts directly on the nonlinear part of the nonlinearity $\boldsymbol{R}$,
such that condition Eq. \eqref{eq:LinearizingAssumption2} is satisfied.
Furthermore, in both cases the coupling matrix $\boldsymbol{\mathcal{B}}\left(\boldsymbol{x}\right)$
is a coupling vector $\boldsymbol{B}\left(\boldsymbol{x}\right)=\left(\begin{array}{cc}
0, & B\left(x,y\right)\end{array}\right)^{T}$ with only one non-vanishing component. This leads to constant projectors
$\boldsymbol{\mathcal{P}}$ and $\boldsymbol{\mathcal{Q}}$, and condition
Eq. \eqref{eq:LinearizingAssumption1} is also satisfied. Another,
less obvious example satisfying the linearizing assumption is the
controlled SIR model.

\begin{example}[Linearizing assumption satisfied by the controlled SIR model]\label{ex:SIRModel1_1}

The controlled state equation for the SIR model was developed in Example
\ref{ex:SIRModel1}. The nonlinearity $\boldsymbol{R}$ is
\begin{align}
\boldsymbol{R}\left(\boldsymbol{x}\left(t\right)\right) & =\left(-\beta\frac{S\left(t\right)I\left(t\right)}{N},\beta\frac{S\left(t\right)I\left(t\right)}{N}-\gamma I\left(t\right),\gamma I\left(t\right)\right)^{T},
\end{align}
while the coupling vector $\boldsymbol{B}$ explicitly depends on
the state, 
\begin{align}
\boldsymbol{B}\left(\boldsymbol{x}\left(t\right)\right) & =\frac{1}{N}\left(-S\left(t\right)I\left(t\right),S\left(t\right)I\left(t\right),0\right)^{T}.
\end{align}
However, the projectors
\begin{align}
\boldsymbol{\mathcal{P}}\left(\boldsymbol{x}\right) & =\boldsymbol{\mathcal{P}}=\boldsymbol{B}\left(\boldsymbol{x}\right)\left(\boldsymbol{B}^{T}\left(\boldsymbol{x}\right)\boldsymbol{B}\left(\boldsymbol{x}\right)\right)^{-1}\boldsymbol{B}^{T}\left(\boldsymbol{x}\right)=\frac{1}{2}\left(\begin{array}{ccc}
1 & -1 & 0\\
-1 & 1 & 0\\
0 & 0 & 0
\end{array}\right),\\
\boldsymbol{\mathcal{Q}}\left(\boldsymbol{x}\right) & =\boldsymbol{\mathcal{Q}}=\frac{1}{2}\left(\begin{array}{ccc}
1 & 1 & 0\\
1 & 1 & 0\\
0 & 0 & 2
\end{array}\right),
\end{align}
are independent of the state. Furthermore, the model also satisfies
the linearizing assumption Eq. \eqref{eq:LinearizingAssumption1}
because
\begin{align}
\boldsymbol{\mathcal{Q}}\boldsymbol{R}\left(\boldsymbol{x}\right) & =\left(\begin{array}{ccc}
0 & -\dfrac{\gamma}{2} & 0\\
0 & -\dfrac{\gamma}{2} & 0\\
0 & \gamma & 0
\end{array}\right)=\boldsymbol{\mathcal{Q}}\boldsymbol{\mathcal{A}}\boldsymbol{x}.
\end{align}
The constraint equation is a linear differential equation with three
components,
\begin{align}
\left(\begin{array}{c}
\frac{1}{2}\left(\gamma I_{d}\left(t\right)+\dot{I}_{d}\left(t\right)+\dot{S}_{d}\left(t\right)\right)\\
\frac{1}{2}\left(\gamma I_{d}\left(t\right)+\dot{I}_{d}\left(t\right)+\dot{S}_{d}\left(t\right)\right)\\
-\gamma I_{d}\left(t\right)+\dot{R}_{d}\left(t\right)
\end{array}\right) & =\left(\begin{array}{c}
0\\
0\\
0
\end{array}\right),
\end{align}
of which one is redundant. Note that because the projectors $\boldsymbol{\mathcal{P}}$
and $\boldsymbol{\mathcal{Q}}$ are non-diagonal, the time derivatives
of $I_{d}\left(t\right)$ and $S_{d}\left(t\right)$ are mixed in
the constraint equation.

\end{example}

\section{\label{sec:Controllability}Controllability}

A system is called controllable or state controllable if it is possible
to achieve a transfer from an initial state $\boldsymbol{x}\left(t_{0}\right)=\boldsymbol{x}_{0}$
at time $t=t_{0}$ to a final state $\boldsymbol{x}\left(t_{1}\right)=\boldsymbol{x}_{1}$
at the terminal time $t=t_{1}$. Controllability is a condition on
the structure of the dynamical system as given by the nonlinearity
$\boldsymbol{R}$ and the coupling matrix $\boldsymbol{\mathcal{B}}$.
In contrast, for a given control system, trajectory realizability
is a condition on the desired trajectory. While for linear control
systems controllability is easily expressed in terms of a rank condition,
the notion is much more difficult for nonlinear control systems. Section
\ref{sec:StateToStateControllability} discusses the Kalman rank condition
for the controllability of LTI systems as introduced by Kalman \cite{kalman1959general,kalman1960contributions}
in the early sixties. Section \ref{sec:ControllabilityForRealizableTrajectories}
derives a similar rank condition in the context of exactly realizable
trajectories. Remarkably, this rank condition also applies to nonlinear
systems satisfying the linearizing assumption from Section \ref{sec:LinearizingAssumption}.

For general nonlinear systems, the notion of controllability must
be refined and it is necessary to distinguish between controllability,
accessibility, and reachability. Different and not necessarily equivalent
notions of controllability exist, and it is said that there are as
many notions of nonlinear controllability as there are researchers
in the field. When applied to LTI systems, all of these notions reduce
to the Kalman rank condition. Here, no attempt is given to generalize
the notion of controllability to nonlinear systems which violate the
linearizing assumption. The reader is referred to the literature \cite{slotine1991applied,isidori1995nonlinear,khalil2002nonlinear,levine2009analysis}.

\subsection{\label{sec:StateToStateControllability}Kalman rank condition for
LTI systems}

Controllability for LTI systems was first introduced by Kalman \cite{kalman1959general,kalman1960contributions}.
An excellent introduction to linear control systems, including controllability,
can be found in \cite{chen1995linear}.

Consider the LTI system with $n$-dimensional state vector $\boldsymbol{x}\left(t\right)$
and $p$-dimensional control signal $\boldsymbol{u}\left(t\right)$,
\begin{align}
\boldsymbol{\dot{x}}\left(t\right) & =\boldsymbol{\mathcal{A}}\boldsymbol{x}\left(t\right)+\boldsymbol{\mathcal{B}}\boldsymbol{u}\left(t\right),\label{eq:LTIStateEquation}
\end{align}
and initial condition
\begin{align}
\boldsymbol{x}\left(t_{0}\right) & =\boldsymbol{x}_{0}.\label{eq:LTIStateInitCond}
\end{align}
Here, $\boldsymbol{\mathcal{A}}$ is an $n\times n$ real constant
matrix and $\boldsymbol{\mathcal{B}}$ an $n\times p$ real constant
matrix. The system Eq. \eqref{eq:LTIStateEquation} is said to be
controllable if, for any initial state $\boldsymbol{x}_{0}$ at the
initial time $t=t_{0}$ and any final state $\boldsymbol{x}_{1}$
at the terminal time $t=t_{1}$, there exists an input that transfers
$\boldsymbol{x}_{0}$ to $\boldsymbol{x}_{1}$. The terminal condition
for the state is
\begin{align}
\boldsymbol{x}\left(t_{1}\right) & =\boldsymbol{x}_{1}.\label{eq:LTIStateTerminalCond}
\end{align}
The definition of controllability requires only that the input $\boldsymbol{u}\left(t\right)$
be capable of moving any state in the state space to any other state
in finite time. The state trajectory $\boldsymbol{x}\left(t\right)$
traced out in state space is not specified. Kalman showed that this
definition of controllability is equivalent to the statement that
the $n\times np$ \textit{controllability matrix} 
\begin{align}
\boldsymbol{\mathcal{K}}= & \left(\boldsymbol{\mathcal{B}}|\boldsymbol{\mathcal{A}}\boldsymbol{\mathcal{B}}|\boldsymbol{\mathcal{A}}^{2}\boldsymbol{\mathcal{B}}|\cdots|\boldsymbol{\mathcal{A}}^{n-1}\boldsymbol{\mathcal{B}}\right)\label{eq:KalmanControllabilityMatrix}
\end{align}
has rank $n$, i.e., it satisfies the \textit{Kalman rank condition}
\begin{align}
\text{rank}\left(\boldsymbol{\mathcal{K}}\right) & =n.\label{eq:KalmanRankCondition}
\end{align}
Since $n\leq np$, this condition states that $\boldsymbol{\mathcal{K}}$
has full row rank. Equation \eqref{eq:KalmanRankCondition} is derived
in the following.

\subsection{\label{sub:DerivationOfKalmanRank}Derivation of the Kalman rank
condition}

The solution $\boldsymbol{x}\left(t\right)$ to Eq. \eqref{eq:LTIStateEquation}
with initial condition Eq. \eqref{eq:LTIStateInitCond} and arbitrary
control signal $\boldsymbol{u}\left(t\right)$ is
\begin{align}
\boldsymbol{x}\left(t\right) & =e^{\boldsymbol{\mathcal{A}}\left(t-t_{0}\right)}\boldsymbol{x}_{0}+\intop_{t_{0}}^{t}d\tau e^{\boldsymbol{\mathcal{A}}\left(t-\tau\right)}\boldsymbol{\mathcal{B}}\boldsymbol{u}\left(\tau\right).\label{eq:LTISolution}
\end{align}
See also Appendix \ref{sec:GeneralSolutionForForcedLinarDynamicalSystem}
for a derivation of the general solution to a forced linear dynamical
system. The system is controllable if a control signal $\boldsymbol{u}$
can be found such that the terminal condition \eqref{eq:LTIStateTerminalCond}
is satisfied. Evaluating Eq. \eqref{eq:LTISolution} at the terminal
time $t=t_{1}$, multiplying by $e^{-\boldsymbol{\mathcal{A}}\left(t_{1}-t_{0}\right)}$,
rearranging, and expanding the matrix exponential under the integral
yields 
\begin{align}
e^{-\boldsymbol{\mathcal{A}}\left(t_{1}-t_{0}\right)}\boldsymbol{x}_{1}-\boldsymbol{x}_{0} & =\sum_{k=0}^{\infty}\boldsymbol{\mathcal{A}}^{k}\boldsymbol{\mathcal{B}}\intop_{t_{0}}^{t_{1}}d\tau\frac{\left(t_{0}-\tau\right)^{k}}{k!}\boldsymbol{u}\left(\tau\right).\label{eq:Eq206}
\end{align}
For the system to be controllable, it must in principle be possible
to solve for the control signal $\boldsymbol{u}$. As a consequence
of the Cayley-Hamilton theorem, the matrix power $\boldsymbol{\mathcal{A}}^{i}$
for any $n\times n$ matrix with $i\geq n$ can be written as a sum
of lower order powers \cite{fischer2008lineare},
\begin{align}
\boldsymbol{\mathcal{A}}^{i} & =\sum_{k=0}^{n-1}c_{ik}\boldsymbol{\mathcal{A}}^{k}.
\end{align}
It follows that the infinite sum in Eq. \eqref{eq:Eq206} can be rearranged
to include only terms with power in $\boldsymbol{\mathcal{A}}$ up
to $\boldsymbol{\mathcal{A}}^{n-1}$. The sum on the right hand side
(r. h. s.) of Eq. \eqref{eq:Eq206} can be simplified as
\begin{align}
 & \sum_{k=0}^{\infty}\boldsymbol{\mathcal{A}}^{k}\boldsymbol{\mathcal{B}}\intop_{t_{0}}^{t_{1}}d\tau\frac{\left(t_{0}-\tau\right)^{k}}{k!}\boldsymbol{u}\left(\tau\right)\nonumber \\
= & \sum_{k=0}^{n-1}\boldsymbol{\mathcal{A}}^{k}\boldsymbol{\mathcal{B}}\intop_{t_{0}}^{t_{1}}d\tau\frac{\left(t_{0}-\tau\right)^{k}}{k!}\boldsymbol{u}\left(\tau\right)+\sum_{i=n}^{\infty}\boldsymbol{\mathcal{A}}^{i}\boldsymbol{\mathcal{B}}\intop_{t_{0}}^{t_{1}}d\tau\frac{\left(t_{0}-\tau\right)^{i}}{i!}\boldsymbol{u}\left(\tau\right)\nonumber \\
= & \sum_{k=0}^{n-1}\boldsymbol{\mathcal{A}}^{k}\boldsymbol{\mathcal{B}}\intop_{t_{0}}^{t_{1}}d\tau\frac{\left(t_{0}-\tau\right)^{k}}{k!}\boldsymbol{u}\left(\tau\right)+\sum_{i=n}^{\infty}\sum_{k=0}^{n-1}c_{ik}\boldsymbol{\mathcal{A}}^{k}\boldsymbol{\mathcal{B}}\intop_{t_{0}}^{t_{1}}d\tau\frac{\left(t_{0}-\tau\right)^{i}}{i!}\boldsymbol{u}\left(\tau\right)\nonumber \\
= & \sum_{k=0}^{n-1}\boldsymbol{\mathcal{A}}^{k}\boldsymbol{\mathcal{B}}\intop_{t_{0}}^{t_{1}}d\tau\left(\frac{\left(t_{0}-\tau\right)^{k}}{k!}+\sum_{i=n}^{\infty}c_{ik}\frac{\left(t_{0}-\tau\right)^{i}}{i!}\right)\boldsymbol{u}\left(\tau\right).\label{eq:Eq297}
\end{align}

It follows that the sum in Eq. \eqref{eq:Eq206} can be truncated
after $n$ terms,
\begin{align}
e^{-\boldsymbol{\mathcal{A}}\left(t_{1}-t_{0}\right)}\boldsymbol{x}_{1}-\boldsymbol{x}_{0} & =\sum_{k=0}^{n-1}\boldsymbol{\mathcal{A}}^{k}\boldsymbol{\mathcal{B}}\boldsymbol{\beta}_{k}\left(t_{1},t_{0}\right).\label{eq:Eq208}
\end{align}
The $\boldsymbol{\beta}_{k}$ are $p\times1$ vectors defined as
\begin{align}
\boldsymbol{\beta}_{k}\left(t_{1},t_{0}\right) & =\intop_{t_{0}}^{t_{1}}d\tau\left(\frac{\left(t_{0}-\tau\right)^{k}}{k!}+\sum_{i=n}^{\infty}c_{ik}\frac{\left(t_{0}-\tau\right)^{i}}{i!}\right)\boldsymbol{u}\left(\tau\right),\label{eq:BetaIntegralEquation}
\end{align}
which depend on the initial and terminal time $t_{0}$ and $t_{1}$,
respectively. These vectors are functionals of the control $\boldsymbol{u}$
and depend on the matrix $\boldsymbol{\mathcal{A}}$ through the expansion
coefficients $c_{ik}$. Defining the $np\times1$ vector
\begin{align}
\boldsymbol{\beta}\left(t_{1},t_{0}\right) & =\left(\begin{array}{c}
\boldsymbol{\beta}_{0}\left(t_{1},t_{0}\right)\\
\vdots\\
\boldsymbol{\beta}_{n-1}\left(t_{1},t_{0}\right)
\end{array}\right),
\end{align}
Eq. \eqref{eq:Eq208} can be written in terms of $\boldsymbol{\beta}$
and Kalman's controllability $n\times np$ matrix $\boldsymbol{\mathcal{K}}$,
Eq. \eqref{eq:KalmanControllabilityMatrix}, as 
\begin{align}
e^{-\boldsymbol{\mathcal{A}}\left(t_{1}-t_{0}\right)}\boldsymbol{x}_{1}-\boldsymbol{x}_{0} & =\sum_{k=0}^{n-1}\boldsymbol{\mathcal{A}}^{k}\boldsymbol{\mathcal{B}}\boldsymbol{\beta}_{k}\left(t_{1},t_{0}\right)=\boldsymbol{\mathcal{K}}\boldsymbol{\beta}\left(t_{1},t_{0}\right).\label{eq:Eq211}
\end{align}
Equation \eqref{eq:Eq211} is a linear equation for the vector $\boldsymbol{\beta}\left(t_{1},t_{0}\right)$
with inhomogeneity $e^{-\boldsymbol{\mathcal{A}}\left(t_{1}-t_{0}\right)}\boldsymbol{x}_{1}-\boldsymbol{x}_{0}$.
For the system \eqref{eq:LTIStateEquation} to be controllable, every
state point $\boldsymbol{x}_{1}$ must have a corresponding vector
$\boldsymbol{\beta}\left(t_{1},t_{0}\right)$. In other words, the
linear map from $\boldsymbol{\beta}\left(t_{1},t_{0}\right)$ to $\boldsymbol{x}_{1}$
must be \textit{surjective}. This is the case if and only if the matrix
$\boldsymbol{\mathcal{K}}$ has full row rank \cite{fischer2008lineare},
i.e.,
\begin{align}
\text{rank}\left(\boldsymbol{\mathcal{K}}\right) & =n.\label{eq:KalmanRankCondition-1}
\end{align}
The Kalman rank condition Eq. \eqref{eq:KalmanRankCondition-1} is
a necessary and sufficient condition for the controllability of an
LTI system.

A slightly different way to arrive at the same result is to solve
\eqref{eq:Eq211} for the vector $\boldsymbol{\beta}\left(t_{1},t_{0}\right)$.
A solution in terms of the $n\times n$ matrix $\boldsymbol{\mathcal{K}}\boldsymbol{\mathcal{K}}^{T}$
is (see also Appendix \ref{sec:OverAndUnderdetSysOfEqs} how to solve
an underdetermined system of equations) 
\begin{align}
\boldsymbol{\beta}\left(t_{1},t_{0}\right) & =\boldsymbol{\mathcal{K}}^{T}\left(\boldsymbol{\mathcal{K}}\boldsymbol{\mathcal{K}}^{T}\right)^{-1}\left(e^{-\boldsymbol{\mathcal{A}}\left(t_{1}-t_{0}\right)}\boldsymbol{x}_{1}-\boldsymbol{x}_{0}\right).\label{eq:BetaSolution}
\end{align}
The inverse of $\boldsymbol{\mathcal{K}}\boldsymbol{\mathcal{K}}^{T}$
does exist only if it has full rank, i.e., $\text{rank}\left(\boldsymbol{\mathcal{K}}\boldsymbol{\mathcal{K}}^{T}\right)=n$.
This is the case if and only if the Kalman rank condition $\text{rank}\left(\boldsymbol{\mathcal{K}}\right)=n$
is satisfied. 

Controllability has a number of interesting and important consequences.
Two examples illustrate the concept and highlight one important consequence.

\begin{example}[Single input diagonal LTI system]\label{ex:DiagonalLTI1}

We consider an LTI system with state and input matrix
\begin{align}
\boldsymbol{\mathcal{A}} & =\left(\begin{array}{cc}
\lambda_{1} & 0\\
0 & \lambda_{2}
\end{array}\right), & \boldsymbol{\mathcal{B}} & =\left(\begin{array}{c}
1\\
1
\end{array}\right).\label{eq:E110}
\end{align}
Kalman's controllability matrix is
\begin{align}
\boldsymbol{\mathcal{K}} & =\left(\boldsymbol{\mathcal{B}}|\boldsymbol{\mathcal{A}}\boldsymbol{\mathcal{B}}\right)=\left(\begin{array}{cc}
1 & \lambda_{1}\\
1 & \lambda_{2}
\end{array}\right).
\end{align}
As long as $\lambda_{1}\neq\lambda_{2}$, $\boldsymbol{\mathcal{K}}$
has rank 2. If $\lambda_{1}=\lambda_{2}$, the second row equals the
first row, and $\boldsymbol{\mathcal{K}}$ has rank 1. The system
Eq. \eqref{eq:E110} is controllable as long as $\lambda_{1}\neq\lambda_{2}$.

\end{example}

\newpage

\begin{example}[Two pendulums mounted on a cart]\label{ex:TwoPendulumsMountedOnACart}

\begin{minipage}{1.0\linewidth}
\begin{center}
\includegraphics[scale=0.65]{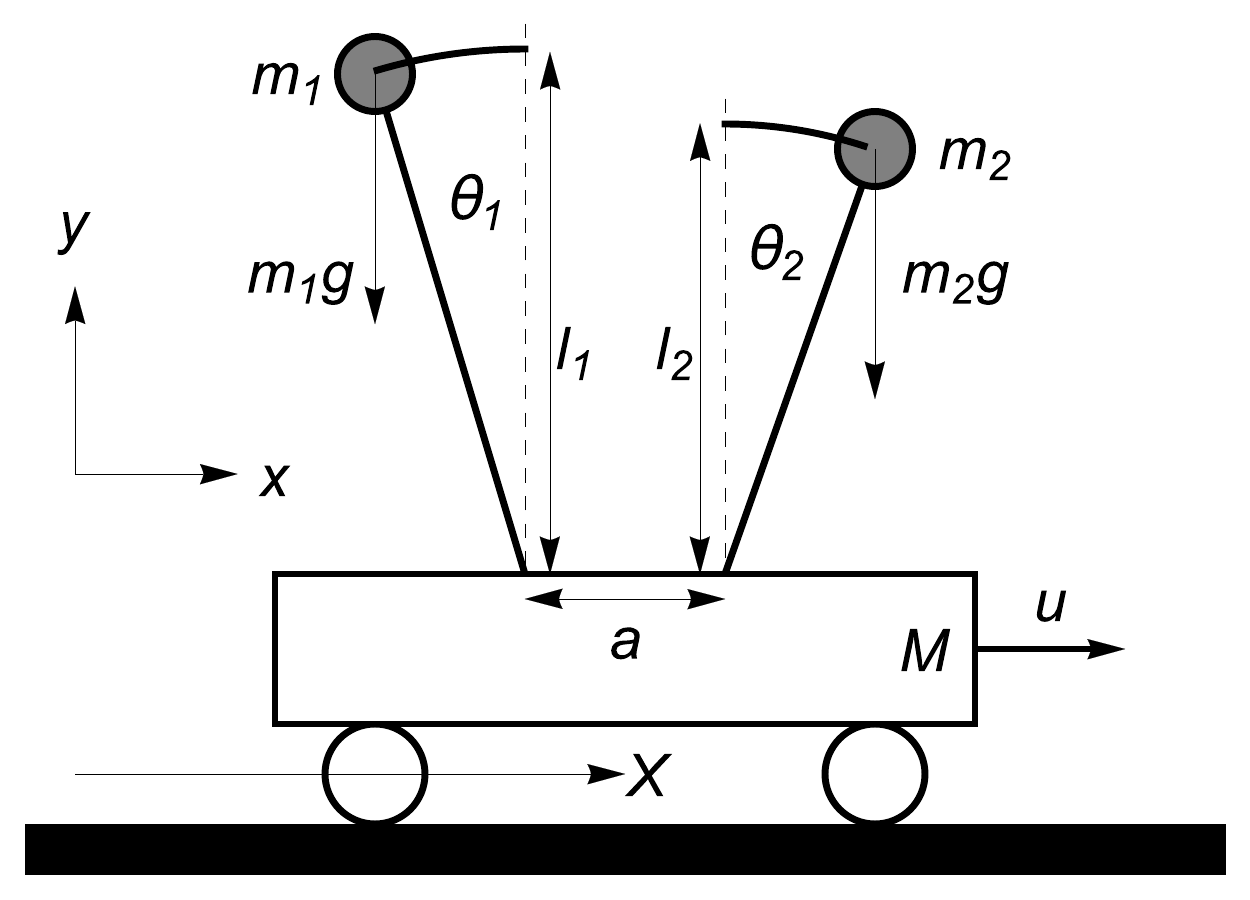}
\captionof{figure}[Two inverted pendulums mounted on a cart]{\label{fig:DoublePendulumOnCart}Two inverted pendulums mounted on a cart. The control task is to keep both pendulums in the upright and unstable equilibrium position. The system is controllable as long as the pendulums are not exactly identical, i.e., as long as either their lengths ($l_1\neq l_2$) or their masses are different ($m_1\neq m_2$).}
%"/home/jakob/svnco/Control/DoublePendulumOnCart/DoublePendulumOnCart.nb"
\end{center}
\end{minipage}

Two pendulums mounted on a cart is a mechanical toy model for linear
control systems \cite{chen1995linear}. As can be seen in Fig. \ref{fig:DoublePendulumOnCart},
the positions of the masses $m_{1}$ and $m_{2}$ given by $\left(\begin{array}{cc}
x_{1}, & y_{1}\end{array}\right)^{T}$ and $\left(\begin{array}{cc}
x_{2}, & y_{2}\end{array}\right)^{T}$, respectively, are 
\begin{align}
x_{1}\left(t\right) & =X\left(t\right)-\frac{a}{2}+l_{1}\sin\left(\theta_{1}\left(t\right)\right), & x_{2}\left(t\right) & =X\left(t\right)+\frac{a}{2}+l_{2}\sin\left(\theta_{2}\left(t\right)\right),\\
y_{1}\left(t\right) & =l_{1}\cos\left(\theta_{1}\left(t\right)\right), & y_{2}\left(t\right) & =l_{2}\cos\left(\theta_{2}\left(t\right)\right).
\end{align}
The cart can only move in the $x$ direction without any motion in
the $y$-direction. Its position is denoted by $X\left(t\right)$.
The Lagrangian $L$ equals the difference between kinetic energy $T$
and potential energy $V$,
\begin{align}
L & =T-V=\frac{1}{2}m_{1}\left(\dot{x}_{1}^{2}\left(t\right)+\dot{y}_{1}^{2}\left(t\right)\right)+\frac{1}{2}m_{2}\left(\dot{x}_{2}^{2}\left(t\right)+\dot{y}_{2}^{2}\left(t\right)\right)\nonumber \\
 & +\frac{1}{2}M\dot{X}^{2}\left(t\right)-m_{1}l_{1}g\cos\left(\theta_{1}\left(t\right)\right)-m_{2}l_{2}g\cos\left(\theta_{1}\left(t\right)\right).
\end{align}
The equations of motion are given by the Euler-Lagrange equations
\begin{align}
\frac{d}{dt}\dfrac{\partial L}{\partial\dot{\theta}_{1}}-\dfrac{\partial L}{\partial\theta_{1}} & =0, & \frac{d}{dt}\dfrac{\partial L}{\partial\dot{\theta}_{2}}-\dfrac{\partial L}{\partial\theta_{2}} & =0, & \frac{d}{dt}\dfrac{\partial L}{\partial\dot{X}}-\dfrac{\partial L}{\partial X} & =u\left(t\right).
\end{align}
The control force $u\left(t\right)$ acts on the cart but not on the
pendulums. Assuming small angles, $0\leq\left|\theta_{1}\left(t\right)\right|\ll1$
and $0\leq\left|\theta_{2}\left(t\right)\right|\ll1$, the equations
of motion are linearized around the stationary point. Rewriting the
second order differential equations as a controlled dynamical system
with $P\left(t\right)=M\dot{X}\left(t\right)$, $p_{1}\left(t\right)=\dot{\theta}_{1}$,
and $p_{2}\left(t\right)=\dot{\theta}_{2}\left(t\right)$ yields
\begin{align}
\boldsymbol{\dot{x}}\left(t\right) & =\boldsymbol{\mathcal{A}}\boldsymbol{x}\left(t\right)+\boldsymbol{B}u\left(t\right),
\end{align}
with
\begin{align}
\boldsymbol{x}\left(t\right) & =\left(\theta_{1}\left(t\right),\theta_{2}\left(t\right),X\left(t\right),p_{1}\left(t\right),p_{2}\left(t\right),P\left(t\right)\right)^{T},\\
\boldsymbol{\mathcal{A}} & =\left(\begin{array}{cccccc}
0 & 0 & 0 & 1 & 0 & 0\\
0 & 0 & 0 & 0 & 1 & 0\\
0 & 0 & 0 & 0 & 0 & \frac{1}{M}\\
\frac{g\left(m_{1}+M\right)}{l_{1}M} & \frac{gm_{2}}{l_{1}M} & 0 & 0 & 0 & 0\\
\frac{gm_{1}}{l_{2}M} & \frac{g\left(m_{2}+M\right)}{l_{2}M} & 0 & 0 & 0 & 0\\
-gm_{1} & -gm_{2} & 0 & 0 & 0 & 0
\end{array}\right),\\
\boldsymbol{B} & =\left(0,0,0,-\frac{1}{l_{1}M},-\frac{1}{l_{2}M},1\right)^{T}.
\end{align}
Kalman's controllability matrix is
\begin{align}
\boldsymbol{\mathcal{K}} & =\left(\boldsymbol{B}|\boldsymbol{\mathcal{A}}\boldsymbol{B}|\boldsymbol{\mathcal{A}}^{2}\boldsymbol{B}|\cdots|\boldsymbol{\mathcal{A}}^{5}\boldsymbol{B}\right)\nonumber \\
 & =\left(\begin{array}{cccccc}
0 & -\frac{1}{l_{1}M} & 0 & -\frac{\alpha_{3}g}{l_{1}^{2}l_{2}M^{2}} & 0 & -\frac{\beta_{1}g^{2}}{l_{1}^{3}l_{2}^{2}M^{3}}\\
0 & -\frac{1}{l_{2}M} & 0 & -\frac{\alpha_{2}g}{l_{1}l_{2}^{2}M^{2}} & 0 & -\frac{\beta_{2}g^{2}}{l_{1}^{2}l_{2}^{3}M^{3}}\\
0 & \frac{1}{M} & 0 & \frac{\alpha_{1}g}{l_{1}l_{2}M^{2}} & 0 & \frac{\beta_{3}g^{2}}{l_{1}^{2}l_{2}^{2}M^{3}}\\
-\frac{1}{l_{1}M} & 0 & -\frac{\alpha_{3}g}{l_{1}^{2}l_{2}M^{2}} & 0 & -\frac{\beta_{1}g^{2}}{l_{1}^{3}l_{2}^{2}M^{3}} & 0\\
-\frac{1}{l_{2}M} & 0 & -\frac{\alpha_{2}g}{l_{1}l_{2}^{2}M^{2}} & 0 & -\frac{\beta_{2}g^{2}}{l_{1}^{2}l_{2}^{3}M^{3}} & 0\\
1 & 0 & \frac{\alpha_{1}g}{l_{1}l_{2}M} & 0 & \frac{\beta_{3}g^{2}}{l_{1}^{2}l_{2}^{2}M^{2}} & 0
\end{array}\right)
\end{align}
with
\begin{align}
\alpha_{1} & =l_{2}m_{1}+l_{1}m_{2}, & \alpha_{2} & =l_{1}\left(m_{2}+M\right)+l_{2}m_{1},\\
\alpha_{3} & =l_{2}\left(m_{1}+M\right)+l_{1}m_{2},
\end{align}
and
\begin{align}
\beta_{1} & =l_{1}^{2}m_{2}\left(m_{2}+M\right)+l_{2}l_{1}m_{2}\left(2m_{1}+M\right)+l_{2}^{2}\left(m_{1}+M\right){}^{2},\\
\beta_{2} & =l_{2}^{2}m_{1}\left(m_{1}+M\right)+l_{1}l_{2}m_{1}\left(2m_{2}+M\right)+l_{1}^{2}\left(m_{2}+M\right){}^{2},\\
\beta_{3} & =l_{1}^{2}m_{2}\left(m_{2}+M\right)+l_{2}^{2}m_{1}\left(m_{1}+M\right)+2l_{2}l_{1}m_{1}m_{2}.
\end{align}
The matrix $\boldsymbol{\mathcal{K}}$ has full row rank,
\begin{align}
\mbox{rank}\left(\boldsymbol{\mathcal{K}}\right) & =6,
\end{align}
as long as the pendulums are not identical. Consequently, the system
is controllable. A small deviation from the equilibrium position can
be counteracted by control. Pendulums with identical mass and lengths,
$m_{2}=m_{1}$ and $l_{2}=l_{1}$, respectively, yield
\begin{align}
\alpha_{2} & =\alpha_{3}, & \beta_{1} & =\beta_{2},
\end{align}
and the first and the second as well as the fourth and the fifth row
of the matrix $\boldsymbol{\mathcal{K}}$ become identical. Consequently,
the rank of $\boldsymbol{\mathcal{K}}$ changes, and two identical
pendulums cannot be controlled.

\end{example}

Both examples show one important consequence of controllability: arbitrary
many, parallel connected identical systems cannot be controlled \cite{kailath1980linear,chen1995linear}.
Expressed in a less rigorous language, controllability renders balancing
two identical brooms with only a single hand mathematically impossible.

\subsection{\label{sec:ControllabilityForRealizableTrajectories}Controllability
for systems satisfying the linearizing assumption}

Kalman's approach to controllability does not allow a direct generalization
to nonlinear systems. Furthermore, nothing is said about the trajectory
along which this transfer is achieved. To some extent, these questions
can be addressed in the framework of exactly realizable trajectories.
Here, a controllability matrix is derived which applies not only to
LTI systems but also to nonlinear systems satisfying the linearizing
assumption from Section \ref{sec:LinearizingAssumption}.

Consider the controlled system 
\begin{align}
\boldsymbol{\dot{x}}\left(t\right) & =\boldsymbol{R}\left(\boldsymbol{x}\left(t\right)\right)+\boldsymbol{\mathcal{B}}\left(\boldsymbol{x}\left(t\right)\right)\boldsymbol{u}\left(t\right)\label{eq:ControlledStateControllability}
\end{align}
together with the linearizing assumption 
\begin{align}
\boldsymbol{\mathcal{Q}}\boldsymbol{R}\left(\boldsymbol{x}\right) & =\boldsymbol{\mathcal{Q}}\boldsymbol{\mathcal{A}}\boldsymbol{x}+\boldsymbol{\mathcal{Q}}\boldsymbol{b}.\label{eq:LinearizingAssumptionControllability}
\end{align}
Equation \eqref{eq:LinearizingAssumptionControllability} implies
a linear constraint equation for an exactly realizable desired trajectory
$\boldsymbol{x}_{d}\left(t\right)$, 
\begin{align}
\boldsymbol{\mathcal{Q}}\boldsymbol{\dot{x}}_{d}\left(t\right) & =\boldsymbol{\mathcal{Q}}\boldsymbol{\mathcal{A}}\boldsymbol{x}_{d}\left(t\right)+\boldsymbol{\mathcal{Q}}\boldsymbol{b}.\label{eq:LinearConstraintEquation}
\end{align}
or, inserting $\boldsymbol{1}=\boldsymbol{\mathcal{P}}+\boldsymbol{\mathcal{Q}}$
between $\boldsymbol{\mathcal{A}}$ and $\boldsymbol{x}_{d}$, 
\begin{align}
\boldsymbol{\mathcal{Q}}\boldsymbol{\dot{x}}_{d}\left(t\right) & =\boldsymbol{\mathcal{Q}}\boldsymbol{\mathcal{A}}\boldsymbol{\mathcal{Q}}\boldsymbol{x}_{d}\left(t\right)+\boldsymbol{\mathcal{Q}}\boldsymbol{\mathcal{A}}\boldsymbol{\mathcal{P}}\boldsymbol{x}_{d}\left(t\right)+\boldsymbol{\mathcal{Q}}\boldsymbol{b}.\label{eq:ConstraintEquationLTI}
\end{align}
From now on, the parts $\boldsymbol{\mathcal{P}}\boldsymbol{x}_{d}\left(t\right)$
and $\boldsymbol{\mathcal{Q}}\boldsymbol{x}_{d}\left(t\right)$ are
considered as independent state components. The part $\boldsymbol{\mathcal{P}}\boldsymbol{x}_{d}\left(t\right)$
is prescribed by the experimenter while the part $\boldsymbol{\mathcal{Q}}\boldsymbol{x}_{d}\left(t\right)$
is governed by Eq. \eqref{eq:LinearConstraintEquation}. Equation
\eqref{eq:ConstraintEquationLTI} is a linear dynamical system for
the variable $\boldsymbol{\mathcal{Q}}\boldsymbol{x}_{d}\left(t\right)$
with inhomogeneity $\boldsymbol{\mathcal{Q}}\boldsymbol{\mathcal{A}}\boldsymbol{\mathcal{P}}\boldsymbol{x}_{d}\left(t\right)+\boldsymbol{\mathcal{Q}}\boldsymbol{b}$.
Achieving a transfer from the initial state $\boldsymbol{x}_{0}$
to the finite state $\boldsymbol{x}_{1}$ means the realizable trajectory
$\boldsymbol{x}_{d}\left(t\right)$ has to satisfy 
\begin{align}
\boldsymbol{x}_{d}\left(t_{0}\right) & =\boldsymbol{x}_{0},\label{eq:InitialConstraint}\\
\boldsymbol{x}_{d}\left(t_{1}\right) & =\boldsymbol{x}_{1}.\label{eq:TerminalConstraint}
\end{align}
Consequently, the prescribed part $\boldsymbol{\mathcal{P}}\boldsymbol{x}_{d}\left(t\right)$
satisfies
\begin{align}
\boldsymbol{\mathcal{P}}\boldsymbol{x}_{d}\left(t_{0}\right) & =\boldsymbol{\mathcal{P}}\boldsymbol{x}_{0}, & \boldsymbol{\mathcal{P}}\boldsymbol{x}_{d}\left(t_{1}\right) & =\boldsymbol{\mathcal{P}}\boldsymbol{x}_{1},\label{eq:Px_dInitialAndTerminalConstraints}
\end{align}
while the part $\boldsymbol{\mathcal{Q}}\boldsymbol{x}_{d}\left(t\right)$
satisfies
\begin{align}
\boldsymbol{\mathcal{Q}}\boldsymbol{x}_{d}\left(t_{0}\right) & =\boldsymbol{\mathcal{Q}}\boldsymbol{x}_{0}, & \boldsymbol{\mathcal{Q}}\boldsymbol{x}_{d}\left(t_{1}\right) & =\boldsymbol{\mathcal{Q}}\boldsymbol{x}_{1}.\label{eq:Qx_dInitialAndTerminalConstraints}
\end{align}
Being a linear equation, the solution $\boldsymbol{\mathcal{Q}}\boldsymbol{x}_{d}\left(t\right)$
to the constraint equation \eqref{eq:ConstraintEquationLTI} can be
expressed as a functional of $\boldsymbol{\mathcal{P}}\boldsymbol{x}_{d}\left(t\right)$,
\begin{align}
\boldsymbol{\mathcal{Q}}\boldsymbol{x}_{d}\left(t\right) & =\exp\left(\boldsymbol{\mathcal{Q}}\boldsymbol{\mathcal{A}}\boldsymbol{\mathcal{Q}}\left(t-t_{0}\right)\right)\boldsymbol{\mathcal{Q}}\boldsymbol{x}_{0}\nonumber \\
 & +\intop_{t_{0}}^{t}d\tau\exp\left(\boldsymbol{\mathcal{Q}}\boldsymbol{\mathcal{A}}\boldsymbol{\mathcal{Q}}\left(t-\tau\right)\right)\boldsymbol{\mathcal{Q}}\left(\boldsymbol{\mathcal{A}}\boldsymbol{\mathcal{P}}\boldsymbol{x}_{d}\left(\tau\right)+\boldsymbol{b}\right).\label{eq:Qx_dSolution}
\end{align}
See also Appendix \ref{sec:GeneralSolutionForForcedLinarDynamicalSystem}
for a derivation of the general solution to a forced linear dynamical
system. The solution Eq. \eqref{eq:Qx_dSolution} satisfies the initial
condition given by Eq. \eqref{eq:Qx_dInitialAndTerminalConstraints}.
Now, all initial and terminal conditions except $\boldsymbol{\mathcal{Q}}\boldsymbol{x}_{d}\left(t_{1}\right)=\boldsymbol{\mathcal{Q}}\boldsymbol{x}_{1}$
are satisfied. Enforcing this remaining terminal condition onto the
solution Eq. \eqref{eq:Qx_dSolution} yields
\begin{align}
\boldsymbol{\mathcal{Q}}\boldsymbol{x}_{1} & =\boldsymbol{\mathcal{Q}}\boldsymbol{x}_{d}\left(t_{1}\right)\nonumber \\
 & =\exp\left(\boldsymbol{\mathcal{Q}}\boldsymbol{\mathcal{A}}\boldsymbol{\mathcal{Q}}\left(t_{1}-t_{0}\right)\right)\boldsymbol{\mathcal{Q}}\boldsymbol{x}_{0}+\intop_{t_{0}}^{t_{1}}d\tau\exp\left(\boldsymbol{\mathcal{Q}}\boldsymbol{\mathcal{A}}\boldsymbol{\mathcal{Q}}\left(t_{1}-\tau\right)\right)\boldsymbol{\mathcal{Q}}\left(\boldsymbol{\mathcal{A}}\boldsymbol{\mathcal{P}}\boldsymbol{x}_{d}\left(\tau\right)+\boldsymbol{b}\right).\label{eq:Px_rConstraint}
\end{align}
This is actually a condition for the part $\boldsymbol{\mathcal{P}}\boldsymbol{x}_{d}\left(t\right)$.
Therefore, the transfer from $\boldsymbol{x}_{0}$ to $\boldsymbol{x}_{1}$
is achieved as long as the part $\boldsymbol{\mathcal{P}}\boldsymbol{x}_{d}$
satisfies Eqs. \eqref{eq:Px_dInitialAndTerminalConstraints} and \eqref{eq:Px_rConstraint}.
In between $t_{0}$ and $t_{1}$, the part $\boldsymbol{\mathcal{P}}\boldsymbol{x}_{d}\left(t\right)$
of the realizable trajectory can be freely chosen by the experimenter.
A system is controllable if at least one exactly realizable trajectory
$\boldsymbol{x}_{d}\left(t\right)$ can be found such that the constraints
Eqs. \eqref{eq:Px_dInitialAndTerminalConstraints} and \eqref{eq:Px_rConstraint}
are satisfied.

Analogously to the derivation of the Kalman rank condition in Section
\ref{sub:DerivationOfKalmanRank}, one can ask for the conditions
on the state matrices $\boldsymbol{\mathcal{A}}$ and projectors $\boldsymbol{\mathcal{P}}$
and $\boldsymbol{\mathcal{Q}}$ such that the constraint Eq. \eqref{eq:Px_rConstraint}
can be satisfied. Equation \eqref{eq:Px_rConstraint} is rearranged
as
\begin{align}
 & \exp\left(-\boldsymbol{\mathcal{Q}}\boldsymbol{\mathcal{A}}\boldsymbol{\mathcal{Q}}\left(t_{1}-t_{0}\right)\right)\boldsymbol{\mathcal{Q}}\boldsymbol{x}_{1}-\boldsymbol{\mathcal{Q}}\boldsymbol{x}_{0}-\intop_{t_{0}}^{t_{1}}d\tau\exp\left(\boldsymbol{\mathcal{Q}}\boldsymbol{\mathcal{A}}\boldsymbol{\mathcal{Q}}\left(t_{0}-\tau\right)\right)\boldsymbol{\mathcal{Q}}\boldsymbol{b}\nonumber \\
= & \intop_{t_{0}}^{t_{1}}d\tau\exp\left(\boldsymbol{\mathcal{Q}}\boldsymbol{\mathcal{A}}\boldsymbol{\mathcal{Q}}\left(t_{0}-\tau\right)\right)\boldsymbol{\mathcal{Q}}\boldsymbol{\mathcal{A}}\boldsymbol{\mathcal{P}}\boldsymbol{x}_{d}\left(\tau\right),\label{eq:Px_rConstraint_1}
\end{align}
and an argument equivalent to the derivation of the Kalman rank condition
\eqref{eq:KalmanRankCondition-1} is applied. Due to the Cayley-Hamilton
theorem, any power of matrices with $i\geq n$ can be expanded in
terms of lower order matrix powers as 
\begin{align}
\left(\boldsymbol{\mathcal{Q}}\boldsymbol{\mathcal{A}}\boldsymbol{\mathcal{Q}}\right)^{i} & =\sum_{k=0}^{n-1}d_{ik}\left(\boldsymbol{\mathcal{Q}}\boldsymbol{\mathcal{A}}\boldsymbol{\mathcal{Q}}\right)^{k}.\label{eq:CayleyHamiltonForQAQ}
\end{align}
The r. h. s. of Eq. \eqref{eq:Px_rConstraint_1} can be simplified
as
\begin{align}
 & \intop_{t_{0}}^{t_{1}}d\tau\exp\left(\boldsymbol{\mathcal{Q}}\boldsymbol{\mathcal{A}}\boldsymbol{\mathcal{Q}}\left(t_{0}-\tau\right)\right)\boldsymbol{\mathcal{Q}}\boldsymbol{\mathcal{A}}\boldsymbol{\mathcal{P}}\boldsymbol{x}_{d}\left(\tau\right)\nonumber \\
= & \sum_{k=0}^{\infty}\left(\boldsymbol{\mathcal{Q}}\boldsymbol{\mathcal{A}}\boldsymbol{\mathcal{Q}}\right)^{k}\boldsymbol{\mathcal{Q}}\boldsymbol{\mathcal{A}}\boldsymbol{\mathcal{P}}\intop_{t_{0}}^{t_{1}}d\tau\frac{\left(t_{0}-\tau\right)^{k}}{k!}\boldsymbol{\mathcal{P}}\boldsymbol{x}_{d}\left(\tau\right)\nonumber \\
= & \sum_{k=0}^{n-1}\left(\boldsymbol{\mathcal{Q}}\boldsymbol{\mathcal{A}}\boldsymbol{\mathcal{Q}}\right)^{k}\boldsymbol{\mathcal{Q}}\boldsymbol{\mathcal{A}}\boldsymbol{\mathcal{P}}\intop_{t_{0}}^{t_{1}}d\tau\frac{\left(t_{0}-\tau\right)^{k}}{k!}\boldsymbol{\mathcal{P}}\boldsymbol{x}_{d}\left(\tau\right)\nonumber \\
 & +\sum_{i=n}^{\infty}\left(\boldsymbol{\mathcal{Q}}\boldsymbol{\mathcal{A}}\boldsymbol{\mathcal{Q}}\right)^{i}\boldsymbol{\mathcal{Q}}\boldsymbol{\mathcal{A}}\boldsymbol{\mathcal{P}}\intop_{t_{0}}^{t_{1}}d\tau\frac{\left(t_{0}-\tau\right)^{i}}{i!}\boldsymbol{\mathcal{P}}\boldsymbol{x}_{d}\left(\tau\right)\nonumber \\
= & \sum_{k=0}^{n-1}\left(\boldsymbol{\mathcal{Q}}\boldsymbol{\mathcal{A}}\boldsymbol{\mathcal{Q}}\right)^{k}\boldsymbol{\mathcal{Q}}\boldsymbol{\mathcal{A}}\boldsymbol{\mathcal{P}}\intop_{t_{0}}^{t_{1}}d\tau\left(\frac{\left(t_{0}-\tau\right)^{k}}{k!}+\sum_{i=n}^{\infty}d_{ik}\frac{\left(t_{0}-\tau\right)^{i}}{i!}\right)\boldsymbol{\mathcal{P}}\boldsymbol{x}_{d}\left(\tau\right),\label{eq:TruncateExpQAQ}
\end{align}
such that Eq. \eqref{eq:Px_rConstraint_1} becomes a truncated sum
\begin{align}
 & \exp\left(-\boldsymbol{\mathcal{Q}}\boldsymbol{\mathcal{A}}\boldsymbol{\mathcal{Q}}\left(t_{1}-t_{0}\right)\right)\boldsymbol{\mathcal{Q}}\boldsymbol{x}_{1}-\boldsymbol{\mathcal{Q}}\boldsymbol{x}_{0}-\intop_{t_{0}}^{t_{1}}d\tau\exp\left(\boldsymbol{\mathcal{Q}}\boldsymbol{\mathcal{A}}\boldsymbol{\mathcal{Q}}\left(t_{0}-\tau\right)\right)\boldsymbol{\mathcal{Q}}\boldsymbol{b}\nonumber \\
= & \sum_{k=0}^{n-1}\left(\boldsymbol{\mathcal{Q}}\boldsymbol{\mathcal{A}}\boldsymbol{\mathcal{Q}}\right)^{k}\boldsymbol{\mathcal{Q}}\boldsymbol{\mathcal{A}}\boldsymbol{\mathcal{P}}\boldsymbol{\alpha}_{k}\left(t_{1},t_{0}\right).\label{eq:Px_rConstraint_2}
\end{align}
Define the $n\times1$ vectors
\begin{align}
\boldsymbol{\alpha}_{k}\left(t_{1},t_{0}\right) & =\intop_{t_{0}}^{t_{1}}d\tau\left(\frac{\left(t_{0}-\tau\right)^{k}}{k!}+\sum_{i=n}^{\infty}d_{ik}\frac{\left(t_{0}-\tau\right)^{i}}{i!}\right)\boldsymbol{\mathcal{P}}\boldsymbol{x}_{d}\left(\tau\right).
\end{align}
The right hand side of Eq. \eqref{eq:Px_rConstraint_2} can be written
with the help of the $n^{2}\times1$ vector 
\begin{align}
\boldsymbol{\alpha}\left(t_{1},t_{0}\right) & =\left(\begin{array}{c}
\boldsymbol{\alpha}_{0}\left(t_{1},t_{0}\right)\\
\boldsymbol{\alpha}_{1}\left(t_{1},t_{0}\right)\\
\vdots\\
\boldsymbol{\alpha}_{n-1}\left(t_{1},t_{0}\right)
\end{array}\right)
\end{align}
as 
\begin{align}
\exp\left(-\boldsymbol{\mathcal{Q}}\boldsymbol{\mathcal{A}}\boldsymbol{\mathcal{Q}}\left(t_{1}-t_{0}\right)\right)\boldsymbol{\mathcal{Q}}\boldsymbol{x}_{1}-\boldsymbol{\mathcal{Q}}\boldsymbol{x}_{0}\nonumber \\
-\intop_{t_{0}}^{t_{1}}d\tau\exp\left(\boldsymbol{\mathcal{Q}}\boldsymbol{\mathcal{A}}\boldsymbol{\mathcal{Q}}\left(t_{0}-\tau\right)\right)\boldsymbol{\mathcal{Q}}\boldsymbol{b} & =\boldsymbol{\mathcal{\tilde{K}}}\boldsymbol{\alpha}\left(t_{1},t_{0}\right).\label{eq:Px_rConstraint_3}
\end{align}
The $n\times n^{2}$ \textit{controllability matrix $\boldsymbol{\mathcal{\tilde{K}}}$
}is defined by
\begin{align}
\boldsymbol{\mathcal{\tilde{K}}} & =\left(\boldsymbol{\mathcal{Q}}\boldsymbol{\mathcal{A}}\boldsymbol{\mathcal{P}}|\boldsymbol{\mathcal{Q}}\boldsymbol{\mathcal{A}}\boldsymbol{\mathcal{Q}}\boldsymbol{\mathcal{A}}\boldsymbol{\mathcal{P}}|\cdots|\left(\boldsymbol{\mathcal{Q}}\boldsymbol{\mathcal{A}}\boldsymbol{\mathcal{Q}}\right)^{n-1}\boldsymbol{\mathcal{Q}}\boldsymbol{\mathcal{A}}\boldsymbol{\mathcal{P}}\right).\label{eq:ControllabilityMatrixForRealizableTrajectories}
\end{align}
The left hand side of Eq. \eqref{eq:Px_rConstraint_3} can be any
point in $\boldsymbol{\mathcal{Q}}\mathbb{R}^{n}=\mathbb{R}^{n-p}$.
The mapping is surjective, i.e., every element on the left hand side
has a corresponding element on the right hand side, if $\boldsymbol{\mathcal{\tilde{K}}}$
has full rank $n-p$. Therefore, the nonlinear affine control system
Eq. \eqref{eq:ControlledStateControllability} satisfying the linearizing
assumption Eq. \eqref{eq:LinearizingAssumptionControllability} is
controllable if
\begin{align}
\text{rank}\left(\boldsymbol{\mathcal{\tilde{K}}}\right) & =n-p.\label{eq:RankConditionRealizableTrajectories}
\end{align}

\begin{example}[Single input diagonal LTI system]\label{ex:DiagonalLTI2}

Consider the LTI system from Example \ref{ex:DiagonalLTI1}. The two
parts of $\boldsymbol{\mathcal{A}}$ necessary for the computation
of the controllability matrix $\boldsymbol{\mathcal{\tilde{K}}}$
are
\begin{align}
\boldsymbol{\mathcal{Q}}\boldsymbol{\mathcal{A}}\boldsymbol{\mathcal{P}} & =\frac{1}{4}\left(\begin{array}{cc}
\lambda_{1}-\lambda_{2} & \lambda_{1}-\lambda_{2}\\
\lambda_{2}-\lambda_{1} & \lambda_{2}-\lambda_{1}
\end{array}\right), & \boldsymbol{\mathcal{Q}}\boldsymbol{\mathcal{A}}\boldsymbol{\mathcal{Q}} & =\frac{1}{4}\left(\begin{array}{cc}
\lambda_{1}+\lambda_{2} & \lambda_{1}-\lambda_{2}\\
\lambda_{2}-\lambda_{1} & \lambda_{1}+\lambda_{2}
\end{array}\right).
\end{align}
The controllability matrix is
\begin{align}
\boldsymbol{\mathcal{\tilde{K}}} & =\left(\boldsymbol{\mathcal{Q}}\boldsymbol{\mathcal{A}}\boldsymbol{\mathcal{P}}|\boldsymbol{\mathcal{Q}}\boldsymbol{\mathcal{A}}\boldsymbol{\mathcal{Q}}\boldsymbol{\mathcal{A}}\boldsymbol{\mathcal{P}}\right)\nonumber \\
 & =\frac{1}{4}\left(\begin{array}{cccc}
\left(\lambda_{1}-\lambda_{2}\right) & \left(\lambda_{1}-\lambda_{2}\right) & \frac{1}{2}\left(\lambda_{1}^{2}-\lambda_{2}^{2}\right) & \frac{1}{2}\left(\lambda_{1}^{2}-\lambda_{2}^{2}\right)\\
\left(\lambda_{2}-\lambda_{1}\right) & \left(\lambda_{2}-\lambda_{1}\right) & \frac{1}{2}\left(\lambda_{2}^{2}-\lambda_{1}^{2}\right) & \frac{1}{2}\left(\lambda_{2}^{2}-\lambda_{1}^{2}\right)
\end{array}\right).
\end{align}
The upper row of $\boldsymbol{\mathcal{\tilde{K}}}$ equals the lower
row times $-1$, i.e., the rows are linearly dependent and so $\boldsymbol{\mathcal{\tilde{K}}}$
has rank 
\begin{align}
\text{rank}\left(\boldsymbol{\mathcal{\tilde{K}}}\right) & =1.
\end{align}
If $\lambda_{1}=\lambda_{2}$, all entries of $\boldsymbol{\mathcal{\tilde{K}}}$
vanish and then $\boldsymbol{\mathcal{\tilde{K}}}$ has zero rank.
The system Eq. \eqref{eq:E110} is controllable as long as $\lambda_{1}\neq\lambda_{2}$.

\end{example}

\begin{example}[Controllability of the activator-controlled FHN model]\label{ex:FHN4}

Controllability in form of a rank condition can be discussed for all
models of the form
\begin{align}
\left(\begin{array}{c}
\dot{x}\left(t\right)\\
\dot{y}\left(t\right)
\end{array}\right) & =\left(\begin{array}{c}
a_{0}+a_{1}x\left(t\right)+a_{2}y\left(t\right)\\
R\left(x\left(t\right),y\left(t\right)\right)
\end{array}\right)+\left(\begin{array}{c}
0\\
B\left(x\left(t\right),y\left(t\right)\right)
\end{array}\right)u\left(t\right).\label{eq:Eq220}
\end{align}
A prominent example is the activator-controlled FHN model. The $\boldsymbol{\mathcal{Q}}$
part of the nonlinearity $\boldsymbol{R}$ is actually a linear function
of the state $\boldsymbol{x}$, 
\begin{align}
\boldsymbol{\mathcal{Q}}\boldsymbol{R}\left(\boldsymbol{x}\left(t\right)\right) & =\left(\begin{array}{c}
a_{1}x\left(t\right)+a_{2}y\left(t\right)\\
0
\end{array}\right)+\left(\begin{array}{c}
a_{0}\\
0
\end{array}\right)=\boldsymbol{\mathcal{Q}}\boldsymbol{\mathcal{A}}\boldsymbol{x}\left(t\right)+\boldsymbol{\mathcal{Q}}\boldsymbol{b},
\end{align}
i.e., this model satisfies the linearizing assumption with the matrix
$\boldsymbol{\mathcal{A}}$ and vector $\boldsymbol{b}$ defined by
\begin{align}
\boldsymbol{\mathcal{Q}}\boldsymbol{\mathcal{A}} & =\left(\begin{array}{cc}
a_{1} & a_{2}\\
0 & 0
\end{array}\right), & \boldsymbol{\mathcal{Q}}\boldsymbol{b} & =\left(\begin{array}{c}
a_{0}\\
0
\end{array}\right).
\end{align}
The controllability matrix $\boldsymbol{\mathcal{\tilde{K}}}$ is
\begin{align}
\boldsymbol{\mathcal{\tilde{K}}} & =\left(\boldsymbol{\mathcal{Q}}\boldsymbol{\mathcal{A}}\boldsymbol{\mathcal{P}}|\boldsymbol{\mathcal{Q}}\boldsymbol{\mathcal{A}}\boldsymbol{\mathcal{Q}}\boldsymbol{\mathcal{A}}\boldsymbol{\mathcal{P}}\right)=\left(\begin{array}{cccc}
0 & a_{2} & 0 & a_{1}a_{2}\\
0 & 0 & 0 & 0
\end{array}\right),
\end{align}
and, for $a_{2}\neq0$, $\boldsymbol{\mathcal{\tilde{K}}}$ has rank
\begin{align}
\mbox{rank}\left(\boldsymbol{\mathcal{\tilde{K}}}\right) & =1=n-p.
\end{align}
The activator-controlled FHN model is controllable as long as $a_{2}\neq0$,
i.e., as long as the equation for the inhibitor $x$ also depends
on the activator $y$. The control directly affects the activator
$y$. If $a_{2}=0$ in Eq. \eqref{eq:Eq220}, the inhibitor evolves
decoupled from the activator, and therefore cannot be affected by
control.

\end{example}

\begin{example}[Controllability of the controlled SIR model]\label{ex:SIRModel2}

With the help of the projectors $\boldsymbol{\mathcal{P}}$ and $\boldsymbol{\mathcal{Q}}$
computed in Example \ref{ex:SIRModel1_1}, the controllability matrix
is obtained as 
\begin{align}
\boldsymbol{\mathcal{\tilde{K}}} & =\left(\boldsymbol{\mathcal{Q}}\boldsymbol{\mathcal{A}}\boldsymbol{\mathcal{P}}|\boldsymbol{\mathcal{Q}}\boldsymbol{\mathcal{A}}\boldsymbol{\mathcal{Q}}\boldsymbol{\mathcal{A}}\boldsymbol{\mathcal{P}}|\boldsymbol{\mathcal{Q}}\boldsymbol{\mathcal{A}}\boldsymbol{\mathcal{Q}}\boldsymbol{\mathcal{A}}\boldsymbol{\mathcal{Q}}\boldsymbol{\mathcal{A}}\boldsymbol{\mathcal{P}}\right)\nonumber \\
 & =\left(\begin{array}{ccccccccc}
\frac{\gamma}{4} & -\frac{\gamma}{4} & 0 & -\frac{\gamma^{2}}{8} & \frac{\gamma^{2}}{8} & 0 & \frac{\gamma^{3}}{16} & -\frac{\gamma^{3}}{16} & 0\\
\frac{\gamma}{4} & -\frac{\gamma}{4} & 0 & -\frac{\gamma^{2}}{8} & \frac{\gamma^{2}}{8} & 0 & \frac{\gamma^{3}}{16} & -\frac{\gamma^{3}}{16} & 0\\
-\frac{\gamma}{2} & \frac{\gamma}{2} & 0 & \frac{\gamma^{2}}{4} & -\frac{\gamma^{2}}{4} & 0 & -\frac{\gamma^{3}}{8} & \frac{\gamma^{3}}{8} & 0
\end{array}\right).
\end{align}
As long as $\gamma\neq0$, the rank of $\boldsymbol{\mathcal{\tilde{K}}}$
is 
\begin{align}
\mbox{rank}\left(\boldsymbol{\mathcal{\tilde{K}}}\right) & =1<n-p=2.
\end{align}
Thus, the rank of $\boldsymbol{\mathcal{\tilde{K}}}$ is smaller than
$n-p$, and consequently the SIR model is \textit{not} controllable.
It is impossible to find a control to reach every final state $\boldsymbol{x}_{1}$
from every other initial state $\boldsymbol{x}_{0}$. Intuitively,
the reason is simple to understand. The controlled SIR model satisfies
a conservation law, see Example \ref{ex:SIRModel1}. Independent of
the actual time dependence of the control signal $u\left(t\right)$,
the total number $N$ of individuals is conserved,
\begin{align}
S\left(t\right)+I\left(t\right)+R\left(t\right) & =N.
\end{align}
The value of $N$ is prescribed by the initial condition $\boldsymbol{x}\left(t_{0}\right)=\boldsymbol{x}_{0}$.
For all times, the dynamics of the controlled SIR model is restricted
to a two-dimensional surface embedded in the three-dimensional state
space. Hence, the system's state vector can only reach points lying
on this surface, and no control can force the system to leave it.

\end{example}

\subsection{Discussion}

A controllability matrix $\boldsymbol{\mathcal{\tilde{K}}}$, Eq.
\eqref{eq:ControllabilityMatrixForRealizableTrajectories}, is derived
in the framework of exactly realizable trajectories. If $\boldsymbol{\mathcal{\tilde{K}}}$
satisfies the rank condition $\text{rank}\left(\boldsymbol{\mathcal{\tilde{K}}}\right)=n-p$,
the system is controllable. At least one control signal exists which
achieves a transfer from an arbitrary initial state $\boldsymbol{x}\left(t_{0}\right)=\boldsymbol{x}_{0}$
to an arbitrary final state $\boldsymbol{x}\left(t_{1}\right)=\boldsymbol{x}_{1}$
within the finite time interval $t_{1}-t_{0}$.

The controllability matrix $\boldsymbol{\mathcal{\tilde{K}}}$ can
be computed for all LTI system. We expect that the rank condition
for controllability, Eq. \eqref{eq:RankConditionRealizableTrajectories}
is fully equivalent to Kalman's rank condition, Eq. \eqref{eq:KalmanRankCondition-1}.
If the system is controllable in terms of $\boldsymbol{\mathcal{\tilde{K}}}$,
it is also controllable in terms of $\boldsymbol{\mathcal{K}}$, and
vice versa. The advantage of controllability in terms of $\boldsymbol{\mathcal{\tilde{K}}}$
is its applicability to a certain class of nonlinear systems. For
affine dynamical systems satisfying the linearizing assumption \eqref{eq:LinearizingAssumption2},
the rank condition for $\boldsymbol{\mathcal{\tilde{K}}}$ remains
a valid check for controllability. This class encompasses a number
of simple nonlinear models which are of interest to physicists. In
particular, it is proven that all mechanical control systems in one
spatial dimension are controllable, see Example \ref{ex:FHN4}. Other
systems satisfying the linearizing assumption are the controlled SIR
model, Example \ref{ex:SIRModel2}, and the activator-controlled FHN
model, see Example \ref{ex:FHN4}. Controllability as proposed here
cannot be applied to the inhibitor-controlled FHN model because the
corresponding constraint equation is nonlinear. Checking its controllability
requires a notion of nonlinear controllability for general nonlinear
systems. Nonlinear controllability cannot be defined in form of a
simple rank condition for a controllability matrix but demands more
difficult concepts.

Exactly realizable trajectories allow a characterization of the entirety
of state trajectories along which a state transfer can be achieved.
Any desired trajectory which satisfies the constraint equation
\begin{align}
\boldsymbol{\mathcal{Q}}\left(\boldsymbol{\dot{x}}_{d}\left(t\right)-\boldsymbol{\mathcal{A}}\boldsymbol{x}_{d}\left(t\right)-\boldsymbol{b}\right) & =\mathbf{0},\label{eq:Eq2153}
\end{align}
and the initial and terminal conditions
\begin{align}
\boldsymbol{x}_{d}\left(t_{0}\right) & =\boldsymbol{x}_{0}, & \boldsymbol{x}_{d}\left(t_{1}\right) & =\boldsymbol{x}_{1}\label{eq:Eq2154}
\end{align}
does the job. For example, a second order differential equation for
$\boldsymbol{x}_{d}\left(t\right)$ can accommodate both initial and
terminal conditions Eqs. \eqref{eq:Eq2154}. A successful transfer
from $\boldsymbol{x}_{0}$ to $\boldsymbol{x}_{1}$ is achieved if
$\boldsymbol{x}_{d}\left(t\right)$ additionally satisfies the constraint
equation \eqref{eq:Eq2153}. The control signal is given by
\begin{align}
\boldsymbol{u}\left(t\right) & =\boldsymbol{\mathcal{B}}^{+}\left(\boldsymbol{\dot{x}}_{d}\left(t\right)-\boldsymbol{\mathcal{A}}\boldsymbol{x}_{d}\left(t\right)-\boldsymbol{b}\right).\label{eq:Eq204}
\end{align}
Equation \eqref{eq:Eq204} can be used to obtain an expression for
the control which depends only on the part $\boldsymbol{\mathcal{P}}\boldsymbol{x}_{d}\left(t\right)$.
According to equation \eqref{eq:Qx_dSolution}, the solution for $\boldsymbol{\mathcal{Q}}\boldsymbol{x}_{d}\left(t\right)$
can be expressed in terms of a functional of $\boldsymbol{\mathcal{P}}\boldsymbol{x}_{d}\left(t\right)$,
\begin{align}
\boldsymbol{\mathcal{Q}}\boldsymbol{x}_{d}\left[\boldsymbol{\mathcal{P}}\boldsymbol{x}_{d}\left(t\right)\right] & =\exp\left(\boldsymbol{\mathcal{Q}}\boldsymbol{\mathcal{A}}\boldsymbol{\mathcal{Q}}\left(t-t_{0}\right)\right)\boldsymbol{\mathcal{Q}}\boldsymbol{x}_{0}\nonumber \\
 & +\intop_{t_{0}}^{t}d\tau\exp\left(\boldsymbol{\mathcal{Q}}\boldsymbol{\mathcal{A}}\boldsymbol{\mathcal{Q}}\left(t-\tau\right)\right)\boldsymbol{\mathcal{Q}}\left(\boldsymbol{\mathcal{A}}\boldsymbol{\mathcal{P}}\boldsymbol{x}_{d}\left(\tau\right)+\boldsymbol{b}\right).
\end{align}
Consequently, Eq. \eqref{eq:Eq204} becomes a functional of $\boldsymbol{\mathcal{P}}\boldsymbol{x}_{d}\left(t\right)$
as well, 
\begin{align}
\boldsymbol{u}\left[\boldsymbol{\mathcal{P}}\boldsymbol{x}_{d}\left(t\right)\right] & =\boldsymbol{\mathcal{B}}^{+}\left(\boldsymbol{\mathcal{P}}\boldsymbol{\dot{x}}_{d}\left(t\right)-\boldsymbol{\mathcal{P}}\boldsymbol{\mathcal{A}}\boldsymbol{\mathcal{P}}\boldsymbol{x}_{d}\left(t\right)-\boldsymbol{\mathcal{P}}\boldsymbol{\mathcal{A}}\exp\left(\boldsymbol{\mathcal{Q}}\boldsymbol{\mathcal{A}}\boldsymbol{\mathcal{Q}}\left(t-t_{0}\right)\right)\boldsymbol{\mathcal{Q}}\boldsymbol{x}_{0}-\boldsymbol{\mathcal{P}}\boldsymbol{b}\right)\nonumber \\
 & -\boldsymbol{\mathcal{B}}^{+}\boldsymbol{\mathcal{A}}\intop_{t_{0}}^{t}d\tau\exp\left(\boldsymbol{\mathcal{Q}}\boldsymbol{\mathcal{A}}\boldsymbol{\mathcal{Q}}\left(t-\tau\right)\right)\boldsymbol{\mathcal{Q}}\left(\boldsymbol{\mathcal{A}}\boldsymbol{\mathcal{P}}\boldsymbol{x}_{d}\left(\tau\right)+\boldsymbol{b}\right).\label{eq:Eq2157}
\end{align}
Thus, any reference to $\boldsymbol{\mathcal{Q}}\boldsymbol{x}_{d}\left(t\right)$
except for the initial condition $\boldsymbol{\mathcal{Q}}\boldsymbol{x}_{0}$
is eliminated from the expression for the control signal. The control
signal is entirely expressed in terms of the part $\boldsymbol{\mathcal{P}}\boldsymbol{x}_{d}\left(t\right)$
prescribed by the experimenter.

Using the complementary projectors $\boldsymbol{\mathcal{P}}$ and
$\boldsymbol{\mathcal{Q}}$, the state matrix $\boldsymbol{\mathcal{A}}$
can be split up in four parts as 
\begin{align}
\boldsymbol{\mathcal{A}} & =\boldsymbol{\mathcal{P}}\boldsymbol{\mathcal{A}}\boldsymbol{\mathcal{P}}+\boldsymbol{\mathcal{P}}\boldsymbol{\mathcal{A}}\boldsymbol{\mathcal{Q}}+\boldsymbol{\mathcal{Q}}\boldsymbol{\mathcal{A}}\boldsymbol{\mathcal{P}}+\boldsymbol{\mathcal{Q}}\boldsymbol{\mathcal{A}}\boldsymbol{\mathcal{Q}}.
\end{align}
Note that the controllability matrix $\boldsymbol{\mathcal{\tilde{K}}}$,
Eq. \eqref{eq:ControllabilityMatrixForRealizableTrajectories}, does
only depend on the parts $\boldsymbol{\mathcal{Q}}\boldsymbol{\mathcal{A}}\boldsymbol{\mathcal{P}}$
and $\boldsymbol{\mathcal{Q}}\boldsymbol{\mathcal{A}}\boldsymbol{\mathcal{Q}}$,
but not on $\boldsymbol{\mathcal{P}}\boldsymbol{\mathcal{A}}\boldsymbol{\mathcal{P}}$
and $\boldsymbol{\mathcal{P}}\boldsymbol{\mathcal{A}}\boldsymbol{\mathcal{Q}}$.
This fact extends the validity of the controllability matrix $\boldsymbol{\mathcal{\tilde{K}}}$
to nonlinear systems satisfying the linearizing assumption. Furthermore,
only the parts $\boldsymbol{\mathcal{Q}}\boldsymbol{\mathcal{A}}\boldsymbol{\mathcal{P}}$
and $\boldsymbol{\mathcal{Q}}\boldsymbol{\mathcal{A}}\boldsymbol{\mathcal{Q}}$
must be known to decide if a system is controllable. Thus, it can
be possible to decide about controllability of a system without knowing
all details of its dynamics. This insight might be useful for experimental
systems with incomplete or approximated model equations.

\section{\label{sec:OutputControllability}Output controllability}

Consider the dynamical system
\begin{align}
\boldsymbol{\dot{x}}\left(t\right) & =\boldsymbol{R}\left(\boldsymbol{x}\left(t\right)\right)+\boldsymbol{\mathcal{B}}\left(\boldsymbol{x}\left(t\right)\right)\boldsymbol{u}\left(t\right),
\end{align}
together with the output 
\begin{align}
\boldsymbol{z}\left(t\right) & =\boldsymbol{h}\left(\boldsymbol{x}\left(t\right)\right).\label{eq:NonlinearOutput}
\end{align}
Here, $\boldsymbol{z}\left(t\right)=\left(z_{1}\left(t\right),\dots,z_{m}\left(t\right)\right)^{T}\in\mathbb{R}^{m}$
with $m\leq n$ components is called the \textit{output vector} and
the \textit{output function} $\boldsymbol{h}$ maps from $\mathbb{R}^{n}$
to $\mathbb{R}^{m}$.

A system is called output controllable if it is possible to achieve
a transfer from an initial output state
\begin{align}
\boldsymbol{z}\left(t_{0}\right) & =\boldsymbol{z}_{0}\label{eq:InitialOutput}
\end{align}
at time $t=t_{0}$ to a terminal output state
\begin{align}
\boldsymbol{z}\left(t_{1}\right) & =\boldsymbol{z}_{1}\label{eq:TerminalOutput}
\end{align}
at the terminal time $t=t_{1}$. In contrast to output controllability,
the notion of controllability discussed in Section \ref{sec:Controllability}
is concerned with the controllability of the state $\boldsymbol{x}\left(t\right)$
and is often referred to as state or full state controllability. Note
that a state controllable system is not necessarily output controllable.
Similarly, an output controllable system is not necessarily state
controllable. For $m=n$ and an output function equal to the identity
function, $\boldsymbol{h}\left(\boldsymbol{x}\right)=\boldsymbol{x}$,
output controllability is equivalent to state controllability.

\subsection{Kalman rank condition for the output controllability of LTI systems}

The notion of state controllability developed in form of a Kalman
rank condition for an output controllability matrix can be adapted
to output controllability \cite{kalman1959general,kalman1960contributions,chen1995linear}.
Consider the LTI system with $n$-dimensional state vector $\boldsymbol{x}\left(t\right)$
and $p$-dimensional control signal $\boldsymbol{u}\left(t\right)$,
\begin{align}
\boldsymbol{\dot{x}}\left(t\right) & =\boldsymbol{\mathcal{A}}\boldsymbol{x}\left(t\right)+\boldsymbol{\mathcal{B}}\boldsymbol{u}\left(t\right).\label{eq:LTIOutput}
\end{align}
The output is assumed to be a linear relation of the form
\begin{align}
\boldsymbol{z}\left(t\right) & =\boldsymbol{\mathcal{C}}\boldsymbol{x}\left(t\right)\label{eq:LinearOutput}
\end{align}
with $m\times n$ \textit{output matrix} $\boldsymbol{\mathcal{C}}$.
The $m\times np$ \textit{output controllability matrix} is defined
as
\begin{align}
\boldsymbol{\mathcal{K}}_{\boldsymbol{\mathcal{C}}}= & \left(\boldsymbol{\mathcal{C}}\boldsymbol{\mathcal{B}}|\boldsymbol{\mathcal{C}}\boldsymbol{\mathcal{A}}\boldsymbol{\mathcal{B}}|\boldsymbol{\mathcal{C}}\boldsymbol{\mathcal{A}}^{2}\boldsymbol{\mathcal{B}}|\cdots|\boldsymbol{\mathcal{C}}\boldsymbol{\mathcal{A}}^{n-1}\boldsymbol{\mathcal{B}}\right).\label{eq:KalmanOutputControllabilityMatrix}
\end{align}
The LTI system Eq. \eqref{eq:LTIOutput} with output \eqref{eq:LinearOutput}
is output controllable if $\boldsymbol{\mathcal{K}}_{\boldsymbol{\mathcal{C}}}$
satisfies the rank condition
\begin{align}
\text{rank}\left(\boldsymbol{\mathcal{K}}_{\boldsymbol{\mathcal{C}}}\right) & =m.\label{eq:KalmanOutputRankCondition}
\end{align}
A proof of Eq. \eqref{eq:KalmanOutputRankCondition} proceeds along
the same lines as the proof for the Kalman rank condition for state
controllability in Section \ref{sub:DerivationOfKalmanRank}. Using
Eq. \eqref{eq:LTISolution}, the solution for $\boldsymbol{z}\left(t\right)$
is
\begin{align}
\boldsymbol{z}\left(t\right) & =\boldsymbol{\mathcal{C}}\boldsymbol{x}\left(t\right)=\boldsymbol{\mathcal{C}}e^{\boldsymbol{\mathcal{A}}\left(t-t_{0}\right)}\boldsymbol{x}_{0}+\boldsymbol{\mathcal{C}}\intop_{t_{0}}^{t}d\tau e^{\boldsymbol{\mathcal{A}}\left(t-\tau\right)}\boldsymbol{\mathcal{B}}\boldsymbol{u}\left(\tau\right).\label{eq:OutputSolution}
\end{align}
Evaluating Eq. \eqref{eq:OutputSolution} at the terminal time $t=t_{1}$
and enforcing the terminal output condition Eq. \eqref{eq:TerminalOutput}
yields a condition for the control signal $\boldsymbol{u}$, 
\begin{align}
\boldsymbol{z}_{1} & =\boldsymbol{z}\left(t_{1}\right)=\boldsymbol{\mathcal{C}}e^{\boldsymbol{\mathcal{A}}\left(t_{1}-t_{0}\right)}\boldsymbol{x}_{0}+\boldsymbol{\mathcal{C}}\intop_{t_{0}}^{t_{1}}d\tau e^{\boldsymbol{\mathcal{A}}\left(t_{1}-\tau\right)}\boldsymbol{\mathcal{B}}\boldsymbol{u}\left(\tau\right).
\end{align}
Exploiting the Cayley-Hamilton theorem and proceeding analogously
to Eq. \eqref{eq:Eq297} yields
\begin{align}
\boldsymbol{z}_{1}-\boldsymbol{\mathcal{C}}e^{\boldsymbol{\mathcal{A}}\left(t_{1}-t_{0}\right)}\boldsymbol{x}_{0} & =\boldsymbol{\mathcal{C}}\sum_{k=0}^{n-1}\boldsymbol{\mathcal{A}}^{k}\boldsymbol{\mathcal{B}}\boldsymbol{\tilde{\beta}}_{k}\left(t_{1},t_{0}\right)=\boldsymbol{\mathcal{K}}_{\boldsymbol{\mathcal{C}}}\boldsymbol{\tilde{\beta}}\left(t_{1},t_{0}\right)
\end{align}
with $p\times1$ vectors $\boldsymbol{\tilde{\beta}}_{k}$ defined
as
\begin{align}
\boldsymbol{\tilde{\beta}}_{k}\left(t_{1},t_{0}\right) & =\intop_{t_{0}}^{t_{1}}d\tau\dfrac{\left(t_{1}-\tau\right)^{k}}{k!}\boldsymbol{u}\left(\tau\right)+\sum_{i=n}^{\infty}c_{ik}\intop_{t_{0}}^{t_{1}}d\tau\dfrac{\left(t_{1}-\tau\right)^{i}}{i!}\boldsymbol{u}\left(\tau\right)
\end{align}
and the $np\times1$ vector 
\begin{align}
\boldsymbol{\tilde{\beta}}\left(t_{1},t_{0}\right) & =\left(\begin{array}{c}
\boldsymbol{\tilde{\beta}}_{0}\left(t_{1},t_{0}\right)\\
\vdots\\
\boldsymbol{\tilde{\beta}}_{n-1}\left(t_{1},t_{0}\right)
\end{array}\right).
\end{align}
The linear map from $\boldsymbol{\tilde{\beta}}\left(t_{1},t_{0}\right)$
to $\boldsymbol{z}_{1}$ must be \textit{surjective}. This is the
case if and only if the matrix $\boldsymbol{\mathcal{K}}_{\boldsymbol{\mathcal{C}}}$
defined in Eq. \eqref{eq:KalmanOutputControllabilityMatrix} has full
row rank, i.e.,
\begin{align}
\text{rank}\left(\boldsymbol{\mathcal{K}}_{\boldsymbol{\mathcal{C}}}\right) & =m.
\end{align}

\subsection{Output controllability for systems satisfying the linearizing assumption}

Using the framework of exactly realizable trajectories, we can generalize
the condition for output controllability in form of a matrix rank
condition to nonlinear affine control systems satisfying the linearizing
assumption from Section \ref{sec:LinearizingAssumption}. Consider
the affine control system
\begin{align}
\boldsymbol{\dot{x}}\left(t\right) & =\boldsymbol{R}\left(\boldsymbol{x}\left(t\right)\right)+\boldsymbol{\mathcal{B}}\left(\boldsymbol{x}\left(t\right)\right)\boldsymbol{u}\left(t\right)
\end{align}
with linear output 
\begin{align}
\boldsymbol{z}\left(t\right) & =\boldsymbol{\mathcal{C}}\boldsymbol{x}\left(t\right).\label{eq:LinearOutput-1}
\end{align}
The constraint equation for exactly realizable desired trajectories
$\boldsymbol{x}_{d}\left(t\right)$ is linear,
\begin{align}
\boldsymbol{\mathcal{Q}}\boldsymbol{\dot{x}}_{d}\left(t\right) & =\boldsymbol{\mathcal{Q}}\boldsymbol{\mathcal{A}}\boldsymbol{x}_{d}\left(t\right)+\boldsymbol{\mathcal{Q}}\boldsymbol{b},
\end{align}
and has the solution 
\begin{align}
\boldsymbol{\mathcal{Q}}\boldsymbol{x}_{d}\left(t\right) & =\exp\left(\boldsymbol{\mathcal{Q}}\boldsymbol{\mathcal{A}}\boldsymbol{\mathcal{Q}}\left(t-t_{0}\right)\right)\boldsymbol{\mathcal{Q}}\boldsymbol{x}_{0}\nonumber \\
 & +\intop_{t_{0}}^{t}d\tau\exp\left(\boldsymbol{\mathcal{Q}}\boldsymbol{\mathcal{A}}\boldsymbol{\mathcal{Q}}\left(t-\tau\right)\right)\boldsymbol{\mathcal{Q}}\left(\boldsymbol{\mathcal{A}}\boldsymbol{\mathcal{P}}\boldsymbol{x}_{d}\left(\tau\right)+\boldsymbol{b}\right).
\end{align}
Enforcing the desired output value at the terminal time $t=t_{1}$
yields 
\begin{align}
\boldsymbol{z}_{1} & =\boldsymbol{z}_{d}\left(t_{1}\right)=\boldsymbol{\mathcal{C}}\boldsymbol{x}_{d}\left(t_{1}\right)=\boldsymbol{\mathcal{C}}\boldsymbol{\mathcal{P}}\boldsymbol{x}_{d}\left(t_{1}\right)+\boldsymbol{\mathcal{C}}\boldsymbol{\mathcal{Q}}\boldsymbol{x}_{d}\left(t_{1}\right)\nonumber \\
 & =\boldsymbol{\mathcal{C}}\boldsymbol{\mathcal{P}}\boldsymbol{x}_{d}\left(t_{1}\right)+\boldsymbol{\mathcal{C}}\boldsymbol{\mathcal{Q}}\exp\left(\boldsymbol{\mathcal{Q}}\boldsymbol{\mathcal{A}}\boldsymbol{\mathcal{Q}}\left(t_{1}-t_{0}\right)\right)\boldsymbol{\mathcal{Q}}\boldsymbol{x}_{0}\nonumber \\
 & +\boldsymbol{\mathcal{C}}\boldsymbol{\mathcal{Q}}\intop_{t_{0}}^{t_{1}}d\tau\exp\left(\boldsymbol{\mathcal{Q}}\boldsymbol{\mathcal{A}}\boldsymbol{\mathcal{Q}}\left(t_{1}-\tau\right)\right)\boldsymbol{\mathcal{Q}}\left(\boldsymbol{\mathcal{A}}\boldsymbol{\mathcal{P}}\boldsymbol{x}_{d}\left(\tau\right)+\boldsymbol{b}\right).
\end{align}
That is a condition for the part $\boldsymbol{\mathcal{P}}\boldsymbol{x}_{d}\left(\tau\right)$.
Exploiting the Cayley-Hamilton theorem and proceeding as in Eq. \eqref{eq:TruncateExpQAQ}
yields
\begin{align}
\boldsymbol{z}_{1}-\boldsymbol{\mathcal{C}}\boldsymbol{\mathcal{Q}}\exp\left(\boldsymbol{\mathcal{Q}}\boldsymbol{\mathcal{A}}\boldsymbol{\mathcal{Q}}\left(t_{1}-t_{0}\right)\right)\boldsymbol{\mathcal{Q}}\boldsymbol{x}_{0}-\boldsymbol{\mathcal{C}}\boldsymbol{\mathcal{Q}}\intop_{t_{0}}^{t_{1}}d\tau\exp\left(\boldsymbol{\mathcal{Q}}\boldsymbol{\mathcal{A}}\boldsymbol{\mathcal{Q}}\left(t_{1}-\tau\right)\right)\boldsymbol{\mathcal{Q}}\boldsymbol{b}\nonumber \\
=\boldsymbol{\mathcal{C}}\boldsymbol{\mathcal{P}}\boldsymbol{x}_{d}\left(t_{1}\right)+\boldsymbol{\mathcal{C}}\boldsymbol{\mathcal{Q}}\sum_{k=0}^{n-1}\left(\boldsymbol{\mathcal{Q}}\boldsymbol{\mathcal{A}}\boldsymbol{\mathcal{Q}}\right)^{k}\boldsymbol{\tilde{\alpha}}_{k}\left(t_{1},t_{0}\right).\label{eq:Eq2175}
\end{align}
In Eq. \eqref{eq:Eq2175} we defined the $n\times1$ vectors
\begin{align}
\boldsymbol{\tilde{\alpha}}_{k}\left(t_{1},t_{0}\right) & =\intop_{t_{0}}^{t_{1}}d\tau\left(\dfrac{\left(t_{1}-\tau\right)^{k}}{k!}+\sum_{i=n}^{\infty}d_{ik}\dfrac{\left(t_{1}-\tau\right)^{i}}{i!}\right)\boldsymbol{\mathcal{Q}}\boldsymbol{\mathcal{A}}\boldsymbol{\mathcal{P}}\boldsymbol{x}_{d}\left(\tau\right).
\end{align}
The right hand side of Eq. \eqref{eq:Eq2175} can be written with
the help of the $n\left(n+1\right)\times1$ vector 
\begin{align}
\boldsymbol{\tilde{\alpha}}\left(t_{1},t_{0}\right) & =\left(\begin{array}{c}
\boldsymbol{\mathcal{P}}\boldsymbol{x}_{d}\left(t_{1}\right)\\
\boldsymbol{\alpha}_{0}\left(t_{1},t_{0}\right)\\
\vdots\\
\boldsymbol{\alpha}_{n-1}\left(t_{1},t_{0}\right)
\end{array}\right)
\end{align}
and the $m\times n\left(n+1\right)$ \textit{output controllability
matrix} 
\begin{align}
\boldsymbol{\mathcal{\tilde{K}}}_{\boldsymbol{\mathcal{C}}} & =\left(\boldsymbol{\mathcal{C}}\boldsymbol{\mathcal{P}}|\boldsymbol{\mathcal{C}}\boldsymbol{\mathcal{Q}}\boldsymbol{\mathcal{A}}\boldsymbol{\mathcal{P}}|\cdots|\boldsymbol{\mathcal{C}}\boldsymbol{\mathcal{Q}}\left(\boldsymbol{\mathcal{Q}}\boldsymbol{\mathcal{A}}\boldsymbol{\mathcal{Q}}\right)^{n-1}\boldsymbol{\mathcal{Q}}\boldsymbol{\mathcal{A}}\boldsymbol{\mathcal{P}}\right)
\end{align}
as
\begin{align}
\boldsymbol{z}_{1}-\boldsymbol{\mathcal{C}}\boldsymbol{\mathcal{Q}}\exp\left(\boldsymbol{\mathcal{Q}}\boldsymbol{\mathcal{A}}\boldsymbol{\mathcal{Q}}\left(t_{1}-t_{0}\right)\right)\boldsymbol{\mathcal{Q}}\boldsymbol{x}_{0}-\boldsymbol{\mathcal{C}}\boldsymbol{\mathcal{Q}}\intop_{t_{0}}^{t_{1}}d\tau\exp\left(\boldsymbol{\mathcal{Q}}\boldsymbol{\mathcal{A}}\boldsymbol{\mathcal{Q}}\left(t_{1}-\tau\right)\right)\boldsymbol{\mathcal{Q}}\boldsymbol{b}\nonumber \\
=\boldsymbol{\mathcal{\tilde{K}}}_{\boldsymbol{\mathcal{C}}}\boldsymbol{\tilde{\alpha}}\left(t_{1},t_{0}\right).
\end{align}
The linear map from $\boldsymbol{\tilde{\alpha}}\left(t_{1},t_{0}\right)$
to $\boldsymbol{z}_{1}$ is surjective if $\boldsymbol{\mathcal{\tilde{K}}}_{\boldsymbol{\mathcal{C}}}$
has full row rank, i.e., if 
\begin{align}
\text{rank}\left(\boldsymbol{\mathcal{\tilde{K}}}_{\boldsymbol{\mathcal{C}}}\right) & =m.\label{eq:OutputControllabilityRankCondition}
\end{align}
Thus, a nonlinear affine control system satisfying the linearizing
assumption is output controllable with linear output Eq. \eqref{eq:LinearOutput-1}
if the matrix $\boldsymbol{\mathcal{K}}_{\boldsymbol{\mathcal{C}}}$
satisfies the \textit{output controllability rank condition} Eq. \eqref{eq:OutputControllabilityRankCondition}.

With $m=n$ and $\boldsymbol{\mathcal{C}}=\boldsymbol{1}$, the notion
of output controllability reduces to the notion of full state controllability.
Indeed, note that for $\boldsymbol{\mathcal{C}}=\boldsymbol{1}$,
$\boldsymbol{\mathcal{\tilde{K}}}_{\boldsymbol{\mathcal{C}}}$ can
be written in terms of the controllability matrix for realizable trajectories
$\boldsymbol{\mathcal{\tilde{K}}}$ given by Eq. \eqref{eq:ControllabilityMatrixForRealizableTrajectories}
 as
\begin{align}
\boldsymbol{\mathcal{\tilde{K}}}_{\boldsymbol{\mathcal{C}}} & =\left(\boldsymbol{\mathcal{P}}|\boldsymbol{\mathcal{\tilde{K}}}\right).\label{eq:KCStateControllability}
\end{align}
Because $\boldsymbol{\mathcal{\tilde{K}}}$ has no components in the
direction of $\boldsymbol{\mathcal{P}}$, i.e., $\boldsymbol{\mathcal{\tilde{K}}}=\boldsymbol{\mathcal{P}}\boldsymbol{\mathcal{\tilde{K}}}+\boldsymbol{\mathcal{Q}}\boldsymbol{\mathcal{\tilde{K}}}=\boldsymbol{\mathcal{Q}}\boldsymbol{\mathcal{\tilde{K}}}$,
the matrix $\boldsymbol{\mathcal{\tilde{K}}}_{\boldsymbol{\mathcal{C}}}$
as given by Eq. \eqref{eq:KCStateControllability} has rank
\begin{align}
\text{rank}\left(\boldsymbol{\mathcal{\tilde{K}}}_{\boldsymbol{\mathcal{C}}}\right) & =p+\text{rank}\left(\boldsymbol{\mathcal{\tilde{K}}}\right)=n.
\end{align}
This proves that the rank condition for output controllability, Eq.
\eqref{eq:OutputControllabilityRankCondition}, indeed reduces, for
$\boldsymbol{\mathcal{C}}=\boldsymbol{1}$, to the rank condition
for full state controllability as given by Eq. \eqref{eq:RankConditionRealizableTrajectories}. 

Output controllability is discussed by means of two examples.

\begin{example}[Output controllability of the activator-controlled\newline FHN model]\label{ex:FHN3_2}

The model
\begin{align}
\left(\begin{array}{c}
\dot{x}\left(t\right)\\
\dot{y}\left(t\right)
\end{array}\right) & =\left(\begin{array}{c}
a_{0}+a_{1}x\left(t\right)+a_{2}y\left(t\right)\\
R\left(x\left(t\right),y\left(t\right)\right)
\end{array}\right)+\left(\begin{array}{c}
0\\
1
\end{array}\right)u\left(t\right)\label{eq:Eq251}
\end{align}
satisfies the linearizing assumption such that the constraint equation
is linear with state matrix
\begin{align}
\boldsymbol{\mathcal{Q}}\boldsymbol{\mathcal{A}} & =\left(\begin{array}{cc}
a_{1} & a_{2}\\
0 & 0
\end{array}\right),
\end{align}
see Examples \ref{ex:FHN3} and \ref{ex:FHN4} for more details. We
check for the controllability of a general desired output with $1\times2$
output matrix $\mbox{\ensuremath{\boldsymbol{\mathcal{C}}}=\ensuremath{\left(\begin{array}{cc}
 c_{1},  &  c_{2},\end{array}\right)^{T}}}$, 
\begin{align}
z_{d}\left(t\right) & =\boldsymbol{\mathcal{C}}\boldsymbol{x}_{d}\left(t\right)=c_{1}x_{d}\left(t\right)+c_{2}y_{d}\left(t\right).\label{eq:FHNOutput}
\end{align}
The $1\times6$ output controllability matrix $\boldsymbol{\mathcal{\tilde{K}}}_{\boldsymbol{\mathcal{C}}}$
becomes
\begin{align}
\boldsymbol{\mathcal{\tilde{K}}}_{\boldsymbol{\mathcal{C}}} & =\left(\boldsymbol{\mathcal{C}}\boldsymbol{\mathcal{P}}|\boldsymbol{\mathcal{C}}\boldsymbol{\mathcal{Q}}\boldsymbol{\mathcal{A}}\boldsymbol{\mathcal{P}}|\boldsymbol{\mathcal{C}}\boldsymbol{\mathcal{Q}}\boldsymbol{\mathcal{A}}\boldsymbol{\mathcal{Q}}\boldsymbol{\mathcal{A}}\boldsymbol{\mathcal{P}}\right)\nonumber \\
 & =\left(\begin{array}{cccccc}
0 & c_{2} & 0 & a_{2}c_{1} & 0 & a_{1}a_{2}c_{1}\end{array}\right).
\end{align}
The rank of $\boldsymbol{\mathcal{\tilde{K}}}_{\boldsymbol{\mathcal{C}}}$
is at most one. Example \ref{ex:FHN4} showed that the system Eq.
\eqref{eq:Eq251} is not controllable if $a_{2}=0$. In this case,
$\boldsymbol{\mathcal{\tilde{K}}}_{\boldsymbol{\mathcal{C}}}$ simplifies
to 
\begin{align}
\boldsymbol{\mathcal{\tilde{K}}}_{\boldsymbol{\mathcal{C}}} & =\left(\begin{array}{cccccc}
0 & c_{2} & 0 & 0 & 0 & 0\end{array}\right).\label{eq:Eq256}
\end{align}
Thus, $\boldsymbol{\mathcal{\tilde{K}}}_{\boldsymbol{\mathcal{C}}}$
still has rank one as long as $c_{2}\neq0$. In conclusion, a model
which is not controllable can nevertheless have a controllable output.
Although for $a_{2}=0$ in Eq. \eqref{eq:Eq251}, the inhibitor $x\left(t\right)$
evolves uncoupled from the activator dynamics, activator and inhibitor
are still coupled in the output $z_{d}\left(t\right)$, Eq. \eqref{eq:FHNOutput}.
In that way the activator $y_{d}\left(t\right)$ can counteract the
inhibitor $x_{d}\left(t\right)$ to control the desired output. If
additionally $c_{2}=0$, this is not possible, and the output is not
controllable. Indeed, for $c_{2}=0$, the output controllability matrix
Eq. \eqref{eq:Eq256} reduces to the zero matrix with vanishing rank.

\end{example}

\begin{example}[Output controllability of the SIR model]\label{ex:SIRModel3}

The controlled state equation for the SIR model was developed in Example
\ref{ex:SIRModel1} and is repeated here for convenience, 
\begin{align}
\dot{S}\left(t\right) & =-\left(\beta+u\left(t\right)\right)\frac{S\left(t\right)I\left(t\right)}{N}, & \dot{I}\left(t\right) & =\left(\beta+u\left(t\right)\right)\frac{S\left(t\right)I\left(t\right)}{N}-\gamma I\left(t\right),\nonumber \\
\dot{R}\left(t\right) & =\gamma I\left(t\right).\label{eq:SIRStateEq}
\end{align}
The controllability of the SIR model was discussed in Example \ref{ex:SIRModel2}.
We check for the controllability of a general single component desired
output with $1\times3$ output matrix $\boldsymbol{\mathcal{C}}=\left(\begin{array}{ccc}
c_{1}, & c_{2}, & c_{3}\end{array}\right)^{T}$, 
\begin{align}
z_{d}\left(t\right) & =\boldsymbol{\mathcal{C}}\boldsymbol{x}_{d}\left(t\right)=c_{1}S_{d}\left(t\right)+c_{2}I_{d}\left(t\right)+c_{3}R_{d}\left(t\right).
\end{align}
The output controllability matrix $\boldsymbol{\mathcal{\tilde{K}}}_{\boldsymbol{\mathcal{C}}}$
becomes
\begin{align}
\boldsymbol{\mathcal{\tilde{K}}}_{\boldsymbol{\mathcal{C}}} & =\left(\boldsymbol{\mathcal{C}}\boldsymbol{\mathcal{P}}|\boldsymbol{\mathcal{C}}\boldsymbol{\mathcal{Q}}\boldsymbol{\mathcal{A}}\boldsymbol{\mathcal{P}}|\boldsymbol{\mathcal{C}}\boldsymbol{\mathcal{Q}}\boldsymbol{\mathcal{A}}\boldsymbol{\mathcal{Q}}\boldsymbol{\mathcal{A}}\boldsymbol{\mathcal{P}}|\boldsymbol{\mathcal{C}}\boldsymbol{\mathcal{Q}}\boldsymbol{\mathcal{A}}\boldsymbol{\mathcal{Q}}\boldsymbol{\mathcal{A}}\boldsymbol{\mathcal{Q}}\boldsymbol{\mathcal{A}}\boldsymbol{\mathcal{P}}\right)\nonumber \\
 & =\left(\begin{array}{cccccccccccc}
\kappa_{1} & -\kappa_{1} & 0 & \gamma\kappa_{2} & -\gamma\kappa_{2} & 0 & -\frac{1}{2}\gamma^{2}\kappa_{2} & \frac{\gamma^{2}}{2}\kappa_{2} & 0 & \frac{\gamma^{3}}{4}\kappa_{2} & -\frac{1}{4}\gamma^{3}\kappa_{2} & 0\end{array}\right),\label{eq:SIROutputControllabilityMatrix}
\end{align}
with $\kappa_{1}=\frac{1}{2}\left(c_{1}-c_{2}\right)$ and $\kappa_{2}=\frac{1}{4}\left(c_{1}+c_{2}-2c_{3}\right)$.
The rank of $\boldsymbol{\mathcal{\tilde{K}}}_{\boldsymbol{\mathcal{C}}}$
is at most one. Depending on the values of the output parameters $c_{1},\,c_{2}$,
$c_{3}$, and the system parameter $\gamma$, the rank of $\boldsymbol{\mathcal{\tilde{K}}}_{\boldsymbol{\mathcal{C}}}$
changes. Two cases are discussed in detail.

First, if $c_{1}=c_{2}=0$ and $c_{3}\neq0$, then $\kappa_{1}=0$
and the output is $z_{d}\left(t\right)=c_{3}R_{d}\left(t\right)$
and prescribes the number of recovered people over time. As can be
seen from Eq. \eqref{eq:SIRStateEq}, $R\left(t\right)$ is decoupled
from the controlled part of the equations if $\gamma=0$. Indeed,
in this case $\boldsymbol{\mathcal{\tilde{K}}}_{\boldsymbol{\mathcal{N}}}$
reduces to the zero matrix with vanishing rank. In conclusion, a desired
output equal to the number $R_{d}\left(t\right)$ of recovered people
cannot be controlled if $\gamma=0$.

Second, for $c_{1}=c_{2}=c_{3}=c$ the desired output becomes
\begin{align}
z_{d}\left(t\right) & =c\left(S_{d}\left(t\right)+I_{d}\left(t\right)+R_{d}\left(t\right)\right)=cN=\text{const}.,\label{eq:SIRConstOutput}
\end{align}
with $N$ being the total number of individuals. This conservation
law can easily be derived from the system dynamics Eq. \eqref{eq:SIRStateEq}
and remains true for the controlled system. We expect that this output
is not controllable because the value of $N$ is fixed by the initial
conditions and cannot be changed by control. Indeed, if $c_{1}=c_{2}=c_{3}$
then $\kappa_{1}=0$ and $\kappa_{2}=0$. The output controllability
matrix $\boldsymbol{\mathcal{\tilde{K}}}_{\boldsymbol{\mathcal{C}}}$
becomes the zero matrix with vanishing rank, and the output Eq. \eqref{eq:SIRConstOutput}
is not controllable.

\end{example}

\section{\label{sec:OutputRealizability}Output realizability}

\subsection{\label{sub:GeneralProcedure}General procedure}

For a desired trajectory $\boldsymbol{x}_{d}\left(t\right)$ to be
exactly realizable, it must satisfy the constraint equation
\begin{align}
\boldsymbol{\mathcal{Q}}\left(\boldsymbol{x}_{d}\left(t\right)\right)\left(\boldsymbol{\dot{x}}_{d}\left(t\right)-\boldsymbol{R}\left(\boldsymbol{x}_{d}\left(t\right)\right)\right) & =\mathbf{0}.\label{eq:ConstraintEquation2}
\end{align}
This equation fixes $n-p$ components of the $n$ components of $\boldsymbol{x}_{d}\left(t\right)$.
Our convention was to choose these $n-p$ independent components as
$\boldsymbol{y}_{d}\left(t\right)=\boldsymbol{\mathcal{Q}}\left(\boldsymbol{x}_{d}\left(t\right)\right)\boldsymbol{x}_{d}\left(t\right)$,
while the $p$ independent components $\boldsymbol{z}_{d}\left(t\right)=\boldsymbol{\mathcal{P}}\left(\boldsymbol{x}_{d}\left(t\right)\right)\boldsymbol{x}_{d}\left(t\right)$
of the desired state trajectory are prescribed by the experimenter.
Equation \eqref{eq:ConstraintEquation2} becomes a non-autonomous
differential equation for $\boldsymbol{y}_{d}\left(t\right)$,
\begin{align}
\boldsymbol{\dot{y}}_{d}\left(t\right) & =\boldsymbol{\mathcal{Q}}\left(\boldsymbol{y}_{d}\left(t\right)+\boldsymbol{z}_{d}\left(t\right)\right)\boldsymbol{R}\left(\boldsymbol{y}_{d}\left(t\right)+\boldsymbol{z}_{d}\left(t\right)\right)\nonumber \\
 & +\boldsymbol{\mathcal{\dot{Q}}}\left(\boldsymbol{y}_{d}\left(t\right)+\boldsymbol{z}_{d}\left(t\right)\right)\left(\boldsymbol{y}_{d}\left(t\right)+\boldsymbol{z}_{d}\left(t\right)\right).\label{eq:ConstraintEquation3}
\end{align}
Here, $\boldsymbol{\mathcal{\dot{Q}}}$ denotes the short hand notation
\begin{align}
\boldsymbol{\mathcal{\dot{Q}}}\left(\boldsymbol{x}\left(t\right)\right) & =\dfrac{d}{dt}\boldsymbol{\mathcal{Q}}\left(\boldsymbol{x}\left(t\right)\right)=\left(\boldsymbol{\dot{x}}^{T}\left(t\right)\nabla\right)\boldsymbol{\mathcal{Q}}\left(\boldsymbol{x}\left(t\right)\right).
\end{align}
The explicit time dependence rendering Eq. \eqref{eq:ConstraintEquation3}
a non-autonomous differential equation comes from the term $\boldsymbol{z}_{d}\left(t\right)$.
The initial condition for Eq. \eqref{eq:ConstraintEquation3} is
\begin{align}
\boldsymbol{y}_{d}\left(t_{0}\right) & =\boldsymbol{\mathcal{Q}}\left(\boldsymbol{x}\left(t_{0}\right)\right)\boldsymbol{x}\left(t_{0}\right),
\end{align}
while $\boldsymbol{z}_{d}\left(t_{0}\right)$ has to satisfy
\begin{align}
\boldsymbol{z}_{d}\left(t_{0}\right) & =\boldsymbol{\mathcal{P}}\left(\boldsymbol{x}\left(t_{0}\right)\right)\boldsymbol{x}\left(t_{0}\right).
\end{align}
Because of $\boldsymbol{\mathcal{B}}^{+}\left(\boldsymbol{x}_{d}\left(t\right)\right)\boldsymbol{\mathcal{Q}}\left(\boldsymbol{x}_{d}\left(t\right)\right)=\boldsymbol{0}$,
the corresponding control signal $\boldsymbol{u}\left(t\right)$ is
given as
\begin{align}
\boldsymbol{u}\left(t\right) & =\boldsymbol{\mathcal{B}}^{+}\left(\boldsymbol{x}_{d}\left(t\right)\right)\left(\boldsymbol{\dot{x}}_{d}\left(t\right)-\boldsymbol{R}\left(\boldsymbol{x}_{d}\left(t\right)\right)\right)\nonumber \\
 & =\boldsymbol{\mathcal{B}}^{+}\left(\boldsymbol{y}_{d}\left(t\right)+\boldsymbol{z}_{d}\left(t\right)\right)\left(\boldsymbol{\dot{z}}_{d}\left(t\right)-\boldsymbol{R}\left(\boldsymbol{y}_{d}\left(t\right)+\boldsymbol{z}_{d}\left(t\right)\right)\right)\nonumber \\
 & -\boldsymbol{\mathcal{B}}^{+}\left(\boldsymbol{y}_{d}\left(t\right)+\boldsymbol{z}_{d}\left(t\right)\right)\boldsymbol{\mathcal{\dot{P}}}\left(\boldsymbol{y}_{d}\left(t\right)+\boldsymbol{z}_{d}\left(t\right)\right)\left(\boldsymbol{y}_{d}\left(t\right)+\boldsymbol{z}_{d}\left(t\right)\right).\label{eq:ControlConvention}
\end{align}
Solving Eq. \eqref{eq:ConstraintEquation3} for $\boldsymbol{y}_{d}\left(t\right)$
in terms of $\boldsymbol{z}_{d}\left(t\right)$, the term $\boldsymbol{y}_{d}\left(t\right)$
is eliminated from Eq. \eqref{eq:ControlConvention}, resulting in
a control signal expressed in terms of $\boldsymbol{z}_{d}\left(t\right)$
only. The dependence of the control signal $\boldsymbol{u}\left(t\right)$
on $\boldsymbol{z}_{d}\left(t\right)$ is in form of a functional,
\begin{align}
\boldsymbol{u}\left(t\right) & =\boldsymbol{u}\left[\boldsymbol{z}_{d}\left(t\right)\right].\label{eq:uFunctional}
\end{align}
However, the choice $\boldsymbol{z}_{d}\left(t\right)=\boldsymbol{\mathcal{P}}\left(\boldsymbol{x}_{d}\left(t\right)\right)\boldsymbol{x}_{d}\left(t\right)$
is not the only possible desired output. A general approach prescribes
an arbitrary $m$-component output
\begin{align}
\boldsymbol{z}_{d}\left(t\right) & =\boldsymbol{h}\left(\boldsymbol{x}_{d}\left(t\right)\right).\label{eq:OutputFunction}
\end{align}
The function $\boldsymbol{h}$ maps the state space $\mathbb{R}^{n}$
to a space $\mathbb{R}^{m}$. Using the constraint equation \eqref{eq:ConstraintEquation2},
one can attempt to eliminate $n-m$ components of $\boldsymbol{x}_{d}\left(t\right)$
in the control signal and obtain a control signal depending on $\boldsymbol{z}_{d}\left(t\right)$
only. If it is possible to do so, also the controlled state trajectory
$\boldsymbol{x}\left(t\right)$ can be expressed in terms of the desired
output $\boldsymbol{z}_{d}\left(t\right)$ only. The output $\boldsymbol{z}_{d}\left(t\right)$
is an \textit{exactly realizable desired output}. Clearly, not all
desired outputs can be realized, and the question arises under which
conditions it is possible to exactly realize a desired output $\boldsymbol{z}_{d}\left(t\right)$.
For example, if the dimension $m$ of the output signals is larger
than the dimension $p$ of the control signals, $m>p$, it should
be impossible to express the control signal in terms of the output.
Here, we are not able to give a definite answer to this question.
We discuss some general aspects of the problem in Section \ref{sub:OutputTrajectoryRealizabilityLinearizingAssumption},
and treat some explicit examples in Sections \ref{sub:OutputRealizationExamples}. 

A remark in order to minimize the confusion: $\boldsymbol{z}_{d}$
as given by $\boldsymbol{z}_{d}=\boldsymbol{\mathcal{P}}\left(\boldsymbol{x}_{d}\right)\boldsymbol{x}_{d}$
is an $n$-component vector, but has only $p$ independent components.
Starting with Eq. \eqref{eq:OutputFunction}, the output $\boldsymbol{z}_{d}\left(t\right)$
is regarded as a $p$-component vector with $p$ independent components,
as it is customary for outputs. Note that the convention of choosing
the part $\boldsymbol{\mathcal{P}}\boldsymbol{x}_{d}\left(t\right)$
as the desired output corresponds to the case $\boldsymbol{\mathcal{M}}=\boldsymbol{\mathcal{P}}$
and $\boldsymbol{\mathcal{N}}=\boldsymbol{\mathcal{Q}}$, which, for
constant coupling matrix $\boldsymbol{\mathcal{B}}\left(\boldsymbol{x}\right)=\boldsymbol{\mathcal{B}}$,
is equivalent to the linear output function $\boldsymbol{z}\left(t\right)=\boldsymbol{\mathcal{B}}^{T}\boldsymbol{x}\left(t\right)$.

\subsection{\label{sub:OutputTrajectoryRealizabilityLinearizingAssumption}Output
trajectory realizability leads to differential-algebraic systems}

Consider a desired output trajectory $\boldsymbol{z}_{d}\left(t\right)$
depending linearly on the desired state trajectory $\boldsymbol{x}_{d}\left(t\right)$,
\begin{align}
\boldsymbol{z}_{d}\left(t\right) & =\boldsymbol{\mathcal{C}}\boldsymbol{x}_{d}\left(t\right).\label{eq:LinearStateOutputRelation_1}
\end{align}
The desired output $\boldsymbol{z}_{d}\left(t\right)$ has $m\leq n$
independent components and $\boldsymbol{\mathcal{C}}$ is assumed
to be a constant $m\times n$ \textit{output matrix} with full rank,
\begin{align}
\text{rank}\left(\boldsymbol{\mathcal{C}}\right) & =m.
\end{align}
Equation \eqref{eq:LinearStateOutputRelation_1} is viewed as an underdetermined
system of linear equations for the desired state $\boldsymbol{x}_{d}\left(t\right)$.
See the Appendix \ref{sec:OverAndUnderdetSysOfEqs} for an introduction
in solving underdetermined systems of equations.

For the linear output given by Eq. \eqref{eq:LinearStateOutputRelation_1},
we can define two complementary projectors $\boldsymbol{\mathcal{M}}$
and $\boldsymbol{\mathcal{N}}$ by 
\begin{align}
\boldsymbol{\mathcal{M}} & =\boldsymbol{\mathcal{C}}^{+}\boldsymbol{\mathcal{C}}, & \boldsymbol{\mathcal{N}} & =\boldsymbol{1}-\boldsymbol{\mathcal{M}}.\label{eq:MNProjectors}
\end{align}
Here, the Moore-Penrose pseudo inverse $\boldsymbol{\mathcal{C}}^{+}$
of $\boldsymbol{\mathcal{C}}$ is given by
\begin{align}
\boldsymbol{\mathcal{C}}^{+} & =\boldsymbol{\mathcal{C}}^{T}\left(\boldsymbol{\mathcal{C}}\boldsymbol{\mathcal{C}}^{T}\right)^{-1}.
\end{align}
The projectors $\boldsymbol{\mathcal{M}}$ and $\boldsymbol{\mathcal{N}}$
are symmetric $n\times n$ matrices. The inverse of the $m\times m$
matrix $\boldsymbol{\mathcal{C}}\boldsymbol{\mathcal{C}}^{T}$ exists
because $\boldsymbol{\mathcal{C}}$ has full rank by assumption. The
ranks of the projectors are
\begin{align}
\text{rank}\left(\boldsymbol{\mathcal{M}}\right) & =m, & \text{rank}\left(\boldsymbol{\mathcal{N}}\right) & =n-m.
\end{align}
Multiplying $\boldsymbol{\mathcal{M}}$ and $\boldsymbol{\mathcal{N}}$
with $\boldsymbol{\mathcal{C}}$ from the left and right yields 
\begin{align}
\boldsymbol{\mathcal{M}}\boldsymbol{\mathcal{C}}^{T} & =\boldsymbol{\mathcal{C}}^{T}, & \boldsymbol{\mathcal{C}}\boldsymbol{\mathcal{M}} & =\boldsymbol{\mathcal{C}}, & \boldsymbol{\mathcal{N}}\boldsymbol{\mathcal{C}}^{T} & =\mathbf{0}, & \boldsymbol{\mathcal{C}}\boldsymbol{\mathcal{N}} & =\mathbf{0}.
\end{align}
Multiplying the state-output relation \eqref{eq:LinearStateOutputRelation_1}
by $\boldsymbol{\mathcal{C}}^{+}$ from the left gives 
\begin{align}
\boldsymbol{\mathcal{M}}\boldsymbol{x}_{d}\left(t\right) & =\boldsymbol{\mathcal{C}}^{+}\boldsymbol{z}_{d}\left(t\right).
\end{align}
Using the last equation, the desired state $\boldsymbol{x}_{d}\left(t\right)$
can be separated in two parts as
\begin{align}
\boldsymbol{x}_{d}\left(t\right) & =\boldsymbol{\mathcal{M}}\boldsymbol{x}_{d}\left(t\right)+\boldsymbol{\mathcal{N}}\boldsymbol{x}_{d}\left(t\right)=\boldsymbol{\mathcal{C}}^{+}\boldsymbol{z}_{d}\left(t\right)+\boldsymbol{\mathcal{N}}\boldsymbol{x}_{d}\left(t\right).\label{eq:StateSeparation}
\end{align}
Thus, the part $\boldsymbol{\mathcal{M}}\boldsymbol{x}_{d}\left(t\right)$
can be expressed in terms of the output $\boldsymbol{z}_{d}\left(t\right)$
while the part $\boldsymbol{\mathcal{N}}\boldsymbol{x}_{d}\left(t\right)$
is left undetermined.

In the following, we enforce the first part of the linearizing assumption,
namely, we assume constant projectors
\begin{align}
\boldsymbol{\mathcal{P}}\left(\boldsymbol{x}\right) & =\boldsymbol{\mathcal{P}}=\text{const.}, & \boldsymbol{\mathcal{Q}}\left(\boldsymbol{x}\right) & =\boldsymbol{1}-\boldsymbol{\mathcal{P}}=\text{const.},
\end{align}
in the constraint equation \eqref{eq:ConstraintEquation2}. The constraint
equation becomes 
\begin{align}
\boldsymbol{\mathcal{Q}}\boldsymbol{\dot{x}}_{d}\left(t\right) & =\boldsymbol{\mathcal{Q}}\boldsymbol{R}\left(\boldsymbol{x}_{d}\left(t\right)\right).
\end{align}
Using the projectors $\boldsymbol{\mathcal{M}}$ and $\boldsymbol{\mathcal{N}}$
introduced in Eq. \eqref{eq:StateSeparation}, the constraint equation
can be written as the \textit{output constraint equation}
\begin{align}
\boldsymbol{\mathcal{Q}}\boldsymbol{\mathcal{N}}\boldsymbol{\dot{x}}_{d}\left(t\right) & =\boldsymbol{\mathcal{Q}}\boldsymbol{R}\left(\boldsymbol{\mathcal{C}}^{+}\boldsymbol{z}_{d}\left(t\right)+\boldsymbol{\mathcal{N}}\boldsymbol{x}_{d}\left(t\right)\right)-\boldsymbol{\mathcal{Q}}\boldsymbol{\mathcal{C}}^{+}\boldsymbol{\dot{z}}_{d}\left(t\right).\label{eq:OutputConstraintEquation}
\end{align}
This is a system of equations for the part $\boldsymbol{\mathcal{N}}\boldsymbol{x}_{d}\left(t\right)$.
However, note that the rank of the matrix product $\boldsymbol{\mathcal{Q}}\boldsymbol{\mathcal{N}}$
is
\begin{align}
r & =\text{rank}\left(\boldsymbol{\mathcal{Q}}\boldsymbol{\mathcal{N}}\right)\leq\min\left(\text{rank}\left(\boldsymbol{\mathcal{Q}}\right),\text{rank}\left(\boldsymbol{\mathcal{N}}\right)\right)\nonumber \\
 & =\min\left(n-p,n-m\right).
\end{align}
In the most extreme case, $\boldsymbol{\mathcal{Q}}\boldsymbol{\mathcal{N}}=\boldsymbol{0}$
and so $r=0$, and Eq. \eqref{eq:OutputConstraintEquation} reduces
to a purely algebraic equation for $\boldsymbol{\mathcal{N}}\boldsymbol{x}_{d}\left(t\right)$,
\begin{align}
\boldsymbol{0} & =\boldsymbol{\mathcal{Q}}\boldsymbol{R}\left(\boldsymbol{\mathcal{C}}^{+}\boldsymbol{z}_{d}\left(t\right)+\boldsymbol{\mathcal{N}}\boldsymbol{x}_{d}\left(t\right)\right)-\boldsymbol{\mathcal{Q}}\boldsymbol{\mathcal{C}}^{+}\boldsymbol{\dot{z}}_{d}\left(t\right).\label{eq:OutputCOnstraintEquationQN0}
\end{align}
In general, Eq. \eqref{eq:OutputConstraintEquation} is a system differential-algebraic
equations for the part $\boldsymbol{\mathcal{N}}\boldsymbol{x}_{d}\left(t\right)$,
and the order of the differential equation depends on the rank of
$\boldsymbol{\mathcal{Q}}\boldsymbol{\mathcal{N}}$. For $m=p$, the
system consists of $r$ independent differential equations and $n-p-r$
algebraic equations. See the books \cite{Campbell1980Singular,Campbell1982Singular,mehrmann2006differentialalgebraic}
for more information about differential-algebraic equations.

Changing the order of differential equations implies consequences
for its initial conditions. For example, evaluating Equation \eqref{eq:OutputCOnstraintEquationQN0}
at the initial time $t=t_{0}$,
\begin{align}
\boldsymbol{0} & =\boldsymbol{\mathcal{Q}}\boldsymbol{R}\left(\boldsymbol{\mathcal{C}}^{+}\boldsymbol{z}_{d}\left(t_{0}\right)+\boldsymbol{\mathcal{N}}\boldsymbol{x}_{d}\left(t_{0}\right)\right)-\boldsymbol{\mathcal{Q}}\boldsymbol{\mathcal{C}}^{+}\boldsymbol{\dot{z}}_{d}\left(t_{0}\right),\label{eq:OutputConstraintEquationt0}
\end{align}
uncovers an additional relation between $\boldsymbol{x}_{d}\left(t_{0}\right)$
and $\boldsymbol{z}_{d}\left(t_{0}\right)$ which also involves the
time derivative $\boldsymbol{\dot{z}}_{d}\left(t_{0}\right)$. If
in an experiment the initial state $\boldsymbol{x}\left(t_{0}\right)=\boldsymbol{x}_{0}$
of the system can be prepared, Eq. \eqref{eq:OutputConstraintEquationt0}
yields the value for the part $\boldsymbol{\mathcal{N}}\boldsymbol{x}_{d}\left(t_{0}\right)$,
while the part $\boldsymbol{\mathcal{M}}\boldsymbol{x}_{d}\left(t_{0}\right)$
is given by
\begin{align}
\boldsymbol{\mathcal{M}}\boldsymbol{x}_{d}\left(t_{0}\right) & =\boldsymbol{\mathcal{C}}^{+}\boldsymbol{z}_{d}\left(t_{0}\right).
\end{align}
On the other hand, if the initial state of the system cannot be prepared,
Eq. \eqref{eq:OutputConstraintEquationt0} enforces an explicit relation
between $\boldsymbol{\mathcal{N}}\boldsymbol{x}_{d}\left(t_{0}\right)$
and $\boldsymbol{\dot{z}}_{d}\left(t_{0}\right)$. In general, for
$m=p$ and $r=\text{rank}\left(\boldsymbol{\mathcal{Q}}\boldsymbol{\mathcal{N}}\right)<n-p$,
$n-p-r$ additional conditions have to be satisfied by the initial
time derivatives of the desired output trajectory $\boldsymbol{z}_{d}\left(t\right)$.
We discuss output realizability with help of several examples.

\subsection{\label{sub:OutputRealizationExamples}Realizing a desired output:
Examples}

\begin{example}[Realizing a desired output for the\newline
activator-controlled FHN model]\label{ex:FHN_4}

Consider the model
\begin{align}
\left(\begin{array}{c}
\dot{x}\left(t\right)\\
\dot{y}\left(t\right)
\end{array}\right) & =\left(\begin{array}{c}
a_{0}+a_{1}x\left(t\right)+a_{2}y\left(t\right)\\
R\left(x\left(t\right),y\left(t\right)\right)
\end{array}\right)+\left(\begin{array}{c}
0\\
1
\end{array}\right)u\left(t\right),
\end{align}
with nonlinearity
\begin{align}
R\left(x,y\right) & =R\left(y\right)-x.
\end{align}
The function $R\left(y\right)=y-\frac{1}{3}y^{3}$ corresponds to
the standard FHN nonlinearity. The constraint equation is linear
\begin{align}
\dot{x}_{d}\left(t\right) & =a_{0}+a_{1}x_{d}\left(t\right)+a_{2}y_{d}\left(t\right).\label{eq:FHNConstraintEquation}
\end{align}
In Example \ref{ex:FHN3}, the conventional choice of prescribing
the activator variable $y_{d}\left(t\right)$ was applied. The constraint
equation \eqref{eq:FHNConstraintEquation} was regarded as a differential
equation for $x_{d}\left(t\right)$. Consequently, by eliminating
$x_{d}\left(t\right)$ in the control $u\left(t\right)=\dot{y}_{d}\left(t\right)-R\left(x_{d}\left(t\right),y_{d}\left(t\right)\right)$,
$u\left(t\right)$ was expressed entirely in terms of $y_{d}\left(t\right)$.
In contrast, here the output $z_{d}\left(t\right)$ is chosen as a
linear combination of activator and inhibitor, 
\begin{align}
z_{d}\left(t\right) & =h\left(x_{d}\left(t\right),y_{d}\left(t\right)\right)=c_{1}x_{d}\left(t\right)+c_{2}y_{d}\left(t\right).\label{eq:FHNMixedOutput}
\end{align}
Rearranging Eq. \eqref{eq:FHNMixedOutput} gives
\begin{align}
y_{d}\left(t\right) & =\dfrac{1}{c_{2}}\left(z_{d}\left(t\right)-c_{1}x_{d}\left(t\right)\right).\label{eq:ydIntermsOfzd}
\end{align}
Using the last relation in the constraint equation \eqref{eq:FHNConstraintEquation}
yields a linear ODE for $x_{d}$ with inhomogeneity $z_{d}$, 
\begin{align}
\dot{x}_{d}\left(t\right) & =\left(a_{1}-\frac{c_{1}a_{2}}{c_{2}}\right)x_{d}\left(t\right)+a_{0}+\frac{a_{2}}{c_{2}}z_{d}\left(t\right).\label{eq:Eq255}
\end{align}
Its solution is, with $\kappa=a_{1}-\frac{a_{2}c_{1}}{c_{2}}$,
\begin{align}
x_{d}\left(t\right) & =x_{d}\left(t_{0}\right)e^{\kappa\left(t-t_{0}\right)}+e^{\kappa t}\frac{a_{2}}{c_{2}}\intop_{t_{0}}^{t}\exp\left(-\kappa\tau\right)z_{d}\left(\tau\right)\,d\tau\nonumber \\
 & +\frac{a_{0}}{\kappa}\left(e^{\kappa\left(t-t_{0}\right)}-1\right),\label{eq:Eq264}
\end{align}
Using relation Eq. \eqref{eq:ydIntermsOfzd} for $y_{d}$ together
with the Eq. \eqref{eq:Eq264}, $x_{d}$ and $y_{d}$ can be eliminated
in terms of $z_{d}\left(t\right)$ from the control signal $u\left(t\right)$.
The result is an expression in terms of the desired output $z_{d}$
only (not shown).

One remark about the initial condition for the desired state $z_{d}\left(t\right)$.
For any desired trajectory $\boldsymbol{x}_{d}\left(t\right)$ to
be exactly realizable, its initial condition $\boldsymbol{x}_{d}\left(t_{0}\right)$
must agree with the initial condition $\boldsymbol{\ensuremath{x}}\left(t_{0}\right)$
of the controlled state $\boldsymbol{\ensuremath{x}}\left(t\right)$.
This naturally restricts the initial value of the desired output to
satisfy $z_{d}\left(t_{0}\right)=c_{1}x_{d}\left(t_{0}\right)+c_{2}y_{d}\left(t_{0}\right)$.

In conclusion, the control as well as the desired state trajectory
$\boldsymbol{x}_{d}\left(t\right)$ is expressed solely in terms of
the desired output $z_{d}\left(t\right)$. A numerical simulation
of the controlled model shown in Fig. \ref{fig:FHN4} demonstrates
the successful realization of the desired output $z_{d}\left(t\right)=4\sin\left(2t\right)$.
Initially at time $t_{0}=0$, the state is set to $x_{0}=y_{0}=0$,
which complies with the initial value of the desired output $z_{d}\left(0\right)=0$.
A comparison of the desired output $z_{d}\left(t\right)$ with the
output $z\left(t\right)$ obtained by numerical simulations of the
controlled system demonstrates perfect agreement (see Fig. \ref{fig:FHN4}
top left), and a plot of $z\left(t\right)-z_{d}\left(t\right)$ reveals
differences within numerical precision (see Fig. \ref{fig:FHN4} top
right). The controlled state trajectories $x\left(t\right)$ and $y\left(t\right)$
are shown in Fig. \ref{fig:FHN4} bottom left, and the control signal
is shown in Fig. \ref{fig:FHN4} bottom right.

\begin{minipage}{1.0\linewidth}%
\begin{center}%
\includegraphics[scale=0.6]{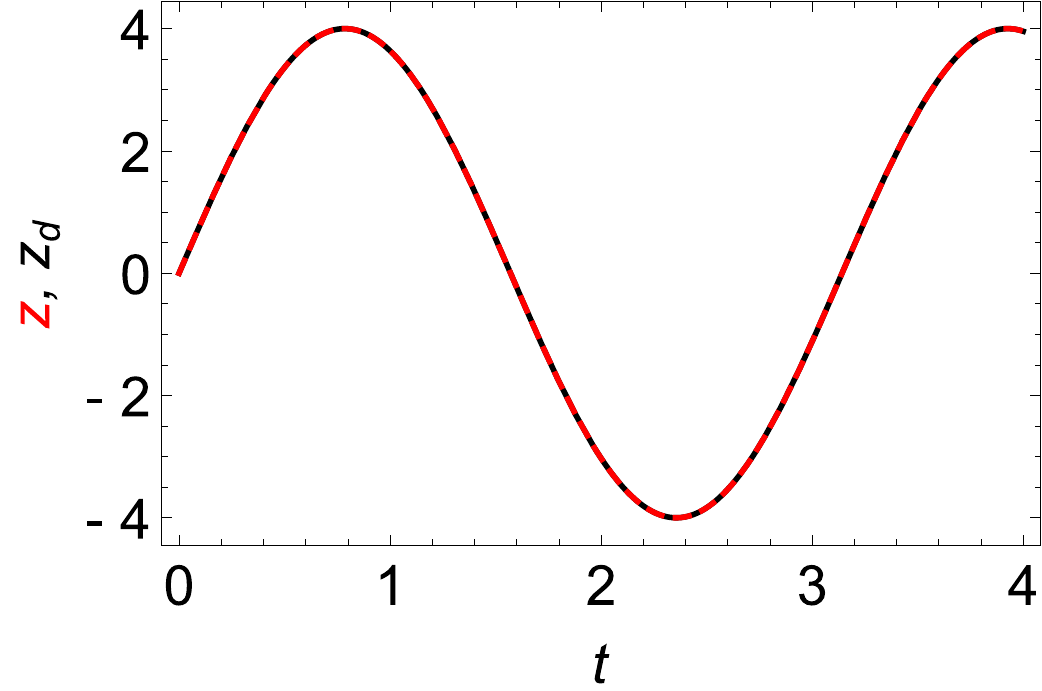}%
\includegraphics[scale=0.62]{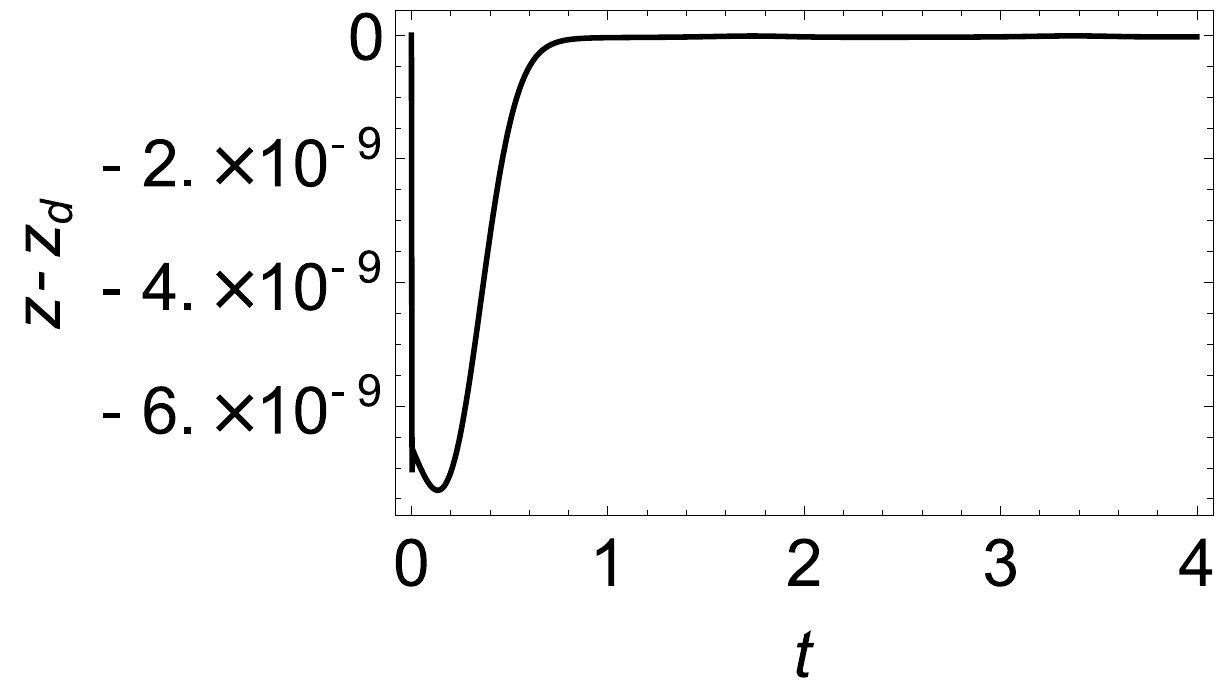}
\includegraphics[scale=0.6]{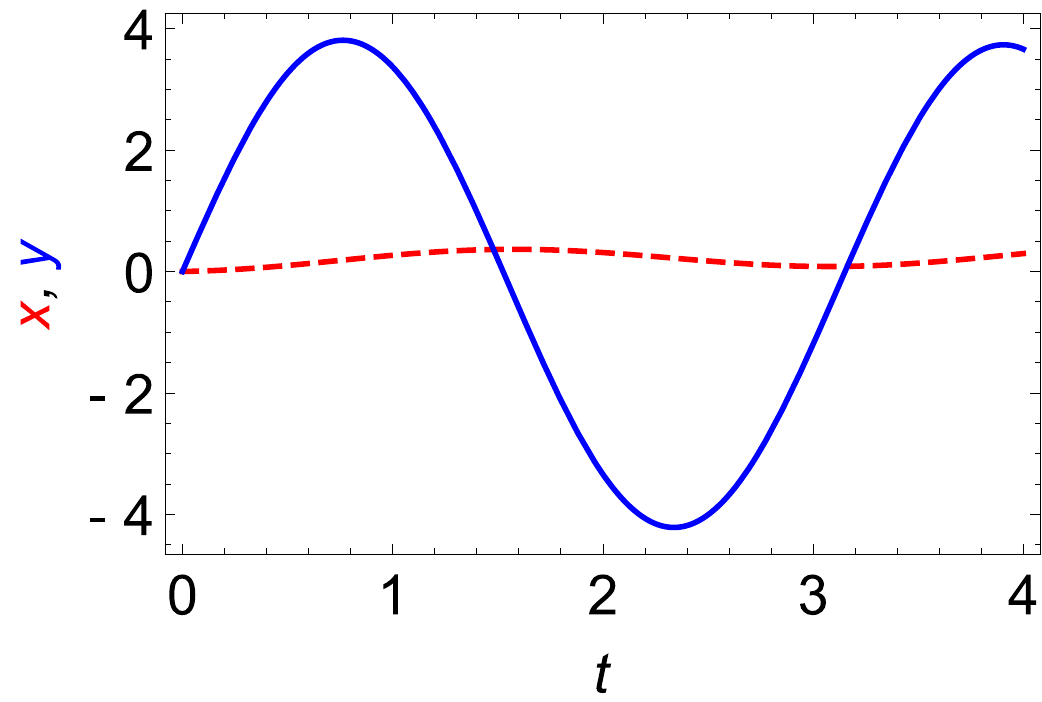}\hspace{1.2cm}%
\includegraphics[scale=0.6]{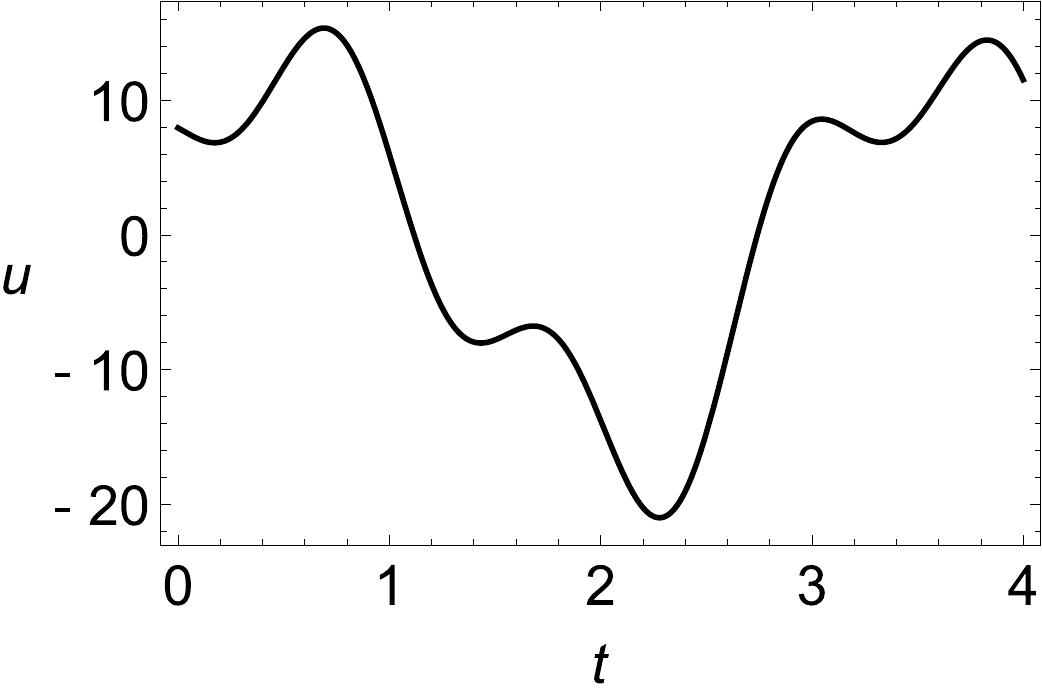}%
\captionof{figure}[Realizing a desired output in the activator-controlled FHN model]{\label{fig:FHN4}Realizing a desired output in the activator-controlled FHN model. The numerical result $z$ (red dashed line) for the output lies on top of the desired output trajectory $z_{d}$ (black line), see top left figure. The difference $z-z_d$ is within the range of numerical precision (top right). The bottom left figure shows the corresponding state trajectories $x$ (red dashed line) and $y$ (blue line) and the control $u$ (bottom right).}%
%"/home/jakob/svnco/Control/FitzHughNagumo/FHNLocalDynamics"%
\end{center}%
\end{minipage}

\end{example}

\begin{example}[Controlling the number of infected individuals in the SIR model]\label{ex:SIRModel4}

The controlled state equation for the SIR model was developed in Example
\ref{ex:SIRModel1}, and its output controllability was discussed
in Example \ref{ex:SIRModel3}. An uncontrolled time evolution is
assumed for all times $t<t_{0}$, upon which the control is switched
on. Starting at time $t=t_{0}$, the number of infected people over
time is prescribed. The desired output is 
\begin{align}
z_{d}\left(t\right) & =I_{d}\left(t\right).\label{eq:Eq2210}
\end{align}
The constraint equation consists of two independent equations 
\begin{align}
\left(\begin{array}{c}
\frac{1}{2}\left(\gamma z_{d}\left(t\right)+\dot{z}_{d}\left(t\right)+\dot{S}_{d}\left(t\right)\right)\\
-\gamma z_{d}\left(t\right)+\dot{R}_{d}\left(t\right)
\end{array}\right) & =\left(\begin{array}{c}
0\\
0
\end{array}\right).\label{eq:SIRConstraint}
\end{align}
The constraint equation is considered as two differential equations
for $S_{d}\left(t\right)$ and $R_{d}\left(t\right)$. Their solutions
are readily obtained as
\begin{align}
S_{d}\left(t\right) & =-\gamma\intop_{t_{0}}^{t}d\tau z_{d}\left(\tau\right)-z_{d}\left(t\right)+S_{d}\left(t_{0}\right)+z_{d}\left(t_{0}\right),\label{eq:Eq2212}\\
R_{d}\left(t\right) & =R_{d}\left(t_{0}\right)+\gamma\intop_{t_{0}}^{t}d\tau z_{d}\left(\tau\right).\label{eq:Eq2213}
\end{align}
Eqs. \eqref{eq:Eq2212} and \eqref{eq:Eq2213} express $S_{d}\left(t\right)$
and $R_{d}\left(t\right)$ solely in terms of the desired output $z_{d}\left(t\right)$
and the initial conditions. For any desired trajectory $\boldsymbol{x}_{d}\left(t\right)$
to be exactly realizable, its initial condition $\boldsymbol{x}_{d}\left(t_{0}\right)$
must comply with the initial condition $\boldsymbol{\ensuremath{x}}\left(t_{0}\right)$.
For the initial conditions of $R_{d}$ and $S_{d}$ follows
\begin{align}
R_{d}\left(t_{0}\right) & =R\left(t_{0}\right), & S_{d}\left(t_{0}\right) & =S\left(t_{0}\right),
\end{align}
while from $z_{d}\left(t\right)=I_{d}\left(t\right)$ follows 
\begin{align}
z_{d}\left(t_{0}\right) & =I\left(t_{0}\right),\label{eq:Eq2211}
\end{align}
with $I\left(t_{0}\right)$ being the number of infected people at
time $t=t_{0}$ when control measures are started.

The solution for the control signal realizing a desired trajectory
$\boldsymbol{x}_{d}\left(t\right)$ is 
\begin{align}
u\left(t\right) & =\left(\boldsymbol{B}^{T}\left(\boldsymbol{x}_{d}\left(t\right)\right)\boldsymbol{B}\left(\boldsymbol{x}_{d}\left(t\right)\right)\right)^{-1}\boldsymbol{B}^{T}\left(\boldsymbol{x}_{d}\left(t\right)\right)\left(\boldsymbol{\dot{x}}_{d}\left(t\right)-\boldsymbol{R}\left(\boldsymbol{x}_{d}\left(t\right)\right)\right)\nonumber \\
 & =N\frac{\gamma I_{d}\left(t\right)+\dot{I}_{d}\left(t\right)-\dot{S}_{d}\left(t\right)}{2I_{d}\left(t\right)S_{d}\left(t\right)}-\beta,
\end{align}
and using the solutions for $S_{d}\left(t\right)$ and $R_{d}\left(t\right)$
in terms of $z_{d}\left(t\right)$, the control signal becomes 
\begin{align}
u\left(t\right) & =N\frac{\gamma z_{d}\left(t\right)+\dot{z}_{d}\left(t\right)}{z_{d}\left(t\right)\left(S\left(t_{0}\right)+I\left(t_{0}\right)-z_{d}\left(t\right)-\gamma\intop_{t_{0}}^{t}d\tau z_{d}\left(\tau\right)\right)}-\beta.\label{eq:SIRControlSolution}
\end{align}
The desired number of infected individuals $z_{d}\left(t\right)$
shall follow a parabolic time evolution,
\begin{align}
z_{d}\left(t\right) & =b_{2}t^{2}+b_{1}t+b_{0}.
\end{align}
Three conditions are necessary to determine the three constants $b_{0},\,b_{1}$,
and $b_{2}$. The first condition follows from Eq. \eqref{eq:Eq2211}.
Second, the number of infected individuals shall vanish at time $t=t_{1}$,
\begin{align}
z_{d}\left(t_{1}\right) & =0,
\end{align}
such that $t_{1}-t_{0}$ is the duration of the epidemic. To obtain
a third relation, we demand that initially, the control signal vanishes.
Evaluating Eq. \eqref{eq:SIRControlSolution} at $t=t_{0}$ yields
\begin{align}
u\left(t_{0}\right) & =N\frac{\gamma I\left(t_{0}\right)+\dot{z}_{d}\left(t_{0}\right)}{I\left(t_{0}\right)S\left(t_{0}\right)}-\beta=0.\label{eq:Eq2221}
\end{align}
This relation can be used to obtain a relation for $\dot{z}_{d}\left(t_{0}\right)$
as
\begin{align}
\dot{z}_{d}\left(t_{0}\right) & =\dfrac{\beta}{N}I\left(t_{0}\right)S\left(t_{0}\right)-\gamma I\left(t_{0}\right).
\end{align}
Equation \eqref{eq:Eq2221} guarantees a smooth transition of the
time-dependent transmission rate $\beta\left(t\right)=\beta+u\left(t\right)$
across $t=t_{0}$.

Figure \ref{fig:ControlledSIR} shows a numerical solution. Up to
time $t=t_{0}$, the system evolves uncontrolled, upon which all initial
state values $S\left(t_{0}\right)$, $I\left(t_{0}\right)$, and $R\left(t_{0}\right)$
are measured. Starting at time $t_{0}=10$, the control signal $u\left(t\right)$,
Eq. \eqref{eq:SIRControlSolution}, acts on the system. To prevent
an unphysical negative transmission rate $\beta\left(t\right)=\beta+u\left(t\right)$,
the control $u\left(t\right)$ is clipped,
\begin{align}
\hat{u}\left(t\right) & =\begin{cases}
u\left(t\right), & u\left(t\right)>-\beta,\\
-\beta, & u\left(t\right)\leq-\beta.
\end{cases}
\end{align}
As can be seen in Fig. \ref{fig:ControlledSIR} bottom right, $\beta+u\left(t\right)$
reaches zero at an approximate time $\tilde{t}_{1}\approx56$, upon
which the system evolves again uncontrolled. At this time, the epidemic
has reached a reproductive number (see Example \ref{ex:SIRModel1})
\begin{align}
R_{0} & =\frac{\beta+u\left(\tilde{t}_{1}\right)}{\gamma}=0<1,
\end{align}
and further spreading of the epidemic is prevented. Comparison of
the controlled output $z\left(t\right)=I\left(t\right)$ with its
desired counterpart $z_{d}\left(t\right)=I_{d}\left(t\right)$ shows
perfect agreement for times $t_{0}<t<\tilde{t}_{1}$ when control
measures are operative, see bottom left of Fig. \ref{fig:ControlledSIR}.
Comparing the left and right top figures of Fig. \ref{fig:ControlledSIR}
reveals a less dramatic epidemic in case of control (top right) than
in case without control (top left), with a lower maximum number of
infected individuals $I\left(t\right)$ (red) and a smaller final
number of recovered individuals $R\left(t\right)$ (black). Note that
$R\left(t\right)$ is equivalent to the cumulative number of peoples
affected by the epidemic.

While no exact analytical solution is known for the uncontrolled SIR
model, we easily managed to find an exact analytical solution for
the control as well as for the controlled state over time. This simple
analytical approach provides statements as ``If the number of infected
individuals $\Delta t$ days from now shall not exceed $I_{\Delta t}$,
the transmission rate has to be lowered by $\Delta\beta$ within the
next $\Delta t_{1}$ days'' without much computational effort. It
is a way to predict the effectiveness versus cost of control measures.
Of course, application of this result to real world systems requires
a model for the cost of quarantine measures or vaccination programs
and their impact on the transmission rate $\beta\left(t\right)$.

\begin{minipage}{1.0\linewidth}
\begin{center}
\includegraphics[scale=0.6]{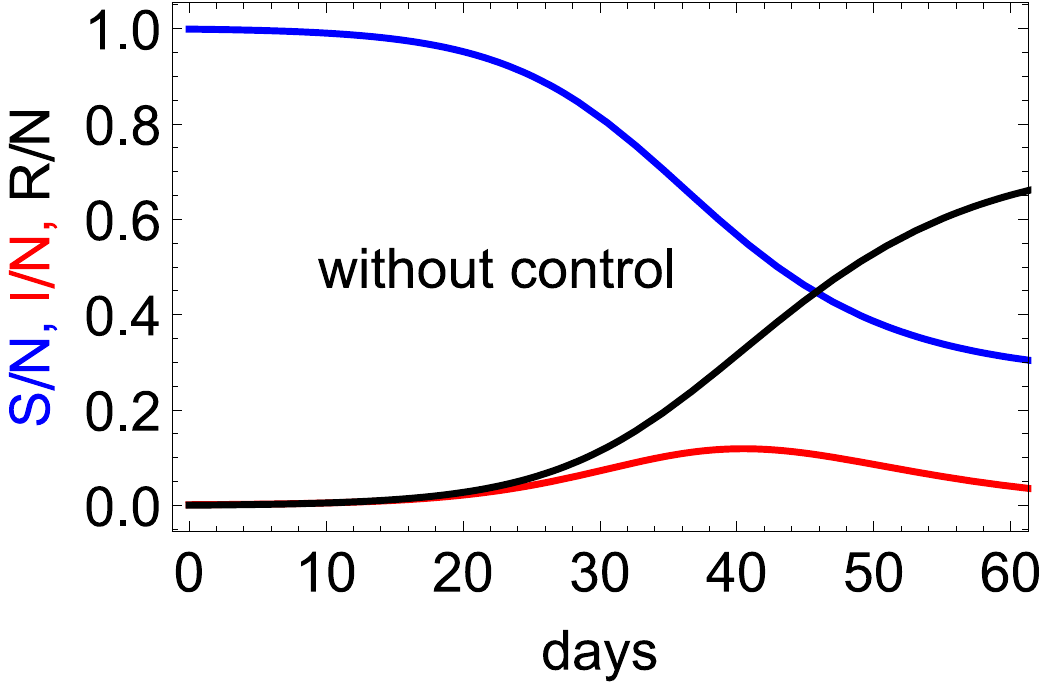}
\includegraphics[scale=0.6]{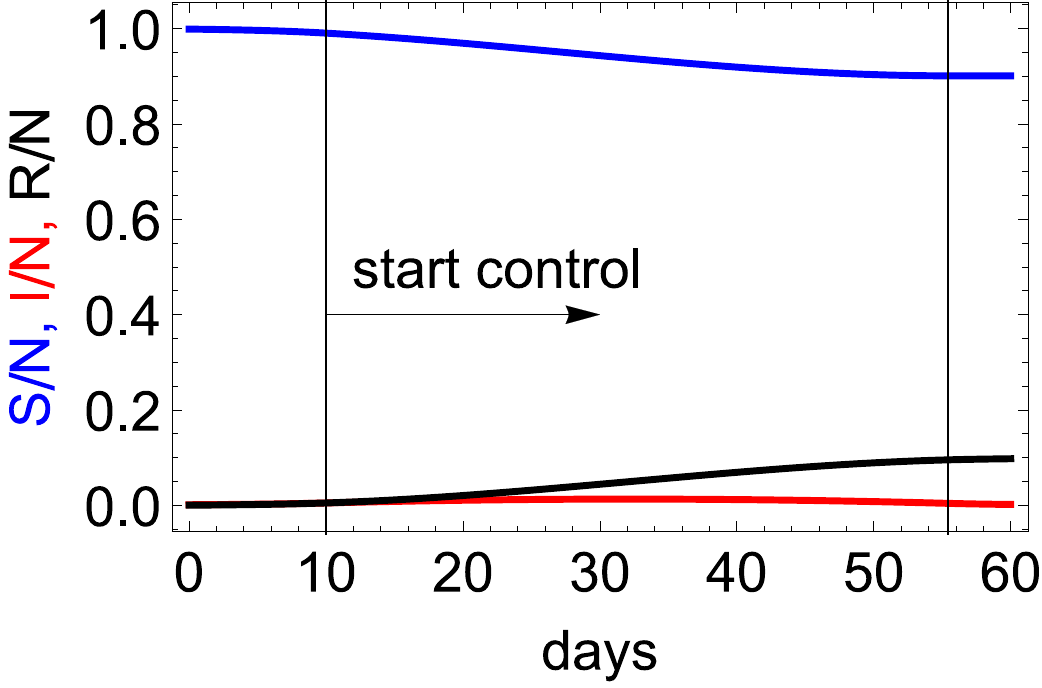}
\includegraphics[scale=0.62]{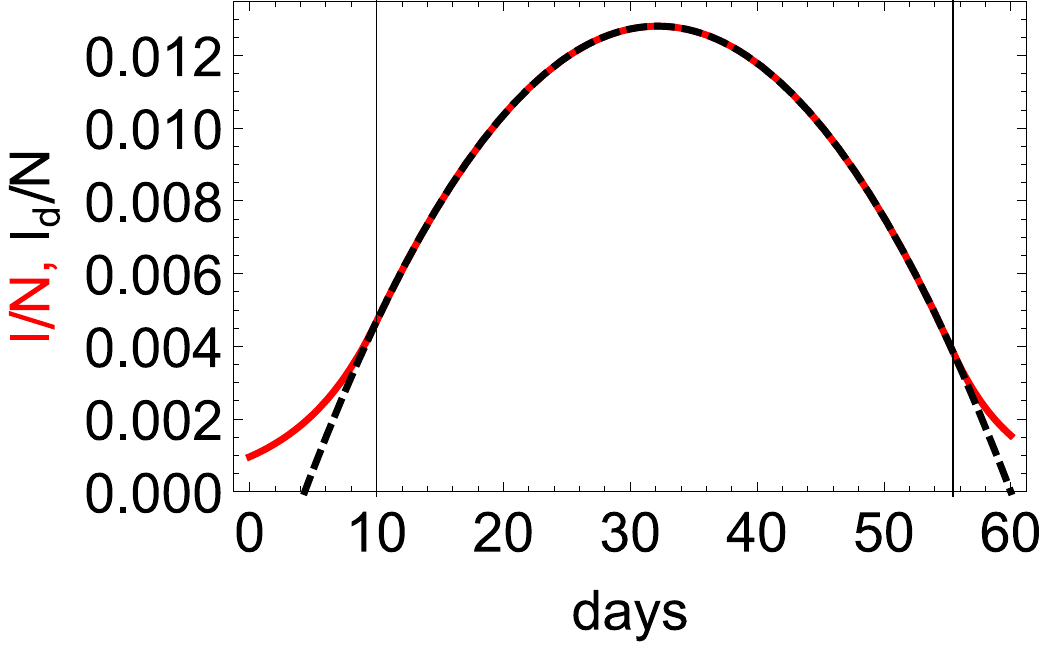}
\includegraphics[scale=0.58]{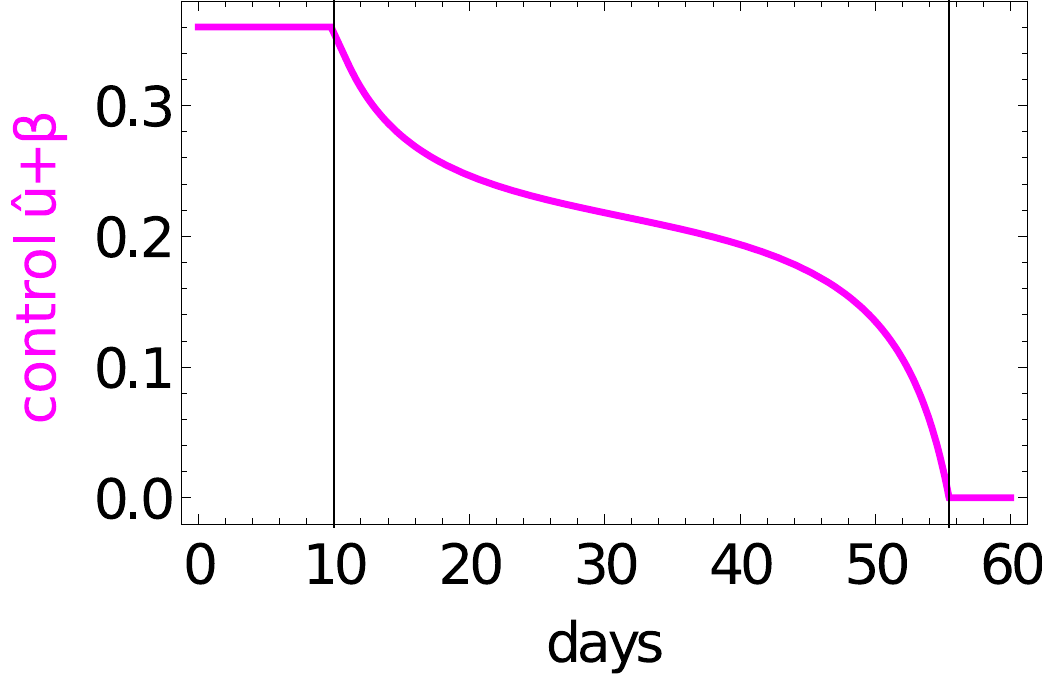}
\captionof{figure}[Control of an epidemic in the SIR model]{\label{fig:ControlledSIR}Control of an epidemic in the SIR model. Without control (top left), much more individuals become infected (red) than with control measures starting at $t=10$ (top right). A comparison of the desired output $z_d \left( t\right) = I_d \left( t \right)$ (black dashed line) of infected individuals with the actual output trajectory $z \left ( t \right) = I \left( t \right)$ (red solid line) of the controlled dynamical system reveals perfect agreement for times $t_{0}<t<\tilde{t}_{1}$ when control measures are operative   (bottom left). The bottom right figure shows the control signal which is clipped such that the time-dependent transmission rate $\beta \left( t \right) = \beta + \hat{u} \left( t \right) > 0$ is always positive.}
%"/home/jakob/svn/Control/SIR/SIRWithVaccination.nb"
\end{center}
\end{minipage}

\end{example}\begin{example}[Activator as output for the inhibitor-controlled\newline
FHN model]\label{ex:FHN_5}

Consider the model with coupling vector $\boldsymbol{B}=\left(\begin{array}{cc}
1, & 0\end{array}\right)^{T}$
\begin{align}
\left(\begin{array}{c}
\dot{x}\left(t\right)\\
\dot{y}\left(t\right)
\end{array}\right) & =\left(\begin{array}{c}
a_{0}+a_{1}x\left(t\right)+a_{2}y\left(t\right)\\
R\left(x\left(t\right),y\left(t\right)\right)
\end{array}\right)+\left(\begin{array}{c}
1\\
0
\end{array}\right)u\left(t\right)
\end{align}
and with nonlinearity
\begin{align}
R\left(x,y\right) & =R\left(y\right)-x.
\end{align}
The function $R\left(y\right)=y-\frac{1}{3}y^{3}$ corresponds to
the standard FHN nonlinearity. Example \ref{ex:FHN3_1} applied the
conventional choice and prescribed the inhibitor variable $x_{d}\left(t\right)$
as the desired output, while $y_{d}\left(t\right)$ was determined
as the solution to the corresponding constraint equation. In contrast,
here the desired output is given by the activator $y_{d}\left(t\right)$
\begin{align}
z_{d}\left(t\right) & =y_{d}\left(t\right).
\end{align}
The control signal in terms of the desired trajectory $\boldsymbol{x}_{d}\left(t\right)=\left(\begin{array}{cc}
x_{d}\left(t\right), & y_{d}\left(t\right)\end{array}\right)^{T}$ is
\begin{align}
u\left(t\right) & =\dot{x}_{d}\left(t\right)-a_{0}-a_{1}x_{d}\left(t\right)-a_{2}y_{d}\left(t\right).\label{eq:FHNControl_1}
\end{align}
The constraint equation for $\boldsymbol{x}_{d}\left(t\right)$ becomes
a nonlinear differential equation for $z_{d}\left(t\right)=y{}_{d}\left(t\right)$,
\begin{align}
\dot{z}_{d}\left(t\right) & =R\left(z_{d}\left(t\right)\right)-x_{d}\left(t\right).\label{eq:FHNConstraintEquation_1}
\end{align}
To realize the desired output $y_{d}\left(t\right)$, any reference
to the inhibitor $x_{d}\left(t\right)$ has to be eliminated from
the control signal Eq. \eqref{eq:FHNControl_1}. To achieve that,
the constraint equation \eqref{eq:FHNConstraintEquation_1} must be
solved for $x_{d}\left(t\right)$ in terms of the desired output $z_{d}\left(t\right)$.
This is a very simple task because Eq. \eqref{eq:FHNConstraintEquation_1}
is a linear algebraic equation for $x_{d}\left(t\right)$. The solution
is
\begin{align}
x_{d}\left(t\right) & =R\left(z_{d}\left(t\right)\right)-\dot{z}_{d}\left(t\right).\label{eq:InhibitorBehavior}
\end{align}
Using the last relation, $x_{d}\left(t\right)$ can be eliminated
from the control signal Eq. \eqref{eq:FHNControl_1} to get
\begin{align}
u\left(t\right) & =\dot{x}_{d}\left(t\right)-a_{0}-a_{1}x_{d}\left(t\right)-a_{2}z_{d}\left(t\right)\nonumber \\
 & =R'\left(z_{d}\left(t\right)\right)\dot{z}_{d}\left(t\right)-\ddot{z}_{d}\left(t\right)-a_{0}-a_{1}z_{d}\left(t\right)-a_{2}\left(R\left(z_{d}\left(t\right)\right)-\dot{z}_{d}\left(t\right)\right).\label{eq:Eq2239}
\end{align}
In conclusion, the control signal $u\left(t\right)$ as well as the
desired state $\boldsymbol{x}_{d}\left(t\right)$ is expressed solely
in terms of the desired output $z_{d}\left(t\right)$. Although the
system does not satisfy the linearizing assumption because the constraint
equation is a nonlinear differential equation, only a linear algebraic
equation had to be solved. Thus, linear structures underlying nonlinear
control systems may exist independently of the linearizing assumption.
Interestingly, the approach of open loop control proposed here yields
a similar result for the control as feedback linearization, see e.g.
\cite{khalil2002nonlinear}. This hints at deep connections between
our approach and feedback linearization. The framework of exactly
realizable trajectories might open up a way to generalize feedback
linearization to open loop control systems.

A remark about the initial conditions. For exactly realizable trajectories
the initial state of the desired trajectory must be equal to the initial
system state, $\boldsymbol{x}_{d}\left(t_{0}\right)=\boldsymbol{x}\left(t_{0}\right)$.
Due to Eq. \eqref{eq:InhibitorBehavior}, the initial value for $x_{d}$
is fully determined by the initial value of the desired output $z_{d}\left(t_{0}\right)$
and its time derivative $\dot{z}_{d}\left(t_{0}\right)$. For a fixed
desired output trajectory $z_{d}\left(t\right)$, the system must
be prepared in the initial state
\begin{align}
x\left(t_{0}\right) & =R\left(z_{d}\left(t_{0}\right)\right)-\dot{z}_{d}\left(t_{0}\right),\label{eq:Eq2226}\\
y\left(t_{0}\right) & =z_{d}\left(t_{0}\right).\label{eq:Eq2227}
\end{align}
On the other hand, if the system cannot be prepared in a certain initial
state, Eq. \eqref{eq:InhibitorBehavior} imposes an additional condition
on the desired output trajectory $z_{d}\left(t\right)$. In fact,
not only is the initial value $z_{d}\left(t_{0}\right)$ prescribed
by Eq. \eqref{eq:Eq2227}, but also the initial value of the time
derivative $\dot{z}_{d}\left(t_{0}\right)$ is fixed by Eq. \eqref{eq:Eq2226}.

Figure \ref{fig:FHN_5} shows the result of a numerical simulation
of the controlled FHN model with the prescribed activator 
\begin{align}
y_{d}\left(t\right) & =4\sin\left(2t\right)\label{eq:DesiredActivator}
\end{align}
as the desired output trajectory. At the initial time $t=t_{0}=0$,
the system is prepared in a state such that Eqs. \eqref{eq:Eq2226}
and \eqref{eq:Eq2227} are satisfied,
\begin{align}
\left(\begin{array}{cc}
x_{0}, & y_{0}\end{array}\right)^{T} & =\left(\begin{array}{cc}
-8, & 0\end{array}\right)^{T}.
\end{align}
Numerically solving the controlled system and comparing the controlled
state trajectories $\boldsymbol{x}\left(t\right)$ with the corresponding
desired reference trajectories reveals a perfect agreement within
numerical precision, see the bottom row of Fig. \ref{fig:FHN_5}.

\begin{minipage}{1.0\linewidth}
\begin{center}
\includegraphics[scale=0.6]{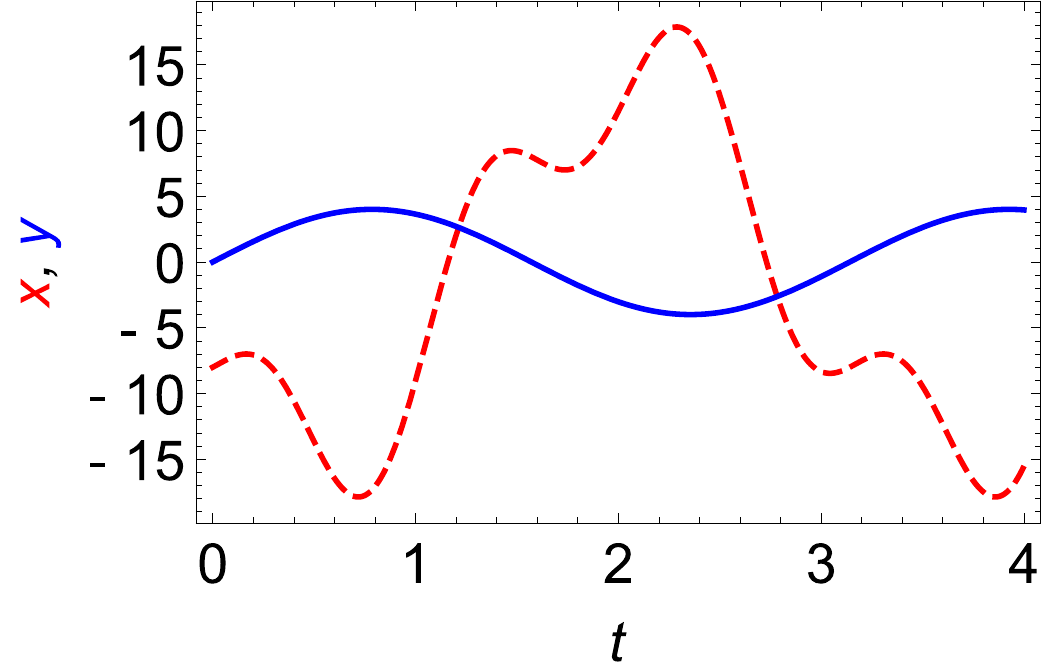}\hspace{1.2cm}
\includegraphics[scale=0.58]{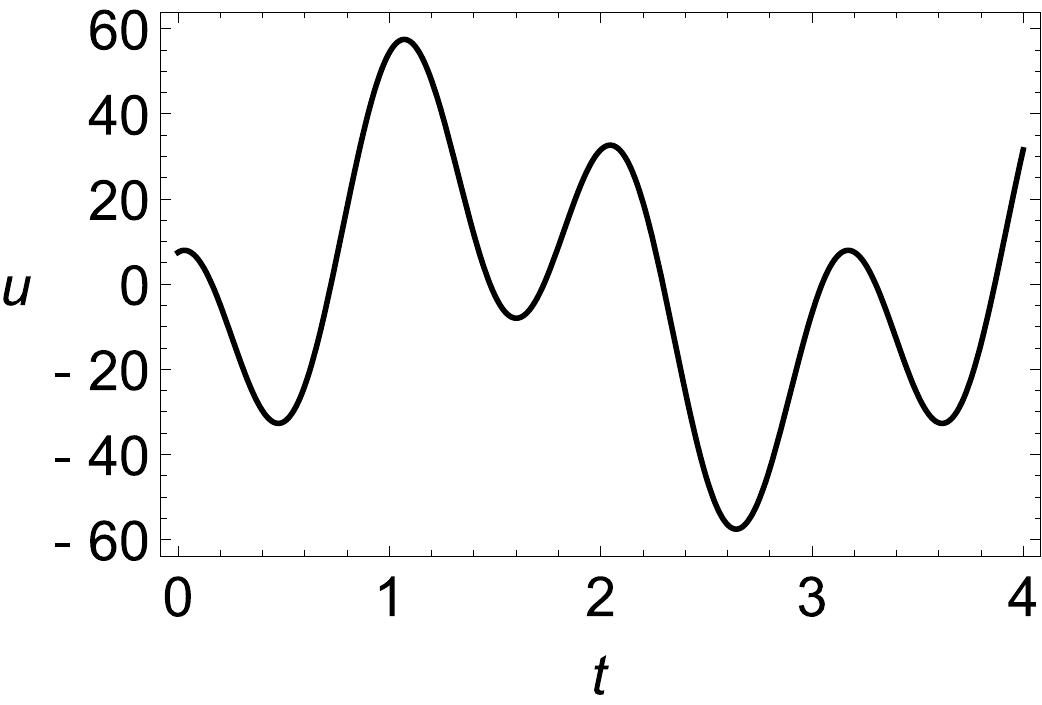}
\includegraphics[scale=0.55]{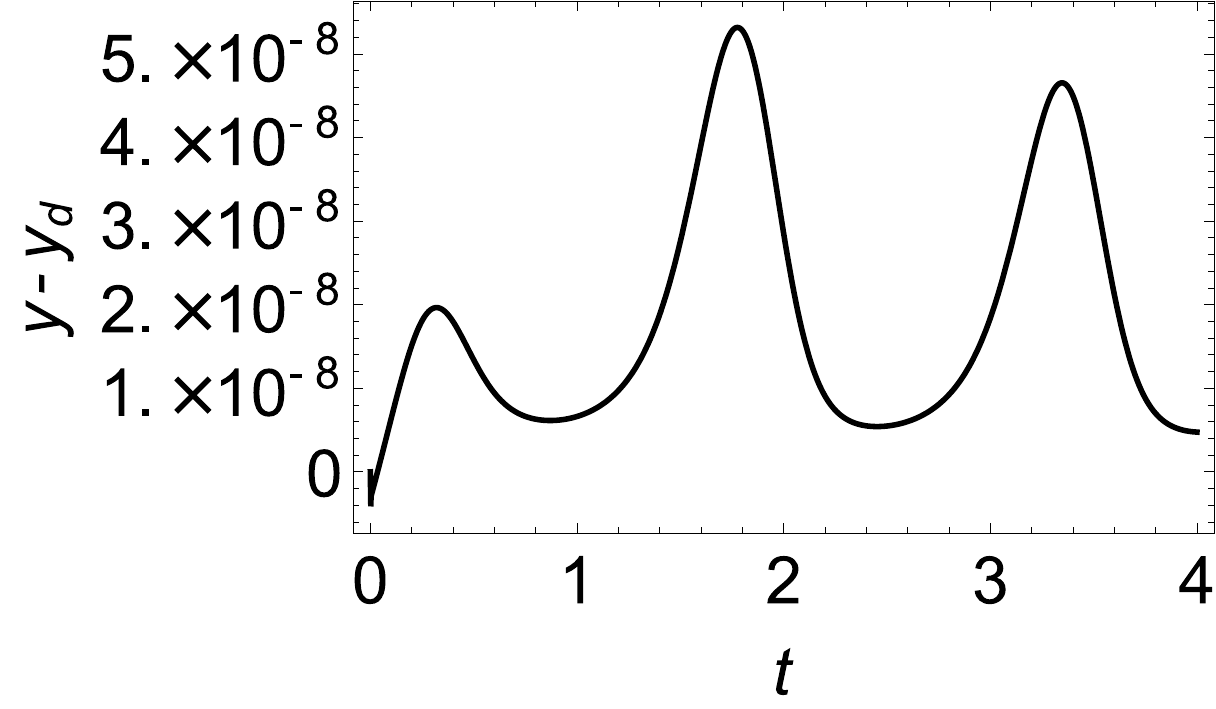}
\includegraphics[scale=0.58]{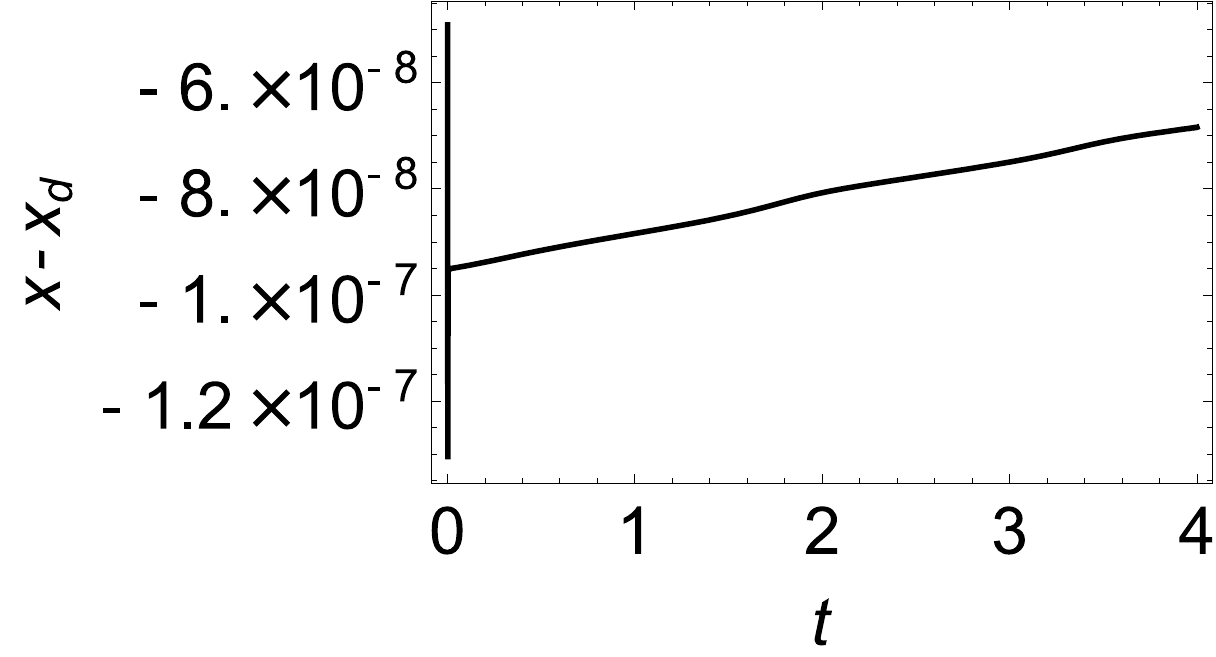}
\captionof{figure}[Inhibitor-controlled FHN model with activator as the desired output]{\label{fig:FHN_5}Inhibitor-controlled FHN model with activator $y_d$ as the desired output trajectory. Top left: the desired activator (blue solid line) is prescribed as in Eq. \eqref{eq:DesiredActivator}, while the desired inhibitor $x_d$ (red dashed line) behaves as given by Eq. \eqref{eq:InhibitorBehavior}. The control $u$ is shown in the top right panel. Comparing the difference between desired activator (bottom left) and inhibitor  (bottom right) with the corresponding controlled time evolution obtained from a numerical solution demonstrates agreement within numerical precision.}
%"/home/jakob/svnco/Control/FitzHughNagumo/FHNLocalDynamics"
\end{center}
\end{minipage}

\end{example}

\begin{example}[Modified Oregonator model]\label{ex:Oregonator}

The modified Oregonator model is a model for the light sensitive Belousov-Zhabotinsky
reaction \cite{krug1990analysis,field1972oscillations,field1974oscillations}.
In experiments, the intensity of illuminated light is used to control
the system. The Belousov-Zhabotinsky reaction has been used as an
experimental play ground for ideas related to the control of complex
systems, see e.g. \cite{mikhailov2006control} for examples. The system
equations for the activator $y$ and inhibitor $x$ read as
\begin{align}
\left(\begin{array}{c}
\dot{x}\left(t\right)\\
\dot{y}\left(t\right)
\end{array}\right) & =\left(\begin{array}{c}
y\left(t\right)-x\left(t\right)\\
\dfrac{1}{\tilde{\epsilon}}\left(y\left(t\right)\left(1-y\left(t\right)\right)+fx\left(t\right)\dfrac{q-y\left(t\right)}{q+y\left(t\right)}\right)
\end{array}\right)\nonumber \\
 & +\left(\begin{array}{c}
0\\
\dfrac{1}{\tilde{\epsilon}}\dfrac{q-y\left(t\right)}{q+y\left(t\right)}
\end{array}\right)u\left(t\right).\label{eq:Oregonator}
\end{align}
The control signal $u\left(t\right)$ is proportional to the applied
light intensity. In experiments, the inhibitor is visible and can
be recorded with a camera. The measured gray scale depends linearly
on the inhibitor and is used as the output $z$,
\begin{align}
z\left(t\right) & =h\left(x\left(t\right)\right)=I_{0}+cx\left(t\right).\label{eq:OregonatorOutput}
\end{align}
For a desired trajectory to be exactly realizable, it has to satisfy
the linear constraint equation,
\begin{align}
\dot{x}_{d}\left(t\right) & =y_{d}\left(t\right)-x_{d}\left(t\right).\label{eq:Eq2229}
\end{align}
Equation \eqref{eq:Eq2229} is solved for $y_{d}\left(t\right)$ to
obtain
\begin{align}
y_{d}\left(t\right) & =\dot{x}_{d}\left(t\right)+x_{d}\left(t\right)=\dfrac{1}{c}\dot{z}_{d}\left(t\right)+\dfrac{1}{c}\left(z_{d}\left(t\right)-I_{0}\right).\label{eq:Eq2230}
\end{align}
The inhibitor $x_{d}\left(t\right)$ was substituted with the desired
output $z_{d}\left(t\right)$ given by Eq. \eqref{eq:OregonatorOutput}.
The control signal $u\left(t\right)$ can be expressed entirely in
terms of the desired output $z_{d}\left(t\right)$ as
\begin{align}
u\left(t\right) & =\dfrac{q+y_{d}\left(t\right)}{q-y_{d}\left(t\right)}\left(\tilde{\epsilon}\dot{y}_{d}\left(t\right)+y_{d}\left(t\right)\left(y_{d}\left(t\right)-1\right)\right)-fx_{d}\left(t\right)\nonumber \\
 & =\dfrac{\tilde{\epsilon}}{c}\dfrac{cq+\dot{z}_{d}\left(t\right)+z_{d}\left(t\right)-I_{0}}{cq-\dot{z}_{d}\left(t\right)-z_{d}\left(t\right)+I_{0}}\left(\ddot{z}_{d}\left(t\right)+\dot{z}_{d}\left(t\right)\right)\nonumber \\
 & +\dfrac{1}{c^{2}}\dfrac{cq+\dot{z}_{d}\left(t\right)+z_{d}\left(t\right)-I_{0}}{cq-\dot{z}_{d}\left(t\right)-z_{d}\left(t\right)+I_{0}}\left(\dot{z}_{d}\left(t\right)+z_{d}\left(t\right)-I_{0}\right)\left(\dot{z}_{d}\left(t\right)+z_{d}\left(t\right)-I_{0}-c\right)\nonumber \\
 & -\dfrac{f}{c}\left(z_{d}\left(t\right)-I_{0}\right).
\end{align}
Since only the output $z\left(t\right)$ can be observed in experiments,
the initial state $\boldsymbol{x}_{0}\left(t\right)=\left(\begin{array}{cc}
x_{0}, & y_{0}\end{array}\right)^{T}$ of the system must be determined from $z\left(t\right)$. Solving
Eq. \eqref{eq:OregonatorOutput} for $x\left(t\right)$ and using
also Eq. \eqref{eq:Eq2230} yields 
\begin{align}
x_{0} & =x\left(t_{0}\right)=\dfrac{1}{c}\left(z\left(t_{0}\right)-I_{0}\right),\label{eq:Eq2232}\\
y_{0} & =y\left(t_{0}\right)=\dfrac{1}{c}\dot{z}\left(t_{0}\right)+\dfrac{1}{c}\left(z\left(t_{0}\right)-I_{0}\right).\label{eq:Eq2233}
\end{align}
Thus, observation of the full initial state requires knowledge of
the output $z\left(t_{0}\right)$ as well as its time derivative $\dot{z}\left(t_{0}\right)$.
A generalization of this fact leads to the notion of observability,
see e.g. \cite{chen1995linear} for the definition of observability
in the context of linear systems. On the other hand, assuming it is
impossible to prepare the system in a desired initial state, the desired
output trajectory $z_{d}\left(t\right)$ has to satisfy specific initial
conditions to comply with the initial system state $\left(\begin{array}{cc}
x_{0}, & y_{0}\end{array}\right)^{T}$. In fact, these conditions are identical in form to Eqs. \eqref{eq:Eq2232}
and \eqref{eq:Eq2233},
\begin{align}
x_{0} & =\dfrac{1}{c}\left(z_{d}\left(t_{0}\right)-I_{0}\right),\label{eq:Eq2234}\\
y_{0} & =\dfrac{1}{c}\dot{z}_{d}\left(t_{0}\right)+\dfrac{1}{c}\left(z_{d}\left(t_{0}\right)-I_{0}\right).\label{eq:Eq2235}
\end{align}
In conclusion, for a successful realization of the desired output
$z_{d}$, not only the initial value of $z_{d}$ but also its time
derivative must be prescribed. This result hints at a connection
between output realizability and observability. A similar connection
between observability and controllability is known as the principle
of duality since the initial work of Kalman \cite{kalman1959general},
see also \cite{chen1995linear}.

\end{example}

\section{\label{sec:2Conclusions}Conclusions}

\subsection{Summary}

A common approach to control, especially in the context of LTI systems,
is concerned with states as the objects to be controlled. Suppose
a controlled system, often called a plant in this context, has a certain
point $\boldsymbol{x}_{1}$ in state space, sometimes called the operating
point, at which the system works efficiently. The control task is
then to bring the system to the operating point $\boldsymbol{x}_{1}$,
and keep it there. This naturally leads to a definition of controllability
as the possibility to achieve a state-to-state transfer from an initial
state $\boldsymbol{x}_{0}$ to the operating point $\boldsymbol{x}_{1}$
within finite time \cite{kalman1959general,chen1995linear}.

In contrast to that, here an approach to control is developed which
centers on the state trajectory $\boldsymbol{x}\left(t\right)$ as
the object of interest. Of course, both approaches to control are
closely related. A single operating point in state space at which
the system is to be kept is nothing more than a degenerate state trajectory.
Equivalently, any state trajectory can be approximated by a succession
of working points.

We distinguish between the controlled state trajectory $\boldsymbol{x}\left(t\right)$
and the desired trajectory $\boldsymbol{x}_{d}\left(t\right)$. The
former is the trajectory which the time-dependent state $\boldsymbol{x}$
traces out in state space under the action of a control signal. The
latter is a fictitious reference trajectory for the state over time.
It is prescribed in analytical or numerical form by the experimenter.
Depending on the choice of the desired trajectory $\boldsymbol{x}_{d}\left(t\right)$,
the controlled state $\boldsymbol{x}\left(t\right)$ may or may not
follow $\boldsymbol{x}_{d}\left(t\right)$.

For affine control systems, the class of exactly realizable desired
trajectories is defined in Section \ref{sec:ExactlyRealizableTrajectories}.
For this subset of desired trajectories, a control signal exists which
enforces the controlled state to follow the desired trajectory exactly,
\begin{align}
\boldsymbol{x}\left(t\right) & =\boldsymbol{x}_{d}\left(t\right),
\end{align}
for all times $t\geq t_{0}$. Exactly realizable desired trajectories
satisfy the constraint equation 
\begin{align}
\boldsymbol{0} & =\boldsymbol{\mathcal{Q}}\left(\boldsymbol{x}_{d}\left(t\right)\right)\left(\boldsymbol{\dot{x}}_{d}\left(t\right)-\boldsymbol{R}\left(\boldsymbol{x}_{d}\left(t\right)\right)\right),\label{eq:Eq2248}
\end{align}
with projector 
\begin{align}
\boldsymbol{\mathcal{Q}}\left(\boldsymbol{x}\right) & =\boldsymbol{1}-\boldsymbol{\mathcal{B}}\left(\boldsymbol{x}\right)\boldsymbol{\mathcal{B}}^{+}\left(\boldsymbol{x}\right)
\end{align}
with rank $n-p$. The vector of control signals $\boldsymbol{u}\left(t\right)$
is expressed in terms of the desired trajectory as
\begin{align}
\boldsymbol{u}\left(t\right) & =\boldsymbol{\mathcal{B}}^{+}\left(\boldsymbol{x}_{d}\left(t\right)\right)\left(\boldsymbol{\dot{x}}_{d}\left(t\right)-\boldsymbol{R}\left(\boldsymbol{x}_{d}\left(t\right)\right)\right).\label{eq:Eq2247}
\end{align}
The matrix $\boldsymbol{\mathcal{B}}^{+}\left(\boldsymbol{x}\right)$
is the Moore-Penrose pseudo inverse of the coupling matrix $\boldsymbol{\mathcal{B}}\left(\boldsymbol{x}\right)$.
Equation \eqref{eq:Eq2247} establishes a one-to-one relationship
between the $p$-dimensional control signal $\boldsymbol{u}\left(t\right)$
and $p$ out of $n$ components of the desired trajectory $\boldsymbol{x}_{d}\left(t\right)$.
The constraint equation \eqref{eq:Eq2248} fixes those $n-p$ components
of the desired trajectory $\boldsymbol{x}_{d}\left(t\right)$ without
a one-to-one relationship to the control signal. The projectors $\boldsymbol{\mathcal{P}}\left(\boldsymbol{x}\right)=\boldsymbol{\mathcal{B}}\left(\boldsymbol{x}\right)\boldsymbol{\mathcal{B}}^{+}\left(\boldsymbol{x}\right)$
and $\boldsymbol{\mathcal{Q}}\left(\boldsymbol{x}\right)=\boldsymbol{1}-\boldsymbol{\mathcal{P}}\left(\boldsymbol{x}\right)$
allow a coordinate-free separation of the state $\boldsymbol{x}$
as well as the controlled state equation in two parts. The part of
the state equation proportional to $\boldsymbol{\mathcal{P}}\left(\boldsymbol{x}\right)$
determines the control signal. This approach allows the elimination
of the control signal Eq. \eqref{eq:Eq2247} from the system. The
remaining part of the state equation, which is proportional to $\boldsymbol{\mathcal{Q}}\left(\boldsymbol{x}\right)$,
is the constraint equation \eqref{eq:Eq2248}. For the control of
exactly realizable trajectories, only the constraint equation must
be solved.

Note that the control signal Eq. \eqref{eq:Eq2247} does not depend
on the state of the system and is therefore an open loop control.
As such, it may suffer from instability. An exactly realizable desired
trajectory might or might not be stable against perturbations of the
initial conditions or external perturbations as e.g. noise.

On the basis of the control signal Eq. \eqref{eq:Eq2247} and constraint
equation \eqref{eq:Eq2248}, a hierarchy of desired trajectories $\boldsymbol{x}_{d}\left(t\right)$
comprising 3 classes is established:

\begin{enumerate}[label=(\Alph*)]

\item desired trajectories $\boldsymbol{x}_{d}\left(t\right)$ which
are solutions to the uncontrolled system,

\item desired trajectories $\boldsymbol{x}_{d}\left(t\right)$ which
are exactly realizable,

\item arbitrary desired trajectories $\boldsymbol{x}_{d}\left(t\right)$.

\end{enumerate}

Desired trajectories of class (A) satisfy the uncontrolled state equation
\begin{align}
\boldsymbol{\dot{x}}_{d}\left(t\right) & =\boldsymbol{R}\left(\boldsymbol{x}_{d}\left(t\right)\right).\label{eq:Eq2249}
\end{align}
This constitutes the most specific class of desired trajectories.
Because of Eq. \eqref{eq:Eq2249}, the constraint equation \eqref{eq:Eq2248}
is trivially satisfied and the control signal as given by Eq. \eqref{eq:Eq2247}
vanishes,
\begin{align}
\boldsymbol{u}\left(t\right) & =\boldsymbol{0}.\label{eq:ControlVanishes}
\end{align}
Equation \eqref{eq:ControlVanishes} implies a non-invasive control
signal, i.e., the control signal vanishes upon achieving the control
target. Because of Eq. \eqref{eq:ControlVanishes}, the open loop
control approach proposed here cannot be employed for desired trajectories
of class (A). Instead, these desired trajectories require feedback
control. Class (A) encompasses several important control tasks, as
e.g. the stabilization of unstable stationary states \cite{sontag2011stability}.
A prominent example extensively studied by the physics community is
the control of chaotic systems by small perturbations \cite{ott1990controlling,shinbrot1993using}.
One of the fundamental aspects of chaos is that many different possible
motions are simultaneously present in the system. In particular, an
infinite number of unstable periodic orbits co-exist with the chaotic
motion. All orbits are solutions to the uncontrolled system dynamics
Eq. \eqref{eq:Eq2249}. Using non-invasive feedback control, a particular
orbit may be stabilized. See also \cite{SchoellSchuster200712,Schimansky-geierFiedlerKurthsScholl200701}
and references therein for more information and examples.

Desired trajectories of class (B) satisfy the constraint equation
\eqref{eq:Eq2248} and yield a non-vanishing control signal $\boldsymbol{u}\left(t\right)\neq\boldsymbol{0}$.
The approach developed in this chapter applies to this class. Several
other techniques developed in mathematical control theory, as e.g.
feedback linearization and differential flatness, also work with this
class of desired trajectories \cite{khalil2002nonlinear,sira2004differentially}.
Class (B) contains the desired trajectories from class (A) as a special
case. For desired trajectories of class (A) and class (B), the solution
of the controlled state trajectory is simply given by $\boldsymbol{x}\left(t\right)=\boldsymbol{x}_{d}\left(t\right)$.

Finally, class (C) is the most general class of desired trajectories
and contains class (A) and (B) as special cases. In general, these
desired trajectories do not satisfy the constraint equation,
\begin{align}
\boldsymbol{0} & \neq\boldsymbol{\mathcal{Q}}\left(\boldsymbol{x}_{d}\left(t\right)\right)\left(\boldsymbol{\dot{x}}_{d}\left(t\right)-\boldsymbol{R}\left(\boldsymbol{x}_{d}\left(t\right)\right)\right),
\end{align}
such that, in general, the approach developed in this chapter cannot
be applied to desired trajectories of class (C). No general expression
for the control signal in terms of the desired trajectory $\boldsymbol{x}_{d}\left(t\right)$
is available. In general, the solution for the controlled state trajectory
$\boldsymbol{x}\left(t\right)$ is not simply given by $\boldsymbol{x}_{d}\left(t\right)$,
$\boldsymbol{x}\left(t\right)\neq\boldsymbol{x}_{d}\left(t\right)$.
Thus, a solution to control problems defined by class (C) does not
only consist in finding an expression for the control signal, but
also involves finding a solution for the controlled state trajectory
$\boldsymbol{x}\left(t\right)$ as well. One possible method to solve
such control problems is optimal control.

The linearizing assumption of Section \ref{sec:LinearizingAssumption}
defines a class of nonlinear control systems which essentially behave
like linear control system. Models satisfying the linearizing assumption
allow exact analytical solutions in closed form even if no analytical
solutions for the uncontrolled system exists, see e.g. the SIR model
in Example \ref{ex:SIRModel4}. The linearizing assumption uncovers
a hidden linear structure underlying nonlinear control systems. Similarly,
feedback linearization defines a huge class of nonlinear control systems
possessing an underlying linear structure. The class of feedback linearizable
systems contains the systems satisfying the linearizing assumption
as a trivial case. However, the linearizing assumption defined here
goes much further than feedback linearization. In fact, while general
nonlinear control systems require a fairly abstract treatment for
the definition of controllability \cite{slotine1991applied,isidori1995nonlinear},
we were able to apply the relatively simple notion of controllability
in terms of a rank condition to systems satisfying the linearizing
assumption, see Section \ref{sec:Controllability}. This is a direct
extension of the properties of linear control systems to a class of
nonlinear control systems. Furthermore, as will be shown in the next
two chapters, the class defined by the linearizing assumption exhibits
a linear structure even in case of optimal control for arbitrary,
not necessarily exactly realizable desired trajectories. This enables
the determination of exact, closed form expressions for optimal trajectory
tracking in Chapter \ref{chap:AnalyticalApproximationsForOptimalTrajectoryTracking}.

The approach to control proposed here shares many similarities to
theories developed in mathematical control theory. We already mentioned
inverse dynamics in the context of mechanical systems in Example \ref{ex:OneDimMechSys3}.
For more information about inverse dynamics, we refer the reader to
the literature about robot control \cite{lewis1993control,deWit1996theory,angeles2002fundamentals}.
In the following, we analyze the similarities and differences of our
approach with differential flatness.

\subsection{Differential flatness}

Similar to the concept of exactly realizable trajectories proposed
in this chapter, differential flatness provides an open loop method
for the control of dynamical systems. We first give a short introduction
to differential flatness to be able to compare the similarities and
differences to our approach. For more information about differential
flatness as well as many examples, we refer the reader to \cite{fliess1995flatness,van1997real,sira2004differentially,levine2009analysis}.
The presentation follows \cite{sira2004differentially}.

Differential flatness relies on the notion of differential functions.
A function $\boldsymbol{\phi}$ is a differential function of $\boldsymbol{x}\left(t\right)$
if it depends on $\boldsymbol{x}\left(t\right)$ and its time derivatives
up to order $\beta$, 
\begin{align}
\boldsymbol{\phi}\left(t\right) & =\boldsymbol{\phi}\left(\boldsymbol{x}\left(t\right),\boldsymbol{\dot{x}}\left(t\right),\dots,\boldsymbol{x}^{\left(\beta\right)}\left(t\right)\right).
\end{align}
The symbol
\begin{align}
\boldsymbol{x}^{\left(\beta\right)}\left(t\right) & =\dfrac{d^{\beta}}{dt^{\beta}}\boldsymbol{x}\left(t\right)
\end{align}
denotes the time derivative of order $\beta$. An affine control system
with $n$-component state vector $\boldsymbol{x}$ and $p$-component
control signal $\boldsymbol{u}$ satisfies 
\begin{align}
\boldsymbol{\dot{x}}\left(t\right) & =\boldsymbol{R}\left(\boldsymbol{x}\left(t\right)\right)+\boldsymbol{\mathcal{B}}\left(\boldsymbol{x}\left(t\right)\right)\boldsymbol{u}\left(t\right).\label{eq:AffineControlSystem}
\end{align}
Applying differentiation with respect to time to Eq. \eqref{eq:AffineControlSystem},
the differential function $\boldsymbol{\tilde{\phi}}\left(\boldsymbol{x},\boldsymbol{\dot{x}}\right)$
can be expressed as a function of $\boldsymbol{x}$ and $\boldsymbol{u}$
\begin{align}
\boldsymbol{\phi}\left(\boldsymbol{x}\left(t\right),\boldsymbol{u}\left(t\right)\right) & =\boldsymbol{\tilde{\phi}}\left(\boldsymbol{x}\left(t\right),\boldsymbol{\dot{x}}\left(t\right)\right)=\boldsymbol{\tilde{\phi}}\left(\boldsymbol{x}\left(t\right),\boldsymbol{R}\left(\boldsymbol{x}\left(t\right)\right)+\boldsymbol{\mathcal{B}}\left(\boldsymbol{x}\left(t\right)\right)\boldsymbol{u}\left(t\right)\right).
\end{align}
Similarly, the differential function $\boldsymbol{\tilde{\phi}}\left(\boldsymbol{x}\left(t\right),\boldsymbol{\dot{x}}\left(t\right),\dots,\boldsymbol{x}^{\left(\beta\right)}\left(t\right)\right)$
can be expressed as a differential function $\boldsymbol{\phi}\left(\boldsymbol{x}\left(t\right),\boldsymbol{u}\left(t\right),\boldsymbol{\dot{u}}\left(t\right),\dots,\boldsymbol{u}^{\left(\beta-1\right)}\left(t\right)\right)$.

The system Eq. \eqref{eq:AffineControlSystem} is called differentially
flat if there exists a $p$-component fictional output $\boldsymbol{z}\left(t\right)=\left(\begin{array}{ccc}
z_{1}\left(t\right), & \dots, & z_{p}\left(t\right)\end{array}\right)^{T}$ such that \cite{sira2004differentially}
\begin{enumerate}
\item the output $\boldsymbol{z}\left(t\right)$ is representable as a differential
function of the state $\boldsymbol{x}\left(t\right)$ and the vector
of control signals $\boldsymbol{u}\left(t\right)$ as
\begin{align}
\boldsymbol{z}\left(t\right) & =\boldsymbol{\phi}\left(\boldsymbol{x}\left(t\right),\boldsymbol{u}\left(t\right),\boldsymbol{\dot{u}}\left(t\right),\boldsymbol{\ddot{u}}\left(t\right),\dots,\boldsymbol{u}^{\left(\beta-1\right)}\left(t\right)\right),
\end{align}

\item the state $\boldsymbol{x}\left(t\right)$ and the vector of control
signals $\boldsymbol{u}\left(t\right)$ are representable as a differential
function of the output $\boldsymbol{z}\left(t\right)$ as (with finite
integer $\alpha$)
\begin{align}
\boldsymbol{x}\left(t\right) & =\boldsymbol{\chi}\left(\boldsymbol{z}\left(t\right),\boldsymbol{\dot{z}}\left(t\right),\boldsymbol{\ddot{z}}\left(t\right),\dots,\boldsymbol{z}^{\left(\alpha\right)}\left(t\right)\right),\\
\boldsymbol{u}\left(t\right) & =\boldsymbol{\psi}\left(\boldsymbol{z}\left(t\right),\boldsymbol{\dot{z}}\left(t\right),\boldsymbol{\ddot{z}}\left(t\right),\dots,\boldsymbol{z}^{\left(\alpha+1\right)}\left(t\right)\right),\label{eq:ControlSignalDifferentialFlatness}
\end{align}

\item the components of the output $\boldsymbol{z}\left(t\right)$ are differentially
independent, i.e., they satisfy no differential equation of the form
\begin{align}
\boldsymbol{\Omega}\left(\boldsymbol{z}\left(t\right),\boldsymbol{\dot{z}}\left(t\right),\boldsymbol{\ddot{z}}\left(t\right),\dots,\boldsymbol{z}^{\left(\beta\right)}\left(t\right)\right) & =\boldsymbol{0}.
\end{align}

\end{enumerate}
For a differentially flat system, the full solution for the state
over time $\boldsymbol{x}\left(t\right)$ as well as the control signal
$\boldsymbol{u}\left(t\right)$ can be expressed in terms of the output
over time $\boldsymbol{z}\left(t\right)$. Mathematically, the relation
between control signal $\boldsymbol{u}\left(t\right)$ and output
$\boldsymbol{z}\left(t\right)$ is a differential function. This has
the great advantage that the determination of the control signal can
be done in real time at time $t$ by computing only a finite number
of time derivatives of $\boldsymbol{z}\left(t\right)$. This would
not be possible if $\boldsymbol{u}\left(t\right)$ also involves time
integrals of the output $\boldsymbol{z}$ because these would require
summation over all previous times as well. Another advantage is that
no differential equations need to be solved to obtain the control
signal and state trajectory. Usually, all expressions are generated
by simply differentiating the controlled state equations with respect
to time. The output $\boldsymbol{z}\left(t\right)$ has the same number
$p$ of components as the number of independent input signals $\boldsymbol{u}\left(t\right)$
available to control the system. If the control signal determined
by $\boldsymbol{\psi}$, Eq. \eqref{eq:ControlSignalDifferentialFlatness},
is applied to the system, then the system's output is $\boldsymbol{z}\left(t\right)$.
Differentially flat systems are not necessarily affine in control
but can be nonlinear in the control as well. However, only certain
systems are differentially flat, and it is not known under which conditions
a controlled dynamical system is differentially flat if the number
of independent control signals is larger than one, $p>1$.

Similar to differential flatness, this chapter proposes an open loop
control method. A solution for control signals exactly realizing a
desired output $\boldsymbol{z}_{d}\left(t\right)$ with $p$ components
is determined. In the discussion of output trajectory realizability
in Section \ref{sec:OutputRealizability}, the control signal is expressed
solely in terms of the desired output $\boldsymbol{z}_{d}\left(t\right)$
and the initial conditions for the state. This implies that the controlled
state trajectory, given as the solution to the controlled state equation,
can also be expressed in terms of the desired output and the initial
conditions for the state. These facts fully agree with the concept
of differential flatness. In contrast to the approach here, the literature
about differential flatness does usually not distinguish explicitly
between desired trajectory $\boldsymbol{x}_{d}\left(t\right)$ and
controlled state trajectory $\boldsymbol{x}\left(t\right)$, but implicitly
assumes this identity from the very beginning.

The most striking difference between the approach here and differential
flatness is the restriction to differential functions. In general,
our approach yields a control signal in terms of a functional of the
desired output, 
\begin{align}
\boldsymbol{u}\left(t\right) & =\boldsymbol{u}\left[\boldsymbol{z}_{d}\left(t\right)\right].
\end{align}
Note that a functional is a more general expression than a differential
function. Using the Dirac delta function $\delta\left(t\right)$,
any time derivative of order $\beta$ can be expressed as a functional,
\begin{align}
\boldsymbol{z}_{d}^{\left(\beta\right)}\left(t\right) & =\intop_{-\infty}^{\infty}d\tau\delta\left(\tau-t\right)\boldsymbol{z}_{d}^{\left(\beta\right)}\left(\tau\right)=-\intop_{-\infty}^{\infty}d\tau\delta'\left(\tau-t\right)\boldsymbol{z}_{d}^{\left(\beta-1\right)}\left(\tau\right)\nonumber \\
 & \vdots\nonumber \\
 & =\left(-1\right)^{\beta}\intop_{-\infty}^{\infty}d\tau\delta^{\left(\beta\right)}\left(\tau-t\right)\boldsymbol{z}_{d}\left(\tau\right)
\end{align}
Therefore, any differential function of $\boldsymbol{z}\left(t\right)$
can be expressed in terms of a function of functionals of $\boldsymbol{z}\left(t\right)$,
while the reverse is not true. The restriction to differential functions
might also explain why only certain systems are differentially flat.
In contrast, the approach proposed here can be applied to any affine
control system. As an advantage, differential flatness yields expressions
for state and control which are computationally more efficient because
they do not require the solution of differential equations or integrals,
which is not the case here.

\subsection{Outlook}

The framework of exactly realizable trajectories is interpreted as
an open loop control method. However, it may be possible to extend
this approach to feedback control. As discussed in Section \ref{sec:OutputRealizability},
a control $\boldsymbol{u}\left(t\right)$ realizing a $p$-component
desired output $\boldsymbol{z}_{d}\left(t\right)=\boldsymbol{h}\left(\boldsymbol{x}_{d}\left(t\right)\right)$
is expressed entirely in terms of the desired output. The dependence
of $\boldsymbol{u}\left(t\right)$ on $\boldsymbol{z}_{d}\left(t\right)$
is typically in form of a functional,
\begin{align}
\boldsymbol{u}\left(t\right) & =\boldsymbol{u}\left[\boldsymbol{z}_{d}\left(t\right)\right].\label{eq:Eq2261}
\end{align}
A generalization to feedback control yields a control signal which
does not only depend on the desired output $\boldsymbol{z}_{d}\left(t\right)$
but also on the monitored state $\boldsymbol{x}\left(t\right)$ of
the controlled system, 
\begin{align}
\boldsymbol{u}\left(t\right) & =\boldsymbol{u}\left[\boldsymbol{z}_{d}\left(t\right),\boldsymbol{x}\left(t\right)\right].
\end{align}
In general, the control signal is allowed to depend on the history
of $\boldsymbol{x}\left(t\right)$ such that the dependence of $\boldsymbol{u}\left(t\right)$
on $\boldsymbol{x}\left(t\right)$ is also in form of a functional.
Such a generalization of the approach to control proposed here certainly
changes the stability properties of the controlled system and may
result in an improved stability of the controlled trajectory.

A fundamental problem affecting not only exactly realizable trajectories
but also feedback linearization, differential flatness, and optimal
control, is the requirement of exactly knowing the system dynamics.
This must be contrasted with the fact that the majority of physical
models are idealizations. Unknown external influences in control systems
can be modeled as noise or structural perturbations, which might both
depend on the system state itself. To ensure a successful control
in experiments, the proposed control methods must not only be stable
against perturbations of the initial conditions, but must be sufficiently
stable against structural perturbations as well. Stability against
structural perturbations is also known as robustness in the context
of control theory \cite{freeman2008robust}. Before applying the control
method developed in this chapter to real world problems, a thorough
investigation of the stability of the control problem at hand must
be conducted. In case of instability, countermeasures as e.g. additional
stabilizing feedback control must be applied \cite{khalil2002nonlinear}.

Section \ref{sec:LinearizingAssumption} introduces the linearizing
assumption. On the one hand, this assumption is restrictive, but on
the other hand it has far reaching consequences and results in significant
simplifications for nonlinear affine control systems. A possible generalization
of the linearizing assumption might be as follows. First, relax condition
Eq. \eqref{eq:LinearizingAssumption1} and allow a state dependent
projector $\boldsymbol{\mathcal{Q}}\left(\boldsymbol{x}\right)$ which,
however, does only depend on the state components $\boldsymbol{\mathcal{P}}\boldsymbol{x}$,
\begin{align}
\boldsymbol{\mathcal{Q}}\left(\boldsymbol{x}\right) & =\boldsymbol{\mathcal{Q}}\left(\boldsymbol{\mathcal{P}}\boldsymbol{x}+\boldsymbol{\mathcal{Q}}\boldsymbol{x}\right)=\boldsymbol{\mathcal{Q}}\left(\boldsymbol{\mathcal{P}}\boldsymbol{x}\right).
\end{align}
Second, also relax condition Eq. \eqref{eq:LinearizingAssumption2}
and assume a nonlinearity $\boldsymbol{R}\left(\boldsymbol{x}\right)$
with the following structure, 
\begin{align}
\boldsymbol{\mathcal{Q}}\left(\boldsymbol{x}\right)\boldsymbol{R}\left(\boldsymbol{x}\right) & =\boldsymbol{\mathcal{Q}}\left(\boldsymbol{\mathcal{P}}\boldsymbol{x}\right)\boldsymbol{\mathcal{A}}\left(\boldsymbol{\mathcal{P}}\boldsymbol{x}\right)\boldsymbol{\mathcal{Q}}\left(\boldsymbol{\mathcal{P}}\boldsymbol{x}\right)\boldsymbol{x}+\boldsymbol{\mathcal{Q}}\left(\boldsymbol{\mathcal{P}}\boldsymbol{x}\right)\boldsymbol{b}\left(\boldsymbol{\mathcal{P}}\boldsymbol{x}\right).
\end{align}
The matrix $\boldsymbol{\mathcal{A}}\left(\boldsymbol{\mathcal{P}}\boldsymbol{x}\right)$,
the projector $\boldsymbol{\mathcal{Q}}\left(\boldsymbol{\mathcal{P}}\boldsymbol{x}\right)$
and the inhomogeneity $\boldsymbol{b}\left(\boldsymbol{\mathcal{P}}\boldsymbol{x}\right)$
may all depend on the state components $\boldsymbol{\mathcal{P}}\boldsymbol{x}$.

Together with
\begin{align}
\dfrac{d}{dt}\left(\boldsymbol{\mathcal{Q}}\left(\boldsymbol{\mathcal{P}}\boldsymbol{x}_{d}\right)\boldsymbol{x}_{d}\right) & =\boldsymbol{\mathcal{\dot{Q}}}\left(\boldsymbol{\mathcal{P}}\boldsymbol{x}_{d}\right)\boldsymbol{\mathcal{P}}\left(\boldsymbol{\mathcal{P}}\boldsymbol{x}_{d}\right)\boldsymbol{x}_{d}\nonumber \\
 & +\boldsymbol{\mathcal{\dot{Q}}}\left(\boldsymbol{\mathcal{P}}\boldsymbol{x}_{d}\right)\boldsymbol{\mathcal{Q}}\left(\boldsymbol{\mathcal{P}}\boldsymbol{x}_{d}\right)\boldsymbol{x}_{d}+\boldsymbol{\mathcal{Q}}\left(\boldsymbol{\mathcal{P}}\boldsymbol{x}_{d}\right)\boldsymbol{\dot{x}}_{d},
\end{align}
the constraint equation becomes
\begin{align}
\dfrac{d}{dt}\left(\boldsymbol{\mathcal{Q}}\boldsymbol{x}_{d}\left(t\right)\right) & =\left(\boldsymbol{\mathcal{\dot{Q}}}+\boldsymbol{\mathcal{Q}}\boldsymbol{\mathcal{A}}\right)\boldsymbol{\mathcal{Q}}\boldsymbol{x}_{d}\left(t\right)+\boldsymbol{\mathcal{\dot{Q}}}\boldsymbol{\mathcal{P}}\boldsymbol{x}_{d}\left(t\right)+\boldsymbol{\mathcal{Q}}\boldsymbol{b}.\label{eq:Eq2276}
\end{align}
The arguments are suppressed and it is understood that $\boldsymbol{\mathcal{\dot{Q}}}$,
$\boldsymbol{\mathcal{Q}}$, $\boldsymbol{\mathcal{P}}$, $\boldsymbol{\mathcal{A}}$,
and $\boldsymbol{b}$ may depend on the part $\boldsymbol{\mathcal{P}}\boldsymbol{x}_{d}\left(t\right)$.
Equation \eqref{eq:Eq2276} is a linear equation for $\boldsymbol{\mathcal{Q}}\boldsymbol{x}_{d}\left(t\right)$
and can thus be solved with the help of its state transition matrix,
see Appendix \ref{sec:GeneralSolutionForForcedLinarDynamicalSystem}.
However, the matrix $\boldsymbol{\mathcal{A}}=\boldsymbol{\mathcal{A}}\left(\boldsymbol{\mathcal{P}}\boldsymbol{x}_{d}\left(t\right)\right)$
exhibits an explicit time dependence through its dependence on $\boldsymbol{\mathcal{P}}\boldsymbol{x}_{d}\left(t\right)$.
This necessitates modifications for the notion of controllability
from Section \ref{sec:Controllability}, see also \cite{chen1995linear}.

A central assumption of the formalism presented in this chapter is
that the $n\times p$ coupling matrix $\boldsymbol{\mathcal{B}}\left(\boldsymbol{x}\right)$
has full rank $p$ for all values of $\boldsymbol{x}$. This assumption
leads to a Moore-Penrose pseudo inverse $\boldsymbol{\mathcal{B}}^{+}\left(\boldsymbol{x}\right)$
of $\boldsymbol{\mathcal{B}}\left(\boldsymbol{x}\right)$ given by
\begin{align}
\boldsymbol{\mathcal{B}}^{+}\left(\boldsymbol{x}\right) & =\left(\boldsymbol{\mathcal{B}}^{T}\left(\boldsymbol{x}\right)\boldsymbol{\mathcal{B}}\left(\boldsymbol{x}\right)\right)^{-1}\boldsymbol{\mathcal{B}}^{T}\left(\boldsymbol{x}\right).
\end{align}
If $\boldsymbol{\mathcal{B}}\left(\boldsymbol{x}\right)$ does not
have full rank for some or all values of $\boldsymbol{x}$, the inverse
of $\boldsymbol{\mathcal{B}}^{T}\left(\boldsymbol{x}\right)\boldsymbol{\mathcal{B}}\left(\boldsymbol{x}\right)$
does not exist. However, a unique Moore-Penrose pseudo inverse $\boldsymbol{\mathcal{B}}^{+}\left(\boldsymbol{x}\right)$
does exist for any matrix $\boldsymbol{\mathcal{B}}\left(\boldsymbol{x}\right)$,
regardless of its rank. No closed form expressions exist for the general
case, but $\boldsymbol{\mathcal{B}}^{+}\left(\boldsymbol{x}\right)$
can nevertheless be computed numerically by singular value decomposition,
for example. Because $\boldsymbol{\mathcal{B}}^{+}\left(\boldsymbol{x}\right)$
exists in any case, the $n\times n$ projector defined by
\begin{align}
\boldsymbol{\mathcal{P}}\left(\boldsymbol{x}\right) & =\boldsymbol{\mathcal{B}}\left(\boldsymbol{x}\right)\boldsymbol{\mathcal{B}}^{+}\left(\boldsymbol{x}\right)
\end{align}
exists as well. Thus, using the general Moore-Penrose pseudo inverse
$\boldsymbol{\mathcal{B}}^{+}\left(\boldsymbol{x}\right)$, the formalism
developed in this chapter can be extended to cases with $\boldsymbol{\mathcal{B}}\left(\boldsymbol{x}\right)$
not having full rank for some values of $\boldsymbol{x}$.

A mathematically more rigorous treatment of the notion of exactly
realizable trajectories is desirable. An important question is the
following. Under which conditions does the constraint equation 
\begin{align}
\boldsymbol{0} & =\boldsymbol{\mathcal{Q}}\left(\boldsymbol{x}_{d}\left(t\right)\right)\left(\boldsymbol{\dot{x}}_{d}\left(t\right)-\boldsymbol{R}\left(\boldsymbol{x}_{d}\left(t\right)\right)\right)\label{eq:Eq2277}
\end{align}
have a unique solution for $\boldsymbol{\mathcal{Q}}\left(\boldsymbol{x}_{d}\left(t\right)\right)\boldsymbol{x}_{d}\left(t\right)$?
Note that Eq. \eqref{eq:Eq2277} is a non-autonomous nonlinear system
of differential equations for $\boldsymbol{\mathcal{Q}}\left(\boldsymbol{x}_{d}\left(t\right)\right)\boldsymbol{x}_{d}\left(t\right)$
with the explicit time dependence caused by the part $\boldsymbol{\mathcal{P}}\left(\boldsymbol{x}_{d}\left(t\right)\right)\boldsymbol{x}_{d}\left(t\right)$.
Therefore, a related question is for conditions on the part $\boldsymbol{\mathcal{P}}\left(\boldsymbol{x}_{d}\left(t\right)\right)\boldsymbol{x}_{d}\left(t\right)$
prescribed by the experimenter. For example, is $\boldsymbol{\mathcal{P}}\left(\boldsymbol{x}_{d}\left(t\right)\right)\boldsymbol{x}_{d}\left(t\right)$
required to be a continuously differentiable function or is it allowed
to have jumps? Although some general answers might be possible, such
questions are simpler to answer for specific control systems.

%% file: chapter-3.tex
\lhead[\chaptername~\thechapter\leftmark]{}

\rhead[]{\rightmark}

\lfoot[\thepage]{}

\cfoot{}

\rfoot[]{\thepage}

\chapter{\label{chap:OptimalControl}Optimal control}

This chapter introduces the standard approach to optimal control theory
in form of the necessary optimality conditions in Section \ref{sec:NecessaryOptimalityConditions}.
Additional necessary optimality conditions for singular optimal control
problems, the so-called Kelly or generalized Legendre-Clebsch conditions,
are presented in Section \ref{sec:SingularOptimalControls}. Section
\ref{sec:NumericalSolutionOfOptimalControlProblems} gives a brief
discussion of the difficulties involved in finding a numerical solution
to an optimal control problem. The conditions under which the control
of exactly realizable trajectories is optimal are clarified in Section
\ref{sec:ExactlyRealizableTrajectoriesOptimalControl}. The last Section
\ref{sec:AnExactlySolvableExample} presents a simple linear optimal
control problem for which an exact but cumbersome analytical solution
can be derived. Assuming a small regularization parameter $0<\epsilon\ll1$,
the exact solution is approximated and a simpler expression is obtained.
Additionally, the impact of different terminal conditions on the solution
is investigated.

\section{\label{sec:NecessaryOptimalityConditions}The necessary optimality
conditions}

The foundations of optimal control theory were laid in the 1950's
and early 1960's. The Russian school of Lev Pontryagin and his students
developed the minimum principle \cite{pontryagin1987mathematical};
also called the maximum principle in the Russian literature. This
principle is based on the calculus of variations and contains the
Euler-Lagrange equations of uncontrolled dynamical systems as a special
case. An American school, led by Richard Bellmann, developed the Dynamical
Programming approach \cite{bellman1957dynamic} based on partial differential
equations (PDEs). Both approaches treat essentially the same problem
and yield equivalent results. Which approach is preferred is, to some
extent, a matter of taste. Here, optimal control is discussed in the
framework of Pontryagin's minimum optimal control. An excellent introduction
for both approaches to optimal control is the book by Bryson and Ho
\cite{bryson1969applied}. A more elementary and technical approach
to the calculus of variations and optimal control is provided by the
readable introduction \cite{liberzon2011calculus}. See also \cite{hull2003optimal}
for applications of optimal control and the mathematically rigorous
treatment \cite{vinter2010optimal}.

\subsection{Statement of the problem}

Optimal control is concerned with minimizing a target functional $\mathcal{J}$,
also called the performance index,
\begin{align}
\mathcal{J}\left[\boldsymbol{x}\left(t\right),\boldsymbol{u}\left(t\right)\right] & =\intop_{t_{0}}^{t_{1}}dtL\left(\boldsymbol{x}\left(t\right),t\right)+M\left(\boldsymbol{x}\left(t_{1}\right),t_{1}\right)+\frac{\epsilon^{2}}{2}\intop_{t_{0}}^{t_{1}}dt\left(\boldsymbol{u}\left(t\right)\right)^{2}.\label{eq:OptimalControlFunctional}
\end{align}
Here, $L\left(\boldsymbol{x},t\right)$ is the cost function and $M\left(\boldsymbol{x},t\right)$
is the terminal cost. The target functional Eq. \eqref{eq:OptimalControlFunctional}
is to be minimized subject to the constraint that $\boldsymbol{x}\left(t\right)$
is given as the solution to the controlled dynamical system 
\begin{align}
\boldsymbol{\dot{x}}\left(t\right) & =\boldsymbol{R}\left(\boldsymbol{x}\left(t\right)\right)+\boldsymbol{\mathcal{B}}\left(\boldsymbol{x}\left(t\right)\right)\boldsymbol{u}\left(t\right),
\end{align}
with initial condition
\begin{align}
\boldsymbol{x}\left(t_{0}\right) & =\boldsymbol{x}_{0}.
\end{align}
The parameter $\epsilon$ is called the regularization parameter.
Furthermore, the state $\boldsymbol{x}$ satisfies the $q\leq n$
end point conditions
\begin{align}
\boldsymbol{\psi}\left(\boldsymbol{x}\left(t_{1}\right)\right) & =\left(\psi_{1}\left(\boldsymbol{x}\left(t_{1}\right)\right),\dots,\psi_{q}\left(\boldsymbol{x}\left(t_{1}\right)\right)\right)^{T}=\boldsymbol{0}.
\end{align}

\subsection{Derivation of the necessary optimality conditions}

Following a standard procedure \cite{bryson1969applied}, the constrained
optimization problem is converted to an unconstrained optimization
problem. Similar to minimizing an ordinary function under constraints,
this is done by introducing Lagrange multipliers. However, in optimal
control, the constraint is a differential equation defined on a certain
time interval. Consequently, the Lagrange multipliers are functions
of time and denoted by $\boldsymbol{\lambda}\left(t\right)=\left(\lambda_{1}\left(t\right),\dots,\lambda_{n}\left(t\right)\right)^{T}\in\mathbb{R}^{n}$.
The vector $\boldsymbol{\lambda}\left(t\right)$ is called the \textit{co-state}
or adjoint state. To accommodate the end point condition $\boldsymbol{\psi}$,
additional constant Lagrange multipliers $\boldsymbol{\nu}=\left(\nu_{1},\dots,\nu_{q}\right)^{T}\in\mathbb{R}^{q}$
are introduced. The constrained optimization problem is reduced to
the minimization of the unconstrained functional 
\begin{align}
\bar{\mathcal{J}}\left[\boldsymbol{x}\left(t\right),\boldsymbol{u}\left(t\right),\boldsymbol{\lambda}\left(t\right),\boldsymbol{\nu}\right] & =\mathcal{J}\left[\boldsymbol{x}\left(t\right),\boldsymbol{u}\left(t\right)\right]+\boldsymbol{\nu}^{T}\boldsymbol{\psi}\left(\boldsymbol{x}\left(t_{1}\right)\right)\nonumber \\
 & +\intop_{t_{0}}^{t_{1}}dt\boldsymbol{\lambda}^{T}\left(t\right)\left(\boldsymbol{R}\left(\boldsymbol{x}\left(t\right)\right)+\boldsymbol{\mathcal{B}}\left(\boldsymbol{x}\left(t\right)\right)\boldsymbol{u}\left(t\right)-\boldsymbol{\dot{x}}\left(t\right)\right).\label{eq:NewOptimalControlFunctional}
\end{align}
Introducing the control Hamiltonian 
\begin{align}
H\left(\boldsymbol{x}\left(t\right),\boldsymbol{u}\left(t\right),t\right) & =L\left(\boldsymbol{x}\left(t\right),t\right)+\frac{\epsilon^{2}}{2}\left(\boldsymbol{u}\left(t\right)\right)^{2}\nonumber \\
 & +\boldsymbol{\lambda}^{T}\left(t\right)\left(\boldsymbol{R}\left(\boldsymbol{x}\left(t\right)\right)+\boldsymbol{\mathcal{B}}\left(\boldsymbol{x}\left(t\right)\right)\boldsymbol{u}\left(t\right)\right),
\end{align}
and applying partial integration for the term involving $\boldsymbol{\lambda}^{T}\left(t\right)\boldsymbol{\dot{x}}\left(t\right)$
in Eq. \eqref{eq:NewOptimalControlFunctional} yields
\begin{align}
\bar{\mathcal{J}}\left[\boldsymbol{x}\left(t\right),\boldsymbol{u}\left(t\right),\boldsymbol{\lambda}\left(t\right),\boldsymbol{\nu}\right] & =\intop_{t_{0}}^{t_{1}}dtH\left(\boldsymbol{x}\left(t\right),\boldsymbol{u}\left(t\right),t\right)+M\left(\boldsymbol{x}\left(t_{1}\right),t_{1}\right)+\boldsymbol{\nu}^{T}\boldsymbol{\psi}\left(\boldsymbol{x}\left(t_{1}\right)\right)\nonumber \\
 & -\boldsymbol{\lambda}^{T}\left(t_{1}\right)\boldsymbol{x}\left(t_{1}\right)+\boldsymbol{\lambda}^{T}\left(t_{0}\right)\boldsymbol{x}\left(t_{0}\right)+\intop_{t_{0}}^{t_{1}}dt\boldsymbol{\dot{\lambda}}^{T}\left(t\right)\boldsymbol{x}\left(t\right).
\end{align}
The functional $\bar{\mathcal{J}}\left[\boldsymbol{x}\left(t\right),\boldsymbol{u}\left(t\right),\boldsymbol{\lambda}\left(t\right),\boldsymbol{\nu}\right]$
must be minimized with respect to $\boldsymbol{\lambda}\left(t\right)$,
$\boldsymbol{\lambda}\left(t_{1}\right)$, $\boldsymbol{u}\left(t\right),\,\boldsymbol{\nu},$
and $\boldsymbol{x}\left(t\right)$. The initial condition $\boldsymbol{x}_{0}$
is prescribed and is therefore kept fixed. For $\bar{\mathcal{J}}$
to be extremal, it has to satisfy the variational equations
\begin{align}
\frac{\delta\bar{\mathcal{J}}}{\delta\boldsymbol{x}\left(t\right)} & =\boldsymbol{0}, & \frac{\delta\bar{\mathcal{J}}}{\delta\boldsymbol{u}\left(t\right)} & =\boldsymbol{0}, & \frac{\delta\bar{\mathcal{J}}}{\delta\boldsymbol{\lambda}\left(t\right)} & =\boldsymbol{0}, & \frac{\delta\bar{\mathcal{J}}}{\delta\boldsymbol{\nu}} & =\boldsymbol{0}, & \frac{\delta\bar{\mathcal{J}}}{\delta\boldsymbol{\lambda}\left(t_{1}\right)} & =\boldsymbol{0}.
\end{align}
The variation $\frac{\delta\bar{\mathcal{J}}}{\delta\boldsymbol{u}\left(t\right)}=\boldsymbol{0}$
with respect to the vector of control signals $\boldsymbol{u}$ leads
to $p$ algebraic equations for the control vector $\boldsymbol{u}\left(t\right)$,
\begin{align}
\epsilon^{2}\boldsymbol{u}^{T}\left(t\right)+\boldsymbol{\lambda}^{T}\left(t\right)\boldsymbol{\mathcal{B}}\left(\boldsymbol{x}\left(t\right)\right)= & \boldsymbol{0}.
\end{align}
Transposing yields
\begin{align}
\epsilon^{2}\boldsymbol{u}\left(t\right)+\boldsymbol{\mathcal{B}}^{T}\left(\boldsymbol{x}\left(t\right)\right)\boldsymbol{\lambda}\left(t\right)= & \boldsymbol{0}.\label{eq:ControlSolutionAdjointEquation}
\end{align}
The variation $\frac{\delta\bar{\mathcal{J}}}{\delta\boldsymbol{\lambda}^{T}\left(t\right)}=\boldsymbol{0}$
with respect to the co-state $\boldsymbol{\lambda}\left(t\right)$
leads to the controlled state equation 
\begin{align}
\boldsymbol{\dot{x}}\left(t\right) & =\boldsymbol{R}\left(\boldsymbol{x}\left(t\right)\right)+\boldsymbol{\mathcal{B}}\left(\boldsymbol{x}\left(t\right)\right)\boldsymbol{u}\left(t\right),\label{eq:StateEquations}\\
\boldsymbol{x}\left(t_{0}\right) & =\boldsymbol{x}_{0}.\label{eq:StateEquationInitCond}
\end{align}
The variation $\frac{\delta\bar{\mathcal{J}}}{\delta\boldsymbol{x}\left(t\right)}=\boldsymbol{0}$
with respect to the state $\boldsymbol{x}\left(t\right)$ leads to
the so-called \textit{adjoint} or \textit{co-state equation}
\begin{align}
-\boldsymbol{\dot{\lambda}}^{T}\left(t\right)= & \boldsymbol{\lambda}^{T}\left(t\right)\nabla\boldsymbol{R}\left(\boldsymbol{x}\left(t\right)\right)+\boldsymbol{\lambda}^{T}\left(t\right)\nabla\boldsymbol{\mathcal{B}}\left(\boldsymbol{x}\left(t\right)\right)\boldsymbol{u}\left(t\right)+\nabla L\left(\boldsymbol{x}\left(t\right),t\right).
\end{align}
Here, the $n\times n$ Jacobi matrix of the $n$-dimensional nonlinear
function $\boldsymbol{R}$ is defined as
\begin{align}
\left(\nabla\boldsymbol{R}\left(\boldsymbol{x}\right)\right)_{ij} & =\frac{\partial}{\partial x_{j}}R_{i}\left(\boldsymbol{x}\right).
\end{align}
The $n\times p\times n$ Jacobi matrix of the $n\times p$ matrix
$\boldsymbol{\mathcal{B}}\left(\boldsymbol{x}\right)$ is given by
\begin{align}
\left(\nabla\boldsymbol{\mathcal{B}}\left(\boldsymbol{x}\right)\right)_{ijk} & =\frac{\partial}{\partial x_{k}}\mathcal{B}_{ij}\left(\boldsymbol{x}\right).
\end{align}
Written component-wise, the inner product of $\boldsymbol{\lambda}$
and $\boldsymbol{\mathcal{B}}\left(\boldsymbol{x}\right)\boldsymbol{u}$
is 
\begin{align}
\boldsymbol{\lambda}^{T}\boldsymbol{\mathcal{B}}\left(\boldsymbol{x}\right)\boldsymbol{u} & =\sum_{i=1}^{n}\sum_{j=1}^{p}\lambda_{i}\mathcal{B}_{ij}\left(\boldsymbol{x}\right)u_{j}.
\end{align}
Consequently, the expression $\boldsymbol{\lambda}^{T}\nabla\boldsymbol{\mathcal{B}}\left(\boldsymbol{x}\right)\boldsymbol{u}$
is an $n$-component row vector defined as
\begin{align}
\boldsymbol{\lambda}^{T}\nabla\boldsymbol{\mathcal{B}}\left(\boldsymbol{x}\right)\boldsymbol{u} & =\sum_{i=1}^{n}\sum_{j=1}^{p}\left(\begin{array}{ccc}
\lambda_{i}\dfrac{\partial}{\partial x_{1}}\mathcal{B}_{ij}\left(\boldsymbol{x}\right)u_{j}, & \dots, & \lambda_{i}\dfrac{\partial}{\partial x_{n}}\mathcal{B}_{ij}\left(\boldsymbol{x}\right)u_{j}\end{array}\right).\label{eq:lambdanablaBu}
\end{align}
The variation $\frac{\delta\bar{\mathcal{J}}}{\delta\boldsymbol{x}\left(t_{1}\right)}=\boldsymbol{0}$
with respect to the terminal state $\boldsymbol{x}\left(t_{1}\right)$
leads to the corresponding boundary condition for the co-state,
\begin{align}
\boldsymbol{\lambda}^{T}\left(t_{1}\right)= & \nabla M\left(\boldsymbol{x}\left(t_{1}\right),t_{1}\right)+\boldsymbol{\nu}^{T}\nabla\boldsymbol{\psi}\left(\boldsymbol{x}\left(t_{1}\right)\right).
\end{align}
Transposing finally gives 
\begin{align}
-\boldsymbol{\dot{\lambda}}\left(t\right) & =\left(\nabla\boldsymbol{R}^{T}\left(\boldsymbol{x}\left(t\right)\right)+\boldsymbol{u}^{T}\left(t\right)\nabla\boldsymbol{\mathcal{B}}^{T}\left(\boldsymbol{x}\left(t\right)\right)\right)\boldsymbol{\lambda}\left(t\right)+\left(\nabla L\left(\boldsymbol{x}\left(t\right),t\right)\right)^{T},\label{eq:AdjointEquation2}\\
\boldsymbol{\lambda}\left(t_{1}\right) & =\nabla M^{T}\left(\boldsymbol{x}\left(t_{1}\right),t_{1}\right)+\nabla\boldsymbol{\psi}^{T}\left(\boldsymbol{x}\left(t_{1}\right)\right)\boldsymbol{\nu}.\label{eq:AdjointEquation2TermCond}
\end{align}
Note that the co-state $\boldsymbol{\lambda}\left(t\right)$ satisfies
a terminal condition rather than an initial condition. The variation
$\frac{\delta\bar{\mathcal{J}}}{\delta\boldsymbol{\nu}^{T}}=\boldsymbol{0}$
with respect to the Lagrange multipliers $\boldsymbol{\nu}^{T}$ finally
yields
\begin{align}
\boldsymbol{\psi}\left(\boldsymbol{x}\left(t_{1}\right)\right) & =\boldsymbol{0}.
\end{align}
The state equation \eqref{eq:StateEquations} and the adjoint equation
\eqref{eq:AdjointEquation2} for the state $\boldsymbol{x}\left(t\right)$
and co-state $\boldsymbol{\lambda}\left(t\right)$, respectively,
as well as the algebraic expression Eq. \eqref{eq:ControlSolutionAdjointEquation}
for the vector of control signals $\boldsymbol{u}\left(t\right)$
constitute the \textit{necessary optimality conditions.}

\subsection{Optimal trajectory tracking}

The optimal control problem considered in this thesis is to steer
the system state $\boldsymbol{x}\left(t\right)$ as closely as possible
along a desired reference trajectory $\boldsymbol{x}_{d}\left(t\right)$.
A common choice of the cost functions $L$ and $M$ is the quadratic
difference between the actual trajectory $\boldsymbol{x}\left(t\right)$
and desired trajectory $\boldsymbol{x}_{d}\left(t\right)$,
\begin{align}
L\left(\boldsymbol{x}\left(t\right),t\right) & =\frac{1}{2}\left(\boldsymbol{x}\left(t\right)-\boldsymbol{x}_{d}\left(t\right)\right)^{T}\boldsymbol{\mathcal{S}}\left(\boldsymbol{x}\left(t\right)-\boldsymbol{x}_{d}\left(t\right)\right),\label{eq:Eq324}\\
M\left(\boldsymbol{x}\left(t_{1}\right),t_{1}\right) & =\frac{1}{2}\left(\boldsymbol{x}\left(t_{1}\right)-\boldsymbol{x}_{1}\right)^{T}\boldsymbol{\mathcal{S}}_{1}\left(\boldsymbol{x}\left(t_{1}\right)-\boldsymbol{x}_{1}\right).\label{eq:Eq331}
\end{align}
The matrices $\boldsymbol{\mathcal{S}}$ and $\boldsymbol{\mathcal{S}}_{1}$
are usually assumed to be symmetric and positive definite matrices
of weighting coefficients. The expression $\intop_{t_{0}}^{t_{1}}dt\boldsymbol{x}^{T}\left(t\right)\boldsymbol{\mathcal{S}}\boldsymbol{x}\left(t\right)$
is called the weighted $L^{2}$ norm of $\boldsymbol{x}\left(t\right)$.
The choice of Eqs. \eqref{eq:Eq324} and \eqref{eq:Eq331} for the
cost functions defines the problem of \textit{optimal trajectory tracking}.

Note that there are different possibilities for the terminal condition
at time $t=t_{1}$. The notions end point and terminal point are used
synonymously. A squared difference term $M\left(\boldsymbol{x}\left(t_{1}\right),t_{1}\right)$
as in Eq. \eqref{eq:Eq331} penalizes a large deviation of the terminal
state $\boldsymbol{x}\left(t_{1}\right)$ from the state space point
$\boldsymbol{x}_{1}$. An end point condition of the form $\boldsymbol{\psi}\left(\boldsymbol{x}\left(t_{1}\right)\right)=\left(\boldsymbol{x}\left(t_{1}\right)-\boldsymbol{x}_{1}\right)^{T}$
insists on $\boldsymbol{x}\left(t_{1}\right)=\boldsymbol{x}_{1}$
at the terminal time. Here, the latter case is called an \textit{exact}
or \textit{sharp} terminal condition. Both possibilities can appear
in the same problem. A sharp terminal condition is much more restrictive
than a squared difference term. If $\boldsymbol{x}\left(t_{1}\right)=\boldsymbol{x}_{1}$
cannot be satisfied, a solution to the optimal control problem does
not exist. Sharp terminal conditions require the controlled dynamical
system to be controllable, i.e., there must exist at least one control
signal enforcing a transfer from the initial state $\boldsymbol{x}\left(t_{0}\right)=\boldsymbol{x}_{0}$
to the terminal state $\boldsymbol{x}\left(t_{1}\right)=\boldsymbol{x}_{1}$.
The case with no terminal conditions, $M\left(\boldsymbol{x}\left(t_{1}\right),t_{1}\right)\equiv0$
and $\boldsymbol{\psi}\left(\boldsymbol{x}\left(t_{1}\right)\right)\equiv\boldsymbol{0}$,
is called a \textit{free} end point condition.

\subsection{\label{sub:DiscussionNecessaryOptimality}Discussion}

The state equation \eqref{eq:StateEquations} and the adjoint equation
\eqref{eq:AdjointEquation2} for the state $\boldsymbol{x}\left(t\right)$
and co-state $\boldsymbol{\lambda}\left(t\right)$, respectively,
as well as the expression Eq. \eqref{eq:ControlSolutionAdjointEquation}
for the vector of control signals $\boldsymbol{u}\left(t\right)$
constitute the \textit{necessary optimality conditions.} Note that
the initial condition for the state equation is specified at the initial
time $t_{0}$, while the initial condition for the adjoint equation
is specified at the terminal time $t_{1}$. These mixed boundary conditions
pose considerable difficulties for a numerical treatment. As will~be
discussed in Section \ref{sec:NumericalSolutionOfOptimalControlProblems},
a straightforward numerical solution is not possible and one usually
has to resort to an iterative scheme.

Even if it is possible to find a solution to the necessary optimality
conditions, this solution can only be considered as a possible candidate
solution for the problem of minimizing the target functional Eq. \eqref{eq:OptimalControlFunctional}.
Similar as for the problem of minimizing an ordinary function, the
necessary optimality conditions only determine an extremum, and sufficient
optimality conditions have to be employed to find out if this candidate
indeed minimizes Eq. \eqref{eq:OptimalControlFunctional}. However,
the question of sufficiency is more subtle than for ordinary functions,
see e.g. \cite{bryson1969applied} and \cite{liberzon2011calculus}.
In this thesis, only necessary optimality conditions are discussed.

Here we discuss only the problem of minimizing a target functional
Eq. \eqref{eq:OptimalControlFunctional} with a constraint in form
of a controlled dynamical system. Other constraints in form of differential
and algebraic equalities and inequalities can be introduced in optimal
control. For technical applications, it is useful to introduce inequality
constraints for control signals such that the control signal is not
allowed to exceed or undershoot certain thresholds. For some problems,
as e.g. unregularized optimal control problems, Eq. \eqref{eq:OptimalControlFunctional}
with $\epsilon=0$, the existence of a solution to the minimization
problem can only be guaranteed if inequality constraints for the control
signals are taken into account.

Other possible constraints are state constraints. For example, the
dynamical system might describe chemical reactions such that the state
components are to be interpreted as the concentrations of some chemical
species. Naturally, these concentrations must be positive quantities,
and a control signal which decreases the value of a concentration
below zero is physically impossible. While the controlled system alone
might violate the condition of positivity, it is possible to enforce
this condition in the context of optimal control.

Other variations of optimal control problems are sparse controls,
which are especially useful for spatio-temporal control systems \cite{Ryll2011,casas_troeltzsch_sparse_state2013}.
Sparse control means a term of the form $\intop_{t_{0}}^{t_{1}}dt\left|\boldsymbol{u}\left(t\right)\right|$,
with $\left|\boldsymbol{u}\left(t\right)\right|=\sqrt{\boldsymbol{u}^{T}\left(t\right)\boldsymbol{u}\left(t\right)}$,
is added to the functional Eq. \eqref{eq:OptimalControlFunctional}.
This has the interesting effect that the control signal vanishes exactly
for some time intervals, while it has a larger amplitude in others
compared to a control without a sparsity term. In this way, it is
possible to find out at which times the application of a control is
most effective.

Finally, we comment on the role of the regularization term $\frac{\epsilon^{2}}{2}\intop_{t_{0}}^{t_{1}}dt\left(\boldsymbol{u}\left(t\right)\right)^{2}$,
sometimes called a Tikhonov regularization, in the context of optimal
trajectory tracking, i.e., for a target functional of the form
\begin{align}
\mathcal{J}\left[\boldsymbol{x}\left(t\right),\boldsymbol{u}\left(t\right)\right] & =\frac{1}{2}\intop_{t_{0}}^{t_{1}}dt\left(\boldsymbol{x}\left(t\right)-\boldsymbol{x}_{d}\left(t\right)\right)^{T}\boldsymbol{\mathcal{S}}\left(\boldsymbol{x}\left(t\right)-\boldsymbol{x}_{d}\left(t\right)\right)\nonumber \\
 & +\frac{1}{2}\left(\boldsymbol{x}\left(t_{1}\right)-\boldsymbol{x}_{1}\right)^{T}\boldsymbol{\mathcal{S}}_{1}\left(\boldsymbol{x}\left(t_{1}\right)-\boldsymbol{x}_{1}\right)+\frac{\epsilon^{2}}{2}\intop_{t_{0}}^{t_{1}}dt\left(\boldsymbol{u}\left(t\right)\right)^{2}.\label{eq:OptimalTrajectoryTrackingFunctional}
\end{align}
The general effect of the regularization term is to penalize large
control values. Without any inequality constraints on the control,
a finite value $\epsilon>0$ is usually necessary to guarantee the
existence of a solution to the minimization problem in terms of bounded
and continuous state trajectories $\boldsymbol{x}\left(t\right)$.
On the other hand, it usually increases the stability and accuracy
of numerical computations of an optimal control. Nevertheless, the
case $\epsilon=0$, called an unregularized optimal control, is of
special interest. For a fixed value of $\epsilon\geq0$, among all
\textit{possible} control signals, the corresponding optimal control
signal is the one which brings the controlled state closest to the
desired trajectory $\boldsymbol{x}_{d}\left(t\right)$ as measured
by Eq. \eqref{eq:OptimalTrajectoryTrackingFunctional}. Furthermore,
among all \textit{optimal} controls, the optimal control for $\epsilon=0$
brings the controlled state closest to the desired trajectory $\boldsymbol{x}_{d}\left(t\right)$.
In other words: the distance measure Eq. \eqref{eq:OptimalTrajectoryTrackingFunctional},
considered as a function of $\epsilon$, has a minimum for $\epsilon=0$.
A proof of this fact is relatively simple. The total derivative of
the augmented functional $\bar{\mathcal{J}}$ with respect to $\epsilon$
is 
\begin{align}
 & \dfrac{d}{d\epsilon}\bar{\mathcal{J}}\left[\boldsymbol{x}\left(t\right),\boldsymbol{u}\left(t\right),\boldsymbol{\lambda}\left(t\right),\boldsymbol{\nu}\right]\nonumber \\
= & \dfrac{\delta}{\delta\boldsymbol{x}\left(t\right)}\bar{\mathcal{J}}\left[\boldsymbol{x}\left(t\right),\boldsymbol{u}\left(t\right),\boldsymbol{\lambda}\left(t\right),\boldsymbol{\nu}\right]\dfrac{d}{d\epsilon}\boldsymbol{x}\left(t\right)+\dfrac{\delta}{\delta\boldsymbol{u}\left(t\right)}\bar{\mathcal{J}}\left[\boldsymbol{x}\left(t\right),\boldsymbol{u}\left(t\right),\boldsymbol{\lambda}\left(t\right),\boldsymbol{\nu}\right]\dfrac{d}{d\epsilon}\boldsymbol{u}\left(t\right)\nonumber \\
 & +\dfrac{\delta}{\delta\boldsymbol{\lambda}\left(t\right)}\bar{\mathcal{J}}\left[\boldsymbol{x}\left(t\right),\boldsymbol{u}\left(t\right),\boldsymbol{\lambda}\left(t\right),\boldsymbol{\nu}\right]\dfrac{d}{d\epsilon}\boldsymbol{\lambda}\left(t\right)+\dfrac{\delta}{\delta\boldsymbol{\nu}}\bar{\mathcal{J}}\left[\boldsymbol{x}\left(t\right),\boldsymbol{u}\left(t\right),\boldsymbol{\lambda}\left(t\right),\boldsymbol{\nu}\right]\dfrac{d}{d\epsilon}\boldsymbol{\nu}\nonumber \\
 & +\epsilon\intop_{t_{0}}^{t_{1}}dt\left(\boldsymbol{u}\left(t\right)\right)^{2}.\label{eq:Eq345}
\end{align}
The last term $\dfrac{\partial}{\partial\epsilon}\bar{\mathcal{J}}=\epsilon\intop_{t_{0}}^{t_{1}}dt\left(\boldsymbol{u}\left(t\right)\right)^{2}$
is due to the explicit dependence of $\bar{\mathcal{J}}$ on $\epsilon$.
If $\boldsymbol{\lambda}\left(t\right)$, $\boldsymbol{\lambda}\left(t_{1}\right)$,
$\boldsymbol{u}\left(t\right),\,\boldsymbol{\nu},$ and $\boldsymbol{x}\left(t\right)$
satisfy the necessary optimality conditions, Eq. \eqref{eq:Eq345}
reduces to 
\begin{align}
\dfrac{d}{d\epsilon}\bar{\mathcal{J}}\left[\boldsymbol{x}\left(t\right),\boldsymbol{u}\left(t\right),\boldsymbol{\lambda}\left(t\right),\boldsymbol{\nu}\right] & =\epsilon\intop_{t_{0}}^{t_{1}}dt\left(\boldsymbol{u}\left(t\right)\right)^{2}.
\end{align}
Then
\begin{align}
\dfrac{d}{d\epsilon}\bar{\mathcal{J}}\left[\boldsymbol{x}\left(t\right),\boldsymbol{u}\left(t\right),\boldsymbol{\lambda}\left(t\right),\boldsymbol{\nu}\right] & =0
\end{align}
for a non-vanishing control signal if and only if $\epsilon=0$. In
other words, $\bar{\mathcal{J}}$ attains an extremum at $\epsilon=0$,
and if $\bar{\mathcal{J}}$ is a minimum with respect to $\boldsymbol{\lambda}\left(t\right)$,
$\boldsymbol{\lambda}\left(t_{1}\right)$, $\boldsymbol{u}\left(t\right),\,\boldsymbol{\nu},$
and $\boldsymbol{x}\left(t\right)$, then it is also a minimum with
respect to $\epsilon$ because of $\intop_{t_{0}}^{t_{1}}dt\left(\boldsymbol{u}\left(t\right)\right)^{2}>0$.

Furthermore, any additional inequality constraints for state or control
can only lead to a value of $\mathcal{J}$ smaller than or equal to
its minimal value attained for $\epsilon=0$. The case with $\epsilon=0$
can be seen as the \textit{limit of realizability} of a certain desired
trajectory $\boldsymbol{x}_{d}\left(t\right)$. No other control,
be it open or closed loop control, can enforce a state trajectory
$\boldsymbol{x}\left(t\right)$ with a smaller distance to the desired
state trajectory $\boldsymbol{x}_{d}\left(t\right)$ than an unregularized
($\epsilon=0$) optimal control. However, assuming $\epsilon=0$ leads
to a singular optimal control problems involving additional difficulties.

\section{\label{sec:SingularOptimalControls}Singular optimal control}

Singular optimal control problems \cite{bell1975singular,bryson1969applied}
are best discussed in terms of the control Hamiltonian $H\left(\boldsymbol{x},\boldsymbol{u},t\right)$.
As long as $H\left(\boldsymbol{x},\boldsymbol{u},t\right)$ depends
only linearly on the vector of control signals $\boldsymbol{u}$,
the optimal control problem is singular. For trajectory tracking tasks
in affine control systems, the control Hamiltonian is defined as 
\begin{align}
H\left(\boldsymbol{x},\boldsymbol{u},t\right) & =\dfrac{1}{2}\left(\boldsymbol{x}\left(t\right)-\boldsymbol{x}_{d}\left(t\right)\right)^{T}\boldsymbol{\mathcal{S}}\left(\boldsymbol{x}\left(t\right)-\boldsymbol{x}_{d}\left(t\right)\right)\nonumber \\
 & +\frac{\epsilon^{2}}{2}\left|\boldsymbol{u}\right|^{2}+\boldsymbol{\lambda}^{T}\left(t\right)\left(\boldsymbol{R}\left(\boldsymbol{x}\right)+\boldsymbol{\mathcal{B}}\left(\boldsymbol{x}\right)\boldsymbol{u}\right).\label{eq:ControlHamiltonian}
\end{align}
The control Hamiltonian is quadratic in the control signal $\boldsymbol{u}\left(t\right)$
as long as $\epsilon>0$. The algebraic relation Eq. \eqref{eq:Eq310}
between control signal and co-state is obtained from the condition
of a stationary Hamiltonian with respect to control,
\begin{align}
\boldsymbol{0} & =\left(\nabla_{\boldsymbol{u}}H\left(\boldsymbol{x},\boldsymbol{u},t\right)\right)^{T}=\epsilon^{2}\boldsymbol{u}+\boldsymbol{\mathcal{B}}^{T}\left(\boldsymbol{x}\right)\boldsymbol{\lambda}.\label{eq:Eq310}
\end{align}
Clearly, if the regularization parameter $\epsilon=0$, $H$ depends
only linearly on the control signal, and Eq. \eqref{eq:Eq310} reduces
to 
\begin{align}
\boldsymbol{0} & =\boldsymbol{\mathcal{B}}^{T}\left(\boldsymbol{x}\right)\boldsymbol{\lambda}.\label{eq:Eq311}
\end{align}
While Eq. \eqref{eq:Eq310} can be used to obtain the control signal
$\boldsymbol{u}\left(t\right)$ in terms of the state $\boldsymbol{x}\left(t\right)$
and co-state $\boldsymbol{\lambda}\left(t\right)$, this is clearly
impossible for Eq. \eqref{eq:Eq311}. Additional necessary optimality
condition, known as the Kelly or generalized Legendre-Clebsch condition,
must be employed to determine an expression for the control signal.

\subsection{\label{sub:KellyConditions}The Kelly condition}

It can be rigorously proven that additional necessary optimality condition
besides the usual necessary optimality conditions have to be satisfied
in case of singular optimal controls \cite{bell1975singular}. This
condition is known as the Kelly condition in case of single-component
control signals and as the generalized Legendre-Clebsch condition
in case of multi-component control signals. For simplicity, only the
case of scalar control signals $u\left(t\right)$ is considered in
this section. The singular control Hamiltonian for optimal trajectory
tracking, Eq. \eqref{eq:ControlHamiltonian} with $\epsilon=0$, becomes
\begin{align}
H\left(\boldsymbol{x},u,t\right) & =\dfrac{1}{2}\left(\boldsymbol{x}\left(t\right)-\boldsymbol{x}_{d}\left(t\right)\right)^{T}\boldsymbol{\mathcal{S}}\left(\boldsymbol{x}\left(t\right)-\boldsymbol{x}_{d}\left(t\right)\right)+\boldsymbol{\lambda}^{T}\left(\boldsymbol{R}\left(\boldsymbol{x}\right)+\boldsymbol{B}\left(\boldsymbol{x}\right)u\right).
\end{align}
The Kelly condition is \cite{bell1975singular}
\begin{align}
\left(-1\right)^{k}\dfrac{\partial}{\partial u}\left[\dfrac{d^{2k}}{dt^{2k}}\dfrac{\partial}{\partial u}H\left(\boldsymbol{x}\left(t\right),u\left(t\right),t\right)\right] & \geq0,\,k=1,2,\dots,\label{eq:KellyCondition}
\end{align}
and is utilized as follows. The stationarity condition $\partial_{u}H=0$,
or, equivalently,
\begin{align}
0 & =\boldsymbol{\lambda}^{T}\left(t\right)\boldsymbol{B}\left(\boldsymbol{x}\left(t\right)\right),\label{eq:Eq334}
\end{align}
is valid for all times $t$ but cannot be used to obtain an expression
for the control signal $u\left(t\right)$. Applying the time derivative
to Eq. \eqref{eq:Eq334} yields $\dfrac{d}{dt}\partial_{u}H=0,$ or
\begin{align}
0 & =\boldsymbol{\dot{\lambda}}^{T}\left(t\right)\boldsymbol{B}\left(\boldsymbol{x}\left(t\right)\right)+\boldsymbol{\lambda}^{T}\left(t\right)\nabla\boldsymbol{B}\left(\boldsymbol{x}\left(t\right)\right)\boldsymbol{\dot{x}}\left(t\right)\nonumber \\
 & =\boldsymbol{\lambda}^{T}\left(t\right)\boldsymbol{q}\left(\boldsymbol{x}\left(t\right)\right)-\left(\boldsymbol{x}\left(t\right)-\boldsymbol{x}_{d}\left(t\right)\right)^{T}\boldsymbol{\mathcal{S}}\boldsymbol{B}\left(\boldsymbol{x}\left(t\right)\right).\label{eq:Eq315}
\end{align}
The controlled state equation as well as the co-state equations was
used and $\boldsymbol{q}\left(\boldsymbol{x}\right)$ denotes the
abbreviation
\begin{align}
\boldsymbol{q}\left(\boldsymbol{x}\right) & =\nabla\boldsymbol{B}\left(\boldsymbol{x}\right)\boldsymbol{R}\left(\boldsymbol{x}\right)-\nabla\boldsymbol{R}\left(\boldsymbol{x}\right)\boldsymbol{B}\left(\boldsymbol{x}\right).\label{eq:Defq}
\end{align}
Equation \eqref{eq:Eq315} yields an additional relation between
state $\boldsymbol{x}$ and co-state $\boldsymbol{\lambda}$, but
does not depend on the control signal $u\left(t\right)$. Therefore,
it cannot be used to obtain an expression for $u\left(t\right)$.
Applying the second time derivative $\dfrac{d^{2}}{dt^{2}}\partial_{u}H=0$
to the stationarity condition yields
\begin{align}
\frac{d^{2}}{dt^{2}}\partial_{u}H & =\boldsymbol{\dot{\lambda}}^{T}\boldsymbol{q}\left(\boldsymbol{x}\right)+\boldsymbol{\lambda}^{T}\nabla\boldsymbol{q}\left(\boldsymbol{x}\right)\boldsymbol{\dot{x}}-\left(\boldsymbol{\dot{x}}-\boldsymbol{\dot{x}}_{d}\right)^{T}\boldsymbol{\mathcal{S}}\boldsymbol{B}\left(\boldsymbol{x}\right)-\left(\boldsymbol{x}-\boldsymbol{x}_{d}\right)^{T}\boldsymbol{\mathcal{S}}\nabla\boldsymbol{B}\left(\boldsymbol{x}\right)\boldsymbol{\dot{x}}\nonumber \\
 & =\boldsymbol{\lambda}^{T}\left(\nabla\boldsymbol{q}\left(\boldsymbol{x}\right)\boldsymbol{R}\left(\boldsymbol{x}\right)-\nabla\boldsymbol{R}\left(\boldsymbol{x}\right)\boldsymbol{q}\left(\boldsymbol{x}\right)\right)+\boldsymbol{B}^{T}\left(\boldsymbol{x}\right)\boldsymbol{\mathcal{S}}\left(\boldsymbol{\dot{x}}_{d}-\boldsymbol{R}\left(\boldsymbol{x}\right)\right)\nonumber \\
 & -\left(\boldsymbol{x}-\boldsymbol{x}_{d}\right)^{T}\boldsymbol{\mathcal{S}}\left(\nabla\boldsymbol{B}\left(\boldsymbol{x}\right)\boldsymbol{R}\left(\boldsymbol{x}\right)+\boldsymbol{q}\left(\boldsymbol{x}\right)\right)+p\left(\boldsymbol{x}\right)u.\label{eq:dtdtduH}
\end{align}
Here, $p\left(\boldsymbol{x}\right)$ denotes the abbreviation
\begin{align}
p\left(\boldsymbol{x}\right) & =\boldsymbol{\lambda}^{T}\left(\nabla\boldsymbol{q}\left(\boldsymbol{x}\right)\boldsymbol{B}\left(\boldsymbol{x}\right)-\nabla\boldsymbol{B}\left(\boldsymbol{x}\right)\boldsymbol{q}\left(\boldsymbol{x}\right)\right)-\boldsymbol{B}^{T}\left(\boldsymbol{x}\right)\boldsymbol{\mathcal{S}}\boldsymbol{B}\left(\boldsymbol{x}\right)\nonumber \\
 & -\left(\boldsymbol{x}-\boldsymbol{x}_{d}\right)^{T}\boldsymbol{\mathcal{S}}\nabla\boldsymbol{B}\left(\boldsymbol{x}\right)\boldsymbol{B}\left(\boldsymbol{x}\right).
\end{align}
Equation \eqref{eq:dtdtduH} does depend on the control, and as long
as $p\left(\boldsymbol{x}\right)\neq0$, it can be solved for the
control signal $u$,
\begin{align}
u & =\dfrac{1}{p\left(\boldsymbol{x}\right)}\left(\boldsymbol{\lambda}^{T}\left(\nabla\boldsymbol{q}\left(\boldsymbol{x}\right)\boldsymbol{R}\left(\boldsymbol{x}\right)-\nabla\boldsymbol{R}\left(\boldsymbol{x}\right)\boldsymbol{q}\left(\boldsymbol{x}\right)\right)+\boldsymbol{B}^{T}\left(\boldsymbol{x}\right)\boldsymbol{\mathcal{S}}\left(\boldsymbol{\dot{x}}_{d}-\boldsymbol{R}\left(\boldsymbol{x}\right)\right)\right)\nonumber \\
 & -\dfrac{1}{p\left(\boldsymbol{x}\right)}\left(\boldsymbol{x}-\boldsymbol{x}_{d}\right)^{T}\boldsymbol{\mathcal{S}}\left(\nabla\boldsymbol{B}\left(\boldsymbol{x}\right)\boldsymbol{R}\left(\boldsymbol{x}\right)+\boldsymbol{q}\left(\boldsymbol{x}\right)\right).\label{eq:KellyControl}
\end{align}
If $\frac{d^{2}}{dt^{2}}\partial_{u}H$ would not depend on $u$,
the time derivative must be applied repeatedly to the stationarity
condition $\partial_{u}H$ until an expression depending on $u$ is
generated. It can be shown that scalar control signals only appear
in even orders of the total time derivative. This is the reason for
the term $\dfrac{d^{2k}}{dt^{2k}}$ in the Kelly conditions Eq. \eqref{eq:KellyCondition}
\cite{bell1975singular}. Finally, one has to check for a generalized
convexity condition, 
\begin{align}
\left(-1\right)^{k}\partial_{u}\left[\frac{d^{2k}}{dt^{2k}}\partial_{u}H\right] & >0.
\end{align}
For Eq. \eqref{eq:dtdtduH} with $k=1$, the generalized convexity
condition
\begin{align}
-\partial_{u}\left[\frac{d^{2}}{dt^{2}}\partial_{u}H\right] & =p\left(\boldsymbol{x}\right)>0.
\end{align}
Thus, as long as $p\left(\boldsymbol{x}\right)>0$, the control signal
given by Eq. \eqref{eq:KellyControl} satisfies all necessary optimality
conditions. The procedure is discussed with a simple example, namely
a mechanical control system in one spatial dimension with vanishing
external force.

\begin{example}[Singular optimal control of a  free particle]\label{ex:FreeParticle1}

Consider the Newton's equation of motion for the position $x$ of
a free point mass in one spatial dimension under the influence of
a control force $u$,
\begin{align}
\ddot{x}\left(t\right) & =u\left(t\right).\label{eq:1DimFreeParticle}
\end{align}
Written as a dynamical system, Eq. \eqref{eq:1DimFreeParticle} becomes
\begin{align}
\dot{x}\left(t\right) & =y\left(t\right),\label{eq:XStateExact}\\
\dot{y}\left(t\right) & =u\left(t\right).\label{eq:YStateExact}
\end{align}
The optimal control task is to minimize the constrained functional
Eq. \eqref{eq:OptimalTrajectoryTrackingFunctional} without any regularization
term. For simplicity, the desired trajectories are chosen to vanish,
\begin{align}
x_{d}\left(t\right) & \equiv0, & y_{d}\left(t\right) & \equiv0.
\end{align}
The initial time is set to $t_{0}=0$. The optimal control task reduces
to the minimization of the functional
\begin{align}
\mathcal{J}\left[\boldsymbol{x}\left(t\right),u\left(t\right)\right] & =\frac{1}{2}\intop_{0}^{t_{1}}\left(\left(x\left(t\right)\right)^{2}+\left(y\left(t\right)\right)^{2}\right)dt,
\end{align}
subject to the dynamics Eqs. \eqref{eq:XStateExact} and \eqref{eq:YStateExact}
with sharp terminal conditions and zero initial conditions,
\begin{align}
x\left(0\right) & =0, & y\left(0\right) & =0,\label{eq:InitCond}\\
x\left(t_{1}\right) & =x_{1}, & y\left(t_{1}\right) & =y_{1}.\label{eq:TermCond}
\end{align}
Thus, the state is required to reach the terminal state $\boldsymbol{x}_{1}=\left(\begin{array}{cc}
x_{1}, & y_{1}\end{array}\right)^{T}$exactly. The necessary optimality conditions are
\begin{align}
\left(\begin{array}{c}
\dot{x}\left(t\right)\\
\dot{y}\left(t\right)
\end{array}\right) & =\left(\begin{array}{c}
y\left(t\right)\\
0
\end{array}\right)+\left(\begin{array}{c}
0\\
u\left(t\right)
\end{array}\right),\\
-\left(\begin{array}{c}
\dot{\lambda}_{x}\left(t\right)\\
\dot{\lambda}_{y}\left(t\right)
\end{array}\right) & =\left(\begin{array}{cc}
0 & 0\\
1 & 0
\end{array}\right)\left(\begin{array}{c}
\lambda_{x}\left(t\right)\\
\lambda_{y}\left(t\right)
\end{array}\right)+\left(\begin{array}{c}
x\left(t\right)\\
y\left(t\right)
\end{array}\right).
\end{align}
The stationarity condition $\partial_{u}H=0$ becomes
\begin{align}
\lambda_{y}\left(t\right) & =0.\label{eq:FreeParticleStationarity}
\end{align}
Applying the time derivative to Eq. \eqref{eq:FreeParticleStationarity},
and using the co-state equation to eliminate $\dot{\lambda}_{y}\left(t\right)$
yields a relation between $\lambda_{x}$ and $y$ as
\begin{align}
0 & =\dfrac{d}{dt}\partial_{u}H=\dot{\lambda}_{y}\left(t\right)=-\lambda_{x}\left(t\right)-y\left(t\right).
\end{align}
Applying the second order time derivative to the stationary condition
Eq. \eqref{eq:FreeParticleStationarity},
\begin{align}
0=\frac{d^{2}}{dt^{2}}\partial_{u}H & =-\dot{\lambda}_{x}\left(t\right)-\dot{y}\left(t\right)=x\left(t\right)-u\left(t\right),
\end{align}
yields an expression for the control signal
\begin{align}
u\left(t\right) & =x\left(t\right).
\end{align}
Finally, the generalized convexity condition yields
\begin{align}
-\partial_{u}\left[\frac{d^{2}}{dt^{2}}\partial_{u}H\right] & =1>0,
\end{align}
and all conditions encoded in the Kelly condition are satisfied. In
summary, the co-states are governed by two algebraic equations
\begin{align}
\lambda_{x}\left(t\right) & =-y\left(t\right), & \lambda_{y}\left(t\right) & =0,
\end{align}
while the control signal is obtained in terms of the controlled state
component $x$ as 
\begin{align}
u\left(t\right) & =x\left(t\right).
\end{align}
The state components $x$ and $y$ are governed by two coupled linear
differential equations of first order,
\begin{align}
\dot{x}\left(t\right) & =y\left(t\right), & \dot{y}\left(t\right) & =x\left(t\right).\label{eq:Eq338}
\end{align}
The general solution to Eqs. \eqref{eq:Eq338}, with two constants
of integration $C_{1}$ and $C_{2}$, is
\begin{align}
x\left(t\right) & =C_{1}\cosh\left(t\right)+C_{2}\sinh\left(t\right),\\
y\left(t\right) & =C_{1}\sinh\left(t\right)+C_{2}\cosh\left(t\right).
\end{align}
However, only two out of four boundary conditions Eqs. \eqref{eq:InitCond}
and \eqref{eq:TermCond} can be satisfied the solution. How to resolve
this problem is demonstrated later on by investigating the same problem
in the limit of small regularization parameter $\epsilon\rightarrow0$.

\end{example}

\section{\label{sec:NumericalSolutionOfOptimalControlProblems}Numerical solution
of optimal control problems}

A numerical solution of an optimal control problem requires the simultaneous
solution of
\begin{enumerate}
\item the $n$ state variables $\boldsymbol{x}\left(t\right)$ as solution
to the controlled state equation \eqref{eq:StateEquations},
\item the $n$ adjoint state variables $\boldsymbol{\lambda}\left(t\right)$
as solution to the adjoint equation \eqref{eq:AdjointEquation2},
\item the $p$ algebraic equations for the control vector $\boldsymbol{u}\left(t\right)$
as given by \eqref{eq:ControlSolutionAdjointEquation}.
\end{enumerate}
A straightforward solution of this coupled system of equations is
rarely possible. The problem is that the boundary conditions for the
controlled state equation, $\boldsymbol{x}\left(t_{0}\right)=\boldsymbol{x}_{0}$,
are given at the initial time $t=t_{0}$, while the adjoint equation
is to be solved with the terminal condition $\boldsymbol{\lambda}\left(t_{1}\right)=\nabla M^{T}\left(\boldsymbol{x}\left(t_{1}\right)\right)+\nabla\boldsymbol{\psi}^{T}\left(\boldsymbol{x}\left(t_{1}\right)\right)\boldsymbol{\nu}$,
Eq. \eqref{eq:AdjointEquation2TermCond}, given at the terminal time
$t=t_{1}$. This typically requires an iterative solution algorithm
similar to the shooting method.

Another problem is that the adjoint equation \eqref{eq:AdjointEquation2}
yields an unstable time evolution, which in turn leads to an unstable
numerical algorithm. This problem can be tackled by solving the time-reversed
adjoint equation. Introducing the new time 
\begin{align}
\tilde{t} & =t_{1}-t, & 0 & \leq\tilde{t}\leq t_{1}-t_{0},
\end{align}
and new adjoint state variables
\begin{align}
\boldsymbol{\tilde{\lambda}}\left(\tilde{t}\right) & =\boldsymbol{\tilde{\lambda}}\left(t_{1}-t\right)=\boldsymbol{\lambda}\left(t\right),
\end{align}
the adjoint equation \eqref{eq:AdjointEquation2} is transformed to
a new equation 
\begin{align}
\boldsymbol{\dot{\tilde{\lambda}}}\left(\tilde{t}\right) & =\left(\nabla\boldsymbol{R}^{T}\left(\boldsymbol{x}\left(t_{1}-\tilde{t}\right)\right)+\boldsymbol{u}^{T}\left(t_{1}-\tilde{t}\right)\nabla\boldsymbol{\mathcal{B}}^{T}\left(\boldsymbol{x}\left(t_{1}-\tilde{t}\right)\right)\right)\boldsymbol{\tilde{\lambda}}\left(\tilde{t}\right)\nonumber \\
 & +\left(\nabla L\left(\boldsymbol{x}\left(t_{1}-\tilde{t}\right),t_{1}-\tilde{t}\right)\right)^{T},\label{eq:TransformedAdjointEquation}
\end{align}
with initial condition
\begin{align}
\boldsymbol{\tilde{\lambda}}\left(0\right) & =\nabla M^{T}\left(\boldsymbol{x}\left(t_{1}\right)\right)+\nabla\boldsymbol{\psi}^{T}\left(\boldsymbol{x}\left(t_{1}\right)\right)\boldsymbol{\nu}.
\end{align}
If the original equation \eqref{eq:AdjointEquation2} yields an unstable
time evolution, then Eq. \eqref{eq:TransformedAdjointEquation} yields
a stable time evolution, and vice versa. Although the transformed
adjoint equation \eqref{eq:TransformedAdjointEquation} as well as
the controlled state equation \eqref{eq:StateEquations} are both
to be solved with initial conditions, they cannot be solved straightforwardly:
while $\boldsymbol{\tilde{\lambda}}$ in Eq. \eqref{eq:TransformedAdjointEquation}
is evaluated at $\tilde{t}$, the state variable $\boldsymbol{x}$
is evaluated at the reversed time $t_{1}-\tilde{t}$. This again illustrates
the problem posed by boundary conditions defined at different times.
These difficulties cannot be resolved by a simple time reversion of
the adjoint equation.

A relatively simple algorithm to solve the optimal control problem
is a first order gradient algorithm. Starting with an initial guess
for the control as e.g. $\boldsymbol{u}\left(t\right)\equiv\boldsymbol{0}$,
the iterative algorithm proceeds as follows (the integer $k$ denotes
the $k$-th iterate):
\begin{enumerate}
\item solve the controlled state equation to obtain the controlled state
$\boldsymbol{x}^{k}\left(t\right)$ from Eq. \eqref{eq:StateEquations},
\item use $\boldsymbol{x}^{k}\left(t\right)$ in Eq. \eqref{eq:AdjointEquation2}
to obtain $\boldsymbol{\lambda}^{k}\left(t\right)$
\item compute a new control $\boldsymbol{u}^{k+1}\left(t\right)$ with the
help of Eq. \eqref{eq:ControlSolutionAdjointEquation} as
\begin{align}
\boldsymbol{u}^{k+1}\left(t\right) & =\boldsymbol{u}^{k}\left(t\right)-s\left(\epsilon^{2}\boldsymbol{u}^{k}\left(t\right)+\boldsymbol{\mathcal{B}}^{T}\left(\boldsymbol{x}^{k}\left(t\right)\right)\boldsymbol{\lambda}^{k}\left(t\right)\right)
\end{align}

\item set $k=k+1$ and go to 1.
\end{enumerate}
The idea of this iterative algorithm is to change the control in the
correct ``direction'' in function space such that the control converges
to the optimal solution. The step width $s$ is an important quantity.
It is often chosen adaptively, with large step widths for the first
couple of iterations and progressively smaller step widths as the
solution for the control converges. Depending on the type of problem,
some hundred up to many hundred thousand of iterations have to be
performed to find a sufficiently correct solution. This renders optimal
control algorithms computationally expensive, and prevents application
of optimizations in real time for processes which are too fast.

The ACADO Toolkit \cite{acadoManual,Houska2011a,Houska2011b} is a
readily available open source package to solve optimal control problems.
If not stated otherwise, this toolkit is used for all numerical solutions
of optimal control throughout the thesis. A typical computation for
a dynamical system with two state components on a time interval of
length $1$ and step width $\Delta t=10^{-3}$, as shown in Example
\ref{ex:ControlledFHN1}, takes about half an hour on a standard laptop.

Many varieties and improvements of the algorithm sketched above can
be found in the literature, as e.g. conjugated gradient method, see
\cite{shewchuk1994introduction} and references therein. Other algorithms
to solve optimal control problems exist, as e.g. the Newton-Raphson
root finding algorithm, see \cite{jorge1999numerical} and \cite{bryson1969applied}
for an overview and examples.

\section{\label{sec:ExactlyRealizableTrajectoriesOptimalControl}Exactly realizable
trajectories and optimal control}

This section discusses the conditions under which exactly realizable
trajectories are optimal. In particular, an exactly realizable desired
trajectory together with its corresponding control signal satisfies
all necessary optimality conditions of a singular optimal control
problem.

\subsection{Usual necessary optimality conditions}

The approach to control in terms of exactly realizable trajectories
from Chapter \ref{chap:ExactlyRealizableTrajectories} is closely
related to an unregularized optimal control problem. Let the target
functional be
\begin{align}
\mathcal{J}\left[\boldsymbol{x}\left(t\right),\boldsymbol{u}\left(t\right)\right]= & \frac{1}{2}\intop_{t_{0}}^{t_{1}}dt\left(\boldsymbol{x}\left(t\right)-\boldsymbol{x}_{d}\left(t\right)\right)^{T}\boldsymbol{\mathcal{S}}\left(\boldsymbol{x}\left(t\right)-\boldsymbol{x}_{d}\left(t\right)\right)\nonumber \\
 & +\frac{1}{2}\left(\boldsymbol{x}\left(t_{1}\right)-\boldsymbol{x}_{1}\right)\boldsymbol{\mathcal{S}}_{1}\left(\boldsymbol{x}\left(t_{1}\right)-\boldsymbol{x}_{1}\right)+\frac{\epsilon^{2}}{2}\intop_{t_{0}}^{t_{1}}dt\left|\boldsymbol{u}\left(t\right)\right|^{2}.\label{eq:JFunctional}
\end{align}
Here, $\boldsymbol{x}_{d}\left(t\right)$ is the desired trajectory,
$\boldsymbol{\mathcal{S}}$ and $\boldsymbol{\mathcal{S}}_{1}$ are
symmetric positive definite matrices, and $\boldsymbol{x}_{1}$ is
a desired terminal state. The necessary optimality conditions comprise
the controlled state equation 
\begin{align}
\boldsymbol{\dot{x}}\left(t\right) & =\boldsymbol{R}\left(\boldsymbol{x}\left(t\right)\right)+\boldsymbol{\mathcal{B}}\left(\boldsymbol{x}\left(t\right)\right)\boldsymbol{u}\left(t\right),\label{eq:KellyControlledStateEquation}\\
\boldsymbol{x}\left(t_{0}\right) & =\boldsymbol{x}_{0},
\end{align}
the adjoint equation 
\begin{align}
-\boldsymbol{\dot{\lambda}}\left(t\right) & =\left(\nabla\boldsymbol{R}^{T}\left(\boldsymbol{x}\left(t\right)\right)+\boldsymbol{u}^{T}\left(t\right)\nabla\boldsymbol{\mathcal{B}}^{T}\left(\boldsymbol{x}\left(t\right)\right)\right)\boldsymbol{\lambda}\left(t\right)+\boldsymbol{\mathcal{S}}\left(\boldsymbol{x}\left(t\right)-\boldsymbol{x}_{d}\left(t\right)\right),\\
\boldsymbol{\lambda}\left(t_{1}\right) & =\boldsymbol{\mathcal{S}}_{1}\left(\boldsymbol{x}\left(t_{1}\right)-\boldsymbol{x}_{1}\right),
\end{align}
and an algebraic relation between co-state $\boldsymbol{\lambda}$
and control $\boldsymbol{u}$,
\begin{align}
\epsilon^{2}\boldsymbol{u}\left(t\right)+\boldsymbol{\mathcal{B}}^{T}\left(\boldsymbol{x}\left(t\right)\right)\boldsymbol{\lambda}\left(t\right)= & 0.\label{eq:KellyStationarityCondition}
\end{align}
Usually, Eq. \eqref{eq:KellyStationarityCondition} is used to determine
the control $\boldsymbol{u}\left(t\right)$. Here, we proceed differently,
and assume an exactly realizable desired trajectory $\boldsymbol{x}_{d}\left(t\right)$
such that the controlled state $\boldsymbol{x}\left(t\right)$ exactly
follows $\boldsymbol{x}_{d}\left(t\right)$ for all times,
\begin{align}
\boldsymbol{x}\left(t\right) & =\boldsymbol{x}_{d}\left(t\right).\label{eq:XEqualsXd}
\end{align}
Starting from this assumption, the necessary optimality conditions
are evaluated to determine conditions on $\boldsymbol{x}_{d}\left(t\right)$
such that Eq. \eqref{eq:XEqualsXd} holds.

First of all, for Eq. \eqref{eq:XEqualsXd} to be valid at all times,
the initial value of the desired trajectory must comply with the initial
state,
\begin{align}
\boldsymbol{x}\left(t_{0}\right) & =\boldsymbol{x}_{0}=\boldsymbol{x}_{d}\left(t_{0}\right).
\end{align}
With assumption Eq. \eqref{eq:XEqualsXd}, the adjoint equation becomes
\begin{align}
-\boldsymbol{\dot{\lambda}}\left(t\right) & =\left(\nabla\boldsymbol{R}^{T}\left(\boldsymbol{x}_{d}\left(t\right)\right)+\boldsymbol{u}^{T}\left(t\right)\nabla\boldsymbol{\mathcal{B}}^{T}\left(\boldsymbol{x}_{d}\left(t\right)\right)\right)\boldsymbol{\lambda}\left(t\right),\\
\boldsymbol{\lambda}\left(t_{1}\right) & =\boldsymbol{\mathcal{S}}_{1}\left(\boldsymbol{x}_{d}\left(t_{1}\right)-\boldsymbol{x}_{1}\right).
\end{align}
If the desired trajectory satisfies 
\begin{align}
\boldsymbol{x}_{d}\left(t_{1}\right) & =\boldsymbol{x}_{1},
\end{align}
the boundary condition for the adjoint equation becomes
\begin{align}
\boldsymbol{\lambda}\left(t_{1}\right)= & \boldsymbol{0}.
\end{align}
Consequently, the co-state $\boldsymbol{\lambda}$ vanishes identically
for all times, 
\begin{align}
\boldsymbol{\lambda}\left(t\right) & \equiv\boldsymbol{0}.\label{eq:LambdaVanish}
\end{align}
With the help of the Moore-Penrose pseudo inverse
\begin{align}
\boldsymbol{\mathcal{B}}^{+}\left(\boldsymbol{x}\right) & =\left(\boldsymbol{\mathcal{B}}^{T}\left(\boldsymbol{x}\right)\boldsymbol{\mathcal{B}}\left(\boldsymbol{x}\right)\right)^{-1}\boldsymbol{\mathcal{B}}^{T}\left(\boldsymbol{x}\right),
\end{align}
the controlled state equation \eqref{eq:KellyControlledStateEquation}
is solved for the control signal 
\begin{align}
\boldsymbol{u}\left(t\right) & =\boldsymbol{\mathcal{B}}^{+}\left(\boldsymbol{x}_{d}\left(t\right)\right)\left(\boldsymbol{\dot{x}}_{d}\left(t\right)-\boldsymbol{R}\left(\boldsymbol{x}_{d}\left(t\right)\right)\right).\label{eq:ControlSolutionOptimal}
\end{align}
Using $\boldsymbol{u}\left(t\right)$ in the controlled state equation
yields the constraint equation
\begin{align}
\boldsymbol{0} & =\boldsymbol{\mathcal{Q}}\left(\boldsymbol{x}_{d}\left(t\right)\right)\left(\boldsymbol{\dot{x}}_{d}\left(t\right)-\boldsymbol{R}\left(\boldsymbol{x}_{d}\left(t\right)\right)\right).\label{eq:ConstraintEquationOptimal}
\end{align}
Finally, because of the vanishing co-state Eq. \eqref{eq:LambdaVanish},
the stationarity condition Eq. \eqref{eq:KellyStationarityCondition}
becomes
\begin{align}
\epsilon^{2}\boldsymbol{u}\left(t\right)= & 0.\label{eq:Eq463}
\end{align}
Clearly, because $\boldsymbol{u}\left(t\right)$ is non-vanishing,
Eq. \eqref{eq:Eq463} can only be satisfied if
\begin{align}
\epsilon & =0.\label{eq:epsilon0}
\end{align}
In conclusion, for the necessary optimality conditions to be valid
under the assumption $\boldsymbol{x}\left(t\right)=\boldsymbol{x}_{d}\left(t\right)$,
the desired trajectory $\boldsymbol{x}_{d}\left(t\right)$ has to
be exactly realizable. Furthermore, $\boldsymbol{x}_{d}\left(t\right)$
must comply with the terminal condition $\boldsymbol{x}_{d}\left(t_{1}\right)=\boldsymbol{x}_{1}$.
This additional condition originates simply from different formulations
of the control task. While the control methods from Chapter \ref{chap:ExactlyRealizableTrajectories}
enforce only initial conditions, optimal control is able to impose
terminal conditions as well. As shown by Eq. \eqref{eq:epsilon0},
exactly realizable trajectories naturally lead to singular optimal
control problems. Hence, additional necessary optimality conditions
in form of the generalized Legendre-Clebsch conditions must be evaluated.

\subsection{\label{sub:TheGeneralizedLegendreClebschConditions}The generalized
Legendre-Clebsch conditions}

The generalized Legendre-Clebsch conditions are \cite{bell1975singular}
\begin{align}
\nabla_{\boldsymbol{u}}\dfrac{d^{k}}{dt^{k}}\left(\nabla_{\boldsymbol{u}}H\right) & =\boldsymbol{0},\,k\in\mathbb{N},\,k\text{ odd},\label{eq:LegendreClebsch1}
\end{align}
and
\begin{align}
\left(-1\right)^{l}\nabla_{\boldsymbol{u}}\dfrac{d^{2l}}{dt^{2l}}\left(\nabla_{\boldsymbol{u}}H\right) & \geq\boldsymbol{0},\,l\in\mathbb{N}.\label{eq:LegendreClebsch2}
\end{align}
These conditions are evaluated in the same manner as the Kelly condition,
see Section \ref{sub:KellyConditions}. Due to the vector character
of the control signal, the computations are more involved, and a different
notation is adopted. Written for the individual state components $x_{i}$,
the controlled state equation is
\begin{align}
\dot{x}_{i}= & R_{i}+\sum_{k=1}^{p}\mathcal{B}_{ik}u_{k}.\label{eq:ControlledStateComponentWise}
\end{align}
For the remainder of this section, the state arguments of $\boldsymbol{R}$
and $\boldsymbol{\mathcal{B}}$ and time arguments of $\boldsymbol{\lambda}$,
$\boldsymbol{x}$, $\boldsymbol{x}_{d}$, and $\boldsymbol{u}$ are
suppressed to shorten the notation. The matrix entries $\mathcal{B}_{ik}$
are assumed to depend on the state $\boldsymbol{x}$. With
\begin{align}
\partial_{j} & =\dfrac{\partial}{\partial x_{j}}, & \left(\nabla\boldsymbol{R}\right)_{ij} & =\partial_{j}R_{i}, & \left(\boldsymbol{\lambda}^{T}\nabla\boldsymbol{\mathcal{B}}\boldsymbol{u}\right)_{k} & =\sum_{i=1}^{n}\sum_{j=1}^{p}\lambda_{i}\partial_{k}\mathcal{B}_{ij}u_{j},
\end{align}
the adjoint equation is written as
\begin{align}
-\dot{\lambda}_{i} & =\sum_{k=1}^{n}\left(\lambda_{k}\partial_{i}R_{k}+\sum_{l=1}^{p}\lambda_{k}\partial_{i}\mathcal{B}_{kl}u_{l}+\left(x_{k}-x_{d,k}\right)\mathcal{S}_{ki}\right).
\end{align}
Here, $x_{d,k}$ denotes the $k$-th component of the desired trajectory
$\boldsymbol{x}_{d}\left(t\right)$. The stationarity condition
\begin{align}
\boldsymbol{0} & =\nabla_{\boldsymbol{u}}H=\boldsymbol{\lambda}^{T}\boldsymbol{\mathcal{B}}
\end{align}
becomes
\begin{align}
0 & =\sum_{i=1}^{n}\lambda_{i}\mathcal{B}_{ij}.\label{eq:StationarityCoords}
\end{align}
The procedure is analogous to the Kelly condition and the time derivative
is applied repeatedly onto Eq. \eqref{eq:StationarityCoords}. The
first condition, Eq. \eqref{eq:LegendreClebsch1} for $k=1$,
\begin{align}
\dfrac{d}{dt}\left(\nabla_{\boldsymbol{u}}H\right) & =\boldsymbol{0},
\end{align}
yields
\begin{align}
0 & =\sum_{i=1}^{n}\left(\dot{\lambda}_{i}\mathcal{B}_{ij}+\sum_{k=1}^{n}\lambda_{i}\partial_{k}\mathcal{B}_{ij}\dot{x}_{k}\right)\nonumber \\
 & =\sum_{l=1}^{p}\sum_{i=1}^{n}\sum_{k=1}^{n}\lambda_{i}\left(\partial_{k}\mathcal{B}_{ij}\mathcal{B}_{kl}-\mathcal{B}_{kj}\partial_{k}\mathcal{B}_{il}\right)u_{l}-\sum_{i=1}^{n}\sum_{k=1}^{n}\left(x_{k}-x_{d,k}\right)\mathcal{S}_{ki}\mathcal{B}_{ij}\nonumber \\
 & +\sum_{i=1}^{n}\sum_{k=1}^{n}\lambda_{i}\left(\partial_{k}\mathcal{B}_{ij}R_{k}-\partial_{k}R_{i}\mathcal{B}_{kj}\right).
\end{align}
Because of $\boldsymbol{\lambda}\left(t\right)\equiv\boldsymbol{0}$
and $\boldsymbol{x}\left(t\right)=\boldsymbol{x}_{d}\left(t\right)$
for all times, this expression is satisfied. Note that the condition
\begin{align}
\nabla_{\boldsymbol{u}}\dfrac{d}{dt}\left(\nabla_{\boldsymbol{u}}H\right) & =\boldsymbol{0},
\end{align}
or
\begin{align}
0 & =\sum_{i=1}^{n}\sum_{k=1}^{n}\lambda_{i}\left(\partial_{k}\mathcal{B}_{ij}\mathcal{B}_{kl}-\mathcal{B}_{kj}\partial_{k}\mathcal{B}_{il}\right),
\end{align}
is only valid under certain symmetry conditions on the coupling matrix
$\boldsymbol{\mathcal{B}}\left(\boldsymbol{x}\right)$ for a finite
co-state $\lambda_{i}\neq0$ \cite{bryson1969applied}. However, here
the co-state vanishes exactly, and it is unnecessary to impose these
symmetry conditions on $\boldsymbol{\mathcal{B}}$. The next Legendre-Clebsch
condition is 
\begin{align}
\dfrac{d^{2}}{dt^{2}}\left(\nabla_{\boldsymbol{u}}H\right) & =\boldsymbol{0},
\end{align}
or
\begin{align}
0 & =\sum_{i=1}^{n}\sum_{k=1}^{n}\left(\sum_{l=1}^{p}\dot{\lambda}_{i}\left(\partial_{k}\mathcal{B}_{ij}\mathcal{B}_{kl}-\mathcal{B}_{kj}\partial_{k}\mathcal{B}_{il}\right)u_{l}+\sum_{l=1}^{p}\lambda_{i}\dfrac{d}{dt}\left(\partial_{k}\mathcal{B}_{ij}\mathcal{B}_{kl}-\mathcal{B}_{kj}\partial_{k}\mathcal{B}_{il}\right)u_{l}\right)\nonumber \\
 & +\sum_{i=1}^{n}\sum_{k=1}^{n}\sum_{l=1}^{p}\lambda_{i}\left(\partial_{k}\mathcal{B}_{ij}\mathcal{B}_{kl}-\mathcal{B}_{kj}\partial_{k}\mathcal{B}_{il}\right)\dot{u}_{l}\nonumber \\
 & -\sum_{i=1}^{n}\sum_{k=1}^{n}\left(\left(R_{k}+\sum_{m=1}^{p}\mathcal{B}_{km}u_{m}-\dot{x}_{d,k}\right)\mathcal{S}_{ki}\mathcal{B}_{ij}+\sum_{m=1}^{n}\left(x_{k}-x_{d,k}\right)\mathcal{S}_{ki}\partial_{m}\mathcal{B}_{ij}\dot{x}_{m}\right)\nonumber \\
 & +\sum_{i=1}^{n}\sum_{k=1}^{n}\left(\dot{\lambda}_{i}\left(\partial_{k}\mathcal{B}_{ij}R_{k}-\partial_{k}R_{i}\mathcal{B}_{kj}\right)+\lambda_{i}\dfrac{d}{dt}\left(\partial_{k}\mathcal{B}_{ij}R_{k}-\partial_{k}R_{i}\mathcal{B}_{kj}\right)\right).\label{eq:Eq395}
\end{align}
The controlled state equation \eqref{eq:ControlledStateComponentWise}
was used to substitute $\dot{x}_{k}$. Because of $\boldsymbol{\lambda}\left(t\right)\equiv\boldsymbol{0}$
and $\boldsymbol{x}\left(t\right)=\boldsymbol{x}_{d}\left(t\right)$
for all times, the Eq. \eqref{eq:Eq395} simplifies considerably and
yields a solution for the control signal
\begin{align}
\sum_{i=1}^{n}\sum_{k=1}^{n}\sum_{m=1}^{p}\mathcal{B}_{ij}\mathcal{S}_{ki}\mathcal{B}_{km}u_{m} & =\sum_{i=1}^{n}\sum_{k=1}^{n}\mathcal{B}_{ij}\mathcal{S}_{ki}\left(\dot{x}_{d,k}-R_{k}\right).\label{eq:Eq428}
\end{align}
Casting Eq. \eqref{eq:Eq428} in terms of vectors and matrices and
exploiting the symmetry of $\boldsymbol{\mathcal{S}}$ gives
\begin{align}
\boldsymbol{\mathcal{B}}^{T}\boldsymbol{\mathcal{S}}\boldsymbol{\mathcal{B}}\boldsymbol{u} & =\boldsymbol{\mathcal{B}}^{T}\boldsymbol{\mathcal{S}}\left(\boldsymbol{\dot{x}}_{d}-\boldsymbol{R}\right).\label{eq:Eq427}
\end{align}
Solving for the control and substituting $\boldsymbol{x}\left(t\right)=\boldsymbol{x}_{d}\left(t\right)$
results in
\begin{align}
\boldsymbol{u}\left(t\right) & =\boldsymbol{\mathcal{B}}_{\boldsymbol{\mathcal{S}}}^{g}\left(\boldsymbol{x}_{d}\left(t\right)\right)\left(\boldsymbol{\dot{x}}_{d}\left(t\right)-\boldsymbol{R}\left(\boldsymbol{x}_{d}\left(t\right)\right)\right).\label{eq:KellyCondControl}
\end{align}
The $p\times n$ matrix $\boldsymbol{\mathcal{B}}_{\boldsymbol{\mathcal{S}}}^{g}\left(\boldsymbol{x}\right)$
is defined by
\begin{align}
\boldsymbol{\mathcal{B}}_{\boldsymbol{\mathcal{S}}}^{g}\left(\boldsymbol{x}\right) & =\left(\boldsymbol{\mathcal{B}}^{T}\left(\boldsymbol{x}\right)\boldsymbol{\mathcal{S}}\boldsymbol{\mathcal{B}}\left(\boldsymbol{x}\right)\right)^{-1}\boldsymbol{\mathcal{B}}^{T}\left(\boldsymbol{x}\right)\boldsymbol{\mathcal{S}}.\label{eq:DefBgS}
\end{align}
Note that the matrix $\boldsymbol{\mathcal{B}}_{\boldsymbol{\mathcal{S}}}^{g}\left(\boldsymbol{x}\right)$
is not the Moore-Penrose pseudo inverse but a generalized reflexive
inverse of $\boldsymbol{\mathcal{B}}\left(\boldsymbol{x}\right)$,
see Appendix \ref{sub:GeneralizedInverseMatrices}. Finally, the generalized
convexity condition
\begin{align}
\nabla_{\boldsymbol{u}}\dfrac{d^{2}}{dt^{2}}\left(\nabla_{\boldsymbol{u}}H\right) & \geq\boldsymbol{0}
\end{align}
is satisfied whenever 
\begin{align}
\boldsymbol{\mathcal{B}}^{T}\left(\boldsymbol{x}\right)\boldsymbol{\mathcal{S}}\boldsymbol{\mathcal{B}}\left(\boldsymbol{x}\right) & >\boldsymbol{0}\label{eq:BSBpositive}
\end{align}
for all $\boldsymbol{x}$. Condition Eq. \eqref{eq:BSBpositive} ensures
that $\left(\boldsymbol{\mathcal{B}}^{T}\left(\boldsymbol{x}\right)\boldsymbol{\mathcal{S}}\boldsymbol{\mathcal{B}}\left(\boldsymbol{x}\right)\right)^{-1}$
in Eq. \eqref{eq:DefBgS} exists. Note that the matrix $\boldsymbol{\mathcal{S}}$
was assumed to be symmetric, but $\boldsymbol{\mathcal{S}}$ does
not need to be positive definite to satisfy Eq. \eqref{eq:BSBpositive}.

The matrix $\boldsymbol{\mathcal{B}}_{\boldsymbol{\mathcal{S}}}^{g}\left(\boldsymbol{x}\right)$
is used to define the two complementary $n\times n$ projectors
\begin{align}
\boldsymbol{\mathcal{P}}_{\boldsymbol{\mathcal{S}}}\left(\boldsymbol{x}\right) & =\boldsymbol{\mathcal{B}}\left(\boldsymbol{x}\right)\boldsymbol{\mathcal{B}}_{\boldsymbol{\mathcal{S}}}^{g}\left(\boldsymbol{x}\right)=\boldsymbol{\mathcal{B}}\left(\boldsymbol{x}\right)\left(\boldsymbol{\mathcal{B}}^{T}\left(\boldsymbol{x}\right)\boldsymbol{\mathcal{S}}\boldsymbol{\mathcal{B}}\left(\boldsymbol{x}\right)\right)^{-1}\boldsymbol{\mathcal{B}}^{T}\left(\boldsymbol{x}\right)\boldsymbol{\mathcal{S}},\label{eq:PSQSDef1}\\
\boldsymbol{\mathcal{Q}}_{\boldsymbol{\mathcal{S}}}\left(\boldsymbol{x}\right) & =\boldsymbol{1}-\boldsymbol{\mathcal{P}}_{\boldsymbol{\mathcal{S}}}\left(\boldsymbol{x}\right).\label{eq:PSQSDef2}
\end{align}
The matrix $\boldsymbol{\mathcal{P}}_{\boldsymbol{\mathcal{S}}}\left(\boldsymbol{x}\right)$
is idempotent but not symmetric,
\begin{align}
\boldsymbol{\mathcal{P}}_{\boldsymbol{\mathcal{S}}}\left(\boldsymbol{x}\right)\boldsymbol{\mathcal{P}}_{\boldsymbol{\mathcal{S}}}\left(\boldsymbol{x}\right) & =\boldsymbol{\mathcal{P}}_{\boldsymbol{\mathcal{S}}}\left(\boldsymbol{x}\right), & \boldsymbol{\mathcal{P}}_{\boldsymbol{\mathcal{S}}}\left(\boldsymbol{x}\right) & \neq\boldsymbol{\mathcal{P}}_{\boldsymbol{\mathcal{S}}}^{T}\left(\boldsymbol{x}\right),
\end{align}
and analogously for $\boldsymbol{\mathcal{Q}}_{\boldsymbol{\mathcal{S}}}\left(\boldsymbol{x}\right)$.
Using the solution Eq. \eqref{eq:KellyCondControl} for the control
in the controlled state equation \eqref{eq:KellyControlledStateEquation}
together with $\boldsymbol{x}\left(t\right)=\boldsymbol{x}_{d}\left(t\right)$
yields
\begin{align}
\boldsymbol{\dot{x}}_{d}\left(t\right)-\boldsymbol{R}\left(\boldsymbol{x}_{d}\left(t\right)\right) & =\boldsymbol{\mathcal{B}}\left(\boldsymbol{x}_{d}\left(t\right)\right)\boldsymbol{u}\left(t\right)\nonumber \\
 & =\boldsymbol{\mathcal{P}}_{\boldsymbol{\mathcal{S}}}\left(\boldsymbol{x}_{d}\left(t\right)\right)\left(\boldsymbol{\dot{x}}_{d}\left(t\right)-\boldsymbol{R}\left(\boldsymbol{x}_{d}\left(t\right)\right)\right),
\end{align}
and finally
\begin{align}
\boldsymbol{\mathcal{Q}}_{\boldsymbol{\mathcal{S}}}\left(\boldsymbol{x}_{d}\left(t\right)\right)\left(\boldsymbol{\dot{x}}_{d}\left(t\right)-\boldsymbol{R}\left(\boldsymbol{x}_{d}\left(t\right)\right)\right) & =\boldsymbol{0}.\label{eq:QSConstraintEquation}
\end{align}
Equation \eqref{eq:QSConstraintEquation} looks very much like the
constraint equation \eqref{eq:ConstraintEquationOptimal} found in
the last section, but with a different projector $\boldsymbol{\mathcal{Q}}_{\boldsymbol{\mathcal{S}}}\left(\boldsymbol{x}\right)$
instead of $\boldsymbol{\mathcal{Q}}\left(\boldsymbol{x}\right)$.
Additionally, the control signal Eq. \eqref{eq:KellyCondControl}
appears unequal from Eq. \eqref{eq:ControlSolutionOptimal} obtained
in the last section. It seems that, for the same problem, two unconnected
control solutions were found. The first control solution in terms
of $\boldsymbol{\mathcal{B}}^{+}$ and $\boldsymbol{\mathcal{Q}}$
is given by
\begin{align}
\boldsymbol{u}_{1}\left(t\right) & =\boldsymbol{\mathcal{B}}^{+}\left(\boldsymbol{x}_{d}\left(t\right)\right)\left(\boldsymbol{\dot{x}}_{d}\left(t\right)-\boldsymbol{R}\left(\boldsymbol{x}_{d}\left(t\right)\right)\right),\\
\boldsymbol{0} & =\boldsymbol{\mathcal{Q}}\left(\boldsymbol{x}_{d}\left(t\right)\right)\left(\boldsymbol{\dot{x}}_{d}\left(t\right)-\boldsymbol{R}\left(\boldsymbol{x}_{d}\left(t\right)\right)\right),\label{eq:ConstEq1}
\end{align}
while the second control solution in terms of $\boldsymbol{\mathcal{B}}_{\boldsymbol{\mathcal{S}}}^{g}$
and $\boldsymbol{\mathcal{Q}}_{\boldsymbol{\mathcal{S}}}$ is 
\begin{align}
\boldsymbol{u}_{2}\left(t\right) & =\boldsymbol{\mathcal{B}}_{\boldsymbol{\mathcal{S}}}^{g}\left(\boldsymbol{x}_{d}\left(t\right)\right)\left(\boldsymbol{\dot{x}}_{d}\left(t\right)-\boldsymbol{R}\left(\boldsymbol{x}_{d}\left(t\right)\right)\right),\\
\boldsymbol{0} & =\boldsymbol{\mathcal{Q}}_{\boldsymbol{\mathcal{S}}}\left(\boldsymbol{x}_{d}\left(t\right)\right)\left(\boldsymbol{\dot{x}}_{d}\left(t\right)-\boldsymbol{R}\left(\boldsymbol{x}_{d}\left(t\right)\right)\right).\label{eq:ConstEq2}
\end{align}
The expressions $\boldsymbol{u}_{1}\left(t\right)$ and $\boldsymbol{u}_{2}\left(t\right)$
are identical for a matrix of weighting coefficients $\boldsymbol{\mathcal{S}}=\boldsymbol{1}$
but seem to disagree for $\boldsymbol{\mathcal{S}}\neq\boldsymbol{1}$.

However, the difference in $\boldsymbol{u}_{2}\left(t\right)$ and
$\boldsymbol{u}_{1}\left(t\right)$ is deceptive. In fact, the expressions
are identical, as is demonstrated in the following. Computing the
difference between $\boldsymbol{u}_{2}\left(t\right)$ and $\boldsymbol{u}_{1}\left(t\right)$,
multiplying by $\boldsymbol{\mathcal{B}}\left(\boldsymbol{x}_{d}\left(t\right)\right)$,
adding $\boldsymbol{0}=\boldsymbol{1}-\boldsymbol{1}$, and exploiting
the constraint equations yields, see also Appendix \ref{sec:OverAndUnderdetSysOfEqs},
\begin{align}
\boldsymbol{\mathcal{B}}\left(\boldsymbol{x}_{d}\left(t\right)\right)\left(\boldsymbol{u}_{1}\left(t\right)-\boldsymbol{u}_{2}\left(t\right)\right) & =\left(\boldsymbol{\mathcal{P}}\left(\boldsymbol{x}_{d}\left(t\right)\right)-\boldsymbol{\mathcal{P}}_{\boldsymbol{\mathcal{S}}}\left(\boldsymbol{x}_{d}\left(t\right)\right)\right)\left(\boldsymbol{\dot{x}}_{d}\left(t\right)-\boldsymbol{R}\left(\boldsymbol{x}_{d}\left(t\right)\right)\right)\nonumber \\
 & =\left(\boldsymbol{\mathcal{Q}}\left(\boldsymbol{x}_{d}\left(t\right)\right)-\boldsymbol{\mathcal{Q}}_{\boldsymbol{\mathcal{S}}}\left(\boldsymbol{x}_{d}\left(t\right)\right)\right)\left(\boldsymbol{\dot{x}}_{d}\left(t\right)-\boldsymbol{R}\left(\boldsymbol{x}_{d}\left(t\right)\right)\right)\nonumber \\
 & =\boldsymbol{0}.\label{eq:Eq3111}
\end{align}
Equation \eqref{eq:Eq3111} implies that either $\boldsymbol{u}_{1}\left(t\right)=\boldsymbol{u}_{2}\left(t\right)$,
or $\boldsymbol{u}_{1}\left(t\right)-\boldsymbol{u}_{2}\left(t\right)$
lies in the null space of $\boldsymbol{\mathcal{B}}\left(\boldsymbol{x}_{d}\left(t\right)\right)$.
Because $\boldsymbol{\mathcal{B}}\left(\boldsymbol{x}\right)$ has
full column rank for all $\boldsymbol{x}$ by assumption, the null
space of $\boldsymbol{\mathcal{B}}\left(\boldsymbol{x}\right)$ contains
only the zero vector. In conclusion,
\begin{align}
\boldsymbol{u}_{1}\left(t\right) & =\boldsymbol{u}_{2}\left(t\right),
\end{align}
and, consequently, the control signal is unique and does not depend
on the matrix of weighting coefficients $\boldsymbol{\mathcal{S}}$.
Because identical control signals enforce identical controlled state
trajectories $\boldsymbol{x}\left(t\right)=\boldsymbol{x}_{d}\left(t\right)$,
the desired trajectories constrained by the different constraint equations
\eqref{eq:ConstEq1} and \eqref{eq:ConstEq1} are identical as well.

In the framework of exactly realizable desired trajectories, the appearance
of alternative projectors $\boldsymbol{\mathcal{P}}_{\boldsymbol{\mathcal{S}}}\left(\boldsymbol{x}\right)$
and $\boldsymbol{\mathcal{Q}}_{\boldsymbol{\mathcal{S}}}\left(\boldsymbol{x}\right)$
plays no role. Neither the control signal nor the controlled state
trajectory depends on the matrix of weighting coefficients $\boldsymbol{\mathcal{S}}$.
However, $\boldsymbol{\mathcal{P}}_{\boldsymbol{\mathcal{S}}}\left(\boldsymbol{x}\right)$
and $\boldsymbol{\mathcal{Q}}_{\boldsymbol{\mathcal{S}}}\left(\boldsymbol{x}\right)$
become important for the perturbative approach to trajectory tracking
of arbitrary desired trajectories in Chapter \ref{chap:AnalyticalApproximationsForOptimalTrajectoryTracking}.
Note that the matrix $\boldsymbol{\mathcal{Q}}_{\boldsymbol{\mathcal{S}}}\left(\boldsymbol{x}\right)$
is similar to the matrix $\boldsymbol{\mathcal{Q}}\left(\boldsymbol{x}\right)$,
i.e., there exists an invertible $n\times n$ matrix $\boldsymbol{\mathcal{T}}\left(\boldsymbol{x}\right)$
such that
\begin{align}
\boldsymbol{\mathcal{Q}}_{\boldsymbol{\mathcal{S}}}\left(\boldsymbol{x}\right) & =\boldsymbol{\mathcal{T}}^{-1}\left(\boldsymbol{x}\right)\boldsymbol{\mathcal{Q}}\left(\boldsymbol{x}\right)\boldsymbol{\mathcal{T}}\left(\boldsymbol{x}\right).
\end{align}
This follows from the fact that both $\boldsymbol{\mathcal{Q}}_{\boldsymbol{\mathcal{S}}}\left(\boldsymbol{x}\right)$
and $\boldsymbol{\mathcal{Q}}\left(\boldsymbol{x}\right)$ are projectors
of rank $n-p$. They have identical eigenvalues and, when diagonalized,
identical diagonal forms $\boldsymbol{\mathcal{Q}}_{D}$. See also
Appendix \ref{sec:DiagonalizingTheConstraintEquation} how to diagonalize
the projectors $\boldsymbol{\mathcal{P}}\left(\boldsymbol{x}\right)$
and $\boldsymbol{\mathcal{Q}}\left(\boldsymbol{x}\right)$.

\section{\label{sec:AnExactlySolvableExample}An exactly solvable example}

A simple linear and exactly solvable example for optimal trajectory
tracking is considered in this section. The rather clumsy exact solution
is simplified by assuming a small regularization parameter $0<\epsilon\ll1$.
A generalized perturbation expansion known as a singular perturbation
expansion is necessary to obtain an approximation which is valid over
the whole time domain $t_{0}\leq t\le t_{1}$. The purpose of analyzing
this exact solution is three-fold. First, it serves as a pedagogical
example displaying similar difficulties as the nonlinear system in
Chapter \ref{chap:AnalyticalApproximationsForOptimalTrajectoryTracking}.
Second, it provides a consistency check for the analytical results
of Chapter \ref{chap:AnalyticalApproximationsForOptimalTrajectoryTracking}.
Third, the impact of different terminal conditions on the exact solution
is analyzed.

\subsection{Problem and exact solution}

Optimal trajectory tracking for a free particle with a finite regularization
coefficient $\epsilon>0$ is considered, see Example \ref{ex:FreeParticle1}.
For simplicity, a vanishing desired trajectory $\boldsymbol{x}_{d}\left(t\right)\equiv\boldsymbol{0}$
and zero initial time $t_{0}=0$ is assumed. The optimal control problem
is to minimize the constrained functional
\begin{align}
\mathcal{J}\left[\boldsymbol{x}\left(t\right),u\left(t\right)\right] & =\frac{1}{2}\intop_{t_{0}}^{t_{1}}\left(\left(x\left(t\right)\right)^{2}+\left(y\left(t\right)\right)^{2}\right)dt+\frac{\epsilon^{2}}{2}\intop_{t_{0}}^{t_{1}}dt\left(u\left(t\right)\right)^{2}\nonumber \\
 & +\dfrac{\beta_{1}}{2}\left(x\left(t_{1}\right)-x_{1}\right)^{2}+\dfrac{\beta_{2}}{2}\left(y\left(t_{1}\right)-y_{1}\right)^{2},\label{eq:FreeParticleFunctional}
\end{align}
subject to the system dynamics and initial conditions 
\begin{align}
\dot{x}\left(t\right) & =y\left(t\right), & \dot{y}\left(t\right) & =u\left(t\right),\label{eq:StateExact}\\
x\left(0\right) & =x_{0}, & y\left(0\right) & =y_{0}.
\end{align}
Note that as long as $\boldsymbol{x}_{0}\neq\boldsymbol{0}$, the
desired trajectory $\boldsymbol{x}_{d}\left(t\right)\equiv\boldsymbol{0}$
does not comply with the initial conditions and is therefore not exactly
realizable. Similarly, if $\boldsymbol{x}_{1}\neq\boldsymbol{0}$,
the desired trajectory $\boldsymbol{x}_{d}\left(t\right)$ does not
comply with the terminal conditions for the state.

The co-state equation becomes
\begin{align}
-\left(\begin{array}{c}
\dot{\lambda}_{x}\left(t\right)\\
\dot{\lambda}_{y}\left(t\right)
\end{array}\right) & =\left(\begin{array}{cc}
0 & 0\\
1 & 0
\end{array}\right)\left(\begin{array}{c}
\lambda_{x}\left(t\right)\\
\lambda_{y}\left(t\right)
\end{array}\right)+\left(\begin{array}{c}
x\left(t\right)\\
y\left(t\right)
\end{array}\right),
\end{align}
with terminal condition
\begin{align}
\lambda_{x}\left(t_{1}\right) & =\beta_{1}\left(x\left(t_{1}\right)-x_{1}\right), & \lambda_{y}\left(t_{1}\right) & =\beta_{2}\left(y\left(t_{1}\right)-y_{1}\right).
\end{align}
The stationarity condition is
\begin{align}
0 & =\epsilon^{2}u\left(t\right)+\lambda_{y}\left(t\right).\label{eq:ControlSolExact}
\end{align}
The state and co-state equations together with Eq. \eqref{eq:ControlSolExact}
can be cast in form of a $4\times4$ linear dynamical system 
\begin{align}
\left(\begin{array}{c}
\dot{x}\left(t\right)\\
\dot{y}\left(t\right)\\
\dot{\lambda}_{x}\left(t\right)\\
\dot{\lambda}_{y}\left(t\right)
\end{array}\right) & =\left(\begin{array}{cccc}
0 & 1 & 0 & 0\\
0 & 0 & 0 & -\epsilon^{-2}\\
-1 & 0 & 0 & 0\\
0 & -1 & -1 & 0
\end{array}\right)\left(\begin{array}{c}
x\left(t\right)\\
y\left(t\right)\\
\lambda_{x}\left(t\right)\\
\lambda_{y}\left(t\right)
\end{array}\right).
\end{align}
Assuming the regularization parameter $\epsilon$ is restricted to
the range $0<\epsilon<\frac{1}{2}$, the four eigenvalues $\sigma$
of the constant state matrix
\begin{align}
\boldsymbol{\mathcal{A}} & =\left(\begin{array}{cccc}
0 & 1 & 0 & 0\\
0 & 0 & 0 & -\epsilon^{-2}\\
-1 & 0 & 0 & 0\\
0 & -1 & -1 & 0
\end{array}\right)
\end{align}
are real and given by
\begin{align}
\sigma_{1,2,3,4} & =\pm\frac{1}{\sqrt{2}\epsilon}\sqrt{1\pm\sqrt{1-4\epsilon^{2}}}.
\end{align}
The exact solution for the state and co-state is given as a superposition
of exponentials $\sim\exp\left(\sigma_{i}t\right)$, which is conveniently
written as
\begin{align}
\left(\begin{array}{c}
x\left(t\right)\\
y\left(t\right)\\
\lambda_{x}\left(t\right)\\
\lambda_{y}\left(t\right)
\end{array}\right) & =\left(\begin{array}{c}
\mathcal{H}_{11}\\
\mathcal{H}_{21}\\
\mathcal{H}_{31}\\
\mathcal{H}_{41}
\end{array}\right)\sinh\left(\kappa_{1}t\right)+\left(\begin{array}{c}
\mathcal{H}_{12}\\
\mathcal{H}_{22}\\
\mathcal{H}_{32}\\
\mathcal{H}_{42}
\end{array}\right)\cosh\left(\kappa_{2}t\right)\nonumber \\
 & +\left(\begin{array}{c}
\mathcal{H}_{13}\\
\mathcal{H}_{23}\\
\mathcal{H}_{33}\\
\mathcal{H}_{43}
\end{array}\right)\sinh\left(\kappa_{1}t\right)+\left(\begin{array}{c}
\mathcal{H}_{14}\\
\mathcal{H}_{24}\\
\mathcal{H}_{34}\\
\mathcal{H}_{44}
\end{array}\right)\sinh\left(\kappa_{2}t\right).\label{eq:ExactSolution}
\end{align}
We introduced the abbreviations
\begin{align}
\kappa_{1} & =\frac{\sqrt{1-\sqrt{1-4\epsilon^{2}}}}{\sqrt{2}\epsilon}, & \kappa_{2} & =\frac{\sqrt{\sqrt{1-4\epsilon^{2}}+1}}{\sqrt{2}\epsilon}, & \kappa_{3} & =\epsilon^{2}\sqrt{1-4\epsilon^{2}},
\end{align}
and the $4\times4$ matrix of coefficients $\left(\boldsymbol{\mathcal{H}}\right)_{ij}=\mathcal{H}_{ij}$
given by 
\begin{align}
\boldsymbol{\mathcal{H}} & =\tiny{\left(\begin{array}{cccc}
\frac{2C_{2}\epsilon^{2}+x_{0}\left(\kappa_{3}+\epsilon^{2}\right)}{2\kappa_{3}} & \frac{1}{2}\left(x_{0}-\frac{\epsilon^{2}\left(2C_{2}+x_{0}\right)}{\kappa_{3}}\right) & \frac{\kappa_{3}y_{0}-\epsilon^{2}\left(2C_{1}+y_{0}\right)}{2\kappa_{1}\kappa_{3}} & \frac{2C_{1}\epsilon^{2}+y_{0}\left(\kappa_{3}+\epsilon^{2}\right)}{2\kappa_{2}\kappa_{3}}\\
\frac{1}{2}\left(y_{0}-\frac{\epsilon^{2}\left(2C_{1}+y_{0}\right)}{\kappa_{3}}\right) & \frac{2C_{1}\epsilon^{2}+y_{0}\left(\kappa_{3}+\epsilon^{2}\right)}{2\kappa_{3}} & \frac{C_{2}\left(\epsilon^{2}-\kappa_{3}\right)+2x_{0}\epsilon^{4}}{2\kappa_{1}\kappa_{3}\epsilon^{2}} & -\frac{C_{2}\left(\kappa_{3}+\epsilon^{2}\right)+2x_{0}\epsilon^{4}}{2\kappa_{2}\kappa_{3}\epsilon^{2}}\\
\frac{C_{1}\left(\kappa_{3}+\epsilon^{2}\right)+2y_{0}\epsilon^{4}}{2\kappa_{3}} & \frac{C_{1}\left(\kappa_{3}-\epsilon^{2}\right)-2y_{0}\epsilon^{4}}{2\kappa_{3}} & -\frac{2C_{2}\epsilon^{2}+x_{0}\left(\kappa_{3}+\epsilon^{2}\right)}{2\kappa_{1}\kappa_{3}} & \frac{2C_{2}\epsilon^{2}+x_{0}\left(\epsilon^{2}-\kappa_{3}\right)}{2\kappa_{2}\kappa_{3}}\\
\frac{C_{2}\left(\kappa_{3}-\epsilon^{2}\right)-2x_{0}\epsilon^{4}}{2\kappa_{3}} & \frac{C_{2}\left(\kappa_{3}+\epsilon^{2}\right)+2x_{0}\epsilon^{4}}{2\kappa_{3}} & \frac{-\kappa_{3}\left(C_{1}+y_{0}\right)+\epsilon^{2}\left(C_{1}+y_{0}\right)-2y_{0}\epsilon^{4}}{2\kappa_{1}\kappa_{3}} & -\frac{2C_{1}\kappa_{2}^{2}\epsilon^{4}+y_{0}\left(\kappa_{3}-2\epsilon^{4}+\epsilon^{2}\right)}{2\kappa_{2}\kappa_{3}}
\end{array}\right)}.
\end{align}
The constants $C_{1}$ and $C_{2}$ are very long expressions and
not important for the subsequent discussion.

Although the solution Eq. \eqref{eq:ExactSolution} is unhandy, it
can be studied relatively easily with the computer algebra system
Mathematica \cite{wolfram2014mathematica}. The next section compares
Eq. \eqref{eq:ExactSolution} with exact solutions to slightly different
optimization problems. These problems differ from Eq. \eqref{eq:FreeParticleFunctional}
only in the terminal conditions.

\subsection{Different terminal conditions}

The last section discussed the case of penalized terminal conditions
leading to terminal conditions for the co-state,
\begin{align}
\lambda_{x}\left(t_{1}\right) & =\beta_{1}\left(x\left(t_{1}\right)-x_{1}\right), & \lambda_{y}\left(t_{1}\right) & =\beta_{2}\left(y\left(t_{1}\right)-y_{1}\right).
\end{align}
The weighting coefficients $\beta_{1/2}\geq0$ quantify the cost for
deviating from the desired terminal state $\boldsymbol{x}_{1}$. With
increasing value of $\beta_{1/2}$, the cost is increasing. In the
limit $\beta_{1/2}\rightarrow\infty$, the terminal co-state can be
finite only for sharp terminal conditions
\begin{align}
x\left(t_{1}\right) & =x_{1}, & y\left(t_{1}\right) & =y_{1}.
\end{align}
This limit can only exist if the system is controllable. Controllability
for mechanical control systems in one spatial dimension, including
the free particle discussed here, was proven in Section \ref{sec:Controllability}.

These considerations lead to the following conjecture. Minimizing
the functional 
\begin{align}
\mathcal{J}\left[\boldsymbol{x}\left(t\right),u\left(t\right)\right] & =\frac{1}{2}\intop_{t_{0}}^{t_{1}}\left(\left(x\left(t\right)\right)^{2}+\left(y\left(t\right)\right)^{2}\right)dt+\frac{\epsilon^{2}}{2}\intop_{t_{0}}^{t_{1}}dt\left(u\left(t\right)\right)^{2}\nonumber \\
 & +\dfrac{\beta_{1}}{2}\left(x\left(t_{1}\right)-x_{1}\right)^{2}+\dfrac{\beta_{2}}{2}\left(y\left(t_{1}\right)-y_{1}\right)^{2}
\end{align}
subject to
\begin{align}
\dot{x}\left(t\right) & =y\left(t\right), & \dot{y}\left(t\right) & =u\left(t\right), & x\left(0\right) & =x_{0}, & y\left(0\right) & =y_{0},
\end{align}
and subsequently applying the limit $\beta_{1/2}\rightarrow\infty$,
is equivalent to the minimization of 
\begin{align}
\mathcal{J}\left[\boldsymbol{x}\left(t\right),u\left(t\right)\right] & =\frac{1}{2}\intop_{t_{0}}^{t_{1}}\left(\left(x\left(t\right)\right)^{2}+\left(y\left(t\right)\right)^{2}\right)dt+\frac{\epsilon^{2}}{2}\intop_{t_{0}}^{t_{1}}dt\left(u\left(t\right)\right)^{2}
\end{align}
subject to
\begin{align}
\dot{x}\left(t\right) & =y\left(t\right), & \dot{y}\left(t\right) & =u\left(t\right),\\
x\left(0\right) & =x_{0}, & y\left(0\right) & =y_{0}, & x\left(t_{1}\right) & =x_{1}, & y\left(t_{1}\right) & =y_{1}.\label{eq:SharpTerminalCond}
\end{align}
This conjecture is confirmed by computing the limit $\beta_{1/2}\rightarrow\infty$
of the exact solution, Eq. \eqref{eq:ExactSolution}, and comparing
it with the exact solution to the corresponding problem with sharp
terminal conditions Eq. \eqref{eq:SharpTerminalCond}.

Similarly, we confirm that the limit $\beta_{1/2}\rightarrow0$ of
Eq. \eqref{eq:ExactSolution} is equivalent to the optimization problem
with free terminal conditions, i.e., to minimizing the functional
\begin{align}
\mathcal{J}\left[\boldsymbol{x}\left(t\right),u\left(t\right)\right] & =\frac{1}{2}\intop_{t_{0}}^{t_{1}}\left(\left(x\left(t\right)\right)^{2}+\left(y\left(t\right)\right)^{2}\right)dt+\frac{\epsilon^{2}}{2}\intop_{t_{0}}^{t_{1}}dt\left(u\left(t\right)\right)^{2}
\end{align}
subject to
\begin{align}
\dot{x}\left(t\right) & =y\left(t\right), & \dot{y}\left(t\right) & =u\left(t\right), & x\left(0\right) & =x_{0}, & y\left(0\right) & =y_{0}.
\end{align}
In conclusion, the task of minimizing Eq. \eqref{eq:FreeParticleFunctional},
with exact solution Eq. \eqref{eq:ExactSolution}, is the most general
way to formulate the terminal condition. All other cases can be generated
by applying the appropriate limits of the weighting coefficients $\beta_{1/2}$.
This insight enables an efficient perturbative treatment of the general
nonlinear problem in Chapter \ref{chap:AnalyticalApproximationsForOptimalTrajectoryTracking}.
However, the perturbative approach uses the regularization parameter
$\epsilon$ as the small parameter, which corresponds to investigating
the limit $\epsilon\rightarrow0$. As demonstrated in the next section,
the limit $\epsilon\rightarrow0$ is not without difficulties. The
subtle question arises if applying the limit $\epsilon\rightarrow0$
commutes with the limit $\beta_{1/2}\rightarrow0$.

\subsection{Approximating the exact solution\label{sub:ApproxExactSol}}

A small regularization parameter $0<\epsilon\ll1$ is assumed to approximate
the exact solution. A regular expansion in form of a power series
in $\epsilon$ results in an approximation which is not uniformly
valid over the entire time interval. Phenomenologically, this non-uniformity
manifests in the appearance of \textit{boundary layers} at the beginning
and end of the time interval for small $\epsilon$. Uniformly valid
approximation are obtained for a series expansion in form of an asymptotic
series. This procedure is known as a singular perturbation expansion.
For a detailed account of asymptotic series, the difference between
singular and regular perturbation theory, and boundary layers see
the excellent book \cite{bender1999advanced} and also \cite{johnson2006singular}.

\subsubsection{Inner and outer limits}

The solution Eq. \eqref{eq:ExactSolution} is rather intimidating
and clumsy and not very useful for subsequent computations. Here,
the leading order approximation to Eq. \eqref{eq:ExactSolution} for
small $\epsilon$ is obtained. The perturbation expansion also involves
expansion of the eigenvalues $\sigma_{i}$. Due to the appearance
of $1/\epsilon$ in the eigenvalues $\sigma_{i}$, care has to be
taken when considering the limit $\epsilon\rightarrow0$. Note that
$\kappa_{1,2}$ behave as
\begin{align}
\kappa_{1} & =\frac{\sqrt{1-\sqrt{1-4\epsilon^{2}}}}{\sqrt{2}\epsilon}=1+\dfrac{\epsilon^{2}}{2}+\mathcal{O}\left(\epsilon^{3}\right),\\
\kappa_{2} & =\frac{\sqrt{\sqrt{1-4\epsilon^{2}}+1}}{\sqrt{2}\epsilon}=\dfrac{1}{\epsilon}-\dfrac{\epsilon}{2}+\mathcal{O}\left(\epsilon^{3}\right),
\end{align}
such that the exponential terms contained in Eq. \eqref{eq:ExactSolution}
are of the form
\begin{align}
\exp\left(-\kappa_{2}t\right) & \approx\exp\left(-\dfrac{t}{\epsilon}\right).\label{eq:ExpKappa1}\\
\exp\left(-\kappa_{1}t\right) & \approx\exp\left(-t\right),\label{eq:ExpKappa2}
\end{align}
While the limit for Eq. \eqref{eq:ExpKappa2} can safely be applied
independent of the value of $t$, the limit for Eq. \eqref{eq:ExpKappa1}
depends on the actual value of $t$. The result is
\begin{align}
\lim_{\epsilon\rightarrow0}\exp\left(-\frac{1}{\epsilon}t\right) & =\begin{cases}
0, & t>0,\\
1, & t=0,
\end{cases}
\end{align}
i.e., the exact solution \eqref{eq:ExactSolution} for $\epsilon=0$
is discontinuous and jumps at the left end of the time domain. For
small but finite $\epsilon$, such a behavior results in the appearance
boundary layers. Analytically, these can be resolved by rescaling
time $t$ appropriately with the small parameter $\epsilon$. Because
the solution contains (for small $\epsilon$) also exponential terms
of the form $\exp\left(-\left(t-t_{1}\right)/\epsilon\right)$, a
similar boundary layer is expected at the right end point $t=t_{1}$
of the time domain.

\begin{figure}[h]
\centering\includegraphics{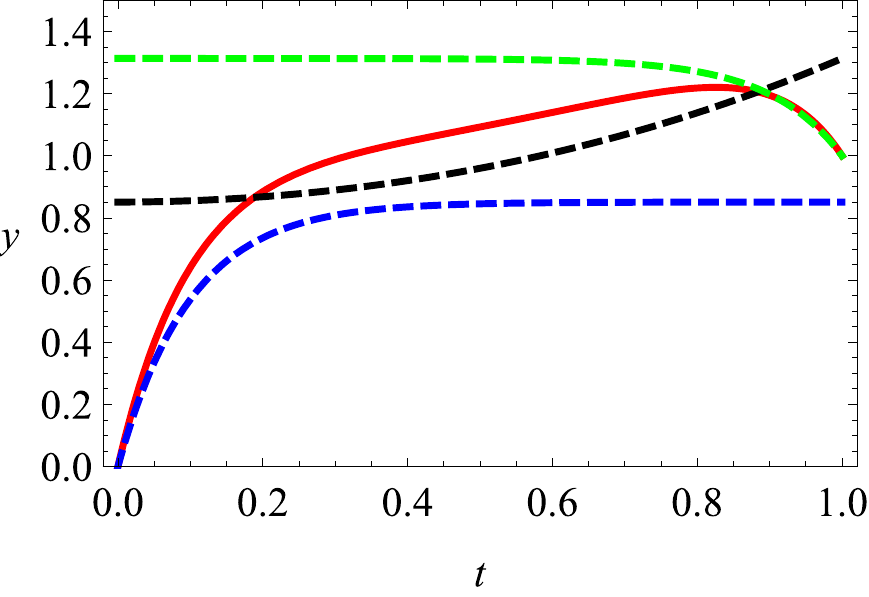}

\caption[Exact and approximate solution for the state component $y\left(t\right)$]{\label{fig:BoundaryLayers}Exact solution for the state component
$y\left(t\right)$ (red solid line) for $\epsilon=0.1$ shows the
appearance of boundary layers at the beginning and end of the time
domain. The black dashed line shows the outer approximation $y_{O}\left(t\right)$
valid in the bulk of the time domain. The blue and green dashed lines
shows that the left and right inner approximations given by $Y_{L}\left(t/\epsilon\right)$
and $Y_{R}\left(\left(t_{1}-t\right)/\epsilon\right)$ are valid close
to the initial and terminal time, respectively, and resolve the boundary
layers. A combination of all three solutions is necessary to yield
an approximation which is uniformly valid over the whole time interval
and also satisfies the initial and terminal conditions.}
\end{figure}

First, consider the limit $\epsilon\rightarrow0$ in an interior point
$0<t<t_{1}$ of the time domain. This limit is called the \textit{outer
limit} and is denoted with index $O$,
\begin{align}
\left(\begin{array}{c}
x_{O}\left(t\right)\\
y_{O}\left(t\right)\\
\lambda_{x,O}\left(t\right)\\
\lambda_{y,O}\left(t\right)
\end{array}\right) & =\lim_{\epsilon\rightarrow0}\left(\begin{array}{c}
x\left(t\right)\\
y\left(t\right)\\
\lambda_{x}\left(t\right)\\
\lambda_{y}\left(t\right)
\end{array}\right).
\end{align}
To obtain this limit by hand is very tedious due to the complexity
of the exact solution. We rely on the capabilities of the computer
algebra system Mathematica and simply state the result,
\begin{align}
\left(\begin{array}{c}
x_{O}\left(t\right)\\
y_{O}\left(t\right)\\
\lambda_{x,O}\left(t\right)\\
\lambda_{y,O}\left(t\right)
\end{array}\right) & =\left(\begin{array}{c}
\dfrac{\beta_{1}x_{1}}{\kappa}\sinh\left(t\right)+\dfrac{x_{0}}{\kappa}\left(\cosh\left(t-t_{1}\right)-\beta_{1}\sinh\left(t-t_{1}\right)\right)\\
\dfrac{\beta_{1}x_{1}}{\kappa}\cosh\left(t\right)+\dfrac{x_{0}}{\kappa}\left(\sinh\left(t-t_{1}\right)-\beta_{1}\cosh\left(t-t_{1}\right)\right)\\
\dfrac{1}{\kappa}\cosh\left(t\right)\left(x_{0}\left(\beta_{1}\cosh\left(t_{1}\right)+\sinh\left(t_{1}\right)\right)-\beta_{1}x_{1}\right)-x_{0}\sinh\left(t\right)\\
0
\end{array}\right),\label{eq:OuterLimit}
\end{align}
with the abbreviation
\begin{align}
\kappa & =\cosh\left(t_{1}\right)+\beta_{1}\sinh\left(t_{1}\right).
\end{align}
Note that this solution does not depend on the initial and terminal
points $y_{0}$ and $y_{1}$! Furthermore, the outer limit $y_{O}\left(t\right)$
does not obey the initial condition $y\left(0\right)=y_{0}$ because
\begin{align}
y_{O}\left(0\right) & =\dfrac{\beta_{1}x_{1}}{\kappa}-\dfrac{x_{0}}{\kappa}\left(\sinh\left(t_{1}\right)+\beta_{1}\cosh\left(t_{1}\right)\right)\neq y_{0}.
\end{align}
To obtain the behavior of the exact solution \eqref{eq:ExactSolution}
near to $t=t_{0}=0$, the time $t$ is rescaled and a new time scale
is introduced as
\begin{align}
\tau_{L} & =\left(t-t_{0}\right)/\epsilon=t/\epsilon.
\end{align}
The corresponding limits of Eq. \eqref{eq:ExactSolution}, valid for
times $t$ on the order of $\epsilon$, $t\sim\epsilon$, are called
\textit{left inner limits} and denoted by
\begin{align}
\left(\begin{array}{c}
X\left(\tau_{L}\right)\\
Y_{L}\left(\tau_{L}\right)\\
\Lambda_{x,L}\left(\tau_{L}\right)\\
\Lambda_{y,L}\left(\tau_{L}\right)
\end{array}\right) & =\lim_{\epsilon\rightarrow0}\left(\begin{array}{c}
x\left(\epsilon\tau_{L}\right)\\
y\left(\epsilon\tau_{L}\right)\\
\lambda_{x}\left(\epsilon\tau_{L}\right)\\
\lambda_{y}\left(\epsilon\tau_{L}\right)
\end{array}\right).
\end{align}
With the help of Mathematica, we obtain
\begin{align}
\left(\begin{array}{c}
X_{L}\left(\tau_{L}\right)\\
Y_{L}\left(\tau_{L}\right)\\
\Lambda_{x,L}\left(\tau_{L}\right)\\
\Lambda_{y,L}\left(\tau_{L}\right)
\end{array}\right) & =\left(\begin{array}{c}
x_{0}\\
\dfrac{\left(1-e^{-\tau_{L}}\right)}{\kappa}\left(\beta_{1}x_{1}-x_{0}\left(\beta_{1}\cosh\left(t_{1}\right)+\sinh\left(t_{1}\right)\right)\right)+y_{0}e^{-\tau_{L}}\\
\dfrac{\beta_{1}\text{csch}\left(t_{1}\right)\left(x_{0}\cosh\left(t_{1}\right)-x_{1}\right)+x_{0}}{\beta_{1}+\coth\left(t_{1}\right)}\\
0
\end{array}\right).\label{eq:InnerLeftLimit}
\end{align}
Note that the solution $Y_{L}\left(\tau_{L}\right)$ satisfies the
appropriate initial condition 
\begin{align}
Y_{L}\left(0\right) & =y\left(0\right)=y_{0}.
\end{align}
A similar procedure is applied for the \textit{right inner limit}
by introducing an appropriate time scale as
\begin{align}
\tau_{R} & =\left(t_{1}-t\right)/\epsilon.
\end{align}
Solutions valid for times $t$ close to the end of the time interval,
$\left(t_{1}-t\right)\sim\epsilon$, are denoted by 
\begin{align}
\left(\begin{array}{c}
X_{R}\left(\tau_{R}\right)\\
Y_{R}\left(\tau_{R}\right)\\
\Lambda_{x,R}\left(\tau_{R}\right)\\
\Lambda_{y,R}\left(\tau_{R}\right)
\end{array}\right) & =\lim_{\epsilon\rightarrow0}\left(\begin{array}{c}
x\left(t_{1}-\epsilon\tau_{R}\right)\\
y\left(t_{1}-\epsilon\tau_{R}\right)\\
\lambda_{x}\left(t_{1}-\epsilon\tau_{R}\right)\\
\lambda_{y}\left(t_{1}-\epsilon\tau_{R}\right)
\end{array}\right).
\end{align}
With the help of Mathematica, we obtain
\begin{align}
\left(\begin{array}{c}
X_{R}\left(\tau_{R}\right)\\
Y_{R}\left(\tau_{R}\right)\\
\Lambda_{x,R}\left(\tau_{R}\right)\\
\Lambda_{y,R}\left(\tau_{R}\right)
\end{array}\right) & =\left(\begin{array}{c}
\dfrac{1}{\kappa}\left(x_{0}+x_{1}\beta_{1}\sinh\left(t_{1}\right)\right)\\
y_{1}e^{-\tau_{R}}-\dfrac{\beta_{1}}{\kappa}\left(1-e^{-\tau_{R}}\right)\left(x_{0}-x_{1}\cosh\left(t_{1}\right)\right)\\
\dfrac{\beta_{1}}{\kappa}\left(x_{0}-x_{1}\cosh\left(t_{1}\right)\right)\\
0
\end{array}\right).\label{eq:InnerRightLimit}
\end{align}
Fig. \ref{fig:BoundaryLayers} shows the exact solution for the state
variable $y\left(t\right)$ (red line) together with the corresponding
outer and both inner limits (dashed lines) for a relatively large
value of $\epsilon=1/10$. While the outer limit provides an approximation
inside the domain but fails at the beginning and end of the time domain,
the left and right inner limits approximate these regions quite well.
The idea is now to combine all limits in a single composite solution
and obtain an approximation which is uniformly valid over the entire
time interval.

\subsubsection{Matching and composite solution}

For a composition of the inner and outer limits to a single and uniformly
valid approximation, certain \textit{matching conditions} must be
satisfied. If the matching conditions are violated, other scalings
of time $t$ with the small parameter $\epsilon$ exist and cannot
be neglected in the composite solution. In general, several scaling
regimes for the inner solutions are possible. A larger variety of
scalings may lead to more complicated structures such as \textit{nested
boundary layers}, \textit{interior boundary layers} or \textit{super
sharp boundary layers} \cite{bender1999advanced}.

The matching conditions at the left end of the time domain are
\begin{align}
O_{L} & =\lim_{t\rightarrow0}\left(\begin{array}{c}
x_{O}\left(t\right)\\
y_{O}\left(t\right)\\
\lambda_{x,O}\left(t\right)\\
\lambda_{y,O}\left(t\right)
\end{array}\right)=\lim_{\tau_{L}\rightarrow\infty}\left(\begin{array}{c}
X\left(\tau_{L}\right)\\
Y_{L}\left(\tau_{L}\right)\\
\Lambda_{x,L}\left(\tau_{L}\right)\\
\Lambda_{y,L}\left(\tau_{L}\right)
\end{array}\right).
\end{align}
Computing both limits yields identical results, and the matching
conditions are satisfied. The \textit{left overlap} $O_{L}$ is obtained
as 
\begin{align}
O_{L} & =\left(\begin{array}{c}
x_{0}\\
\dfrac{1}{\kappa}\left(\beta_{1}x_{1}-x_{0}\left(\beta_{1}\cosh\left(t_{1}\right)+\sinh\left(t_{1}\right)\right)\right)\\
\dfrac{\beta_{1}\text{csch}\left(t_{1}\right)\left(x_{0}\cosh\left(t_{1}\right)-x_{1}\right)+x_{0}}{\beta_{1}+\coth\left(t_{1}\right)}\\
0
\end{array}\right).
\end{align}
The matching conditions at the right end of the time domain are
\begin{align}
O_{R} & =\lim_{t\rightarrow t_{1}}\left(\begin{array}{c}
x_{O}\left(t\right)\\
y_{O}\left(t\right)\\
\lambda_{x,O}\left(t\right)\\
\lambda_{y,O}\left(t\right)
\end{array}\right)=\lim_{\tau_{R}\rightarrow\infty}\left(\begin{array}{c}
X_{R}\left(\tau_{R}\right)\\
Y_{R}\left(\tau_{R}\right)\\
\Lambda_{x,R}\left(\tau_{R}\right)\\
\Lambda_{y,R}\left(\tau_{R}\right)
\end{array}\right).
\end{align}
The \textit{right overlap} $O_{R}$ is obtained as
\begin{align}
O_{R} & =\left(\begin{array}{c}
\dfrac{1}{\kappa}\left(x_{0}+x_{1}\beta_{1}\sinh\left(t_{1}\right)\right)\\
-\dfrac{\beta_{1}}{\kappa}\left(x_{0}-x_{1}\cosh\left(t_{1}\right)\right)\\
\dfrac{\beta_{1}}{\kappa}\left(x_{0}-x_{1}\cosh\left(t_{1}\right)\right)\\
0
\end{array}\right).
\end{align}
In conclusion, both matching conditions are satisfied.

\begin{figure}[h]
\centering\subfloat{\includegraphics[scale=0.79]{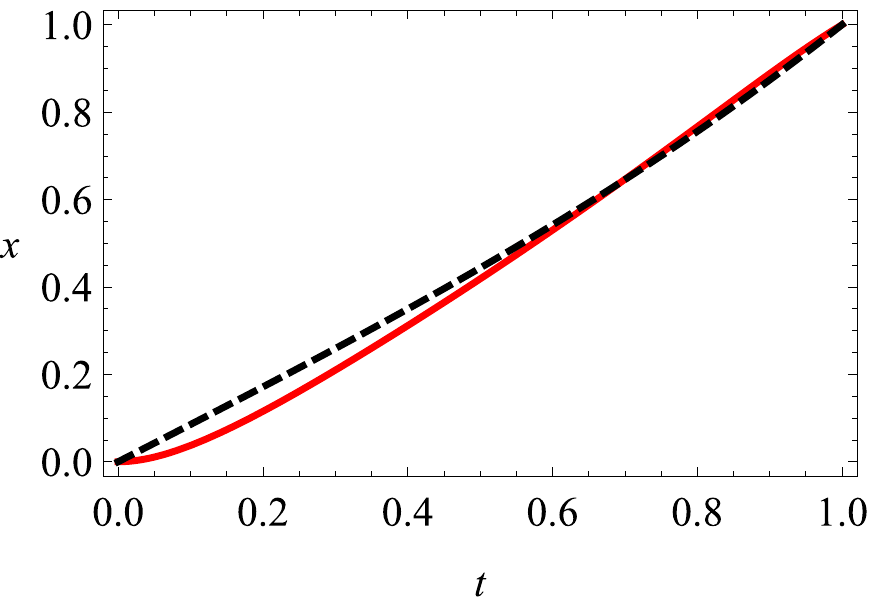}}\hspace{0.7cm}\subfloat{\includegraphics[scale=0.79]{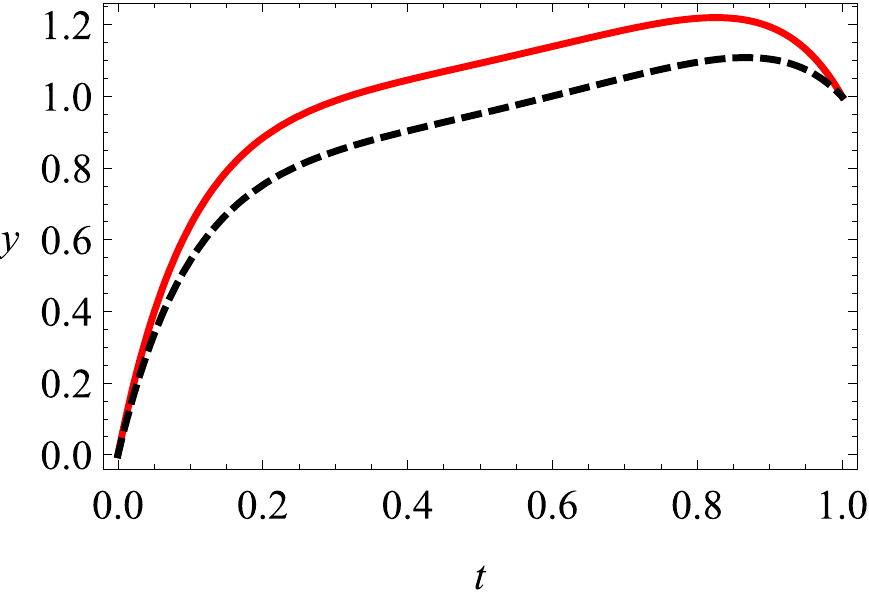}}

\subfloat{\includegraphics[scale=0.84]{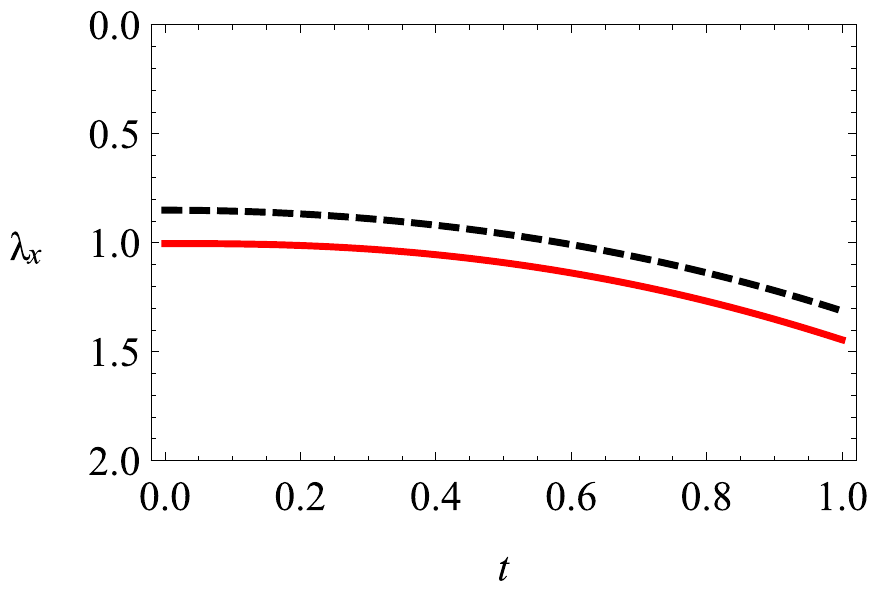}}\hspace{0.1cm}\subfloat{\includegraphics[scale=0.84]{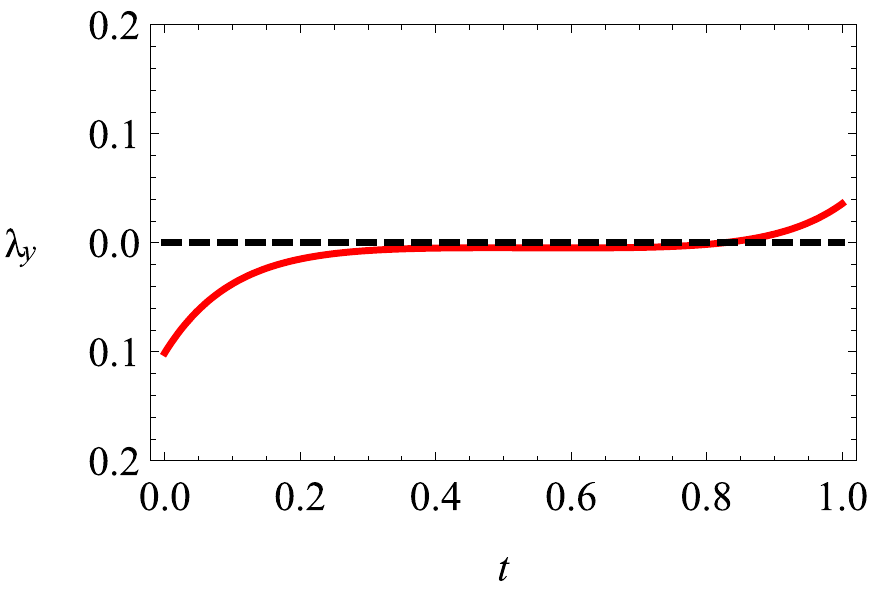}}\caption[Exact and approximate solution for $\epsilon=1/10$ for state and
co-state]{\label{fig:ExactSolution1}The exact solutions (red solid line) for
the states $x\left(t\right)$ (top left), $y\left(t\right)$ (top
right) and co-states $\lambda_{x}\left(t\right)$ (bottom left) and
$\lambda_{y}\left(t\right)$ (bottom right) for $\epsilon=0.1$. The
black dashed lines show the approximate composite solutions obtained
by a singular perturbation expansion.}
\end{figure}
Knowing all inner and outer approximations, they are combined to a
\textit{composite solution} uniformly valid on the entire time domain.
This is done by adding up all inner and outer solutions and subtracting
the overlaps \cite{bender1999advanced},
\begin{align}
\left(\begin{array}{c}
x_{\text{comp}}\left(t\right)\\
y_{\text{comp}}\left(t\right)\\
\lambda_{x,\text{comp}}\left(t\right)\\
\lambda_{y,\text{comp}}\left(t\right)
\end{array}\right) & =\left(\begin{array}{c}
x_{O}\left(t\right)\\
y_{O}\left(t\right)\\
\lambda_{x,O}\left(t\right)\\
\lambda_{y,O}\left(t\right)
\end{array}\right)+\left(\begin{array}{c}
X_{L}\left(t/\epsilon\right)\\
Y_{L}\left(t/\epsilon\right)\\
\Lambda_{x,L}\left(t/\epsilon\right)\\
\Lambda_{y,L}\left(t/\epsilon\right)
\end{array}\right)\nonumber \\
 & +\left(\begin{array}{c}
X_{R}\left(\left(t_{1}-t\right)/\epsilon\right)\\
Y_{R}\left(\left(t_{1}-t\right)/\epsilon\right)\\
\Lambda_{x,R}\left(\left(t_{1}-t\right)/\epsilon\right)\\
\Lambda_{y,R}\left(\left(t_{1}-t\right)/\epsilon\right)
\end{array}\right)-O_{L}-O_{R}.
\end{align}
Finally, the composite solution is given by
\begin{align}
y_{\text{comp}}\left(t\right) & =\dfrac{1}{\kappa}\left(\beta_{1}x_{1}\cosh\left(t\right)+x_{0}\left(\sinh\left(t-t_{1}\right)-\beta_{1}\cosh\left(t-t_{1}\right)\right)\right)\nonumber \\
 & +\dfrac{1}{\kappa}e^{-\left(t_{1}-t\right)/\epsilon}\left(\beta_{1}\left(y_{1}\sinh\left(t_{1}\right)+x_{0}\right)+\cosh\left(t_{1}\right)\left(y_{1}-\beta_{1}x_{1}\right)\right)\nonumber \\
 & +\dfrac{1}{\kappa}e^{-t/\epsilon}\left(\sinh\left(t_{1}\right)\left(x_{0}+\beta_{1}y_{0}\right)+\cosh\left(t_{1}\right)\left(\beta_{1}x_{0}+y_{0}\right)-\beta_{1}x_{1}\right),\label{eq:CompositeApproximation1}
\end{align}
and
\begin{align}
\left(\begin{array}{c}
x_{\text{comp}}\left(t\right)\\
\lambda_{x,\text{comp}}\left(t\right)\\
\lambda_{y,\text{comp}}\left(t\right)
\end{array}\right) & =\left(\begin{array}{c}
\dfrac{\beta_{1}x_{1}}{\kappa}\sinh\left(t\right)+\dfrac{x_{0}}{\kappa}\left(\cosh\left(t-t_{1}\right)-\beta_{1}\sinh\left(t-t_{1}\right)\right)\\
\dfrac{\cosh\left(t\right)}{\kappa}\left(x_{0}\left(\beta_{1}\cosh\left(t_{1}\right)+\sinh\left(t_{1}\right)\right)-\beta_{1}x_{1}\right)-x_{0}\sinh\left(t\right)\\
0
\end{array}\right).\label{eq:CompositeApproximation2}
\end{align}
This is the leading order approximation as $\epsilon\rightarrow0$
of the exact solution Eq. \eqref{eq:ExactSolution}. Note that this
approximate solution depends on the small parameter $\epsilon$ itself,
as it is generally the case for singular perturbation expansions.
Thus, the leading order composite solution is not simply given by
the limit
\begin{align}
\left(\begin{array}{c}
x_{\text{comp}}\left(t\right)\\
y_{\text{comp}}\left(t\right)\\
\lambda_{x,\text{comp}}\left(t\right)\\
\lambda_{y,\text{comp}}\left(t\right)
\end{array}\right) & \neq\lim_{\epsilon\rightarrow0}\left(\begin{array}{c}
x\left(t\right)\\
y\left(t\right)\\
\lambda_{x}\left(t\right)\\
\lambda_{y}\left(t\right)
\end{array}\right),
\end{align}
with $\left(\begin{array}{cccc}
x\left(t\right), & y\left(t\right), & \lambda_{x}\left(t\right), & \lambda_{y}\left(t\right)\end{array}\right)^{T}$ denoting the exact solution. To distinguish the operation ``obtain
the leading order contribution'' from the limit $\epsilon\rightarrow0$,
the notation
\begin{align}
\left(\begin{array}{c}
x_{\text{comp}}\left(t\right)\\
y_{\text{comp}}\left(t\right)\\
\lambda_{x,\text{comp}}\left(t\right)\\
\lambda_{y,\text{comp}}\left(t\right)
\end{array}\right) & =\left(\begin{array}{c}
x\left(t\right)\\
y\left(t\right)\\
\lambda_{x}\left(t\right)\\
\lambda_{y}\left(t\right)
\end{array}\right)+\text{h.o.t.}
\end{align}
is adopted. Here, $\text{h.o.t.}$ stands for ``higher order terms'',
and means that all higher order contributions are neglected as $\epsilon\rightarrow0$.

The composite solutions for all state and co-state variables (black
dashed lines) are compared with the corresponding exact analytical
solutions (red solid lines) in two figures. For a relatively large
value of $\epsilon=1/10$ the solutions agree at least qualitatively,
see Fig. \ref{fig:ExactSolution1}. For the smaller values of $\epsilon=1/40$
the exact solution is approximated quite well. The boundary layers
of width $\epsilon$ displayed by state component $y$ become steeper
with decreasing $\epsilon$ and degenerate to discontinuous jumps
for $\epsilon=0$.
\begin{figure}[h]
\centering\subfloat{\includegraphics[scale=0.79]{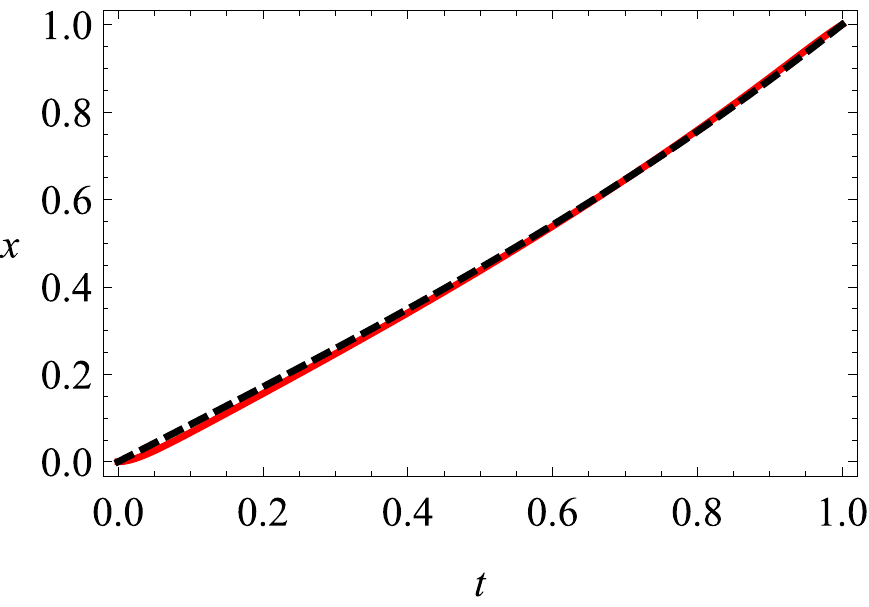}}\hspace{0.8cm}\subfloat{\includegraphics[scale=0.79]{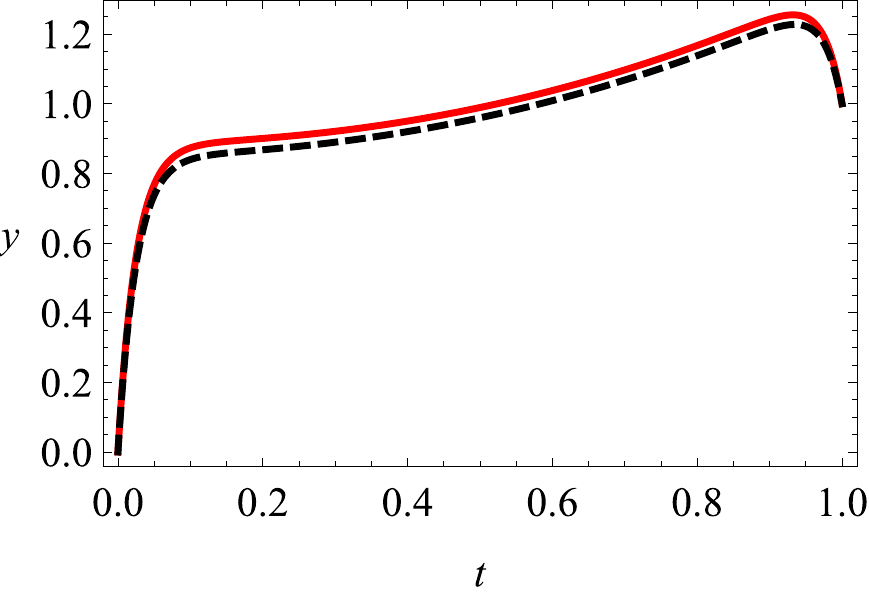}}

\subfloat{\includegraphics[scale=0.84]{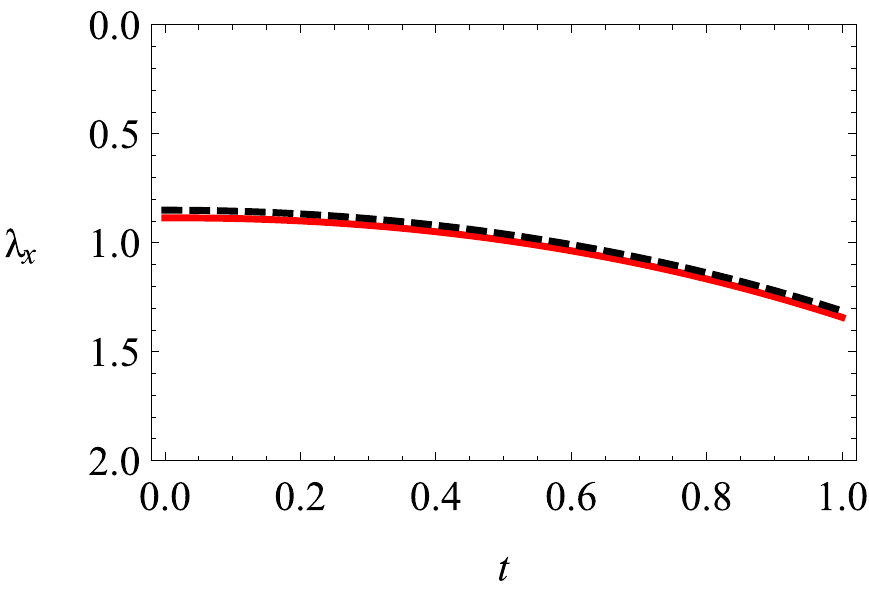}}\hspace{0.1cm}\subfloat{\includegraphics[scale=0.84]{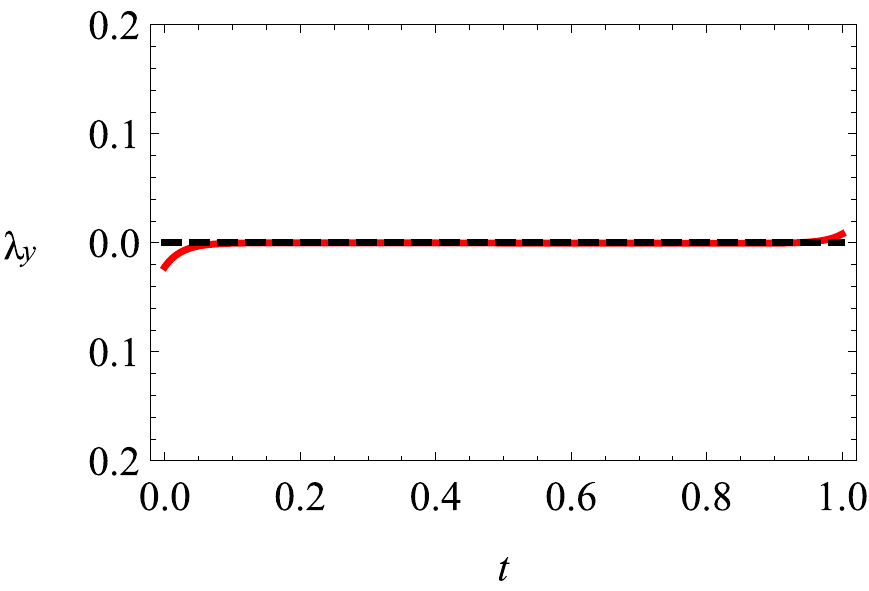}}\caption[Exact solution and approximations for $\epsilon=1/40$ for state and
co-state]{\label{fig:ExactSolution2}Same as in Fig. \ref{fig:ExactSolution1}
but for a smaller value of $\epsilon=1/40=0.025$. The agreement between
exact (red solid line) and approximate solution (black dashed line)
becomes better for smaller values of $\epsilon$ while the slopes
of the boundary layers exhibited by the state component $y\left(t\right)$
steepen.}
\end{figure}

\subsubsection{Different end point conditions}

Interestingly, the composite solution \eqref{eq:CompositeApproximation1}
and \eqref{eq:CompositeApproximation2} does not depend on the weighting
coefficient $\beta_{2}$. This can be explained as follows. Note that
the co-state $\lambda_{y}$ has to satisfy the terminal condition
\begin{align}
\lambda_{y}\left(t_{1}\right) & =\beta_{2}\left(y\left(t_{1}\right)-y_{1}\right),\label{eq:Eq471}
\end{align}
as well as the stationarity condition
\begin{align}
0 & =\epsilon^{2}u\left(t\right)+\lambda_{y}\left(t\right).
\end{align}
Clearly, a finite control signal $u\left(t\right)\neq0$ together
with $\epsilon=0$ implies $\lambda_{y}\left(t\right)\equiv0$. Thus,
the terminal condition Eq. \eqref{eq:Eq471} can only be satisfied
if 
\begin{align}
y\left(t_{1}\right) & =y_{1}.
\end{align}
In the limit $\epsilon\rightarrow0$, the state variable $y$ satisfies
the sharp terminal condition. Consequently, the parameter $\beta_{2}$
drops out of the equations, rendering state, co-state, and control
signal independent of $\beta_{2}$.

A different question is how the composite solutions Eqs. \eqref{eq:CompositeApproximation1}
and \eqref{eq:CompositeApproximation2} behave for $\beta_{1}\rightarrow\infty$
and $\beta_{1}=0$. Due to the subtle nature of the leading order
approximation as $\epsilon\rightarrow0$, it is not immediately clear
if computing the leading order commutes with the limits $\beta_{1}\rightarrow\infty$
and $\beta_{1}\rightarrow0$. Comparison of the leading order approximations
of the exact solutions for $\beta_{1}\rightarrow\infty$ and $\beta_{1}=0$
with the limits $\beta_{1}\rightarrow\infty$ and $\beta_{1}\rightarrow0$
of the leading order composite solution \eqref{eq:CompositeApproximation1}
and \eqref{eq:CompositeApproximation2} confirms that these operations
commute. In other words, computing first
\begin{align}
\left(\begin{array}{c}
x^{\infty}\left(t\right)\\
y^{\infty}\left(t\right)\\
\lambda_{x}^{\infty}\left(t\right)\\
\lambda_{y}^{\infty}\left(t\right)
\end{array}\right) & =\lim_{\beta_{1}\rightarrow\infty}\left(\begin{array}{c}
x\left(t\right)\\
y\left(t\right)\\
\lambda_{x}\left(t\right)\\
\lambda_{y}\left(t\right)
\end{array}\right)
\end{align}
followed by a leading order approximation as $\epsilon\rightarrow0$
leads to the same result as computing the limit $\beta_{1}\rightarrow\infty$
of the composite solution \eqref{eq:CompositeApproximation1} and
\eqref{eq:CompositeApproximation2}. This identity is expressed as
\begin{align}
\lim_{\beta_{1}\rightarrow\infty}\left(\begin{array}{c}
x_{\text{comp}}\left(t\right)\\
\lambda_{x,\text{comp}}\left(t\right)\\
\lambda_{y,\text{comp}}\left(t\right)
\end{array}\right) & =\left(\begin{array}{c}
x^{\infty}\left(t\right)\\
y^{\infty}\left(t\right)\\
\lambda_{x}^{\infty}\left(t\right)\\
\lambda_{y}^{\infty}\left(t\right)
\end{array}\right)+\text{h.o.t.}.
\end{align}
An analogous result is valid for the limit $\beta_{1}\rightarrow0$.
Let the exact result for $\beta_{1}=0$ be defined by 
\begin{align}
\left(\begin{array}{c}
x^{0}\left(t\right)\\
y^{0}\left(t\right)\\
\lambda_{x}^{0}\left(t\right)\\
\lambda_{y}^{0}\left(t\right)
\end{array}\right) & =\lim_{\beta_{1}\rightarrow0}\left(\begin{array}{c}
x\left(t\right)\\
y\left(t\right)\\
\lambda_{x}\left(t\right)\\
\lambda_{y}\left(t\right)
\end{array}\right).
\end{align}
Then the following identity holds
\begin{align}
\lim_{\beta_{1}\rightarrow0}\left(\begin{array}{c}
x_{\text{comp}}\left(t\right)\\
\lambda_{x,\text{comp}}\left(t\right)\\
\lambda_{y,\text{comp}}\left(t\right)
\end{array}\right) & =\left(\begin{array}{c}
x^{0}\left(t\right)\\
y^{0}\left(t\right)\\
\lambda_{x}^{0}\left(t\right)\\
\lambda_{y}^{0}\left(t\right)
\end{array}\right)+\text{h.o.t.}.
\end{align}
In summary, for the leading order approximation as $\epsilon\rightarrow0$,
it is sufficient to study the problem 
\begin{align}
\mathcal{J}\left[\boldsymbol{x}\left(t\right),u\left(t\right)\right] & =\frac{1}{2}\intop_{t_{0}}^{t_{1}}\left(\left(x\left(t\right)\right)^{2}+\left(y\left(t\right)\right)^{2}\right)dt+\frac{\epsilon^{2}}{2}\intop_{t_{0}}^{t_{1}}dt\left(u\left(t\right)\right)^{2}+\dfrac{\beta_{1}}{2}\left(x\left(t_{1}\right)-x_{1}\right)^{2},\label{eq:Eq3167}
\end{align}
subject to 
\begin{align}
\dot{x}\left(t\right) & =y\left(t\right), & \dot{y}\left(t\right) & =u\left(t\right), & x\left(0\right) & =x_{0}, & y\left(0\right) & =y_{0}, & y\left(t_{1}\right) & =y_{1}.\label{eq:Eq3168}
\end{align}
All other leading order approximations for sharp or free terminal
conditions can be generated from the solution to Eqs. \eqref{eq:Eq3167}
and \eqref{eq:Eq3168} by applying the limits $\beta_{1}\rightarrow\infty$
and $\beta_{1}=0$.

\subsubsection{Solution for the control signal}

To obtain simpler expressions, the limit $\beta_{1}\rightarrow\infty$
together with $x_{0}=y_{0}=0$ is assumed in the remainder of this
section.

\begin{figure}[h]
\centering\subfloat{\includegraphics[scale=0.82]{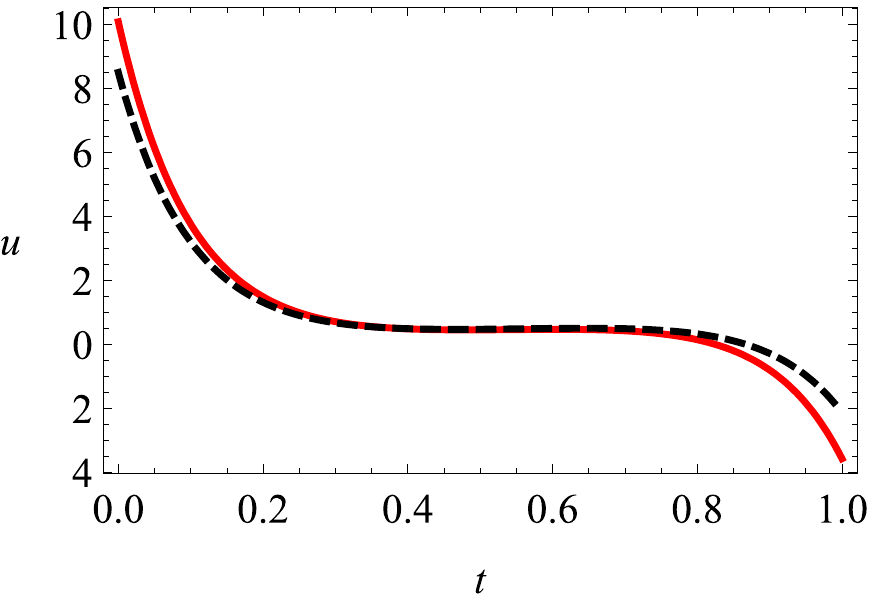}}\hspace{0.2cm}\subfloat{\includegraphics[scale=0.82]{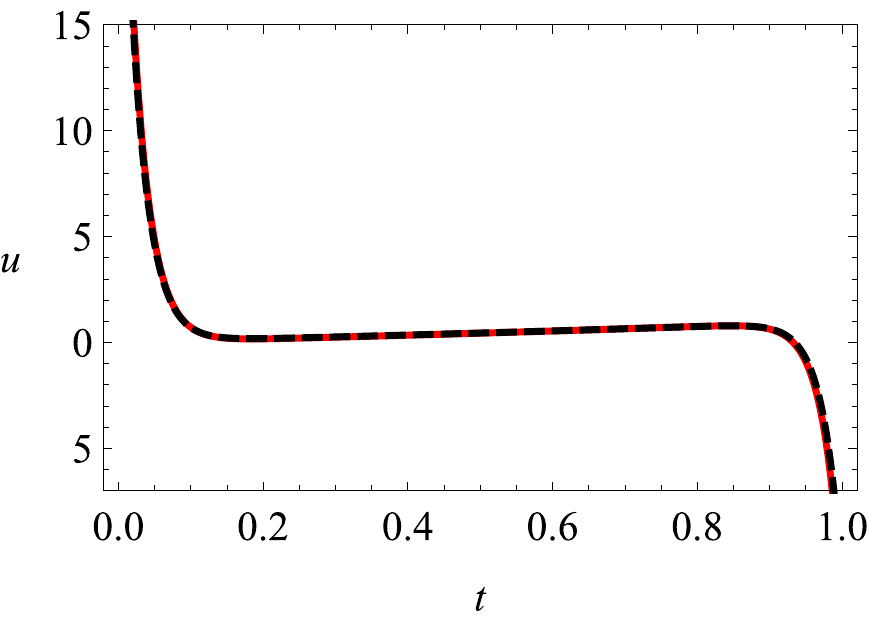}}\caption[Exact solution and approximations for the control signal]{\label{fig:CompareControls}Comparison of the exact control signal
(red solid line) and its approximation (black dashed line) for two
different values of the small parameter $\epsilon$, $\epsilon=1/10$
(left) and $\epsilon=1/40$ (right). For $\epsilon$ approaching zero,
the steeply rising boundary layers will approach infinity, leading
to delta-like kicks located at the beginning and end of the time domain.}
\end{figure}

It remains to find an approximation for the control signal $u\left(t\right)$
as follows. One way is to analyze the exact solution for the co-state
$\lambda_{y}\left(t\right)$ together with Eq. \eqref{eq:ControlSolExact}
and perform the singular perturbation expansion to find the leading
order approximation for $u\left(t\right)$. Instead, here the control
signal is derived from the composite solution for state and co-state.
Note that Eq. \eqref{eq:ControlSolExact} is satisfied to the lowest
order in $\epsilon$ also for the composite solution because $\lambda_{y,\text{comp}}\left(t\right)=0$
vanishes to leading order in $\epsilon$. However, to compute $u\left(t\right)$
from Eq. \eqref{eq:ControlSolExact} requires the knowledge of higher
order contributions to $\lambda_{y,\text{comp}}\left(t\right)$. As
an alternative, Eq. \eqref{eq:StateExact} is utilized to obtain 
\begin{align}
u_{\text{comp}}\left(t\right) & =\dot{y}_{\text{comp}}\left(t\right)\nonumber \\
 & =x_{1}\text{csch}\left(t_{1}\right)\sinh\left(t\right)+x_{1}\text{csch}\left(t_{1}\right)\frac{e^{-\frac{t}{\epsilon}}}{\epsilon}+\frac{e^{-\frac{t_{1}-t}{\epsilon}}}{\epsilon}\left(y_{1}-x_{1}\coth\left(t_{1}\right)\right).\label{eq:UCompSol}
\end{align}
Similar to the composite solution for the state $y\left(t\right)$,
the approximate control signal also contains two boundary layers,
one at $t=0$ and a second at $t=t_{1}$. The outer limit of Eq. \eqref{eq:UCompSol},
valid for times $t_{0}<t<t_{1}$ is
\begin{align}
u_{O}\left(t\right) & =\lim_{\epsilon\rightarrow0}u_{\text{comp}}\left(t\right)=x_{1}\text{csch}\left(t_{1}\right)\sinh\left(t\right).
\end{align}
For the left and right inner limits, rescaled time scales $\tau_{L}=t/\epsilon$
and $\tau_{R}=\left(t_{1}-t\right)/\epsilon$ are introduced. The
corresponding leading order solutions are denoted by $U_{L}\left(\tau_{L}\right)$
and $U_{R}\left(\tau_{R}\right)$, respectively. Because of the factor
$1/\epsilon$ in Eq. \eqref{eq:UCompSol}, it is impossible to apply
the limit $\epsilon\rightarrow0$. Instead, the leading order contribution
is computed as
\begin{align}
u\left(t\right) & =u\left(t_{0}+\epsilon\tau_{L}\right)=u_{\text{comp}}\left(t_{0}+\epsilon\tau_{L}\right)+\text{h.o.t.}\nonumber \\
 & =U_{L}\left(\tau_{L}\right)+\text{h.o.t.},\\
u\left(t\right) & =u\left(t_{1}-\epsilon\tau_{R}\right)=u_{\text{comp}}\left(t_{1}-\epsilon\tau_{R}\right)+\text{h.o.t.}\nonumber \\
 & =U_{R}\left(\tau_{R}\right)+\text{h.o.t.}.
\end{align}
All contributions of higher order in $\epsilon$ are neglected. This
procedure yields 
\begin{align}
U_{L}\left(\tau_{L}\right) & =x_{1}\text{csch}\left(t_{1}\right)\frac{e^{-\tau_{L}}}{\epsilon},\\
U_{R}\left(\tau_{R}\right) & =x_{1}+\frac{e^{-\tau_{R}}}{\epsilon}\left(y_{1}-x_{1}\coth\left(t_{1}\right)\right).
\end{align}
Finally, the matching conditions 
\begin{align}
u_{O}\left(0\right) & =0=\lim_{\tau_{L}\rightarrow\infty}U_{L}\left(\tau_{L}\right),\\
u_{O}\left(t_{1}\right) & =x_{1}=\lim_{\tau_{R}\rightarrow\infty}U_{R}\left(\tau_{R}\right),
\end{align}
are satisfied, confirming that three scaling regions are sufficient
for the composite solution. Hence, Eq. \eqref{eq:UCompSol} is the
leading order approximation to the exact solution for the control
signal. Equation \eqref{eq:UCompSol} (black dashed line) is compared
with the exact solution (red solid line) for two different values
of $\epsilon$ in Fig. \ref{fig:CompareControls}.

\subsubsection{Exact solution for \texorpdfstring{$\epsilon = 0$}{epsilon = 0}}

Having obtained approximate solutions for small but finite $\epsilon>0$,
it becomes possible to understand what exactly happens for $\epsilon=0$.
The composite solution for the state component $y\left(t\right)$
for $\beta_{1}\rightarrow\infty$ and initial conditions $x_{0}=y_{0}=0$
is given by
\begin{align}
y_{\text{comp}}\left(t\right) & =x_{1}\text{csch}\left(t_{1}\right)\cosh\left(t\right)-x_{1}\text{csch}\left(t_{1}\right)e^{-t/\epsilon}\nonumber \\
 & +e^{-\left(t_{1}-t\right)/\epsilon}\left(y_{1}-x_{1}\coth\left(t_{1}\right)\right).
\end{align}
The exact solution for $\epsilon=0$ is obtained by computing the
limit
\begin{align}
\lim_{\epsilon\rightarrow0}y\left(t\right) & =\lim_{\epsilon\rightarrow0}y_{\text{comp}}\left(t\right)=\begin{cases}
0, & t=0,\\
x_{1}\text{csch}\left(t_{1}\right)\cosh\left(t\right), & 0<t<t_{1},\\
y_{1}, & t=t_{1}.
\end{cases}\label{eq:Eq3172}
\end{align}
Using the Kronecker delta defined as
\begin{align}
\delta_{a,b} & =\begin{cases}
1, & a=b,\\
0, & a\neq b,
\end{cases}
\end{align}
Eq. \eqref{eq:Eq3172} can be written in the form
\begin{align}
\lim_{\epsilon\rightarrow0}y\left(t\right) & =x_{1}\text{csch}\left(t_{1}\right)\cosh\left(t\right)-x_{1}\text{csch}\left(t_{1}\right)\delta_{t,0}+\left(y_{1}-x_{1}\coth\left(t_{1}\right)\right)\delta_{t,t_{1}}.
\end{align}
The control signal $u\left(t\right)$ for $\epsilon=0$ is 
\begin{align}
\lim_{\epsilon\rightarrow0}u\left(t\right) & =\lim_{\epsilon\rightarrow0}u_{\text{comp}}\left(t\right)=\begin{cases}
x_{1}\text{csch}\left(t_{1}\right)\sinh\left(t\right), & 0<t<t_{1},\\
\infty, & t=0,\\
\infty, & t=t_{1},
\end{cases}
\end{align}
with $u_{\text{comp}}\left(t\right)$ given by Eq. \eqref{eq:UCompSol}.
To write that in a more enlightening form, the Dirac delta function
$\delta\left(t\right)$ is introduced. The Dirac delta function is
defined by its properties
\begin{align}
\delta\left(t\right) & =\begin{cases}
\infty, & t=0,\\
0, & t\neq0,
\end{cases}
\end{align}
and
\begin{align}
\intop_{-\infty}^{\infty}\delta\left(t\right)dt & =1.
\end{align}
The function 
\begin{align}
g_{\epsilon}\left(t\right) & =\frac{e^{-\frac{\left|t\right|}{\epsilon}}}{2\epsilon}
\end{align}
is a representation of the Dirac delta function in the limit $\epsilon\rightarrow0$,
i.e.,
\begin{align}
\lim_{\epsilon\rightarrow0}g_{\epsilon}\left(t\right) & =\delta\left(t\right).
\end{align}
This can be seen as follows. First, computing the integral yields
\begin{align}
\intop_{-\infty}^{\infty}dtg_{\epsilon}\left(t\right) & =1,
\end{align}
independent of the value of $\epsilon$. Second, computing the limit
gives 
\begin{align}
\lim_{\epsilon\rightarrow0}g_{\epsilon}\left(t\right) & =\begin{cases}
\lim_{\epsilon\rightarrow0}\frac{1}{2\epsilon}=\infty, & t=0,\\
\lim_{\epsilon\rightarrow0}\frac{e^{-\frac{\left|t\right|}{\epsilon}}}{2\epsilon}=0, & t\neq0.
\end{cases}
\end{align}
Noting that $0\leq t\leq t_{1}$ and therefore $g_{\epsilon}\left(t\right)=\frac{e^{-\frac{\left|t\right|}{\epsilon}}}{2\epsilon}=\frac{e^{-\frac{t}{\epsilon}}}{2\epsilon}$,
the control signal can be written in terms of the Dirac delta functions
as
\begin{align}
\lim_{\epsilon\rightarrow0}u_{\text{comp}}\left(t\right) & =x_{1}\text{csch}\left(t_{1}\right)\sinh\left(t\right)+2x_{1}\text{csch}\left(t_{1}\right)\delta\left(t\right)\nonumber \\
 & +2\left(y_{1}-x_{1}\coth\left(t_{1}\right)\right)\delta\left(t-t_{1}\right).
\end{align}
The interpretation is as follows. For $\epsilon=0$, the boundary
layers of the composite solution $y_{\text{comp}}\left(t\right)$
degenerate to jumps located exactly at the beginning and end of the
time interval. At $t=0$, the jump leads from the initial condition
$y_{0}=0$ to the initial value $y_{O}\left(0\right)=x_{1}\text{csch}\left(t_{1}\right)$
of the outer solution. Similarly, at $t=t_{1}$, the jump leads from
the terminal value of the outer solution $y_{O}\left(t_{1}\right)=x_{1}\coth\left(t_{1}\right)$
to the terminal condition $y_{1}$. Correspondingly, at the initial
time $t=0$, the control signal is a delta-like impulse which kicks
the state component $y$ from its initial value $y_{0}$ to $y_{O}\left(0\right)$.
The strength of the kick, given by the coefficient $2x_{1}\text{csch}\left(t_{1}\right)$
of the Dirac delta function, is two times the jump height of $y\left(0\right)$.
Intuitively, the reason is that the delta kick is located right at
the time domain boundary, and regarding the Dirac delta function as
a symmetric function, only half of the kick contributes to the time
evolution. Thus, the strength of the kick must be twice as large.
Analogously, at the terminal time, a reverse kick occurs with strength
$2\left(y_{1}-x_{1}\coth\left(t_{1}\right)\right)$, which is twice
the height of the jump between $y_{O}\left(t_{1}\right)$ and the
terminal state $y_{1}$. This picture remains valid for more complicated
unregularized optimal control problems.

\section{Conclusions}

Unregularized optimal trajectory tracking, defined by the target functional
Eq. \eqref{eq:OptimalTrajectoryTrackingFunctional} with $\epsilon=0$,
is of special interest. Its solution for the controlled state trajectory
$\boldsymbol{x}\left(t\right)$ can be seen as the limit of realizability
for a chosen desired trajectory $\boldsymbol{x}_{d}\left(t\right)$.
No other control, be it open or closed loop control, can enforce a
state trajectory $\boldsymbol{x}\left(t\right)$ with a smaller distance
to the desired state trajectory $\boldsymbol{x}_{d}\left(t\right)$.
Unregularized trajectory tracking results in singular optimal control.
Besides the usual necessary optimality conditions, additional necessary
optimality conditions, called the Kelly- or generalized Legendre-Clebsch
conditions, have to be satisfied in this case.

The concept of an exactly realizable trajectory, introduced in Chapter
\ref{chap:ExactlyRealizableTrajectories}, is closely related to an
unregularized optimal control problem. In fact, if the desired trajectory
$\boldsymbol{x}_{d}\left(t\right)$ complies with the initial and
terminal conditions
\begin{align}
\boldsymbol{x}_{d}\left(t_{0}\right) & =\boldsymbol{x}_{0}, & \boldsymbol{x}_{d}\left(t_{1}\right) & =\boldsymbol{x}_{1},
\end{align}
and satisfies the constraint equation
\begin{align}
\boldsymbol{\mathcal{Q}}\left(\boldsymbol{x}_{d}\left(t\right)\right)\left(\boldsymbol{\dot{x}}_{d}\left(t\right)-\boldsymbol{R}\left(\boldsymbol{x}_{d}\left(t\right)\right)\right) & =\boldsymbol{0},\label{eq:Eq490}
\end{align}
then the unregularized optimal control problem of minimizing 
\begin{align}
\mathcal{J}\left[\boldsymbol{x}\left(t\right),\boldsymbol{u}\left(t\right)\right]= & \frac{1}{2}\intop_{t_{0}}^{t_{1}}dt\left(\boldsymbol{x}\left(t\right)-\boldsymbol{x}_{d}\left(t\right)\right)^{T}\boldsymbol{\mathcal{S}}\left(\boldsymbol{x}\left(t\right)-\boldsymbol{x}_{d}\left(t\right)\right)\nonumber \\
 & +\frac{1}{2}\left(\boldsymbol{x}\left(t_{1}\right)-\boldsymbol{x}_{1}\right)\boldsymbol{\mathcal{S}}_{1}\left(\boldsymbol{x}\left(t_{1}\right)-\boldsymbol{x}_{1}\right)
\end{align}
subject to
\begin{align}
\boldsymbol{\dot{x}}\left(t\right) & =\boldsymbol{R}\left(\boldsymbol{x}\left(t\right)\right)+\boldsymbol{\mathcal{B}}\left(\boldsymbol{x}\left(t\right)\right)\boldsymbol{u}\left(t\right), & \boldsymbol{x}\left(t_{0}\right) & =\boldsymbol{x}_{0}
\end{align}
is solved by
\begin{align}
\boldsymbol{x}\left(t\right) & =\boldsymbol{x}_{d}\left(t\right), & \boldsymbol{u}\left(t\right) & =\boldsymbol{\mathcal{B}}^{+}\left(\boldsymbol{x}_{d}\left(t\right)\right)\left(\boldsymbol{\dot{x}}_{d}\left(t\right)-\boldsymbol{R}\left(\boldsymbol{x}_{d}\left(t\right)\right)\right).
\end{align}
The corresponding co-state $\boldsymbol{\lambda}\left(t\right)$ vanishes
for all times, $\boldsymbol{\lambda}\left(t\right)\equiv\boldsymbol{0}$,
and the functional $\mathcal{J}$ attains its minimally possible value
\begin{align}
\mathcal{J}\left[\boldsymbol{x}\left(t\right),\boldsymbol{u}\left(t\right)\right] & =0.
\end{align}
Furthermore, the controlled state trajectory $\boldsymbol{x}\left(t\right)$
and the control signal $\boldsymbol{u}\left(t\right)$ are independent
of the symmetric matrices of weighting coefficients $\boldsymbol{\mathcal{S}}$
and $\boldsymbol{\mathcal{S}}_{1}$.

If the linearizing assumption
\begin{align}
\boldsymbol{\mathcal{Q}}\boldsymbol{R}\left(\boldsymbol{x}\right) & =\boldsymbol{\mathcal{Q}}\boldsymbol{\mathcal{A}}\boldsymbol{x}+\boldsymbol{\mathcal{Q}}\boldsymbol{b}
\end{align}
holds, the linear constraint equation \eqref{eq:Eq490} is readily
solved. This results in an exact solution to a nonlinear optimization
problem in terms of linear equations. A linear structure underlying
nonlinear unregularized optimal trajectory tracking is uncovered.
Here, this linear structure is restricted to exactly realizable desired
trajectories. This restriction is dropped in the next chapter, which
demonstrates that the same underlying linear structure can be exploited
for arbitrary desired trajectories and small regularization parameter
$\epsilon$.

The exactly solvable linear example in Section \ref{sec:AnExactlySolvableExample}
highlights the difficulties encountered for small $\epsilon$ for
a desired trajectory which is not exactly realizable. Applying the
limit $\epsilon\rightarrow0$ to the exact solution yields an approximation,
called the outer solution, which is not uniformly valid over the entire
time interval $t_{0}\leq t\leq t_{1}$. Phenomenologically, this non-uniformity
manifests in steep transition regions of width $\epsilon$ displayed
by the state component $y$ close to the initial $t\gtrsim t_{0}$
and terminal $t\lesssim t_{1}$ time. These transition regions are
known as boundary layers. In analytical approximation, they are described
by the inner solutions obtained by rescaling time with $\epsilon$
in the exact solution and subsequently applying the limit $\epsilon\rightarrow0$.
The matching conditions relate inner and outer solutions in form of
overlaps. All inner and outer solutions together with their overlaps
are additively combined in a composite solution resulting in a uniformly
valid approximation. This procedure is known as a singular perturbation
expansion \cite{bender1999advanced}.

The analytical approximations obtained by singular perturbation expansion
reveal a dramatically different behavior of the solution for $\epsilon>0$
and $\epsilon=0$. While all state and co-state components are continuous
for $\epsilon>0$, the state component $y$ becomes discontinuous
for $\epsilon=0$. The smooth boundary layers displayed by $y$ degenerate
to jumps situated at the time domain boundaries. Even worse, the control
signal, being given by $u\left(t\right)=\dot{y}\left(t\right)$, diverges
at exactly these instants at which $y\left(t\right)$ becomes discontinuous.
Analytically, these divergences attain the form of Dirac delta functions
located at the beginning and end of the time interval. The strength
of the delta kicks is twice the height of the corresponding jumps
of $y\left(t\right)$.

Regarding the behavior of the exact solution with respect to different
terminal conditions, the target functional 
\begin{align}
\mathcal{J}\left[\boldsymbol{x}\left(t\right),u\left(t\right)\right] & =\frac{1}{2}\intop_{t_{0}}^{t_{1}}\left(\left(x\left(t\right)\right)^{2}+\left(y\left(t\right)\right)^{2}\right)dt\nonumber \\
 & +\frac{\epsilon^{2}}{2}\intop_{t_{0}}^{t_{1}}dt\left(u\left(t\right)\right)^{2}+\dfrac{\beta_{1}}{2}\left(x\left(t_{1}\right)-x_{1}\right)^{2}\label{eq:Eq3195}
\end{align}
subject to the dynamics constraints 
\begin{align}
\dot{x}\left(t\right) & =y\left(t\right), & \dot{y}\left(t\right) & =u\left(t\right), & x\left(0\right) & =x_{0}, & y\left(0\right) & =y_{0}, & y\left(t_{1}\right) & =y_{1},
\end{align}
constitutes the most general one in the limit of small $\epsilon$.
All variants of sharp and free terminal conditions can be generated
from the solution to Eq. \eqref{eq:Eq3195} by applying the limits
$\beta_{1}\rightarrow\infty$ and $\beta_{1}\rightarrow0$. A terminal
term of the form $\beta_{2}\left(y\left(t_{1}\right)-y_{1}\right)^{2}$
in Eq. \eqref{eq:Eq3195} becomes irrelevant as $\epsilon\rightarrow0$,
and the solution becomes independent of $\beta_{2}$. The limits $\beta_{1}\rightarrow\infty$
and $\beta_{1}\rightarrow0$ commute with the determination of the
leading order approximation for small $\epsilon$.

%% file: chapter-4.tex
\lhead[\chaptername~\thechapter\leftmark]{}

\rhead[]{\rightmark}

\lfoot[\thepage]{}

\cfoot{}

\rfoot[]{\thepage}

\chapter{\label{chap:AnalyticalApproximationsForOptimalTrajectoryTracking}Analytical
approximations for optimal trajectory tracking}

Chapters \ref{chap:ExactlyRealizableTrajectories} and \ref{chap:OptimalControl}
uncovered an underlying linear structure of unregularized nonlinear
optimal trajectory tracking for exactly realizable desired trajectories
$\boldsymbol{x}_{d}\left(t\right)$. Clearly, the assumption of $\boldsymbol{x}_{d}\left(t\right)$
to be exactly realizable is overly restrictive. Only $p$ components
of $\boldsymbol{x}_{d}\left(t\right)$ can be prescribed by the experimenter,
while the remaining $n-p$ components are fixed by the constraint
equation. This approach seriously limits the true power of optimal
control, which guarantees the existence of solutions for a huge class
of desired trajectories.

This chapter drops the assumption of exactly realizable trajectories
and allows for arbitrary desired trajectories $\boldsymbol{x}_{d}\left(t\right)$.
The regularization parameter $\epsilon$ is assumed to be small, $\epsilon\ll1$,
and used for a perturbation expansion. Rearranging the necessary optimality
condition leads to a reinterpretation of unregularized optimal control
problems as singularly perturbed differential equations. For systems
satisfying a linearizing assumption, the leading order equations become
linear. The linearity allows the derivation of closed form expressions
for optimal trajectory tracking in a general class of nonlinear systems
affine in control. The perturbative approach yields exact results
for $\epsilon=0$. However, this exact result comes at a price. The
limit $\epsilon\rightarrow0$ leads to a diverging control signal
and a discontinuous state trajectory.

The analytical approach applies to mechanical control systems defined
in Example \ref{ex:OneDimMechSys1} as well as to the activator-controlled
FHN model of Example \ref{ex:FHN1}. Section \ref{sec:TwoDimensionalDynamicalSystem}
presents the straightforward but tedious derivation for a general
model comprising both examples. A comparison with numerical solutions
of optimal control is performed in Section \ref{sec:ComparisonWithNumericalResults}.
Analytical solutions to optimal control problems also apply to optimal
feedback control. Continuous-time and time-delayed feedback are discussed
in Section \ref{sec:OptimalFeedback}. Section \ref{sec:GeneralDynamicalSystem}
tackles the singular perturbation expansion of general dynamical systems,
and Section \ref{sec:4Conclusions} draws conclusions.

\section{\label{sec:TwoDimensionalDynamicalSystem}Two-dimensional dynamical
systems}

\subsection{General procedure}

The analytic approach to optimal trajectory tracking is a straightforward
application of singular perturbation theory. The first step consists
in rearranging the necessary optimality conditions in Section \ref{sub:RearrangeNecessaryOptimalityConditions}.
This rearrangement allows the interpretation of a singular optimal
control problem as a singularly perturbed system of differential equations.
The small regularization parameter $\epsilon$ multiplies the highest
order derivative of the system. Setting $\epsilon=0$ changes the
differential order of the system and results in a violation of initial
conditions. The equations thus obtained are called the \textit{outer
equations}. Their solution is discussed in Section \ref{sub:OuterEquations}.
The outer solutions are not uniformly valid over the whole time interval.
The situation is analogous to the exactly solvable example from Section
\ref{sec:AnExactlySolvableExample}. Similar measures are taken to
resolve the problem. The time is rescaled by $\epsilon$, and a different
set of equations called the \textit{inner equations} is derived. Their
solutions are able to accommodate all initial and terminal conditions.
The inner equations are valid close to the initial and terminal conditions.
Eventually, both sets of solutions have to be connected by the \textit{matching
procedure}. Several free constants of inner and outer solutions must
be determined by matching conditions. The inner equations, their solutions
as well as matching is discussed in Section \ref{sub:InnerEquations}.
It remains to combine all inner and outer solutions in a single \textit{composite
solution.} Only the composite solution yields a uniformly valid approximation
over the whole time domain. The control signal is given in terms of
the composite solution. Both points are discussed in Section \ref{sub:CompositeSolutionsAndControl}.
The last step involves a discussion of the exact solution obtained
for $\epsilon=0$ in Section \ref{sub:TheLimitEpsilon0}.

\subsection{Necessary optimality conditions}

The solution to optimal trajectory tracking is derived for nonlinear
two-dimensional dynamical systems of the form
\begin{align}
\dot{x}\left(t\right) & =a_{0}+a_{1}x\left(t\right)+a_{2}y\left(t\right),\label{eq:xState}\\
\dot{y}\left(t\right) & =R\left(x\left(t\right),y\left(t\right)\right)+b\left(x\left(t\right),y\left(t\right)\right)u\left(t\right).\label{eq:yState}
\end{align}
The parameters $a_{0}$ and $a_{1}$ are arbitrary. The function $b\left(x,y\right)$
is not allowed to vanish for any value of $x$ and $y$. The system
is controllable as long as $a_{2}\neq0$. This was proven in Example
\ref{ex:FHN4} of Section \ref{sec:ControllabilityForRealizableTrajectories},
and $a_{2}\neq0$ is assumed in the following. The optimal control
problem is the minimization of the functional 
\begin{align}
\mathcal{J}\left[\boldsymbol{x}\left(t\right),u\left(t\right)\right] & =\frac{1}{2}\intop_{t_{0}}^{t_{1}}\left(s_{1}\left(x\left(t\right)-x_{d}\left(t\right)\right)^{2}+s_{2}\left(y\left(t\right)-y_{d}\left(t\right)\right)^{2}\right)dt\nonumber \\
 & +\frac{\beta_{1}}{2}\left(x\left(t_{1}\right)-x_{1}\right)^{2}+\frac{\epsilon^{2}}{2}\intop_{t_{0}}^{t_{1}}dt\left(u\left(t\right)-u_{0}\right)^{2}.\label{eq:Functional}
\end{align}
The minimization of Eq. \eqref{eq:Functional} is constrained by the
condition that $x\left(t\right)$ and $y\left(t\right)$ are a solution
of the controlled dynamical system Eqs. \eqref{eq:xState}, \eqref{eq:yState},
together with the initial and terminal conditions
\begin{align}
x\left(t_{0}\right) & =x_{0}, & y\left(t_{0}\right) & =y_{0}, & y\left(t_{1}\right) & =y_{1}.
\end{align}
In contrast to $x\left(t\right)$, the $y$-component satisfies a
sharp terminal condition. This special choice for the terminal conditions
is motivated by the exact solution discussed in Section \ref{sec:AnExactlySolvableExample}.
All relevant cases of terminal conditions are covered by the weighting
coefficient $0\leq\beta_{1}\leq\infty$. The constants $s_{1}>0$
and $s_{2}>0$ are positive weights which correspond to a positive
definite diagonal matrix of weighting coefficients,
\begin{align}
\boldsymbol{\mathcal{S}} & =\left(\begin{array}{cc}
s_{1} & 0\\
0 & s_{2}
\end{array}\right).
\end{align}
Equation \eqref{eq:Functional} takes a constant background value
$u_{0}$ into account. For a vanishing regularization coefficient
$\epsilon=0$, the minimization problem becomes singular. The perturbation
expansion applies in the limit of small $\epsilon>0$. The co-state
is denoted by
\begin{align}
\boldsymbol{\lambda}\left(t\right) & =\left(\begin{array}{c}
\lambda_{x}\left(t\right)\\
\lambda_{y}\left(t\right)
\end{array}\right).
\end{align}
The adjoint or co-state equation involves the Jacobian $\nabla\boldsymbol{R}$
of the nonlinearity $\boldsymbol{R}$,
\begin{align}
\nabla\boldsymbol{R}\left(\boldsymbol{x}\left(t\right)\right) & =\left(\begin{array}{cc}
a_{1} & a_{2}\\
\partial_{x}R\left(x\left(t\right),y\left(t\right)\right) & \partial_{y}R\left(x\left(t\right),y\left(t\right)\right)
\end{array}\right).
\end{align}
The necessary optimality conditions are (see Section \ref{sec:NecessaryOptimalityConditions}
for a discussion of the solution to general linear dynamical systems)
\begin{align}
0 & =\epsilon^{2}\left(u\left(t\right)-u_{0}\right)+b\left(x\left(t\right),y\left(t\right)\right)\lambda_{y}\left(t\right),\label{eq:PontryaginMaximumCondition}\\
\left(\begin{array}{c}
\dot{x}\left(t\right)\\
\dot{y}\left(t\right)
\end{array}\right) & =\left(\begin{array}{c}
y\left(t\right)\\
R\left(x\left(t\right),y\left(t\right)\right)
\end{array}\right)+b\left(x\left(t\right),y\left(t\right)\right)\left(\begin{array}{c}
0\\
u\left(t\right)
\end{array}\right),\label{eq:ControlledTwoDimDynamicalSystem}\\
-\left(\begin{array}{c}
\dot{\lambda}_{x}\left(t\right)\\
\dot{\lambda}_{y}\left(t\right)
\end{array}\right) & =\left(\begin{array}{cc}
a_{1} & \partial_{x}R\left(x\left(t\right),y\left(t\right)\right)+\partial_{x}b\left(x\left(t\right),y\left(t\right)\right)u\left(t\right)\\
a_{2} & \partial_{y}R\left(x\left(t\right),y\left(t\right)\right)+\partial_{y}b\left(x\left(t\right),y\left(t\right)\right)u\left(t\right)
\end{array}\right)\left(\begin{array}{c}
\lambda_{x}\left(t\right)\\
\lambda_{y}\left(t\right)
\end{array}\right)\nonumber \\
 & +\left(\begin{array}{c}
s_{1}\left(x\left(t\right)-x_{d}\left(t\right)\right)\\
s_{2}\left(y\left(t\right)-y_{d}\left(t\right)\right)
\end{array}\right),\label{eq:CostateEquations}\\
\left(\begin{array}{c}
x\left(t_{0}\right)\\
y\left(t_{0}\right)
\end{array}\right) & =\left(\begin{array}{c}
x_{0}\\
y_{0}
\end{array}\right),\label{eq:InitCondxy}\\
y\left(t_{1}\right) & =y_{1},\label{eq:TermCondy}\\
\lambda_{x}\left(t_{1}\right) & =\beta_{1}\left(x\left(t_{1}\right)-x_{1}\right).\label{eq:TermCondlambdax}
\end{align}

\subsection{\label{sub:RearrangeNecessaryOptimalityConditions}Rearranging the
necessary optimality conditions}

The first step is to solve Eq. \eqref{eq:ControlledTwoDimDynamicalSystem}
for the control signal,
\begin{align}
u\left(t\right) & =\frac{1}{b\left(x\left(t\right),y\left(t\right)\right)}\left(\dot{y}\left(t\right)-R\left(x\left(t\right),y\left(t\right)\right)\right).
\end{align}
Using $u\left(t\right)$ in Eq. \eqref{eq:PontryaginMaximumCondition}
yields
\begin{align}
\lambda_{y}\left(t\right) & =-\frac{\epsilon^{2}}{b\left(x\left(t\right),y\left(t\right)\right)^{2}}\left(\dot{y}\left(t\right)-R\left(x\left(t\right),y\left(t\right)\right)-b\left(x\left(t\right),y\left(t\right)\right)u_{0}\right).\label{eq:Lambday1}
\end{align}
This expression holds for all times $t$ such that it is allowed to
apply the time derivative,
\begin{align}
\dot{\lambda}_{y}\left(t\right) & =-\epsilon^{2}\frac{1}{b\left(x\left(t\right),y\left(t\right)\right)^{2}}\dot{x}\left(t\right)\left(u_{0}\partial_{x}b\left(x\left(t\right),y\left(t\right)\right)-\partial_{x}R\left(x\left(t\right),y\left(t\right)\right)\right)\nonumber \\
 & -\epsilon^{2}\frac{1}{b\left(x\left(t\right),y\left(t\right)\right)^{2}}\left(\dot{y}\left(t\right)\left(u_{0}\partial_{y}b\left(x\left(t\right),y\left(t\right)\right)-\partial_{y}R\left(x\left(t\right),y\left(t\right)\right)\right)+\ddot{y}\left(t\right)\right)\nonumber \\
 & -\epsilon^{2}\frac{2\left(R\left(x\left(t\right),y\left(t\right)\right)-\dot{y}\left(t\right)\right)}{b\left(x\left(t\right),y\left(t\right)\right)^{3}}\left(\dot{x}\left(t\right)\partial_{x}b\left(x\left(t\right),y\left(t\right)\right)+\dot{y}\left(t\right)\partial_{y}b\left(x\left(t\right),y\left(t\right)\right)\right).\label{eq:Lambday2}
\end{align}
Using both relations Eqs. \eqref{eq:Lambday1} and Eq. \eqref{eq:Lambday2},
all occurrences of $\lambda_{y}\left(t\right)$ can be eliminated
in the co-state equations \eqref{eq:CostateEquations}. The equation
for $\lambda_{x}$ becomes 
\begin{align}
-\dot{\lambda}_{x}\left(t\right) & =a_{1}\lambda_{x}\left(t\right)+s_{1}\left(x\left(t\right)-x_{d}\left(t\right)\right)\nonumber \\
 & -\frac{\epsilon^{2}}{b\left(x\left(t\right),y\left(t\right)\right)^{2}}\left(\partial_{x}R\left(x\left(t\right),y\left(t\right)\right)+\partial_{x}b\left(x\left(t\right),y\left(t\right)\right)u\left(t\right)\right)\nonumber \\
 & \times\left(\dot{y}\left(t\right)-R\left(x\left(t\right),y\left(t\right)\right)-b\left(x\left(t\right),y\left(t\right)\right)u_{0}\right),\label{eq:LambdaXTransformed}
\end{align}
while the equation for $\lambda_{y}$ transforms to a second order
differential equation for $y$, 
\begin{align}
 & \epsilon^{2}\ddot{y}\left(t\right)+\epsilon^{2}\dot{x}\left(t\right)\left(u_{0}\partial_{x}b\left(x\left(t\right),y\left(t\right)\right)-\partial_{x}R\left(x\left(t\right),y\left(t\right)\right)\right)\nonumber \\
 & +\epsilon^{2}\dot{y}\left(t\right)\left(u_{0}\partial_{y}b\left(x\left(t\right),y\left(t\right)\right)-\partial_{y}R\left(x\left(t\right),y\left(t\right)\right)\right)\nonumber \\
 & +2\epsilon^{2}\frac{\left(R\left(x\left(t\right),y\left(t\right)\right)-\dot{y}\left(t\right)\right)}{b\left(x\left(t\right),y\left(t\right)\right)}w_{1}\left(x\left(t\right),y\left(t\right)\right)\nonumber \\
 & =b\left(x\left(t\right),y\left(t\right)\right)^{2}\left(a_{2}\lambda_{x}\left(t\right)+s_{2}\left(y\left(t\right)-y_{d}\left(t\right)\right)\right)\nonumber \\
 & +\epsilon^{2}w_{2}\left(x\left(t\right),y\left(t\right)\right)\left(b\left(x\left(t\right),y\left(t\right)\right)u_{0}-\left(\dot{y}\left(t\right)-R\left(x\left(t\right),y\left(t\right)\right)\right)\right).\label{eq:YTransformed}
\end{align}
Here, $w_{1}$ and $w_{2}$ denote the abbreviations
\begin{align}
w_{1}\left(x\left(t\right),y\left(t\right)\right) & =\dot{x}\left(t\right)\partial_{x}b\left(x\left(t\right),y\left(t\right)\right)+\dot{y}\left(t\right)\partial_{y}b\left(x\left(t\right),y\left(t\right)\right),\\
w_{2}\left(x\left(t\right),y\left(t\right)\right) & =\partial_{y}R\left(x\left(t\right),y\left(t\right)\right)+\frac{\partial_{y}b\left(x\left(t\right),y\left(t\right)\right)}{b\left(x\left(t\right),y\left(t\right)\right)}\left(\dot{y}\left(t\right)-R\left(x\left(t\right),y\left(t\right)\right)\right).
\end{align}
The equation for $x\left(t\right)$ does not change,
\begin{align}
\dot{x}\left(t\right) & =a_{0}+a_{1}x\left(t\right)+a_{2}y\left(t\right).\label{eq:XTransformed}
\end{align}
The system of equations \eqref{eq:LambdaXTransformed}-\eqref{eq:XTransformed}
must be solved with four initial and terminal conditions Eqs. \eqref{eq:InitCondxy}-\eqref{eq:TermCondlambdax}.
The rearranged equations look horrible and much more difficult than
before. However, the small parameter $\epsilon^{2}$ multiplies the
second order derivative of $y\left(t\right)$ as well as every occurrence
of the nonlinear force term $R$ and its derivatives. Due to $\epsilon$
multiplying the highest order derivative, Eqs. \eqref{eq:LambdaXTransformed}-\eqref{eq:XTransformed}
constitute a singularly perturbed system of differential equations.
Setting $\epsilon=0$ changes the order of the system. Consequently,
not all four boundary conditions Eqs. \eqref{eq:InitCondxy}-\eqref{eq:TermCondlambdax}
can be satisfied. Singular perturbation theory has to be applied to
solve Eqs. \eqref{eq:LambdaXTransformed}-\eqref{eq:XTransformed}
in the limit $\epsilon\rightarrow0$.

\subsection{\label{sub:OuterEquations}Outer equations}

The outer equations are defined for the ordinary time scale $t$.
The outer solutions are denoted with index $O$, 
\begin{align}
x_{O}\left(t\right) & =x\left(t\right), & y_{O}\left(t\right) & =y\left(t\right), & \lambda_{O}\left(t\right) & =\lambda_{x}\left(t\right).
\end{align}
Expanding Eqs. \eqref{eq:LambdaXTransformed}-\eqref{eq:XTransformed}
up to leading order in $\epsilon$ yields two linear differential
equations of first order and an algebraic equation, 
\begin{align}
\dot{\lambda}_{O}\left(t\right) & =-a_{1}\lambda_{x}\left(t\right)+s_{1}\left(x_{d}\left(t\right)-x_{O}\left(t\right)\right),\label{eq:Outer1}\\
\lambda_{O}\left(t\right) & =\frac{s_{2}}{a_{2}}\left(y_{d}\left(t\right)-y_{O}\left(t\right)\right),\label{eq:Outer2}\\
\dot{x}_{O}\left(t\right) & =a_{0}+a_{1}x_{O}\left(t\right)+a_{2}y_{O}\left(t\right).\label{eq:Outer3}
\end{align}
Equations \eqref{eq:Outer1}-\eqref{eq:Outer3} allow for two initial
conditions, 
\begin{align}
x_{O}\left(t_{0}\right) & =x_{\text{init}}, & y_{O}\left(t_{0}\right) & =y_{\text{init}}.
\end{align}
Note that $x_{\text{init}}$ and $y_{\text{init}}$ are not given
by the initial conditions \eqref{eq:InitCondxy} but have to be determined
by matching with the inner solutions. Eliminating $\lambda_{x}\left(t\right)$
from Eqs. \eqref{eq:Outer1}-\eqref{eq:Outer3} yields two coupled
ODEs for $x_{O}\left(t\right)$ and $y_{O}\left(t\right)$,
\begin{align}
\dot{x}_{O}\left(t\right) & =a_{1}x_{O}\left(t\right)+a_{2}y_{O}\left(t\right)+a_{0},\label{eq:xOEquation}\\
\dot{y}_{O}\left(t\right) & =-a_{1}y_{O}\left(t\right)+\frac{a_{2}s_{1}}{s_{2}}x_{O}\left(t\right)+\dot{y}_{d}\left(t\right)+a_{1}y_{d}\left(t\right)-\frac{a_{2}s_{1}}{s_{2}}x_{d}\left(t\right).\label{eq:yOEquation}
\end{align}
It is convenient to express the solutions for $x_{O}$ and $y_{O}$
in terms of the state transition matrix $\boldsymbol{\Phi}\left(t,t_{0}\right)$
(see Appendix \ref{sec:GeneralSolutionForForcedLinarDynamicalSystem}
for a general derivation of state transition matrices)
\begin{align}
\left(\begin{array}{c}
x_{O}\left(t\right)\\
y_{O}\left(t\right)
\end{array}\right) & =\boldsymbol{\Phi}\left(t,t_{0}\right)\left(\begin{array}{c}
x_{\text{init}}\\
y_{\text{init}}
\end{array}\right)+\intop_{t_{0}}^{t}d\tau\boldsymbol{\Phi}\left(t,\tau\right)\boldsymbol{f}\left(\tau\right),
\end{align}
with
\begin{align}
\boldsymbol{\Phi}\left(t,t_{0}\right) & =\left(\begin{array}{cc}
\cosh\left(\Delta t\varphi_{1}\right)+\frac{a_{1}}{\varphi_{1}}\sinh\left(\Delta t\varphi_{1}\right) & \frac{a_{2}}{\varphi_{1}}\sinh\left(\Delta t\varphi_{1}\right)\\
\frac{a_{2}s_{1}}{s_{2}\varphi_{1}}\sinh\left(\Delta t\varphi_{1}\right) & \cosh\left(\Delta t\varphi_{1}\right)-\frac{a_{1}}{\varphi_{1}}\sinh\left(\Delta t\varphi_{1}\right)
\end{array}\right),\\
\varphi_{1} & =\frac{\sqrt{a_{1}^{2}s_{2}+a_{2}^{2}s_{1}}}{\sqrt{s_{2}}},\,\Delta t=t-t_{0},
\end{align}
and inhomogeneity
\begin{align}
\boldsymbol{f}\left(t\right) & =\left(\begin{array}{c}
a_{0}\\
\dot{y}_{d}\left(t\right)+a_{1}y_{d}\left(t\right)-\frac{a_{2}s_{1}}{s_{2}}x_{d}\left(t\right)
\end{array}\right).
\end{align}
For later reference, the solutions for $x_{O}\left(t\right)$ and
$y_{O}\left(t\right)$ are given explicitly,
\begin{align}
x_{O}\left(t\right) & =\frac{a_{1}a_{2}}{\varphi_{1}}\intop_{t_{0}}^{t}y_{d}\left(\tau\right)\sinh\left(\varphi_{1}\left(t-\tau\right)\right)\,d\tau+a_{2}\intop_{t_{0}}^{t}y_{d}\left(\tau\right)\cosh\left(\varphi_{1}\left(t-\tau\right)\right)\,d\tau\nonumber \\
 & -\frac{a_{2}^{2}s_{1}}{s_{2}\varphi_{1}}\intop_{t_{0}}^{t}x_{d}\left(\tau\right)\sinh\left(\varphi_{1}\left(t-\tau\right)\right)\,d\tau+\frac{1}{\varphi_{1}}\sinh\left(\left(t-t_{0}\right)\varphi_{1}\right)\left(a_{0}-a_{2}y_{d}\left(t_{0}\right)\right)\nonumber \\
 & +\frac{a_{0}a_{1}}{\varphi_{1}^{2}}\left(\cosh\left(\left(t-t_{0}\right)\varphi_{1}\right)-1\right)\nonumber \\
 & +x_{\text{init}}\left(\frac{a_{1}}{\varphi_{1}}\sinh\left(\left(t-t_{0}\right)\varphi_{1}\right)+\cosh\left(\left(t-t_{0}\right)\varphi_{1}\right)\right)\nonumber \\
 & +\frac{a_{2}}{\varphi_{1}}y_{\text{init}}\sinh\left(\left(t-t_{0}\right)\varphi_{1}\right),\label{eq:xOSol}
\end{align}
 
\begin{align}
y_{O}\left(t\right) & =\frac{a_{0}a_{2}s_{1}}{s_{2}\varphi_{1}^{2}}\left(\cosh\left(\left(t-t_{0}\right)\varphi_{1}\right)-1\right)+\frac{a_{2}s_{1}}{s_{2}\varphi_{1}}x_{\text{init}}\sinh\left(\left(t-t_{0}\right)\varphi_{1}\right)\nonumber \\
 & +y_{\text{init}}\cosh\left(\left(t-t_{0}\right)\varphi_{1}\right)-y_{\text{init}}\frac{a_{1}}{\varphi_{1}}\sinh\left(\left(t-t_{0}\right)\varphi_{1}\right)\nonumber \\
 & -\frac{a_{2}s_{1}}{s_{2}}\intop_{t_{0}}^{t}x_{d}\left(\tau\right)\cosh\left(\varphi_{1}\left(t-\tau\right)\right)\,d\tau+\frac{a_{2}s_{1}a_{1}}{s_{2}\varphi_{1}}\intop_{t_{0}}^{t}x_{d}\left(\tau\right)\sinh\left(\varphi_{1}\left(t-\tau\right)\right)\,d\tau\nonumber \\
 & +\frac{a_{2}^{2}s_{1}}{s_{2}\varphi_{1}}\intop_{t_{0}}^{t}y_{d}\left(\tau\right)\sinh\left(\varphi_{1}\left(t-\tau\right)\right)\,d\tau\nonumber \\
 & +\frac{a_{1}}{\varphi_{1}}y_{d}\left(t_{0}\right)\sinh\left(\left(t-t_{0}\right)\varphi_{1}\right)+y_{d}\left(t\right)-y_{d}\left(t_{0}\right)\cosh\left(\left(t-t_{0}\right)\varphi_{1}\right).\label{eq:yOSol}
\end{align}
The solution for the co-state $\lambda_{O}\left(t\right)$ reads
\begin{align}
\lambda_{O}\left(t\right) & =\frac{s_{2}}{a_{2}}\left(y_{d}\left(t\right)-y_{O}\left(t\right)\right).\label{eq:lambdaOSol}
\end{align}
Evaluating Eq. \eqref{eq:lambdaOSol} at the initial and terminal
time $t_{0}$ and $t_{1}$ yields relations for $\lambda_{O}\left(t_{0}\right)$
and $\lambda_{O}\left(t_{1}\right)$, respectively,
\begin{align}
\lambda_{O}\left(t_{0}\right) & =\frac{s_{2}}{a_{2}}\left(y_{d}\left(t_{0}\right)-y_{\text{init}}\right),\label{eq:Outer2t0}\\
\lambda_{O}\left(t_{1}\right) & =\frac{s_{2}}{a_{2}}\left(y_{d}\left(t_{1}\right)-y_{\text{end}}\right).\label{eq:Outer2t1}
\end{align}
The abbreviation $y_{\text{end}}$ is
\begin{align}
y_{\text{end}} & =y_{O}\left(t_{1}\right).\label{eq:Abbreviation437}
\end{align}
Equations \eqref{eq:Outer2t0}-\eqref{eq:Abbreviation437} will be
useful for matching.

\subsection{\label{sub:InnerEquations}Inner equations}

\subsubsection{\label{sub:InitialBoundaryLayer}Initial boundary layer}

Boundary layers occur at both ends of the time domain. The initial
boundary layer at the left end of the time domain is resolved using
the time scale $\tau_{L}=\left(t-t_{0}\right)/\epsilon$ and scaled
solutions
\begin{align}
X_{L}\left(\tau_{L}\right) & =X_{L}\left(\left(t-t_{0}\right)/\epsilon\right)=x\left(t\right)=x\left(t_{0}+\epsilon\tau_{L}\right),\label{eq:XLDef}\\
Y_{L}\left(\tau_{L}\right) & =Y_{L}\left(\left(t-t_{0}\right)/\epsilon\right)=y\left(t\right)=y\left(t_{0}+\epsilon\tau_{L}\right),\label{eq:YLDef}\\
\Lambda_{L}\left(\tau_{L}\right) & =\Lambda_{L}\left(\left(t-t_{0}\right)/\epsilon\right)=\lambda_{x}\left(t\right)=\lambda_{x}\left(t_{0}+\epsilon\tau_{L}\right).\label{eq:LambdaLDef}
\end{align}
From the definitions of $X_{L},\,Y_{L},$ and $\Lambda_{L}$ together
with the initial conditions for $x$ and $y$, Eqs. \eqref{eq:InitCondxy},
follow the initial conditions
\begin{align}
x\left(t_{0}\right) & =X_{L}\left(0\right)=x_{0}, & y\left(t_{0}\right) & =Y_{L}\left(0\right)=y_{0}.
\end{align}
The derivatives of $x$ transform as
\begin{align}
\dot{x}\left(t\right) & =\dfrac{1}{\epsilon}X_{L}'\left(\tau_{L}\right), & \ddot{x}\left(t\right) & =\dfrac{1}{\epsilon^{2}}X_{L}''\left(\tau_{L}\right),
\end{align}
and analogously for $y$ and $\lambda_{x}$. The prime $X_{L}'\left(\tau_{L}\right)$
denotes the derivative of $X_{L}$ with respect to its argument. The
matching conditions at the left boundary layer are
\begin{align}
\lim_{t\rightarrow t_{0}}x_{O}\left(t\right) & =\lim_{\tau_{L}\rightarrow\infty}X_{L}\left(\tau_{L}\right),\label{eq:MatchXL}\\
\lim_{t\rightarrow t_{0}}y_{O}\left(t\right) & =\lim_{\tau_{L}\rightarrow\infty}Y_{L}\left(\tau_{L}\right),\label{eq:MatchYL}\\
\lim_{t\rightarrow t_{0}}\lambda_{O}\left(t\right) & =\lim_{\tau_{L}\rightarrow\infty}\Lambda_{L}\left(\tau_{L}\right).\label{eq:MatchLambdaL}
\end{align}
Using the definitions Eqs. \eqref{eq:XLDef}-\eqref{eq:LambdaLDef}
in Eqs. \eqref{eq:LambdaXTransformed}-\eqref{eq:XTransformed} and
expanding in $\epsilon$ yields the left inner equations in leading
order as
\begin{align}
\Lambda_{L}'\left(\tau_{L}\right) & =0,\label{eq:Eq447}\\
Y_{L}''\left(\tau_{L}\right) & =Y_{L}'\left(\tau_{L}\right)^{2}\frac{\partial_{y}b\left(X_{L}\left(\tau_{L}\right),Y_{L}\left(\tau_{L}\right)\right)}{b\left(X_{L}\left(\tau_{L}\right),Y_{L}\left(\tau_{L}\right)\right)}+2X_{L}'\left(\tau_{L}\right)Y_{L}'\left(\tau_{L}\right)\frac{\partial_{x}b\left(X_{L}\left(\tau_{L}\right),Y_{L}\left(\tau_{L}\right)\right)}{b\left(X_{L}\left(\tau_{L}\right),Y_{L}\left(\tau_{L}\right)\right)}\nonumber \\
 & +b\left(X_{L}\left(\tau_{L}\right),Y_{L}\left(\tau_{L}\right)\right){}^{2}\left(s_{2}\left(Y_{L}\left(\tau_{L}\right)-y_{d}\left(t_{0}\right)\right)+a_{2}\Lambda_{L}\left(\tau_{L}\right)\right),\\
X_{L}'\left(\tau_{L}\right) & =0.\label{eq:Eq449}
\end{align}
The differential equations do not involve the nonlinearity $R$. The
solutions for $\Lambda_{L}$ and $X_{L}$ are
\begin{align}
\Lambda_{L}\left(\tau_{L}\right) & =\Lambda_{L,0}=\lambda_{O}\left(t_{0}\right)=\frac{s_{2}}{a_{2}}\left(y_{d}\left(t_{0}\right)-y_{\text{init}}\right),\\
X_{L}\left(\tau_{L}\right) & =x_{0}.
\end{align}
To obtain the value for $\Lambda_{L,0}$, Eq. \eqref{eq:Outer2t0}
was used together with the matching condition Eq. \eqref{eq:MatchLambdaL}.
The matching condition for $x$, Eq. \eqref{eq:MatchXL}, immediately
yields
\begin{align}
x_{\text{init}} & =x_{0},
\end{align}
while $y_{\text{init}}$ will be determined later on. The equation
for $Y_{L}$ simplifies to
\begin{align}
Y_{L}''\left(\tau_{L}\right) & =Y_{L}'\left(\tau_{L}\right)^{2}\frac{\partial_{y}b\left(x_{0},Y_{L}\left(\tau_{L}\right)\right)}{b\left(x_{0},Y_{L}\left(\tau_{L}\right)\right)}+s_{2}b\left(x_{0},Y_{L}\left(\tau_{L}\right)\right){}^{2}\left(Y_{L}\left(\tau_{L}\right)-y_{\text{init}}\right).
\end{align}
As long as $b$ depends on $Y_{L}$, this is a nonlinear equation.
Because it is autonomous, it can be transformed to a first order ODE
by introducing a new function $v_{L}$ defined as
\begin{align}
Y_{L}'\left(\tau_{L}\right) & =v_{L}\left(Y_{L}\left(\tau_{L}\right)\right)b\left(x_{0},Y_{L}\left(\tau_{L}\right)\right)\label{eq:DefvL}
\end{align}
such that the second order time derivative is
\begin{align}
Y_{L}''\left(\tau_{L}\right) & =v_{L}'\left(Y_{L}\left(\tau_{L}\right)\right)v_{L}\left(Y_{L}\left(\tau_{L}\right)\right)b\left(x_{0},Y_{L}\left(\tau_{L}\right)\right)^{2}\nonumber \\
 & +v_{L}\left(Y_{L}\left(\tau_{L}\right)\right)^{2}\partial_{y}b\left(x_{0},Y_{L}\left(\tau_{L}\right)\right)b\left(x_{0},Y_{L}\left(\tau_{L}\right)\right).
\end{align}
This leads to a fairly simple equation for $v_{L}$,
\begin{align}
\frac{1}{2}\partial_{Y_{L}}\left(v_{L}\left(Y_{L}\right)\right)^{2} & =s_{2}\left(Y_{L}-y_{\text{init}}\right).\label{eq:EqvL}
\end{align}
Equation \eqref{eq:EqvL} is to be solved with the matching condition
Eq. \eqref{eq:MatchYL},
\begin{align}
\lim_{\tau_{L}\rightarrow\infty}Y_{L}\left(\tau_{L}\right) & =\lim_{t\rightarrow t_{0}}y_{O}\left(t\right)=y_{\text{init}}.\label{eq:Eq732}
\end{align}
A limit $\lim_{\tau_{L}\rightarrow\infty}Y_{L}\left(\tau_{L}\right)$
exists if $Y_{L}\left(\tau_{L}\right)$ is neither infinite nor oscillatory
as $\tau_{L}$ approaches $\infty$. The existence of both $\lim_{\tau_{L}\rightarrow\infty}Y_{L}\left(\tau_{L}\right)$
and $\lim_{\tau_{L}\rightarrow\infty}Y_{L}'\left(\tau_{L}\right)$
implies
\begin{align}
\lim_{\tau_{L}\rightarrow\infty}Y_{L}'\left(\tau_{L}\right) & =0.\label{eq:Eq733}
\end{align}
A proof is straightforward. From Eq. \eqref{eq:Eq732} follows
\begin{align}
\lim_{\tau_{L}\rightarrow\infty}\frac{Y_{L}\left(\tau_{L}\right)}{\tau_{L}} & =\lim_{\tau_{L}\rightarrow\infty}\frac{y_{\text{init}}}{\tau_{L}}=0,
\end{align}
and applying L'H\^opital's rule to
\begin{align}
\lim_{\tau_{L}\rightarrow\infty}\frac{Y_{L}\left(\tau_{L}\right)-y_{\text{init}}}{\tau_{L}} & =\lim_{\tau_{L}\rightarrow\infty}\frac{Y_{L}\left(\tau_{L}\right)}{\tau_{L}}-\lim_{\tau_{L}\rightarrow\infty}\frac{y_{\text{init}}}{\tau_{L}}=0,
\end{align}
and 
\begin{align}
\lim_{\tau_{L}\rightarrow\infty}\frac{Y_{L}\left(\tau_{L}\right)-y_{\text{init}}}{\tau_{L}} & =\frac{0}{0}=\lim_{\tau_{L}\rightarrow\infty}Y_{L}'\left(\tau_{L}\right)=0
\end{align}
yields the result.

From Eq. \eqref{eq:Eq733} follows the initial condition for $v_{L}$
as 
\begin{align}
v_{L}\left(y_{\text{init}}\right) & =0.
\end{align}
Solving Eq. \eqref{eq:EqvL} with this condition yields 
\begin{align}
v_{L}\left(Y_{L}\right) & =\pm\sqrt{s_{2}}\left|y_{\text{init}}-Y_{L}\right|.\label{eq:vLSol}
\end{align}
The solution Eq. \eqref{eq:vLSol} together with the definition for
$v_{L}$, Eq. \eqref{eq:DefvL}, leads to a first order nonlinear
ODE for $Y_{L}$,
\begin{align}
Y_{L}'\left(\tau_{L}\right) & =\pm\sqrt{s_{2}}\left|y_{\text{init}}-Y_{L}\left(\tau_{L}\right)\right|b\left(x_{0},Y_{L}\left(\tau_{L}\right)\right).\label{eq:EqYLTemp}
\end{align}
The last point is to determine which sign in Eq. \eqref{eq:EqYLTemp}
is the relevant one. Note that $Y_{L}\left(\tau_{L}\right)=y_{\text{init}}$
is a stationary point of Eq. \eqref{eq:EqYLTemp} which cannot be
crossed by the dynamics. Furthermore, this is the only stationary
point. Because of $b\left(x_{0},Y_{L}\left(\tau_{L}\right)\right)\neq0$
by assumption, $b$ cannot change its sign. If initially $y_{\text{init}}>Y_{L}\left(0\right)=y_{0}$
and $b\left(x_{0},Y_{L}\left(\tau_{L}\right)\right)>0$ for all times
$\tau_{L}>0$, $Y_{L}$ must grow and therefore $Y_{L}'=\sqrt{s_{2}}\left(y_{\text{init}}-Y_{L}\right)b\left(x_{0},Y_{L}\right)>0$
is the correct choice. On the other hand, if $y_{\text{init}}>Y_{L}\left(0\right)=y_{0}$
and $b\left(x_{0},Y_{L}\left(\tau_{L}\right)\right)<0$ for all times
$\tau_{L}>0$, $Y_{L}$ must decrease and consequently $Y_{L}$ evolves
according to $Y_{L}'=-\sqrt{s_{2}}\left(y_{\text{init}}-Y_{L}\right)b\left(x_{0},Y_{L}\right)>0$.
These considerations finally lead to
\begin{align}
Y_{L}'\left(\tau_{L}\right) & =\sqrt{s_{2}}\left(y_{\text{init}}-Y_{L}\left(\tau_{L}\right)\right)\left|b\left(x_{0},Y_{L}\left(\tau_{L}\right)\right)\right|,\label{eq:EqYL}\\
Y_{L}\left(0\right) & =y_{0}.
\end{align}
An analytical solution of Eq. \eqref{eq:EqYL} for arbitrary functions
$b$ does not exist in closed form. If $b\left(x,y\right)=b\left(x\right)$
does not depend on $y$, Eq. \eqref{eq:EqYL} is linear and has the
solution
\begin{align}
Y_{L}\left(\tau_{L}\right) & =y_{\text{init}}+\exp\left(-\sqrt{s_{2}}\left|b\left(x_{0}\right)\right|\tau_{L}\right)\left(y_{0}-y_{\text{init}}\right).
\end{align}

\subsubsection{\label{sub:TerminalBoundaryLayer}Terminal boundary layer}

A treatment analogous to Section \ref{sub:InitialBoundaryLayer} is
performed to resolve the boundary layer at the right end of the time
domain. The relevant time scale is $\tau_{R}=\left(t_{1}-t\right)/\epsilon$,
and the scaled solutions are defined as
\begin{align}
X_{R}\left(\tau_{R}\right) & =X_{R}\left(\left(t_{1}-t\right)/\epsilon\right)=x\left(t\right),\\
Y_{R}\left(\tau_{R}\right) & =Y_{R}\left(\left(t_{1}-t\right)/\epsilon\right)=y\left(t\right),\\
\Lambda_{R}\left(\tau_{R}\right) & =\Lambda_{R}\left(\left(t_{1}-t\right)/\epsilon\right)=\lambda_{x}\left(t\right).
\end{align}
The terminal conditions Eqs. \eqref{eq:TermCondlambdax} and \eqref{eq:TermCondy}
lead to the boundary conditions
\begin{align}
\Lambda_{R}\left(0\right) & =\beta_{1}\left(X_{R}\left(0\right)-x_{1}\right), & Y_{R}\left(0\right) & =y_{1}.
\end{align}
Furthermore, $X_{R},\,Y_{R},$ and $\Lambda_{R}$ have to satisfy
the matching conditions
\begin{align}
\lim_{t\rightarrow t_{1}}x_{O}\left(t\right) & =\lim_{\tau_{R}\rightarrow\infty}X_{R}\left(\tau_{R}\right),\label{eq:MatchXR}\\
\lim_{t\rightarrow t_{1}}y_{O}\left(t\right) & =\lim_{\tau_{R}\rightarrow\infty}Y_{R}\left(\tau_{R}\right),\label{eq:MatchYR}\\
\lim_{t\rightarrow t_{1}}\lambda_{O}\left(t\right) & =\lim_{\tau_{R}\rightarrow\infty}\Lambda_{R}\left(\tau_{R}\right).\label{eq:MatchLambdaR}
\end{align}
The derivatives of $x$ transform as
\begin{align}
\dot{x}\left(t\right) & =-\dfrac{1}{\epsilon}X_{R}'\left(\tau_{R}\right), & \ddot{x}\left(t\right) & =\dfrac{1}{\epsilon^{2}}X_{R}''\left(\tau_{R}\right),
\end{align}
and analogously for $y$ and $\lambda_{x}$. Plugging these definitions
in Eqs. \eqref{eq:LambdaXTransformed}-\eqref{eq:XTransformed} and
expanding in $\epsilon$ yields the right inner equations in leading
order,
\begin{align}
\Lambda_{R}'\left(\tau_{R}\right) & =0,\\
Y_{R}''\left(\tau_{R}\right) & =Y_{R}'\left(\tau_{R}\right)^{2}\frac{\partial_{y}b\left(X_{R}\left(\tau_{R}\right),Y_{R}\left(\tau_{R}\right)\right)}{b\left(X_{R}\left(\tau_{R}\right),Y_{R}\left(\tau_{R}\right)\right)}+2X_{R}'\left(\tau_{R}\right)Y_{R}'\left(\tau_{R}\right)\frac{\partial_{x}b\left(X_{R}\left(\tau_{R}\right),Y_{R}\left(\tau_{R}\right)\right)}{b\left(X_{R}\left(\tau_{R}\right),Y_{R}\left(\tau_{R}\right)\right)}\nonumber \\
 & +b\left(X_{R}\left(\tau_{R}\right),Y_{R}\left(\tau_{R}\right)\right){}^{2}\left(s_{2}\left(Y_{R}\left(\tau_{R}\right)-y_{d}\left(t_{1}\right)\right)+a_{2}\Lambda_{R}\left(\tau_{R}\right)\right),\\
X_{R}'\left(\tau_{R}\right) & =0.
\end{align}
These equations are identical in form to the left inner equations
\eqref{eq:Eq447}-\eqref{eq:Eq449}. The solutions for $X_{R}\left(\tau_{R}\right)$
and $\Lambda_{R}\left(\tau_{R}\right)$ are constant and can be written
as
\begin{align}
X_{R}\left(\tau_{R}\right) & =x_{1}+\frac{1}{\beta_{1}}\Lambda_{R,0}, & \Lambda_{R}\left(\tau_{R}\right) & =\Lambda_{R,0}.
\end{align}
Applying the matching condition Eq. \eqref{eq:MatchLambdaR} together
with Eq. \eqref{eq:Outer2t1} yields the solution for $\Lambda_{R}\left(\tau_{R}\right)$
and $X_{R}\left(\tau_{R}\right)$ as 
\begin{align}
\Lambda_{R}\left(\tau_{R}\right) & =\Lambda_{R,0}=\lambda_{O}\left(t_{1}\right)=\frac{s_{2}}{a_{2}}\left(y_{d}\left(t_{1}\right)-y_{\text{end}}\right),\\
X_{R}\left(\tau_{R}\right) & =x_{1}+\frac{s_{2}}{\beta_{1}a_{2}}\left(y_{d}\left(t_{1}\right)-y_{\text{end}}\right).
\end{align}
With the analogous considerations as for the left inner equations,
see Eq. \eqref{eq:EqYL}, the solution to $Y_{R}\left(\tau_{R}\right)$
is given by the first order ODE 
\begin{align}
Y_{R}'\left(\tau_{R}\right) & =\sqrt{s_{2}}\left(y_{\text{end}}-Y_{R}\left(\tau_{R}\right)\right)\left|b\left(x_{1},Y_{R}\left(\tau_{R}\right)\right)\right|,\label{eq:EqYR}\\
Y_{R}\left(0\right) & =y_{1}.
\end{align}
Equation \eqref{eq:EqYR} satisfies the matching condition Eq. \eqref{eq:MatchYR}.
The remaining matching condition Eq. \eqref{eq:MatchXR} gives
\begin{align}
x_{O}\left(t_{1}\right) & =x_{1}+\frac{s_{2}}{\beta_{1}a_{2}}\left(y_{d}\left(t_{1}\right)-y_{\text{end}}\right).\label{eq:Eq484}
\end{align}
The constant $y_{\text{end}}=y_{O}\left(t_{1}\right)$ depends on
$y_{\text{init}}$, which is the last free parameter of the outer
solution. Solving Eq. \eqref{eq:Eq484} for $y_{\text{init}}$ yields
\begin{align}
y_{\text{init}} & =\frac{s_{2}\varphi_{1}}{\kappa}\sinh\left(\left(t_{0}-t_{1}\right)\varphi_{1}\right)\left(\left(a_{1}s_{2}-a_{2}^{2}\beta_{1}\right)y_{d}\left(t_{0}\right)+a_{2}\left(a_{0}\beta_{1}+x_{0}\left(a_{1}\beta_{1}+s_{1}\right)\right)\right)\nonumber \\
 & -\frac{1}{\kappa}\cosh\left(\left(t_{1}-t_{0}\right)\varphi_{1}\right)\left(a_{2}\left(\beta_{1}\left(a_{1}^{2}s_{2}+a_{2}^{2}s_{1}\right)x_{0}+a_{0}s_{2}\left(a_{1}\beta_{1}+s_{1}\right)\right)\right)\nonumber \\
 & +\frac{s_{2}}{\kappa}\left(a_{1}^{2}s_{2}+a_{2}^{2}s_{1}\right)\cosh\left(\left(t_{1}-t_{0}\right)\varphi_{1}\right)y_{d}\left(t_{0}\right)\nonumber \\
 & +\frac{a_{2}\varphi_{1}s_{1}}{\kappa}\left(a_{2}^{2}\beta_{1}-a_{1}s_{2}\right)\intop_{t_{0}}^{t_{1}}x_{d}\left(\tau\right)\sinh\left(\varphi_{1}\left(t_{1}-\tau\right)\right)\,d\tau\nonumber \\
 & -\frac{a_{2}^{2}\varphi_{1}s_{2}}{\kappa}\left(a_{1}\beta_{1}+s_{1}\right)\intop_{t_{0}}^{t_{1}}y_{d}\left(\tau\right)\sinh\left(\varphi_{1}\left(t_{1}-\tau\right)\right)\,d\tau\nonumber \\
 & -\frac{\beta_{1}a_{2}^{2}\varphi_{1}^{2}s_{2}}{\kappa}\intop_{t_{0}}^{t_{1}}y_{d}\left(\tau\right)\cosh\left(\varphi_{1}\left(t_{1}-\tau\right)\right)\,d\tau+\frac{a_{2}a_{0}s_{2}}{\kappa}\left(a_{1}\beta_{1}+s_{1}\right)\nonumber \\
 & +\frac{a_{2}\varphi_{1}^{2}s_{1}s_{2}}{\kappa}\intop_{t_{0}}^{t_{1}}x_{d}\left(\tau\right)\cosh\left(\varphi_{1}\left(t_{1}-\tau\right)\right)\,d\tau+\beta_{1}x_{1}\frac{a_{2}}{\kappa}\left(a_{1}^{2}s_{2}+a_{2}^{2}s_{1}\right).\label{eq:yinit2}
\end{align}
Equation \eqref{eq:yinit2} contains the abbreviation
\begin{align}
\kappa & =s_{2}\varphi_{1}\left(a_{2}^{2}\beta_{1}-a_{1}s_{2}\right)\sinh\left(\left(t_{1}-t_{0}\right)\varphi_{1}\right)+s_{2}^{2}\varphi_{1}^{2}\cosh\left(\left(t_{1}-t_{0}\right)\varphi_{1}\right).
\end{align}
Finally, all inner and outer solutions are determined; and all matching
conditions are satisfied.

\begin{example}[Consistency check]

For appropriate parameter values, the solution derived in this section
must reduce to the corresponding leading order approximation for the
exact solution of Section \ref{sec:AnExactlySolvableExample}. For
simplicity, only the case of vanishing initial conditions, $x_{0}=y_{0}=0$,
is considered. The other parameters have values
\begin{align}
a_{0} & =0, & a_{1} & =0, & a_{2} & =1, & s_{1} & =1, & s_{2} & =1, & t_{0} & =0.
\end{align}
The desired trajectories and the coupling function are
\begin{align}
x_{d}\left(t\right) & \equiv0, & y_{d}\left(t\right) & \equiv0, & b\left(x,y\right) & \equiv1.
\end{align}
From these assumptions follows
\begin{align}
\varphi_{1} & =\frac{\sqrt{a_{1}^{2}s_{2}+a_{2}^{2}s_{1}}}{\sqrt{s_{2}}}=1, & y_{\text{init}} & =\frac{\beta_{1}x_{1}}{\kappa}, & \kappa & =\beta_{1}\sinh\left(t_{1}\right)+\cosh\left(t_{1}\right).
\end{align}
The outer solution of the controlled state is obtained as
\begin{align}
\left(\begin{array}{c}
x_{O}\left(t\right)\\
y_{O}\left(t\right)\\
\lambda_{O}\left(t\right)
\end{array}\right) & =\left(\begin{array}{c}
\dfrac{\beta_{1}x_{1}}{\kappa}\sinh\left(t\right)\\
\dfrac{\beta_{1}x_{1}}{\kappa}\cosh\left(t\right)\\
-\dfrac{\beta_{1}x_{1}}{\kappa}\cosh\left(t\right)
\end{array}\right).
\end{align}
This is indeed the outer limit of the exact solution from Section
\ref{sec:AnExactlySolvableExample}, Eq. \eqref{eq:OuterLimit} with
$x_{0}=y_{0}=0$. The abbreviation $y_{\text{end}}$ simplifies to
\begin{align}
y_{\text{end}} & =y_{O}\left(t_{1}\right)=\dfrac{\beta_{1}x_{1}}{\kappa}\cosh\left(t_{1}\right).
\end{align}
The solution to the left inner equations yields 
\begin{align}
\left(\begin{array}{c}
X_{L}\left(\tau_{L}\right)\\
Y_{L}\left(\tau_{L}\right)\\
\Lambda_{L}\left(\tau_{L}\right)
\end{array}\right) & =\left(\begin{array}{c}
0\\
\dfrac{\beta_{1}x_{1}}{\kappa}\left(1-e^{-\tau_{L}}\right)\\
-\dfrac{\beta_{1}x_{1}}{\kappa}
\end{array}\right),
\end{align}
while the solutions to the right inner equations becomes
\begin{align}
\left(\begin{array}{c}
X_{R}\left(\tau_{R}\right)\\
Y_{R}\left(\tau_{R}\right)\\
\Lambda_{R}\left(\tau_{R}\right)
\end{array}\right) & =\left(\begin{array}{c}
\dfrac{\beta_{1}x_{1}}{\kappa}\sinh\left(t_{1}\right)\\
\dfrac{\beta_{1}x_{1}}{\kappa}\cosh\left(t_{1}\right)+e^{-\tau_{R}}\left(y_{1}-\dfrac{\beta_{1}x_{1}}{\kappa}\cosh\left(t_{1}\right)\right)\\
-\dfrac{\beta_{1}x_{1}}{\kappa}\cosh\left(t_{1}\right)
\end{array}\right).
\end{align}
The left and right inner solutions indeed agree with the left and
right inner limit, Eqs. \eqref{eq:InnerLeftLimit} and \eqref{eq:InnerRightLimit},
respectively. Note that the solution for $\lambda_{y}\left(t\right)$
in leading order of $\epsilon$ vanishes in all three cases of inner
and left and right outer equations.

\end{example}

\subsection{\label{sub:CompositeSolutionsAndControl}Composite solutions and
solution for control}

The composite solutions are the sum of inner and outer solutions minus
the overlaps, 
\begin{align}
x_{\text{comp}}\left(t\right) & =x_{O}\left(t\right),\label{eq:xFinalSol}\\
y_{\text{comp}}\left(t\right) & =y_{O}\left(t\right)+Y_{L}\left(\left(t-t_{0}\right)/\epsilon\right)-y_{\text{init}}+Y_{R}\left(\left(t_{1}-t\right)/\epsilon\right)-y_{\text{end}},\label{eq:yFinalSol}\\
\lambda_{\text{comp}}\left(t\right) & =\lambda_{O}\left(t\right).\label{eq:lambdaFinalSol}
\end{align}
See Section \ref{sec:AnExactlySolvableExample} and \cite{bender1999advanced}
for further information about composite solutions. Here, $x_{O}\left(t\right),\,y_{O}\left(t\right)$,
and $\lambda_{O}\left(t\right)$ are the outer solutions Eqs. \eqref{eq:xOSol}-\eqref{eq:lambdaOSol},
while $Y_{L}\left(\tau_{L}\right)$ and $Y_{R}\left(\tau_{R}\right)$
are the left and right inner solutions given by Eqs. \eqref{eq:EqYL}
and \eqref{eq:EqYR}, respectively. The constant $y_{\text{end}}$
is defined as $y_{\text{end}}=y_{O}\left(t_{1}\right)$ and the expression
for $y_{\text{init}}$ is given by Eq. \eqref{eq:yinit2}. Equations
\eqref{eq:xFinalSol}-\eqref{eq:lambdaFinalSol} are the approximate
solution to leading order in $\epsilon$ for optimal trajectory tracking.
As a result of the singular perturbation expansion, the leading order
solution depends on $\epsilon$ itself. The solution does not depend
on the specific choice of the nonlinear term $R\left(x,y\right)$.
The sole remains left by the nonlinearities $R\left(x,y\right)$ and
$b\left(x,y\right)$ governing the system dynamics are the inner solutions
$Y_{L}$ and $Y_{R}$ which depend on the coupling function $b\left(x,y\right)$.

The general expression for the control in terms of the state components
$x\left(t\right)$ and $y\left(t\right)$ is
\begin{align}
u\left(t\right) & =\frac{1}{b\left(x\left(t\right),y\left(t\right)\right)}\left(\dot{y}\left(t\right)-R\left(x\left(t\right),y\left(t\right)\right)\right).
\end{align}
In terms of the composite solutions $x_{\text{comp}}$ and $y_{\text{comp}}$,
$u\left(t\right)$ is
\begin{align}
u\left(t\right) & =\frac{1}{b\left(x_{\text{comp}}\left(t\right),y_{\text{comp}}\left(t\right)\right)}\left(\dot{y}_{\text{comp}}\left(t\right)-R\left(x_{\text{comp}}\left(t\right),y_{\text{comp}}\left(t\right)\right)\right).\label{eq:controlFinalSol}
\end{align}
To obtain a result consistent with the approximate solution for state
and co-state, Eq. \eqref{eq:controlFinalSol} is expanded up to leading
order in $\epsilon$. The outer limit, valid for times $t_{0}<t<t_{1}$,
yields the identities 
\begin{align}
\lim_{\epsilon\rightarrow0}Y_{L}\left(\left(t-t_{0}\right)/\epsilon\right) & =y_{\text{init}}, & \lim_{\epsilon\rightarrow0}Y_{R}\left(\left(t_{1}-t\right)/\epsilon\right) & =y_{\text{end}},\\
\lim_{\epsilon\rightarrow0}Y_{L}'\left(\left(t-t_{0}\right)/\epsilon\right) & =0, & \lim_{\epsilon\rightarrow0}Y_{R}'\left(\left(t_{1}-t\right)/\epsilon\right) & =0,
\end{align}
and the outer control signal as
\begin{align}
u_{O}\left(t\right) & =\lim_{\epsilon\rightarrow0}u\left(t\right)=\frac{1}{b\left(x_{O}\left(t\right),y_{O}\left(t\right)\right)}\left(\dot{y}_{O}\left(t\right)-R\left(x_{O}\left(t\right),y_{O}\left(t\right)\right)\right).\label{eq:OuterControl}
\end{align}
The inner limits of the control are defined by
\begin{align}
u\left(t_{0}+\epsilon\tau_{L}\right) & =U_{L}\left(\tau_{L}\right)+\text{h.o.t.},\\
u\left(t_{1}-\epsilon\tau_{R}\right) & =U_{R}\left(\tau_{R}\right)+\text{h.o.t.}.
\end{align}
The abbreviation h.o.t. stands for higher order terms which vanish
as $\epsilon\rightarrow0$. The left and right outer control signals
$U_{L}$ and $U_{R}$ depend on the rescaled times $\tau_{L}=\left(t-t_{0}\right)/\epsilon$
and $\tau_{R}=\left(t_{1}-t\right)/\epsilon$, respectively.

To compute $U_{L}$ and $U_{R}$, the derivative $\dot{y}_{\text{comp}}$
must be expanded in $\epsilon$, 
\begin{align}
\dot{y}_{\text{comp}}\left(t_{0}+\epsilon\tau_{L}\right) & =\dot{y}_{O}\left(t_{0}+\epsilon\tau_{L}\right)+\dfrac{1}{\epsilon}Y_{L}'\left(\tau_{L}\right)-\dfrac{1}{\epsilon}Y_{R}'\left(\left(t_{1}-t_{0}\right)/\epsilon-\tau_{L}\right)\nonumber \\
 & =\dot{y}_{O}\left(t_{0}\right)+\dfrac{1}{\epsilon}Y_{L}'\left(\tau_{L}\right)+\text{h.o.t.},\\
\dot{y}_{\text{comp}}\left(t_{1}-\epsilon\tau_{R}\right) & =\dot{y}_{O}\left(t_{1}-\epsilon\tau_{R}\right)+\dfrac{1}{\epsilon}Y_{L}'\left(\left(t_{1}-t_{0}\right)/\epsilon-\tau_{R}\right)-\dfrac{1}{\epsilon}Y_{R}'\left(\tau_{R}\right)\nonumber \\
 & =\dot{y}_{O}\left(t_{1}\right)-\dfrac{1}{\epsilon}Y_{R}'\left(\tau_{R}\right)+\text{h.o.t.}.
\end{align}
The expressions $\dfrac{1}{\epsilon}Y_{R}'\left(\left(t_{1}-t_{0}\right)/\epsilon-\tau_{L}\right)$
and $\dfrac{1}{\epsilon}Y_{L}'\left(\left(t_{1}-t_{0}\right)/\epsilon-\tau_{R}\right)$
are assumed to approach zero sufficiently fast as $\epsilon\rightarrow0$.
The terms proportional to $1/\epsilon$ diverge as $\epsilon\rightarrow0$.
This forbids a straightforward computation of $\lim_{\epsilon\rightarrow0}$,
and is the reason for the $\text{h.o.t.}$ notation. The left and
right inner limits of the control signal are obtained as 
\begin{align}
U_{L}\left(\tau_{L}\right) & =\frac{1}{b\left(x_{0},Y_{L}\left(\tau_{L}\right)\right)}\left(\dot{y}_{O}\left(t_{0}\right)+\dfrac{1}{\epsilon}Y_{L}'\left(\tau_{L}\right)-R\left(x_{0},Y_{L}\left(\tau_{L}\right)\right)\right)\nonumber \\
 & =\frac{1}{b\left(x_{0},Y_{L}\left(\tau_{L}\right)\right)}\left(\dot{y}_{O}\left(t_{0}\right)-R\left(x_{0},Y_{L}\left(\tau_{L}\right)\right)\right)\nonumber \\
 & +\dfrac{\sqrt{s_{2}}}{\epsilon}\text{sign}\left(b\left(x_{0},Y_{L}\left(\tau_{L}\right)\right)\right)\left(y_{\text{init}}-Y_{L}\left(\tau_{L}\right)\right),\label{eq:InnerLeftControl}\\
U_{R}\left(\tau_{R}\right) & =\frac{1}{b\left(x_{1},Y_{R}\left(\tau_{R}\right)\right)}\left(\dot{y}_{O}\left(t_{1}\right)-\dfrac{1}{\epsilon}Y_{R}'\left(\tau_{R}\right)-R\left(x_{1},Y_{R}\left(\tau_{R}\right)\right)\right)\nonumber \\
 & =\frac{1}{b\left(x_{1},Y_{R}\left(\tau_{R}\right)\right)}\left(\dot{y}_{O}\left(t_{1}\right)-R\left(x_{1},Y_{R}\left(\tau_{R}\right)\right)\right)\nonumber \\
 & -\dfrac{\sqrt{s_{2}}}{\epsilon}\text{sign}\left(b\left(x_{1},Y_{R}\left(\tau_{R}\right)\right)\right)\left(y_{\text{end}}-Y_{R}\left(\tau_{R}\right)\right).\label{eq:InnerRightControl}
\end{align}
Equation \eqref{eq:EqYL} is used to substitute $Y_{L}'\left(\tau_{L}\right)$,
and analogously for $Y_{R}'\left(\tau_{R}\right)$.

The left and right matching conditions
\begin{align}
u_{O}\left(t_{0}\right) & =\lim_{\tau_{L}\rightarrow\infty}U_{L}\left(\tau_{L}\right), & u_{O}\left(t_{1}\right) & =\lim_{\tau_{R}\rightarrow\infty}U_{R}\left(\tau_{R}\right),
\end{align}
are satisfied because $Y_{L}'\left(\tau_{L}\right)$ and $Y_{R}'\left(\tau_{R}\right)$
approach zero as $\tau_{L}\rightarrow\infty$ and $\tau_{R}\rightarrow\infty$.
The overlaps are obtained as
\begin{align}
u_{O}\left(t_{0}\right) & =\frac{1}{b\left(x_{O}\left(t_{0}\right),y_{O}\left(t_{0}\right)\right)}\left(\dot{y}_{O}\left(t_{0}\right)-R\left(x_{O}\left(t_{0}\right),y_{O}\left(t_{0}\right)\right)\right)\nonumber \\
 & =\frac{1}{b\left(x_{0},y_{\text{init}}\right)}\left(\dot{y}_{d}\left(t_{0}\right)+\frac{a_{2}s_{1}}{s_{2}}\left(x_{0}-x_{d}\left(t_{0}\right)\right)+a_{1}\left(y_{d}\left(t_{0}\right)-y_{\text{init}}\right)-R\left(x_{0},y_{\text{init}}\right)\right),\\
u_{O}\left(t_{1}\right) & =\frac{1}{b\left(x_{O}\left(t_{1}\right),y_{O}\left(t_{1}\right)\right)}\left(\dot{y}_{O}\left(t_{1}\right)-R\left(x_{O}\left(t_{1}\right),y_{O}\left(t_{1}\right)\right)\right)\nonumber \\
 & =\frac{1}{b\left(x_{1},y_{\text{end}}\right)}\left(\dot{y}_{d}\left(t_{1}\right)+\frac{a_{2}s_{1}}{s_{2}}\left(x_{1}-x_{d}\left(t_{1}\right)\right)+a_{1}\left(y_{d}\left(t_{1}\right)-y_{\text{end}}\right)-R\left(x_{1},y_{\text{end}}\right)\right).
\end{align}
The differential equation for $y_{O}$, Eq. \eqref{eq:yOEquation},
is used to eliminate $\dot{y}\left(t_{0}\right)$ and $\dot{y}\left(t_{1}\right)$.
Finally, the composite solution for the control signal is 
\begin{align}
u_{\text{comp}}\left(t\right) & =u_{O}\left(t\right)+U_{L}\left(\left(t-t_{0}\right)/\epsilon\right)+U_{R}\left(\left(t_{1}-t\right)/\epsilon\right)-u_{O}\left(t_{0}\right)-u_{O}\left(t_{1}\right).\label{eq:uFinalSol}
\end{align}
The outer control signal $u_{O}\left(t\right)$ is given by Eq. \eqref{eq:OuterControl},
while the left and right inner control signals $U_{L}$ and $U_{R}$
are given by Eqs. \eqref{eq:InnerLeftControl} and \eqref{eq:InnerRightControl},
respectively. The long explicit expression for Eq. \eqref{eq:uFinalSol}
is not written down explicitly.

\subsection{\label{sub:TheLimitEpsilon0}The limit \texorpdfstring{$\epsilon\rightarrow0$}{epsilon -> 0}}

For $\epsilon=0$, the analytical approximations derived in this section
become exact. The exact solution displays a discontinuous state trajectory
and a diverging control signal. The boundary layers of the state component
$y\left(t\right)$ degenerates to a jump located at the beginning
and the end of the time interval,
\begin{align}
\lim_{\epsilon\rightarrow0}y\left(t\right) & =\lim_{\epsilon\rightarrow0}y_{\text{comp}}\left(t\right)\nonumber \\
 & =y_{O}\left(t\right)+\lim_{\epsilon\rightarrow0}Y_{L}\left(\left(t-t_{0}\right)/\epsilon\right)-y_{\text{init}}+\lim_{\epsilon\rightarrow0}Y_{R}\left(\left(t_{1}-t\right)/\epsilon\right)-y_{\text{end}}\nonumber \\
 & =\begin{cases}
y_{O}\left(t_{0}\right)+Y_{L}\left(0\right)-y_{\text{init}}, & t=t_{0},\\
y_{O}\left(t\right), & t_{0}<t<t_{1},\\
y_{O}\left(t_{1}\right)+Y_{R}\left(0\right)-y_{\text{end}}, & t=t_{1},
\end{cases}\nonumber \\
 & =\begin{cases}
y_{0}, & t=t_{0},\\
y_{O}\left(t\right), & t_{0}<t<t_{1},\\
y_{1}, & t=t_{1}.
\end{cases}
\end{align}
The state component $x\left(t\right)$ as well as the co-state $\lambda_{x}\left(t\right)$
do not exhibit boundary layers. Their solutions are continuous also
for $\epsilon=0$ and simply given by the outer solutions Eqs. \eqref{eq:xFinalSol}
and \eqref{eq:lambdaFinalSol}. The remaining co-state $\lambda_{y}\left(t\right)$
vanishes identically, $\lambda_{y}\left(t\right)=0$. Although $y_{O}\left(t\right)$
depends on the matching constants $y_{\text{init}}$ and $y_{\text{end}}$,
both constants are given solely in terms of the outer solutions and
the initial and terminal conditions. Thus, to determine the height
and the position of the jumps, it is not necessary to know any details
about the dynamics of the boundary layers. For $\epsilon=0$, no trace
of the boundary layers is left in the composite solution except for
the mere existence of the jumps. In particular, while the form of
the boundary layers depends on the specific choice of the coupling
function $b\left(x,y\right)$, the solution becomes independent of
the coupling function $b\left(x,y\right)$ for $\epsilon=0$. Thus,
for $\epsilon=0$, both possible sources of nonlinear system dynamics,
the nonlinearity $R\left(x,y\right)$ and the coupling function $b\left(x,y\right)$,
do entirely disappear from the solution for the controlled state trajectory.

To obtain the control signal for $\epsilon=0$ from Eq. \eqref{eq:uFinalSol},
the expression $Y_{L}'\left(\left(t-t_{0}\right)/\epsilon\right)/\epsilon$
must be analyzed in the limit $\epsilon\rightarrow0$. To that end,
let the function $\delta_{\epsilon}\left(t\right)$ be defined as
\begin{align}
\delta_{\epsilon}\left(t\right) & =\begin{cases}
\dfrac{1}{2\epsilon}\dfrac{1}{\left(y_{\text{init}}-y_{0}\right)}Y_{L}'\left(t/\epsilon\right), & t\geq0,\\
\dfrac{1}{2\epsilon}\dfrac{1}{\left(y_{\text{init}}-y_{0}\right)}Y_{L}'\left(-t/\epsilon\right), & t<0,
\end{cases}\nonumber \\
 & =\dfrac{1}{2\epsilon}\dfrac{1}{\left(y_{\text{init}}-y_{0}\right)}Y_{L}'\left(\left|t\right|/\epsilon\right).
\end{align}
Note that $\delta_{\epsilon}\left(t\right)$ is continuous across
$t=0$ for all $\epsilon>0$. A few computations prove that $\delta_{\epsilon}\left(t\right)$
is a representation of the Dirac delta function as $\epsilon\rightarrow0$,
\begin{align}
\lim_{\epsilon\rightarrow0}\delta_{\epsilon}\left(t\right) & =\delta\left(t\right).\label{eq:Eq622}
\end{align}
Indeed, from the differential equation for $Y_{L}$, Eq. \eqref{eq:EqYL}
follows 
\begin{align}
\lim_{\epsilon\rightarrow0}\dfrac{1}{\epsilon}Y_{L}'\left(0\right) & =\lim_{\epsilon\rightarrow0}\dfrac{1}{\epsilon}\sqrt{s_{2}}\left(y_{\text{init}}-Y_{L}\left(0\right)\right)\left|b\left(x_{0},Y_{L}\left(0\right)\right)\right|=\text{sign}\left(y_{\text{init}}-y_{0}\right)\infty,
\end{align}
and therefore
\begin{align}
\lim_{\epsilon\rightarrow0}\delta_{\epsilon}\left(t\right) & =\begin{cases}
0, & \left|t\right|>0,\\
\infty, & t=0.
\end{cases}
\end{align}
It remains to show that
\begin{align}
\intop_{-\infty}^{\infty}dt\delta_{\epsilon}\left(t\right) & =1
\end{align}
for all $\epsilon$. Together with the substitutions $t_{1}=-\epsilon\tau_{L}$
and $t_{2}=\epsilon\tau_{L}$ and the initial and matching condition
for $Y_{L}$, the integral over $\delta_{\epsilon}\left(t\right)$
yields
\begin{align}
\intop_{-\infty}^{\infty}dt\delta_{\epsilon}\left(t\right) & =\intop_{-\infty}^{0}dt_{1}\delta_{\epsilon}\left(t_{1}\right)+\intop_{0}^{\infty}dt_{2}\delta_{\epsilon}\left(t_{2}\right)\nonumber \\
 & =\dfrac{1}{2\epsilon}\dfrac{1}{\left(y_{\text{init}}-y_{0}\right)}\left(\intop_{-\infty}^{0}dt_{1}Y_{L}'\left(-t_{1}/\epsilon\right)+\intop_{0}^{\infty}dt_{2}Y_{L}'\left(t_{2}/\epsilon\right)\right)\nonumber \\
 & =\dfrac{1}{2}\dfrac{1}{\left(y_{\text{init}}-y_{0}\right)}\left(\intop_{0}^{\infty}d\tau_{L}Y_{L}'\left(\tau_{L}\right)+\intop_{0}^{\infty}d\tau_{L}Y_{L}'\left(\tau_{L}\right)\right)=1.\label{eq:Eq4118}
\end{align}
The result Eq. \eqref{eq:Eq4118} is independent of the value of $\epsilon$.
Thus, Eq. \eqref{eq:Eq622} can be used to establish the identity
\begin{align}
\lim_{\epsilon\rightarrow0}\dfrac{1}{\epsilon}Y_{L}'\left(\left(t-t_{0}\right)/\epsilon\right) & =2\left(y_{\text{init}}-y_{0}\right)\delta\left(t-t_{0}\right),\,t\geq t_{0}.\label{eq:Eq4119}
\end{align}
An analogous procedure applied to $Y_{R}'\left(\tau_{R}\right)$ yields
the analogous result 
\begin{align}
\lim_{\epsilon\rightarrow0}\dfrac{1}{\epsilon}Y_{R}'\left(\left(t_{1}-t\right)/\epsilon\right) & =2\left(y_{\text{end}}-y_{1}\right)\delta\left(t_{1}-t\right),\,t_{1}\geq t.\label{eq:Eq4120}
\end{align}
Equations \eqref{eq:Eq4119} and \eqref{eq:Eq4120} are used to compute
the limits
\begin{align}
 & \lim_{\epsilon\rightarrow0}U_{L}\left(\left(t-t_{0}\right)/\epsilon\right)\nonumber \\
= & \lim_{\epsilon\rightarrow0}\frac{1}{b\left(x_{0},Y_{L}\left(\left(t-t_{0}\right)/\epsilon\right)\right)}\left(\dot{y}_{O}\left(t_{0}\right)+\dfrac{1}{\epsilon}Y_{L}'\left(\left(t-t_{0}\right)/\epsilon\right)-R\left(x_{0},Y_{L}\left(\left(t-t_{0}\right)/\epsilon\right)\right)\right)\nonumber \\
= & \begin{cases}
\frac{1}{b\left(x_{0},y_{\text{init}}\right)}\left(\dot{y}_{O}\left(t_{0}\right)-R\left(x_{0},y_{\text{init}}\right)\right)=u_{O}\left(t_{0}\right), & t_{0}<t\leq t_{1},\\
\frac{1}{b\left(x_{0},y_{0}\right)}\left(\dot{y}_{O}\left(t_{0}\right)+2\left(y_{\text{init}}-y_{0}\right)\delta\left(t-t_{0}\right)-R\left(x_{0},y_{0}\right)\right), & t=t_{0},
\end{cases}
\end{align}
and
\begin{align}
 & \lim_{\epsilon\rightarrow0}U_{R}\left(\left(t_{1}-t\right)/\epsilon\right)\nonumber \\
= & \lim_{\epsilon\rightarrow0}\frac{1}{b\left(x_{1},Y_{R}\left(\left(t_{1}-t\right)/\epsilon\right)\right)}\left(\dot{y}_{O}\left(t_{1}\right)-\dfrac{1}{\epsilon}Y_{R}'\left(\left(t_{1}-t\right)/\epsilon\right)-R\left(x_{1},Y_{R}\left(\left(t_{1}-t\right)/\epsilon\right)\right)\right)\nonumber \\
 & =\begin{cases}
\frac{1}{b\left(x_{1},y_{\text{end}}\right)}\left(\dot{y}_{O}\left(t_{1}\right)-R\left(x_{1},y_{\text{end}}\right)\right)=u_{O}\left(t_{1}\right), & t_{0}\leq t<t_{1},\\
\frac{1}{b\left(x_{1},y_{1}\right)}\left(\dot{y}_{O}\left(t_{1}\right)-2\left(y_{\text{end}}-y_{1}\right)\delta\left(t_{1}-t\right)-R\left(x_{1},y_{1}\right)\right), & t=t_{1},
\end{cases}
\end{align}
respectively. Finally, the exact solution for the control signal is
\begin{align}
\lim_{\epsilon\rightarrow0}u\left(t\right) & =\lim_{\epsilon\rightarrow0}u_{\text{comp}}\left(t\right)\nonumber \\
 & =u_{O}\left(t\right)+\lim_{\epsilon\rightarrow0}U_{L}\left(\left(t-t_{0}\right)/\epsilon\right)-u_{O}\left(t_{0}\right)+\lim_{\epsilon\rightarrow0}U_{R}\left(\left(t_{1}-t\right)/\epsilon\right)-u_{O}\left(t_{1}\right)\nonumber \\
 & =\begin{cases}
\frac{1}{b\left(x_{0},y_{0}\right)}\left(\dot{y}_{O}\left(t_{0}\right)+2\left(y_{\text{init}}-y_{0}\right)\delta\left(t-t_{0}\right)-R\left(x_{0},y_{0}\right)\right), & t=t_{0},\\
\frac{1}{b\left(x_{O}\left(t\right),y_{O}\left(t\right)\right)}\left(\dot{y}_{O}\left(t\right)-R\left(x_{O}\left(t\right),y_{O}\left(t\right)\right)\right), & t_{0}<t<t_{1},\\
\frac{1}{b\left(x_{1},y_{1}\right)}\left(\dot{y}_{O}\left(t_{1}\right)-2\left(y_{\text{end}}-y_{1}\right)\delta\left(t_{1}-t\right)-R\left(x_{1},y_{1}\right)\right), & t=t_{1}.
\end{cases}
\end{align}
In contrast to the controlled state trajectory, the solution for the
control signal depends on both possible nonlinearities, the coupling
function $b\left(x,y\right)$ and the nonlinearity $R\left(x,y\right)$.

A discussion of the results can be found at the end of the next section.

\section{\label{sec:ComparisonWithNumericalResults}Comparison with numerical
results}

\subsection{Results}

Numerical computations are performed with the \href{http://www.acadotoolkit.org}{ACADO}
Toolkit \cite{Houska2011a,Houska2011b}, an open source program package
for solving optimal control problems. Typically, a problem is solved
on a time interval of length $1$ with a time step width of $\Delta t=10^{-3}$.
The numerical computation of such an example takes about 20-30min
with a standard Laptop. Computation time increases quickly with decreasing
step width or increasing length of the time interval. For comparison
with analytical solutions, the numerical result provided by ACADO
is imported in Mathematica \cite{wolfram2014mathematica} and interpolated.
Two example systems are investigated. Both are covered by the general
analytical result from Section \ref{sec:TwoDimensionalDynamicalSystem}.
The activator-controlled FHN model in is discussed in Example \ref{ex:OptimallyControlledFHN}.
Example \ref{ex:ControlledPendulum} presents results for a mechanical
system, namely the damped mathematical pendulum.

\begin{example}[Activator-controlled FHN model]\label{ex:OptimallyControlledFHN}

See Example \ref{ex:FHN1} for an introduction to the model and parameter
values. The state equations are repeated here for convenience,
\begin{align}
\left(\begin{array}{c}
\dot{x}\left(t\right)\\
\dot{y}\left(t\right)
\end{array}\right) & =\left(\begin{array}{c}
a_{0}+a_{1}x\left(t\right)+a_{2}y\left(t\right)\\
y\left(t\right)-\frac{1}{3}y\left(t\right)^{3}-x\left(t\right)
\end{array}\right)+\left(\begin{array}{c}
0\\
b\left(x\left(t\right),y\left(t\right)\right)
\end{array}\right)u\left(t\right).
\end{align}
In contrast to Example \ref{ex:ControlledFHN1} with $b\left(x,y\right)=1$,
here a state-dependent coupling function $b$ is assumed,
\begin{align}
b\left(x,y\right) & =\dfrac{11}{4}+x^{2}.\label{eq:CouplingFunctionFHN}
\end{align}
The small regularization parameter $\epsilon$ is set to 
\begin{align}
\epsilon & =10^{-3},
\end{align}
which results in a regularization term with coefficient $\frac{1}{2}\epsilon^{2}=0.5\times10^{-6}$,
see Eq. \eqref{eq:Functional}. The desired reference trajectory is
an ellipse,
\begin{align}
x_{d}\left(t\right) & =A_{x}\cos\left(2\pi t/T\right)-\dfrac{1}{2}, & y_{d}\left(t\right) & =A_{y}\sin\left(2\pi t/T\right)+\dfrac{1}{2},\label{eq:DesiredTrajectoryFHN}
\end{align}
with $A_{x}=1,\,A_{y}=15$, and $T=1$. Within the time interval $0=t_{0}\leq t<t_{1}=1$,
the controlled state trajectory shall follow the ellipse as closely
as possible. The initial and terminal state lie on the desired trajectory,
\begin{align}
x\left(t_{0}\right) & =x_{d}\left(t_{0}\right)=\dfrac{1}{2}, & y\left(t_{0}\right) & =y_{d}\left(t_{1}\right)=\dfrac{15}{2},\label{eq:NumInit}\\
x\left(t_{1}\right) & =x_{d}\left(t_{1}\right)=\dfrac{1}{2}, & y\left(t_{1}\right) & =y_{d}\left(t_{1}\right)=\dfrac{15}{2}.\label{eq:NumTerm}
\end{align}
Figure \ref{fig:Fig1} compares the prescribed desired trajectory
$\boldsymbol{x}_{d}\left(t\right)$ as given by Eq. \eqref{eq:DesiredTrajectoryFHN}
(blue solid line) with the numerically obtained optimally controlled
state trajectory $\boldsymbol{x}_{num}\left(t\right)$ (red dashed
line). While the controlled activator (Fig. \ref{fig:Fig1} right)
looks similar to the desired trajectory except for a constant shift,
the controlled inhibitor (Fig. \ref{fig:Fig1} right) is way off.
Although the initial and terminal conditions comply with the desired
trajectory, the solution for the activator component $y\left(t\right)$
exhibits some very steep transients at the beginning and end of the
time interval. These transients can be interpreted as boundary layers
described by the inner solutions. While this example differs from
Example \ref{ex:ControlledFHN1} in the coupling function $b\left(x,y\right)$,
its controlled state trajectories are very similar. Indeed, the analytical
solutions predicts an effect of $b\left(x,y\right)$ restricted to
the boundary layer region of the activator component $y\left(t\right)$.
Consequently, the solutions for different coupling functions $b\left(x,y\right)$
are essentially identical inside the time domain.

\begin{minipage}{1.0\linewidth}
\begin{center}
\includegraphics[scale=0.49]{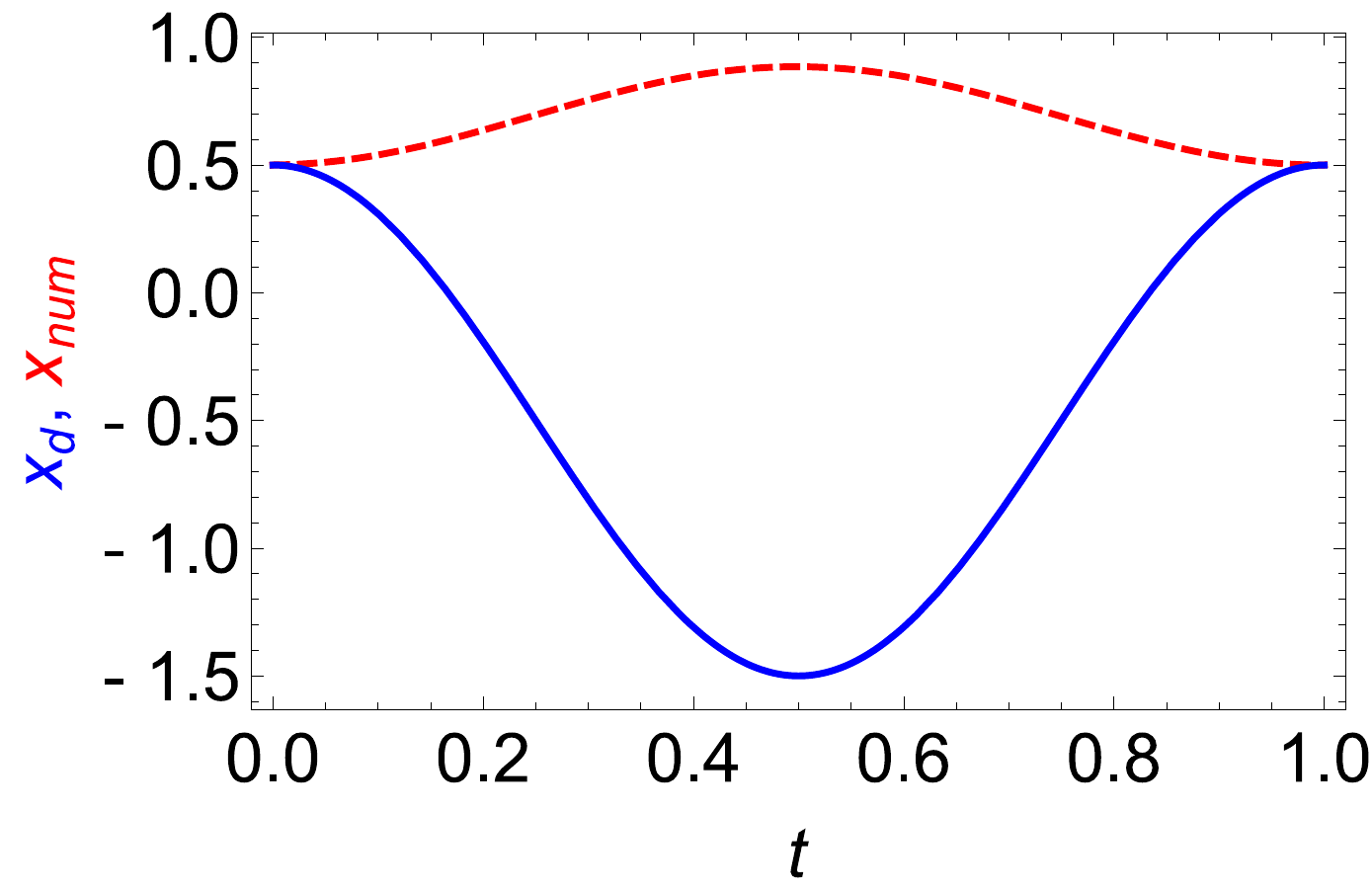}\hspace{0.2cm}\includegraphics[scale=0.485]{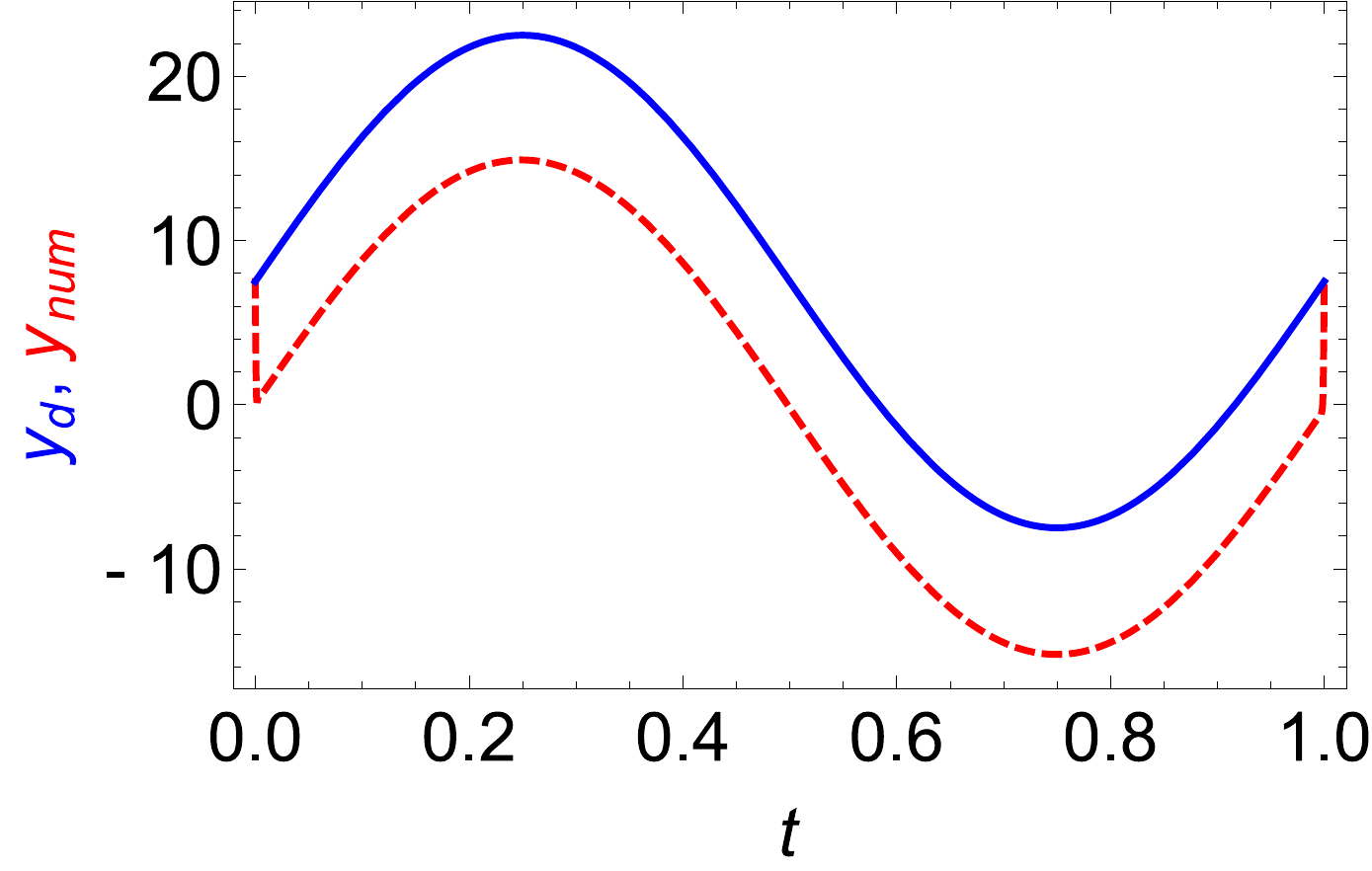}
\captionof{figure}[Comparison of desired and optimal state trajectory in FHN model]{\label{fig:Fig1}Comparison of desired reference trajectory $\boldsymbol{x}_{d}\left(t\right)$ (blue solid line) and numerically obtained optimal trajectory $\boldsymbol{x}_{num}\left(t\right)$ (red dashed line) in the FHN model. The activator over time $y$ (right) is similar to the reference trajectory in shape but shifted by an almost constant value, while the inhibitor over time $x$ (left) is far off.}
%"/home/jakob/ACADOtoolkit/examples/my_examples/FHN/CompareResults_FHN.nb"
\end{center}
\end{minipage}

Figure \ref{fig:Fig2} compares the analytical result for $\boldsymbol{x}\left(t\right)$
from Section \ref{sec:TwoDimensionalDynamicalSystem} with its numerical
counterpart. The agreement is excellent. Although only approximately
valid, this demonstrates an astonishing accuracy of the analytical
result for small values of $\epsilon$. No difference between analytical
and numerical is visible on this scale. Figure \ref{fig:Fig3} visualizes
$\boldsymbol{x}\left(t\right)-\boldsymbol{x}_{num}\left(t\right)$
and reveals relatively small differences in the bulk but somewhat
larger differences close to the initial and terminal time, especially
for the state component $y\left(t\right)$.

\begin{minipage}{1.0\linewidth}
\begin{center}
\includegraphics[scale=0.49]{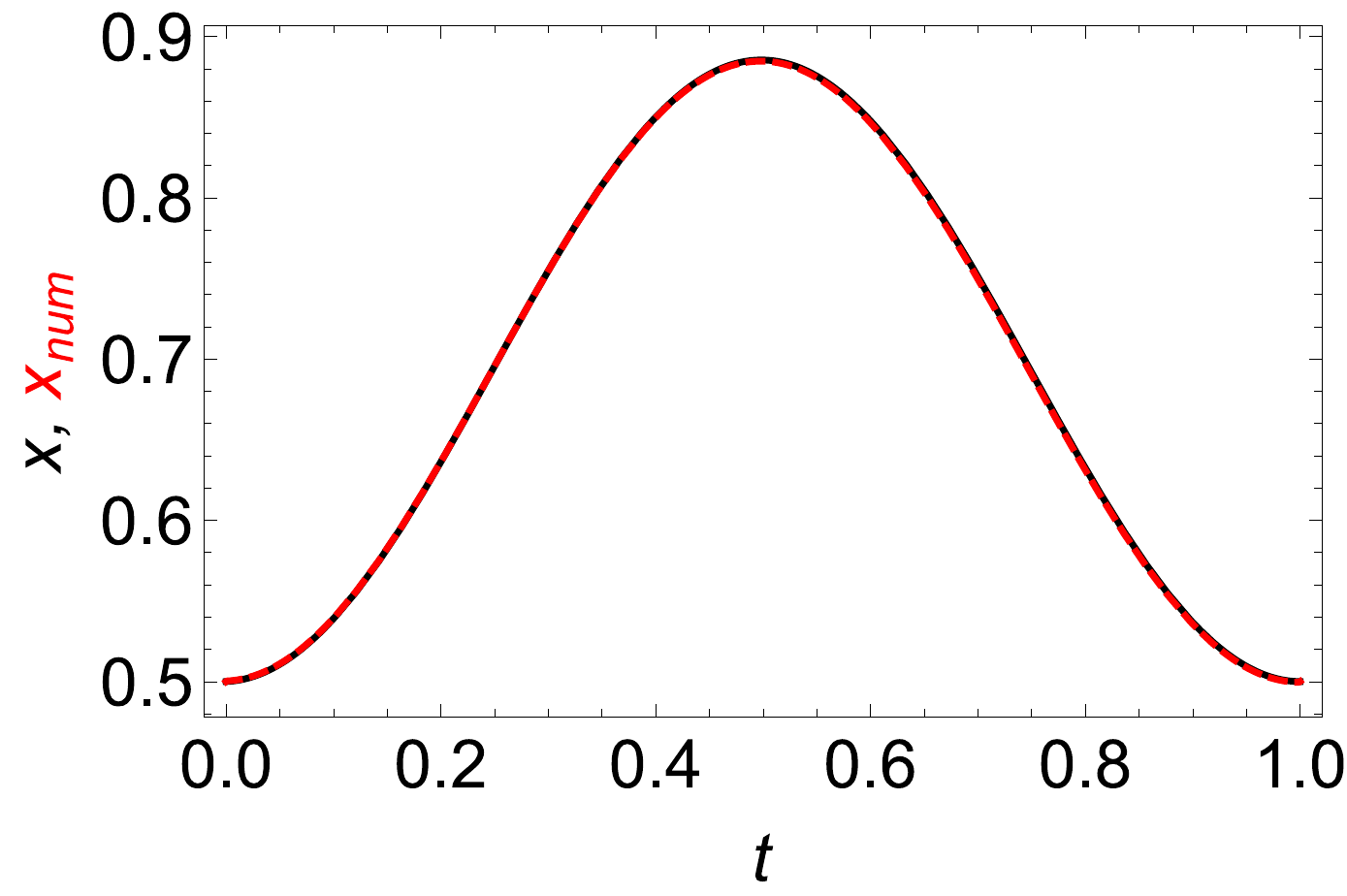}\hspace{0.2cm}\includegraphics[scale=0.4975]{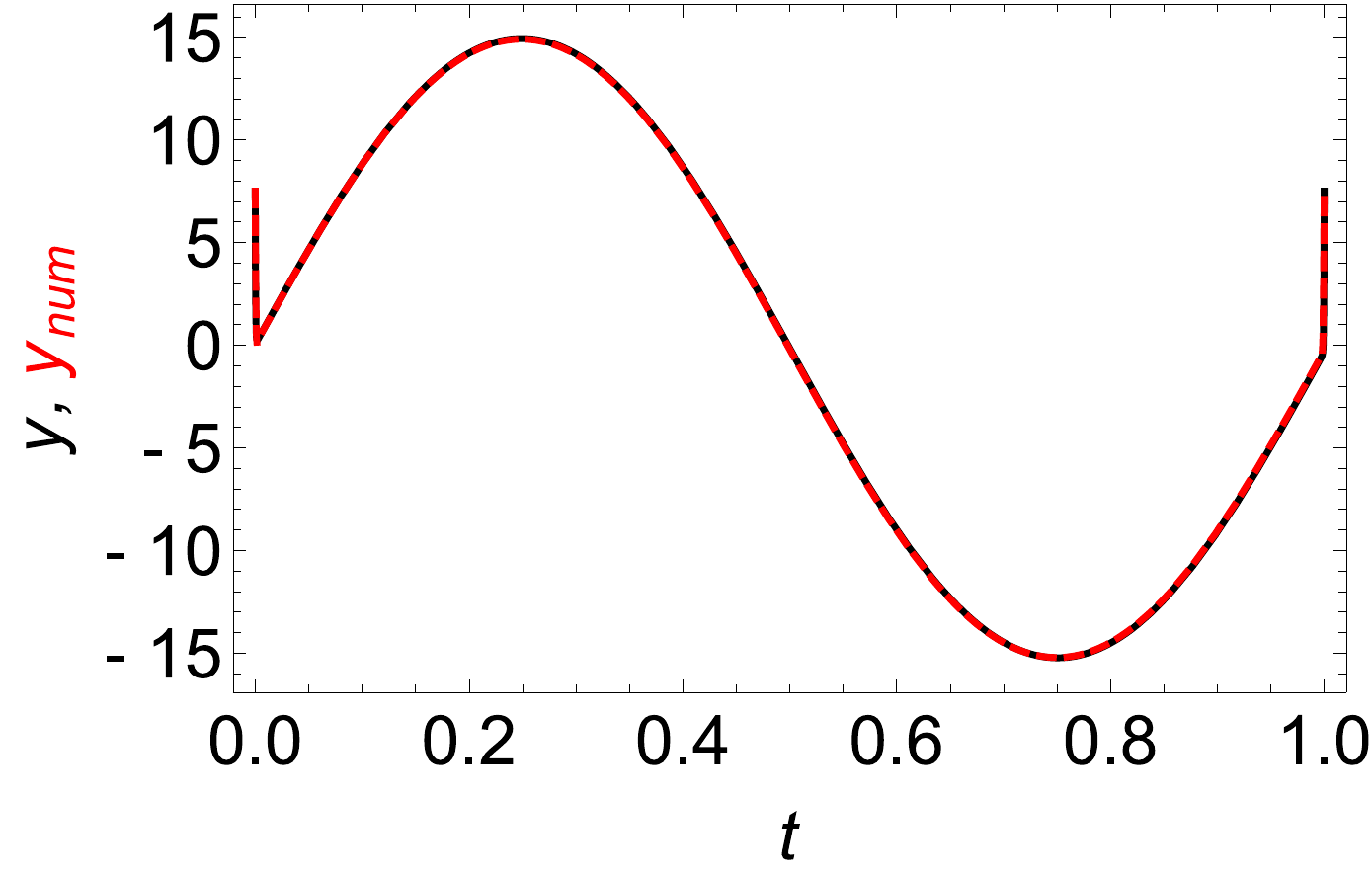}
\captionof{figure}[Comparison of numerical and analytical optimal trajectory]{\label{fig:Fig2}Comparison of numerically obtained optimal trajectory (red dashed line) and analytical approximation (black solid line). On this scale, the agreement is perfect.}
%"/home/jakob/ACADOtoolkit/examples/my_examples/FHN/CompareResults_FHN.nb"
\end{center}
\end{minipage}

\begin{minipage}{1.0\linewidth}
\begin{center}
\includegraphics[scale=0.5]{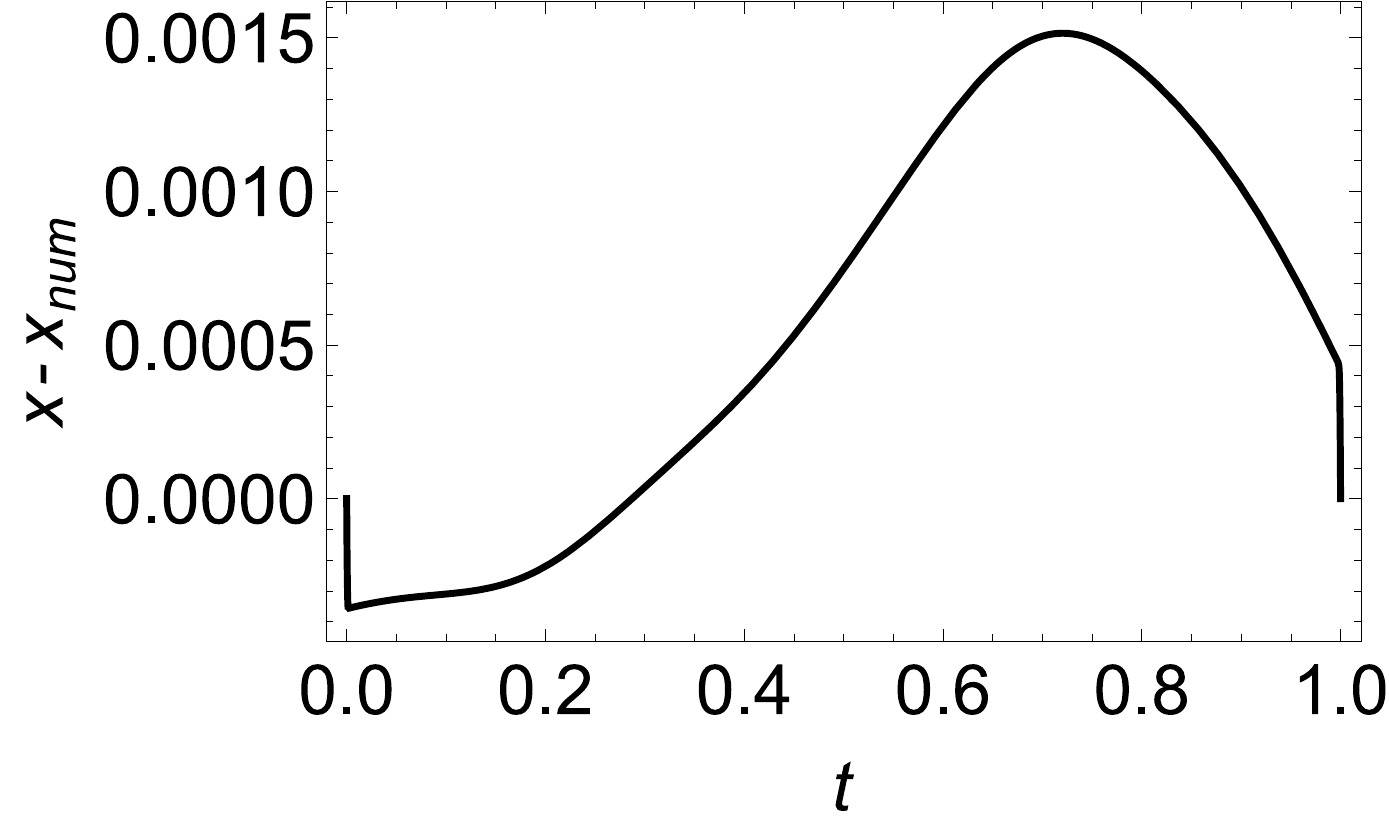}\hspace{0.2cm}\includegraphics[scale=0.475]{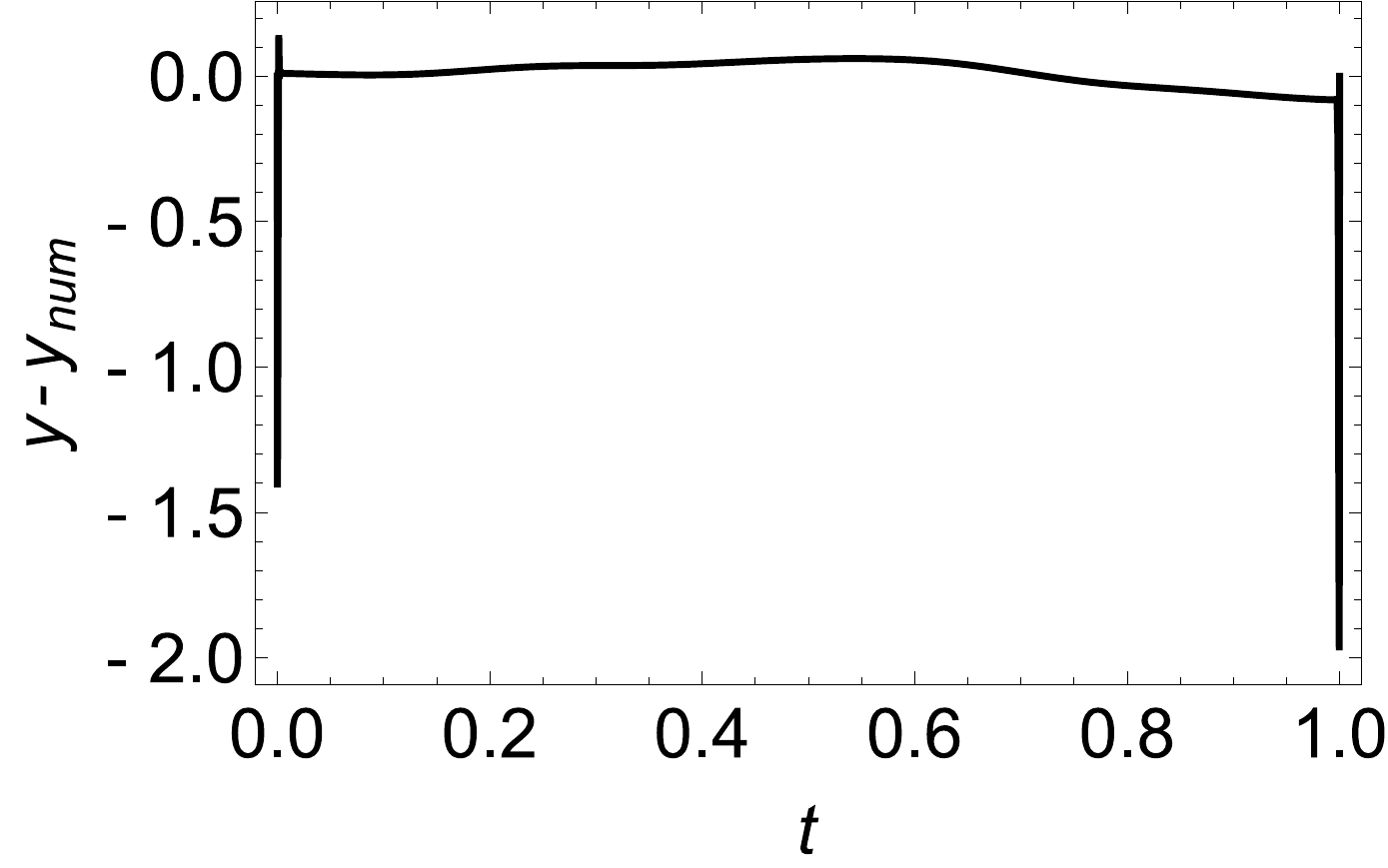}
\captionof{figure}[Difference between analytical and numerical solution]{\label{fig:Fig3}Difference between analytical and numerical solution for inhibitor (left) and activator component (right) over time.}
%"/home/jakob/ACADOtoolkit/examples/my_examples/FHN/CompareResults_FHN.nb"
\end{center}
\end{minipage}

Zooming in on the initial (Fig. \ref{fig:Fig4}) and the terminal
time (Fig. \ref{fig:Fig5}) uncovers the boundary layers. The activator
(right) displays small deviations in the regions with the steepest
slopes. This is certainly due to the limited temporal resolution of
the numerical simulation. Note that the width of the boundary layers
is approximately determined by the value of the regularization parameter
$\epsilon$. The value $\epsilon=10^{-3}$ chosen for numerical simulations
is identical to the temporal resolution of $\Delta t=10^{-3}$. A
result is the relatively large difference between analytical and numerical
result at the initial and terminal times in Fig. \ref{fig:Fig3}.
For values of $\epsilon$ in the range of the temporal resolution,
the boundary layers cannot be resolved numerically with sufficient
accuracy and result in discretization errors. Due to the computational
cost of optimal control algorithms, decreasing the step width $\Delta t$
is not really an option. Figure \ref{fig:Fig4} and \ref{fig:Fig5}
left reveals a small boundary layer displayed only by the numerically
obtained inhibitor component. This is not predicted by the analytical
leading order approximation and results probably from higher order
contributions of the perturbation expansion. Finally, note that the
deviations between analytical and numerical result are slightly larger
close to the terminal time (Fig. \ref{fig:Fig5}) than to the initial
time (Fig. \ref{fig:Fig4}). This hints at the accumulation of numerical
errors in the numerical result.

\begin{minipage}{1.0\linewidth}
\begin{center}
\includegraphics[scale=0.51]{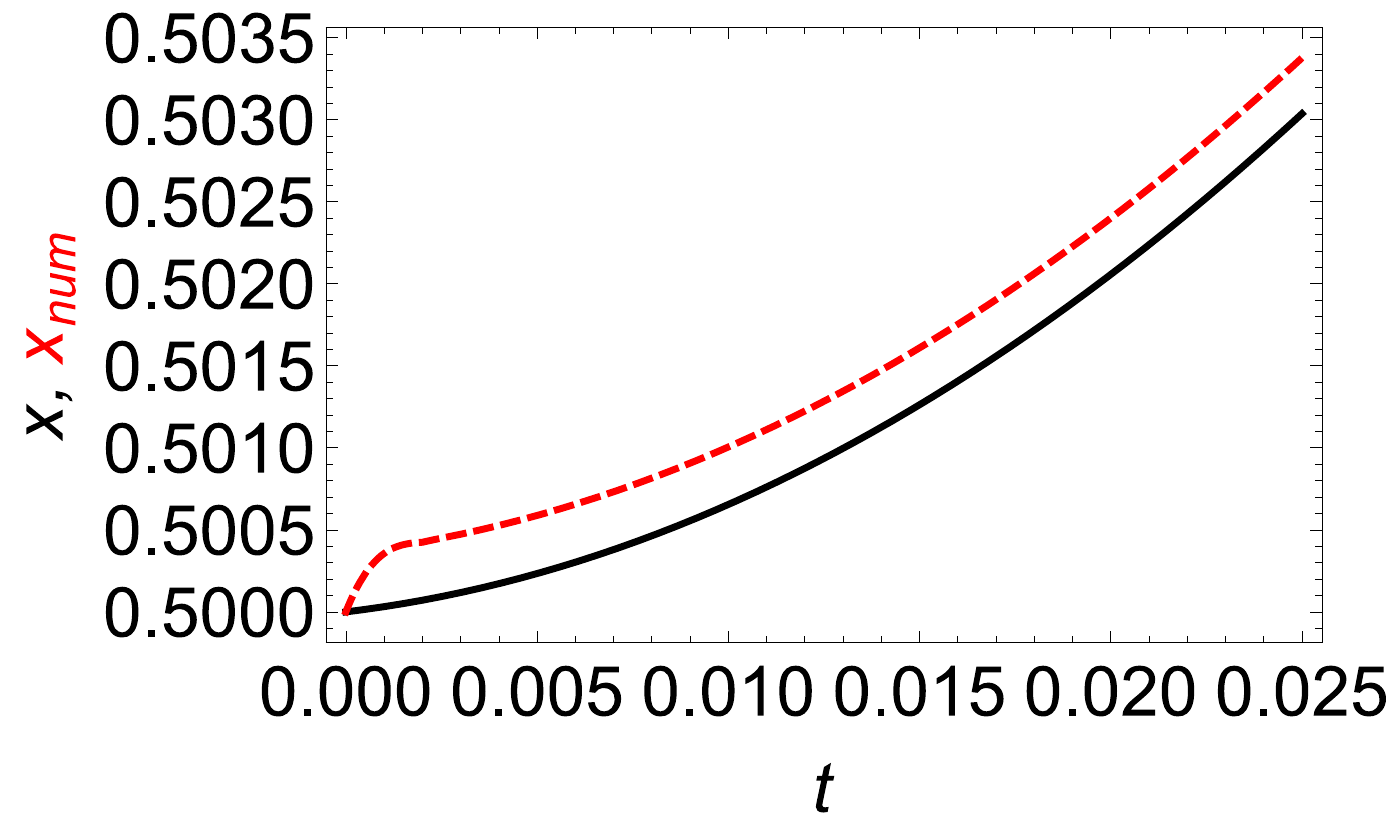}\hspace{0.2cm}\includegraphics[scale=0.45]{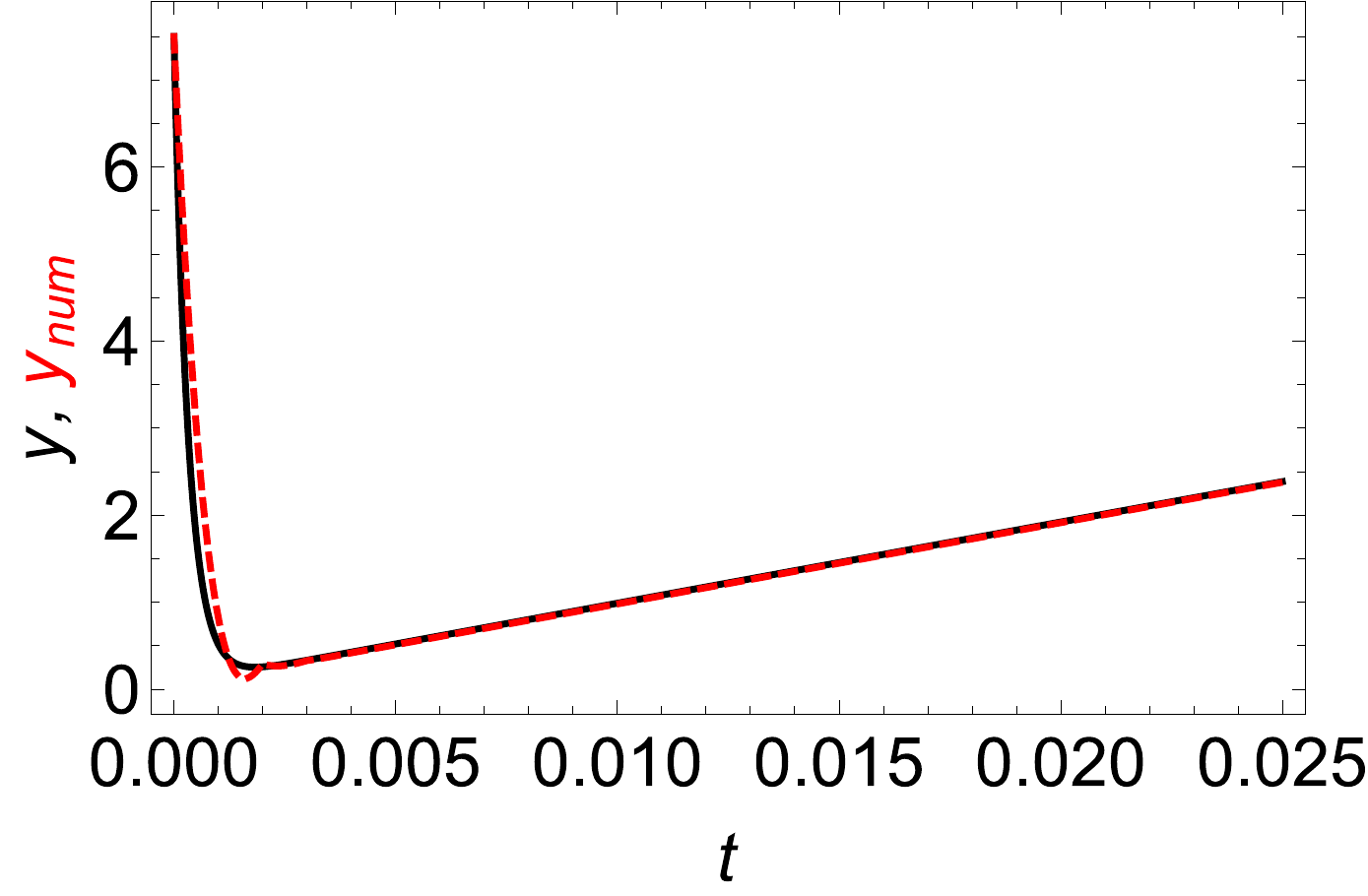}
\captionof{figure}[Closeup of the left boundary layer for the activator]{\label{fig:Fig4}Closeup of the left boundary layer for the activator component (right) shows perfect agreement between numerical (red dashed line) and analytical (black line) result except for the steepest slopes. The numerically obtained inhibitor component exhibits a very small boundary layer as well (left), while the leading order analytical result does not. Analytically, this boundary layer arises probably from contributions of higher order in $\epsilon$.}
%"/home/jakob/ACADOtoolkit/examples/my_examples/FHN/CompareResults_FHN.nb"
\end{center}
\end{minipage}

\begin{minipage}{1.0\linewidth}
\begin{center}
\includegraphics[scale=0.5]{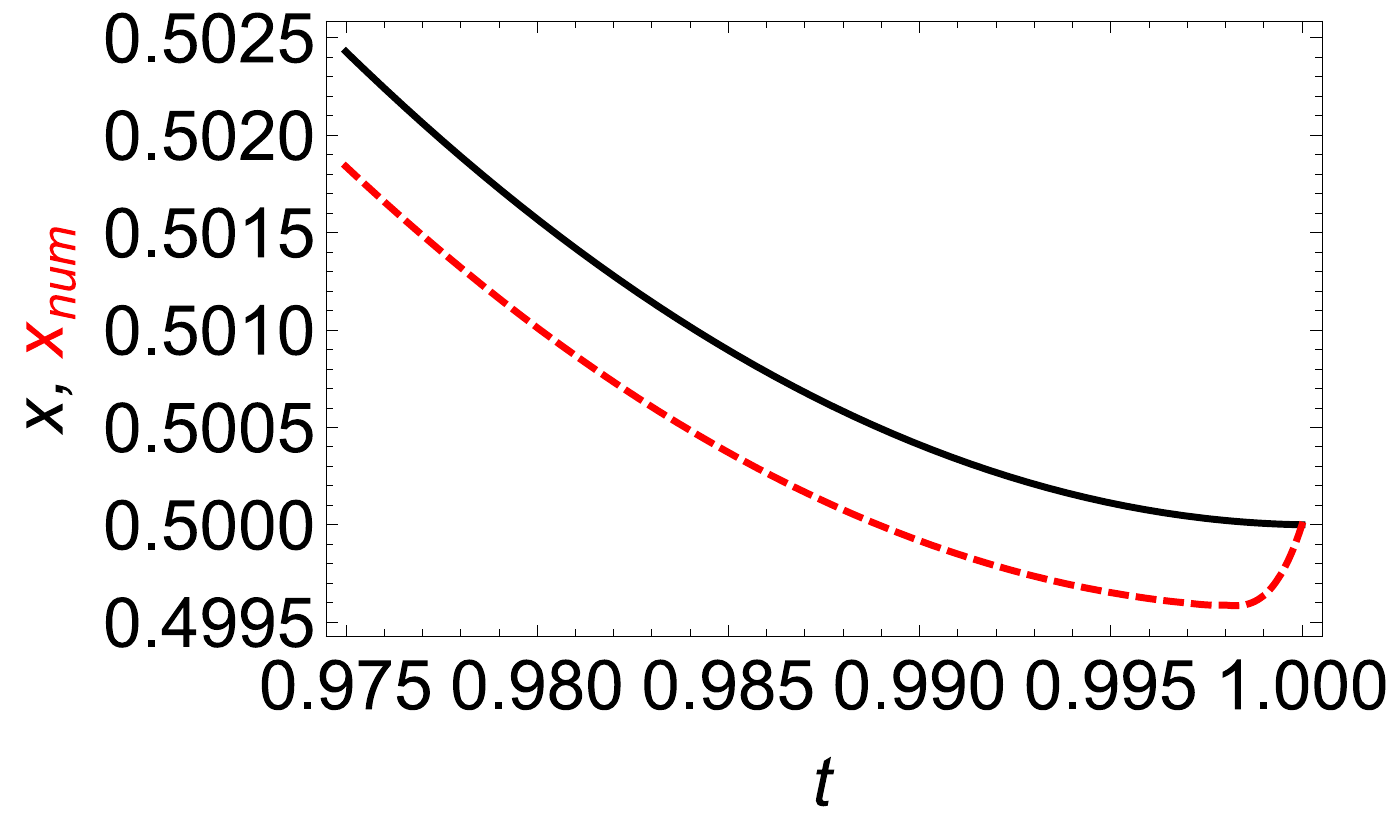}\hspace{0.2cm}\includegraphics[scale=0.455]{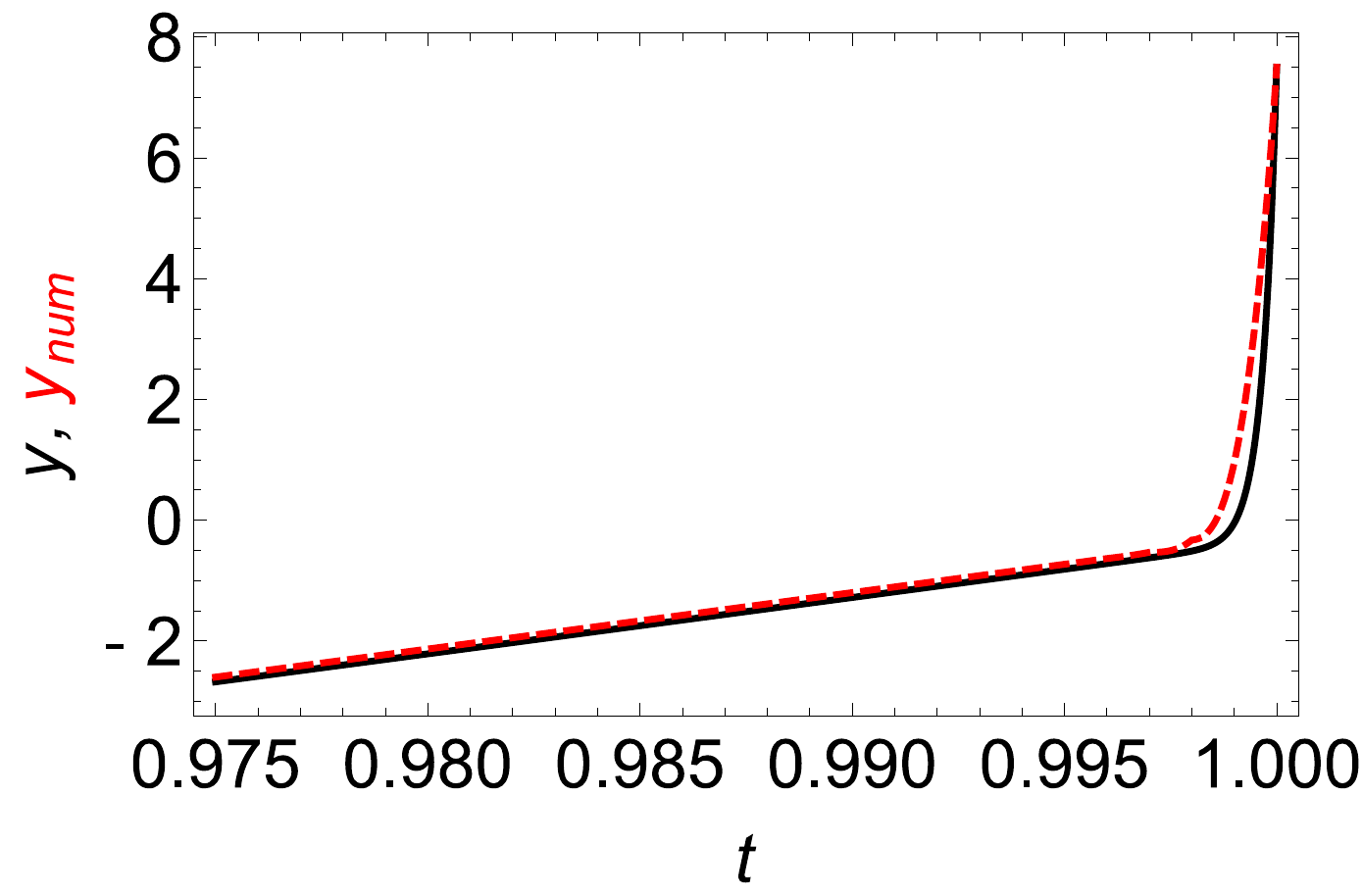}
\captionof{figure}[Closeup of the right boundary layer for the activator]{\label{fig:Fig5}Closeup of the right boundary layer for the activator. The agreement is worse than for the left boundary layer. This hints at the accumulation of numerical errors in the numerical result.}
%"/home/jakob/ACADOtoolkit/examples/my_examples/FHN/CompareResults_FHN.nb"
\end{center}
\end{minipage}

Figure \ref{fig:Fig6} concludes with a comparison between analytical
and numerical solution for the control. The control attains its largest
values at the initial and terminal time. Analytically, these spikes
approach the form of a Dirac delta distribution for a decreasing values
of $\epsilon$.

\begin{minipage}{1.0\linewidth}
\begin{center}
\includegraphics[scale=0.5]{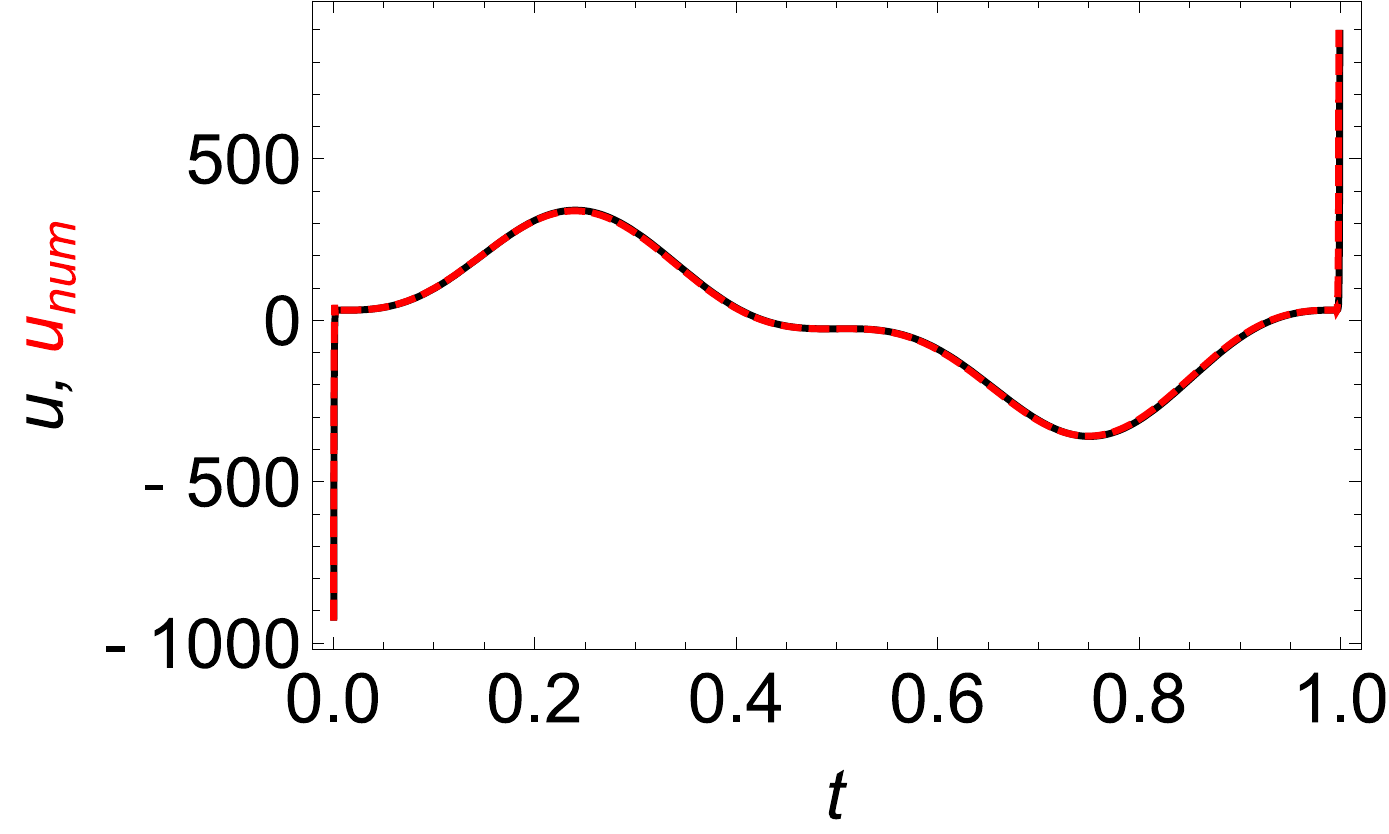}
\captionof{figure}[Comparison of numerical and analytical control for the FHN model]{\label{fig:Fig6}Comparison of the numerically obtained control $u$ (red dashed line) and its analytical approximation (black solid line). The numerical and analytical control solutions attain large values at the boundaries of the time domain due to the appearance of boundary layers in the state component $y\left( t\right)$.}
%"/home/jakob/ACADOtoolkit/examples/my_examples/FHN/CompareResults_FHN.nb"
\end{center}
\end{minipage}

\end{example}

\begin{example}[Mathematical Pendulum]\label{ex:ControlledPendulum}

The mathematical pendulum is an example for a mechanical control system,
see Example \ref{ex:OneDimMechSys1}. The controlled state equation
is
\begin{align}
\dot{x}\left(t\right) & =y\left(t\right),\\
\dot{y}\left(t\right) & =-\gamma y\left(t\right)-\sin\left(x\left(t\right)\right)+u\left(t\right).
\end{align}
For mechanical systems, the general analytical result from Section
\ref{sec:TwoDimensionalDynamicalSystem} simplifies considerably due
to the fixed parameter values 
\begin{align}
a_{0} & =0, & a_{1} & =0, & a_{2} & =1.
\end{align}
The desired trajectory $\boldsymbol{x}_{d}\left(t\right)$ is
\begin{align}
x_{d}\left(t\right) & =\cos\left(2\pi t\right), & y_{d}\left(t\right) & =\cos\left(2\pi t\right)+\sin\left(4\pi t\right).\label{eq:PendulumDesiredTrajectory}
\end{align}
The initial and terminal conditions
\begin{align}
x\left(t_{0}\right) & =-1, & y\left(t_{0}\right) & =-1, & x\left(t_{1}\right) & =-1, & y\left(t_{1}\right) & =-1,\label{eq:InitTermPendulum}
\end{align}
do not comply with the desired trajectory. As before, the small parameter
$\epsilon$ is
\begin{align}
\epsilon & =10^{-3}.
\end{align}
Figure \ref{fig:FigP1} compares the desired trajectory $\boldsymbol{x}_{d}\left(t\right)$
with the numerically obtained optimally controlled state trajectory
$\boldsymbol{x}_{num}\left(t\right)$. The velocity $y\left(t\right)$
is much closer to its desired counterpart than the position over time
$x\left(t\right)$. Initial and terminal boundary layers occur for
the velocity $y\left(t\right)$. Figure \ref{fig:FigP2} compares
the corresponding analytical result with the numerical solution and
reveals almost perfect agreement.

\begin{minipage}{1.0\linewidth}
\begin{center}
\includegraphics[scale=0.49]{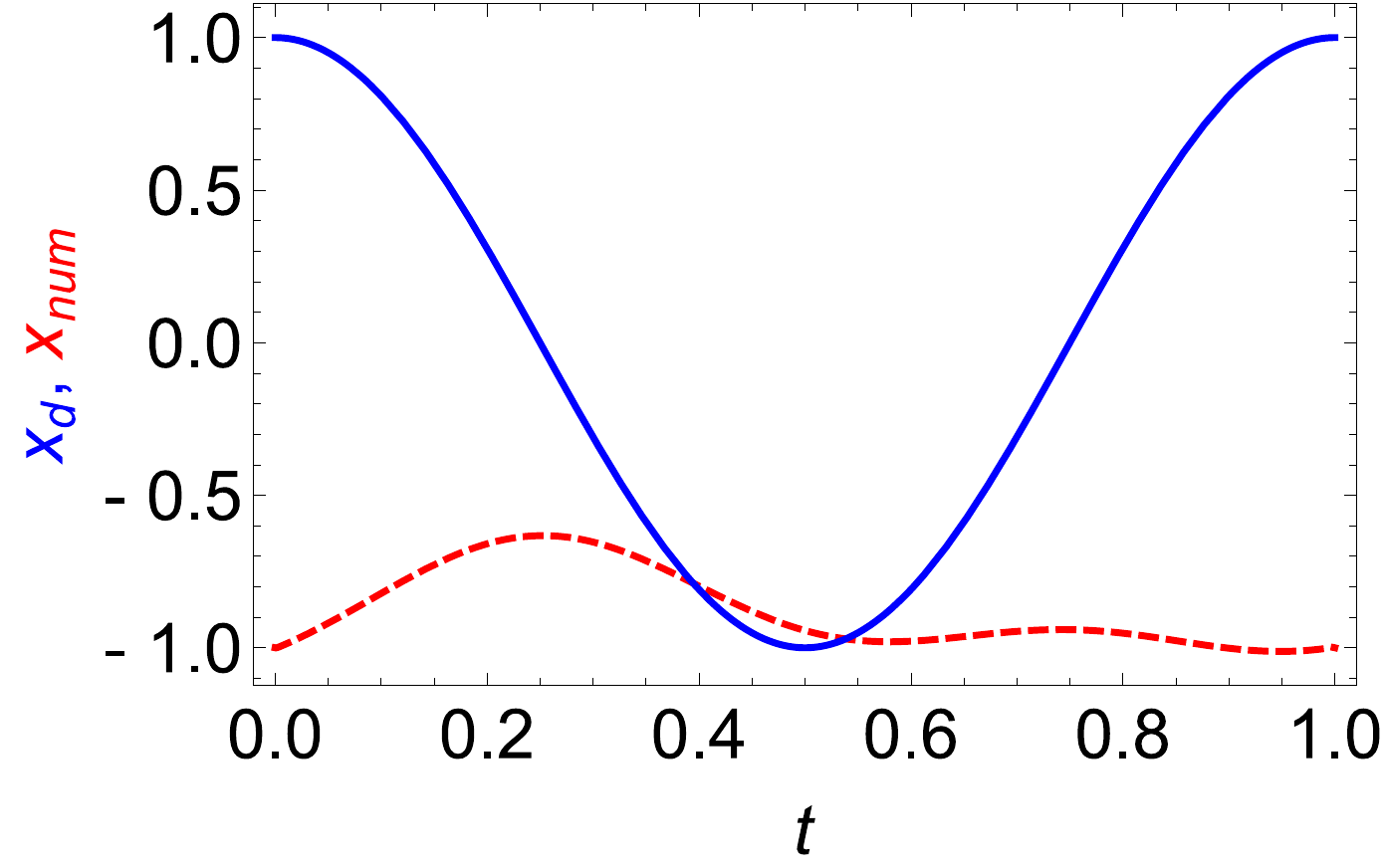}\hspace{0.2cm}\includegraphics[scale=0.475]{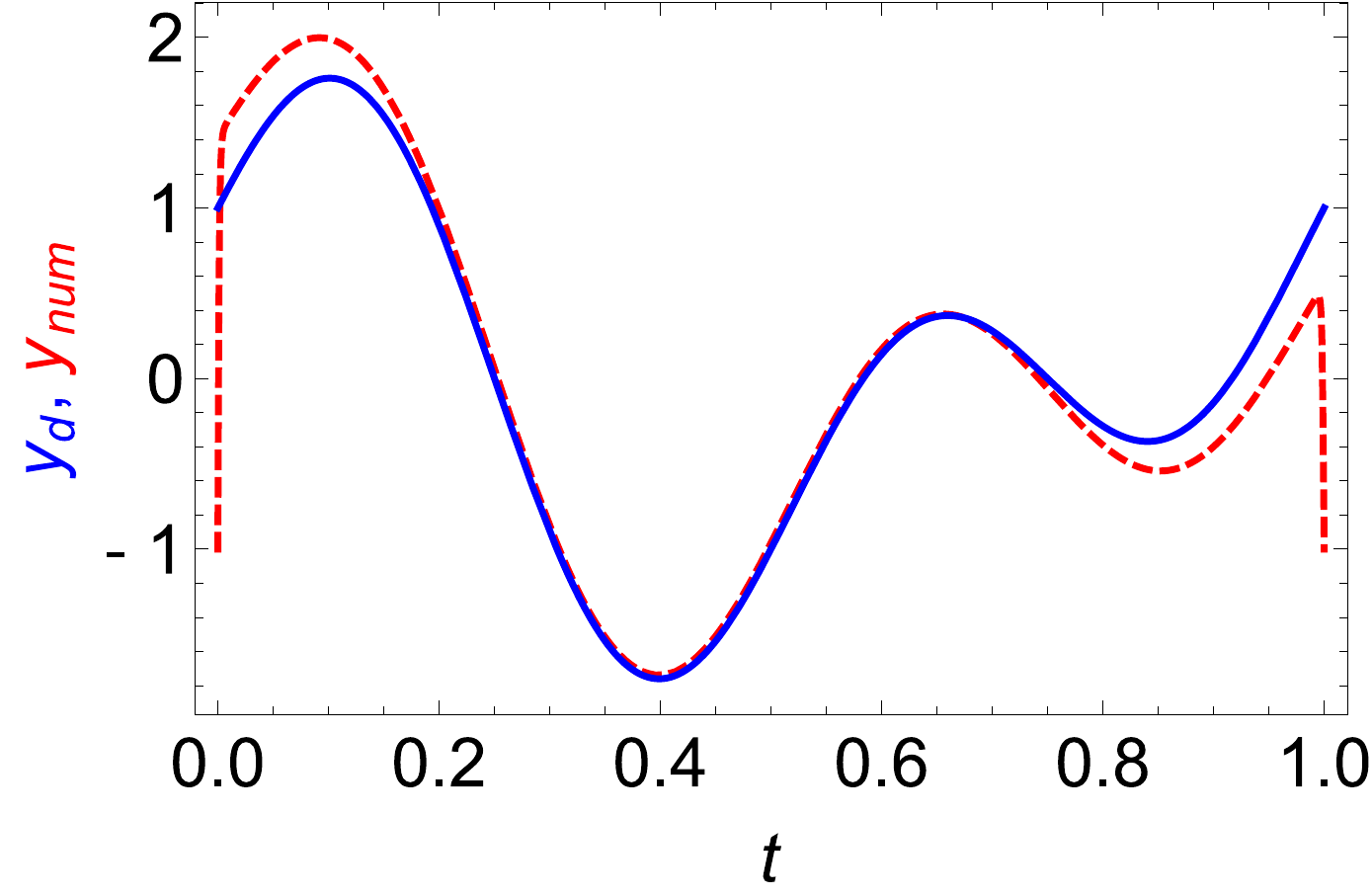}%
\captionof{figure}[Desired and optimal position and velocity  for the damped pendulum]{\label{fig:FigP1}Desired (blue solid line) and optimally controlled position $x$ (left) and velocity $y$ (right) over time for the damped pendulum. While the controlled velocity follows the desired velocity closely, the position is far off. The velocity exhibits an initial and terminal boundary layer.}
%"/home/jakob/ACADOtoolkit/examples/my_examples/Pendulum/MathematicalPendulum.nb"
\end{center}
\end{minipage}

\begin{minipage}{1.0\linewidth}
\begin{center}
\includegraphics[scale=0.49]{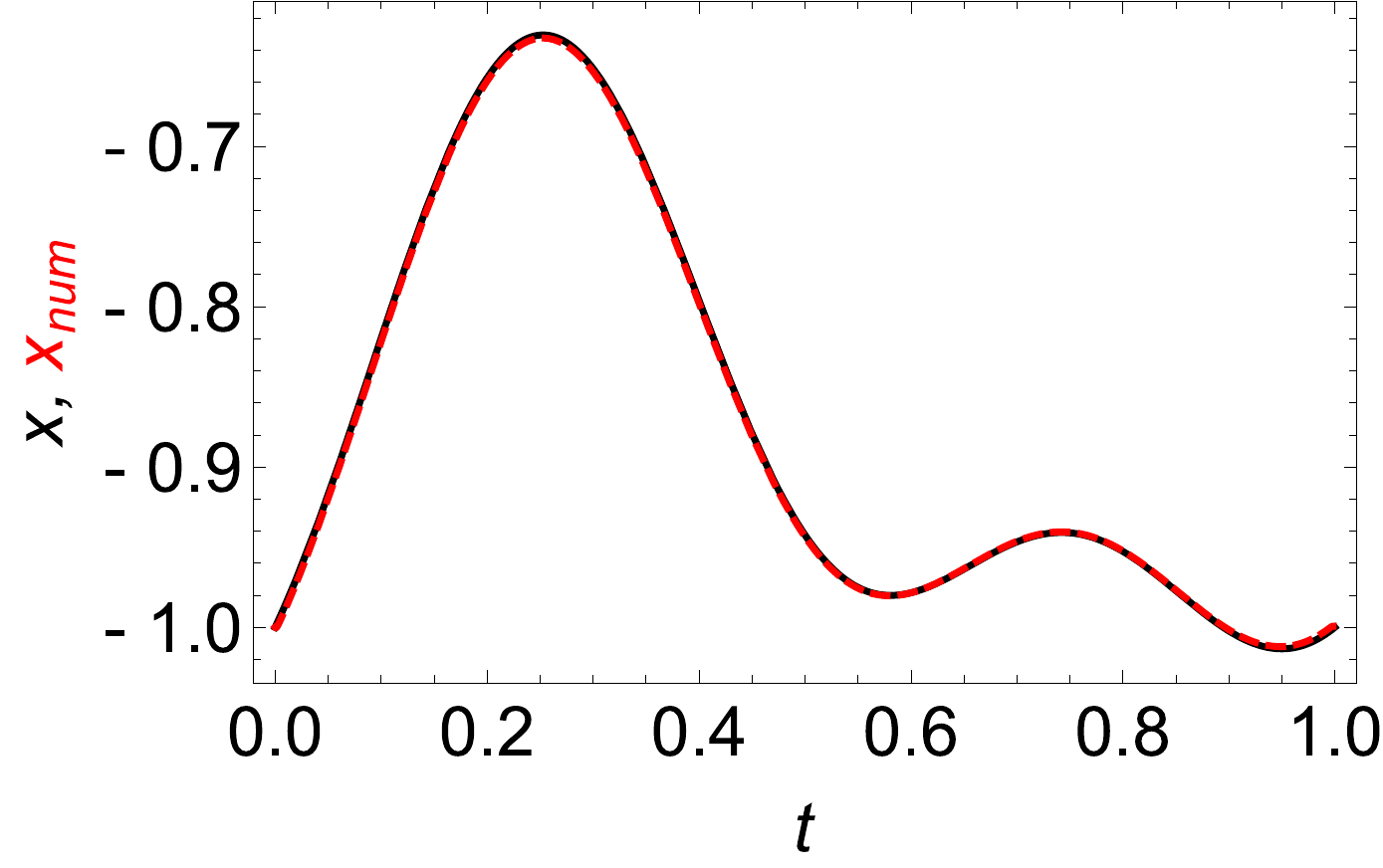}\hspace{0.2cm}\includegraphics[scale=0.475]{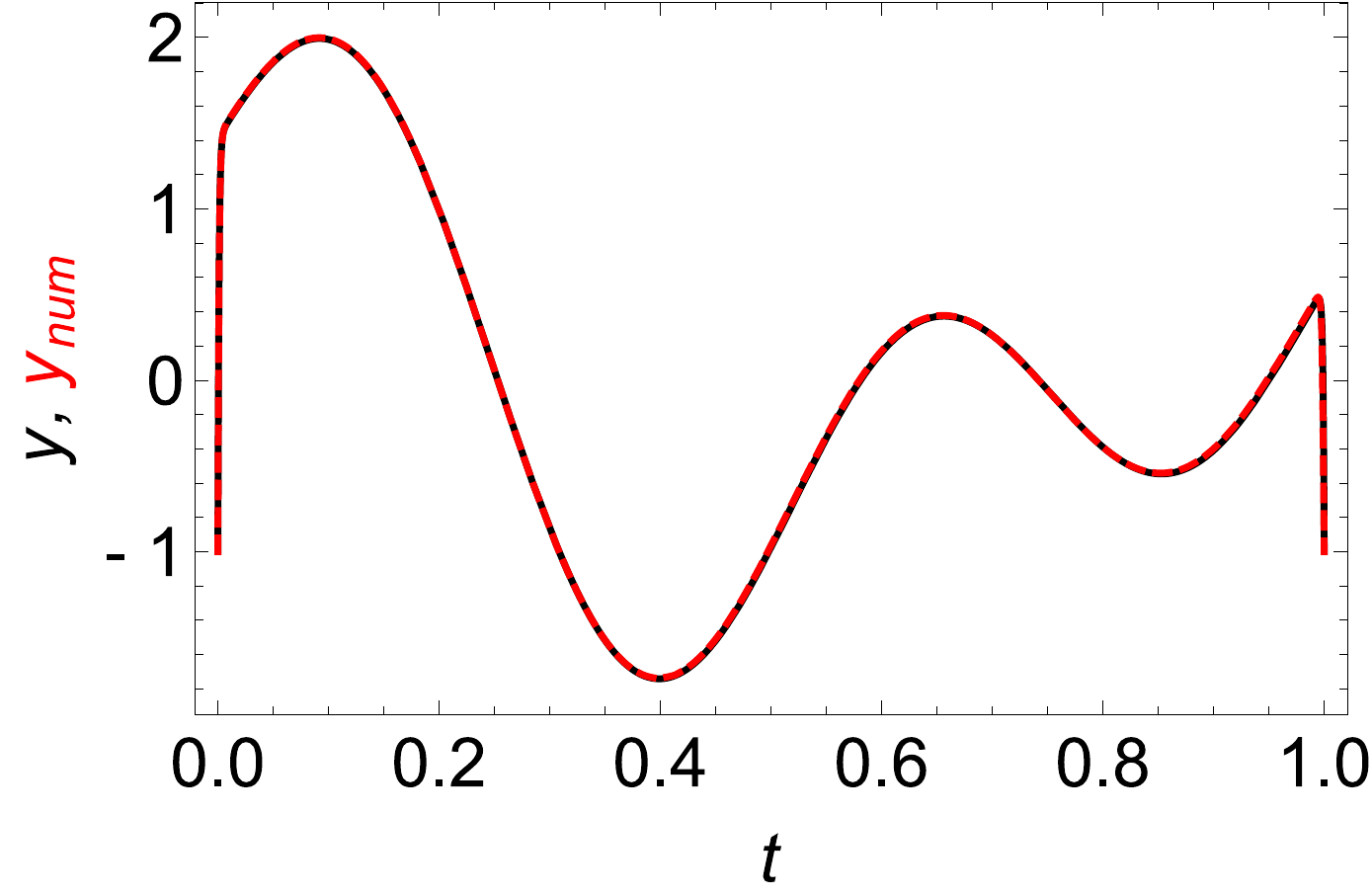}%
\captionof{figure}[Comparison of analytical and numerical position and velocity]{\label{fig:FigP2}Comparison of analytical (black solid line) and numerical (red dashed line) result for the controlled position (left) and velocity (right) over time reveals nearly perfect agreement.}
%"/home/jakob/ACADOtoolkit/examples/my_examples/Pendulum/MathematicalPendulum.nb"
\end{center}
\end{minipage}

As a second example, consider the desired trajectory
\begin{align}
x_{d}\left(t\right) & =\cos\left(2\pi t\right), & y_{d}\left(t\right) & =\cos\left(2\pi t\right)+\sin\left(4\pi t\right)+100,\label{eq:PendulumDesiredTrajectory2}
\end{align}
together with the same initial and terminal conditions Eqs. \eqref{eq:InitTermPendulum}
as above. Equations \eqref{eq:PendulumDesiredTrajectory2} differ
from the desired trajectory from Eqs. \eqref{eq:PendulumDesiredTrajectory}
only in a constant shift of the activator component, $y_{d}\left(t\right)\rightarrow y_{d}\left(t\right)+\alpha$.
Figure \ref{fig:FigP11} compares the numerical solution for the controlled
state trajectory $\boldsymbol{x}_{num}\left(t\right)$ with the desired
trajectory $\boldsymbol{x}_{d}\left(t\right)$ as given by Eq. \eqref{eq:PendulumDesiredTrajectory2}.
In contrast to Fig. \ref{fig:FigP1}, the controlled position over
time (Fig. \ref{fig:FigP11} left) is much closer to its desired counterpart
than the velocity over time (Fig. \ref{fig:FigP11} right). Figure
\ref{fig:FigP12} shows a direct comparison of the controlled state
trajectories for the desired trajectory Eqs. \eqref{eq:PendulumDesiredTrajectory}
(black solid line) and for the desired trajectory Eqs. \eqref{eq:PendulumDesiredTrajectory2}
(yellow dotted line). Surprisingly, the controlled state trajectories
are identical.

\begin{minipage}{1.0\linewidth}
\begin{center}
\includegraphics[scale=0.49]{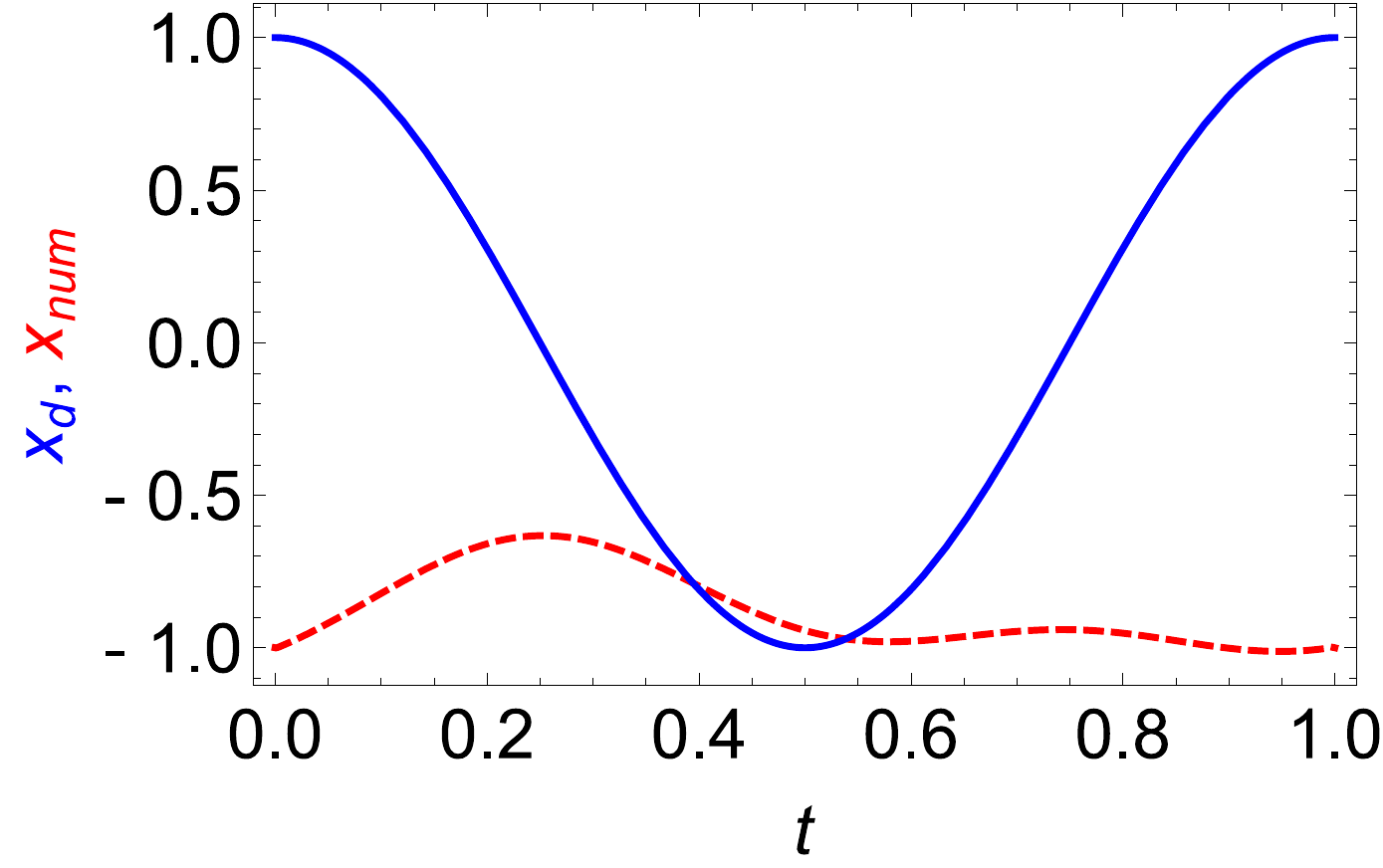}\hspace{0.2cm}\includegraphics[scale=0.475]{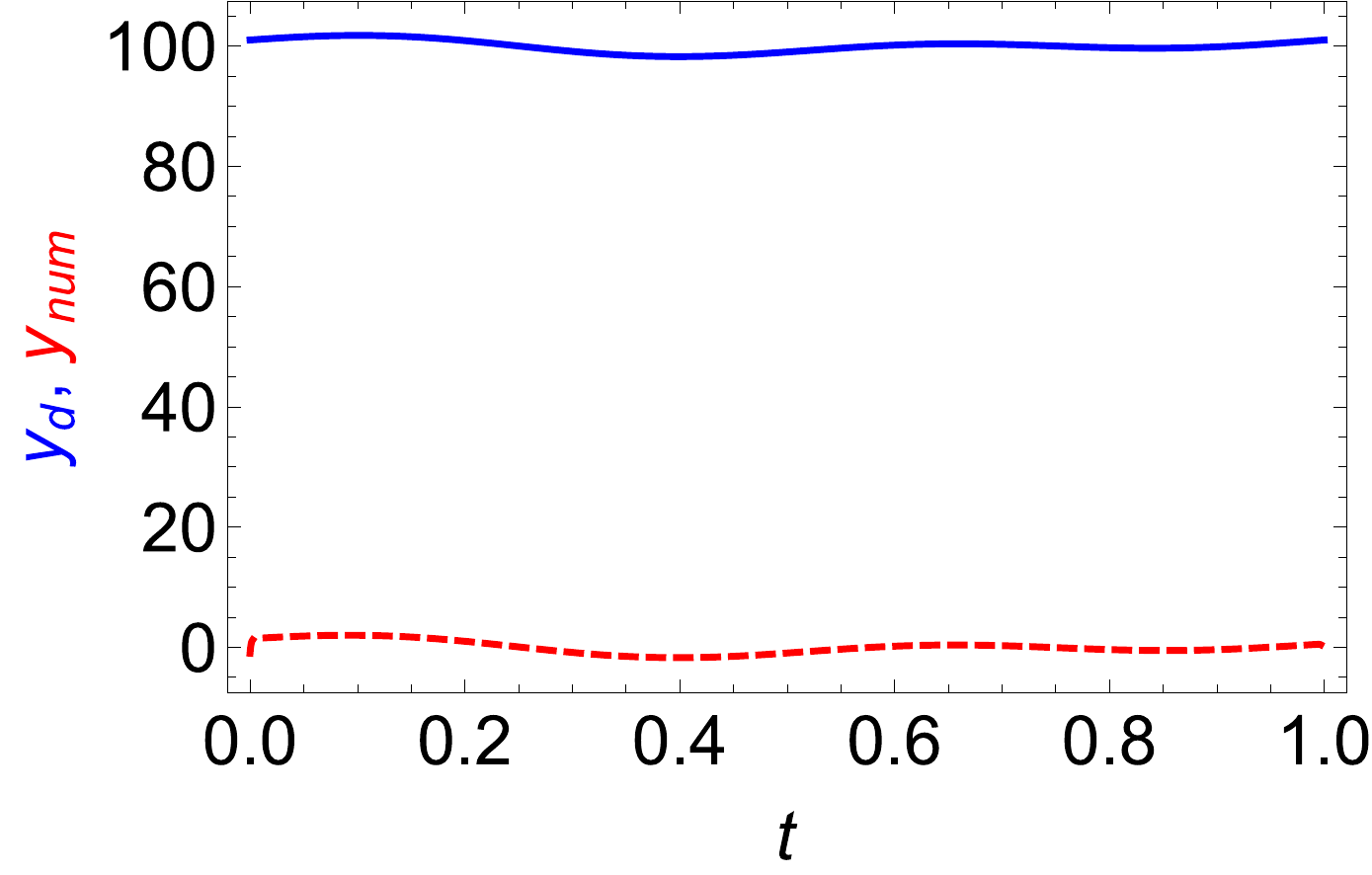}%
\captionof{figure}[Desired and optimal trajectory of the damped pendulum]{\label{fig:FigP11}Desired (blue solid line) and actually realized position $x$ (left) and velocity $y$ (right) over time for the damped pendulum for a desired trajectory given by Eqs. \eqref{eq:PendulumDesiredTrajectory2}. The velocity over time differs significantly from its desired counterpart (right), while the position  over time lies in the correct range of values (left).}
%"/home/jakob/ACADOtoolkit/examples/my_examples/Pendulum/MathematicalPendulum.nb"
\end{center}
\end{minipage}

\begin{minipage}{1.0\linewidth}
\begin{center}
\includegraphics[scale=0.49]{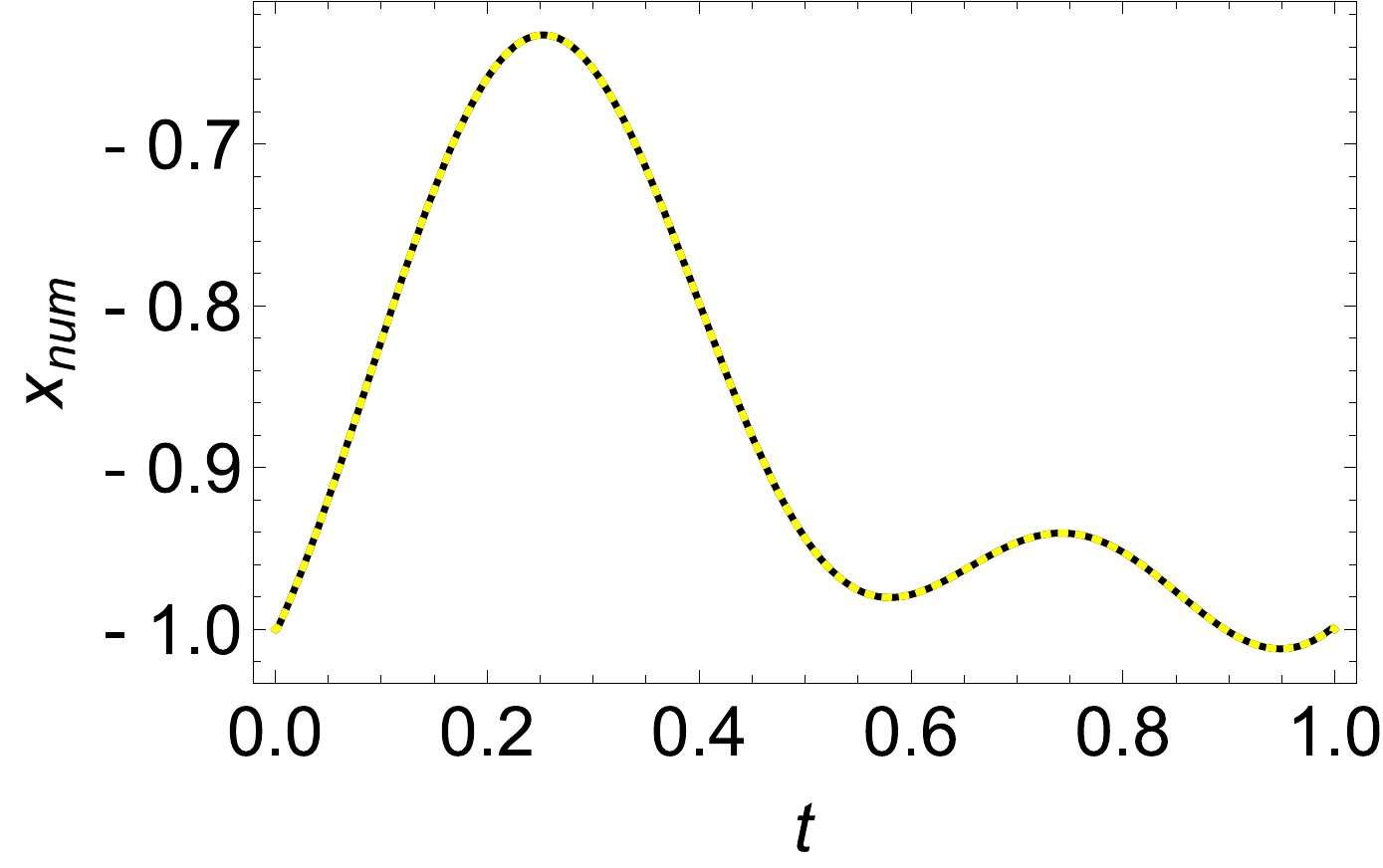}\hspace{0.2cm}\includegraphics[scale=0.475]{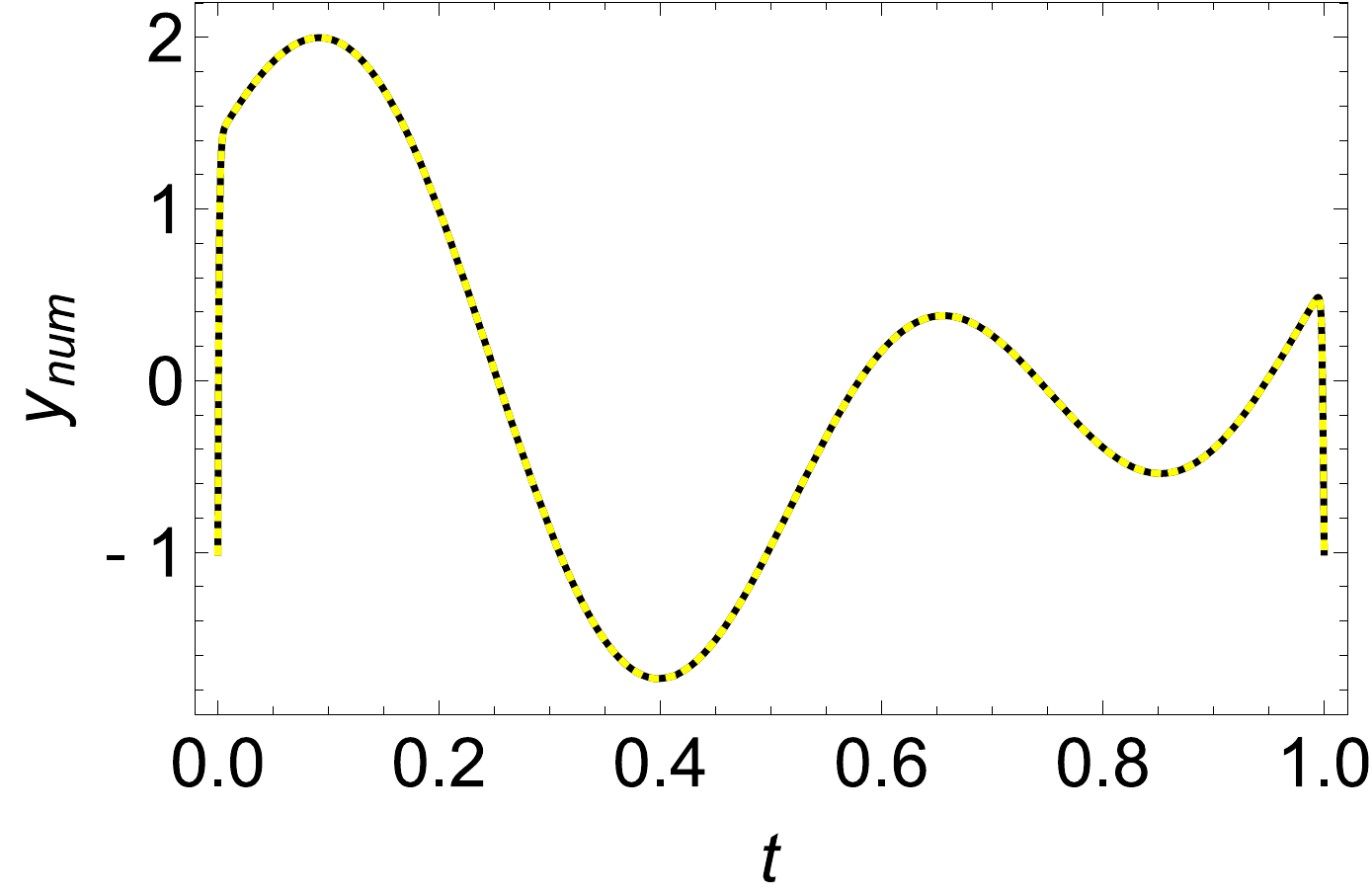}%
\captionof{figure}[Invariance of the optimal state trajectory under a  constant shift]{\label{fig:FigP12}Invariance of the optimal state trajectory under a  constant shift of the desired velocity over time. Position (left) and velocity (right) over time of the controlled state trajectory for two desired trajectories, Eqs. \eqref{eq:PendulumDesiredTrajectory} (black solid line) and Eqs. \eqref{eq:PendulumDesiredTrajectory2} (yellow dotted line). The desired trajectories differ by a constant shift of the desired velocity. Surprisingly, both controlled state trajectories are identical.}
%"/home/jakob/ACADOtoolkit/examples/my_examples/Pendulum/MathematicalPendulum.nb"
\end{center}
\end{minipage}

A more detailed look at the analytical result confirms the findings
of Fig. \ref{fig:FigP12}. For mechanical systems, the outer solution
is given by
\begin{align}
\left(\begin{array}{c}
x_{O}\left(t\right)\\
y_{O}\left(t\right)
\end{array}\right) & =\boldsymbol{\Phi}\left(t,t_{0}\right)\left(\begin{array}{c}
x_{\text{init}}\\
y_{\text{init}}
\end{array}\right)+\intop_{t_{0}}^{t}d\tau\boldsymbol{\Phi}\left(t,\tau\right)\boldsymbol{f}\left(\tau\right),
\end{align}
with state transition matrix $\boldsymbol{\Phi}\left(t,t_{0}\right)$
and inhomogeneity
\begin{align}
\boldsymbol{f}\left(t\right) & =\left(\begin{array}{c}
0\\
\dot{y}_{d}\left(t\right)-\frac{s_{1}}{s_{2}}x_{d}\left(t\right)
\end{array}\right).
\end{align}
While $\boldsymbol{f}\left(t\right)$ depends on $x_{d}\left(t\right)$
and the time derivative $\dot{y}_{d}\left(t\right)$, it is independent
of $y_{d}\left(t\right)$ itself. Thus, the inhomogeneity is invariant
under a constant shift $y_{d}\left(t\right)\rightarrow y_{d}\left(t\right)+\alpha$.
A more careful analysis reveals that the full composite solution for
$\boldsymbol{x}\left(t\right)$, Eqs. \eqref{eq:xFinalSol} and \eqref{eq:yFinalSol},
is independent of a constant shift $\alpha$ as long as $a_{1}=0$.
Consequently, the composite solution for the control signal, Eq. \eqref{eq:uFinalSol},
is independent of a constant shift $\alpha$ as well.

\end{example}

\subsection{Discussion}

Analytical approximations for optimal trajectory tracking for a two-dimensional
dynamical system were derived in Section \ref{sec:TwoDimensionalDynamicalSystem}.
The system includes the mechanical control systems from Example \ref{ex:OneDimMechSys1}
and the activator-controlled FHN model from Example \ref{ex:FHN1}
as special cases. The control acts on the nonlinear equation \eqref{eq:yState}
for $y\left(t\right)$ while the uncontrolled equation \eqref{eq:xState}
for $x\left(t\right)$ is linear.

The necessary optimality conditions are rearranged such that the highest
order time derivative as well as every occurrence of the nonlinearity
$R\left(x,y\right)$ is multiplied by the small parameter $\epsilon$.
This constitutes a system of singularly perturbed differential equations
amenable to a perturbative treatment. The solution reveals that the
$y$-component, i.e., the activator of the FHN model or the velocity
of mechanical systems, exhibits steep transition regions close to
the initial and terminal time. In the context of singular perturbation
theory, these transitions are interpreted as boundary layers with
width $\epsilon$ and arise as solutions to the inner equations. The
inner solutions connect the initial and terminal condition, respectively,
with the outer solution. The outer solution is valid only within the
time domain. Boundary layers occur even if the initial and terminal
conditions lie on the desired trajectory. The outer equations are
linear and their analytical solutions are available in closed form.
The inner equations depend on the nonlinear coupling function $b\left(x,y\right)$.
Neither the outer nor the inner equations depend on the nonlinearity
$R\left(x,y\right)$. The nonlinearity $R\left(x,y\right)$ is entirely
absorbed by the control signal. The control signal depends on $R\left(x,y\right)$
as well as on $b\left(x,y\right)$.

As $\epsilon\rightarrow0$, the boundary layers degenerate to jumps
located at the initial and terminal time. Simultaneously, the control
signal diverges and approaches the form of a Dirac delta function.
The strength of the delta kicks, i.e., the coefficient of the Dirac
delta function, is twice the height of the jumps. Because the delta
kick is located right at the time domain boundaries, only half of
the kick contributes to the time evolution. For $\epsilon=0$, the
composite solution is an exact solution to optimal trajectory tracking.
The analytical form of the composite solution is entirely determined
by the outer equations. No traces of the boundary layers remain except
for the mere existence of the jumps in $y\left(t\right)$. However,
the existence of these jumps, and therefore the existence of the entire
exact solution, relies on the existence of solutions to the inner
equations for the appropriate initial, terminal, and matching conditions.
In contrast to the perturbative result for $\epsilon>0$, the dynamics
is independent of the coupling function $b\left(x,y\right)$. Both
sources $R\left(x,y\right)$ and $b\left(x,y\right)$ of nonlinearity
are irrelevant for the controlled state trajectory. It is in this
sense that we are able to speak about linearity in unregularized nonlinear
optimal control. This result unveils a linear structure underlying
nonlinear optimal trajectory tracking.

A first analysis of the analytical results reveals that the controlled
state $y\left(t\right)$ is invariant under a constant shift $\alpha$
of the desired velocity, $y_{d}\left(t\right)\rightarrow y_{d}\left(t\right)+\alpha$,
as long as $a_{1}=0$. Note that $a_{1}=0$ for all mechanical control
systems, see Eq. \eqref{eq:xState}. This behavior is partially retained
for the FHN model as long as $a_{1}$ is small. Such insights can
hardly be obtained from numerical simulations alone. The impact of
this finding depends on the physical interpretation of the dynamical
system. For mechanical systems, $\dot{x}\left(t\right)=y\left(t\right)$
denotes the velocity of the system. Whereas shifting $y_{d}\left(t\right)$
has no effect on the controlled velocity $y\left(t\right)$, transforming
the desired position over time as $x_{d}\left(t\right)\rightarrow x_{d}\left(t\right)+\alpha t$
also changes $y\left(t\right)$. The controlled velocity can nevertheless
be affected by appropriately designing desired trajectories.

\section{\label{sec:OptimalFeedback}Optimal feedback control}

The approximate solution to the necessary optimality conditions depends
on the initial state $\boldsymbol{x}_{0}=\boldsymbol{x}\left(t_{0}\right)$.
Two possibilities exist to determine $\boldsymbol{x}_{0}$. Either
the system is prepared in state $\boldsymbol{x}_{0}$, or $\boldsymbol{x}_{0}$
is obtained by measurement. Knowing the value of $\boldsymbol{x}_{0}$,
no further information about the controlled system's state is necessary
to compute the control signal and the state's time evolution. The
control is an \textit{open loop control}.

External influences not modeled by the system dynamics can destabilize
the controlled system. A measurement $\tilde{x}_{0}=x\left(\tilde{t}_{0}\right),\,\tilde{y}_{0}=y\left(\tilde{t}_{0}\right)$
performed at a later time $t=\tilde{t}_{0}>t_{0}$ can be used to
update the control with $x_{0}=\tilde{x}_{0},\,y_{0}=\tilde{y}_{0}$
as the new initial condition. By feeding repeated measurements back
into the controlled system, it is possible to counteract unmodeled
perturbations occurring in between measurements. The optimal control
solution derived in Section \ref{sec:TwoDimensionalDynamicalSystem}
is a \textit{sampled-data feedback law }\cite{bryson1969applied}\textit{.}
The initial time $t_{0}$ is the most recent sampling time, and the
initial conditions $x_{0}=x\left(t_{0}\right),\,y_{0}=y\left(t_{0}\right)$
are measurements of the controlled system's state.

If a continuous monitoring of the system's state is possible, one
can set $t_{0}\rightarrow t$ in the composite state solutions Eqs.
\eqref{eq:xFinalSol}-\eqref{eq:controlFinalSol}. The initial conditions
$x_{0}=x\left(t\right),\,y_{0}=y\left(t\right)$ become functions
of the current state of the controlled system itself. This is known
as an optimal \textit{continuous time feedback law,} also called a
\textit{closed loop control}. Similarly, setting $t_{0}\rightarrow t-T$
and $x_{0}=x\left(t-T\right),\,y_{0}=y\left(t-T\right)$ yields a
\textit{continuous time-delayed feedback law.} The state measurements
are fed back to the system after a delay time $T>0$.

Thus, optimal feedback control requires knowledge of the controlled
state trajectory's dependence on its initial state $\boldsymbol{x}_{0}$.
A numerical solution to optimal control, determined for a single specified
value of $\boldsymbol{x}_{0}$, cannot be used for feedback control.
Instead, optimal feedback control is obtained from the Hamilton-Jacobi-Bellman
equation. This PDE is the central object of the Dynamic Programming
approach to optimal control founded by Richard Bellmann and coworkers
\cite{bellman1957dynamic}.

Consider the target functional 
\begin{align}
\mathcal{J}\left[\boldsymbol{x}\left(t\right),\boldsymbol{u}\left(t\right)\right] & =\intop_{t_{0}}^{t_{1}}dt\frac{1}{2}\left(\boldsymbol{x}\left(t\right)-\boldsymbol{x}_{d}\left(t\right)\right)^{T}\boldsymbol{\mathcal{S}}\left(\boldsymbol{x}\left(t\right)-\boldsymbol{x}_{d}\left(t\right)\right)\nonumber \\
 & +\frac{1}{2}\left(\boldsymbol{x}\left(t_{1}\right)-\boldsymbol{x}_{1}\right)^{T}\boldsymbol{\mathcal{S}}_{1}\left(\boldsymbol{x}\left(t_{1}\right)-\boldsymbol{x}_{1}\right)+\frac{\epsilon^{2}}{2}\intop_{t_{0}}^{t_{1}}dt\left(\boldsymbol{u}\left(t\right)\right)^{2}.\label{eq:OptimalTrajectoryTrackingFunctional-1}
\end{align}
The functional $\mathcal{J}\left[\boldsymbol{x}\left(t\right),\boldsymbol{u}\left(t\right)\right]$
is to be minimized subject to the constraints
\begin{align}
\boldsymbol{\dot{x}}\left(t\right) & =\boldsymbol{R}\left(\boldsymbol{x}\left(t\right)\right)+\boldsymbol{\mathcal{B}}\left(\boldsymbol{x}\left(t\right)\right)\boldsymbol{u}\left(t\right), & \boldsymbol{x}\left(t_{0}\right) & =\boldsymbol{x}_{0}.
\end{align}
Denote the minimal value of $\mathcal{J}$ by $\mathcal{J}_{0}$.
$\mathcal{J}_{0}$ is obtained by evaluating $\mathcal{J}$ at the
optimally controlled state trajectory $\boldsymbol{x}\left(t\right)$
and its corresponding control signal $\boldsymbol{u}\left(t\right)$.
$\mathcal{J}_{0}$ can be considered as a function of the initial
state $\boldsymbol{x}_{0}$ and the initial time $t_{0}$, $\mathcal{J}_{0}=\mathcal{J}_{0}\left(\boldsymbol{x}_{0},t_{0}\right)$.
The Hamilton-Jacobi-Bellman equation is a nonlinear evolution equation
for $\mathcal{J}_{0}$ given by \cite{bryson1969applied} 
\begin{align}
0 & =\min_{\boldsymbol{u}}\left\{ \nabla\mathcal{J}_{0}\left(\boldsymbol{x},t\right)\left(\boldsymbol{R}\left(\boldsymbol{x}\right)+\boldsymbol{\mathcal{B}}\left(\boldsymbol{x}\right)\boldsymbol{u}\right)+\frac{1}{2}\left(\boldsymbol{x}-\boldsymbol{x}_{d}\left(t\right)\right)^{T}\boldsymbol{\mathcal{S}}\left(\boldsymbol{x}-\boldsymbol{x}_{d}\left(t\right)\right)+\frac{\epsilon^{2}}{2}\boldsymbol{u}^{2}\right\} \nonumber \\
 & +\dfrac{\partial}{\partial t}\mathcal{J}_{0}\left(\boldsymbol{x},t\right).\label{eq:HJBEquation}
\end{align}
Equation \eqref{eq:HJBEquation} is supplemented with the terminal
condition
\begin{align}
\mathcal{J}_{0}\left(\boldsymbol{x},t_{1}\right) & =\frac{1}{2}\left(\boldsymbol{x}-\boldsymbol{x}_{1}\right)^{T}\boldsymbol{\mathcal{S}}_{1}\left(\boldsymbol{x}-\boldsymbol{x}_{1}\right).\label{eq:HJBInitiCond}
\end{align}
The co-state $\boldsymbol{\lambda}\left(t\right)$, considered as
a function of the initial state $\boldsymbol{x}_{0}$, is given by
the gradient of $\mathcal{J}_{0}$,
\begin{align}
\boldsymbol{\lambda}^{T}\left(t\right) & =\nabla\mathcal{J}_{0}\left(\boldsymbol{x}_{0},t\right).
\end{align}
Determining the minimum on the right hand side of Eq. \eqref{eq:HJBEquation}
yields a relation between the control and the gradient of $\mathcal{J}_{0}$.
This relation is analogous to the stationarity condition, Eq. \eqref{eq:ControlSolutionAdjointEquation},
of the necessary optimality conditions, 
\begin{align}
\nabla\mathcal{J}_{0}\left(\boldsymbol{x},t\right)\boldsymbol{\mathcal{B}}\left(\boldsymbol{x}\right)+\epsilon^{2}\boldsymbol{u}^{T} & =\boldsymbol{0}.
\end{align}
Solving for the control signal $\boldsymbol{u}$ yields 
\begin{align}
\boldsymbol{u} & =-\dfrac{1}{\epsilon^{2}}\boldsymbol{\mathcal{B}}^{T}\left(\boldsymbol{x}\right)\nabla\mathcal{J}_{0}^{T}\left(\boldsymbol{x},t\right).\label{eq:HJBControl}
\end{align}
The Hamilton-Jacobi-Bellman equation for optimal trajectory tracking
becomes
\begin{align}
-\epsilon^{2}\dfrac{\partial}{\partial t}\mathcal{J}_{0}\left(\boldsymbol{x},t\right) & =-\dfrac{1}{2}\nabla\mathcal{J}_{0}\left(\boldsymbol{x},t\right)\boldsymbol{\mathcal{B}}\left(\boldsymbol{x}\right)\boldsymbol{\mathcal{B}}^{T}\left(\boldsymbol{x}\right)\nabla\mathcal{J}_{0}^{T}\left(\boldsymbol{x},t\right)+\epsilon^{2}\nabla\mathcal{J}_{0}\left(\boldsymbol{x},t\right)\boldsymbol{R}\left(\boldsymbol{x}\right)\nonumber \\
 & +\frac{\epsilon^{2}}{2}\left(\boldsymbol{x}-\boldsymbol{x}_{d}\left(t\right)\right)^{T}\boldsymbol{\mathcal{S}}\left(\boldsymbol{x}-\boldsymbol{x}_{d}\left(t\right)\right).\label{eq:HJBTrajectoryTracking}
\end{align}
Knowing the solution to Eq. \eqref{eq:HJBTrajectoryTracking}, the
open loop control signal for a system with initial state $\boldsymbol{x}\left(t_{0}\right)=\boldsymbol{x}_{0}$
is recovered from Eq. \eqref{eq:HJBControl} as 
\begin{align}
\boldsymbol{u}\left(t\right) & =-\dfrac{1}{\epsilon^{2}}\boldsymbol{\mathcal{B}}^{T}\left(\boldsymbol{x}_{0}\right)\nabla\mathcal{J}_{0}^{T}\left(\boldsymbol{x}_{0},t\right).
\end{align}
A continuous time feedback law $\boldsymbol{u}\left(t\right)=\boldsymbol{u}\left(\boldsymbol{x}\left(t\right),t\right)$
depends on the actual state $\boldsymbol{x}\left(t\right)$ of the
controlled system and is given by
\begin{align}
\boldsymbol{u}\left(t\right) & =\boldsymbol{u}\left(\boldsymbol{x}\left(t\right),t\right)=-\dfrac{1}{\epsilon^{2}}\boldsymbol{\mathcal{B}}^{T}\left(\boldsymbol{x}\left(t\right)\right)\nabla\mathcal{J}_{0}^{T}\left(\boldsymbol{x}\left(t\right),t\right).
\end{align}
We do not attempt to solve the nonlinear PDE Eq. \eqref{eq:HJBTrajectoryTracking},
but end with some concluding remarks about the difficulties encountered
in doing so. First of all, for a vanishing regularization parameter
$\epsilon=0$, Eq. \eqref{eq:HJBTrajectoryTracking} suffers from
a similar degeneracy as the necessary optimality conditions. Because
the time derivative $\partial_{t}\mathcal{J}_{0}$ vanishes for $\epsilon=0$,
$\mathcal{J}_{0}$ cannot satisfy the terminal condition \eqref{eq:HJBInitiCond}.
Second, to solve Eq. \eqref{eq:HJBTrajectoryTracking} numerically
for $\epsilon>0$ is a formidable task, especially if the dimension
$n$ of the state space is large. A discretization with $N_{x}$ points
for a single state space dimension results in $N_{x}^{n}$ discretization
points for the full state space. The computational cost increases
exponentially with the state space dimension. This is the ``curse
of dimensionality'', as it was called by Bellman himself \cite{bellman1957dynamic}.

\subsection{\label{sub:ContinousTimeFeedback}Continuous time feedback}

\subsubsection{Derivation of the feedback law}

The approximate solution to the optimal trajectory tracking problem,
Eqs. \eqref{eq:xFinalSol}-\eqref{eq:lambdaFinalSol} and Eq. \eqref{eq:uFinalSol},
is rendered as a continuous time feedback law. First, the initial
conditions $x_{0}$ and $y_{0}$ are given by the controlled state
components $x_{0}=x\left(t_{0}\right)=x\left(t\right)$ and $y_{0}=y\left(t_{0}\right)=y\left(t\right)$,
respectively. Second, every explicit appearance of $t_{0}$ in Eqs.
\eqref{eq:xFinalSol}-\eqref{eq:lambdaFinalSol} and Eq. \eqref{eq:uFinalSol}
is substituted by $t_{0}\rightarrow t$. All constants which depend
on time $t_{0}$, as e.g. $x_{\text{init}}$, become time dependent
on the current time $t$. To minimize the confusion, these constants
are written as
\begin{align}
x_{\text{init}} & =x_{\text{init}}\left(t\right), & y_{\text{init}} & =y_{\text{init}}\left(t\right), & y_{\text{end}} & =y_{\text{end}}\left(t\right).\label{eq:Eq4149}
\end{align}
The outer solutions $x_{O}\left(t\right)$ and $y_{O}\left(t\right)$
given by Eqs. \eqref{eq:xOSol} and \eqref{eq:yOSol} assume a particularly
simple form,
\begin{align}
x_{O}\left(t\right) & =x_{\text{init}}\left(t\right), & y_{O}\left(t\right) & =y_{\text{init}}\left(t\right),
\end{align}
and all integral terms vanish. The composite solutions Eqs. \eqref{eq:xFinalSol}
and \eqref{eq:yFinalSol} reduce to
\begin{align}
x_{\text{comp}}\left(t\right) & =x_{O}\left(t\right)=x_{\text{init}}\left(t\right), & y_{\text{comp}}\left(t\right) & =Y_{L}\left(0\right)+Y_{R}\left(\left(t_{1}-t\right)/\epsilon\right)-y_{\text{end}}\left(t\right).\label{eq:Eq658}
\end{align}
The composite control signal Eq. \eqref{eq:uFinalSol} becomes
\begin{align}
u_{\text{comp}}\left(t\right) & =U_{L}\left(0\right)+U_{R}\left(\left(t_{1}-t\right)/\epsilon\right)-u_{O}\left(t_{1}\right).\label{eq:Eq659}
\end{align}
Note that $Y_{R},\,U_{L},\,U_{R}$ as well as $u_{O}\left(t_{1}\right)$
are just abbreviations and still depend on time $t$ through the time-dependent
parameters from Eq. \eqref{eq:Eq4149}. The terms originating from
the right boundary layer become important only for times $t\lesssim t_{1}$
close to the terminal time. For simplicity, consider the limit $t_{1}\rightarrow\infty$.
Because of
\begin{align}
\lim_{t_{1}\rightarrow\infty}Y_{R}\left(\left(t_{1}-t\right)/\epsilon\right) & =y_{\text{end}},\\
\lim_{t_{1}\rightarrow\infty}U_{R}\left(\left(t_{1}-t\right)/\epsilon\right) & =\lim_{t_{1}\rightarrow\infty}u_{O}\left(t_{1}\right),
\end{align}
the terms originating from the right boundary layer in Eqs. \eqref{eq:Eq658},
\eqref{eq:Eq659} cancel. Together with $x_{\text{init}}\left(t\right)=x_{0}=x\left(t\right)$
and $Y_{L}\left(0\right)=y_{0}=y\left(t\right)$, the composite state
is 
\begin{align}
x_{\text{comp}}\left(t\right) & =x_{0}=x\left(t\right), & y_{\text{comp}}\left(t\right) & =y_{0}=y\left(t\right).
\end{align}
The composite control solution simplifies to 
\begin{align}
u_{\text{comp}}\left(t\right) & =U_{L}\left(0\right).
\end{align}
Note that $U_{L}\left(0\right)$ still depends on time $t$ through
the constants $x_{0}=x\left(t\right),\,y_{0}=y\left(t\right)$ and
$y_{\text{init}}^{\infty}\left(t\right)$. Here, $y_{\text{init}}^{\infty}\left(t\right)$
denotes the constant $y_{\text{init}}\left(t\right)$ in the limit
$t_{1}\rightarrow\infty$ given by 
\begin{align}
y_{\text{init}}^{\infty}\left(t\right) & =\lim_{t_{1}\rightarrow\infty}y_{\text{init}}\left(t\right)=\intop_{t}^{\infty}e^{\left(t-\tau\right)\varphi_{1}}\left(\frac{a_{2}s_{1}}{s_{2}}x_{d}\left(\tau\right)-\left(a_{1}+\varphi_{1}\right)y_{d}\left(\tau\right)\right)d\tau\nonumber \\
 & -\dfrac{1}{a_{2}}\left(\varphi_{1}+a_{1}\right)x\left(t\right)-\frac{a_{0}}{a_{2}}\left(\frac{a_{1}}{\varphi_{1}}+1\right)+y_{d}\left(t\right).\label{eq:yinitinf}
\end{align}
All occurrences of $x_{0}$ and $y_{0}$ are substituted with $x\left(t\right)$
and $y\left(t\right)$, respectively. The feedback control derived
by the outlined procedure is called $u_{\text{fb}}\left(t\right)$
and given by 
\begin{align}
u_{\text{fb}}\left(t\right) & =U_{L}\left(0\right)=\frac{1}{b\left(x_{0},y_{0}\right)}\left(\dot{y}_{O}\left(t_{0}\right)+\dfrac{1}{\epsilon}Y_{L}'\left(0\right)-R\left(x_{0},y_{0}\right)\right)\nonumber \\
 & =\frac{1}{b\left(x\left(t\right),y\left(t\right)\right)}\left(a_{1}\left(y_{d}\left(t\right)-y_{\text{init}}^{\infty}\left(t\right)\right)+\dot{y}_{d}\left(t\right)+\frac{a_{2}s_{1}}{s_{2}}\left(x\left(t\right)-x_{d}\left(t\right)\right)\right)\nonumber \\
 & +\frac{1}{b\left(x\left(t\right),y\left(t\right)\right)}\left(\dfrac{1}{\epsilon}\sqrt{s_{2}}\left(y_{\text{init}}^{\infty}\left(t\right)-y\left(t\right)\right)\left|b\left(x\left(t\right),y\left(t\right)\right)\right|-R\left(x\left(t\right),y\left(t\right)\right)\right).\label{eq:FeedBackControl}
\end{align}
Equations \eqref{eq:EqYL} and \eqref{eq:yOEquation} were used to
substitute the expression $Y_{L}'\left(0\right)$ and $\dot{y}_{O}\left(t\right)$,
respectively. The feedback law Eq. \eqref{eq:FeedBackControl} explicitly
depends on the nonlinearities $R\left(x,y\right)$ and $b\left(x,y\right)$
and is only valid for an infinite terminal time $t_{1}\rightarrow\infty$.
An analogous but more complicated and longer expression can be derived
for finite time intervals $t_{0}\leq t\leq t_{1}$ from Eq. \eqref{eq:Eq659}.

\subsubsection{Feedback-controlled state trajectory}

The time evolution of the controlled state trajectory under feedback
control Eq. \eqref{eq:FeedBackControl} is analyzed. Using Eq. \eqref{eq:FeedBackControl}
in the controlled state equations \eqref{eq:xState} and \eqref{eq:yState}
results in
\begin{align}
\dot{x}\left(t\right) & =y\left(t\right),\label{eq:FeedBackControlledStatex}\\
\dot{y}\left(t\right) & =R\left(x\left(t\right),y\left(t\right)\right)+b\left(x\left(t\right),y\left(t\right)\right)u_{\text{fb}}\left(t\right)\nonumber \\
 & =a_{1}\left(y_{d}\left(t\right)-y_{\text{init}}^{\infty}\left(t\right)\right)+\dot{y}_{d}\left(t\right)+\frac{a_{2}s_{1}}{s_{2}}\left(x\left(t\right)-x_{d}\left(t\right)\right)\nonumber \\
 & +\dfrac{1}{\epsilon}\sqrt{s_{2}}\left(y_{\text{init}}^{\infty}\left(t\right)-y\left(t\right)\right)\left|b\left(x\left(t\right),y\left(t\right)\right)\right|.\label{eq:FeedBackControlledStatey}
\end{align}
Note that the nonlinearity $R\left(x,y\right)$ is eliminated from
Eq. \eqref{eq:FeedBackControlledStatey}, whereas the dependence on
the coupling function $b\left(x,y\right)$ is retained. Note that
$y_{\text{init}}^{\infty}\left(t\right)$ as given by Eq. \eqref{eq:yinitinf}
depends on $x\left(t\right)$ as well. Equations \eqref{eq:FeedBackControlledStatex}
and \eqref{eq:FeedBackControlledStatey} have to be solved with the
initial conditions 
\begin{align}
x\left(t_{0}\right) & =x_{0}^{\text{fb}}, & y\left(t_{0}\right) & =y_{0}^{\text{fb}}.\label{eq:InitialCondFeedback}
\end{align}
Due to the coupling function $b\left(x,y\right)$, Eqs. \eqref{eq:FeedBackControlledStatex}
and \eqref{eq:FeedBackControlledStatey} are nonlinear. No exact analytical
closed form solution exists. However, Eqs. \eqref{eq:FeedBackControlledStatex}
and \eqref{eq:FeedBackControlledStatey} can be solved perturbatively
using the small parameter $\epsilon$ for a perturbation expansion.
Due to the appearance of $1/\epsilon$, Eqs. \eqref{eq:FeedBackControlledStatex}
and \eqref{eq:FeedBackControlledStatey} constitute a singularly perturbed
system of differential equations. 

A procedure analogous to the approach to open loop control is applied.
The inner and outer equations and their solutions must be determined.
Combining them to a composite solution yields an approximate solution
uniformly valid over the while time interval. An initial boundary
layer is expected close to the initial time $t_{0}$. Because of the
assumed infinite terminal time $t_{1}\rightarrow\infty$, no terminal
boundary layer exists. The outer variables are denoted with index
$O$,
\begin{align}
x\left(t\right) & =x_{O}\left(t\right), & y\left(t\right) & =y_{O}\left(t\right).
\end{align}
To leading order in $\epsilon$, the outer equations are
\begin{align}
\dot{x}_{O}\left(t\right) & =y_{O}\left(t\right),\\
0 & =\left(y_{\text{init}}^{\infty}\left(t\right)-y_{O}\left(t\right)\right)\left|b\left(x_{O}\left(t\right),y_{O}\left(t\right)\right)\right|.
\end{align}
Because $b$ does not have a root by assumption, the unique solution
for $y_{O}\left(t\right)$ is
\begin{align}
y_{O}\left(t\right) & =y_{\text{init}}^{\infty}\left(t\right).\label{eq:YOSolFB}
\end{align}
The expression for $y_{\text{init}}^{\infty}\left(t\right)$ still
depends on $x\left(t\right)$, see Eq. \eqref{eq:yinitinf}. The state
$x_{O}\left(t\right)$ is governed by the differential equation
\begin{align}
\dot{x}_{O}\left(t\right) & =-\dfrac{1}{a_{2}}\left(\varphi_{1}+a_{1}\right)x_{O}\left(t\right)-\frac{a_{0}}{a_{2}}\left(\frac{a_{1}}{\varphi_{1}}+1\right)+y_{d}\left(t\right)\nonumber \\
 & +\intop_{t}^{\infty}e^{\left(t-\tau\right)\varphi_{1}}\left(\frac{a_{2}s_{1}}{s_{2}}x_{d}\left(\tau\right)-\left(a_{1}+\varphi_{1}\right)y_{d}\left(\tau\right)\right)d\tau.\label{Eq4166}
\end{align}
Equation \eqref{Eq4166} must be solved with the initial condition
\begin{align}
x\left(t_{0}\right) & =x_{\text{init}}^{\text{fb}}.
\end{align}
The constant $x_{\text{init}}^{\text{fb}}$ has to be determined by
matching the outer solutions with the inner solutions. The solution
for $x_{O}\left(t\right)$ is given by
\begin{align}
x_{O}\left(t\right) & =x_{\text{init}}^{\text{fb}}\exp\left(-\left(t-t_{0}\right)\frac{\left(a_{1}+\varphi_{1}\right)}{a_{2}}\right)+\int_{t_{0}}^{t}g\left(\tilde{t}\right)\exp\left(\frac{\left(a_{1}+\varphi_{1}\right)}{a_{2}}\left(\tilde{t}-t\right)\right)d\tilde{t},\label{eq:xOSolFB}
\end{align}
with abbreviation $g\left(t\right)$ 
\begin{align}
g\left(t\right) & =-\frac{a_{0}}{a_{2}}\left(\frac{a_{1}}{\varphi_{1}}+1\right)+y_{d}\left(t\right)+\intop_{t}^{\infty}e^{\left(t-\tau\right)\varphi_{1}}\left(\frac{a_{2}s_{1}}{s_{2}}x_{d}\left(\tau\right)-\left(a_{1}+\varphi_{1}\right)y_{d}\left(\tau\right)\right)d\tau.
\end{align}
The outer equation are not able to satisfy both initial conditions
Eqs. \eqref{eq:InitialCondFeedback}. The initial boundary layer is
resolved using the time scale
\begin{align}
\tau_{L} & =\left(t-t_{0}\right)/\epsilon
\end{align}
and rescaled inner solutions
\begin{align}
X_{L}\left(\tau_{L}\right) & =X_{L}\left(\left(t-t_{0}\right)/\epsilon\right)=x\left(t\right)=x\left(t_{0}+\epsilon\tau_{L}\right),\\
Y_{L}\left(\tau_{L}\right) & =Y_{L}\left(\left(t-t_{0}\right)/\epsilon\right)=y\left(t\right)=y\left(t_{0}+\epsilon\tau_{L}\right).
\end{align}
Rewritten with the new time scale and rescaled functions, the feedback-controlled
state equations \eqref{eq:FeedBackControlledStatex}, \eqref{eq:FeedBackControlledStatey}
are 
\begin{align}
\dfrac{1}{\epsilon}\dot{X}_{L}\left(\tau_{L}\right) & =Y\left(\tau_{L}\right),\label{eq:FeedBackInnerEqx}\\
\dfrac{1}{\epsilon}\dot{Y}_{L}\left(\tau_{L}\right) & =a_{1}\left(y_{d}\left(t_{0}+\epsilon\tau_{L}\right)-y_{\text{init}}^{\infty}\left(t_{0}+\epsilon\tau_{L}\right)\right)+\dot{y}_{d}\left(t_{0}+\epsilon\tau_{L}\right)\nonumber \\
 & +\frac{a_{2}s_{1}}{s_{2}}\left(X_{L}\left(\tau_{L}\right)-x_{d}\left(t_{0}+\epsilon\tau_{L}\right)\right)\nonumber \\
 & +\dfrac{1}{\epsilon}\sqrt{s_{2}}\left(y_{\text{init}}^{\infty}\left(t_{0}+\epsilon\tau_{L}\right)-Y_{L}\left(\tau_{L}\right)\right)\left|b\left(x_{L}\left(\tau_{L}\right),Y_{L}\left(\tau_{L}\right)\right)\right|.\label{eq:FeedBackInnerEqy}
\end{align}
Note that $y_{\text{init}}^{\infty}\left(t_{0}+\epsilon\tau_{L}\right)$
still depends on $x\left(t\right)$ and becomes
\begin{align}
y_{\text{init}}^{\infty}\left(t_{0}+\epsilon\tau_{L}\right) & =\intop_{t_{0}+\epsilon\tau_{L}}^{\infty}e^{\left(t_{0}+\epsilon\tau_{L}-\tau\right)\varphi_{1}}\left(\frac{a_{2}s_{1}}{s_{2}}x_{d}\left(\tau\right)-\left(a_{1}+\varphi_{1}\right)y_{d}\left(\tau\right)\right)d\tau\nonumber \\
 & -\dfrac{1}{a_{2}}\left(\varphi_{1}+a_{1}\right)X_{L}\left(\tau_{L}\right)-\frac{a_{0}}{a_{2}}\left(\frac{a_{1}}{\varphi_{1}}+1\right)+y_{d}\left(t_{0}+\epsilon\tau_{L}\right).
\end{align}
The inner equations \eqref{eq:FeedBackInnerEqx}, \eqref{eq:FeedBackInnerEqy}
must be solved with the boundary conditions
\begin{align}
X_{L}\left(0\right) & =x_{0}^{\text{fb}}, & Y_{L}\left(0\right) & =y_{0}^{\text{fb}}.
\end{align}
To leading order in $\epsilon$, $y_{\text{init}}^{\infty}\left(t_{0}+\epsilon\tau_{L}\right)$
simplifies to 
\begin{align}
y_{\text{init}}^{\infty}\left(t_{0}\right) & =\intop_{t_{0}}^{\infty}e^{\left(t_{0}-\tau\right)\varphi_{1}}\left(\frac{a_{2}s_{1}}{s_{2}}x_{d}\left(\tau\right)-\left(a_{1}+\varphi_{1}\right)y_{d}\left(\tau\right)\right)d\tau\nonumber \\
 & -\dfrac{1}{a_{2}}\left(\varphi_{1}+a_{1}\right)X_{L}\left(\tau_{L}\right)-\frac{a_{0}}{a_{2}}\left(\frac{a_{1}}{\varphi_{1}}+1\right)+y_{d}\left(t_{0}\right),
\end{align}
and the inner equations simplify to 
\begin{align}
\dot{X}_{L}\left(\tau_{L}\right) & =0,\\
\dot{Y}_{L}\left(\tau_{L}\right) & =\sqrt{s_{2}}\left(y_{\text{init}}^{\infty}\left(t_{0}\right)-Y_{L}\left(\tau_{L}\right)\right)\left|b\left(x_{L}\left(\tau_{L}\right),Y_{L}\left(\tau_{L}\right)\right)\right|.\nonumber 
\end{align}
The solution for $X_{L}$ is
\begin{align}
X_{L}\left(\tau_{L}\right) & =x_{0}^{\text{fb}},
\end{align}
and the equation for $Y_{L}$ reduces to
\begin{align}
\dot{Y}_{L}\left(\tau_{L}\right) & =\sqrt{s_{2}}\left(y_{\text{init}}^{\infty}-Y_{L}\left(\tau_{L}\right)\right)\left|b\left(x_{0},Y_{L}\left(\tau_{L}\right)\right)\right|,\label{eq:YLFeedBack}\\
Y_{L}\left(0\right) & =y_{0}^{\text{fb}},
\end{align}
with constant $y_{\text{init}}^{\infty}$ given by
\begin{align}
y_{\text{init}}^{\infty} & =\intop_{t_{0}}^{\infty}e^{\left(t_{0}-\tau\right)\varphi_{1}}\left(\frac{a_{2}s_{1}}{s_{2}}x_{d}\left(\tau\right)-\left(a_{1}+\varphi_{1}\right)y_{d}\left(\tau\right)\right)d\tau\nonumber \\
 & -\dfrac{1}{a_{2}}\left(\varphi_{1}+a_{1}\right)x_{0}^{\text{fb}}-\frac{a_{0}}{a_{2}}\left(\frac{a_{1}}{\varphi_{1}}+1\right)+y_{d}\left(t_{0}\right).\label{eq:yinitinfConstant}
\end{align}
Equation \eqref{eq:YLFeedBack} is nonlinear and has no analytical
solution in closed form. Remarkably, Eq. \eqref{eq:YLFeedBack} has
the same form as Eq. \eqref{eq:EqYL} for the inner boundary layer
of open loop control. The initial boundary layers are governed by
the same dynamics regardless of open or closed loop control. The only
difference is the value of the constant $y_{\text{init}}^{\infty}$.
Assuming that the coupling function $b\left(x,y\right)$ does not
depend on $y$, Eq. \eqref{eq:YLFeedBack} can immediately be solved,
\begin{align}
Y_{L}\left(\tau_{L}\right) & =\exp\left(-\sqrt{s_{2}}\tau_{L}\left|b\left(x_{0}\right)\right|\right)\left(y_{0}^{\text{fb}}-y_{\text{init}}^{\infty}\right)+y_{\text{init}}^{\infty}.
\end{align}
Finally, the matching procedure must be carried out. The matching
conditions are
\begin{align}
\lim_{\tau_{L}\rightarrow\infty}Y_{L}\left(\tau_{L}\right) & =\lim_{t\rightarrow t_{0}}y_{O}\left(t\right), & \lim_{\tau_{L}\rightarrow\infty}X_{L}\left(\tau_{L}\right) & =\lim_{t\rightarrow t_{0}}x_{O}\left(t\right).\label{eq:Eq4184}
\end{align}
Equation \eqref{eq:Eq4184} immediately yields
\begin{align}
x_{\text{init}}^{\text{fb}} & =x_{0}^{\text{fb}}
\end{align}
for the initial condition of the outer equation. The remaining matching
condition for $y$ is satisfied as well. The overlaps are obtained
as
\begin{align}
y_{O}\left(t_{0}\right) & =y_{\text{init}}^{\infty}, & x_{O}\left(t_{0}\right) & =x_{0}^{\text{fb}}.
\end{align}
Inner and outer solutions are combined in a composite solution as
\begin{align}
x_{\text{comp}}\left(t\right) & =x_{O}\left(t\right)+X_{L}\left(\left(t-t_{0}\right)/\epsilon\right)-x_{0}^{\text{fb}}=x_{O}\left(t\right),\label{eq:FeedBackStateX}\\
y_{\text{comp}}\left(t\right) & =y_{\text{init}}^{\infty}\left(t\right)+Y_{L}\left(\left(t-t_{0}\right)/\epsilon\right)-y_{\text{init}}^{\infty}\nonumber \\
 & =\intop_{t}^{\infty}\left(e^{t\varphi_{1}}-e^{t_{0}\varphi_{1}}\right)e^{-\tau\varphi_{1}}\left(\frac{a_{2}s_{1}}{s_{2}}x_{d}\left(\tau\right)-\left(a_{1}+\varphi_{1}\right)y_{d}\left(\tau\right)\right)d\tau\nonumber \\
 & -\intop_{t_{0}}^{t}e^{\left(t_{0}-\tau\right)\varphi_{1}}\left(\frac{a_{2}s_{1}}{s_{2}}x_{d}\left(\tau\right)-\left(a_{1}+\varphi_{1}\right)y_{d}\left(\tau\right)\right)d\tau\nonumber \\
 & +\dfrac{1}{a_{2}}\left(\varphi_{1}+a_{1}\right)\left(x_{0}^{\text{fb}}-x\left(t\right)\right)+y_{d}\left(t\right)-y_{d}\left(t_{0}\right)+Y_{L}\left(\left(t-t_{0}\right)/\epsilon\right),\label{eq:FeedBackStateY}
\end{align}
Here, the outer solution $x_{O}\left(t\right)$ is given by Eq. \eqref{eq:xOSolFB},
while the left inner solution $Y_{L}\left(\left(t-t_{0}\right)/\epsilon\right)$
is given as the solution to Eq. \eqref{eq:YLFeedBack}.

Comparing the analytical result Eq. \eqref{eq:FeedBackControl} with
a numerical result requires the numerical solution of the Hamilton-Jacobi-Bellman
equation \eqref{eq:HJBTrajectoryTracking}. Unfortunately, this task
is rather difficult and not pursued here. Instead, selected feedback-controlled
state trajectories are compared with their open loop counterparts.
An open loop control is determined for a specified value of the initial
condition $\boldsymbol{x}_{0}$. Applying the same control signal
to another initial condition usually fails. In contrast to that, a
feedback-controlled state trajectory may start at an arbitrary initial
condition $\boldsymbol{x}_{0}^{\text{fb}}$. If $\boldsymbol{x}_{0}=\boldsymbol{x}_{0}^{\text{fb}}$,
open loop and feedback-controlled state trajectories agree. Example
\ref{ex:FBControlledFHN} investigates the impact of selected initial
conditions $\boldsymbol{x}_{0}\neq\boldsymbol{x}_{0}^{\text{fb}}$
on the feedback-controlled state trajectory.

\begin{example}[Feedback-controlled  FHN model]\label{ex:FBControlledFHN}

The feedback control Eq. \eqref{eq:FeedBackControl} is applied to
the activator-controlled FHN model, see Example \ref{ex:FHN1} for
details. The feedback-controlled state equations \eqref{eq:FeedBackControlledStatex}
and \eqref{eq:FeedBackControlledStatey} are solved numerically. The
desired trajectory and all parameters are the same as in Example \ref{ex:OptimallyControlledFHN},
except for the amplitude $A_{x}=5$ of the desired trajectory and
an infinite terminal time $t_{1}\rightarrow\infty$. Figure \ref{fig:FHNFB1}
shows the feedback-controlled state trajectory for inhibitor (left)
and activator (right) for different values of the inhibitor initial
condition $x_{0}^{\text{fb}}$. Figure \ref{fig:FHNFB2} shows state
trajectories for different values of the activator initial condition
$y_{0}^{\text{fb}}$. In all plots, the blue solid line is the controlled
state trajectory enforced by open loop control. During a transient,
the feedback-controlled trajectories converge on the open loop-controlled
state trajectory. However, the time scales of the transients are dramatically
different. A deviation of the inhibitor initial condition decays very
slowly, see Fig. \ref{fig:FHNFB1}. A deviation of the activator initial
condition $y_{0}^{\text{fb}}$ displays relaxation in form of a boundary
layer. It decays on a time scale set by the boundary layer width $\epsilon$,
see Fig. \ref{fig:FHNFB2}.

\begin{minipage}{1.0\linewidth}
\begin{center}
\includegraphics[scale=0.475]{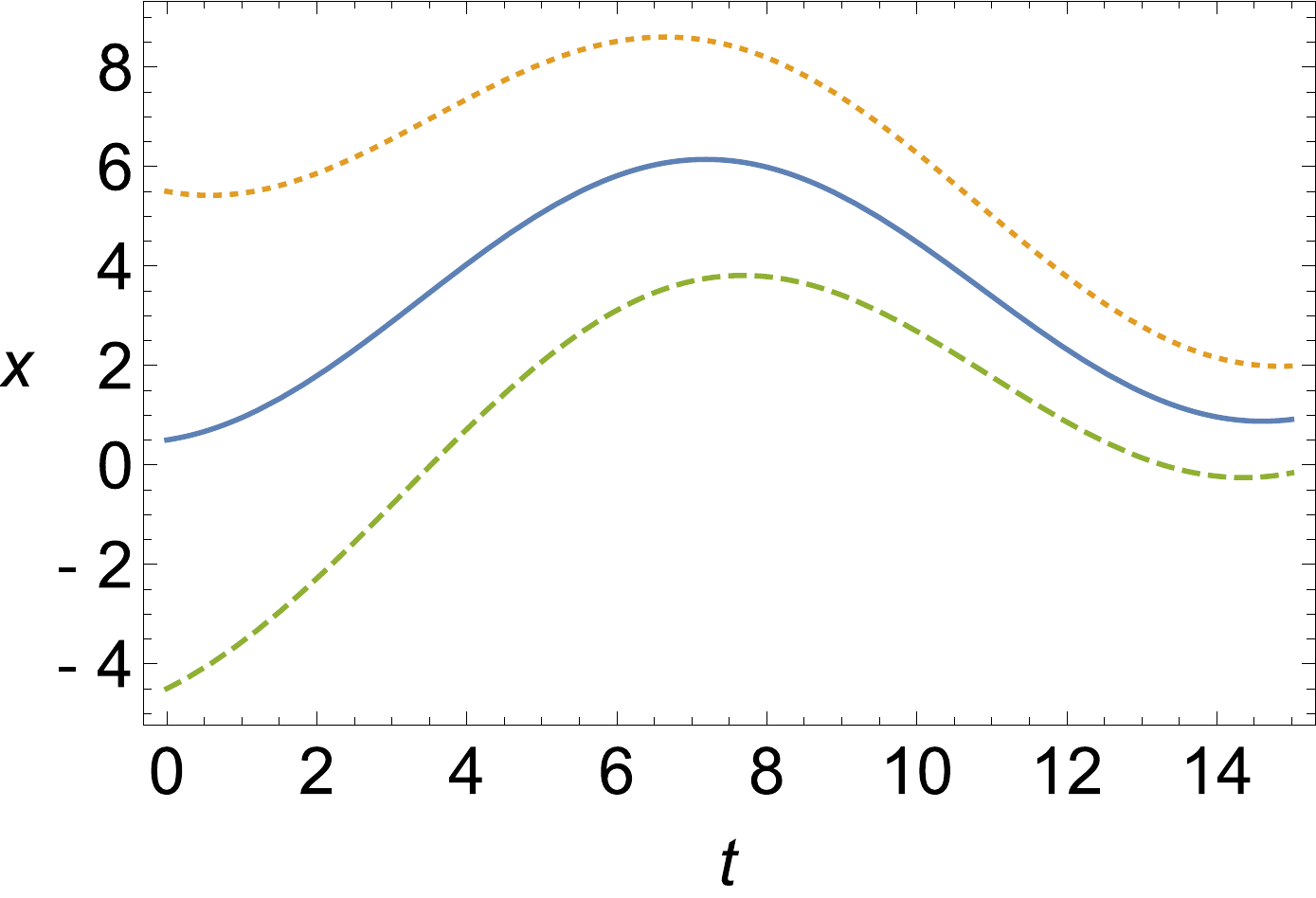}\hspace{0.2cm}\includegraphics[scale=0.49]{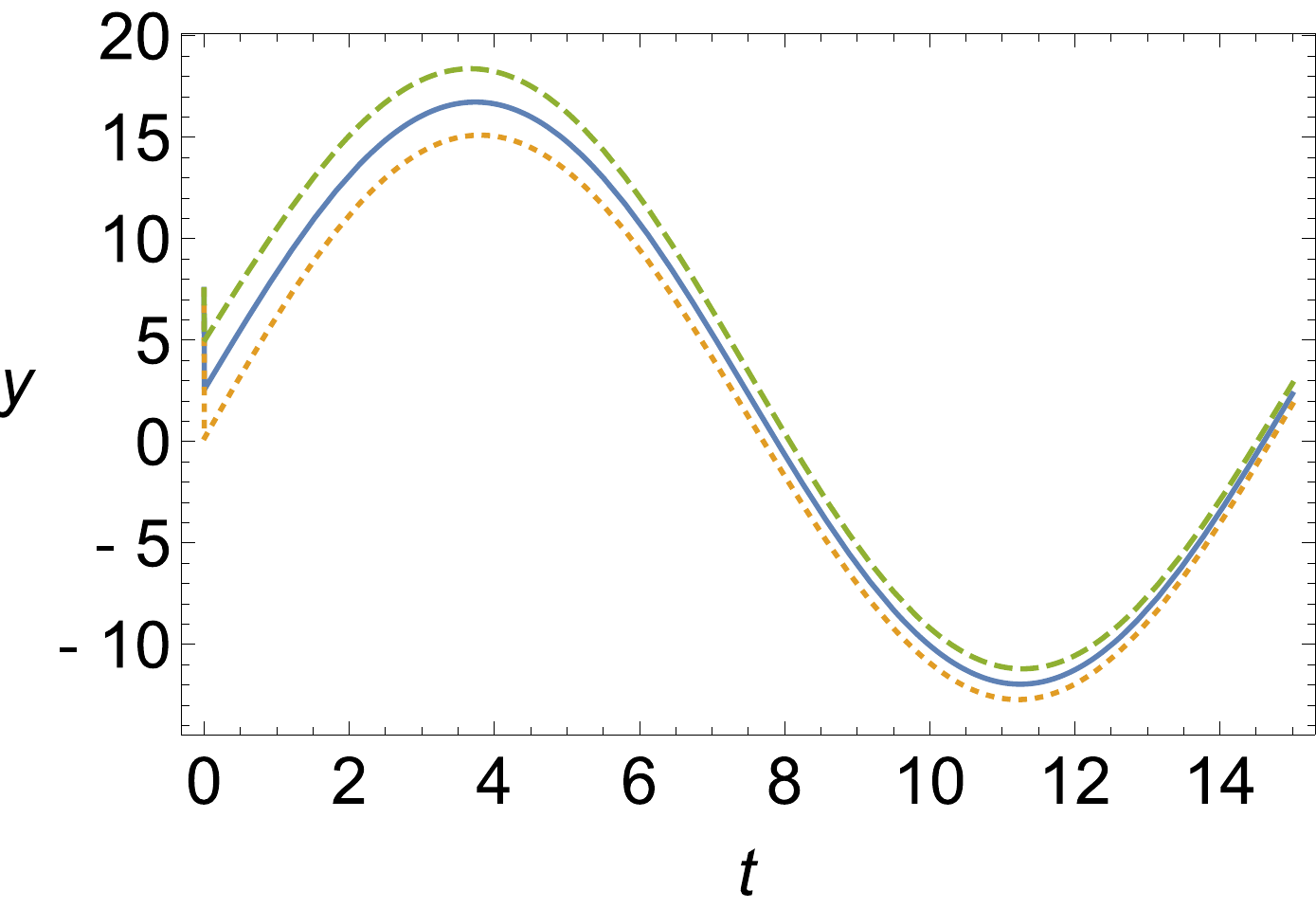}%
\captionof{figure}[Feedback-controlled trajectories for some inhibitor initial conditions]{\label{fig:FHNFB1}Feedback-controlled state trajectories for different inhibitor initial conditions. The blue solid line with initial condition $\boldsymbol{x}^{T}_{0}=\left(x_{0},\, y_{0}\right)=\left(0.5,\, 7.5\right)$ is the state trajectory enforced by open loop control.  During a long transient, the feedback-controlled trajectories for $\left(x^{\text{fb}}_{0},\, y^{\text{fb}}_{0}\right)=\left(x_{0},\, y_{0}\right) + \left(5,\, 0\right)$ (orange dashed line) and $\left(x^{\text{fb}}_{0},\, y^{\text{fb}}_{0}\right) = \left(x_{0},\, y_{0}\right) -  \left(5,\, 0\right)$ (green dashed line) converge on the blue solid trajectory.}
%"/home/jakob/ACADOtoolkit/examples/my_examples/FHN/CompareResults.nb"
\end{center}
\end{minipage}

\begin{minipage}{1.0\linewidth}
\begin{center}%
\includegraphics[scale=0.47]{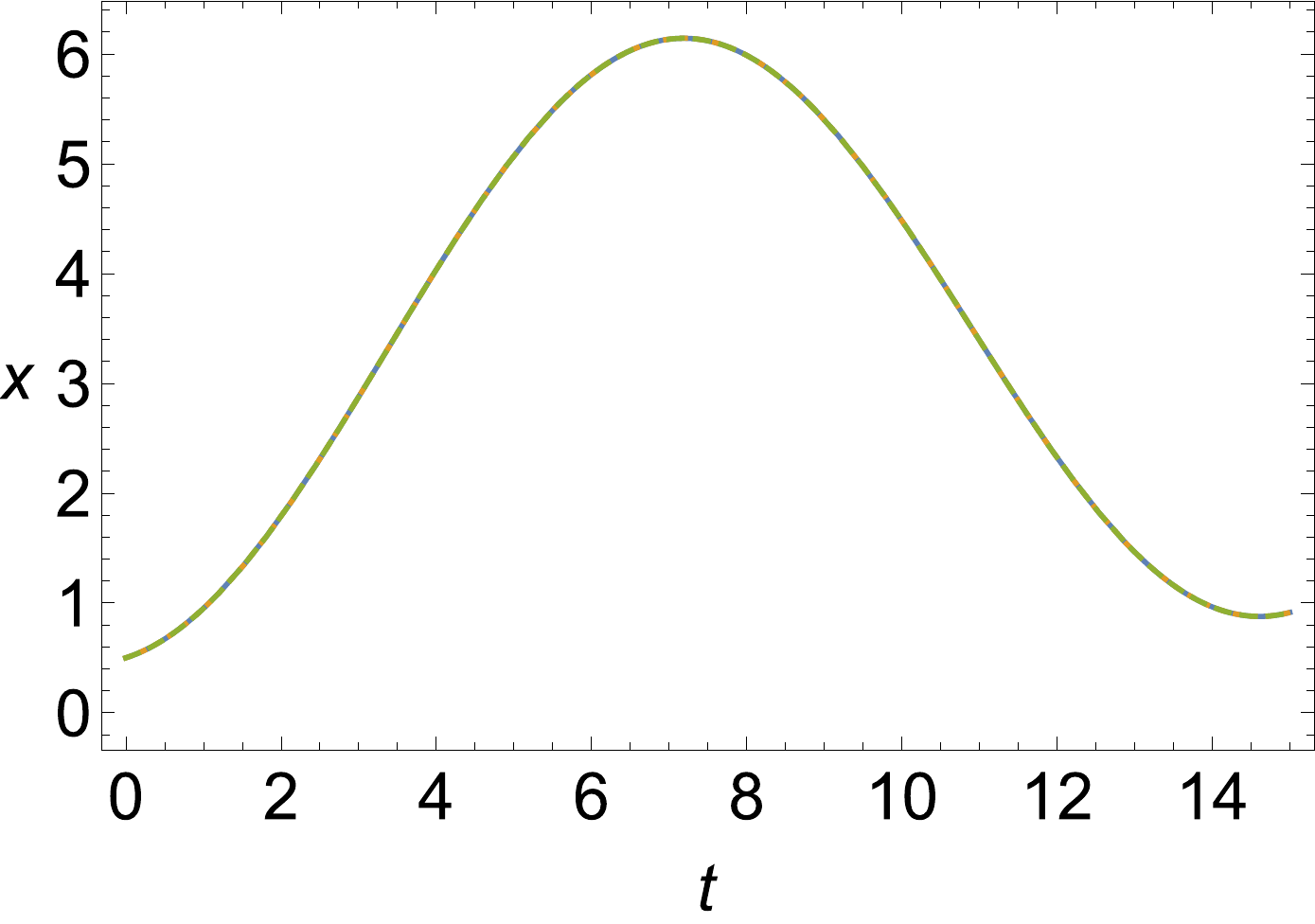}\hspace{0.2cm}\includegraphics[scale=0.5]{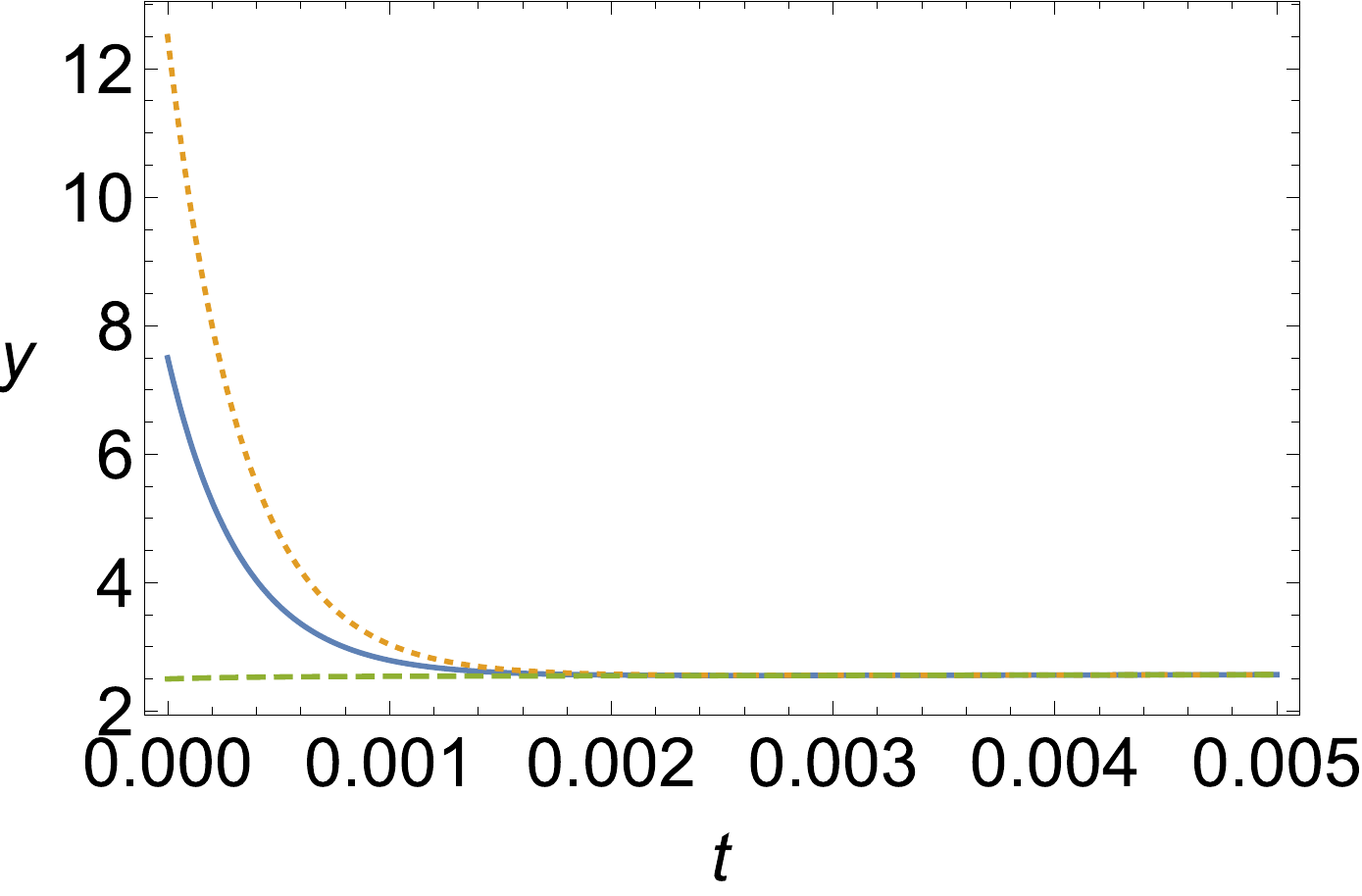}%
\captionof{figure}[Feedback-controlled trajectories for some activator initial conditions]{\label{fig:FHNFB2}Feedback-controlled state trajectories for different activator initial conditions. Compared with Fig. \ref{fig:FHNFB1}, the relaxation back on the state trajectory enforced by open loop control (blue solid line for $\left(x_{0},\, y_{0}\right) = \left(0.5,\, 7.5\right)$) occurs much faster on a time scale set by $\epsilon$. Feedback-controlled state trajectories are shown for $\left(x^{\text{fb}}_{0},\, y^{\text{fb}}_{0}\right) = \left(x_{0},\, y_{0}\right) + \left(0,\, 5\right)$ and $\left(x^{\text{fb}}_{0},\, y^{\text{fb}}_{0}\right) = \left(x_{0},\, y_{0}\right) - \left(0,\, 5\right)$ (green dashed line).}%
%"/home/jakob/ACADOtoolkit/examples/my_examples/FHN/CompareResults.nb"
\end{center}
\end{minipage}

\end{example}

\subsection{\label{sub:ContinuousTimeDelayedFeedback}Continuous time-delayed
feedback}

Continuous time feedback cannot be applied to systems with fast dynamics.
If measurement and processing of the system's state takes place on
a time scale comparable with system dynamics, a delay is induced.
The control signal fed back to the system cannot be assumed to depend
on the current state $\boldsymbol{x}\left(t\right)$ of the system.
Instead, the feedback signal depends on the delayed system state $\boldsymbol{x}\left(t-T\right)$.
The time delay $T>0$ accounts for the duration of measurement and
information processing. Introducing artificial time delays and superpositions
of continuous time and time-delayed feedback signals can be beneficial
for stabilization by feedback. A prominent example is the stabilization
of unstable periodic orbits in chaotic systems developed by Pyragas
\cite{pyragas1992continous,pyragas2006delayed}, see also the theme
issue \cite{just2010delayed}.

The initial conditions $x_{0}$ and $y_{0}$ in Eqs. \eqref{eq:xFinalSol}-\eqref{eq:lambdaFinalSol}
Eq. \eqref{eq:uFinalSol} are assumed to depend on the delayed state
components as $x_{0}=x\left(t_{0}\right)=x\left(t-T\right)$ and $y_{0}=y\left(t_{0}\right)=y\left(t-T\right)$,
respectively. Every explicit appearance of $t_{0}$ in these equations
is substituted by $t_{0}\rightarrow t-T$. All constants which depend
on time $t_{0}$, as e.g. $x_{\text{init}}$, become time dependent,
\begin{align}
x_{\text{init}} & =x_{\text{init}}\left(t\right), & y_{\text{init}} & =y_{\text{init}}\left(t\right), & y_{\text{end}} & =y_{\text{end}}\left(t\right).
\end{align}
The outer solutions $x_{O}\left(t\right)$ and $y_{O}\left(t\right)$
given by Eqs. \eqref{eq:xOSol} and \eqref{eq:yOSol} become
\begin{align}
x_{O}\left(t\right) & =\frac{a_{1}a_{2}}{\varphi_{1}}\intop_{t-T}^{t}y_{d}\left(\tau\right)\sinh\left(\varphi_{1}\left(t-\tau\right)\right)\,d\tau+a_{2}\intop_{t-T}^{t}y_{d}\left(\tau\right)\cosh\left(\varphi_{1}\left(t-\tau\right)\right)\,d\tau\nonumber \\
 & -\frac{a_{2}^{2}s_{1}}{s_{2}\varphi_{1}}\intop_{t-T}^{t}x_{d}\left(\tau\right)\sinh\left(\varphi_{1}\left(t-\tau\right)\right)\,d\tau+\frac{1}{\varphi_{1}}\sinh\left(T\varphi_{1}\right)\left(a_{0}-a_{2}y_{d}\left(t-T\right)\right)\nonumber \\
 & +\frac{a_{0}a_{1}}{\varphi_{1}^{2}}\left(\cosh\left(T\varphi_{1}\right)-1\right)+x_{\text{init}}\left(t\right)\left(\frac{a_{1}}{\varphi_{1}}\sinh\left(T\varphi_{1}\right)+\cosh\left(T\varphi_{1}\right)\right)\nonumber \\
 & +\frac{a_{2}}{\varphi_{1}}y_{\text{init}}\left(t\right)\sinh\left(T\varphi_{1}\right),\label{eq:xOSolFBDelay}
\end{align}
and
\begin{align}
y_{O}\left(t\right) & =\frac{a_{0}a_{2}s_{1}}{s_{2}\varphi_{1}^{2}}\left(\cosh\left(T\varphi_{1}\right)-1\right)+\frac{a_{2}s_{1}}{s_{2}\varphi_{1}}x_{\text{init}}\left(t\right)\sinh\left(T\varphi_{1}\right)\nonumber \\
 & +y_{\text{init}}\left(t\right)\cosh\left(T\varphi_{1}\right)-y_{\text{init}}\left(t\right)\frac{a_{1}}{\varphi_{1}}\sinh\left(T\varphi_{1}\right)\nonumber \\
 & -\frac{a_{2}s_{1}}{s_{2}}\intop_{t-T}^{t}x_{d}\left(\tau\right)\cosh\left(\varphi_{1}\left(t-\tau\right)\right)\,d\tau+\frac{a_{2}s_{1}a_{1}}{s_{2}\varphi_{1}}\intop_{t-T}^{t}x_{d}\left(\tau\right)\sinh\left(\varphi_{1}\left(t-\tau\right)\right)\,d\tau\nonumber \\
 & +\frac{a_{2}^{2}s_{1}}{s_{2}\varphi_{1}}\intop_{t-T}^{t}y_{d}\left(\tau\right)\sinh\left(\varphi_{1}\left(t-\tau\right)\right)\,d\tau\nonumber \\
 & +\frac{a_{1}}{\varphi_{1}}y_{d}\left(t-T\right)\sinh\left(T\varphi_{1}\right)+y_{d}\left(t\right)-y_{d}\left(t-T\right)\cosh\left(T\varphi_{1}\right).\label{eq:yOSolFBDelay}
\end{align}
Due to the time delay, the integral terms in $x_{O}\left(t\right)$
and $y_{O}\left(t\right)$ do not vanish. The time-dependent constants
$x_{\text{init}}\left(t\right),\,y_{\text{init}}\left(t\right)$ are
given by
\begin{align}
x_{\text{init}}\left(t\right) & =x\left(t-T\right),
\end{align}
and\allowdisplaybreaks 
\begin{align}
y_{\text{init}}\left(t\right) & =\frac{s_{2}\varphi_{1}}{\kappa\left(t\right)}\sinh\left(\left(t-T-t_{1}\right)\varphi_{1}\right)\left(a_{1}s_{2}-a_{2}^{2}\beta_{1}\right)y_{d}\left(t-T\right)\nonumber \\
 & +\frac{a_{2}s_{2}\varphi_{1}}{\kappa\left(t\right)}\sinh\left(\left(t-T-t_{1}\right)\varphi_{1}\right)\left(a_{0}\beta_{1}+x\left(t-T\right)\left(a_{1}\beta_{1}+s_{1}\right)\right)\nonumber \\
 & -\frac{a_{2}\beta_{1}}{\kappa\left(t\right)}\left(a_{1}^{2}s_{2}+a_{2}^{2}s_{1}\right)\cosh\left(\left(t_{1}-t+T\right)\varphi_{1}\right)x\left(t-T\right)\nonumber \\
 & -\frac{a_{2}a_{0}s_{2}}{\kappa\left(t\right)}\left(a_{1}\beta_{1}+s_{1}\right)\cosh\left(\left(t_{1}-t+T\right)\varphi_{1}\right)\nonumber \\
 & +\frac{s_{2}}{\kappa\left(t\right)}\left(a_{1}^{2}s_{2}+a_{2}^{2}s_{1}\right)\cosh\left(\left(t_{1}-t+T\right)\varphi_{1}\right)y_{d}\left(t-T\right)\nonumber \\
 & +\frac{a_{2}\varphi_{1}s_{1}}{\kappa\left(t\right)}\left(a_{2}^{2}\beta_{1}-a_{1}s_{2}\right)\intop_{t-T}^{t_{1}}x_{d}\left(\tau\right)\sinh\left(\varphi_{1}\left(t_{1}-\tau\right)\right)\,d\tau\nonumber \\
 & -\frac{a_{2}^{2}\varphi_{1}s_{2}}{\kappa\left(t\right)}\left(a_{1}\beta_{1}+s_{1}\right)\intop_{t-T}^{t_{1}}y_{d}\left(\tau\right)\sinh\left(\varphi_{1}\left(t_{1}-\tau\right)\right)\,d\tau\nonumber \\
 & -\frac{\beta_{1}a_{2}^{2}\varphi_{1}^{2}s_{2}}{\kappa\left(t\right)}\intop_{t-T}^{t_{1}}y_{d}\left(\tau\right)\cosh\left(\varphi_{1}\left(t_{1}-\tau\right)\right)\,d\tau+\frac{a_{2}a_{0}s_{2}}{\kappa}\left(a_{1}\beta_{1}+s_{1}\right)\nonumber \\
 & +\frac{a_{2}\varphi_{1}^{2}s_{1}s_{2}}{\kappa\left(t\right)}\intop_{t-T}^{t_{1}}x_{d}\left(\tau\right)\cosh\left(\varphi_{1}\left(t_{1}-\tau\right)\right)\,d\tau+\beta_{1}x_{1}\frac{a_{2}}{\kappa}\left(a_{1}^{2}s_{2}+a_{2}^{2}s_{1}\right).
\end{align}
\allowdisplaybreaks[0]The abbreviation $\kappa\left(t\right)$ is
defined by
\begin{align}
\kappa\left(t\right) & =s_{2}\varphi_{1}\left(a_{2}^{2}\beta_{1}-a_{1}s_{2}\right)\sinh\left(\left(t_{1}-t+T\right)\varphi_{1}\right)+s_{2}^{2}\varphi_{1}^{2}\cosh\left(\left(t_{1}-t+T\right)\varphi_{1}\right).\label{eq:kappat}
\end{align}
Using Eqs. \eqref{eq:xOSolFBDelay}-\eqref{eq:kappat} together with
the solution Eq. \eqref{eq:uFinalSol} for optimal open loop control,
the optimal time-delayed feedback control signal $u_{\text{comp}}^{\text{fb}}\left(t\right)$
is obtained as 
\begin{align}
u_{\text{comp}}^{\text{fb}}\left(t\right) & =u_{O}\left(t\right)+U_{L}\left(\left(t-t_{0}\right)/\epsilon\right)+U_{R}\left(\left(t_{1}-t\right)/\epsilon\right)-u_{O}\left(t_{0}\right)-u_{O}\left(t_{1}\right)\nonumber \\
 & =u_{O}\left(t\right)+U_{L}\left(T/\epsilon\right)+U_{R}\left(\left(t_{1}-t+T\right)/\epsilon\right)-u_{O}\left(t-T\right)-u_{O}\left(t_{1}\right).
\end{align}
Here, the outer control signal $u_{O}\left(t\right)$ is given by
Eq. \eqref{eq:OuterControl}, and $U_{L}$ and $U_{R}$ are given
by Eqs. \eqref{eq:InnerLeftControl} and \eqref{eq:InnerRightControl},
respectively. Caution has to be taken when evaluating an expression
involving the time derivative $\dot{y}_{O}\left(t\right)$. The dot
$\dot{y}_{O}\left(t\right)$ denotes the time derivative with respect
to the current time $t$. It does not commute with the substitution
$t_{0}\rightarrow t-T$. Consequently, the time derivative has to
be computed before the substitution $t_{0}\rightarrow t-T$.

\subsection{Discussion}

A modification of the analytical results from Section \ref{sec:TwoDimensionalDynamicalSystem}
extends their scope to optimal feedback control. This requires knowledge
about the dependency of the controlled state trajectory on its initial
conditions. The essential idea is to replace the initial state $\boldsymbol{x}_{0}=\boldsymbol{x}\left(t_{0}\right)$
with the monitored state $\boldsymbol{x}\left(t\right)$ of the controlled
dynamical system. Numerically, optimal feedback is obtained by solving
the Hamilton-Jacobi-Bellman equation \eqref{eq:HJBTrajectoryTracking}.

The analytical approach yields a closed form expression for the continuous
time feedback law Eq. \eqref{eq:FeedBackControl}. Remarkably, for
infinite terminal time $t_{1}\rightarrow\infty$, the feedback-controlled
state equation does not depend on the nonlinearity $R\left(x,y\right)$.
The controlled state equations depend on $\epsilon$ and are singularly
perturbed. Separating inner and outer equations reveals an initial
boundary layer for the $y$-component. Its dynamics depends on the
coupling function $b\left(x,y\right)$ and is identical in form to
the boundary layers encountered for open loop control. For identical
initial conditions, a feedback-controlled state trajectory is identical
to a state trajectory enforced by open-loop control. For differing
initial conditions, a feedback-controlled state trajectory relaxes
onto the open loop trajectory. An initial deviation of the $y$-component
converges swiftly. The relaxation is in form of a boundary layer with
a time scale set by $\epsilon$. Initial deviations of the $x$-component
decay at low speed on a time scale independent of $\epsilon$.

Stabilizing feedback control of unstable attractors received plentiful
attention by the physics community, especially in the context of chaos
control \cite{ott1990controlling,SchoellSchuster200712,Schimansky-geierFiedlerKurthsScholl200701}.
However, optimality of these methods is rarely investigated. This
is hardly surprising in view of the fundamental difficulties. Numerically
solving the Hamilton-Jacobi-Bellman equation represents a difficult
task itself. Analytical methods are largely restricted to linear systems
which lack e.g. limit cycles and chaos. The analytical approach outlined
in this section opens up a possibility to study the optimality of
continuous time as well as time-delayed feedback stabilization of
unstable attractors. An interesting option would be an optimal variant
of the Pyragas control to stabilize periodic orbits in chaotic systems
\cite{pyragas1992continous,pyragas2006delayed}. Similarly, it is
possible to investigate the optimality of techniques from mathematical
control theory. An interesting problem concerns the optimality of
feedback linearization. The continuous time feedback control given
by Eq. \eqref{eq:FeedBackControl} shares similarities with the control
signal Eq. \eqref{eq:FeedbackLinearizeControl} obtained by feedback
linearization. In both cases the control signal simply absorbs the
nonlinearity $R\left(x,y\right)$. However, the optimal feedback law
retains a nontrivial dependence on the coupling function $b\left(x,y\right)$
and results in nonlinear evolution equations for the feedback-controlled
state.

\section{\label{sec:GeneralDynamicalSystem}General dynamical system}

This section discusses optimal trajectory tracking for general dynamical
systems. The task is to minimize the target functional 
\begin{align}
\mathcal{J}\left[\boldsymbol{x}\left(t\right),\boldsymbol{u}\left(t\right)\right] & =\intop_{t_{0}}^{t_{1}}dt\frac{1}{2}\left(\boldsymbol{x}\left(t\right)-\boldsymbol{x}_{d}\left(t\right)\right)^{T}\boldsymbol{\mathcal{S}}\left(\boldsymbol{x}\left(t\right)-\boldsymbol{x}_{d}\left(t\right)\right)\nonumber \\
 & +\frac{\epsilon^{2}}{2}\intop_{t_{0}}^{t_{1}}dt\left(\boldsymbol{u}\left(t\right)\right)^{2}\label{eq:OptimalTrajectoryTrackingFunctionalGeneral}
\end{align}
subject to the dynamic constraints 
\begin{align}
\boldsymbol{\dot{x}}\left(t\right) & =\boldsymbol{R}\left(\boldsymbol{x}\left(t\right)\right)+\boldsymbol{\mathcal{B}}\left(\boldsymbol{x}\left(t\right)\right)\boldsymbol{u}\left(t\right), & \boldsymbol{x}\left(t_{0}\right) & =\boldsymbol{x}_{0}, & \boldsymbol{x}\left(t_{1}\right) & =\boldsymbol{x}_{1}.
\end{align}
Here, $\boldsymbol{\mathcal{S}}=\boldsymbol{\mathcal{S}}^{T}$ is
a symmetric $n\times n$ matrix of weights. Only sharp terminal conditions
$\boldsymbol{x}\left(t_{1}\right)=\boldsymbol{x}_{1}$ are discussed.
The starting point for the perturbative treatment are the necessary
optimality conditions,
\begin{align}
\boldsymbol{0} & =\epsilon^{2}\boldsymbol{u}\left(t\right)+\boldsymbol{\mathcal{B}}^{T}\left(\boldsymbol{x}\left(t\right)\right)\boldsymbol{\lambda}\left(t\right),\label{eq:StationarityCondition}\\
\boldsymbol{\dot{x}}\left(t\right) & =\boldsymbol{R}\left(\boldsymbol{x}\left(t\right)\right)+\boldsymbol{\mathcal{B}}\left(\boldsymbol{x}\left(t\right)\right)\boldsymbol{u}\left(t\right),\label{eq:StateEq}\\
-\boldsymbol{\dot{\lambda}}\left(t\right) & =\left(\nabla\boldsymbol{R}^{T}\left(\boldsymbol{x}\left(t\right)\right)+\boldsymbol{u}^{T}\left(t\right)\nabla\boldsymbol{\mathcal{B}}^{T}\left(\boldsymbol{x}\left(t\right)\right)\right)\boldsymbol{\lambda}\left(t\right)+\boldsymbol{\mathcal{S}}\left(\boldsymbol{x}\left(t\right)-\boldsymbol{x}_{d}\left(t\right)\right),\label{eq:AdjointEq}
\end{align}
together with the initial and terminal conditions
\begin{align}
\boldsymbol{x}\left(t_{0}\right) & =\boldsymbol{x}_{0}, & \boldsymbol{x}\left(t_{1}\right) & =\boldsymbol{x}_{1}.\label{eq:GenDynSysInitTermCond}
\end{align}
See Section \ref{sec:NecessaryOptimalityConditions} for a derivation
of Eqs. \eqref{eq:StationarityCondition} and \eqref{eq:AdjointEq}.

The idea of the analytical treatment is to utilize the two projectors
$\boldsymbol{\mathcal{P}}_{\boldsymbol{\mathcal{S}}}\left(\boldsymbol{x}\right)$
and $\boldsymbol{\mathcal{Q}}_{\boldsymbol{\mathcal{S}}}\left(\boldsymbol{x}\right)$
defined by Eqs. \eqref{eq:PSQSDef1} and \eqref{eq:PSQSDef2} to split
up the necessary optimality conditions. While the state projections
$\boldsymbol{\mathcal{P}}_{\boldsymbol{\mathcal{S}}}\left(\boldsymbol{x}\right)\boldsymbol{x}$
exhibit boundary layers, the state projections $\boldsymbol{\mathcal{Q}}_{\boldsymbol{\mathcal{S}}}\left(\boldsymbol{x}\right)\boldsymbol{x}$
do not. The equations are rearranged to obtain a singularly perturbed
system of differential equations. Inner and outer equations are determined
by a perturbation expansion to leading order of the small parameter
$\epsilon$. A linearizing assumption similar to Section \ref{sec:LinearizingAssumption}
results in linear outer equations which can be formally solved.

\subsection{Rearranging the necessary optimality conditions}

To shorten the notation, the time argument of $\boldsymbol{x}\left(t\right)$,
$\boldsymbol{\lambda}\left(t\right)$ and $\boldsymbol{u}\left(t\right)$
is suppressed in this subsection and some abbreviating matrices are
introduced. Let the $n\times n$ matrix $\boldsymbol{\Omega}_{\boldsymbol{\mathcal{S}}}\left(\boldsymbol{x}\right)$
be defined by 
\begin{align}
\boldsymbol{\Omega}_{\boldsymbol{\mathcal{S}}}\left(\boldsymbol{x}\right) & =\boldsymbol{\mathcal{B}}\left(\boldsymbol{x}\right)\left(\boldsymbol{\mathcal{B}}^{T}\left(\boldsymbol{x}\right)\boldsymbol{\mathcal{S}}\boldsymbol{\mathcal{B}}\left(\boldsymbol{x}\right)\right)^{-1}\boldsymbol{\mathcal{B}}^{T}\left(\boldsymbol{x}\right).\label{eq:DefOmega}
\end{align}
A simple calculation shows that $\boldsymbol{\Omega}_{\boldsymbol{\mathcal{S}}}\left(\boldsymbol{x}\right)$
is symmetric, 
\begin{align}
\boldsymbol{\Omega}_{\boldsymbol{\mathcal{S}}}^{T}\left(\boldsymbol{x}\right) & =\boldsymbol{\mathcal{B}}\left(\boldsymbol{x}\right)\left(\boldsymbol{\mathcal{B}}^{T}\left(\boldsymbol{x}\right)\boldsymbol{\mathcal{S}}\boldsymbol{\mathcal{B}}\left(\boldsymbol{x}\right)\right)^{-T}\boldsymbol{\mathcal{B}}^{T}\left(\boldsymbol{x}\right)\nonumber \\
 & =\boldsymbol{\mathcal{B}}\left(\boldsymbol{x}\right)\left(\boldsymbol{\mathcal{B}}^{T}\left(\boldsymbol{x}\right)\boldsymbol{\mathcal{S}}\boldsymbol{\mathcal{B}}\left(\boldsymbol{x}\right)\right)^{-1}\boldsymbol{\mathcal{B}}^{T}\left(\boldsymbol{x}\right)=\boldsymbol{\Omega}_{\boldsymbol{\mathcal{S}}}\left(\boldsymbol{x}\right).
\end{align}
Note that $\boldsymbol{\mathcal{B}}^{T}\left(\boldsymbol{x}\right)\boldsymbol{\mathcal{S}}\boldsymbol{\mathcal{B}}\left(\boldsymbol{x}\right)=\left(\boldsymbol{\mathcal{B}}^{T}\left(\boldsymbol{x}\right)\boldsymbol{\mathcal{S}}\boldsymbol{\mathcal{B}}\left(\boldsymbol{x}\right)\right)^{T}$
is a symmetric $p\times p$ matrix because $\boldsymbol{\mathcal{S}}$
is symmetric by assumption, and the inverse of a symmetric matrix
is symmetric. Let the two $n\times n$ projectors $\boldsymbol{\mathcal{P}}_{\boldsymbol{\mathcal{S}}}\left(\boldsymbol{x}\right)$
and $\boldsymbol{\mathcal{Q}}_{\boldsymbol{\mathcal{S}}}\left(\boldsymbol{x}\right)$
be defined by
\begin{align}
\boldsymbol{\mathcal{P}}_{\boldsymbol{\mathcal{S}}}\left(\boldsymbol{x}\right) & =\boldsymbol{\Omega}_{\boldsymbol{\mathcal{S}}}\left(\boldsymbol{x}\right)\boldsymbol{\mathcal{S}}, & \boldsymbol{\mathcal{Q}}_{\boldsymbol{\mathcal{S}}}\left(\boldsymbol{x}\right) & =\mathbf{1}-\boldsymbol{\mathcal{P}}_{\boldsymbol{\mathcal{S}}}\left(\boldsymbol{x}\right).
\end{align}
These projectors are derived during the discussion of singular optimal
control in Section \ref{sub:TheGeneralizedLegendreClebschConditions}.
$\boldsymbol{\mathcal{P}}_{\boldsymbol{\mathcal{S}}}\left(\boldsymbol{x}\right)$
and $\boldsymbol{\mathcal{Q}}_{\boldsymbol{\mathcal{S}}}\left(\boldsymbol{x}\right)$
are idempotent, $\boldsymbol{\mathcal{P}}_{\boldsymbol{\mathcal{S}}}^{2}\left(\boldsymbol{x}\right)=\boldsymbol{\mathcal{P}}_{\boldsymbol{\mathcal{S}}}\left(\boldsymbol{x}\right)$
and $\boldsymbol{\mathcal{Q}}_{\boldsymbol{\mathcal{S}}}^{2}\left(\boldsymbol{x}\right)=\boldsymbol{\mathcal{Q}}_{\boldsymbol{\mathcal{S}}}\left(\boldsymbol{x}\right)$.
Furthermore, $\boldsymbol{\mathcal{P}}_{\boldsymbol{\mathcal{S}}}\left(\boldsymbol{x}\right)$
and $\boldsymbol{\mathcal{Q}}_{\boldsymbol{\mathcal{S}}}\left(\boldsymbol{x}\right)$
satisfy the relations
\begin{align}
\boldsymbol{\mathcal{P}}_{\boldsymbol{\mathcal{S}}}\left(\boldsymbol{x}\right)\boldsymbol{\mathcal{B}}\left(\boldsymbol{x}\right) & =\boldsymbol{\mathcal{B}}\left(\boldsymbol{x}\right), & \boldsymbol{\mathcal{Q}}_{\boldsymbol{\mathcal{S}}}\left(\boldsymbol{x}\right)\boldsymbol{\mathcal{B}}\left(\boldsymbol{x}\right) & =\boldsymbol{0},\\
\boldsymbol{\mathcal{B}}^{T}\left(\boldsymbol{x}\right)\boldsymbol{\mathcal{S}}\boldsymbol{\mathcal{P}}_{\boldsymbol{\mathcal{S}}}\left(\boldsymbol{x}\right) & =\boldsymbol{\mathcal{B}}^{T}\left(\boldsymbol{x}\right)\boldsymbol{\mathcal{S}}, & \boldsymbol{\mathcal{B}}^{T}\left(\boldsymbol{x}\right)\boldsymbol{\mathcal{S}}\boldsymbol{\mathcal{Q}}_{\boldsymbol{\mathcal{S}}} & =\boldsymbol{0}.
\end{align}
Computing the transposed of $\boldsymbol{\mathcal{P}}_{\boldsymbol{\mathcal{S}}}\left(\boldsymbol{x}\right)$
and $\boldsymbol{\mathcal{Q}}_{\boldsymbol{\mathcal{S}}}\left(\boldsymbol{x}\right)$
yields
\begin{align}
\boldsymbol{\mathcal{P}}_{\boldsymbol{\mathcal{S}}}^{T}\left(\boldsymbol{x}\right) & =\boldsymbol{\mathcal{S}}^{T}\boldsymbol{\Omega}_{\boldsymbol{\mathcal{S}}}^{T}\left(\boldsymbol{x}\right)=\boldsymbol{\mathcal{S}}\boldsymbol{\Omega}_{\boldsymbol{\mathcal{S}}}\left(\boldsymbol{x}\right)\neq\boldsymbol{\mathcal{P}}_{\boldsymbol{\mathcal{S}}}\left(\boldsymbol{x}\right),\label{eq:ProjectorSymmetry}
\end{align}
and analogously for $\boldsymbol{\mathcal{Q}}_{\boldsymbol{\mathcal{S}}}\left(\boldsymbol{x}\right)$.
Equation \eqref{eq:ProjectorSymmetry} shows that $\boldsymbol{\mathcal{P}}_{\boldsymbol{\mathcal{S}}}\left(\boldsymbol{x}\right)$,
and therefore also $\boldsymbol{\mathcal{Q}}_{\boldsymbol{\mathcal{S}}}\left(\boldsymbol{x}\right)$,
is not symmetric. However, $\boldsymbol{\mathcal{P}}_{\boldsymbol{\mathcal{S}}}^{T}\left(\boldsymbol{x}\right)$
satisfies the convenient property 
\begin{align}
\boldsymbol{\mathcal{P}}_{\boldsymbol{\mathcal{S}}}^{T}\left(\boldsymbol{x}\right)\boldsymbol{\mathcal{S}} & =\boldsymbol{\mathcal{S}}\boldsymbol{\Omega}_{\boldsymbol{\mathcal{S}}}\left(\boldsymbol{x}\right)\boldsymbol{\mathcal{S}}=\boldsymbol{\mathcal{S}}\boldsymbol{\mathcal{P}}_{\boldsymbol{\mathcal{S}}}\left(\boldsymbol{x}\right),\label{eq:SPCommute}
\end{align}
which implies
\begin{align}
\boldsymbol{\mathcal{P}}_{\boldsymbol{\mathcal{S}}}^{T}\left(\boldsymbol{x}\right)\boldsymbol{\mathcal{S}} & =\boldsymbol{\mathcal{P}}_{\boldsymbol{\mathcal{S}}}^{T}\left(\boldsymbol{x}\right)\boldsymbol{\mathcal{P}}_{\boldsymbol{\mathcal{S}}}^{T}\left(\boldsymbol{x}\right)\boldsymbol{\mathcal{S}}=\boldsymbol{\mathcal{P}}_{\boldsymbol{\mathcal{S}}}^{T}\left(\boldsymbol{x}\right)\boldsymbol{\mathcal{S}}\boldsymbol{\mathcal{P}}_{\boldsymbol{\mathcal{S}}}\left(\boldsymbol{x}\right),\label{eq:SPCommute2}
\end{align}
and similarly for $\boldsymbol{\mathcal{S}}\boldsymbol{\mathcal{Q}}_{\boldsymbol{\mathcal{S}}}\left(\boldsymbol{x}\right)$.
The product of $\boldsymbol{\Omega}_{\boldsymbol{\mathcal{S}}}\left(\boldsymbol{x}\right)$
with $\boldsymbol{\mathcal{P}}_{\boldsymbol{\mathcal{S}}}^{T}\left(\boldsymbol{x}\right)$
yields 
\begin{align}
\boldsymbol{\Omega}_{\boldsymbol{\mathcal{S}}}\left(\boldsymbol{x}\right)\boldsymbol{\mathcal{P}}_{\boldsymbol{\mathcal{S}}}^{T}\left(\boldsymbol{x}\right) & =\boldsymbol{\Omega}_{\boldsymbol{\mathcal{S}}}\left(\boldsymbol{x}\right)\boldsymbol{\mathcal{S}}\boldsymbol{\Omega}_{\boldsymbol{\mathcal{S}}}\left(\boldsymbol{x}\right)\nonumber \\
 & =\boldsymbol{\mathcal{B}}\left(\boldsymbol{x}\right)\left(\boldsymbol{\mathcal{B}}^{T}\left(\boldsymbol{x}\right)\boldsymbol{\mathcal{S}}\boldsymbol{\mathcal{B}}\left(\boldsymbol{x}\right)\right)^{-1}\boldsymbol{\mathcal{B}}^{T}\left(\boldsymbol{x}\right)=\boldsymbol{\Omega}_{\boldsymbol{\mathcal{S}}}\left(\boldsymbol{x}\right),\label{eq:Eq4210}
\end{align}
and
\begin{align}
\boldsymbol{\Omega}_{\boldsymbol{\mathcal{S}}}\left(\boldsymbol{x}\right)\boldsymbol{\mathcal{P}}_{\boldsymbol{\mathcal{S}}}^{T}\left(\boldsymbol{x}\right)\boldsymbol{\mathcal{S}} & =\boldsymbol{\Omega}_{\boldsymbol{\mathcal{S}}}\left(\boldsymbol{x}\right)\boldsymbol{\mathcal{S}}=\boldsymbol{\mathcal{P}}_{\boldsymbol{\mathcal{S}}}\left(\boldsymbol{x}\right).\label{eq:Eq4211}
\end{align}
Let the $n\times n$ matrix $\boldsymbol{\Gamma}_{\boldsymbol{\mathcal{S}}}\left(\boldsymbol{x}\right)$
be defined by 
\begin{align}
\boldsymbol{\Gamma}_{\boldsymbol{\mathcal{S}}}\left(\boldsymbol{x}\right) & =\boldsymbol{\mathcal{S}}\boldsymbol{\mathcal{B}}\left(\boldsymbol{x}\right)\left(\boldsymbol{\mathcal{B}}^{T}\left(\boldsymbol{x}\right)\boldsymbol{\mathcal{S}}\boldsymbol{\mathcal{B}}\left(\boldsymbol{x}\right)\right)^{-2}\boldsymbol{\mathcal{B}}^{T}\left(\boldsymbol{x}\right)\boldsymbol{\mathcal{S}}.\label{eq:DefGamma}
\end{align}
$\boldsymbol{\Gamma}_{\boldsymbol{\mathcal{S}}}\left(\boldsymbol{x}\right)$
is symmetric, 
\begin{align}
\boldsymbol{\Gamma}_{\boldsymbol{\mathcal{S}}}^{T}\left(\boldsymbol{x}\right) & =\left(\boldsymbol{\mathcal{S}}\boldsymbol{\mathcal{B}}\left(\boldsymbol{x}\right)\left(\boldsymbol{\mathcal{B}}^{T}\left(\boldsymbol{x}\right)\boldsymbol{\mathcal{S}}\boldsymbol{\mathcal{B}}\left(\boldsymbol{x}\right)\right)^{-2}\boldsymbol{\mathcal{B}}^{T}\left(\boldsymbol{x}\right)\boldsymbol{\mathcal{S}}\right)^{T}\nonumber \\
 & =\boldsymbol{\mathcal{S}}\boldsymbol{\mathcal{B}}\left(\boldsymbol{x}\right)\left(\boldsymbol{\mathcal{B}}^{T}\left(\boldsymbol{x}\right)\boldsymbol{\mathcal{S}}\boldsymbol{\mathcal{B}}\left(\boldsymbol{x}\right)\right)^{-2T}\boldsymbol{\mathcal{B}}^{T}\left(\boldsymbol{x}\right)\boldsymbol{\mathcal{S}}=\boldsymbol{\Gamma}_{\boldsymbol{\mathcal{S}}}\left(\boldsymbol{x}\right),\label{eq:GammaSSymmetry}
\end{align}
and satisfies
\begin{align}
 & \boldsymbol{\Gamma}_{\boldsymbol{\mathcal{S}}}\left(\boldsymbol{x}\right)\boldsymbol{\mathcal{P}}_{\boldsymbol{\mathcal{S}}}\left(\boldsymbol{x}\right)\nonumber \\
= & \boldsymbol{\mathcal{S}}\boldsymbol{\mathcal{B}}\left(\boldsymbol{x}\right)\left(\boldsymbol{\mathcal{B}}^{T}\left(\boldsymbol{x}\right)\boldsymbol{\mathcal{S}}\boldsymbol{\mathcal{B}}\left(\boldsymbol{x}\right)\right)^{-2}\boldsymbol{\mathcal{B}}^{T}\left(\boldsymbol{x}\right)\boldsymbol{\mathcal{S}}\boldsymbol{\mathcal{B}}\left(\boldsymbol{x}\right)\left(\boldsymbol{\mathcal{B}}^{T}\left(\boldsymbol{x}\right)\boldsymbol{\mathcal{S}}\boldsymbol{\mathcal{B}}\left(\boldsymbol{x}\right)\right)^{-1}\boldsymbol{\mathcal{B}}^{T}\left(\boldsymbol{x}\right)\boldsymbol{\mathcal{S}}\nonumber \\
= & \boldsymbol{\mathcal{S}}\boldsymbol{\mathcal{B}}\left(\boldsymbol{x}\right)\left(\boldsymbol{\mathcal{B}}^{T}\left(\boldsymbol{x}\right)\boldsymbol{\mathcal{S}}\boldsymbol{\mathcal{B}}\left(\boldsymbol{x}\right)\right)^{-2}\boldsymbol{\mathcal{B}}^{T}\left(\boldsymbol{x}\right)\boldsymbol{\mathcal{S}}=\boldsymbol{\Gamma}_{\boldsymbol{\mathcal{S}}}\left(\boldsymbol{x}\right).
\end{align}
Transposing yields
\begin{align}
\left(\boldsymbol{\Gamma}_{\boldsymbol{\mathcal{S}}}\left(\boldsymbol{x}\right)\boldsymbol{\mathcal{P}}_{\boldsymbol{\mathcal{S}}}\left(\boldsymbol{x}\right)\right)^{T} & =\boldsymbol{\mathcal{P}}_{\boldsymbol{\mathcal{S}}}^{T}\left(\boldsymbol{x}\right)\boldsymbol{\Gamma}_{\boldsymbol{\mathcal{S}}}^{T}\left(\boldsymbol{x}\right)=\boldsymbol{\mathcal{P}}_{\boldsymbol{\mathcal{S}}}^{T}\left(\boldsymbol{x}\right)\boldsymbol{\Gamma}_{\boldsymbol{\mathcal{S}}}\left(\boldsymbol{x}\right)=\boldsymbol{\Gamma}_{\boldsymbol{\mathcal{S}}}\left(\boldsymbol{x}\right).\label{eq:PSGammaS}
\end{align}
The projectors $\boldsymbol{\mathcal{P}}_{\boldsymbol{\mathcal{S}}}\left(\boldsymbol{x}\right)$
and $\boldsymbol{\mathcal{Q}}_{\boldsymbol{\mathcal{S}}}\left(\boldsymbol{x}\right)$
are used to partition the state $\boldsymbol{x}$, 
\begin{align}
\boldsymbol{x} & =\boldsymbol{\mathcal{P}}_{\boldsymbol{\mathcal{S}}}\left(\boldsymbol{x}\right)\boldsymbol{x}+\boldsymbol{\mathcal{Q}}_{\boldsymbol{\mathcal{S}}}\left(\boldsymbol{x}\right)\boldsymbol{x}.\label{eq:PartitionState}
\end{align}
The controlled state equation \eqref{eq:StateEq} is split in two
parts,
\begin{align}
\boldsymbol{\mathcal{P}}_{\boldsymbol{\mathcal{S}}}\left(\boldsymbol{x}\right)\boldsymbol{\dot{x}} & =\boldsymbol{\mathcal{P}}_{\boldsymbol{\mathcal{S}}}\left(\boldsymbol{x}\right)\boldsymbol{R}\left(\boldsymbol{x}\right)+\boldsymbol{\mathcal{B}}\left(\boldsymbol{x}\right)\boldsymbol{u},\label{eq:StatePartitionP}\\
\boldsymbol{\mathcal{Q}}_{\boldsymbol{\mathcal{S}}}\left(\boldsymbol{x}\right)\boldsymbol{\dot{x}} & =\boldsymbol{\mathcal{Q}}_{\boldsymbol{\mathcal{S}}}\left(\boldsymbol{x}\right)\boldsymbol{R}\left(\boldsymbol{x}\right).\label{eq:StatePartitionQ}
\end{align}
The initial and terminal conditions are split up as well,
\begin{align}
\boldsymbol{\mathcal{P}}_{\boldsymbol{\mathcal{S}}}\left(\boldsymbol{x}\left(t_{0}\right)\right)\boldsymbol{x}\left(t_{0}\right) & =\boldsymbol{\mathcal{P}}_{\boldsymbol{\mathcal{S}}}\left(\boldsymbol{x}_{0}\right)\boldsymbol{x}_{0}, & \boldsymbol{\mathcal{Q}}_{\boldsymbol{\mathcal{S}}}\left(\boldsymbol{x}\left(t_{0}\right)\right)\boldsymbol{x}\left(t_{0}\right) & =\boldsymbol{\mathcal{Q}}_{\boldsymbol{\mathcal{S}}}\left(\boldsymbol{x}_{0}\right)\boldsymbol{x}_{0},\label{eq:PSQSInitCond}\\
\boldsymbol{\mathcal{P}}_{\boldsymbol{\mathcal{S}}}\left(\boldsymbol{x}\left(t_{1}\right)\right)\boldsymbol{x}\left(t_{1}\right) & =\boldsymbol{\mathcal{P}}_{\boldsymbol{\mathcal{S}}}\left(\boldsymbol{x}_{1}\right)\boldsymbol{x}_{1}, & \boldsymbol{\mathcal{Q}}_{\boldsymbol{\mathcal{S}}}\left(\boldsymbol{x}\left(t_{1}\right)\right)\boldsymbol{x}\left(t_{1}\right) & =\boldsymbol{\mathcal{Q}}_{\boldsymbol{\mathcal{S}}}\left(\boldsymbol{x}_{1}\right)\boldsymbol{x}_{1}.\label{eq:PSQSTermCond}
\end{align}
With the help of the relation $\boldsymbol{\mathcal{B}}^{T}\left(\boldsymbol{x}\right)\boldsymbol{\mathcal{S}}\boldsymbol{\mathcal{P}}_{\boldsymbol{\mathcal{S}}}\left(\boldsymbol{x}\right)=\boldsymbol{\mathcal{B}}^{T}\left(\boldsymbol{x}\right)\boldsymbol{\mathcal{S}}$,
Eq. \eqref{eq:StatePartitionP} is solved to obtain an expression
for the control $\boldsymbol{u}$ in terms of the controlled state
trajectory $\boldsymbol{x}$, 
\begin{align}
\boldsymbol{u} & =\left(\boldsymbol{\mathcal{B}}^{T}\left(\boldsymbol{x}\right)\boldsymbol{\mathcal{S}}\boldsymbol{\mathcal{B}}\left(\boldsymbol{x}\right)\right)^{-1}\boldsymbol{\mathcal{B}}^{T}\left(\boldsymbol{x}\right)\boldsymbol{\mathcal{S}}\left(\boldsymbol{\dot{x}}-\boldsymbol{R}\left(\boldsymbol{x}\right)\right)\nonumber \\
 & =\boldsymbol{\mathcal{B}}_{\boldsymbol{\mathcal{S}}}^{g}\left(\boldsymbol{x}\right)\left(\boldsymbol{\dot{x}}-\boldsymbol{R}\left(\boldsymbol{x}\right)\right).
\end{align}
The $p\times n$ matrix $\boldsymbol{\mathcal{B}}_{\boldsymbol{\mathcal{S}}}^{g}\left(\boldsymbol{x}\right)$
is a generalized reflexive inverse of $\boldsymbol{\mathcal{B}}\left(\boldsymbol{x}\right)$,
see Appendix \ref{sub:GeneralizedInverseMatrices}, and defined by
\begin{align}
\boldsymbol{\mathcal{B}}_{\boldsymbol{\mathcal{S}}}^{g}\left(\boldsymbol{x}\right) & =\left(\boldsymbol{\mathcal{B}}^{T}\left(\boldsymbol{x}\right)\boldsymbol{\mathcal{S}}\boldsymbol{\mathcal{B}}\left(\boldsymbol{x}\right)\right)^{-1}\boldsymbol{\mathcal{B}}^{T}\left(\boldsymbol{x}\right)\boldsymbol{\mathcal{S}}.
\end{align}
The matrix $\boldsymbol{\mathcal{B}}_{\boldsymbol{\mathcal{S}}}^{g}\left(\boldsymbol{x}\right)$
can be used to rewrite the matrices $\boldsymbol{\mathcal{P}}_{\boldsymbol{\mathcal{S}}}\left(\boldsymbol{x}\right)$
and $\boldsymbol{\Gamma}_{\boldsymbol{\mathcal{S}}}\left(\boldsymbol{x}\right)$
as
\begin{align}
\boldsymbol{\mathcal{P}}_{\boldsymbol{\mathcal{S}}}\left(\boldsymbol{x}\right) & =\boldsymbol{\mathcal{B}}\left(\boldsymbol{x}\right)\boldsymbol{\mathcal{B}}_{\boldsymbol{\mathcal{S}}}^{g}\left(\boldsymbol{x}\right), & \boldsymbol{\Gamma}_{\boldsymbol{\mathcal{S}}}\left(\boldsymbol{x}\right) & =\boldsymbol{\mathcal{B}}_{\boldsymbol{\mathcal{S}}}^{gT}\left(\boldsymbol{x}\right)\boldsymbol{\mathcal{B}}_{\boldsymbol{\mathcal{S}}}^{g}\left(\boldsymbol{x}\right),\label{eq:Eq4222}
\end{align}
respectively. The $n\times p$ matrix $\boldsymbol{\mathcal{B}}_{\boldsymbol{\mathcal{S}}}^{gT}\left(\boldsymbol{x}\right)$
is the transposed of $\boldsymbol{\mathcal{B}}_{\boldsymbol{\mathcal{S}}}^{g}\left(\boldsymbol{x}\right)$.

The solution for $\boldsymbol{u}$ is inserted in the stationarity
condition Eq. \eqref{eq:StationarityCondition} to yield
\begin{align}
\boldsymbol{0} & =\epsilon^{2}\boldsymbol{u}^{T}+\boldsymbol{\lambda}^{T}\boldsymbol{\mathcal{B}}\left(\boldsymbol{x}\right)\nonumber \\
 & =\epsilon^{2}\left(\boldsymbol{\dot{x}}^{T}-\boldsymbol{R}^{T}\left(\boldsymbol{x}\right)\right)\boldsymbol{\mathcal{B}}_{\boldsymbol{\mathcal{S}}}^{gT}\left(\boldsymbol{x}\right)+\boldsymbol{\lambda}^{T}\boldsymbol{\mathcal{B}}\left(\boldsymbol{x}\right).\label{eq:StationarityCondition_1}
\end{align}
Equation \eqref{eq:StationarityCondition_1} is utilized to eliminate
any occurrence of the part $\boldsymbol{\mathcal{P}}_{\boldsymbol{\mathcal{S}}}^{T}\left(\boldsymbol{x}\right)\boldsymbol{\lambda}$
in all equations. In contrast to the state $\boldsymbol{x}$, cf.
Eq. \eqref{eq:PartitionState}, the co-state is split up with the
transposed projectors $\boldsymbol{\mathcal{P}}_{\boldsymbol{\mathcal{S}}}^{T}\left(\boldsymbol{x}\right)$
and $\boldsymbol{\mathcal{Q}}_{\boldsymbol{\mathcal{S}}}^{T}\left(\boldsymbol{x}\right)$,
\begin{align}
\boldsymbol{\lambda} & =\boldsymbol{\mathcal{P}}_{\boldsymbol{\mathcal{S}}}^{T}\left(\boldsymbol{x}\right)\boldsymbol{\lambda}+\boldsymbol{\mathcal{Q}}_{\boldsymbol{\mathcal{S}}}^{T}\left(\boldsymbol{x}\right)\boldsymbol{\lambda}.\label{eq:PartitionCoState}
\end{align}
Multiplying Eq. \eqref{eq:StationarityCondition_1} with $\boldsymbol{\mathcal{B}}_{\boldsymbol{\mathcal{S}}}^{g}\left(\boldsymbol{x}\right)$
from the right and using Eq. \eqref{eq:Eq4222} yields an expression
for $\boldsymbol{\mathcal{P}}_{\boldsymbol{\mathcal{S}}}^{T}\left(\boldsymbol{x}\right)\boldsymbol{\lambda}$,
\begin{align}
\boldsymbol{0} & =\epsilon^{2}\left(\boldsymbol{\dot{x}}^{T}-\boldsymbol{R}^{T}\left(\boldsymbol{x}\right)\right)\boldsymbol{\Gamma}_{\boldsymbol{\mathcal{S}}}\left(\boldsymbol{x}\right)+\boldsymbol{\lambda}^{T}\boldsymbol{\mathcal{P}}_{\boldsymbol{\mathcal{S}}}\left(\boldsymbol{x}\right).
\end{align}
Transposing the last equation and exploiting the symmetry of $\boldsymbol{\Gamma}_{\boldsymbol{\mathcal{S}}}\left(\boldsymbol{x}\right)$,
Eq. \eqref{eq:GammaSSymmetry}, yields 
\begin{align}
\boldsymbol{\mathcal{P}}_{\boldsymbol{\mathcal{S}}}^{T}\left(\boldsymbol{x}\right)\boldsymbol{\lambda} & =-\epsilon^{2}\boldsymbol{\Gamma}_{\boldsymbol{\mathcal{S}}}\left(\boldsymbol{x}\right)\left(\boldsymbol{\dot{x}}-\boldsymbol{R}\left(\boldsymbol{x}\right)\right).\label{eq:PSLambda}
\end{align}
Equation \eqref{eq:PSLambda} is valid for all times $t_{0}\leq t\leq t_{1}$.
Applying the time derivative gives
\begin{align}
\boldsymbol{0} & =\epsilon^{2}\boldsymbol{\Gamma}_{\boldsymbol{\mathcal{S}}}\left(\boldsymbol{x}\right)\left(\boldsymbol{\ddot{x}}-\nabla\boldsymbol{R}\left(\boldsymbol{x}\right)\boldsymbol{\dot{x}}\right)+\epsilon^{2}\boldsymbol{\dot{\Gamma}}_{\boldsymbol{\mathcal{S}}}\left(\boldsymbol{x}\right)\left(\boldsymbol{\dot{x}}-\boldsymbol{R}\left(\boldsymbol{x}\right)\right)\nonumber \\
 & +\boldsymbol{\mathcal{\dot{P}}}_{\boldsymbol{\mathcal{S}}}^{T}\left(\boldsymbol{x}\right)\boldsymbol{\lambda}+\boldsymbol{\mathcal{P}}_{\boldsymbol{\mathcal{S}}}^{T}\left(\boldsymbol{x}\right)\boldsymbol{\dot{\lambda}}.\label{eq:PSLambdaDot_11}
\end{align}
The short hand notations
\begin{align}
\boldsymbol{\dot{\Gamma}}_{\boldsymbol{\mathcal{S}}}\left(\boldsymbol{x}\right) & =\nabla\boldsymbol{\Gamma}_{\boldsymbol{\mathcal{S}}}\left(\boldsymbol{x}\right)\boldsymbol{\dot{x}}, & \boldsymbol{\mathcal{\dot{P}}}_{\boldsymbol{\mathcal{S}}}^{T}\left(\boldsymbol{x}\right) & =\left(\boldsymbol{\dot{x}}^{T}\nabla\right)\boldsymbol{\mathcal{P}}_{\boldsymbol{\mathcal{S}}}^{T}\left(\boldsymbol{x}\right),
\end{align}
were introduced in Eq. \eqref{eq:PSLambdaDot_11}. Splitting the
co-state $\boldsymbol{\lambda}$ as in Eq. \eqref{eq:PartitionCoState}
and using Eq. \eqref{eq:PSLambda} to eliminate $\boldsymbol{\mathcal{P}}_{\boldsymbol{\mathcal{S}}}^{T}\left(\boldsymbol{x}\right)\boldsymbol{\lambda}$
leads to 
\begin{align}
-\boldsymbol{\mathcal{P}}_{\boldsymbol{\mathcal{S}}}^{T}\left(\boldsymbol{x}\right)\boldsymbol{\dot{\lambda}} & =\epsilon^{2}\boldsymbol{\Gamma}_{\boldsymbol{\mathcal{S}}}\left(\boldsymbol{x}\right)\left(\boldsymbol{\ddot{x}}-\nabla\boldsymbol{R}\left(\boldsymbol{x}\right)\boldsymbol{\dot{x}}\right)+\boldsymbol{\mathcal{\dot{P}}}_{\boldsymbol{\mathcal{S}}}^{T}\left(\boldsymbol{x}\right)\boldsymbol{\mathcal{Q}}_{\boldsymbol{\mathcal{S}}}^{T}\left(\boldsymbol{x}\right)\boldsymbol{\lambda}\nonumber \\
 & +\epsilon^{2}\left(\boldsymbol{\dot{\Gamma}}_{\boldsymbol{\mathcal{S}}}\left(\boldsymbol{x}\right)-\boldsymbol{\mathcal{\dot{P}}}_{\boldsymbol{\mathcal{S}}}^{T}\left(\boldsymbol{x}\right)\boldsymbol{\Gamma}_{\boldsymbol{\mathcal{S}}}\left(\boldsymbol{x}\right)\right)\left(\boldsymbol{\dot{x}}-\boldsymbol{R}\left(\boldsymbol{x}\right)\right).\label{eq:PSLambdaDot_12}
\end{align}
Equation is an expression for $\boldsymbol{\mathcal{P}}_{\boldsymbol{\mathcal{S}}}^{T}\left(\boldsymbol{x}\right)\boldsymbol{\dot{\lambda}}$
independent of $\boldsymbol{\mathcal{P}}_{\boldsymbol{\mathcal{S}}}^{T}\left(\boldsymbol{x}\right)\boldsymbol{\lambda}$.

A similar procedure is performed for the adjoint equation \eqref{eq:AdjointEq}.
Eliminating the control signal $\boldsymbol{u}$ from Eq. \eqref{eq:AdjointEq}
gives 
\begin{align}
-\boldsymbol{\dot{\lambda}} & =\left(\nabla\boldsymbol{R}^{T}\left(\boldsymbol{x}\right)+\left(\boldsymbol{\dot{x}}^{T}-\boldsymbol{R}^{T}\left(\boldsymbol{x}\right)\right)\boldsymbol{\mathcal{B}}_{\boldsymbol{\mathcal{S}}}^{gT}\left(\boldsymbol{x}\right)\nabla\boldsymbol{\mathcal{B}}^{T}\left(\boldsymbol{x}\right)\right)\boldsymbol{\lambda}+\boldsymbol{\mathcal{S}}\left(\boldsymbol{x}-\boldsymbol{x}_{d}\left(t\right)\right).\label{eq:Eq4230}
\end{align}
The expression $\boldsymbol{\mathcal{B}}_{\boldsymbol{\mathcal{S}}}^{gT}\left(\boldsymbol{x}\right)\nabla\boldsymbol{\mathcal{B}}^{T}\left(\boldsymbol{x}\right)$
is a third order tensor with $n\times n\times n$ components defined
as (see also Eq. \eqref{eq:lambdanablaBu} for the meaning of $\nabla\boldsymbol{\mathcal{B}}\left(\boldsymbol{x}\right)$)\emph{
}
\begin{align}
\left(\nabla\boldsymbol{\mathcal{B}}\left(\boldsymbol{x}\right)\boldsymbol{\mathcal{B}}_{\boldsymbol{\mathcal{S}}}^{g}\left(\boldsymbol{x}\right)\right)_{ijk} & =\sum_{l=1}^{p}\dfrac{\partial}{\partial x_{j}}\mathcal{B}_{il}\left(\boldsymbol{x}\right)\mathcal{B}_{\boldsymbol{\mathcal{S}},lk}^{g}\left(\boldsymbol{x}\right).
\end{align}
The product $\boldsymbol{\lambda}^{T}\nabla\boldsymbol{\mathcal{B}}\left(\boldsymbol{x}\right)\boldsymbol{\mathcal{B}}_{\boldsymbol{\mathcal{S}}}^{g}\left(\boldsymbol{x}\right)\boldsymbol{x}$
yields an $n$-component row vector with entries
\begin{align}
\left(\boldsymbol{\lambda}^{T}\nabla\boldsymbol{\mathcal{B}}\left(\boldsymbol{x}\right)\boldsymbol{\mathcal{B}}_{\boldsymbol{\mathcal{S}}}^{g}\left(\boldsymbol{x}\right)\boldsymbol{x}\right)_{j} & =\sum_{i=1}^{n}\sum_{k=1}^{n}\sum_{l=1}^{p}\lambda_{i}\dfrac{\partial}{\partial x_{j}}\mathcal{B}_{il}\left(\boldsymbol{x}\right)\mathcal{B}_{\boldsymbol{\mathcal{S}},lk}^{g}\left(\boldsymbol{x}\right)x_{k}.
\end{align}
Tho shorten the notation, the $n\times n$ matrix
\begin{align}
\boldsymbol{\mathcal{W}}\left(\boldsymbol{x},\boldsymbol{y}\right) & =\nabla\boldsymbol{\mathcal{B}}\left(\boldsymbol{x}\right)\boldsymbol{\mathcal{B}}_{\boldsymbol{\mathcal{S}}}^{g}\left(\boldsymbol{x}\right)\left(\boldsymbol{y}-\boldsymbol{R}\left(\boldsymbol{x}\right)\right)
\end{align}
with entries
\begin{align}
\mathcal{W}_{ij}\left(\boldsymbol{x},\boldsymbol{y}\right) & =\sum_{k=1}^{n}\sum_{l=1}^{p}\dfrac{\partial}{\partial x_{j}}\mathcal{B}_{il}\left(\boldsymbol{x}\right)\mathcal{B}_{\boldsymbol{\mathcal{S}},lk}^{g}\left(\boldsymbol{x}\right)\left(y_{k}-R_{k}\left(\boldsymbol{x}\right)\right)
\end{align}
is introduced. Using Eq. \eqref{eq:Eq4222} yields the identity
\begin{align}
\dfrac{\partial}{\partial x_{j}}\boldsymbol{\mathcal{B}}\left(\boldsymbol{x}\right)\boldsymbol{\mathcal{B}}_{\boldsymbol{\mathcal{S}}}^{g}\left(\boldsymbol{x}\right) & =\dfrac{\partial}{\partial x_{j}}\boldsymbol{\mathcal{P}}_{\boldsymbol{\mathcal{S}}}\left(\boldsymbol{x}\right)-\boldsymbol{\mathcal{B}}\left(\boldsymbol{x}\right)\dfrac{\partial}{\partial x_{j}}\boldsymbol{\mathcal{B}}_{\boldsymbol{\mathcal{S}}}^{g}\left(\boldsymbol{x}\right),
\end{align}
and the entries of $\boldsymbol{\mathcal{W}}\left(\boldsymbol{x},\boldsymbol{y}\right)$
can be expressed as 
\begin{align}
\mathcal{W}_{ij}\left(\boldsymbol{x},\boldsymbol{y}\right) & =\sum_{k=1}^{n}\dfrac{\partial}{\partial x_{j}}\mathcal{P}_{\boldsymbol{\mathcal{S}},ik}\left(\boldsymbol{x}\right)\left(y_{k}-R_{k}\left(\boldsymbol{x}\right)\right)\nonumber \\
 & -\sum_{k=1}^{n}\sum_{l=1}^{p}\mathcal{B}_{il}\left(\boldsymbol{x}\right)\dfrac{\partial}{\partial x_{j}}\mathcal{B}_{\boldsymbol{\mathcal{S}},lk}^{g}\left(\boldsymbol{x}\right)\left(y_{k}-R_{k}\left(\boldsymbol{x}\right)\right).\label{eq:Wentries}
\end{align}
Because of $\boldsymbol{\mathcal{Q}}_{\boldsymbol{\mathcal{S}}}\left(\boldsymbol{x}\right)\boldsymbol{\mathcal{B}}\left(\boldsymbol{x}\right)=\boldsymbol{0}$,
the product $\boldsymbol{\mathcal{Q}}_{\boldsymbol{\mathcal{S}}}\left(\boldsymbol{x}\right)\boldsymbol{\mathcal{W}}\left(\boldsymbol{x},\boldsymbol{y}\right)$
is
\begin{align}
\sum_{i=1}^{n}\mathcal{Q}_{\boldsymbol{\mathcal{S}},li}\left(\boldsymbol{x}\right)\mathcal{W}_{ij}\left(\boldsymbol{x},\boldsymbol{y}\right) & =\sum_{i=1}^{n}\sum_{k=1}^{n}\mathcal{Q}_{\boldsymbol{\mathcal{S}},li}\left(\boldsymbol{x}\right)\dfrac{\partial}{\partial x_{j}}\mathcal{P}_{\boldsymbol{\mathcal{S}},ik}\left(\boldsymbol{x}\right)\left(y_{k}-R_{k}\left(\boldsymbol{x}\right)\right).\label{eq:QTimesW}
\end{align}
In terms of the matrix $\boldsymbol{\mathcal{W}}\left(\boldsymbol{x},\boldsymbol{y}\right)$,
Eq. \eqref{eq:Eq4230} assumes the shorter form 
\begin{align}
-\boldsymbol{\dot{\lambda}} & =\left(\nabla\boldsymbol{R}^{T}\left(\boldsymbol{x}\right)+\boldsymbol{\mathcal{W}}^{T}\left(\boldsymbol{x},\boldsymbol{\dot{x}}\right)\right)\boldsymbol{\lambda}+\boldsymbol{\mathcal{S}}\left(\boldsymbol{x}-\boldsymbol{x}_{d}\left(t\right)\right).\label{eq:Eq4230-1}
\end{align}
With the help of the projectors $\boldsymbol{\mathcal{P}}_{\boldsymbol{\mathcal{S}}}^{T}$
and $\boldsymbol{\mathcal{Q}}_{\boldsymbol{\mathcal{S}}}^{T}$, Eq.
\eqref{eq:Eq4230-1} is split up in two parts,
\begin{align}
-\boldsymbol{\mathcal{P}}_{\boldsymbol{\mathcal{S}}}^{T}\left(\boldsymbol{x}\right)\boldsymbol{\dot{\lambda}} & =\boldsymbol{\mathcal{P}}_{\boldsymbol{\mathcal{S}}}^{T}\left(\boldsymbol{x}\right)\left(\nabla\boldsymbol{R}^{T}\left(\boldsymbol{x}\right)+\boldsymbol{\mathcal{W}}^{T}\left(\boldsymbol{x},\boldsymbol{\dot{x}}\right)\right)\boldsymbol{\lambda}+\boldsymbol{\mathcal{P}}_{\boldsymbol{\mathcal{S}}}^{T}\left(\boldsymbol{x}\right)\boldsymbol{\mathcal{S}}\left(\boldsymbol{x}-\boldsymbol{x}_{d}\left(t\right)\right),\label{eq:PSLambdaDot_21}\\
-\boldsymbol{\mathcal{Q}}_{\boldsymbol{\mathcal{S}}}^{T}\left(\boldsymbol{x}\right)\boldsymbol{\dot{\lambda}} & =\boldsymbol{\mathcal{Q}}_{\boldsymbol{\mathcal{S}}}^{T}\left(\boldsymbol{x}\right)\left(\nabla\boldsymbol{R}^{T}\left(\boldsymbol{x}\right)+\boldsymbol{\mathcal{W}}^{T}\left(\boldsymbol{x},\boldsymbol{\dot{x}}\right)\right)\boldsymbol{\lambda}+\boldsymbol{\mathcal{Q}}_{\boldsymbol{\mathcal{S}}}^{T}\left(\boldsymbol{x}\right)\boldsymbol{\mathcal{S}}\left(\boldsymbol{x}-\boldsymbol{x}_{d}\left(t\right)\right).\label{eq:QSLambdaDot_21}
\end{align}
Using Eq. \eqref{eq:PSLambda} to eliminate $\boldsymbol{\mathcal{P}}_{\boldsymbol{\mathcal{S}}}^{T}\left(\boldsymbol{x}\right)\boldsymbol{\lambda}$
in Eqs. \eqref{eq:PSLambdaDot_21} and \eqref{eq:QSLambdaDot_21}
results in 
\begin{align}
-\boldsymbol{\mathcal{P}}_{\boldsymbol{\mathcal{S}}}^{T}\left(\boldsymbol{x}\right)\boldsymbol{\dot{\lambda}} & =\boldsymbol{\mathcal{P}}_{\boldsymbol{\mathcal{S}}}^{T}\left(\boldsymbol{x}\right)\left(\nabla\boldsymbol{R}^{T}\left(\boldsymbol{x}\right)+\boldsymbol{\mathcal{W}}^{T}\left(\boldsymbol{x},\boldsymbol{\dot{x}}\right)\right)\left(\boldsymbol{\mathcal{Q}}_{\boldsymbol{\mathcal{S}}}^{T}\left(\boldsymbol{x}\right)\boldsymbol{\lambda}-\epsilon^{2}\boldsymbol{\Gamma}_{\boldsymbol{\mathcal{S}}}\left(\boldsymbol{x}\right)\left(\boldsymbol{\dot{x}}-\boldsymbol{R}\left(\boldsymbol{x}\right)\right)\right)\nonumber \\
 & +\boldsymbol{\mathcal{P}}_{\boldsymbol{\mathcal{S}}}^{T}\left(\boldsymbol{x}\right)\boldsymbol{\mathcal{S}}\left(\boldsymbol{x}-\boldsymbol{x}_{d}\left(t\right)\right),\label{eq:PSLambdaDot_22}\\
-\boldsymbol{\mathcal{Q}}_{\boldsymbol{\mathcal{S}}}^{T}\left(\boldsymbol{x}\right)\boldsymbol{\dot{\lambda}} & =\boldsymbol{\mathcal{Q}}_{\boldsymbol{\mathcal{S}}}^{T}\left(\boldsymbol{x}\right)\left(\nabla\boldsymbol{R}^{T}\left(\boldsymbol{x}\right)+\boldsymbol{\mathcal{W}}^{T}\left(\boldsymbol{x},\boldsymbol{\dot{x}}\right)\right)\left(\boldsymbol{\mathcal{Q}}_{\boldsymbol{\mathcal{S}}}^{T}\left(\boldsymbol{x}\right)\boldsymbol{\lambda}-\epsilon^{2}\boldsymbol{\Gamma}_{\boldsymbol{\mathcal{S}}}\left(\boldsymbol{x}\right)\left(\boldsymbol{\dot{x}}-\boldsymbol{R}\left(\boldsymbol{x}\right)\right)\right)\nonumber \\
 & +\boldsymbol{\mathcal{Q}}_{\boldsymbol{\mathcal{S}}}^{T}\left(\boldsymbol{x}\right)\boldsymbol{\mathcal{S}}\left(\boldsymbol{x}-\boldsymbol{x}_{d}\left(t\right)\right).\label{eq:eq:QSLambdaDot_22}
\end{align}
Equations \eqref{eq:PSLambdaDot_12} and \eqref{eq:PSLambdaDot_22}
are two independent expressions for $\boldsymbol{\mathcal{P}}_{\boldsymbol{\mathcal{S}}}^{T}\left(\boldsymbol{x}\right)\boldsymbol{\dot{\lambda}}$.
Combining them yields a second order differential equation independent
of $\boldsymbol{\mathcal{P}}_{\boldsymbol{\mathcal{S}}}^{T}\left(\boldsymbol{x}\right)\boldsymbol{\lambda}$,
\begin{align}
\epsilon^{2}\boldsymbol{\Gamma}_{\boldsymbol{\mathcal{S}}}\left(\boldsymbol{x}\right)\boldsymbol{\ddot{x}} & =\epsilon^{2}\boldsymbol{\Gamma}_{\boldsymbol{\mathcal{S}}}\left(\boldsymbol{x}\right)\nabla\boldsymbol{R}\left(\boldsymbol{x}\right)\boldsymbol{\dot{x}}-\epsilon^{2}\left(\boldsymbol{\dot{\Gamma}}_{\boldsymbol{\mathcal{S}}}\left(\boldsymbol{x}\right)-\boldsymbol{\mathcal{\dot{P}}}_{\boldsymbol{\mathcal{S}}}^{T}\left(\boldsymbol{x}\right)\boldsymbol{\Gamma}_{\boldsymbol{\mathcal{S}}}\left(\boldsymbol{x}\right)\right)\left(\boldsymbol{\dot{x}}-\boldsymbol{R}\left(\boldsymbol{x}\right)\right)\nonumber \\
 & +\boldsymbol{\mathcal{P}}_{\boldsymbol{\mathcal{S}}}^{T}\left(\boldsymbol{x}\right)\left(\nabla\boldsymbol{R}^{T}\left(\boldsymbol{x}\right)+\boldsymbol{\mathcal{W}}^{T}\left(\boldsymbol{x},\boldsymbol{\dot{x}}\right)\right)\left(\boldsymbol{\mathcal{Q}}_{\boldsymbol{\mathcal{S}}}^{T}\left(\boldsymbol{x}\right)\boldsymbol{\lambda}-\epsilon^{2}\boldsymbol{\Gamma}_{\boldsymbol{\mathcal{S}}}\left(\boldsymbol{x}\right)\left(\boldsymbol{\dot{x}}-\boldsymbol{R}\left(\boldsymbol{x}\right)\right)\right)\nonumber \\
 & -\boldsymbol{\mathcal{\dot{P}}}_{\boldsymbol{\mathcal{S}}}^{T}\left(\boldsymbol{x}\right)\boldsymbol{\mathcal{Q}}_{\boldsymbol{\mathcal{S}}}^{T}\left(\boldsymbol{x}\right)\boldsymbol{\lambda}+\boldsymbol{\mathcal{P}}_{\boldsymbol{\mathcal{S}}}^{T}\left(\boldsymbol{x}\right)\boldsymbol{\mathcal{S}}\left(\boldsymbol{x}-\boldsymbol{x}_{d}\left(t\right)\right).\label{eq:PSLambdaDot_3}
\end{align}
Equation \eqref{eq:PSLambdaDot_3} contains several time dependent
matrices which can be simplified. From Eq. \eqref{eq:PSGammaS} follows
for the time derivative of $\boldsymbol{\Gamma}_{\boldsymbol{\mathcal{S}}}\left(\boldsymbol{x}\right)$
\begin{align}
\boldsymbol{\dot{\Gamma}}_{\boldsymbol{\mathcal{S}}}\left(\boldsymbol{x}\right) & =\boldsymbol{\mathcal{\dot{P}}}_{\boldsymbol{\mathcal{S}}}^{T}\left(\boldsymbol{x}\right)\boldsymbol{\Gamma}_{\boldsymbol{\mathcal{S}}}\left(\boldsymbol{x}\right)+\boldsymbol{\mathcal{P}}_{\boldsymbol{\mathcal{S}}}^{T}\left(\boldsymbol{x}\right)\boldsymbol{\dot{\Gamma}}_{\boldsymbol{\mathcal{S}}}\left(\boldsymbol{x}\right)
\end{align}
such that
\begin{align}
\boldsymbol{\dot{\Gamma}}_{\boldsymbol{\mathcal{S}}}\left(\boldsymbol{x}\right)-\boldsymbol{\mathcal{\dot{P}}}_{\boldsymbol{\mathcal{S}}}^{T}\left(\boldsymbol{x}\right)\boldsymbol{\Gamma}_{\boldsymbol{\mathcal{S}}}\left(\boldsymbol{x}\right) & =\boldsymbol{\mathcal{P}}_{\boldsymbol{\mathcal{S}}}^{T}\left(\boldsymbol{x}\right)\boldsymbol{\dot{\Gamma}}_{\boldsymbol{\mathcal{S}}}\left(\boldsymbol{x}\right).\label{eq:Eq221}
\end{align}
Furthermore, from
\begin{align}
\boldsymbol{\mathcal{P}}_{\boldsymbol{\mathcal{S}}}^{T}\left(\boldsymbol{x}\right)\boldsymbol{\mathcal{Q}}_{\boldsymbol{\mathcal{S}}}^{T}\left(\boldsymbol{x}\right) & =\boldsymbol{0}
\end{align}
follows
\begin{align}
\boldsymbol{\mathcal{\dot{P}}}_{\boldsymbol{\mathcal{S}}}^{T}\left(\boldsymbol{x}\right)\boldsymbol{\mathcal{Q}}_{\boldsymbol{\mathcal{S}}}^{T}\left(\boldsymbol{x}\right) & =-\boldsymbol{\mathcal{P}}_{\boldsymbol{\mathcal{S}}}^{T}\left(\boldsymbol{x}\right)\boldsymbol{\mathcal{\dot{Q}}}_{\boldsymbol{\mathcal{S}}}^{T}\left(\boldsymbol{x}\right),\label{eq:Eq223}
\end{align}
and also
\begin{align}
\boldsymbol{\mathcal{\dot{P}}}_{\boldsymbol{\mathcal{S}}}^{T}\left(\boldsymbol{x}\right)\boldsymbol{\mathcal{Q}}_{\boldsymbol{\mathcal{S}}}^{T}\left(\boldsymbol{x}\right)\boldsymbol{\mathcal{Q}}_{\boldsymbol{\mathcal{S}}}^{T}\left(\boldsymbol{x}\right) & =-\boldsymbol{\mathcal{P}}_{\boldsymbol{\mathcal{S}}}^{T}\left(\boldsymbol{x}\right)\boldsymbol{\mathcal{\dot{Q}}}_{\boldsymbol{\mathcal{S}}}^{T}\left(\boldsymbol{x}\right)\boldsymbol{\mathcal{Q}}_{\boldsymbol{\mathcal{S}}}^{T}\left(\boldsymbol{x}\right)\label{eq:Eq224}
\end{align}
due to idempotence of projectors. See also Section \ref{sec:PropertiesOfTimeDependentProjectors}
of the Appendix for more relations between time-dependent complementary
projectors. Using Eqs. \eqref{eq:Eq221} and \eqref{eq:Eq223} in
Eq. \eqref{eq:PSLambdaDot_3} yields 
\begin{align}
\epsilon^{2}\boldsymbol{\Gamma}_{\boldsymbol{\mathcal{S}}}\left(\boldsymbol{x}\right)\boldsymbol{\ddot{x}} & =\epsilon^{2}\boldsymbol{\mathcal{P}}_{\boldsymbol{\mathcal{S}}}^{T}\left(\boldsymbol{x}\right)\left(\boldsymbol{\Gamma}_{\boldsymbol{\mathcal{S}}}\left(\boldsymbol{x}\right)\nabla\boldsymbol{R}\left(\boldsymbol{x}\right)\boldsymbol{\dot{x}}-\boldsymbol{\dot{\Gamma}}_{\boldsymbol{\mathcal{S}}}\left(\boldsymbol{x}\right)\left(\boldsymbol{\dot{x}}-\boldsymbol{R}\left(\boldsymbol{x}\right)\right)\right)\nonumber \\
 & +\boldsymbol{\mathcal{P}}_{\boldsymbol{\mathcal{S}}}^{T}\left(\boldsymbol{x}\right)\left(\nabla\boldsymbol{R}^{T}\left(\boldsymbol{x}\right)+\boldsymbol{\mathcal{W}}^{T}\left(\boldsymbol{x},\boldsymbol{\dot{x}}\right)\right)\left(\boldsymbol{\mathcal{Q}}_{\boldsymbol{\mathcal{S}}}^{T}\left(\boldsymbol{x}\right)\boldsymbol{\lambda}-\epsilon^{2}\boldsymbol{\Gamma}_{\boldsymbol{\mathcal{S}}}\left(\boldsymbol{x}\right)\left(\boldsymbol{\dot{x}}-\boldsymbol{R}\left(\boldsymbol{x}\right)\right)\right)\nonumber \\
 & +\boldsymbol{\mathcal{P}}_{\boldsymbol{\mathcal{S}}}^{T}\left(\boldsymbol{x}\right)\boldsymbol{\mathcal{\dot{Q}}}_{\boldsymbol{\mathcal{S}}}^{T}\left(\boldsymbol{x}\right)\boldsymbol{\mathcal{Q}}_{\boldsymbol{\mathcal{S}}}^{T}\left(\boldsymbol{x}\right)\boldsymbol{\lambda}+\boldsymbol{\mathcal{P}}_{\boldsymbol{\mathcal{S}}}^{T}\left(\boldsymbol{x}\right)\boldsymbol{\mathcal{S}}\left(\boldsymbol{x}-\boldsymbol{x}_{d}\left(t\right)\right).\label{eq:PSLambdaDot_4}
\end{align}
The form of Eq. \eqref{eq:PSLambdaDot_4} makes it obvious that
it contains no component in the ''direction'' $\boldsymbol{\mathcal{Q}}_{\boldsymbol{\mathcal{S}}}^{T}\left(\boldsymbol{x}\right)$.
Equation \eqref{eq:PSLambdaDot_4} is a second order differential
equation for $n-p$ independent state components $\boldsymbol{\mathcal{P}}_{\boldsymbol{\mathcal{S}}}\left(\boldsymbol{x}\right)\boldsymbol{x}$.
The $2\left(n-p\right)$ initial or terminal conditions necessary
to solve Eq. \eqref{eq:PSLambdaDot_4} are given by Eqs. \eqref{eq:PSQSInitCond}
and \eqref{eq:PSQSTermCond}.

To summarize the derivation, the rearranged necessary optimality
conditions are 
\begin{align}
-\boldsymbol{\mathcal{Q}}_{\boldsymbol{\mathcal{S}}}^{T}\left(\boldsymbol{x}\right)\boldsymbol{\dot{\lambda}} & =\boldsymbol{\mathcal{Q}}_{\boldsymbol{\mathcal{S}}}^{T}\left(\boldsymbol{x}\right)\left(\nabla\boldsymbol{R}^{T}\left(\boldsymbol{x}\right)+\boldsymbol{\mathcal{W}}^{T}\left(\boldsymbol{x},\boldsymbol{\dot{x}}\right)\right)\left(\boldsymbol{\mathcal{Q}}_{\boldsymbol{\mathcal{S}}}^{T}\left(\boldsymbol{x}\right)\boldsymbol{\lambda}-\epsilon^{2}\boldsymbol{\Gamma}_{\boldsymbol{\mathcal{S}}}\left(\boldsymbol{x}\right)\left(\boldsymbol{\dot{x}}-\boldsymbol{R}\left(\boldsymbol{x}\right)\right)\right)\nonumber \\
 & +\boldsymbol{\mathcal{Q}}_{\boldsymbol{\mathcal{S}}}^{T}\left(\boldsymbol{x}\right)\boldsymbol{\mathcal{S}}\boldsymbol{\mathcal{Q}}_{\boldsymbol{\mathcal{S}}}\left(\boldsymbol{x}\right)\left(\boldsymbol{x}-\boldsymbol{x}_{d}\left(t\right)\right),\label{eq:Rearranged1}\\
\epsilon^{2}\boldsymbol{\Gamma}_{\boldsymbol{\mathcal{S}}}\left(\boldsymbol{x}\right)\boldsymbol{\ddot{x}} & =\epsilon^{2}\boldsymbol{\mathcal{P}}_{\boldsymbol{\mathcal{S}}}^{T}\left(\boldsymbol{x}\right)\left(\boldsymbol{\Gamma}_{\boldsymbol{\mathcal{S}}}\left(\boldsymbol{x}\right)\nabla\boldsymbol{R}\left(\boldsymbol{x}\right)\boldsymbol{\dot{x}}-\boldsymbol{\dot{\Gamma}}_{\boldsymbol{\mathcal{S}}}\left(\boldsymbol{x}\right)\left(\boldsymbol{\dot{x}}-\boldsymbol{R}\left(\boldsymbol{x}\right)\right)\right)\nonumber \\
 & +\boldsymbol{\mathcal{P}}_{\boldsymbol{\mathcal{S}}}^{T}\left(\boldsymbol{x}\right)\left(\nabla\boldsymbol{R}^{T}\left(\boldsymbol{x}\right)+\boldsymbol{\mathcal{W}}^{T}\left(\boldsymbol{x},\boldsymbol{\dot{x}}\right)\right)\left(\boldsymbol{\mathcal{Q}}_{\boldsymbol{\mathcal{S}}}^{T}\left(\boldsymbol{x}\right)\boldsymbol{\lambda}-\epsilon^{2}\boldsymbol{\Gamma}_{\boldsymbol{\mathcal{S}}}\left(\boldsymbol{x}\right)\left(\boldsymbol{\dot{x}}-\boldsymbol{R}\left(\boldsymbol{x}\right)\right)\right)\nonumber \\
 & +\boldsymbol{\mathcal{P}}_{\boldsymbol{\mathcal{S}}}^{T}\left(\boldsymbol{x}\right)\boldsymbol{\mathcal{\dot{Q}}}_{\boldsymbol{\mathcal{S}}}^{T}\left(\boldsymbol{x}\right)\boldsymbol{\mathcal{Q}}_{\boldsymbol{\mathcal{S}}}^{T}\left(\boldsymbol{x}\right)\boldsymbol{\lambda}+\boldsymbol{\mathcal{P}}_{\boldsymbol{\mathcal{S}}}^{T}\left(\boldsymbol{x}\right)\boldsymbol{\mathcal{S}}\boldsymbol{\mathcal{P}}_{\boldsymbol{\mathcal{S}}}\left(\boldsymbol{x}\right)\left(\boldsymbol{x}-\boldsymbol{x}_{d}\left(t\right)\right),\label{eq:Rearranged2}\\
\boldsymbol{\mathcal{Q}}_{\boldsymbol{\mathcal{S}}}\left(\boldsymbol{x}\right)\boldsymbol{\dot{x}} & =\boldsymbol{\mathcal{Q}}_{\boldsymbol{\mathcal{S}}}\left(\boldsymbol{x}\right)\boldsymbol{R}\left(\boldsymbol{x}\right).\label{eq:Rearranged3}
\end{align}
Equation \eqref{eq:SPCommute2} was used for the terms $\boldsymbol{\mathcal{Q}}_{\boldsymbol{\mathcal{S}}}^{T}\boldsymbol{\mathcal{S}}\boldsymbol{\mathcal{Q}}_{\boldsymbol{\mathcal{S}}}$
and $\boldsymbol{\mathcal{P}}_{\boldsymbol{\mathcal{S}}}^{T}\boldsymbol{\mathcal{S}}\boldsymbol{\mathcal{P}}_{\boldsymbol{\mathcal{S}}}$.
We emphasize that these equations are just a rearrangement of the
necessary optimality conditions Eqs. \eqref{eq:StationarityCondition}-\eqref{eq:AdjointEq},
and no approximation is involved. The small regularization parameter
$\epsilon^{2}$ multiplies the highest derivative $\boldsymbol{\ddot{x}}\left(t\right)$
in the system. This fact is exploited to perform a singular perturbation
expansion.

\subsection{Outer equations}

The outer equations are obtained by expanding the rearranged necessary
optimality conditions Eqs. \eqref{eq:Rearranged1}-\eqref{eq:Rearranged3}
in $\epsilon$. They are defined on the original time scale $t$.
For the sake of clarity, the outer solutions are distinguished from
the solutions $\boldsymbol{x}\left(t\right)$ and $\boldsymbol{\lambda}\left(t\right)$
by an index $O$,
\begin{align}
\boldsymbol{x}_{O}\left(t\right) & =\boldsymbol{x}\left(t\right), & \boldsymbol{\lambda}_{O}\left(t\right) & =\boldsymbol{\lambda}\left(t\right).
\end{align}
To shorten the notation, the time argument is suppressed in this
subsection. Setting $\epsilon=0$ in Eqs. \eqref{eq:Rearranged1}-\eqref{eq:Rearranged3}
yields a system of algebraic and first order differential equations,
\begin{align}
-\boldsymbol{\mathcal{Q}}_{\boldsymbol{\mathcal{S}}}^{T}\left(\boldsymbol{x}_{O}\right)\boldsymbol{\dot{\lambda}}_{O} & =\boldsymbol{\mathcal{Q}}_{\boldsymbol{\mathcal{S}}}^{T}\left(\boldsymbol{x}_{O}\right)\left(\nabla\boldsymbol{R}^{T}\left(\boldsymbol{x}_{O}\right)+\boldsymbol{\mathcal{W}}^{T}\left(\boldsymbol{x}_{O},\boldsymbol{\dot{x}}_{O}\right)\right)\boldsymbol{\mathcal{Q}}_{\boldsymbol{\mathcal{S}}}^{T}\left(\boldsymbol{x}_{O}\right)\boldsymbol{\lambda}_{O}\nonumber \\
 & +\boldsymbol{\mathcal{Q}}_{\boldsymbol{\mathcal{S}}}^{T}\left(\boldsymbol{x}_{O}\right)\boldsymbol{\mathcal{S}}\boldsymbol{\mathcal{Q}}_{\boldsymbol{\mathcal{S}}}\left(\boldsymbol{x}_{O}\right)\left(\boldsymbol{x}_{O}-\boldsymbol{x}_{d}\left(t\right)\right),\label{eq:OuterLeadingOrder1}\\
\boldsymbol{\mathcal{P}}_{\boldsymbol{\mathcal{S}}}^{T}\left(\boldsymbol{x}_{O}\right)\boldsymbol{\mathcal{S}}\boldsymbol{\mathcal{P}}_{\boldsymbol{\mathcal{S}}}\left(\boldsymbol{x}_{O}\right)\boldsymbol{x}_{O} & =\boldsymbol{\mathcal{P}}_{\boldsymbol{\mathcal{S}}}^{T}\left(\boldsymbol{x}_{O}\right)\boldsymbol{\mathcal{S}}\boldsymbol{\mathcal{P}}_{\boldsymbol{\mathcal{S}}}\left(\boldsymbol{x}_{O}\right)\boldsymbol{x}_{d}\left(t\right)\nonumber \\
 & -\boldsymbol{\mathcal{P}}_{\boldsymbol{\mathcal{S}}}^{T}\left(\boldsymbol{x}_{O}\right)\left(\nabla\boldsymbol{R}^{T}\left(\boldsymbol{x}_{O}\right)+\boldsymbol{\mathcal{W}}^{T}\left(\boldsymbol{x}_{O},\boldsymbol{\dot{x}}_{O}\right)\right)\boldsymbol{\mathcal{Q}}_{\boldsymbol{\mathcal{S}}}^{T}\left(\boldsymbol{x}_{O}\right)\boldsymbol{\lambda}_{O}\nonumber \\
 & -\boldsymbol{\mathcal{P}}_{\boldsymbol{\mathcal{S}}}^{T}\left(\boldsymbol{x}_{O}\right)\boldsymbol{\mathcal{\dot{Q}}}_{\boldsymbol{\mathcal{S}}}^{T}\left(\boldsymbol{x}_{O}\right)\boldsymbol{\mathcal{Q}}_{\boldsymbol{\mathcal{S}}}^{T}\left(\boldsymbol{x}_{O}\right)\boldsymbol{\lambda}_{O},\label{eq:OuterLeadingOrder2}\\
\boldsymbol{\mathcal{Q}}_{\boldsymbol{\mathcal{S}}}\left(\boldsymbol{x}_{O}\right)\boldsymbol{\dot{x}}_{O} & =\boldsymbol{\mathcal{Q}}_{\boldsymbol{\mathcal{S}}}\left(\boldsymbol{x}_{O}\right)\boldsymbol{R}\left(\boldsymbol{x}_{O}\right).\label{eq:OuterLeadingOrder3}
\end{align}
Equation \eqref{eq:OuterLeadingOrder2} is used to obtain an expression
for $\boldsymbol{\mathcal{P}}_{\boldsymbol{\mathcal{S}}}\left(\boldsymbol{x}_{O}\right)\boldsymbol{x}_{O}$.
Multiplying Eq. \eqref{eq:OuterLeadingOrder2} from the left by $\boldsymbol{\Omega}_{\boldsymbol{\mathcal{S}}}\left(\boldsymbol{x}_{O}\right)$
and exploiting Eqs. \eqref{eq:Eq4210} and \eqref{eq:Eq4211} yields
\begin{align}
\boldsymbol{\mathcal{P}}_{\boldsymbol{\mathcal{S}}}\left(\boldsymbol{x}_{O}\right)\boldsymbol{x}_{O} & =-\boldsymbol{\Omega}_{\boldsymbol{\mathcal{S}}}\left(\boldsymbol{x}_{O}\right)\left(\boldsymbol{\mathcal{\dot{Q}}}_{\boldsymbol{\mathcal{S}}}^{T}\left(\boldsymbol{x}_{O}\right)+\nabla\boldsymbol{R}^{T}\left(\boldsymbol{x}_{O}\right)+\boldsymbol{\mathcal{W}}^{T}\left(\boldsymbol{x}_{O},\boldsymbol{\dot{x}}_{O}\right)\right)\boldsymbol{\mathcal{Q}}_{\boldsymbol{\mathcal{S}}}^{T}\left(\boldsymbol{x}_{O}\right)\boldsymbol{\lambda}_{O}\nonumber \\
 & +\boldsymbol{\mathcal{P}}_{\boldsymbol{\mathcal{S}}}\left(\boldsymbol{x}_{O}\right)\boldsymbol{x}_{d}\left(t\right).\label{eq:OuterLeadingOrder21}
\end{align}
Equation \eqref{eq:OuterLeadingOrder21} is not a closed form expression
for $\boldsymbol{\mathcal{P}}_{\boldsymbol{\mathcal{S}}}\left(\boldsymbol{x}_{O}\right)\boldsymbol{x}_{O}$
as long as $\boldsymbol{\Omega}_{\boldsymbol{\mathcal{S}}}\left(\boldsymbol{x}_{O}\right)$,
$\boldsymbol{\mathcal{W}}^{T}\left(\boldsymbol{x}_{O},\boldsymbol{\dot{x}}_{O}\right)$,
and $\nabla\boldsymbol{R}^{T}\left(\boldsymbol{x}_{O}\right)$ depend
on $\boldsymbol{\mathcal{P}}_{\boldsymbol{\mathcal{S}}}\left(\boldsymbol{x}_{O}\right)\boldsymbol{x}_{O}$
as well. Nevertheless, Eq. \eqref{eq:OuterLeadingOrder21} can be
used to eliminate $\boldsymbol{\mathcal{P}}_{\boldsymbol{\mathcal{S}}}\left(\boldsymbol{x}_{O}\right)\boldsymbol{x}_{O}$
in Eqs. \eqref{eq:OuterLeadingOrder1} and \eqref{eq:OuterLeadingOrder3}.
The derivation of the explicit expression for $\boldsymbol{\mathcal{P}}_{\boldsymbol{\mathcal{S}}}\left(\boldsymbol{x}_{O}\right)\boldsymbol{x}_{O}$
in Eq. \eqref{eq:OuterLeadingOrder21} relies on the usage of the
projectors $\boldsymbol{\mathcal{P}}_{\boldsymbol{\mathcal{S}}}$
and $\boldsymbol{\mathcal{Q}}_{\boldsymbol{\mathcal{S}}}$. It is
impossible to derive an analogous relation using the projectors $\boldsymbol{\mathcal{P}}$
and $\boldsymbol{\mathcal{Q}}$ defined in Chapter \ref{chap:ExactlyRealizableTrajectories}
instead of $\boldsymbol{\mathcal{P}}_{\boldsymbol{\mathcal{S}}}$
and $\boldsymbol{\mathcal{Q}}_{\boldsymbol{\mathcal{S}}}$. In particular,
$\boldsymbol{\mathcal{P}}$ and $\boldsymbol{\mathcal{Q}}$ cannot
satisfy relations analogous to Eqs. \eqref{eq:Eq4210}, \eqref{eq:Eq4211},
and \eqref{eq:SPCommute2}. This motivates the usage of the projectors
$\boldsymbol{\mathcal{P}}_{\boldsymbol{\mathcal{S}}}$ and $\boldsymbol{\mathcal{Q}}_{\boldsymbol{\mathcal{S}}}$
retrospectively.

\subsection{\label{sub:InnerEquationsLeftSide}Inner equations - left side}

\begin{figure}[h]
\centering\includegraphics[scale=0.6]{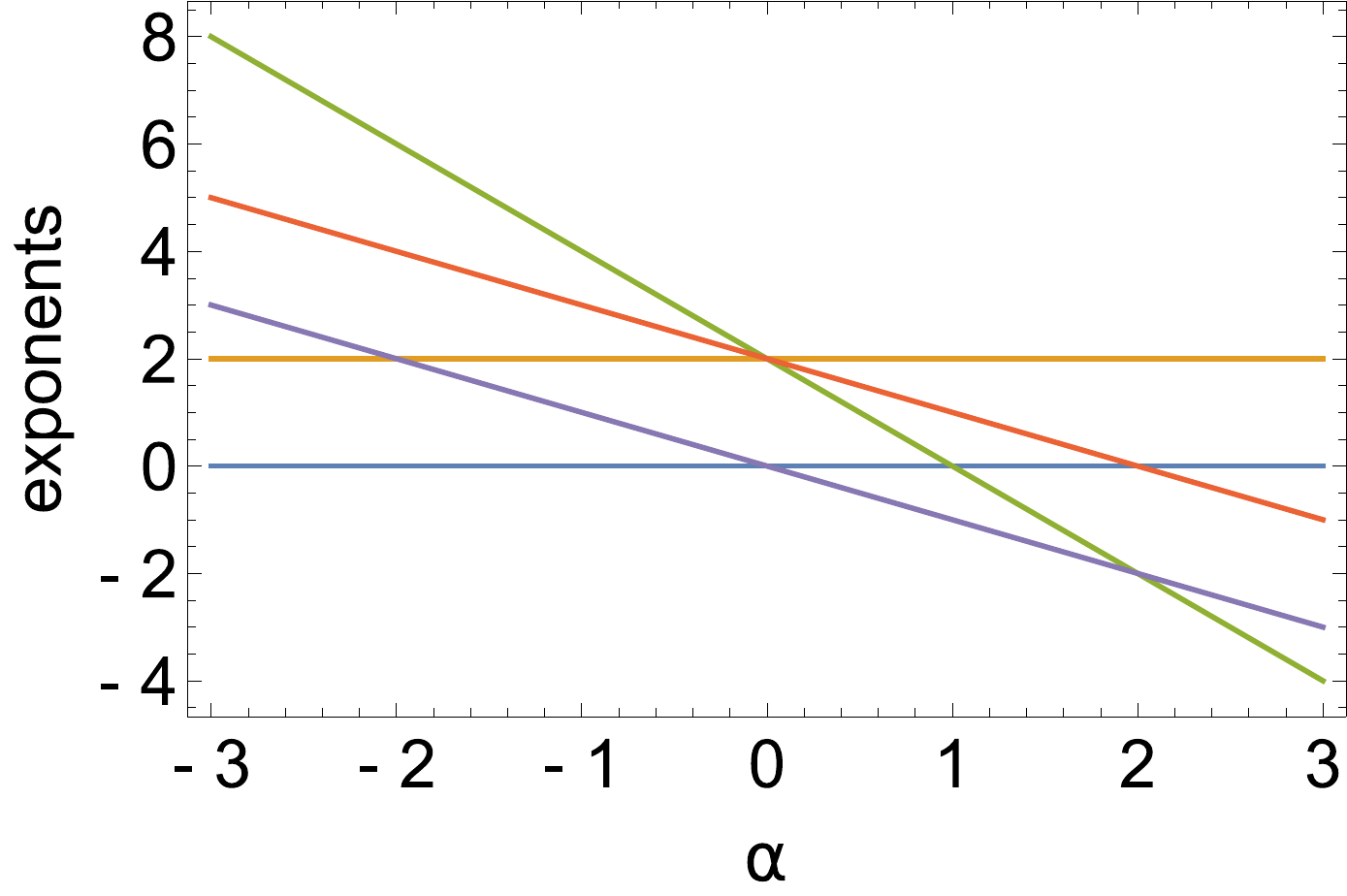}\caption[Values of exponents of $\epsilon$ over $\alpha$]{\label{fig:DominantBalance}The values of exponents of $\epsilon$
are plotted over $\alpha$. An intersection indicates a dominant balance
and determines a possible value for $\alpha$.}
\end{figure}

Boundary layers are expected at the left and right hand side of the
time domain. An initial boundary layer at the left end is resolved
by the time scale $\tau_{L}$ defined as
\begin{align}
\tau_{L} & =\left(t-t_{0}\right)/\epsilon^{\alpha}.
\end{align}
The exponent $\alpha$ has to be determined by dominant balance of
the leading order terms as $\epsilon\rightarrow0$ \cite{bender1999advanced,johnson2006singular}.
The left inner solutions are denoted by capital letters with index
$L$, 
\begin{align}
\boldsymbol{X}_{L}\left(\tau_{L}\right) & =\boldsymbol{X}_{L}\left(\left(t-t_{0}\right)/\epsilon^{\alpha}\right)=\boldsymbol{x}\left(t\right), & \boldsymbol{\Lambda}_{L}\left(\tau_{L}\right) & =\boldsymbol{\Lambda}_{L}\left(\left(t-t_{0}\right)/\epsilon^{\alpha}\right)=\boldsymbol{\lambda}\left(t\right).
\end{align}
Expressed in terms of the inner solutions, the time derivatives become
\begin{align}
\boldsymbol{\dot{x}}\left(t\right) & =\frac{d}{dt}\boldsymbol{X}_{L}\left(\left(t-t_{0}\right)/\epsilon^{\alpha}\right)=\epsilon^{-\alpha}\boldsymbol{X}_{L}'\left(\tau_{L}\right), & \boldsymbol{\dot{\lambda}}\left(t\right) & =\epsilon^{-\alpha}\boldsymbol{\Lambda}_{L}'\left(\tau_{L}\right),\\
\boldsymbol{\ddot{x}}\left(t\right) & =\frac{d^{2}}{dt^{2}}\boldsymbol{X}_{L}\left(\left(t-t_{0}\right)/\epsilon^{\alpha}\right)=\epsilon^{-2\alpha}\boldsymbol{X}_{L}''\left(\tau_{L}\right).
\end{align}
The prime $\boldsymbol{X}_{L}'\left(\tau_{L}\right)$ denotes the
derivative of $\boldsymbol{X}_{L}$ with respect to its argument $\tau_{L}$.
The time derivatives of $\boldsymbol{\mathcal{Q}}_{\boldsymbol{\mathcal{S}}}\left(\boldsymbol{x}\left(t\right)\right)$
and $\boldsymbol{\Gamma}_{\boldsymbol{\mathcal{S}}}\left(\boldsymbol{x}\left(t\right)\right)$
transform as 
\begin{align}
\boldsymbol{\mathcal{\dot{Q}}}_{\boldsymbol{\mathcal{S}}}^{T}\left(\boldsymbol{x}\left(t\right)\right) & =\epsilon^{-\alpha}\nabla\boldsymbol{\mathcal{Q}}_{\boldsymbol{\mathcal{S}}}^{T}\left(\boldsymbol{X}_{L}\left(\tau_{L}\right)\right)\boldsymbol{X}_{L}'\left(\tau_{L}\right)=\epsilon^{-\alpha}\boldsymbol{\mathcal{Q}}_{\boldsymbol{\mathcal{S}}}^{T}\vspace{0cm}'\left(\boldsymbol{X}_{L}\left(\tau_{L}\right)\right),\\
\boldsymbol{\dot{\Gamma}}_{\boldsymbol{\mathcal{S}}}\left(\boldsymbol{x}\left(t\right)\right) & =\epsilon^{-\alpha}\nabla\boldsymbol{\Gamma}_{\boldsymbol{\mathcal{S}}}\left(\boldsymbol{X}_{L}\left(\tau_{L}\right)\right)\boldsymbol{X}_{L}'\left(\tau_{L}\right)=\epsilon^{-\alpha}\boldsymbol{\Gamma}_{\boldsymbol{\mathcal{S}}}'\left(\boldsymbol{X}_{L}\left(\tau_{L}\right)\right).
\end{align}
To shorten the notation, the prime on the matrix is defined as
\begin{align}
\boldsymbol{\Gamma}_{\boldsymbol{\mathcal{S}}}'\left(\boldsymbol{X}_{L}\left(\tau_{L}\right)\right) & =\nabla\boldsymbol{\Gamma}_{\boldsymbol{\mathcal{S}}}\left(\boldsymbol{X}_{L}\left(\tau_{L}\right)\right)\boldsymbol{X}_{L}'\left(\tau_{L}\right).
\end{align}
The matrix $\boldsymbol{\mathcal{W}}\left(\boldsymbol{x}\left(t\right),\boldsymbol{\dot{x}}\left(t\right)\right)$
transforms as 
\begin{align}
\boldsymbol{\mathcal{W}}\left(\boldsymbol{x}\left(t\right),\boldsymbol{\dot{x}}\left(t\right)\right) & =\nabla\boldsymbol{\mathcal{B}}\left(\boldsymbol{x}\left(t\right)\right)\boldsymbol{\mathcal{B}}_{\boldsymbol{\mathcal{S}}}^{g}\left(\boldsymbol{x}\left(t\right)\right)\boldsymbol{\dot{x}}\left(t\right)-\nabla\boldsymbol{\mathcal{B}}\left(\boldsymbol{x}\left(t\right)\right)\boldsymbol{\mathcal{B}}_{\boldsymbol{\mathcal{S}}}^{g}\left(\boldsymbol{x}\left(t\right)\right)\boldsymbol{R}\left(\boldsymbol{x}\left(t\right)\right)\nonumber \\
 & =\epsilon^{-\alpha}\boldsymbol{\mathcal{V}}\left(\boldsymbol{X}_{L}\left(\tau_{L}\right),\boldsymbol{X}_{L}'\left(\tau_{L}\right)\right)+\boldsymbol{\mathcal{U}}\left(\boldsymbol{X}_{L}\left(\tau_{L}\right)\right),
\end{align}
with $n\times n$ matrices $\boldsymbol{\mathcal{U}}$ and $\boldsymbol{\mathcal{V}}$
defined by 
\begin{align}
\boldsymbol{\mathcal{U}}\left(\boldsymbol{x}\right) & =-\nabla\boldsymbol{\mathcal{B}}\left(\boldsymbol{x}\right)\boldsymbol{\mathcal{B}}_{\boldsymbol{\mathcal{S}}}^{g}\left(\boldsymbol{x}\right)\boldsymbol{R}\left(\boldsymbol{x}\right),\label{eq:DefU}\\
\boldsymbol{\mathcal{V}}\left(\boldsymbol{x},\boldsymbol{y}\right) & =\nabla\boldsymbol{\mathcal{B}}\left(\boldsymbol{x}\right)\boldsymbol{\mathcal{B}}_{\boldsymbol{\mathcal{S}}}^{g}\left(\boldsymbol{x}\right)\boldsymbol{y}.\label{eq:DefV}
\end{align}
The entries of $\boldsymbol{\mathcal{U}}$ and $\boldsymbol{\mathcal{V}}$
are 
\begin{align}
\mathcal{U}_{ij}\left(\boldsymbol{x}\right) & =\sum_{k=1}^{n}\sum_{l=1}^{p}\dfrac{\partial}{\partial x_{j}}\mathcal{B}_{il}\left(\boldsymbol{x}\right)\mathcal{B}_{\boldsymbol{\mathcal{S}},lk}^{g}\left(\boldsymbol{x}\right)R_{k}\left(\boldsymbol{x}\right),\\
\mathcal{V}_{ij}\left(\boldsymbol{x},\boldsymbol{y}\right) & =\sum_{k=1}^{n}\sum_{l=1}^{p}\dfrac{\partial}{\partial x_{j}}\mathcal{B}_{il}\left(\boldsymbol{x}\right)\mathcal{B}_{\boldsymbol{\mathcal{S}},lk}^{g}\left(\boldsymbol{x}\right)y_{k}.
\end{align}
From the initial conditions Eq. \eqref{eq:GenDynSysInitTermCond}
follow the initial conditions for $\boldsymbol{X}_{L}\left(\tau_{L}\right)$
as
\begin{align}
\boldsymbol{X}_{L}\left(0\right) & =\boldsymbol{x}_{0}.\label{eq:InitCondL}
\end{align}
Transforming the necessary optimality conditions Eqs. \eqref{eq:Rearranged1}-\eqref{eq:Rearranged3}
yields
\begin{align}
-\epsilon^{-\alpha}\boldsymbol{\mathcal{Q}}_{\boldsymbol{\mathcal{S}}}^{T}\left(\boldsymbol{X}_{L}\right)\boldsymbol{\Lambda}_{L}' & =-\epsilon^{2-2\alpha}\boldsymbol{\mathcal{Q}}_{\boldsymbol{\mathcal{S}}}^{T}\left(\boldsymbol{X}_{L}\right)\boldsymbol{\mathcal{V}}^{T}\left(\boldsymbol{X}_{L},\boldsymbol{X}_{L}'\right)\boldsymbol{\Gamma}_{\boldsymbol{\mathcal{S}}}\left(\boldsymbol{X}_{L}\right)\boldsymbol{X}_{L}'\nonumber \\
 & +\epsilon^{2-\alpha}\boldsymbol{\mathcal{Q}}_{\boldsymbol{\mathcal{S}}}^{T}\left(\boldsymbol{X}_{L}\right)\boldsymbol{\mathcal{V}}^{T}\left(\boldsymbol{X}_{L},\boldsymbol{X}_{L}'\right)\boldsymbol{\Gamma}_{\boldsymbol{\mathcal{S}}}\left(\boldsymbol{X}_{L}\right)\boldsymbol{R}\left(\boldsymbol{X}_{L}\right)\nonumber \\
 & -\epsilon^{2-\alpha}\boldsymbol{\mathcal{Q}}_{\boldsymbol{\mathcal{S}}}^{T}\left(\boldsymbol{X}_{L}\right)\left(\nabla\boldsymbol{R}^{T}\left(\boldsymbol{X}_{L}\right)+\boldsymbol{\mathcal{U}}^{T}\left(\boldsymbol{X}_{L}\right)\right)\boldsymbol{\Gamma}_{\boldsymbol{\mathcal{S}}}\left(\boldsymbol{X}_{L}\right)\boldsymbol{X}_{L}'\nonumber \\
 & +\epsilon^{-\alpha}\boldsymbol{\mathcal{Q}}_{\boldsymbol{\mathcal{S}}}^{T}\left(\boldsymbol{X}_{L}\right)\boldsymbol{\mathcal{V}}^{T}\left(\boldsymbol{X}_{L},\boldsymbol{X}_{L}'\right)\boldsymbol{\mathcal{Q}}_{\boldsymbol{\mathcal{S}}}^{T}\left(\boldsymbol{X}_{L}\right)\boldsymbol{\Lambda}_{L}\nonumber \\
 & +\epsilon^{2}\boldsymbol{\mathcal{Q}}_{\boldsymbol{\mathcal{S}}}^{T}\left(\boldsymbol{X}_{L}\right)\left(\nabla\boldsymbol{R}^{T}\left(\boldsymbol{X}_{L}\right)+\boldsymbol{\mathcal{U}}^{T}\left(\boldsymbol{X}_{L}\right)\right)\boldsymbol{\Gamma}_{\boldsymbol{\mathcal{S}}}\left(\boldsymbol{X}_{L}\right)\boldsymbol{R}\left(\boldsymbol{X}_{L}\right)\nonumber \\
 & +\boldsymbol{\mathcal{Q}}_{\boldsymbol{\mathcal{S}}}^{T}\left(\boldsymbol{X}_{L}\right)\left(\nabla\boldsymbol{R}^{T}\left(\boldsymbol{X}_{L}\right)+\boldsymbol{\mathcal{U}}^{T}\left(\boldsymbol{X}_{L}\right)\right)\boldsymbol{\mathcal{Q}}_{\boldsymbol{\mathcal{S}}}^{T}\left(\boldsymbol{X}_{L}\right)\boldsymbol{\Lambda}_{L}\nonumber \\
 & +\boldsymbol{\mathcal{Q}}_{\boldsymbol{\mathcal{S}}}^{T}\left(\boldsymbol{X}_{L}\right)\boldsymbol{\mathcal{S}}\boldsymbol{\mathcal{Q}}_{\boldsymbol{\mathcal{S}}}\left(\boldsymbol{X}_{L}\right)\left(\boldsymbol{X}_{L}-\boldsymbol{x}_{d}\left(t_{0}+\epsilon^{\alpha}\tau_{L}\right)\right),\label{eq:InnerLEq1General}
\end{align}
\begin{align}
\epsilon^{2-2\alpha}\boldsymbol{\Gamma}_{\boldsymbol{\mathcal{S}}}\left(\boldsymbol{X}_{L}\right)\boldsymbol{X}_{L}'' & =-\epsilon^{2-2\alpha}\boldsymbol{\mathcal{P}}_{\boldsymbol{\mathcal{S}}}^{T}\left(\boldsymbol{X}_{L}\right)\left(\boldsymbol{\mathcal{V}}^{T}\left(\boldsymbol{X}_{L},\boldsymbol{X}_{L}'\right)\boldsymbol{\Gamma}_{\boldsymbol{\mathcal{S}}}\left(\boldsymbol{X}_{L}\right)+\boldsymbol{\Gamma}_{\boldsymbol{\mathcal{S}}}'\left(\boldsymbol{X}_{L}\right)\right)\boldsymbol{X}_{L}'\nonumber \\
 & -\epsilon^{2-\alpha}\boldsymbol{\mathcal{P}}_{\boldsymbol{\mathcal{S}}}^{T}\left(\boldsymbol{X}_{L}\right)\boldsymbol{\Gamma}_{\boldsymbol{\mathcal{S}}}\left(\boldsymbol{X}_{L}\right)\nabla\boldsymbol{R}\left(\boldsymbol{X}_{L}\right)\boldsymbol{X}_{L}'\left(\tau_{L}\right)\nonumber \\
 & +\epsilon^{2-\alpha}\boldsymbol{\mathcal{P}}_{\boldsymbol{\mathcal{S}}}^{T}\left(\boldsymbol{X}_{L}\right)\boldsymbol{\Gamma}_{\boldsymbol{\mathcal{S}}}'\left(\boldsymbol{X}_{L}\right)\boldsymbol{R}\left(\boldsymbol{X}_{L}\right)\nonumber \\
 & -\epsilon^{2-\alpha}\boldsymbol{\mathcal{P}}_{\boldsymbol{\mathcal{S}}}^{T}\left(\boldsymbol{X}_{L}\right)\left(\nabla\boldsymbol{R}^{T}\left(\boldsymbol{X}_{L}\right)+\boldsymbol{\mathcal{U}}^{T}\left(\boldsymbol{X}_{L}\right)\right)\boldsymbol{\Gamma}_{\boldsymbol{\mathcal{S}}}\left(\boldsymbol{X}_{L}\right)\boldsymbol{X}_{L}'\nonumber \\
 & +\epsilon^{2-\alpha}\boldsymbol{\mathcal{P}}_{\boldsymbol{\mathcal{S}}}^{T}\left(\boldsymbol{X}_{L}\right)\boldsymbol{\mathcal{V}}^{T}\left(\boldsymbol{X}_{L},\boldsymbol{X}_{L}'\right)\boldsymbol{\Gamma}_{\boldsymbol{\mathcal{S}}}\left(\boldsymbol{X}_{L}\right)\boldsymbol{R}\left(\boldsymbol{X}_{L}\right)\nonumber \\
 & +\epsilon^{-\alpha}\boldsymbol{\mathcal{P}}_{\boldsymbol{\mathcal{S}}}^{T}\left(\boldsymbol{X}_{L}\right)\left(\boldsymbol{\mathcal{V}}^{T}\left(\boldsymbol{X}_{L},\boldsymbol{X}_{L}'\right)+\boldsymbol{\mathcal{Q}}_{\boldsymbol{\mathcal{S}}}^{T}\vspace{0cm}'\left(\boldsymbol{X}_{L}\right)\right)\boldsymbol{\mathcal{Q}}_{\boldsymbol{\mathcal{S}}}^{T}\left(\boldsymbol{X}_{L}\right)\boldsymbol{\Lambda}_{L}\nonumber \\
 & +\epsilon^{2}\boldsymbol{\mathcal{P}}_{\boldsymbol{\mathcal{S}}}^{T}\left(\boldsymbol{X}_{L}\right)\left(\nabla\boldsymbol{R}^{T}\left(\boldsymbol{X}_{L}\right)+\boldsymbol{\mathcal{U}}^{T}\left(\boldsymbol{X}_{L}\right)\right)\boldsymbol{\Gamma}_{\boldsymbol{\mathcal{S}}}\left(\boldsymbol{X}_{L}\right)\boldsymbol{R}\left(\boldsymbol{X}_{L}\right)\nonumber \\
 & +\boldsymbol{\mathcal{P}}_{\boldsymbol{\mathcal{S}}}^{T}\left(\boldsymbol{X}_{L}\right)\left(\nabla\boldsymbol{R}^{T}\left(\boldsymbol{X}_{L}\right)+\boldsymbol{\mathcal{U}}^{T}\left(\boldsymbol{X}_{L}\right)\right)\boldsymbol{\mathcal{Q}}_{\boldsymbol{\mathcal{S}}}^{T}\left(\boldsymbol{X}_{L}\right)\boldsymbol{\Lambda}_{L}\nonumber \\
 & +\boldsymbol{\mathcal{P}}_{\boldsymbol{\mathcal{S}}}^{T}\left(\boldsymbol{X}_{L}\right)\boldsymbol{\mathcal{S}}\boldsymbol{\mathcal{P}}_{\boldsymbol{\mathcal{S}}}\left(\boldsymbol{X}_{L}\right)\left(\boldsymbol{X}_{L}-\boldsymbol{x}_{d}\left(t_{0}+\epsilon^{\alpha}\tau_{L}\right)\right),\label{eq:InnerLEq2General}
\end{align}
\begin{align}
\epsilon^{-\alpha}\boldsymbol{\mathcal{Q}}_{\boldsymbol{\mathcal{S}}}\left(\boldsymbol{X}_{L}\right)\boldsymbol{X}_{L}' & =\boldsymbol{\mathcal{Q}}_{\boldsymbol{\mathcal{S}}}\left(\boldsymbol{X}_{L}\right)\boldsymbol{R}\left(\boldsymbol{X}_{L}\right).\label{eq:InnerLEq3General}
\end{align}
A dominant balance argument is applied to determine the possible
values of the exponent $\alpha$. Collecting the exponents of $\epsilon$
yields a list 
\begin{equation}
\begin{array}{ccccc}
2-2\alpha, & 2-\alpha, & -\alpha, & 2, & 0.\end{array}\label{eq:Eq4275}
\end{equation}
A dominant balance occurs if at least one pair of equal exponents
exists \cite{bender1999advanced}. Pairs of equal exponents appear
as intersections of straight lines in a plot of the exponent values
\eqref{eq:Eq4275} over $\alpha$, see Fig. \ref{fig:DominantBalance}.

The value $\alpha=0$ leads to the outer equations and is discarded.
The relevant values of $\alpha$ as given by dominant balance are
$\alpha=2$, $\alpha=1$, and $\alpha=-2$. A case by case analysis
is performed in the following.

\subsubsection{\label{sub:Case1L}Case 1 with \texorpdfstring{$\alpha = 2$}{alpha = 2}}

The leading order equations are $2\left(n-p\right)$ first order differential
equations and $p$ second order differential equations,
\begin{align}
\boldsymbol{\mathcal{Q}}_{\boldsymbol{\mathcal{S}}}^{T}\left(\boldsymbol{X}_{L}\right)\boldsymbol{\Lambda}_{L}' & =\boldsymbol{\mathcal{Q}}_{\boldsymbol{\mathcal{S}}}^{T}\left(\boldsymbol{X}_{L}\right)\boldsymbol{\mathcal{V}}^{T}\left(\boldsymbol{X}_{L},\boldsymbol{X}_{L}'\right)\left(\boldsymbol{\Gamma}_{\boldsymbol{\mathcal{S}}}\left(\boldsymbol{X}_{L}\right)\boldsymbol{X}_{L}'-\boldsymbol{\mathcal{Q}}_{\boldsymbol{\mathcal{S}}}^{T}\left(\boldsymbol{X}_{L}\right)\boldsymbol{\Lambda}_{L}\right)\label{eq:Case1LEq1}\\
\dfrac{\partial}{\partial\tau_{L}}\left(\boldsymbol{\Gamma}_{\boldsymbol{\mathcal{S}}}\left(\boldsymbol{X}_{L}\right)\boldsymbol{X}_{L}'\right) & =\boldsymbol{\mathcal{P}}_{\boldsymbol{\mathcal{S}}}^{T}\left(\boldsymbol{X}_{L}\right)\left(\boldsymbol{\mathcal{V}}^{T}\left(\boldsymbol{X}_{L},\boldsymbol{X}_{L}'\right)+\boldsymbol{\mathcal{Q}}_{\boldsymbol{\mathcal{S}}}^{T}\vspace{0cm}'\left(\boldsymbol{X}_{L}\right)\right)\boldsymbol{\mathcal{Q}}_{\boldsymbol{\mathcal{S}}}^{T}\left(\boldsymbol{X}_{L}\right)\boldsymbol{\Lambda}_{L}\nonumber \\
 & -\boldsymbol{\mathcal{P}}_{\boldsymbol{\mathcal{S}}}^{T}\left(\boldsymbol{X}_{L}\right)\boldsymbol{\mathcal{V}}^{T}\left(\boldsymbol{X}_{L},\boldsymbol{X}_{L}'\right)\boldsymbol{\Gamma}_{\boldsymbol{\mathcal{S}}}\left(\boldsymbol{X}_{L}\right)\boldsymbol{X}_{L}',\label{eq:Case1LEq2}\\
\boldsymbol{\mathcal{Q}}_{\boldsymbol{\mathcal{S}}}\left(\boldsymbol{X}_{L}\right)\boldsymbol{X}_{L}' & =\boldsymbol{0}.\label{eq:Case1LEq3}
\end{align}
Equation \eqref{eq:Case1LEq2} was simplified by using the identity
\begin{align}
\dfrac{\partial}{\partial\tau_{L}}\left(\boldsymbol{\Gamma}_{\boldsymbol{\mathcal{S}}}\left(\boldsymbol{X}_{L}\right)\boldsymbol{X}_{L}'\right) & =\boldsymbol{\Gamma}_{\boldsymbol{\mathcal{S}}}'\left(\boldsymbol{X}_{L}\right)\boldsymbol{X}_{L}'+\boldsymbol{\mathcal{P}}_{\boldsymbol{\mathcal{S}}}^{T}\left(\boldsymbol{X}_{L}\right)\boldsymbol{\Gamma}_{\boldsymbol{\mathcal{S}}}\left(\boldsymbol{X}_{L}\right)\boldsymbol{X}_{L}''\label{eq:Eq4280}
\end{align}
together with Eq. \eqref{eq:PSGammaS}.

\subsubsection{\label{sub:Case21L}Case 2.1 with \texorpdfstring{$\alpha = 1$}{alpha = 1}}

As long as 
\begin{align}
\boldsymbol{\mathcal{P}}_{\boldsymbol{\mathcal{S}}}^{T}\left(\boldsymbol{X}_{L}\right)\left(\boldsymbol{\mathcal{V}}^{T}\left(\boldsymbol{X}_{L},\boldsymbol{X}_{L}'\right)+\boldsymbol{\mathcal{Q}}_{\boldsymbol{\mathcal{S}}}^{T}\vspace{0cm}'\left(\boldsymbol{X}_{L}\right)\right)\boldsymbol{\mathcal{Q}}_{\boldsymbol{\mathcal{S}}}^{T}\left(\boldsymbol{X}_{L}\right)\boldsymbol{\Lambda}_{L} & \neq\boldsymbol{0},
\end{align}
the leading order equations are $2\left(n-p\right)+p$ first order
differential equations,
\begin{align}
\boldsymbol{\mathcal{Q}}_{\boldsymbol{\mathcal{S}}}^{T}\left(\boldsymbol{X}_{L}\right)\boldsymbol{\Lambda}_{L}' & =-\boldsymbol{\mathcal{Q}}_{\boldsymbol{\mathcal{S}}}^{T}\left(\boldsymbol{X}_{L}\right)\boldsymbol{\mathcal{V}}^{T}\left(\boldsymbol{X}_{L},\boldsymbol{X}_{L}'\right)\boldsymbol{\mathcal{Q}}_{\boldsymbol{\mathcal{S}}}^{T}\left(\boldsymbol{X}_{L}\right)\boldsymbol{\Lambda}_{L},\\
\boldsymbol{0} & =\boldsymbol{\mathcal{P}}_{\boldsymbol{\mathcal{S}}}^{T}\left(\boldsymbol{X}_{L}\right)\left(\boldsymbol{\mathcal{V}}^{T}\left(\boldsymbol{X}_{L},\boldsymbol{X}_{L}'\right)+\boldsymbol{\mathcal{Q}}_{\boldsymbol{\mathcal{S}}}^{T}\vspace{0cm}'\left(\boldsymbol{X}_{L}\right)\right)\boldsymbol{\mathcal{Q}}_{\boldsymbol{\mathcal{S}}}^{T}\left(\boldsymbol{X}_{L}\right)\boldsymbol{\Lambda}_{L},\\
\boldsymbol{\mathcal{Q}}_{\boldsymbol{\mathcal{S}}}\left(\boldsymbol{X}_{L}\right)\boldsymbol{X}_{L}' & =\boldsymbol{0}.
\end{align}

\subsubsection{\label{sub:Case22L}Case 2.2 with \texorpdfstring{$\alpha = 1$}{alpha = 1}}

The leading order equations for $\alpha=1$ change if 
\begin{align}
\boldsymbol{\mathcal{P}}_{\boldsymbol{\mathcal{S}}}^{T}\left(\boldsymbol{X}_{L}\right)\left(\boldsymbol{\mathcal{V}}^{T}\left(\boldsymbol{X}_{L},\boldsymbol{X}_{L}'\right)+\boldsymbol{\mathcal{Q}}_{\boldsymbol{\mathcal{S}}}^{T}\vspace{0cm}'\left(\boldsymbol{X}_{L}\right)\right)\boldsymbol{\mathcal{Q}}_{\boldsymbol{\mathcal{S}}}^{T}\left(\boldsymbol{X}_{L}\right)\boldsymbol{\Lambda}_{L} & =\boldsymbol{0}.
\end{align}
In this case, the leading order equations are $2\left(n-p\right)$
first order differential equations and $p$ second order differential
equations, 
\begin{align}
\boldsymbol{\mathcal{Q}}_{\boldsymbol{\mathcal{S}}}^{T}\left(\boldsymbol{X}_{L}\right)\boldsymbol{\Lambda}_{L}' & =-\boldsymbol{\mathcal{Q}}_{\boldsymbol{\mathcal{S}}}^{T}\left(\boldsymbol{X}_{L}\right)\boldsymbol{\mathcal{V}}^{T}\left(\boldsymbol{X}_{L},\boldsymbol{X}_{L}'\right)\boldsymbol{\mathcal{Q}}_{\boldsymbol{\mathcal{S}}}^{T}\left(\boldsymbol{X}_{L}\right)\boldsymbol{\Lambda}_{L},\label{eq:Case22LEq1}\\
\dfrac{\partial}{\partial\tau_{L}}\left(\boldsymbol{\Gamma}_{\boldsymbol{\mathcal{S}}}\left(\boldsymbol{X}_{L}\right)\boldsymbol{X}_{L}'\right) & =-\boldsymbol{\mathcal{P}}_{\boldsymbol{\mathcal{S}}}^{T}\left(\boldsymbol{X}_{L}\right)\boldsymbol{\mathcal{V}}^{T}\left(\boldsymbol{X}_{L},\boldsymbol{X}_{L}'\right)\boldsymbol{\Gamma}_{\boldsymbol{\mathcal{S}}}\left(\boldsymbol{X}_{L}\right)\boldsymbol{X}_{L}'\nonumber \\
 & +\boldsymbol{\mathcal{P}}_{\boldsymbol{\mathcal{S}}}^{T}\left(\boldsymbol{X}_{L}\right)\left(\nabla\boldsymbol{R}^{T}\left(\boldsymbol{X}_{L}\right)+\boldsymbol{\mathcal{U}}^{T}\left(\boldsymbol{X}_{L}\right)\right)\boldsymbol{\mathcal{Q}}_{\boldsymbol{\mathcal{S}}}^{T}\left(\boldsymbol{X}_{L}\right)\boldsymbol{\Lambda}_{L}\nonumber \\
 & +\boldsymbol{\mathcal{P}}_{\boldsymbol{\mathcal{S}}}^{T}\left(\boldsymbol{X}_{L}\right)\boldsymbol{\mathcal{S}}\boldsymbol{\mathcal{P}}_{\boldsymbol{\mathcal{S}}}\left(\boldsymbol{X}_{L}\right)\left(\boldsymbol{X}_{L}-\boldsymbol{x}_{d}\left(t_{0}\right)\right),\label{eq:Case22LEq2}\\
\boldsymbol{\mathcal{Q}}_{\boldsymbol{\mathcal{S}}}\left(\boldsymbol{X}_{L}\right)\boldsymbol{X}_{L}' & =\boldsymbol{0}.\label{eq:Case22LEq3}
\end{align}
This case corresponds to the left inner equations of the two-dimensional
dynamical system from Section \ref{sec:TwoDimensionalDynamicalSystem}.
Equations \eqref{eq:Case22LEq1}-\eqref{eq:Case22LEq3} reduce to
the left inner equations derived in Section \ref{sub:InitialBoundaryLayer}.

\subsubsection{\label{sub:Case3L}Case 3 with \texorpdfstring{$\alpha = -2$}{alpha = -2}}

Let $\boldsymbol{x}_{d,0}^{\infty}$ be defined by the limit for $\tau_{L}>0$
\begin{align}
\boldsymbol{x}_{d,0}^{\infty} & =\lim_{\epsilon\rightarrow0}\boldsymbol{x}_{d}\left(t_{0}+\epsilon^{-2}\tau_{L}\right).\label{eq:Eq4288}
\end{align}
The outer equations for $\alpha=-2$ rely on the existence of the
limit Eq. \eqref{eq:Eq4288}. The leading order equations are $2\left(n-p\right)+p$
algebraic equations, 
\begin{align}
\boldsymbol{0} & =\boldsymbol{\mathcal{Q}}_{\boldsymbol{\mathcal{S}}}^{T}\left(\boldsymbol{X}_{L}\right)\left(\nabla\boldsymbol{R}^{T}\left(\boldsymbol{X}_{L}\right)+\boldsymbol{\mathcal{U}}^{T}\left(\boldsymbol{X}_{L}\right)\right)\boldsymbol{\mathcal{Q}}_{\boldsymbol{\mathcal{S}}}^{T}\left(\boldsymbol{X}_{L}\right)\boldsymbol{\Lambda}_{L}\nonumber \\
 & +\boldsymbol{\mathcal{Q}}_{\boldsymbol{\mathcal{S}}}^{T}\left(\boldsymbol{X}_{L}\right)\boldsymbol{\mathcal{S}}\boldsymbol{\mathcal{Q}}_{\boldsymbol{\mathcal{S}}}\left(\boldsymbol{X}_{L}\right)\left(\boldsymbol{X}_{L}-\boldsymbol{x}_{d,0}^{\infty}\right),\\
\boldsymbol{0} & =\boldsymbol{\mathcal{P}}_{\boldsymbol{\mathcal{S}}}^{T}\left(\boldsymbol{X}_{L}\right)\left(\nabla\boldsymbol{R}^{T}\left(\boldsymbol{X}_{L}\right)+\boldsymbol{\mathcal{U}}^{T}\left(\boldsymbol{X}_{L}\right)\right)\boldsymbol{\mathcal{Q}}_{\boldsymbol{\mathcal{S}}}^{T}\left(\boldsymbol{X}_{L}\right)\boldsymbol{\Lambda}_{L}\nonumber \\
 & +\boldsymbol{\mathcal{P}}_{\boldsymbol{\mathcal{S}}}^{T}\left(\boldsymbol{X}_{L}\right)\boldsymbol{\mathcal{S}}\boldsymbol{\mathcal{P}}_{\boldsymbol{\mathcal{S}}}\left(\boldsymbol{X}_{L}\right)\left(\boldsymbol{X}_{L}-\boldsymbol{x}_{d,0}^{\infty}\right),\\
\boldsymbol{0} & =\boldsymbol{\mathcal{Q}}_{\boldsymbol{\mathcal{S}}}\left(\boldsymbol{X}_{L}\right)\boldsymbol{R}\left(\boldsymbol{X}_{L}\right).
\end{align}

\subsection{\label{sub:InnerEquationsRightSide}Inner equations - right side}

The boundary layer at the right hand is similarly dealt with as the
boundary layer at the left hand side. The new time scale is
\begin{align}
\tau_{R} & =\left(t_{1}-t\right)/\epsilon^{\alpha},
\end{align}
which vanishes at $t=t_{1}$. The inner solutions are denoted by capital
letters with an index $R$,
\begin{align}
\boldsymbol{X}_{R}\left(\tau_{R}\right) & =\boldsymbol{X}_{R}\left(\left(t_{1}-t\right)/\epsilon^{\alpha}\right)=\boldsymbol{x}\left(t\right), & \boldsymbol{\Lambda}_{R}\left(\tau_{R}\right) & =\boldsymbol{\Lambda}_{R}\left(\left(t_{1}-t\right)/\epsilon^{\alpha}\right)=\boldsymbol{\lambda}\left(t\right).
\end{align}
The terminal conditions Eq. \eqref{eq:GenDynSysInitTermCond} imply
initial conditions for $\boldsymbol{X}_{R}\left(\tau_{R}\right)$
as
\begin{align}
\boldsymbol{X}_{R}\left(0\right) & =\boldsymbol{x}_{1}.\label{eq:InitCondR}
\end{align}
The determination of $\alpha$ by dominant balance and the derivation
of the leading order equations proceeds analogous to the inner equations
on the left side. The only difference is that a minus sign appears
for time derivatives of odd order. Note that $\boldsymbol{\mathcal{V}}\left(\boldsymbol{x},-\boldsymbol{y}\right)=-\boldsymbol{\mathcal{V}}\left(\boldsymbol{x},\boldsymbol{y}\right)$.

\subsubsection{\label{sub:Case1R}Case 1 with \texorpdfstring{$\alpha = 2$}{alpha = 2}}

The leading order equations are $2\left(n-p\right)$ first order
differential equations and $p$ second order differential equations,
\begin{align}
\boldsymbol{\mathcal{Q}}_{\boldsymbol{\mathcal{S}}}^{T}\left(\boldsymbol{X}_{R}\right)\boldsymbol{\Lambda}_{R}' & =-\boldsymbol{\mathcal{Q}}_{\boldsymbol{\mathcal{S}}}^{T}\left(\boldsymbol{X}_{R}\right)\boldsymbol{\mathcal{V}}^{T}\left(\boldsymbol{X}_{R},\boldsymbol{X}_{R}'\right)\left(\boldsymbol{\Gamma}_{\boldsymbol{\mathcal{S}}}\left(\boldsymbol{X}_{R}\right)\boldsymbol{X}_{R}'+\boldsymbol{\mathcal{Q}}_{\boldsymbol{\mathcal{S}}}^{T}\left(\boldsymbol{X}_{R}\right)\boldsymbol{\Lambda}_{R}\right)\label{eq:Case1REq1}\\
\dfrac{\partial}{\partial\tau_{R}}\left(\boldsymbol{\Gamma}_{\boldsymbol{\mathcal{S}}}\left(\boldsymbol{X}_{R}\right)\boldsymbol{X}_{R}'\right) & =-\boldsymbol{\mathcal{P}}_{\boldsymbol{\mathcal{S}}}^{T}\left(\boldsymbol{X}_{R}\right)\left(\boldsymbol{\mathcal{V}}^{T}\left(\boldsymbol{X}_{R},\boldsymbol{X}_{R}'\right)+\boldsymbol{\mathcal{Q}}_{\boldsymbol{\mathcal{S}}}^{T}\vspace{0cm}'\left(\boldsymbol{X}_{R}\right)\right)\boldsymbol{\mathcal{Q}}_{\boldsymbol{\mathcal{S}}}^{T}\left(\boldsymbol{X}_{R}\right)\boldsymbol{\Lambda}_{R}\nonumber \\
 & -\boldsymbol{\mathcal{P}}_{\boldsymbol{\mathcal{S}}}^{T}\left(\boldsymbol{X}_{R}\right)\boldsymbol{\mathcal{V}}^{T}\left(\boldsymbol{X}_{R},\boldsymbol{X}_{R}'\right)\boldsymbol{\Gamma}_{\boldsymbol{\mathcal{S}}}\left(\boldsymbol{X}_{R}\right)\boldsymbol{X}_{R}',\\
\boldsymbol{\mathcal{Q}}_{\boldsymbol{\mathcal{S}}}\left(\boldsymbol{X}_{R}\right)\boldsymbol{X}_{R}' & =\boldsymbol{0}.\label{eq:Case1REq3}
\end{align}
Eqs. \eqref{eq:Case1REq1}-\eqref{eq:Case1REq3} are not identical
in form to their counterparts for the left boundary layer Eqs. \eqref{eq:Case1LEq1}-\eqref{eq:Case1LEq3}.

\subsubsection{\label{sub:Case21R}Case 2.1 with \texorpdfstring{$\alpha = 1$}{alpha = 1}}

As long as 
\begin{align}
\boldsymbol{\mathcal{P}}_{\boldsymbol{\mathcal{S}}}^{T}\left(\boldsymbol{X}_{R}\right)\left(\boldsymbol{\mathcal{V}}^{T}\left(\boldsymbol{X}_{R},\boldsymbol{X}_{R}'\right)+\boldsymbol{\mathcal{Q}}_{\boldsymbol{\mathcal{S}}}^{T}\vspace{0cm}'\left(\boldsymbol{X}_{R}\right)\right)\boldsymbol{\mathcal{Q}}_{\boldsymbol{\mathcal{S}}}^{T}\left(\boldsymbol{X}_{R}\right)\boldsymbol{\Lambda}_{R} & \neq\boldsymbol{0},
\end{align}
the leading order equations are $2\left(n-p\right)+p$ first order
differential equations
\begin{align}
\boldsymbol{\mathcal{Q}}_{\boldsymbol{\mathcal{S}}}^{T}\left(\boldsymbol{X}_{R}\right)\boldsymbol{\Lambda}_{R}' & =\boldsymbol{\mathcal{Q}}_{\boldsymbol{\mathcal{S}}}^{T}\left(\boldsymbol{X}_{R}\right)\boldsymbol{\mathcal{V}}^{T}\left(\boldsymbol{X}_{R},\boldsymbol{X}_{R}'\right)\boldsymbol{\mathcal{Q}}_{\boldsymbol{\mathcal{S}}}^{T}\left(\boldsymbol{X}_{R}\right)\boldsymbol{\Lambda}_{R},\\
\boldsymbol{0} & =\boldsymbol{\mathcal{P}}_{\boldsymbol{\mathcal{S}}}^{T}\left(\boldsymbol{X}_{R}\right)\left(\boldsymbol{\mathcal{V}}^{T}\left(\boldsymbol{X}_{R},\boldsymbol{X}_{R}'\right)+\boldsymbol{\mathcal{Q}}_{\boldsymbol{\mathcal{S}}}^{T}\vspace{0cm}'\left(\boldsymbol{X}_{R}\right)\right)\boldsymbol{\mathcal{Q}}_{\boldsymbol{\mathcal{S}}}^{T}\left(\boldsymbol{X}_{R}\right)\boldsymbol{\Lambda}_{R},\\
\boldsymbol{\mathcal{Q}}_{\boldsymbol{\mathcal{S}}}\left(\boldsymbol{X}_{R}\right)\boldsymbol{X}_{R}' & =\boldsymbol{0}.
\end{align}
The leading order equations are identical in form to their counterparts
for the left boundary layer.

\subsubsection{\label{sub:Case22R}Case 2.2 with \texorpdfstring{$\alpha = 1$}{alpha = 1}}

The leading order equations for $\alpha=1$ change if 
\begin{align}
\boldsymbol{\mathcal{P}}_{\boldsymbol{\mathcal{S}}}^{T}\left(\boldsymbol{X}_{R}\right)\left(\boldsymbol{\mathcal{V}}^{T}\left(\boldsymbol{X}_{R},\boldsymbol{X}_{R}'\right)+\boldsymbol{\mathcal{Q}}_{\boldsymbol{\mathcal{S}}}^{T}\vspace{0cm}'\left(\boldsymbol{X}_{R}\right)\right)\boldsymbol{\mathcal{Q}}_{\boldsymbol{\mathcal{S}}}^{T}\left(\boldsymbol{X}_{R}\right)\boldsymbol{\Lambda}_{R} & =\boldsymbol{0}.
\end{align}
In this case, the leading order equations are $2\left(n-p\right)$
first order differential equations and $p$ second order differential
equations, 
\begin{align}
\boldsymbol{\mathcal{Q}}_{\boldsymbol{\mathcal{S}}}^{T}\left(\boldsymbol{X}_{R}\right)\boldsymbol{\Lambda}_{R}' & =-\boldsymbol{\mathcal{Q}}_{\boldsymbol{\mathcal{S}}}^{T}\left(\boldsymbol{X}_{R}\right)\boldsymbol{\mathcal{V}}^{T}\left(\boldsymbol{X}_{R},\boldsymbol{X}_{R}'\right)\boldsymbol{\mathcal{Q}}_{\boldsymbol{\mathcal{S}}}^{T}\left(\boldsymbol{X}_{R}\right)\boldsymbol{\Lambda}_{R},\label{eq:Case22REq1}\\
\dfrac{\partial}{\partial\tau_{R}}\left(\boldsymbol{\Gamma}_{\boldsymbol{\mathcal{S}}}\left(\boldsymbol{X}_{R}\right)\boldsymbol{X}_{R}'\right) & =-\boldsymbol{\mathcal{P}}_{\boldsymbol{\mathcal{S}}}^{T}\left(\boldsymbol{X}_{R}\right)\boldsymbol{\mathcal{V}}^{T}\left(\boldsymbol{X}_{R},\boldsymbol{X}_{R}'\right)\boldsymbol{\Gamma}_{\boldsymbol{\mathcal{S}}}\left(\boldsymbol{X}_{R}\right)\boldsymbol{X}_{R}'\nonumber \\
 & +\boldsymbol{\mathcal{P}}_{\boldsymbol{\mathcal{S}}}^{T}\left(\boldsymbol{X}_{R}\right)\left(\nabla\boldsymbol{R}^{T}\left(\boldsymbol{X}_{R}\right)+\boldsymbol{\mathcal{U}}^{T}\left(\boldsymbol{X}_{R}\right)\right)\boldsymbol{\mathcal{Q}}_{\boldsymbol{\mathcal{S}}}^{T}\left(\boldsymbol{X}_{R}\right)\boldsymbol{\Lambda}_{R}\nonumber \\
 & +\boldsymbol{\mathcal{P}}_{\boldsymbol{\mathcal{S}}}^{T}\left(\boldsymbol{X}_{R}\right)\boldsymbol{\mathcal{S}}\boldsymbol{\mathcal{P}}_{\boldsymbol{\mathcal{S}}}\left(\boldsymbol{X}_{R}\right)\left(\boldsymbol{X}_{R}-\boldsymbol{x}_{d}\left(t_{1}\right)\right),\label{eq:Case22REq2}\\
\boldsymbol{\mathcal{Q}}_{\boldsymbol{\mathcal{S}}}\left(\boldsymbol{X}_{R}\right)\boldsymbol{X}_{R}' & =\boldsymbol{0}.\label{eq:Case22REq3}
\end{align}
The leading order equations are identical in form to their counterparts
for the left boundary layer. This case corresponds to the left inner
equations of the two-dimensional dynamical system from Section \ref{sec:TwoDimensionalDynamicalSystem}.
Equations \eqref{eq:Case22REq1}-\eqref{eq:Case22REq3} reduce to
the left inner equations derived in Section \ref{sub:TerminalBoundaryLayer}.

\subsubsection{\label{sub:Case3R}Case 3 with \texorpdfstring{$\alpha = -2$}{alpha = -2}}

Let $\boldsymbol{x}_{d,1}^{\infty}$ be defined by the limit for $\tau_{R}>0$
\begin{align}
\boldsymbol{x}_{d,1}^{\infty} & =\lim_{\epsilon\rightarrow0}\boldsymbol{x}_{d}\left(t_{1}-\epsilon^{-2}\tau_{R}\right).\label{eq:Eq4306}
\end{align}
The outer equations for $\alpha=-2$ rely on the existence of the
limit Eq. \eqref{eq:Eq4306}. The leading order equations are $2\left(n-p\right)+p$
algebraic equations, 
\begin{align}
\boldsymbol{0} & =\boldsymbol{\mathcal{Q}}_{\boldsymbol{\mathcal{S}}}^{T}\left(\boldsymbol{X}_{R}\right)\left(\nabla\boldsymbol{R}^{T}\left(\boldsymbol{X}_{R}\right)+\boldsymbol{\mathcal{U}}^{T}\left(\boldsymbol{X}_{R}\right)\right)\boldsymbol{\mathcal{Q}}_{\boldsymbol{\mathcal{S}}}^{T}\left(\boldsymbol{X}_{R}\right)\boldsymbol{\Lambda}_{R}\nonumber \\
 & +\boldsymbol{\mathcal{Q}}_{\boldsymbol{\mathcal{S}}}^{T}\left(\boldsymbol{X}_{R}\right)\boldsymbol{\mathcal{S}}\boldsymbol{\mathcal{Q}}_{\boldsymbol{\mathcal{S}}}\left(\boldsymbol{X}_{R}\right)\left(\boldsymbol{X}_{R}-\boldsymbol{x}_{d,1}^{\infty}\right),\\
\boldsymbol{0} & =\boldsymbol{\mathcal{P}}_{\boldsymbol{\mathcal{S}}}^{T}\left(\boldsymbol{X}_{R}\right)\left(\nabla\boldsymbol{R}^{T}\left(\boldsymbol{X}_{R}\right)+\boldsymbol{\mathcal{U}}^{T}\left(\boldsymbol{X}_{R}\right)\right)\boldsymbol{\mathcal{Q}}_{\boldsymbol{\mathcal{S}}}^{T}\left(\boldsymbol{X}_{R}\right)\boldsymbol{\Lambda}_{R}\nonumber \\
 & +\boldsymbol{\mathcal{P}}_{\boldsymbol{\mathcal{S}}}^{T}\left(\boldsymbol{X}_{R}\right)\boldsymbol{\mathcal{S}}\boldsymbol{\mathcal{P}}_{\boldsymbol{\mathcal{S}}}\left(\boldsymbol{X}_{R}\right)\left(\boldsymbol{X}_{R}-\boldsymbol{x}_{d,1}^{\infty}\right),\\
\boldsymbol{0} & =\boldsymbol{\mathcal{Q}}_{\boldsymbol{\mathcal{S}}}\left(\boldsymbol{X}_{R}\right)\boldsymbol{R}\left(\boldsymbol{X}_{R}\right).
\end{align}
These equations are identical in form to their counterparts for the
left boundary layers.

\subsection{Discussion of inner equations}

Several cases of inner equations are possible for general dynamical
systems. Different cases lead to different numbers of differential
and algebraic equations. Consequently, the number of boundary conditions
which can be accommodated by the inner equations differ from case
to case. In the following, we focus on cases providing the maximum
number of $2n$ boundary conditions. These are the cases with $2\left(n-p\right)$
first order and $p$ second order differential equations given by
Case 1, see Sections \ref{sub:Case1L} and \ref{sub:Case1R}, and
Case 2.2, see Sections \ref{sub:Case22L} and \ref{sub:Case22R}.
This choice is also motivated by the fact that Case 2.2 corresponds
to the left and right inner equations of the two-dimensional dynamical
system from Section \ref{sec:TwoDimensionalDynamicalSystem}. Both
cases imply a constant state projection $\boldsymbol{\mathcal{Q}}_{\boldsymbol{\mathcal{S}}}\left(\boldsymbol{X}_{L/R}\right)\boldsymbol{X}_{L/R}$,
see Eqs. \eqref{eq:Case1LEq3} and \eqref{eq:Case22LEq3}. With the
initial and terminal conditions for the state, Eqs. \eqref{eq:InitCondL}
and Eqs. \eqref{eq:InitCondR}, respectively, follows 
\begin{align}
\boldsymbol{\mathcal{Q}}_{\boldsymbol{\mathcal{S}}}\left(\boldsymbol{X}_{L}\left(\tau_{L}\right)\right)\boldsymbol{X}_{L}\left(\tau_{L}\right) & =\boldsymbol{\mathcal{Q}}_{\boldsymbol{\mathcal{S}}}\left(\boldsymbol{X}_{L}\left(\tau_{L}\right)\right)\boldsymbol{x}_{0},\label{eq:QSXLSol}\\
\boldsymbol{\mathcal{Q}}_{\boldsymbol{\mathcal{S}}}\left(\boldsymbol{X}_{R}\left(\tau_{R}\right)\right)\boldsymbol{X}_{R}\left(\tau_{R}\right) & =\boldsymbol{\mathcal{Q}}_{\boldsymbol{\mathcal{S}}}\left(\boldsymbol{X}_{R}\left(\tau_{R}\right)\right)\boldsymbol{x}_{1}.\label{eq:QSXRSol}
\end{align}
In principle, all cases of inner equations listed in Sections \ref{sub:InnerEquationsLeftSide}
and \ref{sub:InnerEquationsRightSide} can play a role for the perturbative
solution. Furthermore, more general scalings are possible for nonlinear
state and co-state equations. The scaling might not only involve a
rescaled time but also rescaled states and co-states as
\begin{align}
\boldsymbol{X}_{L}\left(\tau_{L}\right) & =\epsilon^{\beta}\boldsymbol{x}\left(t_{0}+\epsilon^{\alpha}\tau_{L}\right), & \boldsymbol{\Lambda}_{L}\left(\tau_{L}\right) & =\epsilon^{\gamma}\boldsymbol{x}\left(t_{0}+\epsilon^{\alpha}\tau_{L}\right).
\end{align}
Dominant balance arguments have to be applied to determine all possible
combinations of exponents $\alpha$, $\beta$, and $\gamma$. Usually,
the values of $\beta$ and $\gamma$ depend explicitly on the form
of all nonlinearities of the necessary optimality conditions. A larger
variety of scalings can lead to much more difficult boundary layer
structures than the simple boundary layers encountered for the two-dimensional
system of Section \ref{sec:TwoDimensionalDynamicalSystem}. Multiple
boundary layers are successions of boundary layers connecting the
initial conditions with the outer solutions by two or more scaling
regimes. Other possibilities are nested boundary layers, also called
inner-inner boundary layers, or interior boundary layers located inside
the time domain \cite{bender1999advanced}. All cases of inner equations
have to satisfy appropriate matching conditions connecting them to
their neighboring inner or outer equations. A perturbative solution
uniformly valid over the entire time interval is guaranteed only if
a combination of inner and outer solutions satisfying their appropriate
initial, terminal, and matching conditions exists. An exhaustive treatment
including proofs for the existence and uniqueness of solutions for
all possible combinations of scaling regimes is restricted to specific
control systems and not performed here.

\subsection{Matching}

Here, the matching conditions for Case 1, see Sections \ref{sub:Case1L}
and \ref{sub:Case1R}, and Case 2.2, see Sections \ref{sub:Case22L}
and \ref{sub:Case22R} are discussed. On the left side, the matching
conditions are
\begin{align}
\lim_{\tau_{L}\rightarrow\infty}\boldsymbol{\mathcal{Q}}_{\boldsymbol{\mathcal{S}}}^{T}\left(\boldsymbol{X}_{L}\left(\tau_{L}\right)\right)\boldsymbol{\Lambda}_{L}\left(\tau_{L}\right) & =\lim_{t\rightarrow t{}_{0}}\boldsymbol{\mathcal{Q}}_{\boldsymbol{\mathcal{S}}}^{T}\left(\boldsymbol{x}_{O}\left(t\right)\right)\boldsymbol{\lambda}_{O}\left(t\right),\label{eq:LeftMatching1}\\
\lim_{\tau_{L}\rightarrow\infty}\boldsymbol{\mathcal{Q}}_{\boldsymbol{\mathcal{S}}}\left(\boldsymbol{X}_{L}\left(\tau_{L}\right)\right)\boldsymbol{X}_{L}\left(\tau_{L}\right) & =\lim_{t\rightarrow t{}_{0}}\boldsymbol{\mathcal{Q}}_{\boldsymbol{\mathcal{S}}}\left(\boldsymbol{x}_{O}\left(t\right)\right)\boldsymbol{x}_{O}\left(t\right),\label{eq:LeftMatching2}\\
\lim_{\tau_{L}\rightarrow\infty}\boldsymbol{\mathcal{P}}_{\boldsymbol{\mathcal{S}}}\left(\boldsymbol{X}_{L}\left(\tau_{L}\right)\right)\boldsymbol{X}_{L}\left(\tau_{L}\right) & =\lim_{t\rightarrow t{}_{0}}\boldsymbol{\mathcal{P}}_{\boldsymbol{\mathcal{S}}}\left(\boldsymbol{x}_{O}\left(t\right)\right)\boldsymbol{x}_{O}\left(t\right).\label{eq:LeftMatching3}
\end{align}
Analogously, the matching conditions at the right side are
\begin{align}
\lim_{\tau_{R}\rightarrow\infty}\boldsymbol{\mathcal{Q}}_{\boldsymbol{\mathcal{S}}}^{T}\left(\boldsymbol{X}_{R}\left(\tau_{R}\right)\right)\boldsymbol{\Lambda}_{R}\left(\tau_{R}\right) & =\lim_{t\rightarrow t{}_{1}}\boldsymbol{\mathcal{Q}}_{\boldsymbol{\mathcal{S}}}^{T}\left(\boldsymbol{x}_{O}\left(t\right)\right)\boldsymbol{\lambda}_{O}\left(t\right),\label{eq:RightMatching1}\\
\lim_{\tau_{R}\rightarrow\infty}\boldsymbol{\mathcal{Q}}_{\boldsymbol{\mathcal{S}}}\left(\boldsymbol{X}_{R}\left(\tau_{R}\right)\right)\boldsymbol{X}_{R}\left(\tau_{R}\right) & =\lim_{t\rightarrow t{}_{1}}\boldsymbol{\mathcal{Q}}_{\boldsymbol{\mathcal{S}}}\left(\boldsymbol{x}_{O}\left(t\right)\right)\boldsymbol{x}_{O}\left(t\right),\label{eq:RightMatching2}\\
\lim_{\tau_{R}\rightarrow\infty}\boldsymbol{\mathcal{P}}_{\boldsymbol{\mathcal{S}}}\left(\boldsymbol{X}_{R}\left(\tau_{R}\right)\right)\boldsymbol{X}_{R}\left(\tau_{R}\right) & =\lim_{t\rightarrow t{}_{1}}\boldsymbol{\mathcal{P}}_{\boldsymbol{\mathcal{S}}}\left(\boldsymbol{x}_{O}\left(t\right)\right)\boldsymbol{x}_{O}\left(t\right).\label{eq:RightMatching3}
\end{align}
Adding Eq. \eqref{eq:LeftMatching2} and \eqref{eq:LeftMatching3}
yields 
\begin{align}
\lim_{\tau_{L}\rightarrow\infty}\boldsymbol{X}_{L}\left(\tau_{L}\right) & =\boldsymbol{x}_{O}\left(t_{0}\right),\label{eq:XLinf}
\end{align}
and similarly for the right matching conditions
\begin{align}
\lim_{\tau_{R}\rightarrow\infty}\boldsymbol{X}_{R}\left(\tau_{R}\right) & =\boldsymbol{x}_{O}\left(t_{1}\right).\label{eq:XRinf}
\end{align}
The last two equations yield an identity for the projector,
\begin{align}
\lim_{\tau_{L}\rightarrow\infty}\boldsymbol{\mathcal{Q}}_{\boldsymbol{\mathcal{S}}}\left(\boldsymbol{X}_{L}\left(\tau_{L}\right)\right) & =\boldsymbol{\mathcal{Q}}_{\boldsymbol{\mathcal{S}}}\left(\boldsymbol{x}_{O}\left(t_{0}\right)\right),\\
\lim_{\tau_{R}\rightarrow\infty}\boldsymbol{\mathcal{Q}}_{\boldsymbol{\mathcal{S}}}\left(\boldsymbol{X}_{R}\left(\tau_{R}\right)\right) & =\boldsymbol{\mathcal{Q}}_{\boldsymbol{\mathcal{S}}}\left(\boldsymbol{x}_{O}\left(t_{1}\right)\right),
\end{align}
and similarly for $\boldsymbol{\mathcal{P}}_{\boldsymbol{\mathcal{S}}}$.

The conditions Eqs. \eqref{eq:LeftMatching1} and \eqref{eq:RightMatching1}
yield the boundary conditions for the inner co-states $\boldsymbol{\mathcal{Q}}_{\boldsymbol{\mathcal{S}}}^{T}\left(\boldsymbol{X}_{L}\left(\tau_{L}\right)\right)\boldsymbol{\Lambda}_{L}\left(\tau_{L}\right)$
and $\boldsymbol{\mathcal{Q}}_{\boldsymbol{\mathcal{S}}}^{T}\left(\boldsymbol{X}_{R}\left(\tau_{R}\right)\right)\boldsymbol{\Lambda}_{R}\left(\tau_{R}\right)$,
respectively.

Evaluating the algebraic outer Eq. \eqref{eq:OuterLeadingOrder21}
at the initial time $t_{0}$ yields
\begin{align}
 & \boldsymbol{\mathcal{P}}_{\boldsymbol{\mathcal{S}}}\left(\boldsymbol{x}_{O}\left(t_{0}\right)\right)\boldsymbol{x}_{O}\left(t_{0}\right)=\boldsymbol{\mathcal{P}}_{\boldsymbol{\mathcal{S}}}\left(\boldsymbol{x}_{O}\left(t_{0}\right)\right)\boldsymbol{x}_{d}\left(t_{0}\right)\nonumber \\
 & -\boldsymbol{\Omega}_{\boldsymbol{\mathcal{S}}}\left(\boldsymbol{x}_{O}\left(t_{0}\right)\right)\left(\boldsymbol{\mathcal{\dot{Q}}}_{\boldsymbol{\mathcal{S}}}^{T}\left(\boldsymbol{x}_{O}\left(t_{0}\right)\right)+\nabla\boldsymbol{R}^{T}\left(\boldsymbol{x}_{O}\left(t_{0}\right)\right)\right)\boldsymbol{\mathcal{Q}}_{\boldsymbol{\mathcal{S}}}^{T}\left(\boldsymbol{x}_{O}\left(t_{0}\right)\right)\boldsymbol{\lambda}_{O}\left(t_{0}\right)\nonumber \\
 & -\boldsymbol{\Omega}_{\boldsymbol{\mathcal{S}}}\left(\boldsymbol{x}_{O}\left(t_{0}\right)\right)\boldsymbol{\mathcal{W}}^{T}\left(\boldsymbol{x}_{O}\left(t_{0}\right),\boldsymbol{\dot{x}}_{O}\left(t_{0}\right)\right)\boldsymbol{\mathcal{Q}}_{\boldsymbol{\mathcal{S}}}^{T}\left(\boldsymbol{x}_{O}\left(t_{0}\right)\right)\boldsymbol{\lambda}_{O}\left(t_{0}\right).\label{eq:Eq4323}
\end{align}
Together with the matching condition Eq. \eqref{eq:LeftMatching3},
Eq. \eqref{eq:Eq4323} results in an additional boundary condition
for $\boldsymbol{\mathcal{P}}_{\boldsymbol{\mathcal{S}}}\left(\boldsymbol{X}_{L}\right)\boldsymbol{X}_{L}$,
\begin{align}
 & \lim_{\tau_{L}\rightarrow\infty}\boldsymbol{\mathcal{P}}_{\boldsymbol{\mathcal{S}}}\left(\boldsymbol{x}_{O}\left(t_{0}\right)\right)\boldsymbol{X}_{L}\left(\tau_{L}\right)=\boldsymbol{\mathcal{P}}_{\boldsymbol{\mathcal{S}}}\left(\boldsymbol{x}_{O}\left(t_{0}\right)\right)\boldsymbol{x}_{d}\left(t_{0}\right)\nonumber \\
 & -\boldsymbol{\Omega}_{\boldsymbol{\mathcal{S}}}\left(\boldsymbol{x}_{O}\left(t_{0}\right)\right)\left(\boldsymbol{\mathcal{\dot{Q}}}_{\boldsymbol{\mathcal{S}}}^{T}\left(\boldsymbol{x}_{O}\left(t_{0}\right)\right)+\nabla\boldsymbol{R}^{T}\left(\boldsymbol{x}_{O}\left(t_{0}\right)\right)\right)\boldsymbol{\mathcal{Q}}_{\boldsymbol{\mathcal{S}}}^{T}\left(\boldsymbol{x}_{O}\left(t_{0}\right)\right)\boldsymbol{\lambda}_{O}\left(t_{0}\right)\nonumber \\
 & -\boldsymbol{\Omega}_{\boldsymbol{\mathcal{S}}}\left(\boldsymbol{x}_{O}\left(t_{0}\right)\right)\boldsymbol{\mathcal{W}}^{T}\left(\boldsymbol{x}_{O}\left(t_{0}\right),\boldsymbol{\dot{x}}_{O}\left(t_{0}\right)\right)\boldsymbol{\mathcal{Q}}_{\boldsymbol{\mathcal{S}}}^{T}\left(\boldsymbol{x}_{O}\left(t_{0}\right)\right)\boldsymbol{\lambda}_{O}\left(t_{0}\right).\label{eq:InnerLTermCond}
\end{align}
The existence of this limit, together with the result that $\boldsymbol{\mathcal{Q}}_{\boldsymbol{\mathcal{S}}}\left(\boldsymbol{X}_{L}\right)\boldsymbol{X}_{L}$
is constant, see Eq. \eqref{eq:QSXLSol}, implies 
\begin{align}
\lim_{\tau_{L}\rightarrow\infty}\boldsymbol{X}_{L}'\left(\tau_{L}\right) & =\boldsymbol{0}.\label{eq:XLTauLInfty}
\end{align}
On the other hand, evaluating the algebraic outer Eq. \eqref{eq:OuterLeadingOrder21}
at the terminal time $t_{1}$ yields
\begin{align}
 & \boldsymbol{\mathcal{P}}_{\boldsymbol{\mathcal{S}}}\left(\boldsymbol{x}_{O}\left(t_{1}\right)\right)\boldsymbol{x}_{O}\left(t_{1}\right)=\boldsymbol{\mathcal{P}}_{\boldsymbol{\mathcal{S}}}\left(\boldsymbol{x}_{O}\left(t_{1}\right)\right)\boldsymbol{x}_{d}\left(t_{1}\right)\nonumber \\
 & -\boldsymbol{\Omega}_{\boldsymbol{\mathcal{S}}}\left(\boldsymbol{x}_{O}\left(t_{1}\right)\right)\left(\boldsymbol{\mathcal{\dot{Q}}}_{\boldsymbol{\mathcal{S}}}^{T}\left(\boldsymbol{x}_{O}\left(t_{1}\right)\right)+\nabla\boldsymbol{R}^{T}\left(\boldsymbol{x}_{O}\left(t_{1}\right)\right)\right)\boldsymbol{\mathcal{Q}}_{\boldsymbol{\mathcal{S}}}^{T}\left(\boldsymbol{x}_{O}\left(t_{1}\right)\right)\boldsymbol{\lambda}_{O}\left(t_{1}\right)\nonumber \\
 & -\boldsymbol{\Omega}_{\boldsymbol{\mathcal{S}}}\left(\boldsymbol{x}_{O}\left(t_{1}\right)\right)\boldsymbol{\mathcal{W}}^{T}\left(\boldsymbol{x}_{O}\left(t_{1}\right),\boldsymbol{\dot{x}}_{O}\left(t_{1}\right)\right)\boldsymbol{\mathcal{Q}}_{\boldsymbol{\mathcal{S}}}^{T}\left(\boldsymbol{x}_{O}\left(t_{1}\right)\right)\boldsymbol{\lambda}_{O}\left(t_{1}\right).\label{eq:Eq4326}
\end{align}
Together with the matching condition Eq. \eqref{eq:RightMatching3},
Eq. \eqref{eq:Eq4326} results in an additional boundary condition
for $\boldsymbol{\mathcal{P}}_{\boldsymbol{\mathcal{S}}}\left(\boldsymbol{X}_{R}\right)\boldsymbol{X}_{R}$
in the form
\begin{align}
 & \lim_{\tau_{R}\rightarrow\infty}\boldsymbol{\mathcal{P}}_{\boldsymbol{\mathcal{S}}}\left(\boldsymbol{X}_{R}\left(\tau_{R}\right)\right)\boldsymbol{X}_{R}\left(\tau_{R}\right)=\boldsymbol{\mathcal{P}}_{\boldsymbol{\mathcal{S}}}\left(\boldsymbol{x}_{O}\left(t_{1}\right)\right)\boldsymbol{x}_{d}\left(t_{1}\right)\nonumber \\
 & -\boldsymbol{\Omega}_{\boldsymbol{\mathcal{S}}}\left(\boldsymbol{x}_{O}\left(t_{1}\right)\right)\left(\boldsymbol{\mathcal{\dot{Q}}}_{\boldsymbol{\mathcal{S}}}^{T}\left(\boldsymbol{x}_{O}\left(t_{1}\right)\right)+\nabla\boldsymbol{R}^{T}\left(\boldsymbol{x}_{O}\left(t_{1}\right)\right)\right)\boldsymbol{\mathcal{Q}}_{\boldsymbol{\mathcal{S}}}^{T}\left(\boldsymbol{x}_{O}\left(t_{1}\right)\right)\boldsymbol{\lambda}_{O}\left(t_{1}\right)\nonumber \\
 & -\boldsymbol{\Omega}_{\boldsymbol{\mathcal{S}}}\left(\boldsymbol{x}_{O}\left(t_{1}\right)\right)\boldsymbol{\mathcal{W}}^{T}\left(\boldsymbol{x}_{O}\left(t_{1}\right),\boldsymbol{\dot{x}}_{O}\left(t_{1}\right)\right)\boldsymbol{\mathcal{Q}}_{\boldsymbol{\mathcal{S}}}^{T}\left(\boldsymbol{x}_{O}\left(t_{1}\right)\right)\boldsymbol{\lambda}_{O}\left(t_{1}\right).\label{eq:InnerRTermCond}
\end{align}
Similar to above, the existence of this limit, together with the result
that $\boldsymbol{\mathcal{Q}}_{\boldsymbol{\mathcal{S}}}\left(\boldsymbol{X}_{R}\right)\boldsymbol{X}_{R}$
is constant, see Eq. \eqref{eq:QSXRSol}, implies
\begin{align}
\lim_{\tau_{R}\rightarrow\infty}\boldsymbol{X}_{R}'\left(\tau_{R}\right) & =\boldsymbol{0}.
\end{align}
Finally, two matching conditions Eqs. \eqref{eq:LeftMatching2} and
\eqref{eq:RightMatching2} remain. Because of the constancy of $\boldsymbol{\mathcal{Q}}_{\boldsymbol{\mathcal{S}}}\left(\boldsymbol{X}_{L/R}\right)\boldsymbol{X}_{L/R}$,
Eqs. \eqref{eq:QSXLSol} and \eqref{eq:QSXRSol}, together with Eqs.
\eqref{eq:XLinf} and \eqref{eq:XRinf}, respectively, these can be
written as 
\begin{align}
\boldsymbol{\mathcal{Q}}_{\boldsymbol{\mathcal{S}}}\left(\boldsymbol{x}_{O}\left(t_{0}\right)\right)\boldsymbol{x}_{0} & =\boldsymbol{\mathcal{Q}}_{\boldsymbol{\mathcal{S}}}\left(\boldsymbol{x}_{O}\left(t_{0}\right)\right)\boldsymbol{x}_{O}\left(t_{0}\right),\\
\boldsymbol{\mathcal{Q}}_{\boldsymbol{\mathcal{S}}}\left(\boldsymbol{x}_{O}\left(t_{1}\right)\right)\boldsymbol{x}_{1} & =\boldsymbol{\mathcal{Q}}_{\boldsymbol{\mathcal{S}}}\left(\boldsymbol{x}_{O}\left(t_{1}\right)\right)\boldsymbol{x}_{O}\left(t_{1}\right).
\end{align}
These are the boundary conditions for the outer equations, Eqs. \eqref{eq:OuterLeadingOrder1}
and \eqref{eq:OuterLeadingOrder3}. They depend only on the initial
and terminal conditions $\boldsymbol{x}_{0}$ and $\boldsymbol{x}_{1}$,
respectively. Hence, the solutions to the outer equations are independent
of any details of the inner equations. In particular, the outer solutions
are identical for both cases, Case 1 and Case 2.2, of inner equations
discussed here.

Finally, it is possible to formally write down the composite solutions
for the problem. The parts $\boldsymbol{\mathcal{Q}}_{\boldsymbol{\mathcal{S}}}\left(\boldsymbol{x}_{\text{comp}}\left(t\right)\right)\boldsymbol{x}_{\text{comp}}\left(t\right)$
and $\boldsymbol{\mathcal{Q}}_{\boldsymbol{\mathcal{S}}}\left(\boldsymbol{\lambda}_{\text{comp}}\left(t\right)\right)\boldsymbol{\lambda}_{\text{comp}}\left(t\right)$
do not exhibit boundary layers and are simply given by the solution
to the outer equations,
\begin{align}
\boldsymbol{\mathcal{Q}}_{\boldsymbol{\mathcal{S}}}\left(\boldsymbol{x}_{\text{comp}}\left(t\right)\right)\boldsymbol{x}_{\text{comp}}\left(t\right) & =\boldsymbol{\mathcal{Q}}_{\boldsymbol{\mathcal{S}}}\left(\boldsymbol{x}_{O}\left(t\right)\right)\boldsymbol{x}_{O}\left(t\right),\label{eq:CompSolGen1}\\
\boldsymbol{\mathcal{Q}}_{\boldsymbol{\mathcal{S}}}\left(\boldsymbol{x}_{\text{comp}}\left(t\right)\right)\boldsymbol{\lambda}_{\text{comp}}\left(t\right) & =\boldsymbol{\mathcal{Q}}_{\boldsymbol{\mathcal{S}}}\left(\boldsymbol{x}_{O}\left(t\right)\right)\boldsymbol{\lambda}_{O}\left(t\right).\label{eq:CompSolGen2}
\end{align}
The part $\boldsymbol{\mathcal{P}}_{\boldsymbol{\mathcal{S}}}\left(\boldsymbol{x}_{\text{comp}}\left(t\right)\right)\boldsymbol{x}_{\text{comp}}\left(t\right)$
contains boundary layers and is given by the sum of outer, left inner
and right inner solution minus the overlaps $\boldsymbol{\mathcal{P}}_{\boldsymbol{\mathcal{S}}}\left(\boldsymbol{x}_{O}\left(t_{0}\right)\right)\boldsymbol{x}_{O}\left(t_{0}\right)$
and $\boldsymbol{\mathcal{P}}_{\boldsymbol{\mathcal{S}}}\left(\boldsymbol{x}_{O}\left(t_{1}\right)\right)\boldsymbol{x}_{O}\left(t_{1}\right)$,
\begin{align}
 & \boldsymbol{\mathcal{P}}_{\boldsymbol{\mathcal{S}}}\left(\boldsymbol{x}_{\text{comp}}\left(t\right)\right)\boldsymbol{x}_{\text{comp}}\left(t\right)=\boldsymbol{\mathcal{P}}_{\boldsymbol{\mathcal{S}}}\left(\boldsymbol{x}_{O}\left(t\right)\right)\boldsymbol{x}_{O}\left(t\right)\nonumber \\
 & +\boldsymbol{\mathcal{P}}_{\boldsymbol{\mathcal{S}}}\left(\boldsymbol{X}_{L}\left(\epsilon^{-\alpha}\left(t-t_{0}\right)\right)\right)\boldsymbol{X}_{L}\left(\epsilon^{-\alpha}\left(t-t_{0}\right)\right)-\boldsymbol{\mathcal{P}}_{\boldsymbol{\mathcal{S}}}\left(\boldsymbol{x}_{O}\left(t_{0}\right)\right)\boldsymbol{x}_{O}\left(t_{0}\right)\nonumber \\
 & +\boldsymbol{\mathcal{P}}_{\boldsymbol{\mathcal{S}}}\left(\boldsymbol{X}_{R}\left(\epsilon^{-\alpha}\left(t_{1}-t\right)\right)\right)\boldsymbol{X}_{R}\left(\epsilon^{-\alpha}\left(t_{1}-t\right)\right)-\boldsymbol{\mathcal{P}}_{\boldsymbol{\mathcal{S}}}\left(\boldsymbol{x}_{O}\left(t_{1}\right)\right)\boldsymbol{x}_{O}\left(t_{1}\right).\label{eq:CompSolGen3}
\end{align}
Finally, the controlled state reads as
\begin{align}
\boldsymbol{x}_{\text{comp}}\left(t\right) & =\boldsymbol{\mathcal{P}}_{\boldsymbol{\mathcal{S}}}\left(\boldsymbol{x}_{\text{comp}}\left(t\right)\right)\boldsymbol{x}_{\text{comp}}\left(t\right)+\boldsymbol{\mathcal{Q}}_{\boldsymbol{\mathcal{S}}}\left(\boldsymbol{x}_{\text{comp}}\left(t\right)\right)\boldsymbol{x}_{\text{comp}}\left(t\right)\nonumber \\
 & =\boldsymbol{x}_{O}\left(t\right)-\boldsymbol{\mathcal{P}}_{\boldsymbol{\mathcal{S}}}\left(\boldsymbol{x}_{O}\left(t_{0}\right)\right)\boldsymbol{x}_{O}\left(t_{0}\right)-\boldsymbol{\mathcal{P}}_{\boldsymbol{\mathcal{S}}}\left(\boldsymbol{x}_{O}\left(t_{1}\right)\right)\boldsymbol{x}_{O}\left(t_{1}\right)\nonumber \\
 & +\boldsymbol{\mathcal{P}}_{\boldsymbol{\mathcal{S}}}\left(\boldsymbol{X}_{L}\left(\epsilon^{-\alpha}\left(t-t_{0}\right)\right)\right)\boldsymbol{X}_{L}\left(\epsilon^{-\alpha}\left(t-t_{0}\right)\right)\nonumber \\
 & +\boldsymbol{\mathcal{P}}_{\boldsymbol{\mathcal{S}}}\left(\boldsymbol{X}_{R}\left(\epsilon^{-\alpha}\left(t_{1}-t\right)\right)\right)\boldsymbol{X}_{R}\left(\epsilon^{-\alpha}\left(t_{1}-t\right)\right).
\end{align}
The composite control signal is given in terms of the composite solutions
as
\begin{align}
\boldsymbol{u}_{\text{comp}}\left(t\right) & =\boldsymbol{\mathcal{B}}_{\boldsymbol{\mathcal{S}}}^{g}\left(\boldsymbol{x}_{\text{comp}}\left(t\right)\right)\left(\boldsymbol{\dot{x}}_{\text{comp}}\left(t\right)-\boldsymbol{R}\left(\boldsymbol{x}_{\text{comp}}\left(t\right)\right)\right).
\end{align}

\subsection{\label{sub:ExactSolutionForEpsilon0}Exact state solution for \texorpdfstring{$\epsilon = 0$}{epsilon = 0}}

For a vanishing value of the regularization parameter $\epsilon$,
the inner solutions degenerate to jumps located at the time domain
boundaries. The exact solution for the controlled state trajectory
$\boldsymbol{x}\left(t\right)$ is entirely determined by the outer
equations supplemented with appropriate boundary conditions and jumps.
The time evolution of the state $\boldsymbol{x}\left(t\right)$ and
co-state $\boldsymbol{\lambda}\left(t\right)$ is governed by $2\left(n-p\right)$
first order differential equations and $2p$ algebraic equations.

First, the parts $\boldsymbol{\mathcal{P}}_{\boldsymbol{\mathcal{S}}}^{T}\left(\boldsymbol{x}\right)\boldsymbol{\lambda}$
and $\boldsymbol{\mathcal{P}}_{\boldsymbol{\mathcal{S}}}\left(\boldsymbol{x}\right)\boldsymbol{x}$
are given by algebraic equations. The part $\boldsymbol{\mathcal{P}}_{\boldsymbol{\mathcal{S}}}^{T}\left(\boldsymbol{x}\right)\boldsymbol{\lambda}$
vanishes identically for all times, 
\begin{align}
\boldsymbol{\mathcal{P}}_{\boldsymbol{\mathcal{S}}}^{T}\left(\boldsymbol{x}\left(t\right)\right)\boldsymbol{\lambda}\left(t\right) & =\boldsymbol{0}.
\end{align}
The part $\boldsymbol{\mathcal{P}}_{\boldsymbol{\mathcal{S}}}\left(\boldsymbol{x}\right)\boldsymbol{x}$
behaves discontinuously at the domain boundaries, 
\begin{align}
\boldsymbol{\mathcal{P}}_{\boldsymbol{\mathcal{S}}}\left(\boldsymbol{x}\left(t\right)\right)\boldsymbol{x}\left(t\right) & =\lim_{\epsilon\rightarrow0}\boldsymbol{\mathcal{P}}_{\boldsymbol{\mathcal{S}}}\left(\boldsymbol{x}_{\text{comp}}\left(t\right)\right)\boldsymbol{x}_{\text{comp}}\left(t\right)\nonumber \\
 & =\begin{cases}
\boldsymbol{\mathcal{P}}_{\boldsymbol{\mathcal{S}}}\left(\boldsymbol{x}_{0}\right)\boldsymbol{x}_{0}, & t=t_{0},\\
\boldsymbol{\mathcal{P}}_{\boldsymbol{\mathcal{S}}}\left(\boldsymbol{x}_{O}\left(t\right)\right)\boldsymbol{x}_{O}\left(t\right), & t_{0}<t<t_{1},\\
\boldsymbol{\mathcal{P}}_{\boldsymbol{\mathcal{S}}}\left(\boldsymbol{x}_{1}\right)\boldsymbol{x}_{1}, & t=t_{1}.
\end{cases}
\end{align}
Inside the time domain, $\boldsymbol{\mathcal{P}}_{\boldsymbol{\mathcal{S}}}\left(\boldsymbol{x}\right)\boldsymbol{x}$
behaves continuously and is given in terms of $\boldsymbol{\mathcal{Q}}_{\boldsymbol{\mathcal{S}}}^{T}\left(\boldsymbol{x}_{O}\right)\boldsymbol{\lambda}_{O}$
and $\boldsymbol{\mathcal{Q}}_{\boldsymbol{\mathcal{S}}}\left(\boldsymbol{x}_{O}\right)\boldsymbol{x}_{O}$
as 
\begin{align}
 & \boldsymbol{\mathcal{P}}_{\boldsymbol{\mathcal{S}}}\left(\boldsymbol{x}_{O}\left(t\right)\right)\boldsymbol{x}_{O}\left(t\right)=\boldsymbol{\mathcal{P}}_{\boldsymbol{\mathcal{S}}}\left(\boldsymbol{x}_{O}\left(t\right)\right)\boldsymbol{x}_{d}\left(t\right)\nonumber \\
 & -\boldsymbol{\Omega}_{\boldsymbol{\mathcal{S}}}\left(\boldsymbol{x}_{O}\left(t\right)\right)\left(\boldsymbol{\mathcal{\dot{Q}}}_{\boldsymbol{\mathcal{S}}}^{T}\left(\boldsymbol{x}_{O}\left(t\right)\right)+\nabla\boldsymbol{R}^{T}\left(\boldsymbol{x}_{O}\left(t\right)\right)\right)\boldsymbol{\mathcal{Q}}_{\boldsymbol{\mathcal{S}}}^{T}\left(\boldsymbol{x}_{O}\left(t\right)\right)\boldsymbol{\lambda}_{O}\left(t\right)\nonumber \\
 & -\boldsymbol{\Omega}_{\boldsymbol{\mathcal{S}}}\left(\boldsymbol{x}_{O}\left(t\right)\right)\boldsymbol{\mathcal{W}}^{T}\left(\boldsymbol{x}_{O}\left(t\right),\boldsymbol{\dot{x}}_{O}\left(t\right)\right)\boldsymbol{\mathcal{Q}}_{\boldsymbol{\mathcal{S}}}^{T}\left(\boldsymbol{x}_{O}\left(t\right)\right)\boldsymbol{\lambda}_{O}\left(t\right).\label{eq:ExactSol3}
\end{align}
The parts $\boldsymbol{\mathcal{Q}}_{\boldsymbol{\mathcal{S}}}^{T}\left(\boldsymbol{x}\right)\boldsymbol{\lambda}$
and $\boldsymbol{\mathcal{Q}}_{\boldsymbol{\mathcal{S}}}\left(\boldsymbol{x}\right)\boldsymbol{x}$
are given as the solution to the outer equations
\begin{alignat}{1}
\boldsymbol{\mathcal{Q}}_{\boldsymbol{\mathcal{S}}}^{T}\left(\boldsymbol{x}\left(t\right)\right)\boldsymbol{\lambda}\left(t\right) & =\boldsymbol{\mathcal{Q}}_{\boldsymbol{\mathcal{S}}}^{T}\left(\boldsymbol{x}_{O}\left(t\right)\right)\boldsymbol{\lambda}_{O}\left(t\right),\\
\boldsymbol{\mathcal{Q}}_{\boldsymbol{\mathcal{S}}}\left(\boldsymbol{x}\left(t\right)\right)\boldsymbol{x}\left(t\right) & =\boldsymbol{\mathcal{Q}}_{\boldsymbol{\mathcal{S}}}\left(\boldsymbol{x}_{O}\left(t\right)\right)\boldsymbol{x}_{O}\left(t\right),
\end{alignat}
which satisfy
\begin{align}
-\boldsymbol{\mathcal{Q}}_{\boldsymbol{\mathcal{S}}}^{T}\left(\boldsymbol{x}_{O}\left(t\right)\right)\boldsymbol{\dot{\lambda}}_{O}\left(t\right) & =\boldsymbol{\mathcal{Q}}_{\boldsymbol{\mathcal{S}}}^{T}\left(\boldsymbol{x}_{O}\left(t\right)\right)\nabla\boldsymbol{R}^{T}\left(\boldsymbol{x}_{O}\left(t\right)\right)\boldsymbol{\mathcal{Q}}_{\boldsymbol{\mathcal{S}}}^{T}\left(\boldsymbol{x}_{O}\left(t\right)\right)\boldsymbol{\lambda}_{O}\left(t\right)\nonumber \\
 & +\boldsymbol{\mathcal{Q}}_{\boldsymbol{\mathcal{S}}}^{T}\left(\boldsymbol{x}_{O}\left(t\right)\right)\boldsymbol{\mathcal{W}}^{T}\left(\boldsymbol{x}_{O}\left(t\right),\boldsymbol{\dot{x}}_{O}\left(t\right)\right)\boldsymbol{\mathcal{Q}}_{\boldsymbol{\mathcal{S}}}^{T}\left(\boldsymbol{x}_{O}\left(t\right)\right)\boldsymbol{\lambda}_{O}\left(t\right)\nonumber \\
 & +\boldsymbol{\mathcal{Q}}_{\boldsymbol{\mathcal{S}}}^{T}\left(\boldsymbol{x}_{O}\left(t\right)\right)\boldsymbol{\mathcal{S}}\boldsymbol{\mathcal{Q}}_{\boldsymbol{\mathcal{S}}}\left(\boldsymbol{x}_{O}\left(t\right)\right)\left(\boldsymbol{x}_{O}\left(t\right)-\boldsymbol{x}_{d}\left(t\right)\right),\label{eq:ExactSol1}\\
\boldsymbol{\mathcal{Q}}_{\boldsymbol{\mathcal{S}}}\left(\boldsymbol{x}_{O}\left(t\right)\right)\boldsymbol{\dot{x}}_{O}\left(t\right) & =\boldsymbol{\mathcal{Q}}_{\boldsymbol{\mathcal{S}}}\left(\boldsymbol{x}_{O}\left(t\right)\right)\boldsymbol{R}\left(\boldsymbol{x}_{O}\left(t\right)\right).\label{eq:ExactSol2}
\end{align}
These equations have to satisfy the boundary conditions
\begin{align}
\boldsymbol{\mathcal{Q}}_{\boldsymbol{\mathcal{S}}}\left(\boldsymbol{x}_{O}\left(t_{0}\right)\right)\boldsymbol{x}_{0} & =\boldsymbol{\mathcal{Q}}_{\boldsymbol{\mathcal{S}}}\left(\boldsymbol{x}_{O}\left(t_{0}\right)\right)\boldsymbol{x}_{O}\left(t_{0}\right),\label{eq:OuterInitCond0}\\
\boldsymbol{\mathcal{Q}}_{\boldsymbol{\mathcal{S}}}\left(\boldsymbol{x}_{O}\left(t_{1}\right)\right)\boldsymbol{x}_{1} & =\boldsymbol{\mathcal{Q}}_{\boldsymbol{\mathcal{S}}}\left(\boldsymbol{x}_{O}\left(t_{1}\right)\right)\boldsymbol{x}_{O}\left(t_{1}\right).\label{eq:OuterInitCond1}
\end{align}
The full state can be expressed as
\begin{align}
\boldsymbol{x}\left(t\right) & =\lim_{\epsilon\rightarrow0}\boldsymbol{x}_{\text{comp}}\left(t\right)=\begin{cases}
\boldsymbol{x}_{0}, & t=t_{0},\\
\boldsymbol{x}_{O}\left(t\right), & t_{0}<t<t_{1},\\
\boldsymbol{x}_{1}, & t=t_{1}.
\end{cases}
\end{align}
The jumps exhibited by $\boldsymbol{\mathcal{P}}_{\boldsymbol{\mathcal{S}}}\left(\boldsymbol{x}\right)\boldsymbol{x}$
at the beginning and the end of the time domain are remnants of the
boundary layers. Together with Eq. \eqref{eq:ExactSol3}, their heights
are given by
\begin{align}
 & \boldsymbol{\mathcal{P}}_{\boldsymbol{\mathcal{S}}}\left(\boldsymbol{x}_{O}\left(t_{0}\right)\right)\boldsymbol{x}_{O}\left(t_{0}\right)-\boldsymbol{\mathcal{P}}_{\boldsymbol{\mathcal{S}}}\left(\boldsymbol{x}_{0}\right)\boldsymbol{x}_{0}=\boldsymbol{\mathcal{P}}_{\boldsymbol{\mathcal{S}}}\left(\boldsymbol{x}_{O}\left(t_{0}\right)\right)\boldsymbol{x}_{d}\left(t_{0}\right)-\boldsymbol{\mathcal{P}}_{\boldsymbol{\mathcal{S}}}\left(\boldsymbol{x}_{0}\right)\boldsymbol{x}_{0}\nonumber \\
 & -\boldsymbol{\Omega}_{\boldsymbol{\mathcal{S}}}\left(\boldsymbol{x}_{O}\left(t_{0}\right)\right)\left(\boldsymbol{\mathcal{\dot{Q}}}_{\boldsymbol{\mathcal{S}}}^{T}\left(\boldsymbol{x}_{O}\left(t_{0}\right)\right)+\nabla\boldsymbol{R}^{T}\left(\boldsymbol{x}_{O}\left(t_{0}\right)\right)\right)\boldsymbol{\mathcal{Q}}_{\boldsymbol{\mathcal{S}}}^{T}\left(\boldsymbol{x}_{O}\left(t_{0}\right)\right)\boldsymbol{\lambda}_{O}\left(t_{0}\right)\nonumber \\
 & -\boldsymbol{\Omega}_{\boldsymbol{\mathcal{S}}}\left(\boldsymbol{x}_{O}\left(t_{0}\right)\right)\boldsymbol{\mathcal{W}}^{T}\left(\boldsymbol{x}_{O}\left(t_{0}\right),\boldsymbol{\dot{x}}_{O}\left(t_{0}\right)\right)\boldsymbol{\mathcal{Q}}_{\boldsymbol{\mathcal{S}}}^{T}\left(\boldsymbol{x}_{O}\left(t_{0}\right)\right)\boldsymbol{\lambda}_{O}\left(t_{0}\right),\label{eq:HeightL}
\end{align}
and
\begin{align}
 & \boldsymbol{\mathcal{P}}_{\boldsymbol{\mathcal{S}}}\left(\boldsymbol{x}_{O}\left(t_{1}\right)\right)\boldsymbol{x}_{O}\left(t_{1}\right)-\boldsymbol{\mathcal{P}}_{\boldsymbol{\mathcal{S}}}\left(\boldsymbol{x}_{1}\right)\boldsymbol{x}_{1}=\boldsymbol{\mathcal{P}}_{\boldsymbol{\mathcal{S}}}\left(\boldsymbol{x}_{O}\left(t_{1}\right)\right)\boldsymbol{x}_{d}\left(t_{1}\right)-\boldsymbol{\mathcal{P}}_{\boldsymbol{\mathcal{S}}}\left(\boldsymbol{x}_{1}\right)\boldsymbol{x}_{1}\nonumber \\
 & -\boldsymbol{\Omega}_{\boldsymbol{\mathcal{S}}}\left(\boldsymbol{x}_{O}\left(t_{1}\right)\right)\left(\boldsymbol{\mathcal{\dot{Q}}}_{\boldsymbol{\mathcal{S}}}^{T}\left(\boldsymbol{x}_{O}\left(t_{1}\right)\right)+\nabla\boldsymbol{R}^{T}\left(\boldsymbol{x}_{O}\left(t_{1}\right)\right)\right)\boldsymbol{\mathcal{Q}}_{\boldsymbol{\mathcal{S}}}^{T}\left(\boldsymbol{x}_{O}\left(t_{1}\right)\right)\boldsymbol{\lambda}_{O}\left(t_{1}\right)\nonumber \\
 & -\boldsymbol{\Omega}_{\boldsymbol{\mathcal{S}}}\left(\boldsymbol{x}_{O}\left(t_{1}\right)\right)\boldsymbol{\mathcal{W}}^{T}\left(\boldsymbol{x}_{O}\left(t_{1}\right),\boldsymbol{\dot{x}}_{O}\left(t_{1}\right)\right)\boldsymbol{\mathcal{Q}}_{\boldsymbol{\mathcal{S}}}^{T}\left(\boldsymbol{x}_{O}\left(t_{1}\right)\right)\boldsymbol{\lambda}_{O}\left(t_{1}\right),\label{eq:HeightR}
\end{align}
respectively.

In conclusion, the jump heights are entirely determined in terms of
the solutions to the outer equation together with the initial and
terminal conditions for the state. Thus, for $\epsilon=0$, no traces
of the boundary layers survive except their mere existence and location.
The inner equations play no role for the \textit{form} of the exact
state and co-state trajectory for $\epsilon=0$. However, the inner
equations play a role for the \textit{existence} of the exact solution.
The existence of jumps of appropriate height can only be guaranteed
if inner solutions satisfying appropriate initial, terminal, and matching
conditions exist.

Here, we assumed that solutions exist for the inner equations either
given by Case 1, see Sections \ref{sub:Case1L} and \ref{sub:Case1R},
or Case 2.2, see Sections \ref{sub:Case22L} and \ref{sub:Case22R}.
If other scaling regimes not given by Case 1 or Case 2.2 play a role,
the initial conditions Eqs. \eqref{eq:OuterInitCond0} and \eqref{eq:OuterInitCond1}
for the outer equations might change.

We emphasize that the exact solution stated in this section is highly
formal. The expressions for the parts $\boldsymbol{\mathcal{P}}_{\boldsymbol{\mathcal{S}}}\left(\boldsymbol{x}_{O}\right)\boldsymbol{x}_{O}$
and $\boldsymbol{\mathcal{Q}}_{\boldsymbol{\mathcal{S}}}\left(\boldsymbol{x}_{O}\right)\boldsymbol{x}_{O}$
are not closed form expressions as long as also the projectors depend
on $\boldsymbol{x}_{O}$. In general, the expression Eq. \eqref{eq:ExactSol3}
for $\boldsymbol{\mathcal{P}}_{\boldsymbol{\mathcal{S}}}\left(\boldsymbol{x}_{O}\right)\boldsymbol{x}_{O}$
is a nonlinear equation for $\boldsymbol{\mathcal{P}}_{\boldsymbol{\mathcal{S}}}\left(\boldsymbol{x}_{O}\right)\boldsymbol{x}_{O}$.
Closed form expressions can be obtained by transforming the state
$\boldsymbol{x}$ such that the projectors $\boldsymbol{\mathcal{P}}_{\boldsymbol{\mathcal{S}}}$
and $\boldsymbol{\mathcal{Q}}_{\boldsymbol{\mathcal{S}}}$ are diagonal.

\subsection{Exact control solution for \texorpdfstring{$\epsilon = 0$}{epsilon = 0}}

Formally, the control signal is given in terms of the controlled
state trajectory by the expression
\begin{align}
\boldsymbol{u}\left(t\right) & =\boldsymbol{\mathcal{B}}_{\boldsymbol{\mathcal{S}}}^{g}\left(\boldsymbol{x}\left(t\right)\right)\left(\boldsymbol{\dot{x}}\left(t\right)-\boldsymbol{R}\left(\boldsymbol{x}\left(t\right)\right)\right).
\end{align}
However, care has to be taken when evaluating the time derivative
$\boldsymbol{\dot{x}}\left(t\right)$ at the time domain boundaries.
To determine the control signal at these points, it is necessary to
analyze the expression
\begin{align}
\boldsymbol{u}_{\text{comp}}\left(t\right) & =\boldsymbol{\mathcal{B}}_{\boldsymbol{\mathcal{S}}}^{g}\left(\boldsymbol{x}_{\text{comp}}\left(t\right)\right)\left(\boldsymbol{\dot{x}}_{\text{comp}}\left(t\right)-\boldsymbol{R}\left(\boldsymbol{x}_{\text{comp}}\left(t\right)\right)\right)
\end{align}
in the limit $\epsilon\rightarrow0$. All terms except $\boldsymbol{\dot{x}}_{\text{comp}}\left(t\right)$
are well behaved. The term $\boldsymbol{\dot{x}}_{\text{comp}}\left(t\right)$
requires the investigation of the limit
\begin{align}
\lim_{\epsilon\rightarrow0}\boldsymbol{\dot{x}}_{\text{comp}}\left(t\right) & =\boldsymbol{\dot{x}}_{O}\left(t\right)\nonumber \\
 & +\lim_{\epsilon\rightarrow0}\dfrac{d}{dt}\left(\boldsymbol{\mathcal{P}}_{\boldsymbol{\mathcal{S}}}\left(\boldsymbol{X}_{L}\left(\epsilon^{-\alpha}\left(t-t_{0}\right)\right)\right)\boldsymbol{X}_{L}\left(\epsilon^{-\alpha}\left(t-t_{0}\right)\right)\right)\nonumber \\
 & +\lim_{\epsilon\rightarrow0}\dfrac{d}{dt}\left(\boldsymbol{\mathcal{P}}_{\boldsymbol{\mathcal{S}}}\left(\boldsymbol{X}_{R}\left(\epsilon^{-\alpha}\left(t_{1}-t\right)\right)\right)\boldsymbol{X}_{R}\left(\epsilon^{-\alpha}\left(t_{1}-t\right)\right)\right).
\end{align}
Similar as for two-dimensional dynamical systems in Section \ref{sec:TwoDimensionalDynamicalSystem},
it is possible to prove that $\dfrac{d}{dt}\left(\boldsymbol{\mathcal{P}}_{\boldsymbol{\mathcal{S}}}\left(\boldsymbol{X}_{L}\left(\epsilon^{-\alpha}t\right)\right)\boldsymbol{X}_{L}\left(\epsilon^{-\alpha}t\right)\right)$
yields a term proportional to the Dirac delta function in the limit
$\epsilon\rightarrow0$.

Define the $n$-dimensional vector of functions 
\begin{align}
\boldsymbol{\delta}_{L,\epsilon}\left(t\right) & =\begin{cases}
\dfrac{d}{dt}\left(\boldsymbol{\mathcal{P}}_{\boldsymbol{\mathcal{S}}}\left(\boldsymbol{X}_{L}\left(\epsilon^{-\alpha}t\right)\right)\boldsymbol{X}_{L}\left(\epsilon^{-\alpha}t\right)\right), & t\geq0,\\
\dfrac{d}{d\tilde{t}}\left(\boldsymbol{\mathcal{P}}_{\boldsymbol{\mathcal{S}}}\left(\boldsymbol{X}_{L}\left(\epsilon^{-\alpha}\tilde{t}\right)\right)\boldsymbol{X}_{L}\left(\epsilon^{-\alpha}\tilde{t}\right)\right)\Bigg|_{_{\tilde{t}=-t}}, & t<0.
\end{cases}
\end{align}
The function $\boldsymbol{\delta}_{L,\epsilon}\left(t\right)$ is
continuous for $t=0$ in every component. It can also be expressed
as
\begin{align}
\boldsymbol{\delta}_{L,\epsilon}\left(t\right) & =\epsilon^{-\alpha}\left(\nabla\boldsymbol{\mathcal{P}}_{\boldsymbol{\mathcal{S}}}\left(\boldsymbol{X}_{L}\left(\epsilon^{-\alpha}\left|t\right|\right)\right)\boldsymbol{X}_{L}'\left(\epsilon^{-\alpha}\left|t\right|\right)\right)\boldsymbol{X}_{L}\left(\epsilon^{-\alpha}\left|t\right|\right)\nonumber \\
 & +\epsilon^{-\alpha}\boldsymbol{\mathcal{P}}_{\boldsymbol{\mathcal{S}}}\left(\boldsymbol{X}_{L}\left(\epsilon^{-\alpha}\left|t\right|\right)\right)\boldsymbol{X}_{L}'\left(\epsilon^{-\alpha}\left|t\right|\right).\label{eq:DefDelta}
\end{align}
First, evaluating $\boldsymbol{\delta}_{L,\epsilon}\left(t\right)$
at $t=0$ yields 
\begin{align}
\boldsymbol{\delta}_{L,\epsilon}\left(0\right) & =\epsilon^{-\alpha}\left(\left(\nabla\boldsymbol{\mathcal{P}}_{\boldsymbol{\mathcal{S}}}\left(\boldsymbol{x}_{0}\right)\boldsymbol{X}_{L}'\left(0\right)\right)\boldsymbol{x}_{0}+\boldsymbol{\mathcal{P}}_{\boldsymbol{\mathcal{S}}}\left(\boldsymbol{x}_{0}\right)\boldsymbol{X}_{L}'\left(0\right)\right),
\end{align}
and because $\boldsymbol{X}_{L}'\left(0\right)$ is finite and does
not depend on $\epsilon$, this expression clearly diverges in the
limit $\epsilon\rightarrow0$, 
\begin{align}
\lim_{\epsilon\rightarrow0}\boldsymbol{\delta}_{L,\epsilon}\left(0\right) & =\infty\left(\left(\nabla\boldsymbol{\mathcal{P}}_{\boldsymbol{\mathcal{S}}}\left(\boldsymbol{x}_{0}\right)\boldsymbol{X}_{L}'\left(0\right)\right)\boldsymbol{x}_{0}+\boldsymbol{\mathcal{P}}_{\boldsymbol{\mathcal{S}}}\left(\boldsymbol{x}_{0}\right)\boldsymbol{X}_{L}'\left(0\right)\right).
\end{align}
Second, for $\left|t\right|>0$, $\lim_{\epsilon\rightarrow0}\boldsymbol{\delta}_{L,\epsilon}\left(t\right)$
behaves as 
\begin{align}
\lim_{\epsilon\rightarrow0}\boldsymbol{\delta}_{L,\epsilon}\left(t\right) & =\boldsymbol{0},\,t\neq0,
\end{align}
because $\boldsymbol{X}_{L}'\left(\epsilon^{-\alpha}\left|t\right|\right)$
appears in both terms of Eq. \eqref{eq:DefDelta} and behaves as (see
also Eq. \eqref{eq:XLTauLInfty})
\begin{align}
\lim_{\epsilon\rightarrow0}\boldsymbol{X}_{L}'\left(\epsilon^{-\alpha}\left|t\right|\right) & =0,\,t\neq0.
\end{align}
Third, the integral of $\boldsymbol{\delta}_{L,\epsilon}\left(t\right)$
over time $t$ must be determined. The integral can be split up in
two integrals,
\begin{align}
\intop_{-\infty}^{\infty}d\tilde{t}\boldsymbol{\delta}_{L,\epsilon}\left(\tilde{t}\right) & =\intop_{-\infty}^{0}d\tilde{t}\boldsymbol{\delta}_{L,\epsilon}\left(\tilde{t}\right)+\intop_{0}^{\infty}d\tilde{t}\boldsymbol{\delta}_{L,\epsilon}\left(\tilde{t}\right)\nonumber \\
 & =\epsilon^{-\alpha}\intop_{-\infty}^{0}d\tilde{t}\left(\nabla\boldsymbol{\mathcal{P}}_{\boldsymbol{\mathcal{S}}}\left(\boldsymbol{X}_{L}\left(-\epsilon^{-\alpha}\tilde{t}\right)\right)\boldsymbol{X}_{L}'\left(-\epsilon^{-\alpha}\tilde{t}\right)\right)\boldsymbol{X}_{L}\left(-\epsilon^{-\alpha}\tilde{t}\right)\nonumber \\
 & +\epsilon^{-\alpha}\intop_{-\infty}^{0}d\tilde{t}\boldsymbol{\mathcal{P}}_{\boldsymbol{\mathcal{S}}}\left(\boldsymbol{X}_{L}\left(-\epsilon^{-\alpha}\tilde{t}\right)\right)\boldsymbol{X}_{L}'\left(-\epsilon^{-\alpha}\tilde{t}\right)\nonumber \\
 & +\epsilon^{-\alpha}\intop_{0}^{\infty}d\tilde{t}\left(\nabla\boldsymbol{\mathcal{P}}_{\boldsymbol{\mathcal{S}}}\left(\boldsymbol{X}_{L}\left(\epsilon^{-\alpha}\tilde{t}\right)\right)\boldsymbol{X}_{L}'\left(\epsilon^{-\alpha}\tilde{t}\right)\right)\boldsymbol{X}_{L}\left(\epsilon^{-\alpha}\tilde{t}\right)\nonumber \\
 & +\epsilon^{-\alpha}\intop_{0}^{\infty}d\tilde{t}\boldsymbol{\mathcal{P}}_{\boldsymbol{\mathcal{S}}}\left(\boldsymbol{X}_{L}\left(\epsilon^{-\alpha}\tilde{t}\right)\right)\boldsymbol{X}_{L}'\left(\epsilon^{-\alpha}\tilde{t}\right).
\end{align}
Substituting $\tau=-\epsilon^{-\alpha}\tilde{t}$ in the first and
$\tau=\epsilon^{-\alpha}\tilde{t}$ in the second integral yields
\begin{align}
\intop_{-\infty}^{\infty}d\tilde{t}\boldsymbol{\delta}_{L,\epsilon}\left(\tilde{t}\right) & =2\intop_{0}^{\infty}d\tau\left(\left(\nabla\boldsymbol{\mathcal{P}}_{\boldsymbol{\mathcal{S}}}\left(\boldsymbol{X}_{L}\left(\tau\right)\right)\boldsymbol{X}_{L}'\left(\tau\right)\right)\boldsymbol{X}_{L}\left(\tau\right)+\boldsymbol{\mathcal{P}}_{\boldsymbol{\mathcal{S}}}\left(\boldsymbol{X}_{L}\left(\tau\right)\right)\boldsymbol{X}_{L}'\left(\tau\right)\right)\nonumber \\
 & =2\left(\boldsymbol{\mathcal{P}}_{\boldsymbol{\mathcal{S}}}\left(\boldsymbol{x}_{O}\left(t_{0}\right)\right)\boldsymbol{x}_{O}\left(t_{0}\right)-\boldsymbol{\mathcal{P}}_{\boldsymbol{\mathcal{S}}}\left(\boldsymbol{x}_{0}\right)\boldsymbol{x}_{0}\right).
\end{align}
Thus, we proved that
\begin{align}
\lim_{\epsilon\rightarrow0}\boldsymbol{\delta}_{L,\epsilon}\left(t\right) & =2\left(\boldsymbol{\mathcal{P}}_{\boldsymbol{\mathcal{S}}}\left(\boldsymbol{x}_{O}\left(t_{0}\right)\right)\boldsymbol{x}_{O}\left(t_{0}\right)-\boldsymbol{\mathcal{P}}_{\boldsymbol{\mathcal{S}}}\left(\boldsymbol{x}_{0}\right)\boldsymbol{x}_{0}\right)\delta\left(t\right).
\end{align}
Expressing the time derivative of $\boldsymbol{\mathcal{P}}_{\boldsymbol{\mathcal{S}}}\boldsymbol{X}_{L}$
as 
\begin{align}
\dfrac{d}{dt}\left(\boldsymbol{\mathcal{P}}_{\boldsymbol{\mathcal{S}}}\left(\boldsymbol{X}_{L}\left(\epsilon^{-\alpha}\left(t-t_{0}\right)\right)\right)\boldsymbol{X}_{L}\left(\epsilon^{-\alpha}\left(t-t_{0}\right)\right)\right) & =\boldsymbol{\delta}_{L,\epsilon}\left(t-t_{0}\right),\,t\geq t_{0},
\end{align}
finally gives 
\begin{align}
 & \lim_{\epsilon\rightarrow0}\dfrac{d}{dt}\left(\boldsymbol{\mathcal{P}}_{\boldsymbol{\mathcal{S}}}\left(\boldsymbol{X}_{L}\left(\epsilon^{-\alpha}\left(t-t_{0}\right)\right)\right)\boldsymbol{X}_{L}\left(\epsilon^{-\alpha}\left(t-t_{0}\right)\right)\right)\nonumber \\
 & =\lim_{\epsilon\rightarrow0}\boldsymbol{\delta}_{L,\epsilon}\left(t-t_{0}\right)\nonumber \\
 & =2\left(\boldsymbol{\mathcal{P}}_{\boldsymbol{\mathcal{S}}}\left(\boldsymbol{x}_{O}\left(t_{0}\right)\right)\boldsymbol{x}_{O}\left(t_{0}\right)-\boldsymbol{\mathcal{P}}_{\boldsymbol{\mathcal{S}}}\left(\boldsymbol{x}_{0}\right)\boldsymbol{x}_{0}\right)\delta\left(t-t_{0}\right).
\end{align}
A similar discussion for the right inner equation yields the equivalent
result
\begin{align}
 & \lim_{\epsilon\rightarrow0}\dfrac{d}{dt}\left(\boldsymbol{\mathcal{P}}_{\boldsymbol{\mathcal{S}}}\left(\boldsymbol{X}_{R}\left(\epsilon^{-\alpha}\left(t_{1}-t\right)\right)\right)\boldsymbol{X}_{L}\left(\epsilon^{-\alpha}\left(t_{1}-t\right)\right)\right)\nonumber \\
 & =-2\left(\boldsymbol{\mathcal{P}}_{\boldsymbol{\mathcal{S}}}\left(\boldsymbol{x}_{O}\left(t_{1}\right)\right)\boldsymbol{x}_{O}\left(t_{1}\right)-\boldsymbol{\mathcal{P}}_{\boldsymbol{\mathcal{S}}}\left(\boldsymbol{x}_{1}\right)\boldsymbol{x}_{1}\right)\delta\left(t_{1}-t\right).
\end{align}
Finally, the exact solution for the control signal for $\epsilon=0$
reads as 
\begin{align}
\boldsymbol{u}\left(t\right) & =\begin{cases}
\boldsymbol{\mathcal{B}}_{\boldsymbol{\mathcal{S}}}^{g}\left(\boldsymbol{x}_{0}\right)\left(\boldsymbol{\dot{x}}_{O}\left(t_{0}\right)-\boldsymbol{R}\left(\boldsymbol{x}_{0}\right)\right)\\
+2\boldsymbol{\mathcal{B}}_{\boldsymbol{\mathcal{S}}}^{g}\left(\boldsymbol{x}_{0}\right)\left(\boldsymbol{\mathcal{P}}_{\boldsymbol{\mathcal{S}}}\left(\boldsymbol{x}_{O}\left(t_{0}\right)\right)\boldsymbol{x}_{O}\left(t_{0}\right)-\boldsymbol{\mathcal{P}}_{\boldsymbol{\mathcal{S}}}\left(\boldsymbol{x}_{0}\right)\boldsymbol{x}_{0}\right)\delta\left(t-t_{0}\right), & t=t_{0},\\
\boldsymbol{\mathcal{B}}_{\boldsymbol{\mathcal{S}}}^{g}\left(\boldsymbol{x}_{O}\left(t\right)\right)\left(\boldsymbol{\dot{x}}_{O}\left(t\right)-\boldsymbol{R}\left(\boldsymbol{x}_{O}\left(t\right)\right)\right), & t_{0}<t<t_{1},\\
\boldsymbol{\mathcal{B}}_{\boldsymbol{\mathcal{S}}}^{g}\left(\boldsymbol{x}_{1}\right)\left(\boldsymbol{\dot{x}}_{O}\left(t_{1}\right)-\boldsymbol{R}\left(\boldsymbol{x}_{1}\right)\right)\\
-2\boldsymbol{\mathcal{B}}_{\boldsymbol{\mathcal{S}}}^{g}\left(\boldsymbol{x}_{1}\right)\left(\boldsymbol{\mathcal{P}}_{\boldsymbol{\mathcal{S}}}\left(\boldsymbol{x}_{O}\left(t_{1}\right)\right)\boldsymbol{x}_{O}\left(t_{1}\right)-\boldsymbol{\mathcal{P}}_{\boldsymbol{\mathcal{S}}}\left(\boldsymbol{x}_{1}\right)\boldsymbol{x}_{1}\right)\delta\left(t_{1}-t\right) & t=t_{1}.
\end{cases}
\end{align}
In conclusion, the control diverges at the initial and terminal time,
$t=t_{0}$ and $t=t_{1}$, respectively. The divergence is in form
of a Dirac delta function. The delta kick has a direction in state
space parallel to the jump of the discontinuous state components.
The strength of the delta kick is twice the height of the jump. Inside
the time domain, the control signal is continuous and finite and entirely
given in terms of the outer solution $\boldsymbol{x}_{O}\left(t\right)$.

\subsection{\label{sub:LinearizingAssumptionGeneral}Linearizing assumption}

The exact state solution for optimal trajectory tracking for $\epsilon=0$
is given solely in terms of the outer equations \eqref{eq:ExactSol1}
and \eqref{eq:ExactSol2}. Although these equations are simpler than
the full necessary optimality conditions, they are nevertheless nonlinear
and cannot be solved easily. However, Eqs. \eqref{eq:ExactSol1} and
\eqref{eq:ExactSol2} become linear if a linearizing assumption holds,
and a solution in closed form can be given.

First, the matrix $\boldsymbol{\Omega}_{\boldsymbol{\mathcal{S}}}\left(\boldsymbol{x}\right)$
is assumed to be constant,
\begin{align}
\boldsymbol{\Omega}_{\boldsymbol{\mathcal{S}}}\left(\boldsymbol{x}\right) & =\boldsymbol{\mathcal{B}}\left(\boldsymbol{x}\right)\left(\boldsymbol{\mathcal{B}}^{T}\left(\boldsymbol{x}\right)\boldsymbol{\mathcal{S}}\boldsymbol{\mathcal{B}}\left(\boldsymbol{x}\right)\right)^{-1}\boldsymbol{\mathcal{B}}^{T}\left(\boldsymbol{x}\right)=\text{const.}=\boldsymbol{\Omega}_{\boldsymbol{\mathcal{S}}}.\label{eq:LinearizingAssumptionS_1}
\end{align}
Note that this assumption does neither imply a constant coupling matrix
$\boldsymbol{\mathcal{B}}\left(\boldsymbol{x}\right)$ nor a constant
matrix $\boldsymbol{\Gamma}_{\boldsymbol{\mathcal{S}}}\left(\boldsymbol{x}\right)$
defined in Eq. \eqref{eq:DefGamma}. Equation \eqref{eq:LinearizingAssumptionS_1}
implies constant projectors $\boldsymbol{\mathcal{P}}_{\boldsymbol{\mathcal{S}}}\left(\boldsymbol{x}\right)$
and $\boldsymbol{\mathcal{Q}}_{\boldsymbol{\mathcal{S}}}\left(\boldsymbol{x}\right)$,
\begin{align}
\boldsymbol{\mathcal{P}}_{\boldsymbol{\mathcal{S}}}\left(\boldsymbol{x}\right) & =\boldsymbol{\Omega}_{\boldsymbol{\mathcal{S}}}\left(\boldsymbol{x}\right)\boldsymbol{\mathcal{S}}=\text{const.}=\boldsymbol{\mathcal{P}}_{\boldsymbol{\mathcal{S}}},\\
\boldsymbol{\mathcal{Q}}_{\boldsymbol{\mathcal{S}}}\left(\boldsymbol{x}\right) & =\mathbf{1}-\boldsymbol{\mathcal{P}}_{\boldsymbol{\mathcal{S}}}\left(\boldsymbol{x}\right)=\text{const.}=\boldsymbol{\mathcal{Q}}_{\boldsymbol{\mathcal{S}}},
\end{align}
and analogously constant transposed projectors $\boldsymbol{\mathcal{P}}_{\boldsymbol{\mathcal{S}}}^{T}\left(\boldsymbol{x}\right)$
and $\boldsymbol{\mathcal{Q}}_{\boldsymbol{\mathcal{S}}}^{T}\left(\boldsymbol{x}\right)$.

Second, the nonlinearity $\boldsymbol{R}\left(\boldsymbol{x}\right)$
is assumed to have the following structure with respect to the control,
\begin{align}
\boldsymbol{\mathcal{Q}}_{\boldsymbol{\mathcal{S}}}\boldsymbol{R}\left(\boldsymbol{x}\right) & =\boldsymbol{\mathcal{Q}}_{\boldsymbol{\mathcal{S}}}\boldsymbol{\mathcal{A}}\boldsymbol{x}+\boldsymbol{\mathcal{Q}}_{\boldsymbol{\mathcal{S}}}\boldsymbol{b},\label{eq:LinearizingAssumptionS_2}
\end{align}
with constant $n\times n$ matrix $\boldsymbol{\mathcal{A}}$ and
constant $n$-component vector $\boldsymbol{b}$. Cast into a single
sentence, assumption Eq. \eqref{eq:LinearizingAssumptionS_2} states
that the control signals act on the nonlinear equations of the system,
and all other equations are linear. Note that the linearizing assumption
Eqs. \eqref{eq:LinearizingAssumptionS_1} and \eqref{eq:LinearizingAssumptionS_2}
differs from the linearizing assumption of Section \ref{sec:LinearizingAssumption}
in that it involves the matrix of weights $\boldsymbol{\mathcal{S}}$.

The linearizing assumption implies that the part $\boldsymbol{\mathcal{Q}}_{\boldsymbol{\mathcal{S}}}\nabla\boldsymbol{R}$
of the Jacobi matrix is independent of the state and given by
\begin{align}
\boldsymbol{\mathcal{Q}}_{\boldsymbol{\mathcal{S}}}\nabla\boldsymbol{R}\left(\boldsymbol{x}\right) & =\boldsymbol{\mathcal{Q}}_{\boldsymbol{\mathcal{S}}}\boldsymbol{\mathcal{A}}.
\end{align}
Transposing yields
\begin{align}
\nabla\boldsymbol{R}^{T}\left(\boldsymbol{x}\right)\boldsymbol{\mathcal{Q}}_{\boldsymbol{\mathcal{S}}}^{T} & =\boldsymbol{\mathcal{A}}^{T}\boldsymbol{\mathcal{Q}}_{\boldsymbol{\mathcal{S}}}^{T}.
\end{align}
Furthermore, Eqs. \eqref{eq:LinearizingAssumptionS_1} and \eqref{eq:LinearizingAssumptionS_2}
imply 
\begin{align}
\boldsymbol{\mathcal{Q}}_{\boldsymbol{\mathcal{S}}}\boldsymbol{\mathcal{W}}\left(\boldsymbol{x},\boldsymbol{y}\right) & =\boldsymbol{0},
\end{align}
because of Eq. \eqref{eq:QTimesW}. Analogously, it follows that 
\begin{align}
\boldsymbol{\mathcal{Q}}_{\boldsymbol{\mathcal{S}}}\boldsymbol{\mathcal{U}}\left(\boldsymbol{x}\right) & =\boldsymbol{0}, & \boldsymbol{\mathcal{Q}}_{\boldsymbol{\mathcal{S}}}\boldsymbol{\mathcal{V}}\left(\boldsymbol{x},\boldsymbol{y}\right) & =\boldsymbol{0},
\end{align}
for the matrices $\boldsymbol{\mathcal{U}}$ and $\boldsymbol{\mathcal{V}}$
defined in Eqs. \eqref{eq:DefU} and \eqref{eq:DefV}, respectively.

Under the linearizing assumption, the outer equations \eqref{eq:OuterLeadingOrder1},
\eqref{eq:OuterLeadingOrder3}, and \eqref{eq:OuterLeadingOrder21}
become linear, 
\begin{align}
-\boldsymbol{\mathcal{Q}}_{\boldsymbol{\mathcal{S}}}^{T}\boldsymbol{\dot{\lambda}}_{O}\left(t\right) & =\boldsymbol{\mathcal{Q}}_{\boldsymbol{\mathcal{S}}}^{T}\boldsymbol{\mathcal{A}}^{T}\boldsymbol{\mathcal{Q}}_{\boldsymbol{\mathcal{S}}}^{T}\boldsymbol{\lambda}_{O}\left(t\right)+\boldsymbol{\mathcal{Q}}_{\boldsymbol{\mathcal{S}}}^{T}\boldsymbol{\mathcal{S}}\boldsymbol{\mathcal{Q}}_{\boldsymbol{\mathcal{S}}}\left(\boldsymbol{x}_{O}\left(t\right)-\boldsymbol{x}_{d}\left(t\right)\right),\label{eq:TransformedSystem11}\\
\boldsymbol{\mathcal{P}}_{\boldsymbol{\mathcal{S}}}^{T}\boldsymbol{\lambda}_{O}\left(t\right) & =\boldsymbol{0},\\
\boldsymbol{\mathcal{P}}_{\boldsymbol{\mathcal{S}}}\boldsymbol{x}_{O}\left(t\right) & =\boldsymbol{\mathcal{P}}_{\boldsymbol{\mathcal{S}}}\boldsymbol{x}_{d}\left(t\right)-\boldsymbol{\Omega}_{\boldsymbol{\mathcal{S}}}\boldsymbol{\mathcal{A}}^{T}\boldsymbol{\mathcal{Q}}_{\boldsymbol{\mathcal{S}}}^{T}\boldsymbol{\lambda}_{O}\left(t\right),\label{eq:TransformedSystem21}\\
\boldsymbol{\mathcal{Q}}_{\boldsymbol{\mathcal{S}}}\boldsymbol{\dot{x}}_{O}\left(t\right) & =\boldsymbol{\mathcal{Q}}_{\boldsymbol{\mathcal{S}}}\boldsymbol{\mathcal{A}}\boldsymbol{x}_{O}\left(t\right).\label{eq:TransformedSystem31}
\end{align}
The system of Eqs. \eqref{eq:TransformedSystem11}-\eqref{eq:TransformedSystem31}
is linear and can be solved. Using Eq. \eqref{eq:TransformedSystem21},
$\boldsymbol{\mathcal{P}}_{\boldsymbol{\mathcal{S}}}\boldsymbol{x}_{O}\left(t\right)$
can be eliminated from Eq. \eqref{eq:TransformedSystem31}. This
yields a system of $2\left(n-p\right)$ inhomogeneous first order
ODEs for $\boldsymbol{\mathcal{Q}}_{\boldsymbol{\mathcal{S}}}^{T}\boldsymbol{\lambda}_{O}\left(t\right)$
and $\boldsymbol{\mathcal{Q}}_{\boldsymbol{\mathcal{S}}}\boldsymbol{x}_{O}\left(t\right)$,
\begin{align}
\left(\begin{array}{c}
\boldsymbol{\mathcal{Q}}_{\boldsymbol{\mathcal{S}}}^{T}\boldsymbol{\dot{\lambda}}_{O}\left(t\right)\\
\boldsymbol{\mathcal{Q}}_{\boldsymbol{\mathcal{S}}}\boldsymbol{\dot{x}}_{O}\left(t\right)
\end{array}\right) & =\left(\begin{array}{cc}
-\boldsymbol{\mathcal{Q}}_{\boldsymbol{\mathcal{S}}}^{T}\boldsymbol{\mathcal{A}}^{T}\boldsymbol{\mathcal{Q}}_{\boldsymbol{\mathcal{S}}}^{T} & -\boldsymbol{\mathcal{Q}}_{\boldsymbol{\mathcal{S}}}^{T}\boldsymbol{\mathcal{S}}\boldsymbol{\mathcal{Q}}_{\boldsymbol{\mathcal{S}}}\\
-\boldsymbol{\mathcal{Q}}_{\boldsymbol{\mathcal{S}}}\boldsymbol{\mathcal{A}}\boldsymbol{\Omega}_{\boldsymbol{\mathcal{S}}}\boldsymbol{\mathcal{A}}^{T}\boldsymbol{\mathcal{Q}}_{\boldsymbol{\mathcal{S}}}^{T} & \boldsymbol{\mathcal{Q}}_{\boldsymbol{\mathcal{S}}}\boldsymbol{\mathcal{A}}\boldsymbol{\mathcal{Q}}_{\boldsymbol{\mathcal{S}}}
\end{array}\right)\left(\begin{array}{c}
\boldsymbol{\mathcal{Q}}_{\boldsymbol{\mathcal{S}}}^{T}\boldsymbol{\lambda}_{O}\left(t\right)\\
\boldsymbol{\mathcal{Q}}_{\boldsymbol{\mathcal{S}}}\boldsymbol{x}_{O}\left(t\right)
\end{array}\right)\nonumber \\
 & +\left(\begin{array}{c}
\boldsymbol{\mathcal{Q}}_{\boldsymbol{\mathcal{S}}}^{T}\boldsymbol{\mathcal{S}}\boldsymbol{\mathcal{Q}}_{\boldsymbol{\mathcal{S}}}\boldsymbol{x}_{d}\left(t\right)\\
\boldsymbol{\mathcal{Q}}_{\boldsymbol{\mathcal{S}}}\boldsymbol{\mathcal{A}}\boldsymbol{\mathcal{P}}_{\boldsymbol{\mathcal{S}}}\boldsymbol{x}_{d}\left(t\right)+\boldsymbol{\mathcal{Q}}_{\boldsymbol{\mathcal{S}}}\boldsymbol{b}
\end{array}\right).\label{eq:OuterEquations}
\end{align}
Equation \eqref{eq:OuterEquations} has to be solved with the initial
and terminal conditions
\begin{align}
\boldsymbol{\mathcal{Q}}_{\boldsymbol{\mathcal{S}}}\boldsymbol{x}\left(t_{0}\right) & =\boldsymbol{\mathcal{Q}}_{\boldsymbol{\mathcal{S}}}\boldsymbol{x}_{0}, & \boldsymbol{\mathcal{Q}}_{\boldsymbol{\mathcal{S}}}\boldsymbol{x}\left(t_{1}\right) & =\boldsymbol{\mathcal{Q}}_{\boldsymbol{\mathcal{S}}}\boldsymbol{x}_{1}.
\end{align}
The solution to Eq. \eqref{eq:OuterEquations} can be expressed in
closed form in terms of the state transition matrix $\boldsymbol{\Phi}\left(t,t_{0}\right)$,
\begin{align}
\left(\begin{array}{c}
\boldsymbol{\mathcal{Q}}_{\boldsymbol{\mathcal{S}}}^{T}\boldsymbol{\lambda}_{O}\left(t\right)\\
\boldsymbol{\mathcal{Q}}_{\boldsymbol{\mathcal{S}}}\boldsymbol{x}_{O}\left(t\right)
\end{array}\right) & =\boldsymbol{\Phi}\left(t,t_{0}\right)\left(\begin{array}{c}
\boldsymbol{\mathcal{Q}}_{\boldsymbol{\mathcal{S}}}^{T}\boldsymbol{\lambda}_{\text{init}}\\
\boldsymbol{\mathcal{Q}}_{\boldsymbol{\mathcal{S}}}\boldsymbol{x}_{0}
\end{array}\right)\nonumber \\
 & +\intop_{t_{0}}^{t}d\tau\boldsymbol{\Phi}\left(t,\tau\right)\left(\begin{array}{c}
\boldsymbol{\mathcal{Q}}_{\boldsymbol{\mathcal{S}}}^{T}\boldsymbol{\mathcal{S}}\boldsymbol{\mathcal{Q}}_{\boldsymbol{\mathcal{S}}}\boldsymbol{x}_{d}\left(\tau\right)\\
\boldsymbol{\mathcal{Q}}_{\boldsymbol{\mathcal{S}}}\boldsymbol{\mathcal{A}}\boldsymbol{\mathcal{P}}_{\boldsymbol{\mathcal{S}}}\boldsymbol{x}_{d}\left(\tau\right)+\boldsymbol{\mathcal{Q}}_{\boldsymbol{\mathcal{S}}}\boldsymbol{b}
\end{array}\right),
\end{align}
see Appendix \ref{sec:GeneralSolutionForForcedLinarDynamicalSystem}.
The term $\boldsymbol{\mathcal{Q}}_{\boldsymbol{\mathcal{S}}}^{T}\boldsymbol{\lambda}_{\text{init}}$
must be determined by the terminal condition $\boldsymbol{\mathcal{Q}}_{\boldsymbol{\mathcal{S}}}\boldsymbol{x}\left(t_{1}\right)=\boldsymbol{\mathcal{Q}}_{\boldsymbol{\mathcal{S}}}\boldsymbol{x}_{1}$.
Because the state matrix of Eq. \eqref{eq:OuterEquations} is constant
in time, the state transition matrix is given by the matrix exponential
and can be formally written as
\begin{align}
\boldsymbol{\Phi}\left(t,t_{0}\right) & =\exp\left(\left(\begin{array}{cc}
-\boldsymbol{\mathcal{Q}}_{\boldsymbol{\mathcal{S}}}^{T}\boldsymbol{\mathcal{A}}^{T}\boldsymbol{\mathcal{Q}}_{\boldsymbol{\mathcal{S}}}^{T} & -\boldsymbol{\mathcal{Q}}_{\boldsymbol{\mathcal{S}}}^{T}\boldsymbol{\mathcal{S}}\boldsymbol{\mathcal{Q}}_{\boldsymbol{\mathcal{S}}}\\
-\boldsymbol{\mathcal{Q}}_{\boldsymbol{\mathcal{S}}}\boldsymbol{\mathcal{A}}\boldsymbol{\Omega}_{\boldsymbol{\mathcal{S}}}\boldsymbol{\mathcal{A}}^{T}\boldsymbol{\mathcal{Q}}_{\boldsymbol{\mathcal{S}}}^{T} & \boldsymbol{\mathcal{Q}}_{\boldsymbol{\mathcal{S}}}\boldsymbol{\mathcal{A}}\boldsymbol{\mathcal{Q}}_{\boldsymbol{\mathcal{S}}}
\end{array}\right)\left(t-t_{0}\right)\right).
\end{align}
We emphasize that the linearizing assumption only leads to linear
outer equations. In general, the inner equations are nonlinear even
if the linearizing assumption holds. This is demonstrated for the
Case 2.2 of inner equations, see Sections \ref{sub:Case22L} and \ref{sub:Case22R}.
Case 2.2 corresponds to the left and right inner equations of the
two-dimensional dynamical system from Section \ref{sec:TwoDimensionalDynamicalSystem}.
The left inner equations \eqref{eq:Case22LEq1} and \eqref{eq:Case22LEq3}
become 
\begin{align}
\boldsymbol{\mathcal{Q}}_{\boldsymbol{\mathcal{S}}}^{T}\boldsymbol{\Lambda}_{L}'\left(\tau_{L}\right) & =\boldsymbol{0}, & \boldsymbol{\mathcal{Q}}_{\boldsymbol{\mathcal{S}}}\boldsymbol{X}_{L}'\left(\tau_{L}\right) & =\boldsymbol{0}.
\end{align}
The initial conditions Eqs. \eqref{eq:QSXLSol} and \eqref{eq:LeftMatching1}
lead to 
\begin{align}
\boldsymbol{\mathcal{Q}}_{\boldsymbol{\mathcal{S}}}^{T}\boldsymbol{\Lambda}_{L}\left(\tau_{L}\right) & =\boldsymbol{\mathcal{Q}}_{\boldsymbol{\mathcal{S}}}^{T}\boldsymbol{\lambda}_{O}\left(t_{0}\right), & \boldsymbol{\mathcal{Q}}_{\boldsymbol{\mathcal{S}}}\boldsymbol{X}_{L}\left(\tau_{L}\right) & =\boldsymbol{\mathcal{Q}}_{\boldsymbol{\mathcal{S}}}\boldsymbol{x}_{0},
\end{align}
and the remaining left inner equation \eqref{eq:Case22LEq2} becomes
\begin{align}
\dfrac{\partial}{\partial\tau_{L}}\left(\boldsymbol{\Gamma}_{\boldsymbol{\mathcal{S}}}\left(\boldsymbol{X}_{L}\left(\tau_{L}\right)\right)\boldsymbol{X}_{L}'\left(\tau_{L}\right)\right) & =-\boldsymbol{\mathcal{P}}_{\boldsymbol{\mathcal{S}}}^{T}\boldsymbol{\mathcal{V}}^{T}\left(\boldsymbol{X}_{L}\left(\tau_{L}\right),\boldsymbol{X}_{L}'\left(\tau_{L}\right)\right)\boldsymbol{\Gamma}_{\boldsymbol{\mathcal{S}}}\left(\boldsymbol{X}_{L}\left(\tau_{L}\right)\right)\boldsymbol{X}_{L}'\left(\tau_{L}\right)\nonumber \\
 & +\boldsymbol{\mathcal{P}}_{\boldsymbol{\mathcal{S}}}^{T}\boldsymbol{\mathcal{A}}^{T}\boldsymbol{\mathcal{Q}}_{\boldsymbol{\mathcal{S}}}^{T}\boldsymbol{\lambda}_{O}\left(t_{0}\right)+\boldsymbol{\mathcal{P}}_{\boldsymbol{\mathcal{S}}}^{T}\boldsymbol{\mathcal{S}}\boldsymbol{\mathcal{P}}_{\boldsymbol{\mathcal{S}}}\left(\boldsymbol{X}_{L}\left(\tau_{L}\right)-\boldsymbol{x}_{d}\left(t_{0}\right)\right).\label{eq:EqCase22LEq2LinAss}
\end{align}
In general, the matrices $\boldsymbol{\Gamma}_{\boldsymbol{\mathcal{S}}}$
and $\boldsymbol{\mathcal{V}}$ depend nonlinearly on the state. In
both cases, the nonlinearity originates from the coupling matrix $\boldsymbol{\mathcal{B}}\left(\boldsymbol{x}\right)$,
see Eqs. \eqref{eq:DefGamma} and \eqref{eq:DefV}, and no trace is
left by the nonlinearity $\boldsymbol{R}\left(\boldsymbol{x}\right)$.
Equation \eqref{eq:EqCase22LEq2LinAss} has to be solved with the
initial condition 
\begin{align}
\boldsymbol{\mathcal{P}}_{\boldsymbol{\mathcal{S}}}\boldsymbol{X}_{L}\left(0\right) & =\boldsymbol{\mathcal{P}}_{\boldsymbol{\mathcal{S}}}\boldsymbol{x}_{0}.
\end{align}
The terminal condition, Eq. \eqref{eq:InnerLTermCond}, becomes 
\begin{align}
\lim_{\tau_{L}\rightarrow\infty}\boldsymbol{\mathcal{P}}_{\boldsymbol{\mathcal{S}}}\boldsymbol{X}_{L}\left(\tau_{L}\right) & =\boldsymbol{\mathcal{P}}_{\boldsymbol{\mathcal{S}}}\boldsymbol{x}_{d}\left(t_{0}\right)-\boldsymbol{\Omega}_{\boldsymbol{\mathcal{S}}}\boldsymbol{\mathcal{A}}^{T}\boldsymbol{\mathcal{Q}}_{\boldsymbol{\mathcal{S}}}^{T}\boldsymbol{\lambda}_{O}\left(t_{0}\right).\label{eq:Eq4382}
\end{align}
The existence of the limit Eq. \eqref{eq:Eq4382} implies
\begin{align}
\lim_{\tau_{L}\rightarrow\infty}\boldsymbol{\mathcal{P}}_{\boldsymbol{\mathcal{S}}}\boldsymbol{X}_{L}'\left(\tau_{L}\right) & =\boldsymbol{0}.
\end{align}
Comparing this limit with the limit $\tau_{L}\rightarrow\infty$ of
Eq. \eqref{eq:EqCase22LEq2LinAss}, multiplied by $\boldsymbol{\Omega}_{\boldsymbol{\mathcal{S}}}$
from the left, yields indeed Eq. \eqref{eq:Eq4382}.

The right inner equations \eqref{eq:Case22REq1} and \eqref{eq:Case22REq3}
become
\begin{align}
\boldsymbol{\mathcal{Q}}_{\boldsymbol{\mathcal{S}}}^{T}\boldsymbol{\Lambda}_{R}'\left(\tau_{R}\right) & =\boldsymbol{0}, & \boldsymbol{\mathcal{Q}}_{\boldsymbol{\mathcal{S}}}\boldsymbol{X}_{R}'\left(\tau_{R}\right) & =\boldsymbol{0}.
\end{align}
The initial conditions Eqs. \eqref{eq:QSXRSol} and \eqref{eq:RightMatching1}
lead to
\begin{align}
\boldsymbol{\mathcal{Q}}_{\boldsymbol{\mathcal{S}}}^{T}\boldsymbol{\Lambda}_{R}\left(\tau_{R}\right) & =\boldsymbol{\mathcal{Q}}_{\boldsymbol{\mathcal{S}}}^{T}\boldsymbol{\lambda}_{O}\left(t_{1}\right), & \boldsymbol{\mathcal{Q}}_{\boldsymbol{\mathcal{S}}}\boldsymbol{X}_{R}\left(\tau_{R}\right) & =\boldsymbol{\mathcal{Q}}_{\boldsymbol{\mathcal{S}}}\boldsymbol{x}_{1}.
\end{align}
The remaining right inner equation \eqref{eq:Case22REq2} is 
\begin{align}
\dfrac{\partial}{\partial\tau_{R}}\left(\boldsymbol{\Gamma}_{\boldsymbol{\mathcal{S}}}\left(\boldsymbol{X}_{R}\left(\tau_{R}\right)\right)\boldsymbol{X}_{R}'\left(\tau_{R}\right)\right) & =-\boldsymbol{\mathcal{P}}_{\boldsymbol{\mathcal{S}}}^{T}\boldsymbol{\mathcal{V}}^{T}\left(\boldsymbol{X}_{R}\left(\tau_{R}\right),\boldsymbol{X}_{R}'\left(\tau_{R}\right)\right)\boldsymbol{\Gamma}_{\boldsymbol{\mathcal{S}}}\left(\boldsymbol{X}_{R}\left(\tau_{R}\right)\right)\boldsymbol{X}_{R}'\left(\tau_{R}\right)\nonumber \\
 & +\boldsymbol{\mathcal{P}}_{\boldsymbol{\mathcal{S}}}^{T}\boldsymbol{\mathcal{A}}^{T}\boldsymbol{\mathcal{Q}}_{\boldsymbol{\mathcal{S}}}^{T}\boldsymbol{\lambda}_{O}\left(t_{1}\right)+\boldsymbol{\mathcal{P}}_{\boldsymbol{\mathcal{S}}}^{T}\boldsymbol{\mathcal{S}}\boldsymbol{\mathcal{P}}_{\boldsymbol{\mathcal{S}}}\left(\boldsymbol{X}_{R}\left(\tau_{R}\right)-\boldsymbol{x}_{d}\left(t_{1}\right)\right),\label{eq:EqCase22REq2LinAss}
\end{align}
which is to be solved together with the initial condition 
\begin{align}
\boldsymbol{\mathcal{P}}_{\boldsymbol{\mathcal{S}}}\boldsymbol{X}_{R}\left(0\right) & =\boldsymbol{\mathcal{P}}_{\boldsymbol{\mathcal{S}}}\boldsymbol{x}_{1}.
\end{align}
Similar as for the left side, the terminal condition, Eq. \eqref{eq:InnerRTermCond},
\begin{align}
\lim_{\tau_{R}\rightarrow\infty}\boldsymbol{\mathcal{P}}_{\boldsymbol{\mathcal{S}}}\boldsymbol{X}_{R}\left(\tau_{R}\right) & =\boldsymbol{\mathcal{P}}_{\boldsymbol{\mathcal{S}}}\boldsymbol{x}_{d}\left(t_{1}\right)-\boldsymbol{\Omega}_{\boldsymbol{\mathcal{S}}}\boldsymbol{\mathcal{A}}^{T}\boldsymbol{\mathcal{Q}}_{\boldsymbol{\mathcal{S}}}^{T}\boldsymbol{\lambda}_{O}\left(t_{1}\right)\label{eq:Eq4388}
\end{align}
is already satisfied because Eq. \eqref{eq:Eq4388} implies
\begin{align}
\lim_{\tau_{R}\rightarrow\infty}\boldsymbol{\mathcal{P}}_{\boldsymbol{\mathcal{S}}}\boldsymbol{X}_{R}'\left(\tau_{R}\right) & =\boldsymbol{0}.
\end{align}
The solution for the control signal is given by
\begin{align}
\boldsymbol{u}\left(t\right) & =\begin{cases}
\boldsymbol{\mathcal{B}}_{\boldsymbol{\mathcal{S}}}^{g}\left(\boldsymbol{x}_{0}\right)\left(\boldsymbol{\dot{x}}_{O}\left(t_{0}\right)-\boldsymbol{R}\left(\boldsymbol{x}_{0}\right)+2\left(\boldsymbol{x}_{O}\left(t_{0}\right)-\boldsymbol{x}_{0}\right)\delta\left(t-t_{0}\right)\right), & t=t_{0},\\
\boldsymbol{\mathcal{B}}_{\boldsymbol{\mathcal{S}}}^{g}\left(\boldsymbol{x}_{O}\left(t\right)\right)\left(\boldsymbol{\dot{x}}_{O}\left(t\right)-\boldsymbol{R}\left(\boldsymbol{x}_{O}\left(t\right)\right)\right), & t_{0}<t<t_{1},\\
\boldsymbol{\mathcal{B}}_{\boldsymbol{\mathcal{S}}}^{g}\left(\boldsymbol{x}_{1}\right)\left(\boldsymbol{\dot{x}}_{O}\left(t_{1}\right)-\boldsymbol{R}\left(\boldsymbol{x}_{1}\right)-2\left(\boldsymbol{x}_{O}\left(t_{1}\right)-\boldsymbol{x}_{1}\right)\delta\left(t_{1}-t\right)\right) & t=t_{1}.
\end{cases}
\end{align}
In conclusion, for $\epsilon=0$ and valid linearizing assumption,
the exact state solution is given by the linear outer equations \eqref{eq:OuterEquations}
accompanied by jumps at the time domain boundaries. Thus, the analytical
\textit{form} of the exact solution is solely determined by linear
equations. However, the \textit{existence} of the exact solution relies
not only on the existence of outer solutions to Eq. \eqref{eq:OuterEquations},
but also on the existence of solutions to the generally nonlinear
inner equations \eqref{eq:EqCase22LEq2LinAss} and \eqref{eq:EqCase22REq2LinAss}.
Only the existence of inner solutions guarantees the existence of
jumps connecting the initial and terminal conditions with the outer
solution. It is in this sense that we are able to speak about an
underlying linear structure of nonlinear optimal trajectory tracking.

\subsection{Discussion}

Analytical approximations for optimal trajectory tracking of nonlinear
affine dynamical systems are developed in this section. In contrast
to Chapter \ref{chap:ExactlyRealizableTrajectories}, which discusses
only exactly realizable trajectories, the results given here are valid
for arbitrary desired trajectories. The general structure of the solution
for small regularization parameter $0\leq\epsilon\ll1$ is unveiled.
The $n$ state components $\boldsymbol{x}\left(t\right)$ are separated
by the two complementary projectors $\boldsymbol{\mathcal{P}}_{\boldsymbol{\mathcal{S}}}\left(\boldsymbol{x}\right)$
and $\boldsymbol{\mathcal{Q}}_{\boldsymbol{\mathcal{S}}}\left(\boldsymbol{x}\right)$,
while the $n$ co-state components $\boldsymbol{\lambda}\left(t\right)$
are separated by the transposed projectors $\boldsymbol{\mathcal{P}}_{\boldsymbol{\mathcal{S}}}^{T}\left(\boldsymbol{x}\right)$
and $\boldsymbol{\mathcal{Q}}_{\boldsymbol{\mathcal{S}}}^{T}\left(\boldsymbol{x}\right)$.

For all $\epsilon>0$, the dynamics of an optimal control system takes
place in the combined state space of dimension $2n$ of state $\boldsymbol{x}$
and co-state $\boldsymbol{\lambda}$ and is governed by $2n$ first
order ODEs. The exact solution for $\epsilon=0$ is governed by $2\left(n-p\right)$
first order ODEs for the state components $\boldsymbol{\mathcal{Q}}_{\boldsymbol{\mathcal{S}}}\left(\boldsymbol{x}\left(t\right)\right)\boldsymbol{x}\left(t\right)$
and the co-state components $\boldsymbol{\mathcal{Q}}_{\boldsymbol{\mathcal{S}}}^{T}\left(\boldsymbol{x}\left(t\right)\right)\boldsymbol{\lambda}\left(t\right)$,
respectively. These equations are called the outer equations and given
by \eqref{eq:ExactSol1} and \eqref{eq:ExactSol2}, respectively.
The $2p$ state components $\boldsymbol{\mathcal{P}}_{\boldsymbol{\mathcal{S}}}\left(\boldsymbol{x}\left(t\right)\right)\boldsymbol{x}\left(t\right)$
and co-state components $\boldsymbol{\mathcal{P}}_{\boldsymbol{\mathcal{S}}}^{T}\left(\boldsymbol{x}\left(t\right)\right)\boldsymbol{\lambda}\left(t\right)$
are given by algebraic equations. The part $\allowbreak\boldsymbol{\mathcal{P}}_{\boldsymbol{\mathcal{S}}}^{T}\left(\boldsymbol{x}\left(t\right)\right)\boldsymbol{\lambda}\left(t\right)=\boldsymbol{0}$
vanishes for all times, while the part $\boldsymbol{\mathcal{P}}_{\boldsymbol{\mathcal{S}}}\left(\boldsymbol{x}\left(t\right)\right)\boldsymbol{x}\left(t\right)$
is given by the algebraic equation \eqref{eq:ExactSol3} inside the
time domain, $t_{0}<t<t_{1}$. Th $2p$ algebraic equations restrict
the dynamics to a hypersurface of dimension $2\left(n-p\right)$ embedded
in the extended phase space of dimension $2n$. For all times except
for the beginning, $t=t_{0}$, and end of the time interval, $t=t_{1}$,
the system is evolving on the so-called \textit{singular surface}
\cite{bryson1969applied}. Note that the singular surface is time
dependent if the desired trajectory $\boldsymbol{x}_{d}\left(t\right)$
is time dependent.

The outer equations \eqref{eq:ExactSol1} and \eqref{eq:ExactSol2}
are $2\left(n-p\right)$ first order ODEs which allow for $2\left(n-p\right)$
initial conditions. This is not enough to accommodate all $2n$ initial
and terminal conditions given by Eq. \eqref{eq:GenDynSysInitTermCond}.
For $\epsilon=0$, this results in instantaneous and discontinuous
transitions, or jumps, at the time domain boundaries. At $t=t_{0}$,
a jump from the initial condition $\boldsymbol{x}_{0}$ onto the singular
surface occurs. Similarly, at $t=t_{1}$, a jump from the singular
surface onto the terminal condition $\boldsymbol{x}_{1}$ takes place.
These jumps manifest as discontinuities in the $2p$ state components
$\boldsymbol{\mathcal{P}}_{\boldsymbol{\mathcal{S}}}\left(\boldsymbol{x}\left(t\right)\right)\boldsymbol{x}\left(t\right)$.
The heights and directions of the jumps, see Eqs. \eqref{eq:HeightL}
and \eqref{eq:HeightR}, are given by differences between the initial
and terminal conditions and the initial and terminal values of the
outer solutions, respectively. The jumps are mediated by control impulses
in form of Dirac delta functions located at the beginning and the
end of the time interval. The direction of a delta kick, given by
the coefficient of the Dirac delta function, is parallel to the direction
of the jump occurring at the same instant. The strength of the kick
is twice the height of the jump. Intuitively, the reason is that the
delta kicks are located right at the time domain boundaries such that
only half of the kicks contribute to the time evolution.

Section \ref{sub:ExactSolutionForEpsilon0} demonstrates that the
exact solution for $\epsilon=0$ is entirely expressed in terms of
the outer solutions given by Eqs. \eqref{eq:ExactSol1}, \eqref{eq:ExactSol2},
and \eqref{eq:ExactSol3}. No trace remains of the inner solutions
except the mere existence and location of the jumps at the time domain
boundaries. Section \ref{sub:LinearizingAssumptionGeneral} draws
the final conclusion and states sufficient conditions for the linearity
of the outer equations. Because of their linearity, a formal closed
form solution valid for arbitrary desired trajectories $\boldsymbol{x}_{d}\left(t\right)$
can be given. This establishes an underlying linear structure of nonlinear
unregularized optimal trajectory tracking for affine control systems
satisfying the linearizing assumption. This finding constitutes the
major result of this thesis.

While the form of the exact state trajectory for $\epsilon=0$ is
entirely given by the outer equations, its existence relies on the
existence of solutions to the inner equations for Case 1 and Case
2.2. The existence of inner solutions for appropriate initial, terminal,
and matching conditions ensures the existence of jumps connecting
initial and terminal conditions with the singular surface. If inner
solutions do not exist for Case 1 and Case 2.2, different or additional
scaling regimes may exist which ensure the existence of jumps. The
linearizing assumption eliminates all nonlinear terms originating
from the nonlinearity $\boldsymbol{R}\left(\boldsymbol{x}\right)$
in the outer equations and the inner equations. However, the inner
equations are generally nonlinear, with nonlinear terms originating
from a state-dependent coupling matrix $\boldsymbol{\mathcal{B}}\left(\boldsymbol{x}\right)$.
For a constant coupling matrix, the inner equations are linear as
well. In conclusion, if $\epsilon=0$ and the linearizing assumption
holds, the form of the exact state and co-state trajectory is given
by linear ODEs, but their existence relies on additional, generally
nonlinear ODEs.

Having obtained analytical results for optimal open loop control,
it is in principle possible to extend this result to continuous time
and continuous time-delayed feedback control. The computations proceed
along the same lines as in Section \ref{sec:OptimalFeedback} by promoting
the initial state to a functional of the controlled state. For the
perturbation expansion in this section, sharp terminal conditions
$\boldsymbol{x}\left(t_{1}\right)=\boldsymbol{x}_{1}$ were assumed.
Analytically, this is the simplest choice, but requires the system
to be controllable. An extension to more general terminal conditions
is desirable but not straightforward.

The projectors $\boldsymbol{\mathcal{P}}_{\boldsymbol{\mathcal{S}}}$
and $\boldsymbol{\mathcal{Q}}_{\boldsymbol{\mathcal{S}}}$ play an
essential role for the solution. Both projectors are derived with
the help of the generalized Legendre-Clebsch condition for singular
optimal control in Section \ref{sub:TheGeneralizedLegendreClebschConditions}.
While $\boldsymbol{\mathcal{P}}_{\boldsymbol{\mathcal{S}}}$ depends
on the matrix of weighting coefficients $\boldsymbol{\mathcal{S}}$,
the projector $\boldsymbol{\mathcal{P}}$ defined in Chapter \ref{chap:ExactlyRealizableTrajectories}
is independent of $\boldsymbol{\mathcal{S}}$ and $\boldsymbol{\mathcal{P}}_{\boldsymbol{\mathcal{S}}}$
reduces to $\boldsymbol{\mathcal{P}}$ for $\boldsymbol{\mathcal{S}}=\boldsymbol{1}$.
The necessity to use $\boldsymbol{\mathcal{P}}_{\boldsymbol{\mathcal{S}}}$
instead of $\boldsymbol{\mathcal{P}}$ becomes obvious in the derivation
of Eq. \eqref{eq:OuterLeadingOrder21}, which cannot be obtained with
projector $\boldsymbol{\mathcal{P}}$. For exactly realizable desired
trajectories $\boldsymbol{x}_{d}\left(t\right)$, it is irrelevant
which projector is used. The control and controlled state obtained
with $\boldsymbol{\mathcal{P}}_{\boldsymbol{\mathcal{S}}}$ are independent
of $\boldsymbol{\mathcal{S}}$ and identical to results obtained with
$\boldsymbol{\mathcal{P}}$, see Section \ref{sub:TheGeneralizedLegendreClebschConditions}
for a proof. In contrast, optimal trajectory tracking for arbitrary
desired trajectories $\boldsymbol{x}_{d}\left(t\right)$ yields a
control signal and controlled state trajectory which depends explicitly
on $\boldsymbol{\mathcal{S}}$. Analogously, the linearizing assumption
introduced in Section \ref{sub:LinearizingAssumptionGeneral} depends
on $\boldsymbol{\mathcal{S}}$ and is different from the linearizing
assumption in Section \ref{sec:LinearizingAssumption}. However, the
matrices $\boldsymbol{\mathcal{Q}}_{\boldsymbol{\mathcal{S}}}$ and
$\boldsymbol{\mathcal{Q}}$ are similar and have identical diagonal
representations. This suggests that if a linearizing assumption holds
in terms of $\boldsymbol{\mathcal{Q}}_{\boldsymbol{\mathcal{S}}}$,
it also holds in terms of $\boldsymbol{\mathcal{Q}}$, and vice versa.
A rigorous proof of this conjecture is desirable.

\section{\label{sec:4Conclusions}Conclusions}

\subsection{Analytical results for \texorpdfstring{$\epsilon\rightarrow 0$}{epsilon -> 0}}

Analytical approximations for optimal trajectory tracking in nonlinear
affine control systems were derived in this chapter. The regularization
parameter $\epsilon$ is used as the small parameter for a perturbation
expansion, and the solutions become exact for $\epsilon=0$. As discussed
in Section \ref{sub:DiscussionNecessaryOptimality}, the case $\epsilon=0$
can be seen as the limit of realizability of a certain desired trajectory
$\boldsymbol{x}_{d}\left(t\right)$. No other control, be it open
or closed loop control, can enforce a state trajectory $\boldsymbol{x}\left(t\right)$
with a smaller distance to the desired state trajectory $\boldsymbol{x}_{d}\left(t\right)$.
Importantly, the regularization parameter originates solely from
the formulation of the control problem. The system dynamics is exactly
taken into account. The analytical approximations do neither require
any simplifying assumptions about the strength of nonlinearities,
as e.g. weak nonlinearities, nor about the separation of time scales
between different state components, or similar. To solve the equations
derived by the perturbative treatment in closed form, however, the
nonlinearity $\boldsymbol{R}\left(\boldsymbol{x}\right)$ has to have
a simple structure with respect to the coupling matrix $\boldsymbol{\mathcal{B}}\left(\boldsymbol{x}\right)$.
This structure is defined in abstract notation valid for a general
affine control system in Eqs. \eqref{eq:LinearizingAssumptionS_1}
and \eqref{eq:LinearizingAssumptionS_2} and called the linearizing
assumption. Cast in words, this assumption becomes ``the control
signals act on the state components governed by nonlinear equations,
and all other components are governed by linear equations''. The
linearizing assumption results in linear equations for unregularized
nonlinear optimal trajectory tracking. Only due to this linearity
it is possible to derive solutions in closed form valid for arbitrary
desired trajectories and arbitrary initial and terminal conditions.
Even if no general analytical solution is known for the uncontrolled
dynamics, the optimally controlled system can be solved analytically.
While the linearizing assumption introduced in Section \ref{sec:LinearizingAssumption}
applies only to exactly realizable desired trajectories, its applicability
is extended here to arbitrary desired trajectories. Thus, we proved
that linear structures underlying nonlinear optimal trajectory tracking
are possible.

The analytical treatment is based on a reinterpretation of a singular
optimal control problem as a singularly perturbed system of differential
equations. This reinterpretation is valid for all optimal control
problems with affine control signals. In the light of this reinterpretation,
it is now possible to understand the role of $\epsilon$ more clearly.
In particular, the behavior of unregularized optimal trajectory tracking
in many affine control systems can be outlined, even if they do not
satisfy the linearizing assumption from Section \ref{sub:LinearizingAssumptionGeneral}.

For all $\epsilon>0$, the dynamics of an optimal control system takes
place in the extended state space, i.e., in the combined space of
state $\boldsymbol{x}$ and co-state $\boldsymbol{\lambda}$ of dimension
$2n$, with $n$ being the number of state space components. For $\epsilon=0$,
the dynamics is restricted by $2p$ algebraic equations to a hypersurface
of dimension $2\left(n-p\right)$, with $p$ being the number of independent
control signals. For all times except at the beginning and the end
of the time interval, the system is evolving on the so-called \textit{singular
surface} \cite{bryson1969applied}. At the initial time, a kick in
form of a Dirac delta function mediated by the control signal induces
an instantaneous transition from the initial state onto the singular
surface. Similarly, at the terminal time, a delta-like kick induces
an instantaneous transition from the singular surface to the terminal
state. These instantaneous transitions render certain state components
discontinuous at the initial and terminal time, respectively. For
$\epsilon>0$, the discontinuities of the state are smoothed out in
form of boundary layers, i.e., continuous transition regions with
a slope controlled by the value of $\epsilon$. The control signals
are finite and exhibit a sharp peak at the time domain boundaries
with an amplitude inversely proportional to $\epsilon$.

The general picture of the behavior of unregularized optimal control
problems clearly explains the necessity of a regularization term in
the cost functional $\mathcal{J}$. While the behavior for $\epsilon=0$
is relatively easy to understand and determined by simpler equations
than for $\epsilon>0$, the result is mathematically inconvenient.
For $\epsilon=0$, it is not possible to find a solution for the optimal
controlled state trajectory in terms of continuous functions. Even
worse, the solution for the control signal must be expressed in terms
of the Dirac delta function, i.e., in terms of distributions. Throughout
this thesis, no attention is paid to the function spaces to which
the controlled state trajectory and the control signal belongs. Everything
is assumed to be sufficiently well behaved. However, the analytical
treatment for $\epsilon=0$ leads right to the importance of such
questions. A mathematically more precise characterization of the different
function spaces involved in the problems of optimal trajectory tracking
for $\epsilon=0$ and $\epsilon>0$ is desirable.

Here, we derived perturbative solutions for small $\epsilon$. Note,
however, that finding such solutions is only a first step in a mathematically
rigorous perturbative treatment. In a second step, existence and uniqueness
of the outer and relevant inner equations together with their initial,
terminal, and matching conditions must be established. Third, the
reliability of the approximate result must be demonstrated. This is
usually done by estimates in form of rigorous inequalities which determine
how much the approximate solution deviates from the exact result for
a given value of $\epsilon$. Another point deserving more mathematical
rigor concerns the linearizing assumption. Here, we showed that the
linearizing assumption is a sufficient condition for a linear structure
of optimal trajectory tracking. The question arises if it is also
necessary. Other classes of affine control systems which violate the
linearizing assumption but exhibit an underlying linear structure
may exist.

Due to the limited resolution in numerical simulations, it is at least
difficult, if not impossible, to find a faithful numerical representation
of the solution to optimal trajectory tracking for $\epsilon=0$.
To ensure a solution to an optimal control problem in terms of numerically
treatable functions, a finite value of $\epsilon$ is indispensable.
For the two-dimensional dynamical systems of Section \ref{sec:TwoDimensionalDynamicalSystem},
the width of the boundary layers is directly proportional to the value
of $\epsilon$. Thus, a temporal resolution $\Delta t$ smaller than
$\epsilon$, $\Delta t<\epsilon$, will not be able to numerically
resolve these boundary layers, and is likely to lead to large numerical
errors. Indeed, comparing a numerical result obtained for $\Delta t=\epsilon$
with its analytical counterpart reveals that the largest differences
occur in the boundary layer regions, see Fig. \ref{fig:Fig3} of Example
\ref{ex:OptimallyControlledFHN} in Section \ref{sec:ComparisonWithNumericalResults}.
Note that the initial boundary layer plays an important role for the
future time evolution of the system. An erroneous computation of this
transition region leads to a perturbed initial value on the singular
surface. This is a problem if the system is sensitive with respect
to perturbations of the initial conditions. However, here the problem
can be more severe due to the time dependence of the singular surface.
A desired trajectory $\boldsymbol{x}_{d}\left(t\right)$ changing
rapidly during the initial transient is likely to cause a rapidly
changing singular surface, and an erroneous computation of the initial
boundary layer might lead to a different singular surface altogether.

Combining analytical approximations with numerical methods can be
fruitful for many applications. Analytical solutions for optimal control,
even if only approximately valid, can provide a suitable initial guess
for iterative optimal control algorithms, and result in a considerable
decrease of computational cost. Imaginable are applications to real
time computations of optimal control, which is still an ambitious
task even with modern-day fast computers. Furthermore, analytical
approximations can be used to test the accuracy of numerical optimal
control algorithms and estimate errors caused by discretization.

In technical application and experiments, it is impossible to generate
diverging control signals, and in general, an experimental realization
of unregularized optimal control systems is impossible. Nevertheless,
understanding the behavior of the control system in the limit $\epsilon\rightarrow0$
can be very useful for applications. For example, to avoid any steep
transitions and large control amplitudes, one can exploit the knowledge
about the initial conditions for the singular surface. If the initial
state of the system can be prepared, the initial state could be chosen
to lie on the singular surface. Thereby, any initial steep transitions
can be prevented, or at least minimized. Furthermore, if the initial
state cannot be prepared, it might still be possible to design the
desired trajectory $\boldsymbol{x}_{d}\left(t\right)$ such that the
initial state lies on the singular surface. This is only one example
how an analytical solution can be utilized for the planning of desired
trajectories. Another example is the discovery from Example \ref{ex:ControlledPendulum}
that the desired velocity over time can only be controlled up to a
constant shift for mechanical control systems in one spatial dimension.
In general, analytical solutions of optimal trajectory tracking for
arbitrary desired trajectory $\boldsymbol{x}_{d}\left(t\right)$ enable
to compare the performance of controlled state trajectories for different
choices of desired trajectories. This is useful if the desired trajectory
is not entirely fixed by the problem setting but exhibits some degrees
of freedom. In a second step, these can be optimized with respect
to other aspects as e.g. the control amplitude. Such a procedure is
nearly impossible for numerical optimal control due to the computational
cost of numerical algorithms.

It is clear that the class of optimal control systems with underlying
linear structure is much smaller than the class of feedback linearizable
systems. The reason is that the coupled state and co-state equations
are more complex, and there are many more sources for nonlinearity.
Consider the necessary optimality conditions for optimal trajectory
tracking, 
\begin{align}
\boldsymbol{0} & =\epsilon^{2}\boldsymbol{u}\left(t\right)+\boldsymbol{\mathcal{B}}^{T}\boldsymbol{\lambda}\left(t\right),\\
\boldsymbol{\dot{x}}\left(t\right) & =\boldsymbol{R}\left(\boldsymbol{x}\left(t\right)\right)+\boldsymbol{\mathcal{B}}\boldsymbol{u}\left(t\right),\\
-\boldsymbol{\dot{\lambda}}\left(t\right) & =\nabla\boldsymbol{R}^{T}\left(\boldsymbol{x}\left(t\right)\right)\boldsymbol{\lambda}\left(t\right)+\boldsymbol{\mathcal{S}}\left(\boldsymbol{x}\left(t\right)-\boldsymbol{x}_{d}\left(t\right)\right).
\end{align}
The controlled state equation is coupled with the adjoint equation
via the transposed Jacobian of $\boldsymbol{R}$. Additionally, the
inhomogeneity $\boldsymbol{\mathcal{S}}\left(\boldsymbol{x}\left(t\right)-\boldsymbol{x}_{d}\left(t\right)\right)$
of the adjoint equation depends linearly on $\boldsymbol{\ensuremath{x}}$,
and any nonlinear transformation of the state results in an inhomogeneity
depending nonlinearly on $\boldsymbol{x}$. Applying a nonlinear state
transformation, as it is often required by feedback linearization,
leads to new nonlinearities in the adjoint equation. Thus, optimal
control systems cannot fully benefit from the linear structure underlying
feedback linearizable systems. The class of exactly linear optimal
control systems is certainly much smaller than the class of feedback
linearizable systems.

An important problem is the impact of noise on optimal control. Fundamental
results exist for linear optimal control. The standard problem of
linear optimal feedback control is the so-called linear-quadratic
regulator, a linear controlled state equation together with a cost
function quadratic in the state. The linear-quadratic-Gaussian control
problem considers the linear-quadratic regulator together with additive
Gaussian white noise in the state equation as well as for state measurements
\cite{bryson1969applied}. The discovery of linear structures underlying
nonlinear trajectory tracking might enable a similar investigation
for control systems satisfying the linearizing assumption. In this
context, we mention Ref. \cite{kappen2005linear} which presents a
linear theory for the control of nonlinear stochastic systems. The
approach in \cite{kappen2005linear} relies on an exact linearization
of the Hamilton-Jacobi-Bellman equation for stochastic systems by
a Cole-Hopf transform. However, the method in \cite{kappen2005linear}
is restricted to systems with identical numbers of control signals
and state components, $n=p$.

\subsection{Weak and strong coupling}

Here, optimal trajectory tracking is characterized for the whole range
of the regularization parameter $\epsilon\geq0$. Consider a system
which satisfies the linearizing assumption,
\begin{align}
\boldsymbol{\mathcal{Q}}\boldsymbol{R}\left(\boldsymbol{x}\right) & =\boldsymbol{\mathcal{Q}}\boldsymbol{\mathcal{A}}\boldsymbol{x}+\boldsymbol{\mathcal{Q}}\boldsymbol{b},\label{eq:Eq816}
\end{align}
such that the control signals $\boldsymbol{u}\left(t\right)$ act
on nonlinear state equations. Due to the special structure defined
by Eq. \eqref{eq:Eq816}, the control is able to counteract the nonlinearity.
In general, it is able to do so only if it is allowed to have an arbitrarily
large amplitude. The amplitude of the control signal is closely related
to the value of the regularization parameter $\epsilon$. In the cost
functional for optimal trajectory tracking,
\begin{align}
\mathcal{J}= & \frac{\alpha}{2}\intop_{t_{0}}^{t_{1}}dt\left(\boldsymbol{x}\left(t\right)-\boldsymbol{x}_{d}\left(t\right)\right)^{2}+\frac{\epsilon^{2}}{2}\intop_{t_{0}}^{t_{1}}dt\left(\boldsymbol{u}\left(t\right)\right)^{2},
\end{align}
the regularization term $\sim\epsilon^{2}$ penalizes large control
signals. Depending on the value of  $\epsilon$, different regimes
can be identified. Clearly, in the limit of $\epsilon\rightarrow\infty$,
any non-vanishing control signal leads to a diverging value of $\mathcal{J}$.
Thus, the limit $\epsilon\rightarrow\infty$ implies a vanishing control
signal, $\boldsymbol{u}\left(t\right)\equiv\boldsymbol{0}$, and corresponds
to the uncontrolled system. If $\epsilon$ is much larger than $1$
but finite, $\epsilon\gg1$, the control is allowed to have a small
maximum amplitude. This regime can be regarded as the weak coupling
limit. The nonlinearity $\boldsymbol{R}\left(\boldsymbol{x}\right)$
dominates the system dynamics even if the system satisfies the linearizing
assumption Eq. \eqref{eq:Eq816}. For decreasing values of $\epsilon$,
the control exerts a growing influence on the system, until the regime
with $1\gg\epsilon>0$ is reached where the control dominates over
the nonlinearity. If $\epsilon$ vanishes identically, $\epsilon=0$,
the control is given by 
\begin{align}
\boldsymbol{u}\left(t\right) & =\boldsymbol{\mathcal{B}}^{+}\left(\boldsymbol{x}\left(t\right)\right)\left(\boldsymbol{\dot{x}}\left(t\right)-\boldsymbol{R}\left(\boldsymbol{x}\left(t\right)\right)\right).
\end{align}
In general, without any constraints on the state $\boldsymbol{x}$,
the nonlinearity $\boldsymbol{R}\left(\boldsymbol{x}\right)$ can
attain arbitrarily large values. Only if the control is allowed to
attain arbitrarily large values as well, it is able to counteract
an arbitrary nonlinearity $\boldsymbol{R}$ evaluated at an arbitrary
state value $\boldsymbol{x}$. Only for $\epsilon=0$, one can expect
an exactly linear behavior of nonlinear optimal control systems independent
of the nonlinearity $\boldsymbol{R}$.

Indeed, the analytical results indicate that the control signal scales
as $1/\epsilon$. Although the analytical results are only valid for
small $\epsilon$, this corroborates the above considerations for
the entire range of values of $\epsilon$. In view of the underlying
linear structure of a nonlinear optimal control for $\epsilon=0$,
one might ask if something can be learned about the nonlinear uncontrolled
system by analyzing the linear controlled problem? The answer is clearly
no, because the limit of an uncontrolled system is the opposite limit
of an arbitrarily strongly controlled system assumed for the perturbative
treatment.

This foregoing reasoning might explain why methods like feedback linearization,
which exploit an underlying linear structure of controlled systems,
are relatively unfamiliar in the nonlinear dynamics community and
among physicists in general. Physicists tend to approach controlled
systems from the viewpoint of uncontrolled systems. Having understood
the manifold of solutions to the uncontrolled system, which is the
traditional topic of nonlinear dynamics, the natural approach to controlled
systems is to regard the control as a perturbation. This corresponds
to the weak coupling limit mentioned above. Treating controlled systems
in this limit is often sufficient to discuss stabilization of unstable
attractors and similar topics. Such control tasks can often be achieved
with non-invasive control signals. Usually, small control amplitudes
are technically more feasible, and generally preferred over large
control amplitudes. On the downside, concentrating solely on the weak
coupling limit misses the fact that many nonlinear control systems,
as e.g. feedback-linearizable systems, have an underlying linear structure.
Because basically all exact linearizations of control systems work
by transforming the control such that it cancels the nonlinearity,
this underlying linear structure can only be exploited in the strong
coupling limit and in the absence of constraints for the control.
Any a priori assumptions about the maximum value of the control amplitude,
enforced by a regularization term or inequality constraints in case
of optimal control, destroy the underlying linear structure.

However, the strong coupling limit does not always imply very large
or even diverging control amplitudes. After having obtained the solutions
for the control signal as well as the controlled state trajectory,
it is possible to give a posteriori estimates on the maximum control
amplitude. Depending on the desired trajectories, initial conditions,
and system dynamics, these a posteriori estimates can be comparatively
small. The weak coupling limit only refers to a priori assumptions
about the maximum control amplitude, and excludes or penalizes large
control amplitudes from the very beginning.

%% file: chapter-5.tex
\lhead[\chaptername~\thechapter\leftmark]{}

\rhead[]{\rightmark}

\lfoot[\thepage]{}

\cfoot{}

\rfoot[]{\thepage}

\chapter{\label{chap:ControlOfReactionDiffusion}Control of reaction-diffusion
systems}

Reaction-diffusion systems model phenomena from a large variety of
fields. Examples are chemical systems \cite{KapralShowalter201210,EpsteinPojman199810},
action potential propagation in the heart and neurons \cite{keener1998mathematical1,keener1998mathematical2},
population dynamics \cite{murray1993mathematicalbiology1,murray1993mathematicalbiology2},
vegetation patterns \cite{Hardenberg2001Diversity}, and the motility
of crawling cells \cite{ziebert2011model,ziebert2013effects,lober2014modeling,aranson2014phasefield,lober2015collisions},
to name only a few. These systems possess a rich phenomenology of
solutions, ranging from homogeneous stable steady states, phase waves,
Turing patterns, stationary localized and labyrinthine patterns, traveling,
rotating and scroll waves to fully developed spatio-temporal turbulence
\cite{turing1952chemical,cross1993pfo,Hagberg1994From,kuramoto2003chemical,vanag2007localized}.
Due to the complexity and the nonlinearity of the underlying evolution
equations, their theoretical investigation relies heavily on numerical
simulations. However, more complex patterns can often be understood
as being assembled of simple ``building blocks'' as traveling fronts
and pulses. A solitary pulse in the FHN model in one spatial dimension
can be considered as being built of two propagating interfaces separating
the excited from the refractory state. These interfaces are front
solutions to a simpler reaction-diffusion system. Similarly, many
two-dimensional shapes as e.g. spiral waves can be approximated as
consisting of appropriately shifted one-dimensional pulse profiles
\cite{tyson1988singular,zykov1987simulation,pismen2006patterns,lober2013analytical,Mikhailov2011Synergetics1}.
In many cases, the simplified equations allow the inclusion of additional
effects as e.g. spatial heterogeneities \cite{lober2009nonlinear,alonso2010wave,lober2012front},
noise \cite{schimanskygeier1983efp,engel1985nif}, or curved boundaries
\cite{engel1987interaction,martens2015front}.

The control of patterns in reaction-diffusion system has received
the attention of many researchers in the past \cite{mikhailov2006control,vanag2008design}.
Due to their complexity, it makes sense to develop first a detailed
understanding of the control of simple solutions as e.g. solitary
excitation pulses. A particularly simple but still general control
task is position control of traveling waves. The position of a traveling
wave is shifted according to a prescribed trajectory in position space,
called the protocol of motion. Simultaneously, the wave profile is
kept as close as possible to the uncontrolled wave profile. An example
of open loop control in this spirit is the dragging of chemical pulses
of adsorbed CO during heterogeneous catalysis on platinum single crystal
surfaces \cite{wolff2003gentle}. In experiments with an addressable
catalyst surface, the pulse velocity was controlled by a laser beam
creating a movable localized temperature heterogeneity, resulting
in a V-shaped wave pattern \cite{wolff2001spatiotemporal,wolff2003wave}.
Theoretical studies of dragging one-dimensional chemical fronts or
phase interfaces by anchoring it to a movable parameter heterogeneity
can be found in \cite{nistazakis2002targeted,malomed2002pulled,kevrekidis2004dragging}.
While these approaches assume a fixed spatial profile of the control
signal and vary only its location, the method developed in \cite{lober2014controlling}
determines the profile, amplitude, and location of the control signal
by solving an inverse problem for the position over time of controlled
traveling waves. This control solution is close to the solution of
an appropriately formulated optimal control problem. Furthermore,
an extension allows the investigation of the stability of controlled
traveling waves \cite{lober2014stability}, which can never be taken
for granted in open loop control systems. A modification of the method
provides shaping of wave patterns in two-dimensional reaction-diffusion
systems by shifting the position of a trajectory outlining the pattern
\cite{lober2014shaping}. See also \cite{lober2014control} for a
discussion of experimental realizations. An approach similarly aiming
at the position of wave patterns is the forcing of spiral waves with
temporally periodic and spatially homogeneous control signals. This
can be utilized to guide a meandering spiral wave tip along a wide
range of open and closed hypocycloidal trajectories \cite{steinbock1993control,zykov1994external}.

Position control can be tackled by feedback control as well. In experiments
with spiral waves in the photosensitive Belousov-Zhabotinsky reaction
\cite{krug1990analysis}, the spiral wave core is steered around obstacles
using feedback signals obtained from wave activity measured at detector
points, along detector lines, or in a spatially extended control domain
\cite{Zykov2004global,zykov2004feedback,schlesner2008efficient}.
Two feedback loops were used to guide wave segments along pre-given
trajectories \cite{sakurai2002design}. Furthermore, feedback-mediated
control loops are employed in order to stabilize unstable patterns
such as plane waves undergoing transversal instabilities \cite{molnos2015control},
unstable traveling wave segments \cite{mihaliuk2002feedback}, or
rigidly rotating unstable spiral waves in the regime of stable meandering
spiral waves \cite{schlesner2006stabilization}. Another strategy
is control by imposed geometric constraints such as no-flux boundaries
\cite{paulau2014stabilization} or heterogeneities \cite{luther2011low}.

While feedback control and external forcing of reaction-diffusion
system has received much attention, optimal control of these systems
remains largely unexplored, at least within the physics community.
One reason lies in the computational cost involved in numerical approaches
to optimal open loop control of PDEs which restricts numerical investigations
to relatively small spatial domains and short time intervals. Even
worse, optimal feedback control of PDEs becomes almost intractable.
The reason is the curse of dimensionality \cite{bellman1957dynamic}.
For an $n$-dimensional dynamical system, the Hamilton-Jacobi-Bellman
equation for optimal feedback control is a PDE on an $n$-dimensional
domain, see the discussion at the beginning of Section \ref{sec:OptimalFeedback}.
For a controlled PDE, which can be regarded as a dynamical system
with $n\rightarrow\infty$ dimensions, the corresponding Hamilton-Jacobi-Bellman
equation is a PDE on a domain with $n\rightarrow\infty$ dimensions.

In view of the numerical difficulties, the analytical approach pursued
in this thesis has some benefits when compared with purely numerical
methods, and can be used to obtain solutions to optimal control for
a number of systems with relative ease. In Section \ref{sec:RDSNotation},
the formalism based on projectors is modified and applied to spatio-temporal
systems. The controlled state equation is split up in two equations
in Section \ref{sec:RDSSplitUpTheStateEquation}, and exactly realizable
distributions are introduced as the spatio-temporal analogue of exactly
realizable trajectories in Section \ref{sec:RDSExactlyRealizableDistributions}.
As an important application, the position control of traveling waves
is discussed in Section \ref{sec:PositionControlOfTravelingWaves}.
This chapter concludes with a discussion and outlook in Section \ref{sec:DiscussionAndOutlook}.

\section{\label{sec:RDSNotation}Formalism}

In this section, the formalism developed in Chapter \ref{chap:ExactlyRealizableTrajectories}
is modified and applied to spatio-temporal systems. The emphasis lies
on distributed controls, i.e., the control signal is allowed to depend
on space and time. Often, the control cannot act everywhere in position
space. For example, it might be possible to let a control act at or
close to the boundaries, but it is impossible to reach the interior
of the domain. These restrictions can be accounted for by a modification
of the projectors $\boldsymbol{\mathcal{P}}$ and $\boldsymbol{\mathcal{Q}}$.

Let the position vector $\boldsymbol{r}$ in $N$ spatial dimensions
be
\begin{align}
\boldsymbol{r} & =\left(\begin{array}{cc}
r_{1}, & r_{2},\dots,r_{N}\end{array}\right)^{T}.
\end{align}
The spatial domain is denoted by $\Omega\subset\mathbb{R}^{N}$, and
its boundary is $\Gamma=\partial\Omega\subset\mathbb{R}^{N}$. Let
$\boldsymbol{\chi}\left(\boldsymbol{r}\right)$ be a diagonal $n\times n$
matrix of characteristic functions,
\begin{align}
\boldsymbol{\chi}\left(\boldsymbol{r}\right) & =\left(\begin{array}{ccc}
\chi_{1}\left(\boldsymbol{r}\right) & \cdots & 0\\
\vdots & \ddots & \vdots\\
0 & \cdots & \chi_{n}\left(\boldsymbol{r}\right)
\end{array}\right).
\end{align}
The characteristic function $\chi_{i}\left(\boldsymbol{r}\right)$
only attains the values zero or one, 
\begin{align}
\chi_{i}\left(\boldsymbol{r}\right)= & \begin{cases}
1, & \boldsymbol{r}\in A_{i},\\
0, & \boldsymbol{r}\in\Omega\setminus A_{i},
\end{cases}i\in\left\{ 1,\dots,n\right\} .
\end{align}
The characteristic functions $\chi_{i}\left(\boldsymbol{r}\right)$
do account for the case that a control signal might only act in a
restricted region of space and not in the full domain $\Omega$. The
region in which a control signal acts on the $i$-th state component
is denoted by $A_{i}$. The total spatial region affected by control
is $A=\bigcup_{i=1}^{N}A_{i},$ and no control acts in region $\Omega\setminus A$
outside of $A$. If all state components are controlled in the same
region $A_{1}=A_{2}=\dots=A_{n}$ of space, then $\boldsymbol{\chi}\left(\boldsymbol{r}\right)$
simplifies to a multiple of the identity matrix $\boldsymbol{\chi}\left(\boldsymbol{r}\right)=\chi\left(\boldsymbol{r}\right)\boldsymbol{1}$,
with a scalar characteristic function $\chi$. If additionally, the
control acts everywhere in the spatial domain $\Omega$, then $\boldsymbol{\chi}\left(\boldsymbol{r}\right)=\boldsymbol{1}$.
The matrix $\boldsymbol{\chi}\left(\boldsymbol{r}\right)$ is a space
dependent projector on the state space. It is idempotent,
\begin{align}
\boldsymbol{\chi}\left(\boldsymbol{r}\right)\boldsymbol{\chi}\left(\boldsymbol{r}\right) & =\boldsymbol{\chi}\left(\boldsymbol{r}\right),
\end{align}
and symmetric
\begin{align}
\boldsymbol{\chi}^{T}\left(\boldsymbol{r}\right) & =\boldsymbol{\chi}\left(\boldsymbol{r}\right).
\end{align}
The projector $\boldsymbol{\psi}$ complementary to $\boldsymbol{\chi}$
is defined as
\begin{align}
\boldsymbol{\psi}\left(\boldsymbol{r}\right) & =\boldsymbol{1}-\boldsymbol{\chi}\left(\boldsymbol{r}\right),
\end{align}
such that $\boldsymbol{\psi}\left(\boldsymbol{r}\right)\boldsymbol{\chi}\left(\boldsymbol{r}\right)=\boldsymbol{0}$.

An affine controlled reaction-diffusion system for the $n$-component
state vector
\begin{align}
\boldsymbol{x}\left(\boldsymbol{r},t\right) & =\left(\begin{array}{ccc}
x_{1}\left(\boldsymbol{r},t\right), & \dots, & x_{n}\left(\boldsymbol{r},t\right)\end{array}\right)
\end{align}
with the $p$-component vector of distributed control signals
\begin{align}
\boldsymbol{u}\left(\boldsymbol{r},t\right) & =\left(\begin{array}{ccc}
u_{1}\left(\boldsymbol{r},t\right), & \dots, & u_{p}\left(\boldsymbol{r},t\right)\end{array}\right)
\end{align}
is 
\begin{align}
\partial_{t}\boldsymbol{x}\left(\boldsymbol{r},t\right) & =\boldsymbol{\mathcal{D}}\triangle\boldsymbol{x}\left(\boldsymbol{r},t\right)+\boldsymbol{R}\left(\boldsymbol{x}\left(\boldsymbol{r},t\right)\right)+\boldsymbol{\chi}\left(\boldsymbol{r}\right)\boldsymbol{\mathcal{B}}\left(\boldsymbol{x}\left(\boldsymbol{r},t\right)\right)\boldsymbol{u}\left(\boldsymbol{r},t\right).\label{eq:ControlledRDS}
\end{align}
Here, $\boldsymbol{\mathcal{D}}$ is an $n\times n$ diagonal matrix
of constant diffusion coefficients and $\triangle$ denotes the Laplacian
which, in Cartesian coordinates, assumes the form $\triangle=\sum_{i=1}^{N}\dfrac{\partial^{2}}{\partial r_{i}^{2}}$.
For simplicity, an isotropic medium is considered. The $n\times p$
coupling matrix $\boldsymbol{\mathcal{B}}\left(\boldsymbol{x}\right)$
is assumed to have full column rank, $\text{rank}\left(\boldsymbol{\mathcal{B}}\left(\boldsymbol{x}\right)\right)=p$,
for all $\boldsymbol{x}$. To shorten the notation, the $n\times p$
matrix
\begin{align}
\boldsymbol{\mathcal{B}}\left(\boldsymbol{x},\boldsymbol{r}\right) & =\boldsymbol{\chi}\left(\boldsymbol{r}\right)\boldsymbol{\mathcal{B}}\left(\boldsymbol{x}\right)
\end{align}
is introduced. For reaction-diffusion systems in finite domains, Eq.
\eqref{eq:ControlledRDS} is supplemented with appropriate boundary
conditions. A common choice are homogeneous Neumann or no flux boundary
conditions
\begin{align}
\boldsymbol{n}^{T}\left(\boldsymbol{r}\right)\left(\boldsymbol{\mathcal{D}}\nabla\boldsymbol{x}\left(\boldsymbol{r},t\right)\right) & =\boldsymbol{0},\,\boldsymbol{r}\in\Gamma.\label{eq:NeumannBoundary}
\end{align}
Here, the $N$-component vector $\boldsymbol{n}\left(\boldsymbol{r}\right)$
is the vector normal to the boundary $\Gamma$. If the diffusion coefficient
for a certain component vanishes, the boundary condition for this
component is trivially satisfied.

In principle, additional control signals acting on the domain boundary
$\Gamma$ can be introduced. In case of Neumann boundary conditions,
this corresponds to an inhomogeneity on the right hand side of Eq.
\eqref{eq:NeumannBoundary} prescribing the flux of state components
across $\Gamma$ \cite{theissen2006optimale}. See \cite{lebiedz2003manipulation}
how a desired stationary concentration profile is enforced in a reaction-diffusion
system by boundary control. Although such control schemes are important
for applications, the discussion here is restricted to distributed
controls, i.e., spatio-temporal control signals acting inside the
spatial domain. Other possible boundary conditions for reaction-diffusion
systems are Dirichlet or periodic boundary conditions. Finally, the
initial condition for Eq. \eqref{eq:ControlledRDS} is given by
\begin{align}
\boldsymbol{x}\left(\boldsymbol{r},t_{0}\right) & =\boldsymbol{x}_{0}\left(\boldsymbol{r}\right).
\end{align}

\section{\label{sec:RDSSplitUpTheStateEquation}Split up the state equation}

Similar as in earlier chapters, the control signal can be expressed
in terms of the controlled state $\boldsymbol{x}\left(\boldsymbol{r},t\right)$
as 
\begin{align}
\boldsymbol{u}\left(\boldsymbol{r},t\right) & =\boldsymbol{\mathcal{B}}^{+}\left(\boldsymbol{x}\left(\boldsymbol{r},t\right),\boldsymbol{r}\right)\left(\partial_{t}\boldsymbol{x}\left(\boldsymbol{r},t\right)-\boldsymbol{\mathcal{D}}\triangle\boldsymbol{x}\left(\boldsymbol{r},t\right)-\boldsymbol{R}\left(\boldsymbol{x}\left(\boldsymbol{r},t\right)\right)\right).\label{eq:RDSControl}
\end{align}
Here, $\boldsymbol{\mathcal{B}}^{+}\left(\boldsymbol{x},\boldsymbol{r}\right)$
denotes the $p\times n$ Moore-Penrose pseudo inverse of the matrix
$\boldsymbol{\mathcal{B}}\left(\boldsymbol{x},\boldsymbol{r}\right)$,
\begin{align}
\boldsymbol{\mathcal{B}}^{+}\left(\boldsymbol{x},\boldsymbol{r}\right) & =\left(\boldsymbol{\chi}\left(\boldsymbol{r}\right)\boldsymbol{\mathcal{B}}\left(\boldsymbol{x}\right)\right)^{+}.
\end{align}
In contrast to earlier chapters, no closed form expression can be
given for the pseudo inverse $\boldsymbol{\mathcal{B}}^{+}\left(\boldsymbol{x},\boldsymbol{r}\right)$
in the general case. The reason is that the $p\times p$ matrix $\boldsymbol{\mathcal{B}}^{T}\left(\boldsymbol{x},\boldsymbol{r}\right)\boldsymbol{\mathcal{B}}\left(\boldsymbol{x},\boldsymbol{r}\right)=\boldsymbol{\mathcal{B}}^{T}\left(\boldsymbol{x}\right)\boldsymbol{\chi}\left(\boldsymbol{r}\right)\boldsymbol{\mathcal{B}}\left(\boldsymbol{x}\right)$
does not necessarily have rank $p$ and can therefore not be inverted.
In general, the rank of $\boldsymbol{\mathcal{B}}^{T}\left(\boldsymbol{x}\right)\boldsymbol{\chi}\left(\boldsymbol{r}\right)\boldsymbol{\mathcal{B}}\left(\boldsymbol{x}\right)$
depends on the spatial coordinate $\boldsymbol{r}$. Nevertheless,
the Moore-Penrose pseudo inverse of $\boldsymbol{\mathcal{B}}^{+}\left(\boldsymbol{x},\boldsymbol{r}\right)$
does always exist and can be computed numerically with the help of
singular value decomposition, for example. In some special but important
cases, explicit expressions for the pseudo inverse $\boldsymbol{\mathcal{B}}^{+}\left(\boldsymbol{x},\boldsymbol{r}\right)$
can be given. Note that $\boldsymbol{\chi}\left(\boldsymbol{r}\right)$,
being a symmetric projector, is its own pseudo inverse, i.e.,
\begin{align}
\boldsymbol{\chi}^{+}\left(\boldsymbol{r}\right) & =\boldsymbol{\chi}\left(\boldsymbol{r}\right).
\end{align}
If all state components are affected in the same region of space such
that $\boldsymbol{\chi}\left(\boldsymbol{r}\right)=\chi\left(\boldsymbol{r}\right)\boldsymbol{1}$,
then the pseudo inverse is given by
\begin{align}
\boldsymbol{\mathcal{B}}^{+}\left(\boldsymbol{x},\boldsymbol{r}\right) & =\chi\left(\boldsymbol{r}\right)\boldsymbol{\mathcal{B}}^{+}\left(\boldsymbol{x}\right)=\chi\left(\boldsymbol{r}\right)\left(\boldsymbol{\mathcal{B}}^{T}\left(\boldsymbol{x}\right)\boldsymbol{\mathcal{B}}\left(\boldsymbol{x}\right)\right)^{-1}\boldsymbol{\mathcal{B}}^{T}\left(\boldsymbol{x}\right).
\end{align}
If the number of independent control signals equals the number of
state components, $n=p$, such that $\boldsymbol{\mathcal{B}}\left(\boldsymbol{x}\right)$
is invertible, then the pseudo inverse is
\begin{align}
\boldsymbol{\mathcal{B}}^{+}\left(\boldsymbol{x},\boldsymbol{r}\right) & =\boldsymbol{\mathcal{B}}^{-1}\left(\boldsymbol{x}\right)\boldsymbol{\chi}^{+}\left(\boldsymbol{r}\right)=\boldsymbol{\mathcal{B}}^{-1}\left(\boldsymbol{x}\right)\boldsymbol{\chi}\left(\boldsymbol{r}\right).
\end{align}
Finally, if the $p\times p$ matrix $\boldsymbol{\mathcal{B}}^{T}\left(\boldsymbol{x}\right)\boldsymbol{\chi}\left(\boldsymbol{r}\right)\boldsymbol{\mathcal{B}}\left(\boldsymbol{x}\right)$
has full rank $p$ for all values of $\boldsymbol{r}\in\Omega$ and
for all states $\boldsymbol{x}$, the pseudo inverse is given by 
\begin{align}
\boldsymbol{\mathcal{B}}^{+}\left(\boldsymbol{x},\boldsymbol{r}\right) & =\left(\boldsymbol{\mathcal{B}}^{T}\left(\boldsymbol{x}\right)\boldsymbol{\chi}\left(\boldsymbol{r}\right)\boldsymbol{\mathcal{B}}\left(\boldsymbol{x}\right)\right)^{-1}\boldsymbol{\mathcal{B}}^{T}\left(\boldsymbol{x}\right)\boldsymbol{\chi}\left(\boldsymbol{r}\right).
\end{align}
Note that for any matrix $\boldsymbol{\mathcal{A}}$, its Moore-Penrose
pseudo inverse can also be expressed as \cite{albert1972regression}
\begin{align}
\boldsymbol{\mathcal{A}}^{+} & =\left(\boldsymbol{\mathcal{A}}^{T}\boldsymbol{\mathcal{A}}\right)^{+}\boldsymbol{\mathcal{A}}^{T},
\end{align}
such that $\boldsymbol{\mathcal{B}}^{+}\left(\boldsymbol{x},\boldsymbol{r}\right)$
can be written in the form
\begin{align}
\boldsymbol{\mathcal{B}}^{+}\left(\boldsymbol{x},\boldsymbol{r}\right) & =\left(\boldsymbol{\mathcal{B}}^{T}\left(\boldsymbol{x},\boldsymbol{r}\right)\boldsymbol{\mathcal{B}}\left(\boldsymbol{x},\boldsymbol{r}\right)\right)^{+}\boldsymbol{\mathcal{B}}^{T}\left(\boldsymbol{x},\boldsymbol{r}\right)\nonumber \\
 & =\left(\boldsymbol{\mathcal{B}}^{T}\left(\boldsymbol{x}\right)\boldsymbol{\chi}\left(\boldsymbol{r}\right)\boldsymbol{\mathcal{B}}\left(\boldsymbol{x}\right)\right)^{+}\boldsymbol{\mathcal{B}}^{T}\left(\boldsymbol{x}\right)\boldsymbol{\chi}\left(\boldsymbol{r}\right).\label{eq:BDaggerAlternative}
\end{align}
Together with the projector property of $\boldsymbol{\chi}\left(\boldsymbol{r}\right)$,
Eq. \eqref{eq:BDaggerAlternative} yields the identity
\begin{align}
\boldsymbol{\mathcal{B}}^{+}\left(\boldsymbol{x},\boldsymbol{r}\right)\boldsymbol{\chi}\left(\boldsymbol{r}\right) & =\boldsymbol{\mathcal{B}}^{+}\left(\boldsymbol{x},\boldsymbol{r}\right), & \boldsymbol{\mathcal{B}}^{+}\left(\boldsymbol{x},\boldsymbol{r}\right)\boldsymbol{\psi}\left(\boldsymbol{r}\right) & =\boldsymbol{0}.
\end{align}
Consequently, the control can be written as
\begin{align}
\boldsymbol{u}\left(\boldsymbol{r},t\right) & =\boldsymbol{\mathcal{B}}^{+}\left(\boldsymbol{x}\left(\boldsymbol{r},t\right),\boldsymbol{r}\right)\boldsymbol{\chi}\left(\boldsymbol{r}\right)\left(\partial_{t}\boldsymbol{x}\left(\boldsymbol{r},t\right)-\boldsymbol{\mathcal{D}}\triangle\boldsymbol{x}\left(\boldsymbol{r},t\right)-\boldsymbol{R}\left(\boldsymbol{x}\left(\boldsymbol{r},t\right)\right)\right).
\end{align}
As could be expected intuitively, evaluating this expression at a
position $\boldsymbol{r}\in\Omega\setminus A$ outside the region
$A$ affected by control yields a vanishing control signal,
\begin{align}
\boldsymbol{u}\left(\boldsymbol{r},t\right) & =\boldsymbol{0},\,\boldsymbol{r}\in\Omega\setminus A.
\end{align}
Using expression \eqref{eq:RDSControl} for the control signal in
the controlled reaction-diffusion system, Eq. \eqref{eq:ControlledRDS},
yields 
\begin{align}
\partial_{t}\boldsymbol{x}\left(\boldsymbol{r},t\right) & =\boldsymbol{\mathcal{Q}}\left(\boldsymbol{x}\left(\boldsymbol{r},t\right),\boldsymbol{r}\right)\left(\boldsymbol{\mathcal{D}}\triangle\boldsymbol{x}\left(\boldsymbol{r},t\right)+\boldsymbol{R}\left(\boldsymbol{x}\left(\boldsymbol{r},t\right)\right)\right)\nonumber \\
 & +\boldsymbol{\mathcal{P}}\left(\boldsymbol{x}\left(\boldsymbol{r},t\right),\boldsymbol{r}\right)\partial_{t}\boldsymbol{x}\left(\boldsymbol{r},t\right),
\end{align}
or
\begin{align}
\boldsymbol{0} & =\boldsymbol{\mathcal{Q}}\left(\boldsymbol{x}\left(\boldsymbol{r},t\right),\boldsymbol{r}\right)\left(\partial_{t}\boldsymbol{x}\left(\boldsymbol{r},t\right)-\boldsymbol{\mathcal{D}}\triangle\boldsymbol{x}\left(\boldsymbol{r},t\right)-\boldsymbol{R}\left(\boldsymbol{x}\left(\boldsymbol{r},t\right)\right)\right).\label{eq:RDSConstraintEquation}
\end{align}
Equation \eqref{eq:RDSConstraintEquation} is the constraint equation
for reaction-diffusion systems. The two complementary $n\times n$
projectors $\boldsymbol{\mathcal{P}}$ and $\boldsymbol{\mathcal{Q}}$
are given by
\begin{align}
\boldsymbol{\mathcal{P}}\left(\boldsymbol{x},\boldsymbol{r}\right) & =\boldsymbol{\mathcal{B}}\left(\boldsymbol{x},\boldsymbol{r}\right)\boldsymbol{\mathcal{B}}^{+}\left(\boldsymbol{x},\boldsymbol{r}\right)=\boldsymbol{\chi}\left(\boldsymbol{r}\right)\boldsymbol{\mathcal{B}}\left(\boldsymbol{x}\right)\boldsymbol{\mathcal{B}}^{+}\left(\boldsymbol{x},\boldsymbol{r}\right)\boldsymbol{\chi}\left(\boldsymbol{r}\right),\\
\boldsymbol{\mathcal{Q}}\left(\boldsymbol{x},\boldsymbol{r}\right) & =\boldsymbol{1}-\boldsymbol{\mathcal{P}}\left(\boldsymbol{x},\boldsymbol{r}\right).
\end{align}
For general spatio-temporal control systems, the projectors $\boldsymbol{\mathcal{P}}$
and $\boldsymbol{\mathcal{Q}}$ do not only depend the state $\boldsymbol{x}$,
but also on the position $\boldsymbol{r}$. Acting with projectors
$\boldsymbol{\mathcal{\psi}}$ and $\boldsymbol{\mathcal{\chi}}$
on $\boldsymbol{\mathcal{P}}$ and $\boldsymbol{\mathcal{Q}}$ and
using also Eq. \eqref{eq:BDaggerAlternative} yields the following
relations, 
\begin{align}
\boldsymbol{\mathcal{\psi}}\left(\boldsymbol{r}\right)\boldsymbol{\mathcal{P}}\left(\boldsymbol{x},\boldsymbol{r}\right) & =\boldsymbol{0}, & \boldsymbol{\mathcal{\chi}}\left(\boldsymbol{r}\right)\boldsymbol{\mathcal{P}}\left(\boldsymbol{x},\boldsymbol{r}\right) & =\boldsymbol{\mathcal{P}}\left(\boldsymbol{x},\boldsymbol{r}\right),\\
\boldsymbol{\mathcal{\psi}}\left(\boldsymbol{r}\right)\boldsymbol{\mathcal{Q}}\left(\boldsymbol{x},\boldsymbol{r}\right) & =\boldsymbol{\mathcal{\psi}}\left(\boldsymbol{r}\right), & \boldsymbol{\mathcal{\chi}}\left(\boldsymbol{r}\right)\boldsymbol{\mathcal{Q}}\left(\boldsymbol{x},\boldsymbol{r}\right) & =\boldsymbol{\mathcal{\chi}}\left(\boldsymbol{r}\right)-\boldsymbol{\mathcal{P}}\left(\boldsymbol{x},\boldsymbol{r}\right)=\boldsymbol{\mathcal{Q}}\left(\boldsymbol{x},\boldsymbol{r}\right)-\boldsymbol{\mathcal{\psi}}\left(\boldsymbol{r}\right),\\
\boldsymbol{\mathcal{P}}\left(\boldsymbol{x},\boldsymbol{r}\right)\boldsymbol{\mathcal{\psi}}\left(\boldsymbol{r}\right) & =\boldsymbol{0}, & \boldsymbol{\mathcal{P}}\left(\boldsymbol{x},\boldsymbol{r}\right)\boldsymbol{\mathcal{\chi}}\left(\boldsymbol{r}\right) & =\boldsymbol{\mathcal{P}}\left(\boldsymbol{x},\boldsymbol{r}\right),\\
\boldsymbol{\mathcal{Q}}\left(\boldsymbol{x},\boldsymbol{r}\right)\boldsymbol{\mathcal{\psi}}\left(\boldsymbol{r}\right) & =\boldsymbol{\mathcal{\psi}}\left(\boldsymbol{r}\right), & \boldsymbol{\mathcal{Q}}\left(\boldsymbol{x},\boldsymbol{r}\right)\boldsymbol{\mathcal{\chi}}\left(\boldsymbol{r}\right) & =\boldsymbol{\mathcal{\chi}}\left(\boldsymbol{r}\right)-\boldsymbol{\mathcal{P}}\left(\boldsymbol{x},\boldsymbol{r}\right)=\boldsymbol{\mathcal{Q}}\left(\boldsymbol{x},\boldsymbol{r}\right)-\boldsymbol{\mathcal{\psi}}\left(\boldsymbol{r}\right).
\end{align}
The constraint equation \eqref{eq:RDSConstraintEquation} can be interpreted
as follows. Acting with $\boldsymbol{\mathcal{\psi}}\left(\boldsymbol{r}\right)$
from the left on the constraint equation \eqref{eq:RDSConstraintEquation}
yields
\begin{align}
\boldsymbol{0} & =\boldsymbol{\mathcal{\psi}}\left(\boldsymbol{r}\right)\left(\partial_{t}\boldsymbol{x}\left(\boldsymbol{r},t\right)-\boldsymbol{\mathcal{D}}\triangle\boldsymbol{x}\left(\boldsymbol{r},t\right)-\boldsymbol{R}\left(\boldsymbol{x}\left(\boldsymbol{r},t\right)\right)\right).
\end{align}
Thus, outside the spatial region $A$ affected by control, the state
$\boldsymbol{x}\left(\boldsymbol{r},t\right)$ satisfies the uncontrolled
reaction-diffusion equation. Acting with $\boldsymbol{\mathcal{\chi}}\left(\boldsymbol{r}\right)$
from the left on Eq. \eqref{eq:RDSConstraintEquation} yields an equation
which is equivalent to the constraint equation for dynamical systems,
\begin{align}
\boldsymbol{0} & =\boldsymbol{\mathcal{\chi}}\left(\boldsymbol{r}\right)\boldsymbol{\mathcal{Q}}\left(\boldsymbol{x},\boldsymbol{r}\right)\left(\partial_{t}\boldsymbol{x}\left(\boldsymbol{r},t\right)-\boldsymbol{\mathcal{D}}\triangle\boldsymbol{x}\left(\boldsymbol{r},t\right)-\boldsymbol{R}\left(\boldsymbol{x}\left(\boldsymbol{r},t\right)\right)\right).\label{eq:RDSConstraintEq2}
\end{align}
Inside the region $A$ affected by control, $p$ state components
determine the vector of control signals $\boldsymbol{u}\left(t\right)$
while the remaining $n-p$ components are fixed by Eq. \eqref{eq:RDSConstraintEq2}.

\section{\label{sec:RDSExactlyRealizableDistributions}Exactly realizable
distributions}

The spatio-temporal analogues of desired trajectories in dynamical
systems are called desired distributions. Exactly realizable distributions
are desired distributions for which a control signal can be found
such that the state $\boldsymbol{x}\left(\boldsymbol{r},t\right)$
equals the desired distribution $\boldsymbol{x}_{d}\left(\boldsymbol{r},t\right)$
everywhere and for all times,
\begin{align}
\boldsymbol{x}\left(\boldsymbol{r},t\right) & =\boldsymbol{x}_{d}\left(\boldsymbol{r},t\right).
\end{align}
For a desired distribution $\boldsymbol{x}_{d}\left(\boldsymbol{r},t\right)$
to be exactly realizable, it has to satisfy the constraint equation
\begin{align}
\boldsymbol{0} & =\boldsymbol{\mathcal{Q}}\left(\boldsymbol{x}_{d}\left(\boldsymbol{r},t\right),\boldsymbol{r}\right)\left(\partial_{t}\boldsymbol{x}_{d}\left(\boldsymbol{r},t\right)-\boldsymbol{\mathcal{D}}\triangle\boldsymbol{x}_{d}\left(\boldsymbol{r},t\right)-\boldsymbol{R}\left(\boldsymbol{x}_{d}\left(\boldsymbol{r},t\right)\right)\right).\label{eq:RDSConstraintEq3}
\end{align}
Furthermore, the desired distribution $\boldsymbol{x}_{d}\left(\boldsymbol{r},t\right)$
must comply with the initial and boundary conditions for the state,
\begin{align}
\boldsymbol{x}_{d}\left(\boldsymbol{r},t_{0}\right) & =\boldsymbol{x}\left(\boldsymbol{r},t_{0}\right), & \boldsymbol{n}^{T}\left(\boldsymbol{r}\right)\left(\boldsymbol{\mathcal{D}}\nabla\boldsymbol{x}_{d}\left(\boldsymbol{r},t\right)\right) & =\boldsymbol{0},\,\boldsymbol{r}\in\Gamma.
\end{align}
The control signal enforcing the exactly realizable distribution $\boldsymbol{x}_{d}\left(\boldsymbol{r},t\right)$
is given by 
\begin{align}
\boldsymbol{u}\left(\boldsymbol{r},t\right) & =\boldsymbol{\mathcal{B}}^{+}\left(\boldsymbol{x}_{d}\left(\boldsymbol{r},t\right),\boldsymbol{r}\right)\boldsymbol{\chi}\left(\boldsymbol{r}\right)\left(\partial_{t}\boldsymbol{x}_{d}\left(\boldsymbol{r},t\right)-\boldsymbol{\mathcal{D}}\triangle\boldsymbol{x}_{d}\left(\boldsymbol{r},t\right)-\boldsymbol{R}\left(\boldsymbol{x}_{d}\left(\boldsymbol{r},t\right)\right)\right).\label{eq:RDSControlSignal}
\end{align}
The proof that these assumptions lead to the desired distribution
being an exact solution for the state, $\boldsymbol{x}\left(\boldsymbol{r},t\right)=\boldsymbol{x}_{d}\left(\boldsymbol{r},t\right)$,
is analogous to the proof for dynamical systems from Section \ref{sec:ExactlyRealizableTrajectories}.
In short, introducing $\boldsymbol{y}$ as
\begin{align}
\boldsymbol{x}\left(\boldsymbol{r},t\right) & =\boldsymbol{x}_{d}\left(\boldsymbol{r},t\right)+\boldsymbol{y}\left(\boldsymbol{r},t\right),
\end{align}
and using the control signal Eq. \eqref{eq:RDSControlSignal} in the
controlled state equation \eqref{eq:ControlledRDS} yields, together
with the constraint equation \eqref{eq:RDSConstraintEq3} and after
linearization in $\boldsymbol{y}$, a linear homogeneous partial differential
equation for $\boldsymbol{y}$,
\begin{align}
\partial_{t}\boldsymbol{y}\left(\boldsymbol{r},t\right) & =\boldsymbol{\mathcal{D}}\triangle\boldsymbol{y}\left(\boldsymbol{r},t\right)+\left(\nabla\boldsymbol{R}\left(\boldsymbol{x}_{d}\left(\boldsymbol{r},t\right)\right)+\boldsymbol{\mathcal{T}}\left(\boldsymbol{x}_{d}\left(\boldsymbol{r},t\right),\boldsymbol{r}\right)\right)\boldsymbol{y}\left(\boldsymbol{r},t\right).\label{eq:LinearPDEy}
\end{align}
The $n\times n$ matrix $\boldsymbol{\mathcal{T}}\left(\boldsymbol{x},\boldsymbol{r}\right)$
is defined by 
\begin{align}
\boldsymbol{\mathcal{T}}\left(\boldsymbol{x},\boldsymbol{r}\right)\boldsymbol{y} & =\boldsymbol{\chi}\left(\boldsymbol{r}\right)\left(\nabla\boldsymbol{\mathcal{B}}\left(\boldsymbol{x}\right)\boldsymbol{y}\right)\boldsymbol{\mathcal{B}}^{+}\left(\boldsymbol{x},\boldsymbol{r}\right)\boldsymbol{\chi}\left(\boldsymbol{r}\right)\left(\partial_{t}\boldsymbol{x}-\boldsymbol{\mathcal{D}}\triangle\boldsymbol{x}-\boldsymbol{R}\left(\boldsymbol{x}\right)\right).
\end{align}
Equation \eqref{eq:LinearPDEy} is to be solved with the initial condition
\begin{align}
\boldsymbol{y}\left(\boldsymbol{r},t_{0}\right) & =\boldsymbol{x}_{0}\left(\boldsymbol{r}\right)-\boldsymbol{x}_{d}\left(\boldsymbol{r},t\right).
\end{align}
If the initial state $\boldsymbol{x}_{0}\left(\boldsymbol{r}\right)$
complies with the initial desired distribution, then $\boldsymbol{y}\left(\boldsymbol{r},t_{0}\right)=\boldsymbol{0}$
initially. Furthermore, if $\boldsymbol{x}\left(\boldsymbol{r},t\right)$
as well as $\boldsymbol{x}_{d}\left(\boldsymbol{r},t\right)$ satisfy
homogeneous Neumann boundary conditions, then
\begin{align}
\boldsymbol{n}^{T}\left(\boldsymbol{r}\right)\left(\boldsymbol{\mathcal{D}}\nabla\boldsymbol{y}\left(\boldsymbol{r},t\right)\right) & =\boldsymbol{0},\,\boldsymbol{r}\in\Gamma.
\end{align}
Consequently, $\boldsymbol{y}$ vanishes everywhere and for all times,
\begin{align}
\boldsymbol{y}\left(\boldsymbol{r},t\right) & \equiv\boldsymbol{0}.
\end{align}
Equation \eqref{eq:LinearPDEy} determines the stability of exactly
realizable trajectories against perturbations $\boldsymbol{y}\left(\boldsymbol{r},t_{0}\right)=\boldsymbol{y}_{0}$
of the initial conditions.

Similar as for dynamical systems, a linearizing assumption can be
introduced. First of all, the projectors $\boldsymbol{\mathcal{P}}$
and $\boldsymbol{\mathcal{Q}}$ must be independent of the state $\boldsymbol{x}$,
\begin{align}
\boldsymbol{\mathcal{Q}}\left(\boldsymbol{x},\boldsymbol{r}\right) & =\boldsymbol{\mathcal{Q}}\left(\boldsymbol{r}\right)=\text{const.}\label{eq:RDSLinearizingAssumption1}
\end{align}
Second, the nonlinearity $\boldsymbol{R}\left(\boldsymbol{x}\right)$
must satisfy the condition
\begin{align}
\boldsymbol{\mathcal{Q}}\left(\boldsymbol{r}\right)\boldsymbol{R}\left(\boldsymbol{x}\right) & =\boldsymbol{\mathcal{Q}}\left(\boldsymbol{r}\right)\boldsymbol{\mathcal{A}}\boldsymbol{x}+\boldsymbol{\mathcal{Q}}\left(\boldsymbol{r}\right)\boldsymbol{b},\label{eq:RDSLinearizingAssumption2}
\end{align}
such that the constraint equation \eqref{eq:RDSConstraintEq3} becomes
linear. In principle, the projector $\boldsymbol{\mathcal{Q}}\left(\boldsymbol{r}\right)$
may still depend on the position $\boldsymbol{r}$ to yield a linear
constraint equation. This might be useful if the nonlinearity $\boldsymbol{R}\left(\boldsymbol{x}\right)$
exhibits an explicit dependence on space $\boldsymbol{r}$.

Equation \eqref{eq:RDSConstraintEq3} becomes a linear PDE for $\boldsymbol{\mathcal{Q}}\left(\boldsymbol{r}\right)\boldsymbol{x}_{d}\left(\boldsymbol{r},t\right)$
with $\boldsymbol{\mathcal{P}}\left(\boldsymbol{r}\right)\boldsymbol{x}_{d}\left(\boldsymbol{r},t\right)$
serving as an inhomogeneity, 
\begin{align}
\boldsymbol{\mathcal{Q}}\left(\boldsymbol{r}\right)\partial_{t}\boldsymbol{x}_{d}\left(\boldsymbol{r},t\right) & =\boldsymbol{\mathcal{Q}}\left(\boldsymbol{r}\right)\boldsymbol{\mathcal{D}}\triangle\boldsymbol{x}_{d}\left(\boldsymbol{r},t\right)+\boldsymbol{\mathcal{Q}}\left(\boldsymbol{r}\right)\boldsymbol{\mathcal{A}}\boldsymbol{x}_{d}\left(\boldsymbol{r},t\right)+\boldsymbol{\mathcal{Q}}\left(\boldsymbol{r}\right)\boldsymbol{b},
\end{align}
or, inserting $\boldsymbol{1}=\boldsymbol{\mathcal{P}}\left(\boldsymbol{r}\right)+\boldsymbol{\mathcal{Q}}\left(\boldsymbol{r}\right)$
between $\boldsymbol{\mathcal{A}}$ and $\boldsymbol{x}_{d}$, 
\begin{align}
 & \partial_{t}\left(\boldsymbol{\mathcal{Q}}\left(\boldsymbol{r}\right)\boldsymbol{x}_{d}\left(\boldsymbol{r},t\right)\right)-\boldsymbol{\mathcal{Q}}\left(\boldsymbol{r}\right)\boldsymbol{\mathcal{D}}\triangle\left(\boldsymbol{\mathcal{Q}}\left(\boldsymbol{r}\right)\boldsymbol{x}_{d}\left(\boldsymbol{r},t\right)\right)-\boldsymbol{\mathcal{Q}}\left(\boldsymbol{r}\right)\boldsymbol{\mathcal{A}}\boldsymbol{\mathcal{Q}}\left(\boldsymbol{r}\right)\boldsymbol{x}_{d}\left(\boldsymbol{r},t\right)\nonumber \\
= & \boldsymbol{\mathcal{Q}}\left(\boldsymbol{r}\right)\boldsymbol{\mathcal{D}}\triangle\left(\boldsymbol{\mathcal{P}}\left(\boldsymbol{r}\right)\boldsymbol{x}_{d}\left(\boldsymbol{r},t\right)\right)+\boldsymbol{\mathcal{Q}}\left(\boldsymbol{r}\right)\boldsymbol{\mathcal{A}}\boldsymbol{\mathcal{P}}\left(\boldsymbol{r}\right)\boldsymbol{x}_{d}\left(\boldsymbol{r},t\right)+\boldsymbol{\mathcal{Q}}\left(\boldsymbol{r}\right)\boldsymbol{b}.\label{eq:LinearizingAssumptionPDE}
\end{align}
Being a linear partial differential equation for $\boldsymbol{\mathcal{Q}}\left(\boldsymbol{r}\right)\boldsymbol{x}_{d}\left(\boldsymbol{r},t\right)$
with inhomogeneity on the right hand side, Eq. \eqref{eq:LinearizingAssumptionPDE}
can formally be solved with the help of Green's functions. The explicit
form of the Green's function depends on the form and dimension of
the spatial domain.

Having identified a linearizing assumption, the next step would be
to discuss concepts of controllability for linear PDEs, as e.g. the
linear diffusion equation, and apply these concepts to Eq. \eqref{eq:LinearizingAssumptionPDE}.
However, in contrast to dynamical systems, no condition for controllability
in terms of a rank condition for a controllability matrix can be formulated.
The reason is that PDEs are essentially dynamical systems with an
infinite-dimensional state space. The Cayley-Hamilton theorem cannot
be applied to truncate the exponential of a linear operator after
a finite number of terms, see Section \ref{sub:DerivationOfKalmanRank}.
We omit a discussion of controllability and present position control
of traveling waves as an application of exactly realizable distributions.

\section{\label{sec:PositionControlOfTravelingWaves}Position control of traveling
waves}

For simplicity, the projectors $\boldsymbol{\mathcal{Q}}$ and $\boldsymbol{\mathcal{P}}$
are assumed to be constant in space and independent of the state $\boldsymbol{x}$,
\begin{align}
\boldsymbol{\mathcal{Q}}\left(\boldsymbol{x},\boldsymbol{r}\right) & =\boldsymbol{\mathcal{Q}}=\text{const.}, & \boldsymbol{\mathcal{P}}\left(\boldsymbol{x},\boldsymbol{r}\right) & =\boldsymbol{\mathcal{P}}=\text{const.}.
\end{align}
Note that this implies that the control acts everywhere in position
space,
\begin{align}
\boldsymbol{\chi}\left(\boldsymbol{r}\right) & =\boldsymbol{1}.
\end{align}
Consider the uncontrolled reaction-diffusion system in an unbounded
domain $\Omega=\mathbb{R}^{N}$, 
\begin{align}
\partial_{t}\boldsymbol{x}\left(\boldsymbol{r},t\right) & =\boldsymbol{\mathcal{D}}\triangle\boldsymbol{x}\left(\boldsymbol{r},t\right)+\boldsymbol{R}\left(\boldsymbol{x}\left(\boldsymbol{r},t\right)\right).
\end{align}
Many reaction-diffusion systems exhibit plane traveling wave solutions
propagating with constant velocity $c$ in a constant direction $\boldsymbol{\hat{c}}$,
$\left|\boldsymbol{\hat{c}}\right|=1$. A traveling wave is characterized
by a wave profile $\boldsymbol{X}_{c}$ depending only on a single
coordinate as 
\begin{align}
\boldsymbol{x}\left(\boldsymbol{r},t\right) & =\boldsymbol{X}_{c}\left(\boldsymbol{\hat{c}}^{T}\boldsymbol{r}-ct\right)=\boldsymbol{X}_{c}\left(\sum_{i=1}^{n}\hat{c}_{i}r_{i}-ct\right).\label{eq:TravelingWave}
\end{align}
In a frame of reference $\xi=\boldsymbol{\hat{c}}^{T}\boldsymbol{r}-ct$
comoving with the traveling wave, the wave profile $\boldsymbol{X}_{c}$
appears stationary and satisfies the profile equation
\begin{align}
\boldsymbol{0} & =\boldsymbol{\mathcal{D}}\boldsymbol{X}_{c}''\left(\xi\right)+c\boldsymbol{X}_{c}'\left(\xi\right)+\boldsymbol{R}\left(\boldsymbol{X}_{c}\left(\xi\right)\right).\label{eq:ProfileEquation}
\end{align}
The ODE for the wave profile, Eq. \eqref{eq:ProfileEquation}, can
exhibit one or more homogeneous steady states. Typically, for $\xi\rightarrow\pm\infty$,
the wave profile $\boldsymbol{X}_{c}$ approaches either two different
steady states or the same steady state. This fact can be used to classify
traveling wave profiles. Front profiles connect different steady states
for $\xi\rightarrow\pm\infty$ and are found to be heteroclinic orbits
of Eq. \eqref{eq:ProfileEquation}, while pulse profiles join the
same steady state and are found to be homoclinic orbits. Pulse profiles
are naturally localized and usually every component exhibits one or
several extrema. Fronts are not localized but typically exhibit a
narrow region where the transition from one to the other steady state
occurs. Therefore, all traveling wave solutions are localized in the
sense that the derivatives of any order $m\geq1$ of the wave profile
$\boldsymbol{X}_{c}\left(\xi\right)$ with respect to the traveling
wave coordinate $\xi$ decays to zero,
\begin{align}
\lim_{\xi\rightarrow\pm\infty}\partial_{\xi}^{m}\boldsymbol{X}_{c}\left(\xi\right) & =\boldsymbol{0}.\label{eq:LocalizationOfTravelingWaves}
\end{align}
Note that Eq. \eqref{eq:LocalizationOfTravelingWaves} is not a boundary
condition for $\boldsymbol{X}_{c}$ but characterizes the solution
to Eq. \eqref{eq:ProfileEquation}.

We assume that before control is switched on at time $t=t_{0}$, the
traveling wave moves unperturbed. Thus, the initial condition for
the controlled reaction-diffusion system Eq. \eqref{eq:ControlledRDS}
is
\begin{align}
\boldsymbol{x}\left(\boldsymbol{r},t_{0}\right) & =\boldsymbol{X}_{c}\left(\boldsymbol{\hat{c}}^{T}\boldsymbol{r}-ct_{0}\right).\label{eq:PositionControlInitCond}
\end{align}
The idea of position control is to choose the desired distribution
$\boldsymbol{x}_{d}\left(\boldsymbol{r},t\right)$ in form of a traveling
wave profile $\boldsymbol{X}_{c}$ shifted according to a protocol
of motion $\phi\left(t\right)$. The function $\phi\left(t\right)$
encodes the desired position over time of the controlled traveling
wave along the spatial direction $\boldsymbol{\hat{c}}$. The position
of a traveling wave is defined by a distinguishing point of the wave
profile. For pulse solutions, the extremum of a certain component
of the wave profile defines its position. The position of a front
solution is defined by a characteristic point in the transition region
as e.g. the point of the steepest slope. A problem arises because
for an exactly realizable distribution, only $p$ out of all $n$
components of $\boldsymbol{x}_{d}\left(\boldsymbol{r},t\right)$ can
be prescribed, while the remaining $n-p$ components have to satisfy
the constraint equation \eqref{eq:RDSConstraintEq3}. Here, the convention
is that the part $\boldsymbol{\mathcal{P}}\boldsymbol{x}_{d}\left(\boldsymbol{r},t\right)$
is the traveling wave profile $\boldsymbol{X}_{c}$ shifted according
to the protocol of motion $\phi\left(t\right)$,
\begin{align}
\boldsymbol{\mathcal{P}}\boldsymbol{x}_{d}\left(\boldsymbol{r},t\right) & =\boldsymbol{\mathcal{P}}\boldsymbol{X}_{c}\left(\boldsymbol{\hat{c}}^{T}\boldsymbol{r}-\phi\left(t\right)\right).
\end{align}
The remaining part $\boldsymbol{\mathcal{Q}}\boldsymbol{x}_{d}\left(\boldsymbol{r},t\right)$
has to satisfy the constraint equation \eqref{eq:RDSConstraintEq3},
\begin{align}
\partial_{t}\left(\boldsymbol{\mathcal{Q}}\boldsymbol{x}_{d}\left(\boldsymbol{r},t\right)\right) & =\boldsymbol{\mathcal{Q}}\boldsymbol{\mathcal{D}}\triangle\left(\boldsymbol{\mathcal{Q}}\boldsymbol{x}_{d}\left(\boldsymbol{r},t\right)\right)+\boldsymbol{\mathcal{Q}}\boldsymbol{R}\left(\boldsymbol{x}_{d}\left(\boldsymbol{r},t\right)\right)+\boldsymbol{\mathcal{Q}}\boldsymbol{\mathcal{D}}\boldsymbol{\mathcal{P}}\boldsymbol{X}_{c}''\left(\boldsymbol{\hat{c}}^{T}\boldsymbol{r}-\phi\left(t\right)\right),
\end{align}
with initial condition
\begin{align}
\boldsymbol{\mathcal{Q}}\boldsymbol{x}_{d}\left(\boldsymbol{r},t_{0}\right) & =\boldsymbol{\mathcal{Q}}\boldsymbol{X}_{c}\left(\boldsymbol{\hat{c}}^{T}\boldsymbol{r}-ct_{0}\right).
\end{align}
From the initial condition Eq. \eqref{eq:PositionControlInitCond}
follows an initial condition for $\phi$ as
\begin{align}
\phi\left(t_{0}\right) & =ct_{0}.
\end{align}
The profile equation \eqref{eq:ProfileEquation} for $\boldsymbol{X}_{c}$
is exploited to obtain
\begin{align}
\partial_{t}\left(\boldsymbol{\mathcal{Q}}\boldsymbol{x}_{d}\left(\boldsymbol{r},t\right)\right) & =\boldsymbol{\mathcal{Q}}\boldsymbol{\mathcal{D}}\left(\triangle\boldsymbol{\mathcal{Q}}\boldsymbol{x}_{d}\left(\boldsymbol{r},t\right)-\boldsymbol{\mathcal{Q}}\boldsymbol{X}_{c}''\left(\boldsymbol{\hat{c}}^{T}\boldsymbol{r}-\phi\left(t\right)\right)\right)\nonumber \\
 & +\boldsymbol{\mathcal{Q}}\left(\boldsymbol{R}\left(\boldsymbol{x}_{d}\left(\boldsymbol{r},t\right)\right)-\boldsymbol{R}\left(\boldsymbol{X}_{c}\left(\boldsymbol{\hat{c}}^{T}\boldsymbol{r}-\phi\left(t\right)\right)\right)\right)-c\boldsymbol{\mathcal{Q}}\boldsymbol{X}_{c}'\left(\boldsymbol{\hat{c}}^{T}\boldsymbol{r}-\phi\left(t\right)\right).\label{eq:RDSConstEqTravelingWaves}
\end{align}
Let $\boldsymbol{\mathcal{Q}}\boldsymbol{y}_{d}\left(\boldsymbol{r},t\right)$
be defined by 
\begin{align}
\boldsymbol{x}_{d}\left(\boldsymbol{r},t\right) & =\boldsymbol{X}_{c}\left(\boldsymbol{\hat{c}}^{T}\boldsymbol{r}-\phi\left(t\right)\right)+\boldsymbol{\mathcal{Q}}\boldsymbol{y}_{d}\left(\boldsymbol{r},t\right).\label{eq:xdInTermsOfyd}
\end{align}
The constraint equation \eqref{eq:RDSConstEqTravelingWaves} can be
written as a PDE for $\boldsymbol{\mathcal{Q}}\boldsymbol{y}_{d}\left(\boldsymbol{r},t\right)$,
\begin{align}
\partial_{t}\left(\boldsymbol{\mathcal{Q}}\boldsymbol{y}_{d}\left(\boldsymbol{r},t\right)\right) & =\boldsymbol{\mathcal{Q}}\boldsymbol{\mathcal{D}}\triangle\boldsymbol{\mathcal{Q}}\boldsymbol{y}_{d}\left(\boldsymbol{r},t\right)+\boldsymbol{\mathcal{Q}}\boldsymbol{R}\left(\boldsymbol{X}_{c}\left(\boldsymbol{\hat{c}}^{T}\boldsymbol{r}-\phi\left(t\right)\right)+\boldsymbol{\mathcal{Q}}\boldsymbol{y}_{d}\left(\boldsymbol{r},t\right)\right)\nonumber \\
 & +\left(\dot{\phi}\left(t\right)-c\right)\boldsymbol{\mathcal{Q}}\boldsymbol{X}_{c}'\left(\boldsymbol{\hat{c}}^{T}\boldsymbol{r}-\phi\left(t\right)\right)-\boldsymbol{\mathcal{Q}}\boldsymbol{R}\left(\boldsymbol{X}_{c}\left(\boldsymbol{\hat{c}}^{T}\boldsymbol{r}-\phi\left(t\right)\right)\right).\label{eq:TravelingWaveConstraintEquation}
\end{align}
The next step is the determination of the control signal $\boldsymbol{u}\left(t\right)$.
Due to $\boldsymbol{\mathcal{B}}^{+}\left(\boldsymbol{x}\right)\boldsymbol{\mathcal{P}}=\boldsymbol{\mathcal{B}}^{+}\left(\boldsymbol{x}\right)$,
$\boldsymbol{\mathcal{P}}$ being constant in time and space, and
Eq. \eqref{eq:xdInTermsOfyd}, the control signal Eq. \eqref{eq:RDSControlSignal}
can be cast in the form
\begin{align}
\boldsymbol{u}\left(\boldsymbol{r},t\right) & =-\boldsymbol{\mathcal{B}}^{+}\left(\boldsymbol{X}_{c}\left(\boldsymbol{\hat{c}}^{T}\boldsymbol{r}-\phi\left(t\right)\right)+\boldsymbol{\mathcal{Q}}\boldsymbol{y}_{d}\left(\boldsymbol{r},t\right)\right)\left(\boldsymbol{\mathcal{D}}\boldsymbol{X}_{c}''\left(\boldsymbol{\hat{c}}^{T}\boldsymbol{r}-\phi\left(t\right)\right)\right.\nonumber \\
 & \left.+\boldsymbol{\mathcal{D}}\triangle\boldsymbol{\mathcal{Q}}\boldsymbol{y}_{d}\left(\boldsymbol{r},t\right)+\dot{\phi}\left(t\right)\boldsymbol{X}_{c}'\left(\boldsymbol{\hat{c}}^{T}\boldsymbol{r}-\phi\left(t\right)\right)\right.\nonumber \\
 & \left.+\boldsymbol{R}\left(\boldsymbol{X}_{c}\left(\boldsymbol{\hat{c}}^{T}\boldsymbol{r}-\phi\left(t\right)\right)+\boldsymbol{\mathcal{Q}}\boldsymbol{y}_{d}\left(\boldsymbol{r},t\right)\right)\right).
\end{align}
Exploiting again the profile equation \eqref{eq:ProfileEquation},
the last expression becomes
\begin{align}
\boldsymbol{u}\left(\boldsymbol{r},t\right) & =-\boldsymbol{\mathcal{B}}^{+}\left(\boldsymbol{X}_{c}\left(\boldsymbol{\hat{c}}^{T}\boldsymbol{r}-\phi\left(t\right)\right)+\boldsymbol{\mathcal{Q}}\boldsymbol{y}_{d}\left(\boldsymbol{r},t\right)\right)\left(\left(\dot{\phi}\left(t\right)-c\right)\boldsymbol{X}_{c}'\left(\boldsymbol{\hat{c}}^{T}\boldsymbol{r}-\phi\left(t\right)\right)\right)\nonumber \\
 & \left.+\boldsymbol{R}\left(\boldsymbol{X}_{c}\left(\boldsymbol{\hat{c}}^{T}\boldsymbol{r}-\phi\left(t\right)\right)+\boldsymbol{\mathcal{Q}}\boldsymbol{y}_{d}\left(\boldsymbol{r},t\right)\right)-\boldsymbol{R}\left(\boldsymbol{X}_{c}\left(\boldsymbol{\hat{c}}^{T}\boldsymbol{r}-\phi\left(t\right)\right)\right)\right.\nonumber \\
 & \left.+\boldsymbol{\mathcal{D}}\triangle\boldsymbol{\mathcal{Q}}\boldsymbol{y}_{d}\left(\boldsymbol{r},t\right)\right).\label{eq:TravelingWaveControl}
\end{align}
Equation \eqref{eq:TravelingWaveControl} for the control signal together
with the constraint equation in the form of Eq. \eqref{eq:TravelingWaveConstraintEquation}
is the starting point for the position control of traveling waves
in general reaction-diffusion systems. Several special cases leading
to simpler expressions can be identified.

An invertible coupling matrix yields $\boldsymbol{\mathcal{Q}}=\boldsymbol{0}$
and $\boldsymbol{\mathcal{B}}^{+}\left(\boldsymbol{x}\right)=\boldsymbol{\mathcal{B}}^{-1}\left(\boldsymbol{x}\right)$.
The constraint equation \eqref{eq:TravelingWaveConstraintEquation}
is trivially satisfied, and the control simplifies to 
\begin{align}
\boldsymbol{u}\left(\boldsymbol{r},t\right) & =\left(c-\dot{\phi}\left(t\right)\right)\boldsymbol{\mathcal{B}}^{-1}\left(\boldsymbol{X}_{c}\left(\boldsymbol{\hat{c}}^{T}\boldsymbol{r}-\phi\left(t\right)\right)\right)\boldsymbol{X}_{c}'\left(\boldsymbol{\hat{c}}^{T}\boldsymbol{r}-\phi\left(t\right)\right).
\end{align}
In this case, the control can expressed solely in terms of the traveling
wave profile $\boldsymbol{X}_{c}$ and its velocity $c$. Any reference
to the nonlinearity $\boldsymbol{R}\left(\boldsymbol{x}\right)$ vanishes.
This can be useful if $\boldsymbol{R}\left(\boldsymbol{x}\right)$
is only approximately known but the wave profile $\boldsymbol{X}_{c}$
and its velocity $c$ can be measured with sufficient accuracy in
experiments. However, the assumption of an equal number of control
signals and state components is restrictive and valid only for a limited
number of systems. Nevertheless, it is satisfied for all single-component
reaction-diffusion systems with a scalar distributed control signal
$u\left(\boldsymbol{r},t\right)$. As an example, we discuss position
control of traveling fronts in the Schlögl model.

\begin{example}[Position control of fronts in the Schl\"ogl model]\label{ex:PositionControlSchoeglModel}

Consider an autocatalytic chemical reaction mechanism proposed by
Schlögl \cite{schlogl1972crm}
\begin{align}
A_{1}+2X & \overset{k_{1}^{+}}{\underset{k_{1}^{-}}{\rightleftharpoons}}3X, & X & \overset{k_{2}^{+}}{\underset{k_{2}^{-}}{\rightleftharpoons}}A_{2}.
\end{align}
Under the assumption that the concentrations $a_{1/2}=\left[A_{1/2}\right]$
of the chemical species $A_{1/2}$ are kept constant in space and
time, a nonlinearity $R\left(x\right)$ in form of a cubic polynomial
\begin{align}
R\left(x\right) & =k_{1}^{+}a_{1}x{}^{2}-k_{1}^{-}x{}^{3}-k_{2}^{+}x+k_{2}^{-}a_{2}\nonumber \\
 & =-k\left(x-x_{0}\right)\left(x-x_{1}\right)\left(x-x_{2}\right)
\end{align}
dictates the time evolution of the concentration $x=\left[X\right]$.
For a certain range of parameters, $R\left(x\right)$ possesses three
real positive roots $0<x_{0}<x_{1}<x_{2}$. In one spatial dimension
$r$, the uncontrolled reaction-diffusion system known as the Schlögl
model becomes
\begin{align}
\partial_{t}x\left(r,t\right) & =D\partial_{r}^{2}x\left(r,t\right)+R\left(x\left(r,t\right)\right).\label{eq:SchloeglModelUncontrolled}
\end{align}
Although the Schlögl model is introduced in the context of chemical
reactions, a scalar reaction-diffusion equation of the form \eqref{eq:SchloeglModelUncontrolled}
with cubic nonlinearity $R$ can be seen as a paradigmatic model for
a bistable medium. Such models have found widespread application far
beyond chemical reactions. An important example is the phase field,
which is used to model phenomena as diverse as cell motility \cite{lober2015collisions},
free boundary problems in fluid mechanics \cite{anderson1998diffuse},
and solidification \cite{boettinger2002phase}. Initially, Eq. \eqref{eq:SchloeglModelUncontrolled}
has been discussed in 1938 by Zeldovich and Frank-Kamenetsky in connection
with flame propagation \cite{zeldovich1938theory}.

The roots $x_{0},\,x_{1}$, and $x_{2}$ are homogeneous steady states
of the system, with the upper ($x_{2}$) and lower ($x_{0}$) being
stable steady states while the root $x_{1}$ is unstable. The Schlögl
model exhibits a variety of traveling front solutions, 
\begin{align}
x\left(r,t\right) & =X_{c}\left(r-ct\right),
\end{align}
propagating with velocity $c$. The front profile $X_{c}$ satisfies
the profile equation with $\xi=x-ct$
\begin{align}
0 & =DX_{c}''\left(\xi\right)+cX_{c}'\left(\xi\right)+R\left(X_{c}\left(\xi\right)\right).
\end{align}
Front solutions connect the homogeneous steady states as $\lim_{\xi\rightarrow\pm\infty}$.
A stable traveling front solution connecting the lower and upper stable
states $x_{0}$ and $x_{2}$, respectively, is known analytically
and given by
\begin{align}
X_{c}\left(\xi\right) & =\dfrac{1}{2}\left(x_{0}+x_{2}\right)+\dfrac{1}{2}\left(x_{0}-x_{2}\right)\tanh\left(\dfrac{1}{2\sqrt{2}}\sqrt{\dfrac{k}{D}}\left(x_{2}-x_{0}\right)\xi\right),\\
c & =\sqrt{\dfrac{Dk}{2}}\left(x_{0}+x_{2}-2x_{1}\right).
\end{align}
Assuming that the concentrations $a_{1/2}$ can be controlled spatio-temporally
by the distributed control signal $u\left(r,t\right)$ amounts to
the substitution 
\begin{align}
a_{1/2} & \rightarrow a_{1/2}+u\left(x,t\right)
\end{align}
in Eq. \eqref{eq:SchloeglModelUncontrolled}. The controlled reaction-diffusion
system is
\begin{align}
\partial_{t}x\left(r,t\right) & =D\partial_{r}^{2}x\left(r,t\right)+R\left(x\left(r,t\right)\right)+B\left(x\left(r,t\right)\right)u\left(r,t\right).
\end{align}
Control by $a_{2}$ will be additive with constant coupling function
$B\left(x\right)=k_{2}^{-}$, while for control via $a_{1}$ the spatio-temporal
forcing couples multiplicatively to the RD kinetics and the coupling
function $B\left(x\right)=k_{1}^{+}x^{2}$ becomes state dependent.
See \cite{lober2014control} for a discussion of experimental realizations
of the controlled Schlögl model.

In the following, we assume control by parameter $a_{1}$ such that
the coupling function is
\begin{align}
B\left(x\right) & =k_{1}^{+}x^{2}.\label{eq:SchloeglMultCoupling}
\end{align}
In the context of chemical systems, $x$ is interpreted as a concentration
which only attains positive values, $x\geq0$. As long as $x>0$,
the coupling function as given by Eq. \eqref{eq:SchloeglMultCoupling}
is positive and $B\left(x\right)$ does not change its rank. Because
the Schlögl model is a single component reaction-diffusion system,
a single distributed control signal $u\left(r,t\right)$ is sufficient
to realize any desired distribution which complies with the initial
and boundary conditions of the system. With the desired distribution
$x_{d}\left(r,t\right)$ given in terms of the traveling wave solution
as
\begin{align}
x_{d}\left(r,t\right) & =X_{c}\left(r-\phi\left(t\right)\right),
\end{align}
the solution for the control signal becomes 
\begin{align}
u\left(r,t\right) & =\left(c-\dot{\phi}\left(t\right)\right)\dfrac{1}{B\left(X_{c}\left(r-\phi\left(t\right)\right)\right)}X_{c}'\left(r-\phi\left(t\right)\right).\label{eq:SchloeglMultControl}
\end{align}
The protocol of motion $\phi\left(t\right)$ is chosen to move the
front back and forth sinusoidally as
\begin{align}
\phi\left(t\right) & =A_{0}+A\sin\left(2\pi t/T+A_{1}\right).
\end{align}
The control is applied starting at time $t=t_{0}$, upon which the
front moves unperturbed with velocity $c$. To achieve a smooth transition
of the position $\phi$ and velocity $\dot{\phi}$ across $t=t_{0}$,
the constants $A_{0}$ and $A_{1}$ are determined by the conditions
\begin{align}
\phi\left(t_{0}\right) & =\phi_{0}, & \dot{\phi}\left(t_{0}\right) & =c.
\end{align}
Figure \ref{fig:SchoeglFront} shows a snapshot of the controlled
front solution (black solid line) and the control signal as given
by Eq. \eqref{eq:SchloeglMultControl} (red dashed line).

\begin{minipage}{1.0\linewidth}
\begin{center}
\includegraphics[scale=0.6]{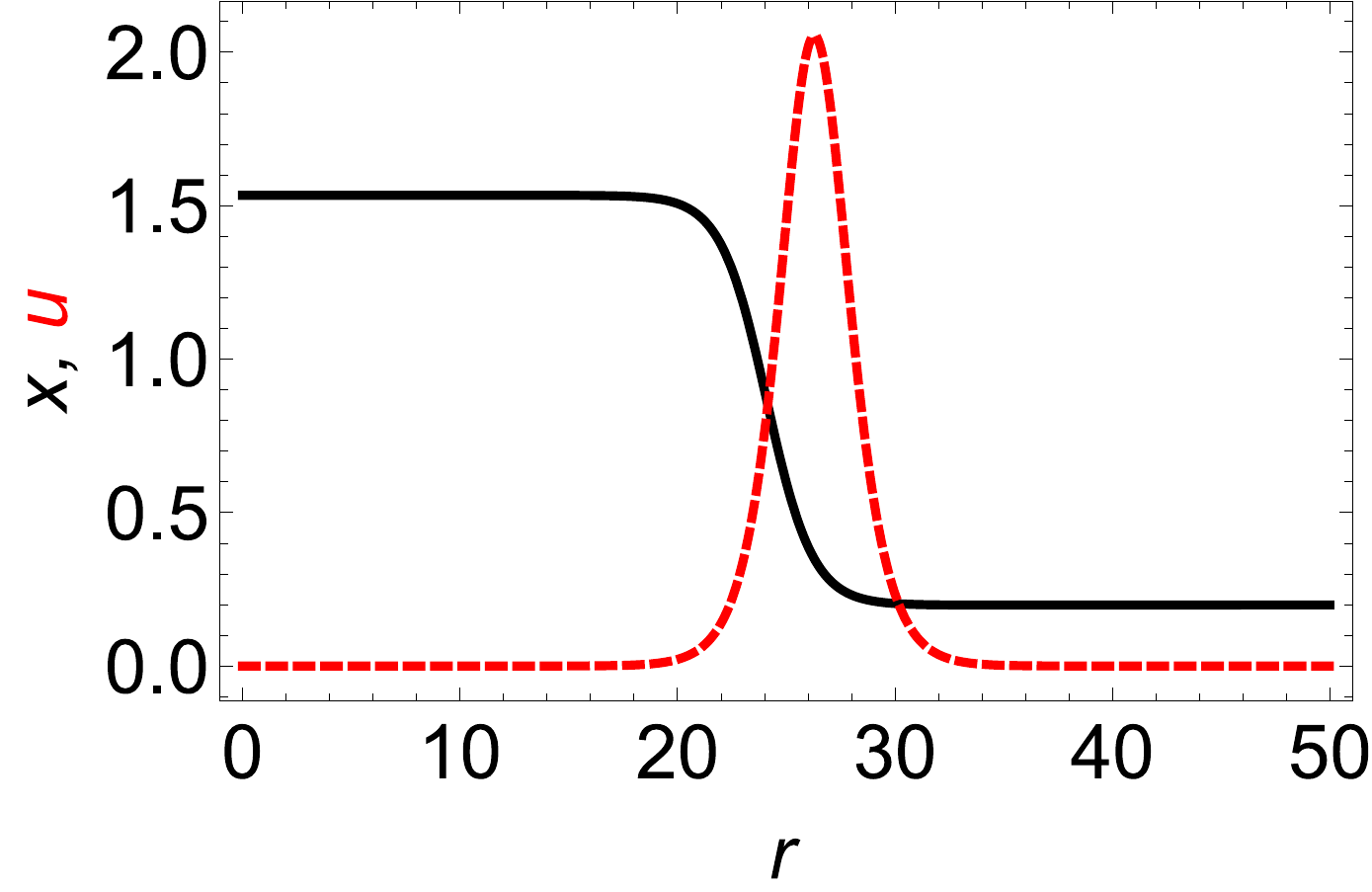}
\captionof{figure}[Snapshot of front solution and control in the Schl\"ogl model]{\label{fig:SchoeglFront}Snapshot of the controlled front solution $x\left( r,t \right) = X_c \left( r-ct \right)$ to the  Schl\"ogl model (black solid line) and distributed control signal (red solid line) as given by Eq. \eqref{eq:SchloeglMultControl}.}
\end{center}
\end{minipage}

To validate the performance of the control, the protocol of motion
$\phi\left(t\right)$ is compared with the position over time recorded
from numerical simulations of the controlled front. In numerical simulations,
the position of the front is defined as the point of the steepest
slope of the transition region. Figure \ref{fig:PositionControlSchoeglModel}
left demonstrates perfect agreement between prescribed (black solid
line) and recorded (red dashed line) position over time. The analogous
comparison for the velocity over time shown in Fig. \ref{fig:PositionControlSchoeglModel}
right reveals an overall perfect agreement but small deviations at
the points of maximum and minimum velocity. Such deviations can be
understood to arise from an underlying instability \cite{lober2014stability}.
This instability manifests as a finite, non-increasing shift between
the positions of control signal and controlled front.

\begin{minipage}{1.0\linewidth}
\begin{center}
\includegraphics[scale=0.475]{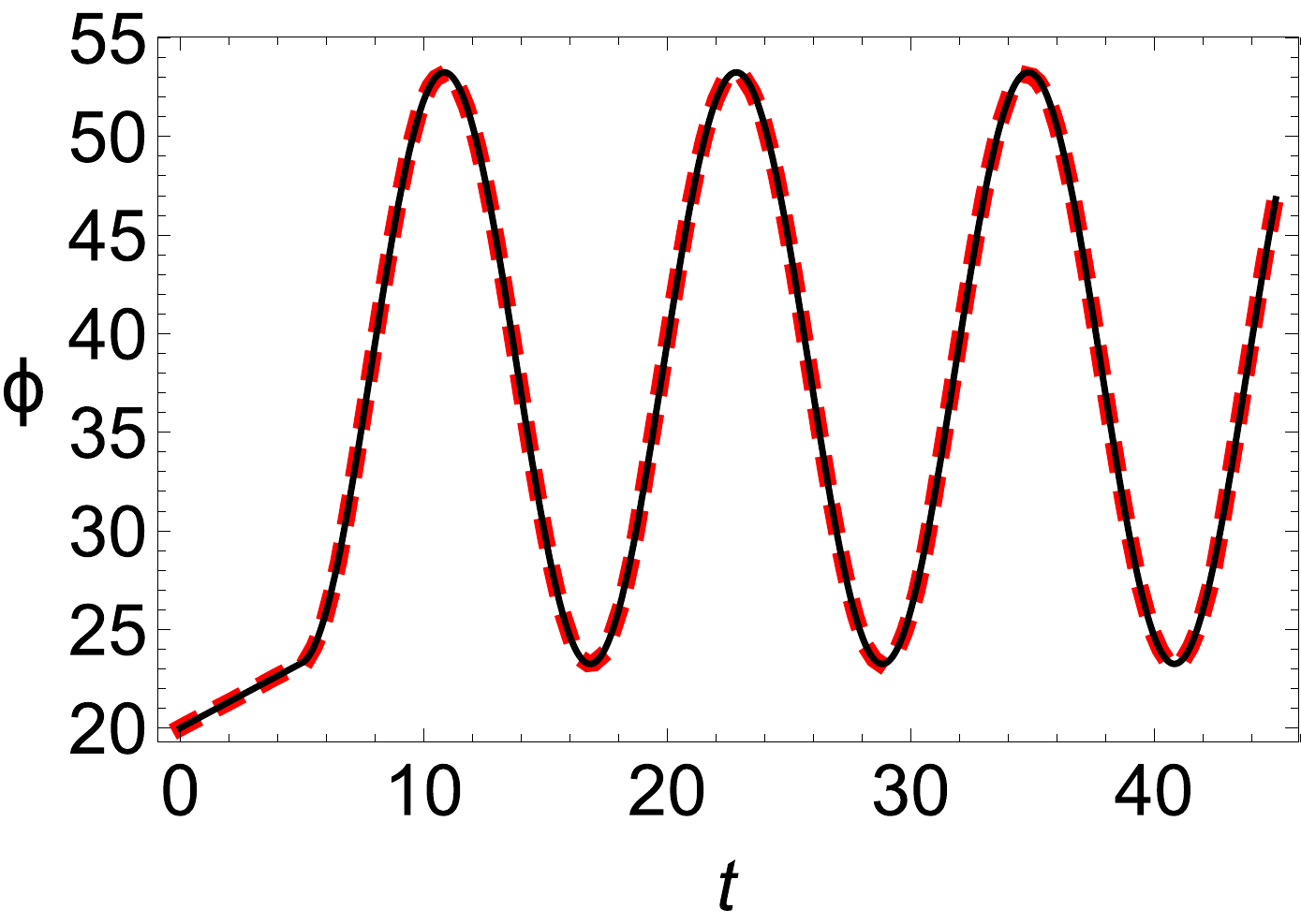}\hspace{0.2cm}\includegraphics[scale=0.475]{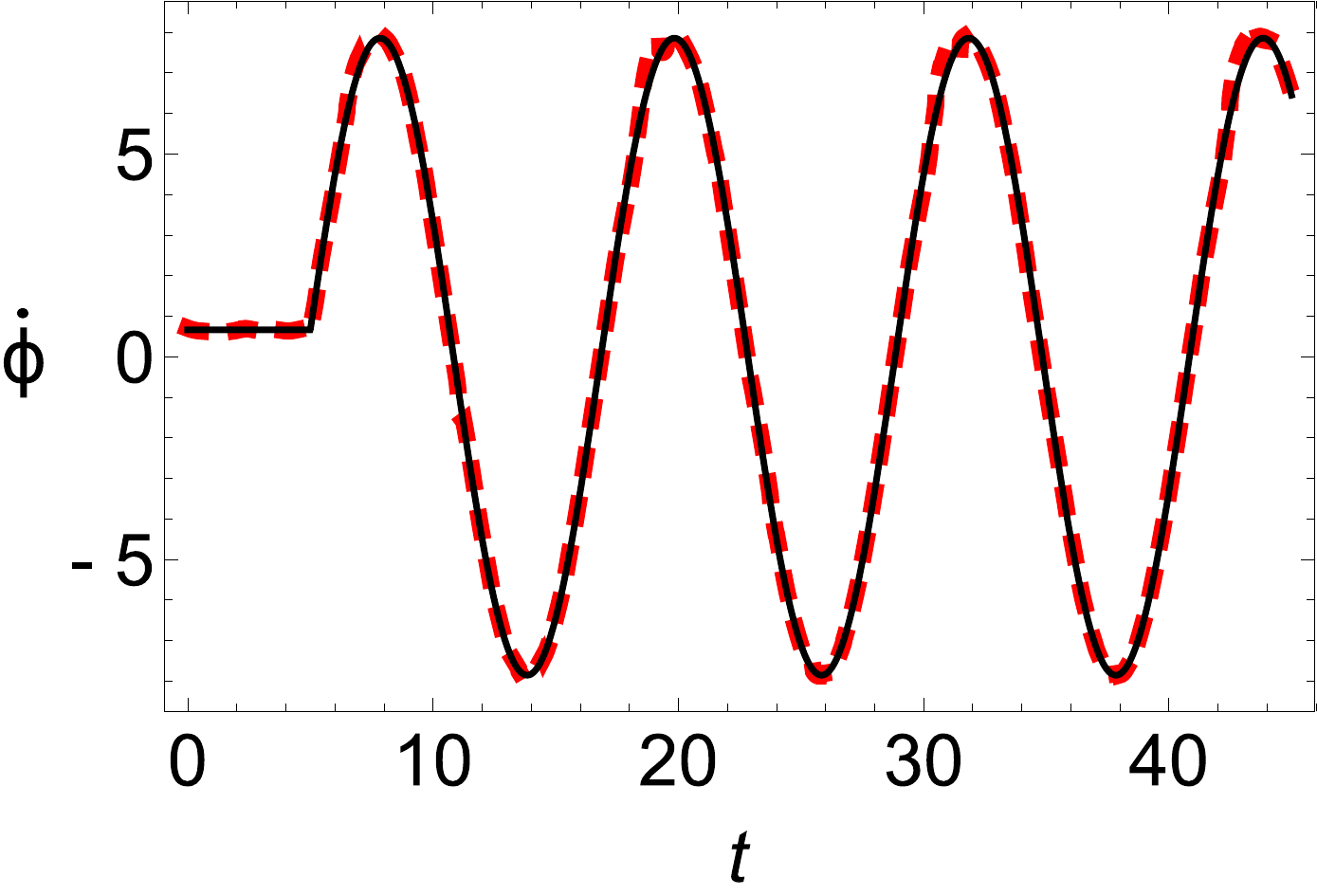}
\captionof{figure}[Position control of fronts in the Schl\"ogl model]{\label{fig:PositionControlSchoeglModel}Position control of fronts in the Schl\"ogl model. Left: position over time of the desired protocol of motion (red dashed line) and the actual position over time of the controlled front (black solid line). Right: Velocity over time. Agreement is nearly perfect in both cases.}
\end{center}
\end{minipage}

\end{example}

A second example of position control with a number of control signals
smaller than the number of state components is discussed in the following.

\begin{example}[Position control of traveling waves in the activator-controlled FHN model]\label{ex:FHNPositionControl}

Consider the one-dimensional spatial domain $0\leq r<L=150$ with
periodic boundary conditions. Apart from an additional diffusion term,
the model equations are the same as in Example \ref{ex:FHN1}, 
\begin{align}
\left(\begin{array}{c}
\partial_{t}x\left(r,t\right)\\
\partial_{t}y\left(r,t\right)
\end{array}\right) & =\left(\begin{array}{c}
D_{x}\partial_{r}^{2}x\left(r,t\right)\\
D_{y}\partial_{r}^{2}y\left(r,t\right)
\end{array}\right)\nonumber \\
 & +\left(\begin{array}{c}
a_{0}+a_{1}x\left(r,t\right)+a_{2}y\left(r,t\right)\\
R\left(x\left(r,t\right),y\left(r,t\right)\right)
\end{array}\right)+\left(\begin{array}{c}
0\\
1
\end{array}\right)u\left(r,t\right).\label{eq:RDSFHN}
\end{align}
The nonlinearity is linear in the inhibitor but nonlinear in the activator,
and $R$ is given by 
\begin{align}
R\left(x,y\right) & =R\left(y\right)-x.
\end{align}
The function $R\left(y\right)$ is a cubic polynomial of the form
\begin{align}
R\left(y\right) & =3y-y^{3}.
\end{align}
The parameter values are set to
\begin{align}
a_{0} & =0.429, & a_{1} & =0, & a_{2} & =0.33, & D_{y} & =1, & D_{x} & =0.3.
\end{align}
The traveling wave profile $\boldsymbol{X}_{c}\left(\xi\right)=\left(\begin{array}{cc}
X_{c}\left(\xi\right), & Y_{c}\left(\xi\right)\end{array}\right)^{T}$ satisfies, with $\xi=x-ct$,
\begin{align}
D_{x}X_{c}''\left(\xi\right)+cX_{c}'\left(\xi\right)+a_{0}+a_{1}X_{c}\left(\xi\right)+a_{2}Y_{c}\left(\xi\right) & =0,\label{eq:FHNEq1}\\
D_{y}Y_{c}''\left(\xi\right)+cY_{c}'\left(\xi\right)+R\left(Y_{c}\left(\xi\right)\right)-X_{c}\left(\xi\right) & =0.\label{eq:FHNEq2}
\end{align}
Exact analytical solutions to Eqs. \eqref{eq:FHNEq1} and \eqref{eq:FHNEq2}
are neither known for finite nor infinite or periodic domains. The
wave profile $\boldsymbol{X}_{c}$ and its velocity $c$ are determined
numerically and $\boldsymbol{X}_{c}$ is interpolated with Mathematica.
See Fig. \ref{fig:FHNWaveProfile} for a snapshot of $\boldsymbol{X}_{c}$.

\begin{minipage}{1.0\linewidth}
\begin{center}
\includegraphics[scale=0.55]{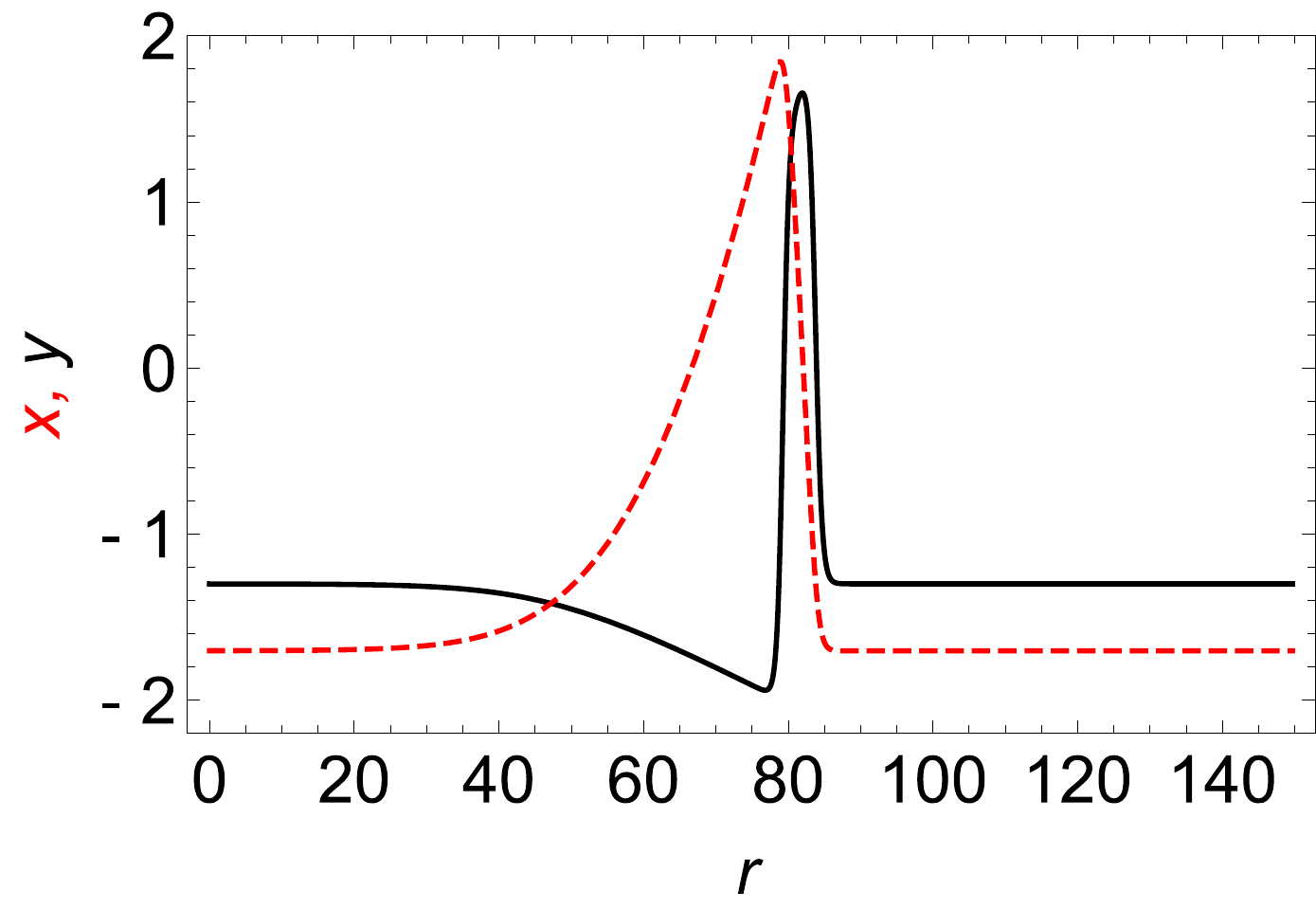}
\captionof{figure}[Wave profile of the uncontrolled FHN model]{\label{fig:FHNWaveProfile}Wave profile of the uncontrolled FHN model moving to the right. Shown is the activator component $y$ (black solid line) and the inhibitor component $x$ (red dashed line). The solution is obtained by numerically solving the uncontrolled reaction-diffusion system, Eq. \eqref{eq:RDSFHN} with $u\left(t \right) =0$, for periodic boundary conditions and then interpolated.}
\end{center}
%"/home/jakob/svnco/Control/FitzHughNagumo/ControlFitzHugNagumo2.nb"
\end{minipage}

The activator component $y_{d}\left(r,t\right)$ of the desired distribution
is the traveling wave profile $Y_{c}$ shifted according to the protocol
$\phi\left(t\right)$,
\begin{align}
y_{d}\left(r,t\right) & =Y_{c}\left(r-\phi\left(t\right)\right),
\end{align}
while the inhibitor component $x_{d}\left(r,t\right)$ has to satisfy
the partial differential equation 
\begin{align}
\partial_{t}x_{d}\left(r,t\right)-D_{x}\partial_{r}^{2}x_{d}\left(r,t\right)-a_{1}x_{d}\left(r,t\right) & =a_{0}+a_{2}Y_{c}\left(r-\phi\left(t\right)\right)\label{eq:RDSInhibitorConstraintEq}
\end{align}
with initial condition
\begin{align}
x_{d}\left(r,t_{0}\right) & =X_{c}\left(r-\phi\left(t_{0}\right)\right).
\end{align}
To simplify Eq. \eqref{eq:RDSInhibitorConstraintEq}, $\hat{x}_{d}$
is defined by the relation 
\begin{align}
x_{d}\left(r,t\right) & =\hat{x}_{d}\left(r,t\right)+X_{c}\left(r-\phi\left(t\right)\right).
\end{align}
After using Eq. \eqref{eq:FHNEq1}, the evolution equation for $\hat{x}_{d}$
becomes
\begin{align}
\partial_{t}\hat{x}_{d}\left(r,t\right)-D_{x}\partial_{r}^{2}\hat{x}_{d}\left(r,t\right)-a_{1}\hat{x}_{d}\left(r,t\right) & =-\left(c-\dot{\phi}\left(t\right)\right)X_{c}'\left(r-\phi\left(t\right)\right),\label{eq:Eq591}\\
\hat{x}_{d}\left(r,t_{0}\right) & =0.
\end{align}
Together with Eq. \eqref{eq:FHNEq2}, the control signal is given
by the relatively simple expression
\begin{align}
u\left(r,t\right) & =\left(c-\dot{\phi}\left(t\right)\right)Y_{c}'\left(r-\phi\left(t\right)\right)+\hat{x}_{d}\left(r,t\right).\label{eq:RDSFHNControlSignal}
\end{align}
Exactly the same result as Eq. \eqref{eq:RDSFHNControlSignal} was
derived in \cite{lober2014controlling} with an approach that focused
exclusively on position control of traveling waves.

In numerical simulations, the interpolated result for the wave profile
is used to formulate the control signal \eqref{eq:RDSFHNControlSignal},
and the resulting controlled reaction-diffusion system is solved numerically.
Although governed by a linear ODE, the solution for $\hat{x}_{d}\left(r,t\right)$
is determined numerically by solving Eq. \eqref{eq:Eq591} with periodic
boundary conditions for simplicity. The protocol of motion is,
\begin{align}
\phi\left(t\right) & =c\left(t-t_{0}\right)+A\sin\left(2\pi\left(t-t_{0}\right)/T\right)
\end{align}
and results in a controlled traveling wave moving sinusoidally back
and forth. The values for amplitude and period are $A=80$ and $T=20$,
respectively. Figure \ref{fig:FHNPositionControl1} compares the desired
activator (left) and inhibitor (right, black solid line) with the
numerically obtained result of the controlled reaction-diffusion system
(red dashed line). On this scale, the agreement is very good. Note
that while the controlled activator profile is identical to its uncontrolled
profile (see black solid line in Fig. \ref{fig:FHNWaveProfile}),
the inhibitor wave profile is largely deformed and very different
from its uncontrolled counterpart. The reason is simply that only
a single state component of the desired distribution can be prescribed,
which was chosen to be the activator component, while the inhibitor
component is determined by the constraint equation \eqref{eq:RDSInhibitorConstraintEq}.

\begin{minipage}{1.0\linewidth}
\begin{center}
\includegraphics[scale=0.46]{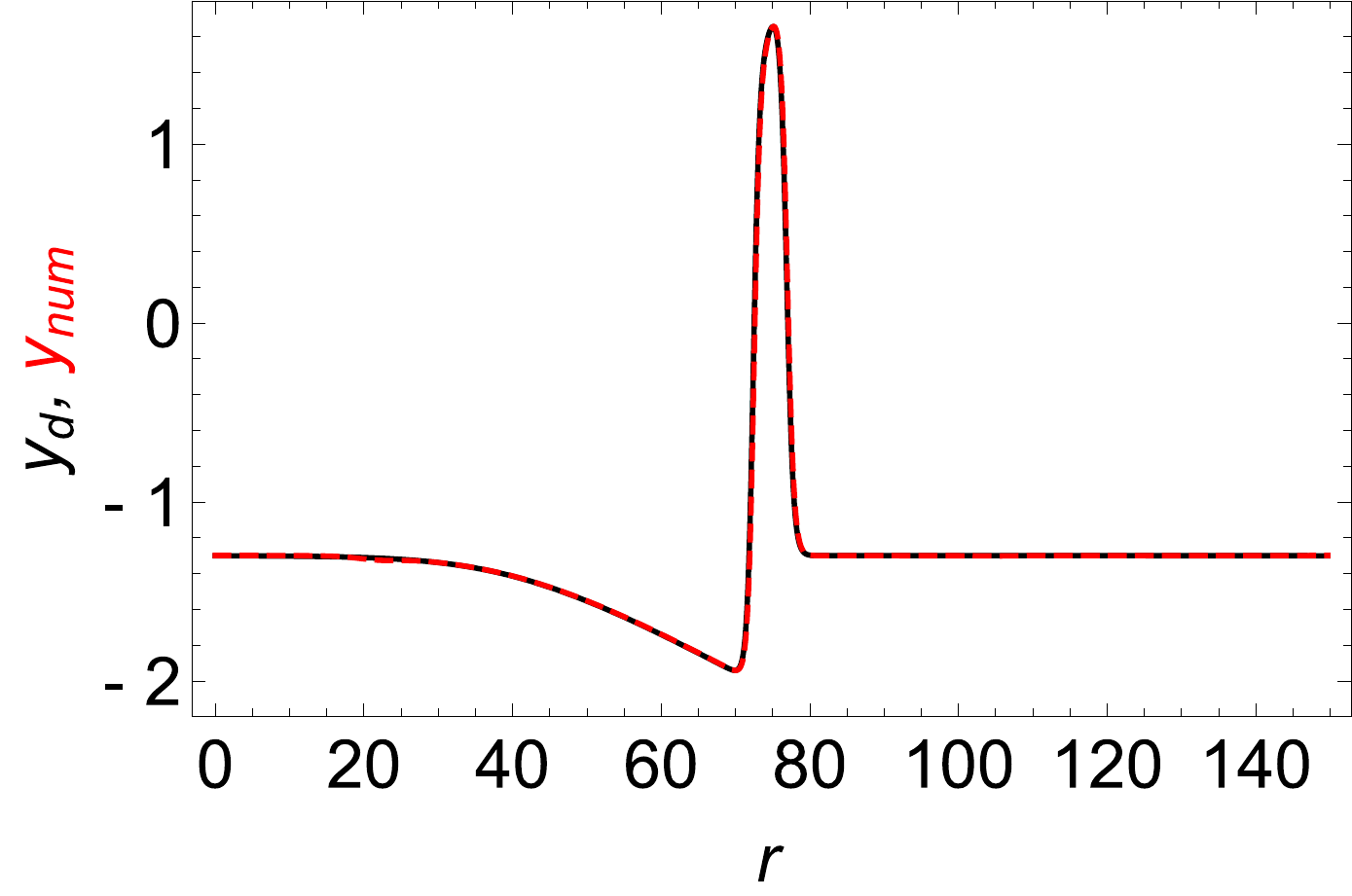}\hspace{0.5cm}\includegraphics[scale=0.46]{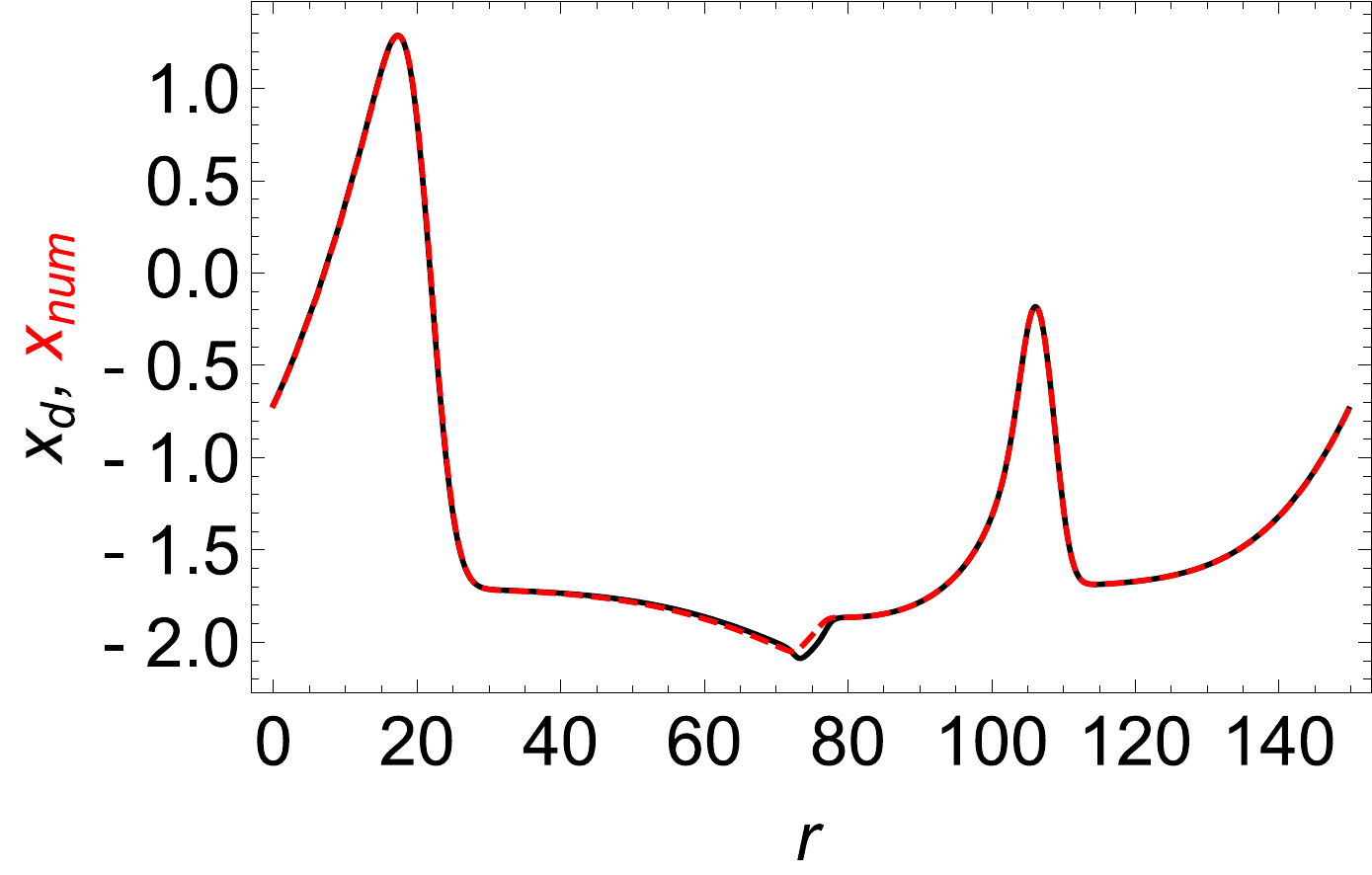}\vspace{0.2cm}
\captionof{figure}[Position control of traveling pulses in the FHN model]{\label{fig:FHNPositionControl1}Position control of a traveling pulse in the activator-controlled FHN model. The numerically obtained result for the controlled traveling pulse (black solid line) is very close to the desired distribution (red dashed line).  Left: Snapshot of controlled activator $y$ over space $r$. Right: Snapshot of controlled inhibitor $x$ over space $r$.}
\end{center}
%"/home/jakob/svnco/Control/FitzHughNagumo/ControlFitzHugNagumo2.nb"
\end{minipage}

Figure \ref{fig:FHNPositionControl2} shows the difference between
the desired and controlled wave profile for activator (left) and inhibitor
(right). The differences are \textit{not} in the range of numerical
accuracy, and are likely due to an instability. See \cite{lober2014stability}
for a discussion of one possible instability. Because of the steep
slopes exhibited by the pulse profile, a small difference in the position
between desired and controlled wave has a large effect on the difference
between the profiles. 

\begin{minipage}{1.0\linewidth}
\begin{center}
\includegraphics[scale=0.46]{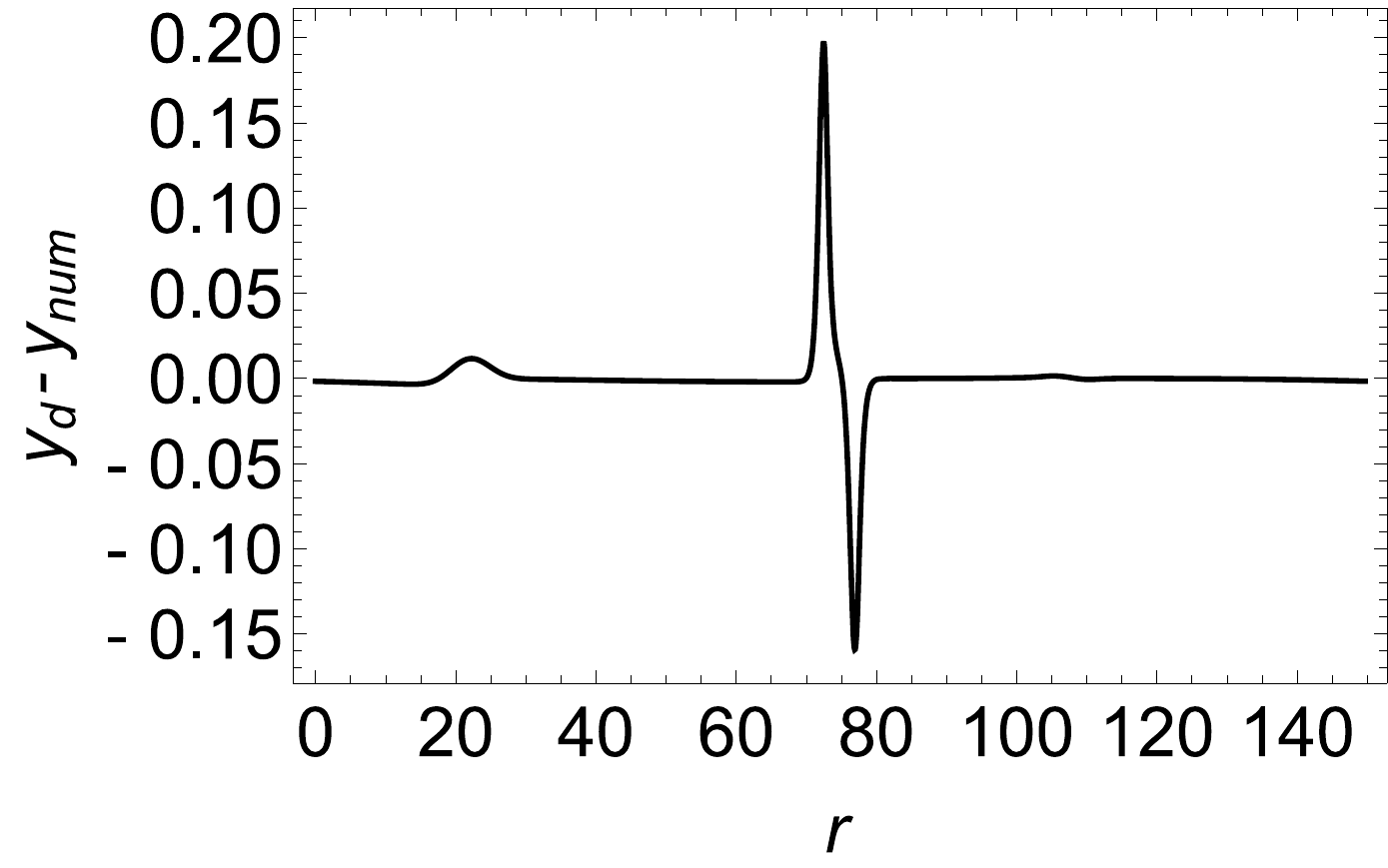}\hspace{0.5cm}\includegraphics[scale=0.46]{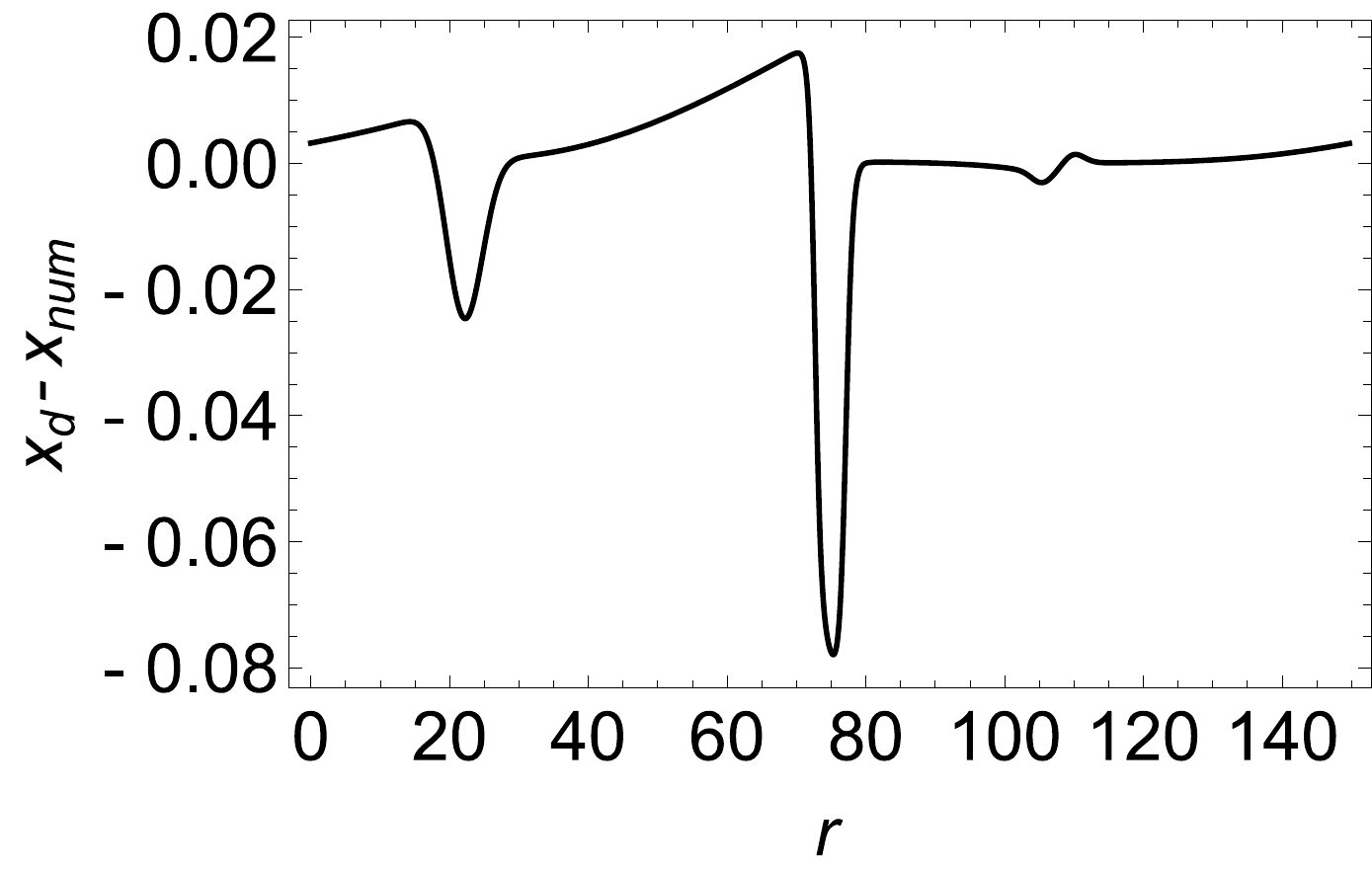}\vspace{0.2cm}
\captionof{figure}[Difference between controlled and desired traveling pulse]{\label{fig:FHNPositionControl2}Difference between controlled and desired traveling pulse for the activator (left) and inhibitor (right) component.}
\end{center}
%"/home/jakob/svnco/Control/FitzHughNagumo/ControlFitzHugNagumo2.nb"
\end{minipage}

Finally, Fig. \ref{fig:FHNPositionControl3} left shows the control
signal as given by Eq. \eqref{eq:RDSFHNControlSignal}. Being proportional
to the derivative of the activator pulse profile $Y_{c}$, the control
signal has its largest amplitude at the points of the steepest slope
of the activator pulse profile. Figure \ref{fig:FHNPositionControl3}
right compares the position over time as prescribed by the protocol
of motion $\phi$ (black solid line) with the position over time obtained
from numerical simulations (red dashed line). The agreement is well
within the range of numerical accuracy. Numerically, the position
of the pulse is defined as the position of the maximum of the controlled
activator pulse profile. 

\begin{minipage}{1.0\linewidth}
\begin{center}
\includegraphics[scale=0.475]{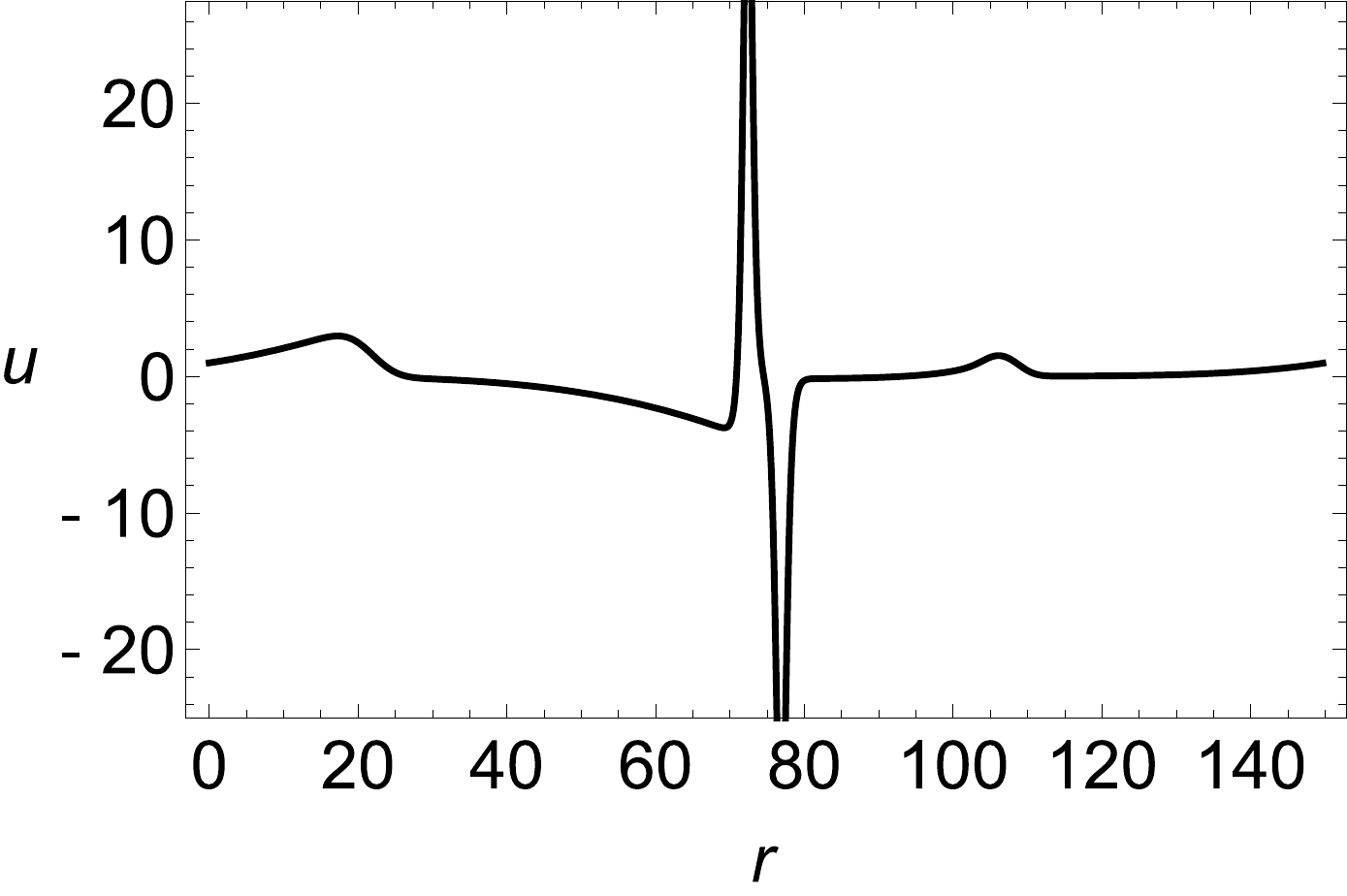}\hspace{0.2cm}\includegraphics[scale=0.475]{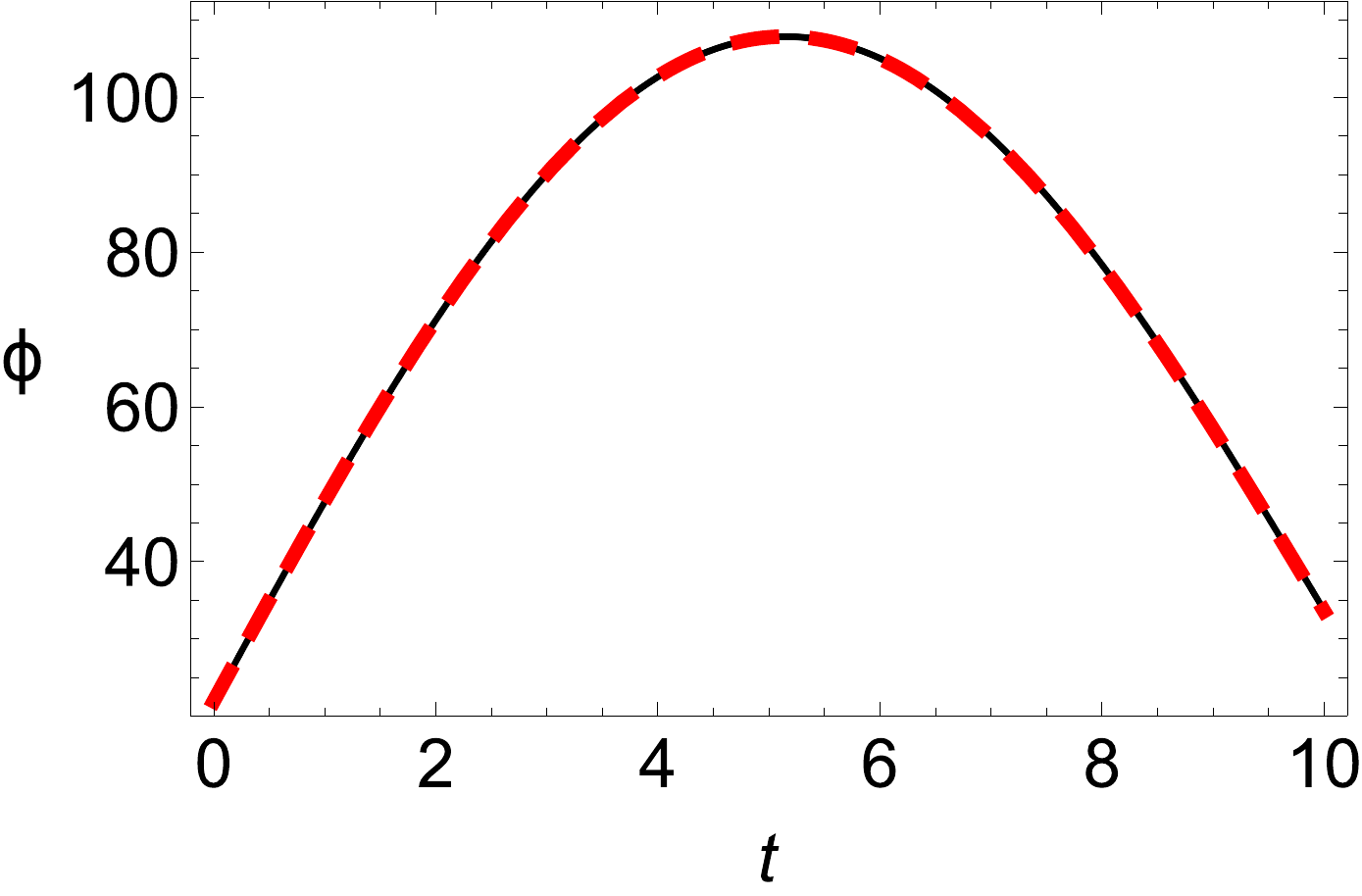}
\captionof{figure}[Control signal and controlled position over time for the FHN model ]{\label{fig:FHNPositionControl3}Control signal and controlled position over time for the FHN model. Left: Snapshot of control signal $u$ over space $r$. Right: Comparison of prescribed protocol $\phi$ over time (black solid line) and numerically recorded position over time of the controlled pulse   (red dashed line).}
\end{center}
%"/home/jakob/svnco/Control/FitzHughNagumo/ControlFitzHugNagumo2.nb"
\end{minipage}

\end{example}

\section{\label{sec:DiscussionAndOutlook}Discussion and outlook}

\subsection{Optimal control of reaction-diffusion systems}

This section briefly discusses optimal trajectory tracking for reaction-diffusion
systems. The mathematical theory of optimal control of PDEs is well
developed. The reader is referred to the book \cite{troltzsch2010optimal}
for a mathematically rigorous treatment. Applications of optimal control
to reaction-diffusion systems can be found in \cite{buchholz2013on,theissen2006optimale,Ryll2011}.

The target functional for optimal trajectory tracking in reaction-diffusion
systems is
\begin{align}
\mathcal{J}\left[\boldsymbol{x}\left(\boldsymbol{r},t\right),\boldsymbol{u}\left(\boldsymbol{r},t\right)\right]= & \frac{1}{2}\intop_{t_{0}}^{t_{1}}dt\intop_{\Omega}d\boldsymbol{r}\left(\boldsymbol{x}\left(\boldsymbol{r},t\right)-\boldsymbol{x}_{d}\left(\boldsymbol{r},t\right)\right)^{T}\boldsymbol{\mathcal{S}}\left(\boldsymbol{x}\left(\boldsymbol{r},t\right)-\boldsymbol{x}_{d}\left(\boldsymbol{r},t\right)\right)\nonumber \\
 & +\frac{1}{2}\intop_{\Omega}d\boldsymbol{r}\left(\boldsymbol{x}\left(\boldsymbol{r},t_{1}\right)-\boldsymbol{x}_{1}\left(\boldsymbol{r}\right)\right)^{T}\boldsymbol{\mathcal{S}}_{1}\left(\boldsymbol{x}\left(\boldsymbol{r},t_{1}\right)-\boldsymbol{x}_{1}\left(\boldsymbol{r}\right)\right)\nonumber \\
 & +\frac{\epsilon^{2}}{2}\intop_{t_{0}}^{t_{1}}dt\intop_{\Omega}d\boldsymbol{r}\left(\boldsymbol{u}\left(\boldsymbol{r},t\right)\right)^{2}.\label{eq:RDSJ}
\end{align}
Here, $\Omega$ denotes the $N$-dimensional spatial domain, $\boldsymbol{\mathcal{S}}$
and $\boldsymbol{\mathcal{S}}_{1}$ are symmetric matrices of weights,
and $\epsilon$ is the regularization parameter. Apart from the integration
over the spatial domain $\Omega$, the functional Eq. \eqref{eq:RDSJ}
is identical to the functional Eq. \eqref{eq:OptimalTrajectoryTrackingFunctional}
for optimal trajectory tracking in dynamical systems from Chapter
\ref{chap:OptimalControl}. Equation \eqref{eq:RDSJ} must be minimized
under the constraint that $\boldsymbol{x}\left(\boldsymbol{r},t\right)$
is governed by the controlled reaction-diffusion equation 
\begin{align}
\partial_{t}\boldsymbol{x}\left(\boldsymbol{r},t\right) & =\boldsymbol{\mathcal{D}}\triangle\boldsymbol{x}\left(\boldsymbol{r},t\right)+\boldsymbol{R}\left(\boldsymbol{x}\left(\boldsymbol{r},t\right)\right)+\boldsymbol{\mathcal{B}}\left(\boldsymbol{x}\left(\boldsymbol{r},t\right),\boldsymbol{r}\right)\boldsymbol{u}\left(\boldsymbol{r},t\right),\label{eq:RDSNecessaryOptimalityConditions1}
\end{align}
supplemented with the boundary and initial conditions
\begin{align}
\boldsymbol{0} & =\boldsymbol{n}^{T}\left(\boldsymbol{r}\right)\left(\boldsymbol{\mathcal{D}}\nabla\boldsymbol{x}\left(\boldsymbol{r},t\right)\right),\,\boldsymbol{r}\in\Gamma, & \boldsymbol{x}\left(\boldsymbol{r},t_{0}\right) & =\boldsymbol{x}_{0}\left(\boldsymbol{r}\right).
\end{align}
Similar as in Section \ref{sec:NecessaryOptimalityConditions}, the
constrained minimization problem can be transformed to an unconstrained
minimization problem by introducing the vector of Lagrange multipliers
$\boldsymbol{\lambda}\left(\boldsymbol{r},t\right)$, also called
adjoint state or co-state. This leads to the adjoint or co-state equation
for $\boldsymbol{\lambda}$ \cite{troltzsch2010optimal,theissen2006optimale,buchholz2013on},
\begin{align}
-\partial_{t}\boldsymbol{\lambda}\left(\boldsymbol{r},t\right) & =\boldsymbol{\mathcal{D}}\triangle\boldsymbol{\lambda}\left(\boldsymbol{r},t\right)+\left(\nabla\boldsymbol{R}^{T}\left(\boldsymbol{x}\left(\boldsymbol{r},t\right)\right)+\boldsymbol{u}^{T}\left(\boldsymbol{r},t\right)\nabla\boldsymbol{\mathcal{B}}^{T}\left(\boldsymbol{x}\left(\boldsymbol{r},t\right),\boldsymbol{r}\right)\right)\boldsymbol{\lambda}\left(\boldsymbol{r},t\right)\nonumber \\
 & +\boldsymbol{\mathcal{S}}\left(\boldsymbol{x}\left(\boldsymbol{r},t\right)-\boldsymbol{x}_{d}\left(\boldsymbol{r},t\right)\right).\label{eq:RDSAdjointEquation}
\end{align}
The co-state $\boldsymbol{\lambda}$ satisfies the same homogeneous
Neumann boundary conditions as the state $\boldsymbol{x}$, 
\begin{align}
\boldsymbol{0} & =\boldsymbol{n}^{T}\left(\boldsymbol{r}\right)\left(\boldsymbol{\mathcal{D}}\nabla\boldsymbol{\lambda}\left(\boldsymbol{r},t\right)\right),\,\boldsymbol{r}\in\Gamma,
\end{align}
and the terminal conditions
\begin{align}
\boldsymbol{\lambda}\left(\boldsymbol{r},t_{1}\right) & =\boldsymbol{\mathcal{S}}_{1}\left(\boldsymbol{x}\left(\boldsymbol{r},t_{1}\right)-\boldsymbol{x}_{1}\left(\boldsymbol{r}\right)\right).
\end{align}
Finally, the relation between control signal $\boldsymbol{u}$ and
co-state $\boldsymbol{\lambda}$ is obtained as 
\begin{align}
\epsilon^{2}\boldsymbol{u}\left(\boldsymbol{r},t\right)+\boldsymbol{\mathcal{B}}^{T}\left(\boldsymbol{x}\left(\boldsymbol{r},t\right),\boldsymbol{r}\right)\boldsymbol{\lambda}\left(\boldsymbol{r},t\right)= & \boldsymbol{0}.\label{eq:RDSNecessaryOptimalityConditionsLast}
\end{align}
Equations \eqref{eq:RDSNecessaryOptimalityConditions1}-\eqref{eq:RDSNecessaryOptimalityConditionsLast}
constitute the \textit{necessary optimality conditions for optimal
trajectory tracking in reaction-diffusion systems}.

For dynamical systems, it was found that the control signal obtained
within the framework of exactly realizable trajectories arises as
the solution to an unregularized optimal control problem. We expect
a similar identity for exactly realizable distributions of spatio-temporal
systems. Indeed, if the state equals the desired trajectory everywhere
and for all times $t_{0}\leq t\leq t_{1}$, $\boldsymbol{x}\left(\boldsymbol{r},t\right)=\boldsymbol{x}_{d}\left(\boldsymbol{r},t\right)$,
Eq. \eqref{eq:RDSAdjointEquation} becomes a homogeneous linear partial
differential equation. If additionally the desired distribution complies
with the terminal state, $\boldsymbol{x}_{d}\left(\boldsymbol{r},t_{1}\right)=\boldsymbol{x}_{1}\left(\boldsymbol{r}\right)$,
the co-state $\boldsymbol{\lambda}\left(\boldsymbol{r},t\right)$
vanishes identically everywhere and for all times,
\begin{align}
\boldsymbol{\lambda}\left(\boldsymbol{r},t\right) & \equiv\boldsymbol{0}.
\end{align}
It follows that for a non-vanishing control signal $\boldsymbol{u}\left(t\right)$,
Eq. \eqref{eq:RDSNecessaryOptimalityConditionsLast} can only be satisfied
if $\epsilon=0$. In conclusion, all necessary optimality conditions
Eqs. \eqref{eq:RDSNecessaryOptimalityConditions1}-\eqref{eq:RDSNecessaryOptimalityConditionsLast}
are satisfied.

However, analogously to the generalized Legendre-Clebsch conditions
for dynamical systems, Eqs. \eqref{eq:LegendreClebsch1} and \eqref{eq:LegendreClebsch2},
we expect that there are additional necessary optimality conditions
for singular optimal control problems, see Section \ref{sub:TheGeneralizedLegendreClebschConditions}.
While the necessity of the generalized Legendre-Clebsch conditions
for dynamical systems is rigorously proven in \cite{bell1975singular},
there seems to be no rigorous proof available for singular optimal
control of PDEs. We omit a discussion of additional necessary optimality
conditions.

\subsection{Outlook}

A possible next step is the application of the singular perturbation
expansion developed in Chapter \ref{chap:AnalyticalApproximationsForOptimalTrajectoryTracking}
to the necessary optimality conditions Eqs. \eqref{eq:RDSNecessaryOptimalityConditions1}-\eqref{eq:RDSNecessaryOptimalityConditionsLast}.
First, the necessary optimality conditions must be rearranged and
split up in equations for the parts $\boldsymbol{\mathcal{P}}\boldsymbol{x}$,
$\boldsymbol{\mathcal{Q}}\boldsymbol{x}$, $\boldsymbol{\mathcal{P}}\boldsymbol{\lambda}$,
and $\boldsymbol{\mathcal{Q}}\boldsymbol{\lambda}$. Second, the inner
and outer equations must be derived. In general, not only time but
also space can be rescaled with the small parameter $\epsilon$. This
might lead to a larger variety of inner equations and combinations
of spatial and temporal boundary layers. However, at least for problems
as simple as the activator-controlled FHN model and a control acting
everywhere within the spatial domain $\Omega$, it seems reasonable
to expect a simple structure of inner and outer equations analogously
to the two-dimensional dynamical system from Section \ref{sec:TwoDimensionalDynamicalSystem}.

The essential difference in the evolution equations between dynamical
systems and reaction-diffusion systems is the diffusion term. Being
a linear differential operator, we anticipate that if the outer equations
of a dynamical system reduce to linear ODEs as the result of a linearizing
assumption, the outer equations for a corresponding reaction-diffusion
system reduce to linear PDEs. This opens up the interesting possibility
to obtain analytical approximations for the optimal control of reaction-diffusion
systems. The arising equations will be linear reaction-diffusion equations
with inhomogeneities which involve the desired distribution $\boldsymbol{x}_{d}$.
Such equations can in principle be solved analytically with the help
of Green's functions. These solutions would not only provide analytical
approximations for open loop control, but would also yield optimal
feedback controls for nonlinear reaction-diffusion systems. As was
discussed above, it is virtually impossible to numerically deal with
optimal feedback control of spatio-temporal systems due to the curse
of dimensionality. The approach outlined here would enable an almost
exclusive approach to optimal feedback control of nonlinear spatio-temporal
systems.

Stability of open loop control methods can never be taken for granted
but requires further investigations. The control of exactly realizable
distributions might or might not be stable with respect to perturbations
of the initial conditions. According to Eq. \eqref{eq:LinearPDEy}
governing the stability of an exactly realizable distribution $\boldsymbol{x}_{d}\left(\boldsymbol{r},t\right)$,
the stability of $\boldsymbol{x}_{d}\left(\boldsymbol{r},t\right)$
depends on $\boldsymbol{x}_{d}\left(\boldsymbol{r},t\right)$ itself.
This observation opens up the investigation of the stability of position
control of traveling waves. Assuming for simplicity a desired distribution
of the form
\begin{align}
\boldsymbol{x}_{d}\left(\boldsymbol{r},t\right) & =\boldsymbol{X}_{c}\left(\boldsymbol{\hat{c}}^{T}\boldsymbol{r}-\phi\left(t\right)\right),
\end{align}
and a constant coupling matrix $\boldsymbol{\mathcal{B}}\left(\boldsymbol{x},\boldsymbol{r}\right)=\boldsymbol{\mathcal{B}}=\text{const.}$,
Eq. \eqref{eq:LinearPDEy}, becomes 
\begin{align}
\partial_{t}\boldsymbol{y}\left(\boldsymbol{r},t\right) & =\boldsymbol{\mathcal{D}}\triangle\boldsymbol{y}\left(\boldsymbol{r},t\right)+\nabla\boldsymbol{R}\left(\boldsymbol{X}_{c}\left(\boldsymbol{\hat{c}}^{T}\boldsymbol{r}-\phi\left(t\right)\right)\right)\boldsymbol{y}\left(\boldsymbol{r},t\right).\label{eq:Eq5103}
\end{align}
If additionally the protocol velocity is close to the velocity $c$
of the uncontrolled traveling wave, $\dot{\phi}\left(t\right)=c+\gamma$
with $\left|\gamma\right|\ll1$, Eq. \eqref{eq:Eq5103} reduces to
the equation which determines the linear stability of the traveling
wave $\boldsymbol{X}_{c}$ \cite{Sandstede2002983}. As long as the
exactly realizable desired distribution $\boldsymbol{x}_{d}\left(\boldsymbol{r},t\right)$
is sufficiently close to a stable traveling wave solution $\boldsymbol{X}_{c}$,
the controlled wave may be stable. In this way, the controlled wave
may benefit from the stability of the uncontrolled traveling wave.
A rigorous discussion of stability must take into account the fact
that only $p$ out of $n$ components of a desired distribution can
be prescribed, and should take into account a state dependent coupling
matrix $\boldsymbol{\mathcal{B}}\left(\boldsymbol{x}\right)$. Additional
problems arise because every stable traveling wave possesses at least
one eigenvalue with vanishing real part. This fact requires a nonlinear
stability analysis. A popular method for that is a multiple scale
perturbation expansion \cite{lober2009nonlinear,lober2012front}.
Some aspects of this nonlinear stability analysis for position control
of traveling waves are presented in \cite{lober2014stability}, see
also the discussion at the end of \cite{lober2014shaping}. Generally
speaking, one can expect a stable open loop control of an exactly
realizable desired distribution as long as the desired distribution
is sufficiently close to a stable solution of the uncontrolled problem.
A thorough understanding of the solutions to an uncontrolled system,
including their stability properties, can be very useful for the design
of exactly realizable desired distributions which do not require stabilization
by additional feedback. Note that the stability analysis of desired
distributions which are not exactly realizable is much more difficult.
The reason is that the controlled state might be very different from
the desired distribution. In general, the stability properties of
controlled and desired state are unrelated.

%% file: Appendix.tex
\lhead[\leftmark]{}

\rhead[]{\rightmark}

\lfoot[\thepage]{}

\cfoot{}

\rfoot[]{\thepage}

\chapter{Appendix}

\section{\label{sec:GeneralSolutionForForcedLinarDynamicalSystem}General
solution for a forced linear dynamical system}

Consider a linear $n$-dimensional dynamical system 
\begin{align}
\boldsymbol{\dot{x}}\left(t\right) & =\boldsymbol{\mathcal{A}}\left(t\right)\boldsymbol{x}\left(t\right)+\boldsymbol{f}\left(t\right),\label{eq:InhomogeneousLinearSystem}
\end{align}
and initial conditions 
\begin{align}
\boldsymbol{x}\left(t_{0}\right) & =\boldsymbol{x}_{0},\label{eq:InitialCondition}
\end{align}
for the state $\boldsymbol{x}$ 
\begin{align}
\boldsymbol{x}\left(t\right) & =\left(x_{1}\left(t\right),\dots,x_{n}\left(t\right)\right)^{T}
\end{align}
with forcing or inhomogeneity $\boldsymbol{f}$
\begin{align}
\boldsymbol{f}\left(t\right) & =\left(f_{1}\left(t\right),\dots,f_{n}\left(t\right)\right)^{T}.
\end{align}
Dynamical systems of the form Eq. \eqref{eq:InhomogeneousLinearSystem}
are called linear time-variant (LTV) in control theory. If $\boldsymbol{\mathcal{A}}\left(t\right)=\boldsymbol{\mathcal{A}}=\text{const.}$
does not depend on time, Eq. \eqref{eq:InhomogeneousLinearSystem}
is called a linear time invariant (LTI) system. See the excellent
book \cite{chen1995linear} and also \cite{katsuhiko2010modern} for
an exhaustive treatment of LTV and LTI systems. The general solution
of Eq. \eqref{eq:InhomogeneousLinearSystem} can be expressed in terms
of the principal fundamental $n\times n$ matrix $\boldsymbol{\Phi}\left(t,t_{0}\right)$,
also called state transition matrix, which satisfies
\begin{align}
\partial_{t}\boldsymbol{\Phi}\left(t,t_{0}\right) & =\boldsymbol{\mathcal{A}}\left(t\right)\boldsymbol{\Phi}\left(t,t_{0}\right), & \boldsymbol{\Phi}\left(t_{0},t_{0}\right) & =\mathbf{1}.\label{eq:LinSysStateTransitionMatrix}
\end{align}
$\boldsymbol{\Phi}$ must be a nonsingular matrix such that its inverse
$\boldsymbol{\Phi}^{-1}\left(t_{2},t_{1}\right)=\boldsymbol{\Phi}\left(t_{1},t_{2}\right)$
exists. This implies
\begin{align}
\boldsymbol{\Phi}^{-1}\left(t_{0},t_{0}\right) & =\mathbf{1},\label{eq:InitialInverse}
\end{align}
and
\begin{align}
\boldsymbol{\Phi}\left(t_{0},t\right)\boldsymbol{\Phi}\left(t,t_{0}\right) & =\boldsymbol{\Phi}\left(t,t_{0}\right)\boldsymbol{\Phi}\left(t_{0},t\right)=\mathbf{1}.\label{eq:InverseIdentity}
\end{align}
Applying the derivative with respect to time to Eq. \eqref{eq:InverseIdentity}
yields
\begin{align}
\partial_{t}\boldsymbol{\Phi}\left(t_{0},t\right)\boldsymbol{\Phi}\left(t,t_{0}\right) & =-\boldsymbol{\Phi}\left(t_{0},t\right)\partial_{t}\boldsymbol{\Phi}\left(t,t_{0}\right)=-\boldsymbol{\Phi}\left(t_{0},t\right)\boldsymbol{\boldsymbol{\mathcal{A}}}\left(t\right)\boldsymbol{\Phi}\left(t,t_{0}\right).\label{eq:DerivativeInverseIdentity}
\end{align}
From Eq. \eqref{eq:DerivativeInverseIdentity} follows the useful
relation
\begin{align}
\boldsymbol{\Phi}\left(t_{0},t\right)\boldsymbol{\boldsymbol{\mathcal{A}}}\left(t\right) & =-\partial_{t}\boldsymbol{\Phi}\left(t_{0},t\right).\label{eq:UsefulRelation}
\end{align}
Transposing  Eq. \eqref{eq:UsefulRelation} yields the so-called adjoint
equation
\begin{align}
\partial_{t}\boldsymbol{\Phi}^{T}\left(t_{0},t\right)= & -\boldsymbol{\boldsymbol{\mathcal{A}}}^{T}\left(t\right)\boldsymbol{\Phi}^{T}\left(t_{0},t\right).\label{eq:FundamentalMatrixAdjointEquation}
\end{align}
Hence, if $\boldsymbol{\Phi}\left(t,t_{0}\right)$ is the fundamental
matrix to the original system, then the inverse and transposed matrix
$\boldsymbol{\Phi}^{-T}\left(t,t_{0}\right)=\boldsymbol{\Phi}^{T}\left(t_{0},t\right)$
is the fundamental matrix to the adjoint system. The general solution
$\boldsymbol{x}\left(t\right)$ to the inhomogeneous linear system
Eq. \eqref{eq:InhomogeneousLinearSystem} is a superposition
\begin{align}
\boldsymbol{x}\left(t\right) & =\boldsymbol{y}\left(t\right)+\boldsymbol{z}\left(t\right)
\end{align}
of the solution $\boldsymbol{y}\left(t\right)$ to the homogeneous
system 
\begin{align}
\boldsymbol{\dot{y}}\left(t\right) & =\boldsymbol{\boldsymbol{\mathcal{A}}}\left(t\right)\boldsymbol{y}\left(t\right), & \boldsymbol{y}\left(t_{0}\right) & =\boldsymbol{x}_{0},
\end{align}
and a solution $\boldsymbol{z}\left(t\right)$ of the inhomogeneous
system. The homogeneous solution $\boldsymbol{y}\left(t\right)$ can
be written in terms of the fundamental matrix as
\begin{align}
\boldsymbol{y}\left(t\right) & =\boldsymbol{\Phi}\left(t,t_{0}\right)\boldsymbol{x}_{0}.
\end{align}
The proof is very simple
\begin{align}
\boldsymbol{\dot{y}}\left(t\right) & =\partial_{t}\boldsymbol{\Phi}\left(t,t_{0}\right)\boldsymbol{x}_{0}=\boldsymbol{\boldsymbol{\mathcal{A}}}\left(t\right)\boldsymbol{\Phi}\left(t,t_{0}\right)\boldsymbol{x}_{0}=\boldsymbol{\boldsymbol{\mathcal{A}}}\left(t\right)\boldsymbol{y}\left(t\right),\\
\boldsymbol{y}\left(t_{0}\right) & =\boldsymbol{\Phi}\left(t_{0},t_{0}\right)\boldsymbol{x}_{0}=\boldsymbol{x}_{0}.
\end{align}
The ansatz for the solution of the inhomogeneous linear system Eq.
\eqref{eq:InhomogeneousLinearSystem} is
\begin{align}
\boldsymbol{z}\left(t\right) & =\boldsymbol{\Phi}\left(t,t_{0}\right)\boldsymbol{v}\left(t\right).\label{eq:ZAnsatz}
\end{align}
Using the ansatz Eq. \eqref{eq:ZAnsatz} in the inhomogeneous linear
system \eqref{eq:InhomogeneousLinearSystem} yields
\begin{align}
\boldsymbol{\dot{z}}\left(t\right) & =\partial_{t}\boldsymbol{\Phi}\left(t,t_{0}\right)\boldsymbol{v}\left(t\right)+\boldsymbol{\Phi}\left(t,t_{0}\right)\boldsymbol{\dot{v}}\left(t\right)=\boldsymbol{\boldsymbol{\mathcal{A}}}\left(t\right)\boldsymbol{\Phi}\left(t,t_{0}\right)\boldsymbol{v}\left(t\right)+\boldsymbol{\Phi}\left(t,t_{0}\right)\boldsymbol{\dot{v}}\left(t\right)\nonumber \\
 & =\boldsymbol{\boldsymbol{\mathcal{A}}}\left(t\right)\boldsymbol{\Phi}\left(t,t_{0}\right)\boldsymbol{v}\left(t\right)+\boldsymbol{f}\left(t\right).
\end{align}
It follows that
\begin{align}
\boldsymbol{\Phi}\left(t,t_{0}\right)\boldsymbol{\dot{v}}\left(t\right) & =\boldsymbol{f}\left(t\right)
\end{align}
and, after rearranging and integrating over time,
\begin{align}
\boldsymbol{v}\left(t\right) & =\intop_{t_{0}}^{t}d\tau\boldsymbol{\Phi}^{-1}\left(\tau,t_{0}\right)\boldsymbol{f}\left(\tau\right)=\intop_{t_{0}}^{t}d\tau\boldsymbol{\Phi}\left(t_{0},\tau\right)\boldsymbol{f}\left(\tau\right).
\end{align}
The solution for $\boldsymbol{z}\left(t\right)$ is thus
\begin{align}
\boldsymbol{z}\left(t\right) & =\boldsymbol{\Phi}\left(t,t_{0}\right)\boldsymbol{v}\left(t\right)=\boldsymbol{\Phi}\left(t,t_{0}\right)\intop_{t_{0}}^{t}d\tau\boldsymbol{\Phi}\left(t_{0},\tau\right)\boldsymbol{f}\left(\tau\right)=\intop_{t_{0}}^{t}d\tau\boldsymbol{\Phi}\left(t,\tau\right)\boldsymbol{f}\left(\tau\right).
\end{align}
The general solution $\boldsymbol{x}\left(t\right)=\boldsymbol{y}\left(t\right)+\boldsymbol{z}\left(t\right)$
to Eq. \eqref{eq:InhomogeneousLinearSystem} is then
\begin{align}
\boldsymbol{x}\left(t\right) & =\boldsymbol{\Phi}\left(t,t_{0}\right)\boldsymbol{x}_{0}+\intop_{t_{0}}^{t}d\tau\boldsymbol{\Phi}\left(t,\tau\right)\boldsymbol{f}\left(\tau\right).\label{eq:GeneralSolutionLinearSystem}
\end{align}
For an LTI system with constant state matrix $\boldsymbol{\mathcal{A}}\left(t\right)=\boldsymbol{\mathcal{A}}=\text{const.}$,
the solution for the state transition matrix $\boldsymbol{\Phi}$
is
\begin{align}
\boldsymbol{\Phi}\left(t,t_{0}\right) & =\exp\left(\boldsymbol{\mathcal{A}}\left(t-t_{0}\right)\right).\label{eq:MatrixExponential}
\end{align}
The matrix exponential is defined by the power series
\begin{align}
\exp\left(\boldsymbol{\mathcal{A}}\right) & =\sum_{k=0}^{\infty}\dfrac{1}{k!}\boldsymbol{\mathcal{A}}^{k}.
\end{align}
A proof of Eq. \eqref{eq:MatrixExponential} reads as follows. The
derivative of $\boldsymbol{\Phi}\left(t,t_{0}\right)$ with respect
to time $t$ is
\begin{align}
\partial_{t}\boldsymbol{\Phi}\left(t,t_{0}\right) & =\partial_{t}\exp\left(\boldsymbol{\mathcal{A}}\left(t-t_{0}\right)\right)=\partial_{t}\left(\sum_{k=0}^{\infty}\dfrac{1}{k!}\boldsymbol{\mathcal{A}}^{k}\left(t-t_{0}\right)^{k}\right)=\sum_{k=0}^{\infty}\dfrac{1}{k!}\boldsymbol{\mathcal{A}}^{k}\partial_{t}\left(t-t_{0}\right)^{k}\nonumber \\
 & =\sum_{k=1}^{\infty}\dfrac{k}{k!}\boldsymbol{\mathcal{A}}^{k}\left(t-t_{0}\right)^{k-1}=\sum_{k=1}^{\infty}\dfrac{1}{\left(k-1\right)!}\boldsymbol{\mathcal{A}}^{k}\left(t-t_{0}\right)^{k-1}=\sum_{\tilde{k}=0}^{\infty}\dfrac{1}{\tilde{k}!}\boldsymbol{\mathcal{A}}^{\tilde{k}+1}\left(t-t_{0}\right)^{\tilde{k}}\nonumber \\
 & =\boldsymbol{\mathcal{A}}\sum_{\tilde{k}=0}^{\infty}\dfrac{1}{\tilde{k}!}\boldsymbol{\mathcal{A}}^{\tilde{k}}\left(t-t_{0}\right)^{\tilde{k}}=\boldsymbol{\mathcal{A}}\exp\left(\boldsymbol{\mathcal{A}}\left(t-t_{0}\right)\right)=\boldsymbol{\mathcal{A}}\boldsymbol{\Phi}\left(t,t_{0}\right).
\end{align}
The index shift $\tilde{k}=k-1$ was introduced in the second line.

\section{\label{sec:OverAndUnderdetSysOfEqs}Over- and underdetermined systems
of linear equations}

The solutions of over- and underdetermined systems of linear equations
are discussed.

\subsection{\label{sub:GeneralizedInverseMatrices}Generalized inverse matrices}

The inverse $\boldsymbol{\mathcal{A}}^{-1}$ of a matrix $\boldsymbol{\mathcal{A}}$
with real or complex entries satisfies $\boldsymbol{\mathcal{A}}\boldsymbol{\mathcal{A}}^{-1}=\boldsymbol{\mathcal{A}}^{-1}\boldsymbol{\mathcal{A}}=\boldsymbol{1}$.
An $n\times m$ matrix $\boldsymbol{\mathcal{A}}$ has an inverse
only if it is square, i.e., $m=n$, and full rank, i.e., $\text{rank}\left(\boldsymbol{\mathcal{A}}\right)=n$.
For other matrices, a generalized inverse can be defined.

A generalized inverse $\boldsymbol{\mathcal{A}}^{g}$ of the $n\times m$
matrix $\boldsymbol{\mathcal{A}}$ with real entries has to satisfy
the condition
\begin{align}
\boldsymbol{\mathcal{A}}\boldsymbol{\mathcal{A}}^{g}\boldsymbol{\mathcal{A}} & =\boldsymbol{\mathcal{A}}.\label{eq:MP1}
\end{align}
If $\boldsymbol{\mathcal{A}}^{g}$ additionally satisfies the condition
\begin{align}
\boldsymbol{\mathcal{A}}^{g}\boldsymbol{\mathcal{A}}\boldsymbol{\mathcal{A}}^{g} & =\boldsymbol{\mathcal{A}}^{g},\label{eq:MP2}
\end{align}
$\boldsymbol{\mathcal{A}}^{g}$ is called a generalized reflexive
inverse. Furthermore, if $\boldsymbol{\mathcal{A}}^{g}$ satisfies
additionally the conditions
\begin{align}
\left(\boldsymbol{\mathcal{A}}\boldsymbol{\mathcal{A}}^{g}\right)^{T} & =\boldsymbol{\mathcal{A}}\boldsymbol{\mathcal{A}}^{g},\label{eq:MP3}
\end{align}
and
\begin{align}
\left(\boldsymbol{\mathcal{A}}^{g}\boldsymbol{\mathcal{A}}\right)^{T} & =\boldsymbol{\mathcal{A}}^{g}\boldsymbol{\mathcal{A}},\label{eq:MP4}
\end{align}
$\boldsymbol{\mathcal{A}}^{g}$ is called the Moore-Penrose pseudo
inverse matrix and denoted by $\boldsymbol{\mathcal{A}}^{+}$. For
any matrix $\boldsymbol{\mathcal{A}}$ with real or complex entries,
the Moore-Penrose pseudo inverse $\boldsymbol{\mathcal{A}}^{+}$ exists
and is unique. A generalized inverse satisfying only condition \eqref{eq:MP1}
is usually not unique \cite{CampbellJr.1991Generalized}.

\subsection{Solving an overdetermined system of linear equations}

An overdetermined system of equations has more equations than unknowns.
Let $\boldsymbol{x}\in\mathbb{R}^{p}$ and $\boldsymbol{b}\in\mathbb{R}^{n}$
with $p<n$, and let $\boldsymbol{\mathcal{A}}$ be an $n\times p$
matrix. The aim is to solve the system of $n$ equations
\begin{align}
\boldsymbol{\mathcal{A}}\boldsymbol{x} & =\boldsymbol{b}\label{eq:AxEqualsb}
\end{align}
for $\boldsymbol{x}$. Such overdetermined equations regularly occur
in data fitting problems. Because $\boldsymbol{\mathcal{A}}$ is not
a quadratic matrix, an exact solution cannot exist. However, a useful
expression for $\boldsymbol{x}$ can be derived as follows. Multiplying
Eq. \eqref{eq:AxEqualsb} with $\boldsymbol{\mathcal{A}}^{T}$ yields
\begin{align}
\boldsymbol{\mathcal{A}}^{T}\boldsymbol{\mathcal{A}}\boldsymbol{x} & =\boldsymbol{\mathcal{A}}^{T}\boldsymbol{b}.\label{eq:EqA30}
\end{align}
To solve for $\boldsymbol{x}$, Eq. \eqref{eq:EqA30} is multiplied
with the inverse of the $p\times p$ matrix $\boldsymbol{\mathcal{A}}^{T}\boldsymbol{\mathcal{A}}$
from the left to get
\begin{align}
\boldsymbol{x} & =\left(\boldsymbol{\mathcal{A}}^{T}\boldsymbol{\mathcal{A}}\right)^{-1}\boldsymbol{\mathcal{A}}^{T}\boldsymbol{b}=\boldsymbol{\mathcal{A}}^{+}\boldsymbol{b}.\label{eq:SolutionForx}
\end{align}
The $p\times n$ matrix $\boldsymbol{\mathcal{A}}^{+}$ is defined
as 
\begin{align}
\boldsymbol{\mathcal{A}}^{+} & =\left(\boldsymbol{\mathcal{A}}^{T}\boldsymbol{\mathcal{A}}\right)^{-1}\boldsymbol{\mathcal{A}}^{T}.
\end{align}
The matrix $\boldsymbol{\mathcal{A}}^{+}$ is the Moore-Penrose pseudo
inverse of matrix $\boldsymbol{\mathcal{A}}$, which can be proven
by checking all four conditions Eqs. \eqref{eq:MP1}-\eqref{eq:MP4}.
The inverse of $\boldsymbol{\mathcal{A}}^{T}\boldsymbol{\mathcal{A}}$
exists whenever $\boldsymbol{\mathcal{A}}$ has full column rank $p$,
\begin{align}
\text{rank}\left(\boldsymbol{\mathcal{A}}\right) & =p.
\end{align}
If $p=n$ and $\boldsymbol{\mathcal{A}}$ has full rank, the inverse
of $\boldsymbol{\mathcal{A}}$ exists and
\begin{align}
\boldsymbol{\mathcal{A}}^{+} & =\left(\boldsymbol{\mathcal{A}}^{T}\boldsymbol{\mathcal{A}}\right)^{-1}\boldsymbol{\mathcal{A}}^{T}=\boldsymbol{\mathcal{A}}^{-1}\boldsymbol{\mathcal{A}}^{-T}\boldsymbol{\mathcal{A}}^{T}=\boldsymbol{\mathcal{A}}^{-1}.
\end{align}
Multiplying the expression \eqref{eq:SolutionForx} for $\boldsymbol{x}$
from the left by $\boldsymbol{\mathcal{A}}$ as on the l. h. s. of
Eq. \eqref{eq:AxEqualsb} gives
\begin{align}
\boldsymbol{\mathcal{A}}\boldsymbol{x} & =\boldsymbol{\mathcal{A}}\boldsymbol{\mathcal{A}}^{+}\boldsymbol{b}=\boldsymbol{\mathcal{A}}\left(\boldsymbol{\mathcal{A}}^{T}\boldsymbol{\mathcal{A}}\right)^{-1}\boldsymbol{\mathcal{A}}^{T}\boldsymbol{b}.\label{eq:Eq41}
\end{align}
Note that
\begin{align}
\boldsymbol{\mathcal{P}} & =\boldsymbol{\mathcal{A}}\boldsymbol{\mathcal{A}}^{+}=\boldsymbol{\mathcal{A}}\left(\boldsymbol{\mathcal{A}}^{T}\boldsymbol{\mathcal{A}}\right)^{-1}\boldsymbol{\mathcal{A}}^{T}
\end{align}
is a projector, i. e., it is an idempotent $n\times n$ matrix, 
\begin{align}
\boldsymbol{\mathcal{P}}^{2} & =\boldsymbol{\mathcal{A}}\left(\boldsymbol{\mathcal{A}}^{T}\boldsymbol{\mathcal{A}}\right)^{-1}\boldsymbol{\mathcal{A}}^{T}\boldsymbol{\mathcal{A}}\left(\boldsymbol{\mathcal{A}}^{T}\boldsymbol{\mathcal{A}}\right)^{-1}\boldsymbol{\mathcal{A}}^{T}=\boldsymbol{\mathcal{A}}\left(\boldsymbol{\mathcal{A}}^{T}\boldsymbol{\mathcal{A}}\right)^{-1}\boldsymbol{\mathcal{A}}^{T}=\boldsymbol{\mathcal{P}}.
\end{align}
Furthermore, $\boldsymbol{\mathcal{P}}$ is symmetric
\begin{align}
\boldsymbol{\mathcal{P}}^{T} & =\left(\boldsymbol{\mathcal{A}}\left(\boldsymbol{\mathcal{A}}^{T}\boldsymbol{\mathcal{A}}\right)^{-1}\boldsymbol{\mathcal{A}}^{T}\right)^{T}=\boldsymbol{\mathcal{A}}\left(\boldsymbol{\mathcal{A}}^{T}\boldsymbol{\mathcal{A}}\right)^{-T}\boldsymbol{\mathcal{A}}^{T}=\boldsymbol{\mathcal{A}}\left(\boldsymbol{\mathcal{A}}^{T}\boldsymbol{\mathcal{A}}\right)^{-1}\boldsymbol{\mathcal{A}}^{T}=\boldsymbol{\mathcal{P}}.
\end{align}
Note that the inverse of the symmetric matrix $\boldsymbol{\mathcal{A}}^{T}\boldsymbol{\mathcal{A}}$
is also symmetric. The projector $\boldsymbol{\mathcal{P}}$ has rank
\begin{align}
\text{rank}\left(\boldsymbol{\mathcal{P}}\right) & =p.
\end{align}
A projector $\boldsymbol{\mathcal{Q}}$ complementary to $\boldsymbol{\mathcal{P}}$
can be defined as
\begin{align}
\boldsymbol{\mathcal{Q}} & =\mathbf{1}-\boldsymbol{\mathcal{P}},
\end{align}
which is also idempotent and symmetric and has rank
\begin{align}
\text{rank}\left(\boldsymbol{\mathcal{Q}}\right) & =n-p.
\end{align}
With the help of these projectors, the l. h. s. of Eq. \eqref{eq:AxEqualsb}
can be written as
\begin{align}
\boldsymbol{\mathcal{A}}\boldsymbol{x} & =\boldsymbol{\mathcal{P}}\boldsymbol{b}=\boldsymbol{b}-\boldsymbol{\mathcal{Q}}\boldsymbol{b}.
\end{align}
According to \eqref{eq:AxEqualsb}, this should be equal to $\boldsymbol{b}$,
which, of course, can only be true if 
\begin{align}
\boldsymbol{\mathcal{Q}}\boldsymbol{b} & =\mathbf{0}.\label{eq:QbEqualsZero}
\end{align}
In general, Eq. \eqref{eq:QbEqualsZero} is not true, and therefore
the ``solution'' Eq. \eqref{eq:SolutionForx} cannot be an exact
solution. In fact, Eq. \eqref{eq:QbEqualsZero} is the condition for
an exact solution to exist. That means that either $\boldsymbol{b}$
is the null vector, or the matrix $\boldsymbol{\mathcal{Q}}$ is the
null matrix. The third possibility is that $\boldsymbol{b}$ lies
in the null space of $\boldsymbol{\mathcal{Q}}$.

The expression Eq. \eqref{eq:SolutionForx} can be understood to give
an optimal approximate solution in the least square sense. In the
following, we demonstrate that $\boldsymbol{x}=\boldsymbol{\mathcal{A}}^{+}\boldsymbol{b}$
is the solution to the minimization problem
\begin{align}
\min_{\boldsymbol{x}} & \;\frac{1}{2}\left(\boldsymbol{\mathcal{A}}\boldsymbol{x}-\boldsymbol{b}\right)^{2}.
\end{align}
Define the scalar function $\mathcal{J}$ as 
\begin{align}
\mathcal{J}\left(\boldsymbol{x}\right) & =\frac{1}{2}\left(\boldsymbol{\mathcal{A}}\boldsymbol{x}-\boldsymbol{b}\right)^{2}=\frac{1}{2}\left(\boldsymbol{x}^{T}\boldsymbol{\mathcal{A}}^{T}\boldsymbol{\mathcal{A}}\boldsymbol{x}-2\boldsymbol{b}^{T}\boldsymbol{\mathcal{A}}\boldsymbol{x}+\boldsymbol{b}^{T}\boldsymbol{b}\right).
\end{align}
The Jacobian $\nabla\mathcal{J}$ of $\mathcal{J}$ with respect to
$\boldsymbol{x}$ is given by 
\begin{align}
\nabla\mathcal{J}\left(\boldsymbol{x}\right) & =\boldsymbol{x}^{T}\boldsymbol{\mathcal{A}}^{T}\boldsymbol{\mathcal{A}}-\boldsymbol{b}^{T}\boldsymbol{\mathcal{A}}.
\end{align}
The function $\mathcal{J}$ attains its extremum whenever
\begin{align}
\nabla\mathcal{J}\left(\boldsymbol{x}\right) & =\mathbf{0}.
\end{align}
Consequently, the vector $\boldsymbol{x}$ for which $\mathcal{J}$
attains its extremum must satisfy the equation 
\begin{align}
\boldsymbol{x}^{T}\boldsymbol{\mathcal{A}}^{T}\boldsymbol{\mathcal{A}} & =\boldsymbol{b}^{T}\boldsymbol{\mathcal{A}},
\end{align}
or, after transposing,
\begin{align}
\boldsymbol{\mathcal{A}}^{T}\boldsymbol{\mathcal{A}}\boldsymbol{x} & =\boldsymbol{\mathcal{A}}^{T}\boldsymbol{b}.
\end{align}
Solving for $\boldsymbol{x}$ indeed yields the expression Eq. \eqref{eq:SolutionForx},
\begin{align}
\boldsymbol{x} & =\left(\boldsymbol{\mathcal{A}}^{T}\boldsymbol{\mathcal{A}}\right)^{-1}\boldsymbol{\mathcal{A}}^{T}\boldsymbol{b}=\boldsymbol{\mathcal{A}}^{+}\boldsymbol{b}.
\end{align}
It remains to check if the extremum is indeed a minimum. Computing
the Hessian matrix $\nabla^{2}\mathcal{J}$ of $\mathcal{J}$ yields
\begin{align}
\nabla^{2}\mathcal{J}\left(\boldsymbol{x}\right) & =\boldsymbol{\mathcal{A}}^{T}\boldsymbol{\mathcal{A}}.
\end{align}
For an arbitrary matrix $\boldsymbol{\mathcal{A}}$, $\boldsymbol{\mathcal{A}}^{T}\boldsymbol{\mathcal{A}}$
is a positive semidefinite matrix. It becomes a positive definite
matrix if $\boldsymbol{\mathcal{A}}^{T}\boldsymbol{\mathcal{A}}$
is nonsingular, or, equivalently, if $\boldsymbol{\mathcal{A}}$ has
full rank, $\text{rank}\left(\boldsymbol{\mathcal{A}}\right)=p$ \cite{chen1995linear}.
Therefore, $\boldsymbol{x}=\boldsymbol{\mathcal{A}}^{+}\boldsymbol{b}$
indeed minimizes $\mathcal{J}$.

In conclusion, the linear equation $\boldsymbol{\mathcal{A}}\boldsymbol{x}=\boldsymbol{b}$
is discussed. An optimal solution for $\boldsymbol{x}$, which minimizes
the squared difference $\left(\boldsymbol{\mathcal{A}}\boldsymbol{x}-\boldsymbol{b}\right)^{2}$,
exists as long as $\boldsymbol{\mathcal{A}}^{T}\boldsymbol{\mathcal{A}}$
is positive definite and is given by $\boldsymbol{x}=\left(\boldsymbol{\mathcal{A}}^{T}\boldsymbol{\mathcal{A}}\right)^{-1}\boldsymbol{\mathcal{A}}^{T}\boldsymbol{b}$.
An exact solution for $\boldsymbol{x}$ can only exist if additionally,
the vector $\boldsymbol{b}$ satisfies the constraint $\left(\boldsymbol{1}-\boldsymbol{\mathcal{A}}\left(\boldsymbol{\mathcal{A}}^{T}\boldsymbol{\mathcal{A}}\right)^{-1}\boldsymbol{\mathcal{A}}^{T}\right)\boldsymbol{b}=\boldsymbol{0}$.

Finally, a slightly more general minimization problem is discussed.
The problem is to minimize 
\begin{align}
\min_{\boldsymbol{x}} & \;\frac{1}{2}\left(\boldsymbol{\mathcal{A}}\boldsymbol{x}-\boldsymbol{b}\right)^{T}\boldsymbol{\mathcal{S}}\left(\boldsymbol{\mathcal{A}}\boldsymbol{x}-\boldsymbol{b}\right),
\end{align}
with the symmetric $n\times n$ matrix $\boldsymbol{\mathcal{S}}^{T}=\boldsymbol{\mathcal{S}}$
of weighting coefficients. In the same manner as before, the scalar
function $\mathcal{J}_{\boldsymbol{\mathcal{S}}}$ is defined as 
\begin{align}
\mathcal{J}_{\boldsymbol{\mathcal{S}}}\left(\boldsymbol{x}\right) & =\frac{1}{2}\left(\boldsymbol{\mathcal{A}}\boldsymbol{x}-\boldsymbol{b}\right)^{T}\boldsymbol{\mathcal{S}}\left(\boldsymbol{\mathcal{A}}\boldsymbol{x}-\boldsymbol{b}\right)=\frac{1}{2}\left(\boldsymbol{x}^{T}\boldsymbol{\mathcal{A}}^{T}\boldsymbol{\mathcal{S}}\boldsymbol{\mathcal{A}}\boldsymbol{x}-2\boldsymbol{b}^{T}\boldsymbol{\mathcal{S}}\boldsymbol{\mathcal{A}}\boldsymbol{x}+\boldsymbol{b}^{T}\boldsymbol{\mathcal{S}}\boldsymbol{b}\right).
\end{align}
The function $\mathcal{J}_{\boldsymbol{\mathcal{S}}}$ attains its
extremum if
\begin{align}
\nabla\mathcal{J}_{\boldsymbol{\mathcal{S}}}\left(\boldsymbol{x}\right) & =\mathbf{0},
\end{align}
which gives
\begin{align}
\boldsymbol{\mathcal{A}}^{T}\boldsymbol{\mathcal{S}}\boldsymbol{\mathcal{A}}\boldsymbol{x} & =\boldsymbol{\mathcal{A}}^{T}\boldsymbol{\mathcal{S}}\boldsymbol{b}.\label{eq:EqA55}
\end{align}
As long as the $p\times p$ matrix $\boldsymbol{\mathcal{A}}^{T}\boldsymbol{\mathcal{S}}\boldsymbol{\mathcal{A}}$
has full rank,
\begin{align}
\text{rank}\left(\boldsymbol{\mathcal{A}}^{T}\boldsymbol{\mathcal{S}}\boldsymbol{\mathcal{A}}\right) & =p,
\end{align}
Eq. \eqref{eq:EqA55} can be solved for $\boldsymbol{x}$ to get 
\begin{align}
\boldsymbol{x} & =\left(\boldsymbol{\mathcal{A}}^{T}\boldsymbol{\mathcal{S}}\boldsymbol{\mathcal{A}}\right)^{-1}\boldsymbol{\mathcal{A}}^{T}\boldsymbol{\mathcal{S}}\boldsymbol{b}=\boldsymbol{\mathcal{A}}_{\boldsymbol{\mathcal{S}}}^{+}\boldsymbol{b}.\label{eq:xSSolution}
\end{align}
The generalized inverse $p\times n$ matrix $\boldsymbol{\mathcal{A}}_{\boldsymbol{\mathcal{S}}}^{+}$
is defined as
\begin{align}
\boldsymbol{\mathcal{A}}_{\boldsymbol{\mathcal{S}}}^{+} & =\left(\boldsymbol{\mathcal{A}}^{T}\boldsymbol{\mathcal{S}}\boldsymbol{\mathcal{A}}\right)^{-1}\boldsymbol{\mathcal{A}}^{T}\boldsymbol{\mathcal{S}}.
\end{align}
The question arises if $\boldsymbol{\mathcal{A}}_{\boldsymbol{\mathcal{S}}}^{+}$
is the Moore-Penrose pseudo inverse. Checking the four conditions
Eqs. \eqref{eq:MP1}-\eqref{eq:MP4} reveals that all conditions except
Eq. \eqref{eq:MP3} are satisfied. Consequently, $\boldsymbol{\mathcal{A}}_{\boldsymbol{\mathcal{S}}}^{+}$
is not a Moore-Penrose pseudo inverse but a generalized reflexive
inverse. What remains to check is if the extremum is indeed a minimum.
Computing the Hessian matrix $\nabla\nabla\mathcal{J}_{\boldsymbol{\mathcal{S}}}$
of $\mathcal{J}_{\boldsymbol{\mathcal{S}}}$ yields 
\begin{align}
\nabla\nabla\mathcal{J}_{\boldsymbol{\mathcal{S}}}\left(\boldsymbol{x}\right) & =\boldsymbol{\mathcal{A}}^{T}\boldsymbol{\mathcal{S}}\boldsymbol{\mathcal{A}}.
\end{align}
Consequently, as long as $\boldsymbol{\mathcal{A}}^{T}\boldsymbol{\mathcal{S}}\boldsymbol{\mathcal{A}}$
is positive definite, the solution Eq. \eqref{eq:xSSolution} minimizes
$\mathcal{J}_{\boldsymbol{\mathcal{S}}}$. Note that a positive definite
matrix has always full rank and is invertible, such that the solution
Eq. \eqref{eq:xSSolution} exists. Similar as above, two complementary
projectors $\boldsymbol{\mathcal{P}}_{\boldsymbol{\mathcal{S}}}$
and $\boldsymbol{\mathcal{Q}}_{\boldsymbol{\mathcal{S}}}$ can be
defined as 
\begin{align}
\boldsymbol{\mathcal{P}}_{\boldsymbol{\mathcal{S}}} & =\boldsymbol{\mathcal{A}}\boldsymbol{\mathcal{A}}_{\boldsymbol{\mathcal{S}}}^{+}=\boldsymbol{\mathcal{A}}\left(\boldsymbol{\mathcal{A}}^{T}\boldsymbol{\mathcal{S}}\boldsymbol{\mathcal{A}}\right)^{-1}\boldsymbol{\mathcal{A}}^{T}\boldsymbol{\mathcal{S}}, & \boldsymbol{\mathcal{Q}}_{\boldsymbol{\mathcal{S}}} & =\boldsymbol{1}-\boldsymbol{\mathcal{P}}_{\boldsymbol{\mathcal{S}}}.
\end{align}
In contrast to the projectors $\boldsymbol{\mathcal{P}}$ and $\boldsymbol{\mathcal{Q}}$,
these projectors are not symmetric. For the optimal solution $\boldsymbol{x}$
to be an exact solution to $\boldsymbol{\mathcal{A}}\boldsymbol{x}=\boldsymbol{b}$,
the vector $\boldsymbol{b}$ has to satisfy an additional condition,
\begin{align}
\boldsymbol{\mathcal{Q}}_{\boldsymbol{\mathcal{S}}}\boldsymbol{b} & =\boldsymbol{0}.
\end{align}
Thus, two exact solutions to the linear equation $\boldsymbol{\mathcal{A}}\boldsymbol{x}=\boldsymbol{b}$
were found. The first solution is given by
\begin{align}
\boldsymbol{x} & =\boldsymbol{x}_{1}=\left(\boldsymbol{\mathcal{A}}^{T}\boldsymbol{\mathcal{A}}\right)^{-1}\boldsymbol{\mathcal{A}}^{T}\boldsymbol{b}, & \boldsymbol{\mathcal{Q}}\boldsymbol{b} & =\boldsymbol{0},
\end{align}
while the second solution is
\begin{align}
\boldsymbol{x} & =\boldsymbol{x}_{2}=\left(\boldsymbol{\mathcal{A}}^{T}\boldsymbol{\mathcal{S}}\boldsymbol{\mathcal{A}}\right)^{-1}\boldsymbol{\mathcal{A}}^{T}\boldsymbol{\mathcal{S}}\boldsymbol{b}, & \boldsymbol{\mathcal{Q}}_{\boldsymbol{\mathcal{S}}}\boldsymbol{b} & =\boldsymbol{0}.
\end{align}
The exact solution to $\boldsymbol{\mathcal{A}}\boldsymbol{x}=\boldsymbol{b}$
should be unique such that
\begin{align}
\boldsymbol{x}_{1} & =\boldsymbol{x}_{2}.
\end{align}
Indeed, computing their difference, multiplying by $\boldsymbol{\mathcal{A}}$,
and exploiting the relations for $\boldsymbol{b}$ yields
\begin{align}
\boldsymbol{\mathcal{A}}\left(\boldsymbol{x}_{1}-\boldsymbol{x}_{2}\right) & =\boldsymbol{\mathcal{A}}\left(\boldsymbol{\mathcal{A}}^{T}\boldsymbol{\mathcal{A}}\right)^{-1}\boldsymbol{\mathcal{A}}^{T}\boldsymbol{b}-\boldsymbol{\mathcal{A}}\left(\boldsymbol{\mathcal{A}}^{T}\boldsymbol{\mathcal{S}}\boldsymbol{\mathcal{A}}\right)^{-1}\boldsymbol{\mathcal{A}}^{T}\boldsymbol{\mathcal{S}}\boldsymbol{b}\nonumber \\
 & =\left(\boldsymbol{\mathcal{P}}-\boldsymbol{\mathcal{P}}_{\boldsymbol{\mathcal{S}}}\right)\boldsymbol{b}=\left(\boldsymbol{1}-\boldsymbol{\mathcal{Q}}-\boldsymbol{1}+\boldsymbol{\mathcal{Q}}_{\boldsymbol{\mathcal{S}}}\right)\boldsymbol{b}=\boldsymbol{0}.
\end{align}
This relation is true if either $\boldsymbol{x}_{1}=\boldsymbol{x}_{2}$,
or $\boldsymbol{x}_{1}-\boldsymbol{x}_{2}$ lies in the null space
of $\boldsymbol{\mathcal{A}}$. However, due to the assumption that
$\boldsymbol{\mathcal{A}}$ has full rank and as a consequence of
the rank-nullity theorem, the null space of $\boldsymbol{\mathcal{A}}$
has zero dimension and contains only the zero vector. Consequently,
the solutions are identical, $\boldsymbol{x}_{1}=\boldsymbol{x}_{2}$.
As expected, the exact solution is unique, and does not depend on
the matrix of weighting coefficients $\boldsymbol{\mathcal{S}}$.
The relations $\boldsymbol{\mathcal{Q}}\boldsymbol{b}=\boldsymbol{0}$
and $\boldsymbol{\mathcal{Q}}_{\boldsymbol{\mathcal{S}}}\boldsymbol{b}=\boldsymbol{0}$
are the analogues of the constraint equations for exactly realizable
desired trajectories introduced in Section \ref{sec:ExactlyRealizableTrajectories}
and Section \ref{sec:ExactlyRealizableTrajectoriesOptimalControl},
respectively.

\subsection{Solving an underdetermined system of equations}

The opposite problem is considered. An underdetermined system is a
system with fewer equations than unknowns. Let $\boldsymbol{x}\in\mathbb{R}^{p}$
and $\boldsymbol{b}\in\mathbb{R}^{n}$ with $p>n$, and let $\boldsymbol{\mathcal{A}}$
be an $n\times p$ matrix. The system of $n$ equations
\begin{align}
\boldsymbol{\mathcal{A}}\boldsymbol{x} & =\boldsymbol{b}\label{eq:AxEqualsb-1}
\end{align}
is to be solved for $\boldsymbol{x}$. Because there are fewer equations
than components of $\boldsymbol{x}$, not all components of $\boldsymbol{x}$
can be determined. Motivated by the example above, two complementary
projectors $\boldsymbol{\mathcal{M}}$ and $\boldsymbol{\mathcal{N}}$
are introduced as 
\begin{align}
\boldsymbol{\mathcal{M}} & =\boldsymbol{\mathcal{A}}^{T}\left(\boldsymbol{\mathcal{A}}\boldsymbol{\mathcal{A}}^{T}\right)^{-1}\boldsymbol{\mathcal{A}},\\
\boldsymbol{\mathcal{N}} & =\mathbf{1}-\boldsymbol{\mathcal{M}}.
\end{align}
These projectors are symmetric $p\times p$ matrices. Note that the
$n\times n$ matrix $\boldsymbol{\mathcal{A}}\boldsymbol{\mathcal{A}}^{T}$
has rank $n$ whenever $\boldsymbol{\mathcal{A}}$ has full row rank,
\begin{align}
\text{rank}\left(\boldsymbol{\mathcal{A}}\right) & =n.
\end{align}
The projectors $\boldsymbol{\mathcal{M}}$ and $\boldsymbol{\mathcal{N}}$
have rank
\begin{align}
\text{rank}\left(\boldsymbol{\mathcal{M}}\right) & =n, & \text{rank}\left(\boldsymbol{\mathcal{N}}\right) & =p-n.
\end{align}
Multiplying Eq. \eqref{eq:AxEqualsb-1} by $\boldsymbol{\mathcal{A}}^{T}\left(\boldsymbol{\mathcal{A}}\boldsymbol{\mathcal{A}}^{T}\right)^{-1}$
from the left yields
\begin{align}
\boldsymbol{\mathcal{A}}^{T}\left(\boldsymbol{\mathcal{A}}\boldsymbol{\mathcal{A}}^{T}\right)^{-1}\boldsymbol{b} & =\boldsymbol{\mathcal{A}}^{T}\left(\boldsymbol{\mathcal{A}}\boldsymbol{\mathcal{A}}^{T}\right)^{-1}\boldsymbol{\mathcal{A}}\boldsymbol{x}=\boldsymbol{\mathcal{M}}\boldsymbol{x}.
\end{align}
Thus, the part $\boldsymbol{\mathcal{M}}\boldsymbol{x}$ can be determined
in terms of $\boldsymbol{b}$, while the part $\boldsymbol{\mathcal{N}}\boldsymbol{x}$
must be left undetermined.

\section{\label{sec:PropertiesOfTimeDependentProjectors}Properties of time-dependent
projectors}

Some relations for the projectors $\boldsymbol{\mathcal{P}}\left(\boldsymbol{x}\right)$
and $\boldsymbol{\mathcal{Q}}\left(\boldsymbol{x}\right)$ are listed.
The projectors may depend on time though its argument $\boldsymbol{x}$.
First, the projectors are idempotent, 
\begin{align}
\boldsymbol{\mathcal{Q}}\left(\boldsymbol{x}\right)\boldsymbol{\mathcal{Q}}\left(\boldsymbol{x}\right) & =\boldsymbol{\mathcal{Q}}\left(\boldsymbol{x}\right), & \boldsymbol{\mathcal{P}}\left(\boldsymbol{x}\right)\boldsymbol{\mathcal{P}}\left(\boldsymbol{x}\right) & =\boldsymbol{\mathcal{P}}\left(\boldsymbol{x}\right),
\end{align}
and complementary,
\begin{align}
\boldsymbol{\mathcal{P}}\left(\boldsymbol{x}\right)+\boldsymbol{\mathcal{Q}}\left(\boldsymbol{x}\right) & =\boldsymbol{1}.
\end{align}
Applying the time derivative $\frac{d}{dt}$ to the last relation
yields

\begin{align}
\frac{d}{dt}\boldsymbol{\mathcal{P}}\left(\boldsymbol{x}\left(t\right)\right) & =-\frac{d}{dt}\boldsymbol{\mathcal{Q}}\left(\boldsymbol{x}\left(t\right)\right)\label{eq:E7}
\end{align}
or
\begin{align}
\nabla\boldsymbol{\mathcal{P}}\left(\boldsymbol{x}\left(t\right)\right)\boldsymbol{\dot{x}}\left(t\right) & =-\nabla\boldsymbol{\mathcal{Q}}\left(\boldsymbol{x}\left(t\right)\right)\boldsymbol{\dot{x}}\left(t\right)
\end{align}
or
\begin{align}
\nabla\boldsymbol{\mathcal{P}}\left(\boldsymbol{x}\left(t\right)\right) & =-\nabla\boldsymbol{\mathcal{Q}}\left(\boldsymbol{x}\left(t\right)\right).
\end{align}
Here, $\nabla\boldsymbol{\mathcal{P}}\left(\boldsymbol{x}\right)$
denotes the Jacobian of $\boldsymbol{\mathcal{P}}\left(\boldsymbol{x}\right)$
with respect to $\boldsymbol{x}$. Note that $\nabla\boldsymbol{\mathcal{P}}\left(\boldsymbol{x}\right)$
is a third order tensor. Some more relations for the time derivatives
of the projectors are given. To shorten the notation, the time-dependent
projectors are rewritten as 
\begin{align}
\boldsymbol{\mathcal{P}}\left(t\right) & =\boldsymbol{\mathcal{P}}\left(\boldsymbol{x}\left(t\right)\right), & \boldsymbol{\mathcal{Q}}\left(t\right) & =\boldsymbol{\mathcal{Q}}\left(\boldsymbol{x}\left(t\right)\right).
\end{align}
The time derivative is denoted as
\begin{align}
\boldsymbol{\mathcal{\dot{P}}}\left(t\right) & =\frac{d}{dt}\boldsymbol{\mathcal{P}}\left(\boldsymbol{x}\left(t\right)\right)=\nabla\boldsymbol{\mathcal{P}}\left(\boldsymbol{x}\left(t\right)\right)\boldsymbol{\dot{x}}\left(t\right).
\end{align}
From the complementarity property Eq. \eqref{eq:PTimesQ} follows
\begin{align}
\boldsymbol{\mathcal{\dot{P}}}\left(t\right)\boldsymbol{\mathcal{Q}}\left(t\right)+\boldsymbol{\mathcal{P}}\left(t\right)\boldsymbol{\mathcal{\dot{Q}}}\left(t\right) & =\boldsymbol{0}, & \boldsymbol{\mathcal{\dot{Q}}}\left(t\right)\boldsymbol{\mathcal{P}}\left(t\right)+\boldsymbol{\mathcal{Q}}\left(t\right)\boldsymbol{\mathcal{\dot{P}}}\left(t\right) & =\boldsymbol{0},\\
\boldsymbol{\mathcal{\dot{Q}}}\left(t\right)\boldsymbol{\mathcal{Q}}\left(t\right)+\boldsymbol{\mathcal{Q}}\left(t\right)\boldsymbol{\mathcal{\dot{Q}}}\left(t\right) & =\boldsymbol{\mathcal{\dot{Q}}}\left(t\right), & \boldsymbol{\mathcal{\dot{P}}}\left(t\right)\boldsymbol{\mathcal{P}}\left(t\right)+\boldsymbol{\mathcal{P}}\left(t\right)\boldsymbol{\mathcal{\dot{P}}}\left(t\right) & =\boldsymbol{\mathcal{\dot{P}}}\left(t\right).
\end{align}
The last line yields
\begin{align}
\boldsymbol{\mathcal{Q}}\left(t\right)\boldsymbol{\mathcal{\dot{Q}}}\left(t\right)\boldsymbol{\mathcal{Q}}\left(t\right)+\boldsymbol{\mathcal{Q}}\left(t\right)\boldsymbol{\mathcal{Q}}\left(t\right)\boldsymbol{\mathcal{\dot{Q}}}\left(t\right) & =\boldsymbol{\mathcal{Q}}\left(t\right)\boldsymbol{\mathcal{\dot{Q}}}\left(t\right),
\end{align}
or
\begin{align}
\boldsymbol{\mathcal{Q}}\left(t\right)\boldsymbol{\mathcal{\dot{Q}}}\left(t\right)\boldsymbol{\mathcal{Q}}\left(t\right) & =\mathbf{0}.
\end{align}
A similar computation results in
\begin{align}
\boldsymbol{\mathcal{P}}\left(t\right)\boldsymbol{\mathcal{\dot{Q}}}\left(t\right)\boldsymbol{\mathcal{Q}}\left(t\right)+\boldsymbol{\mathcal{P}}\left(t\right)\boldsymbol{\mathcal{Q}}\left(t\right)\boldsymbol{\mathcal{\dot{Q}}}\left(t\right) & =\boldsymbol{\mathcal{P}}\left(t\right)\boldsymbol{\mathcal{\dot{Q}}}\left(t\right),
\end{align}
or
\begin{align}
\boldsymbol{\mathcal{P}}\left(t\right)\boldsymbol{\mathcal{\dot{Q}}}\left(t\right)\boldsymbol{\mathcal{Q}}\left(t\right) & =\boldsymbol{\mathcal{P}}\left(t\right)\boldsymbol{\mathcal{\dot{Q}}}\left(t\right).
\end{align}
With the help of the projectors Eqs. \eqref{eq:TimeDependentP}, \eqref{eq:TimeDependentQ},
the time derivative of $\boldsymbol{x}\left(t\right)$ can be written
as
\begin{align}
\frac{d}{dt}\boldsymbol{x}\left(t\right) & =\frac{d}{dt}\left(\boldsymbol{\mathcal{P}}\left(t\right)\boldsymbol{x}\left(t\right)+\boldsymbol{\mathcal{Q}}\left(t\right)\boldsymbol{x}\left(t\right)\right)\nonumber \\
 & =\boldsymbol{\mathcal{\dot{P}}}\left(t\right)\boldsymbol{x}\left(t\right)+\boldsymbol{\mathcal{\dot{Q}}}\left(t\right)\boldsymbol{x}\left(t\right)+\boldsymbol{\mathcal{P}}\left(t\right)\boldsymbol{\dot{x}}\left(t\right)+\boldsymbol{\mathcal{Q}}\left(t\right)\boldsymbol{\dot{x}}\left(t\right).
\end{align}
Applying $\boldsymbol{\mathcal{Q}}\left(t\right)$ from the left and
using Eq. \eqref{eq:E7} yields
\begin{align}
\boldsymbol{\mathcal{Q}}\left(t\right)\boldsymbol{\mathcal{\dot{P}}}\left(t\right)\boldsymbol{x}\left(t\right)+\boldsymbol{\mathcal{Q}}\left(t\right)\boldsymbol{\mathcal{\dot{Q}}}\left(t\right)\boldsymbol{x}\left(t\right)+\boldsymbol{\mathcal{Q}}\left(t\right)\boldsymbol{\dot{x}}\left(t\right) & =\boldsymbol{\mathcal{Q}}\left(t\right)\boldsymbol{\dot{x}}\left(t\right).
\end{align}
Similarly, applying $\boldsymbol{\mathcal{P}}\left(t\right)$ from
the left gives 
\begin{align}
\boldsymbol{\mathcal{P}}\left(t\right)\boldsymbol{\mathcal{\dot{P}}}\left(t\right)\boldsymbol{x}\left(t\right)+\boldsymbol{\mathcal{P}}\left(t\right)\boldsymbol{\mathcal{\dot{Q}}}\left(t\right)\boldsymbol{x}\left(t\right)+\boldsymbol{\mathcal{P}}\left(t\right)\boldsymbol{\dot{x}}\left(t\right) & =\boldsymbol{\mathcal{P}}\left(t\right)\boldsymbol{\dot{x}}\left(t\right).
\end{align}

\section{\label{sec:DiagonalizingTheConstraintEquation}Diagonalizing the
projectors \texorpdfstring{$\boldsymbol{\mathcal{P}}\left(\boldsymbol{x}\right)$}{P(x)}
and \texorpdfstring{$\boldsymbol{\mathcal{Q}}\left(\boldsymbol{x}\right)$}{Q(x)}}

Let the $n\times n$ matrix $\boldsymbol{\mathcal{Q}}\left(\boldsymbol{x}\right)$
be the projector
\begin{align}
\boldsymbol{\mathcal{Q}}\left(\boldsymbol{x}\right) & =\boldsymbol{1}-\boldsymbol{\mathcal{B}}\left(\boldsymbol{x}\right)\boldsymbol{\mathcal{B}}^{+}\left(\boldsymbol{x}\right),
\end{align}
with $\boldsymbol{\mathcal{B}}^{+}\left(\boldsymbol{x}\right)$ the
Moore-Penrose pseudo inverse of the coupling matrix $\boldsymbol{\mathcal{B}}\left(\boldsymbol{x}\right)$.
$\boldsymbol{\mathcal{Q}}\left(\boldsymbol{x}\right)$ is idempotent,
\begin{align}
\boldsymbol{\mathcal{Q}}\left(\boldsymbol{x}\right)\boldsymbol{\mathcal{Q}}\left(\boldsymbol{x}\right) & =\boldsymbol{\mathcal{Q}}\left(\boldsymbol{x}\right),
\end{align}
and has a complementary projector defined by
\begin{align}
\boldsymbol{\mathcal{P}}\left(\boldsymbol{x}\right) & =\boldsymbol{1}-\boldsymbol{\mathcal{Q}}\left(\boldsymbol{x}\right).
\end{align}
Furthermore, $\boldsymbol{\mathcal{Q}}\left(\boldsymbol{x}\right)$
satisfies
\begin{align}
\boldsymbol{\mathcal{Q}}\left(\boldsymbol{x}\right)\boldsymbol{\mathcal{B}}\left(\boldsymbol{x}\right) & =\boldsymbol{\mathcal{B}}\left(\boldsymbol{x}\right)-\boldsymbol{\mathcal{B}}\left(\boldsymbol{x}\right)\boldsymbol{\mathcal{B}}^{+}\left(\boldsymbol{x}\right)\boldsymbol{\mathcal{B}}\left(\boldsymbol{x}\right)=\boldsymbol{0}.
\end{align}
Assume that the rank of $\boldsymbol{\mathcal{Q}}\left(\boldsymbol{x}\right)$
is, with $p\leq n$,
\begin{align}
\text{rank}\left(\boldsymbol{\mathcal{Q}}\left(\boldsymbol{x}\right)\right) & =n-p,
\end{align}
for all $\boldsymbol{x}$.

Any projector $\boldsymbol{\mathcal{Q}}\left(\boldsymbol{x}\right)$
can be diagonalized with zeros and ones as the diagonal entries \cite{fischer2008lineare,liesen2011lineare}.
Let $\boldsymbol{\mathcal{T}}\left(\boldsymbol{x}\right)$ be the
$n\times n$ matrix which diagonalizes the projector $\boldsymbol{\mathcal{Q}}\left(\boldsymbol{x}\right)$,
\begin{align}
\boldsymbol{\mathcal{Q}}_{D} & =\boldsymbol{\mathcal{T}}{}^{-1}\left(\boldsymbol{x}\right)\boldsymbol{\mathcal{Q}}\left(\boldsymbol{x}\right)\boldsymbol{\mathcal{T}}\left(\boldsymbol{x}\right)=\left(\begin{array}{ccccc}
0 & \cdots & 0 & \cdots & 0\\
\vdots & \ddots & \vdots & \vdots & \vdots\\
0 & \cdots & 1 & \cdots & \vdots\\
\vdots & \vdots & \vdots & \ddots & \vdots\\
0 & \cdots & 0 & \cdots & 1
\end{array}\right),
\end{align}
such that the first $p$ diagonal elements are zero, while the last
$n-p$ diagonal elements are one. The same matrix $\boldsymbol{\mathcal{T}}\left(\boldsymbol{x}\right)$
diagonalizes the projector $\boldsymbol{\mathcal{P}}$ as well,
\begin{align}
\boldsymbol{\mathcal{P}}_{D} & =\boldsymbol{\mathcal{T}}{}^{-1}\left(\boldsymbol{x}\right)\boldsymbol{\mathcal{P}}\boldsymbol{\mathcal{T}}\left(\boldsymbol{x}\right)\nonumber \\
 & =\boldsymbol{\mathcal{T}}{}^{-1}\left(\boldsymbol{x}\right)\left(\boldsymbol{1}-\boldsymbol{\mathcal{Q}}\right)\boldsymbol{\mathcal{T}}\left(\boldsymbol{x}\right)=\boldsymbol{1}-\boldsymbol{\mathcal{Q}}_{D}=\left(\begin{array}{ccccc}
1 & \cdots & 0 & \cdots & 0\\
\vdots & \ddots & \vdots & \vdots & \vdots\\
0 & \cdots & 1 & \cdots & \vdots\\
\vdots & \vdots & \vdots & \ddots & \vdots\\
0 & \cdots & 0 & \cdots & 0
\end{array}\right),
\end{align}
such that the first $p$ diagonal elements are one, while the last
$n-p$ diagonal elements are zero. The matrix $\boldsymbol{\mathcal{T}}\left(\boldsymbol{x}\right)$
contains $n$ linearly independent eigenvectors of $\boldsymbol{\mathcal{Q}}\left(\boldsymbol{x}\right)$
and can be constructed as follows. Let $\boldsymbol{q}_{i}\left(\boldsymbol{x}\right)$
denote the $i$-th column of the matrix $\boldsymbol{\mathcal{Q}}\left(\boldsymbol{x}\right)$,
\begin{align}
\boldsymbol{\mathcal{Q}}\left(\boldsymbol{x}\right) & =\left(\begin{array}{ccc}
\boldsymbol{q}_{1}\left(\boldsymbol{x}\right), & \dots, & \boldsymbol{q}_{n}\left(\boldsymbol{x}\right)\end{array}\right),
\end{align}
or, written component wise,
\begin{align}
\mathcal{Q}_{ij}\left(\boldsymbol{x}\right) & =\left(\boldsymbol{q}_{j}\right)_{i}\left(\boldsymbol{x}\right)=q_{j,i}\left(\boldsymbol{x}\right).
\end{align}
The vectors $\boldsymbol{q}_{i}\left(\boldsymbol{x}\right)$ are eigenvectors
of $\boldsymbol{\mathcal{Q}}\left(\boldsymbol{x}\right)$. Indeed,
$\boldsymbol{\mathcal{Q}}\left(\boldsymbol{x}\right)$ is idempotent,
\begin{align}
\boldsymbol{\mathcal{Q}}\left(\boldsymbol{x}\right)\boldsymbol{\mathcal{Q}}\left(\boldsymbol{x}\right) & =\boldsymbol{\mathcal{Q}}\left(\boldsymbol{x}\right),
\end{align}
or, written component wise,
\begin{align}
\sum_{j=1}^{n}\mathcal{Q}_{ij}\left(\boldsymbol{x}\right)\mathcal{Q}_{jk}\left(\boldsymbol{x}\right) & =\mathcal{Q}_{ik}\left(\boldsymbol{x}\right).
\end{align}
Expressed in terms of the vectors $\boldsymbol{q}_{i}\left(\boldsymbol{x}\right)$,
the last relation becomes 
\begin{align}
\sum_{j=1}^{n}\mathcal{Q}_{ij}\left(\boldsymbol{x}\right)q_{k,j}\left(\boldsymbol{x}\right) & =q_{k,i}\left(\boldsymbol{x}\right),
\end{align}
or
\begin{align}
\boldsymbol{\mathcal{Q}}\left(\boldsymbol{x}\right)\boldsymbol{q}_{k}\left(\boldsymbol{x}\right) & =\boldsymbol{q}_{k}\left(\boldsymbol{x}\right),
\end{align}
which shows that the vectors $\boldsymbol{q}_{k}\left(\boldsymbol{x}\right)$
are eigenvectors of $\boldsymbol{\mathcal{Q}}\left(\boldsymbol{x}\right)$
to eigenvalue one. However, because $\boldsymbol{\mathcal{Q}}\left(\boldsymbol{x}\right)$
has $\text{rank}\left(\boldsymbol{\mathcal{Q}}\left(\boldsymbol{x}\right)\right)=n-p$,
only $n-p$ vectors out of $i=1,\dots,n$ vectors $\boldsymbol{q}_{i}\left(\boldsymbol{x}\right)$
are linearly independent. By appropriately ordering the eigenvectors,
one can ensure that the first $n-p$ eigenvectors $\boldsymbol{q}_{1}\left(\boldsymbol{x}\right),\dots,\boldsymbol{q}_{n-p}\left(\boldsymbol{x}\right)$
are linearly independent. The remaining eigenvectors can be constructed
from the coupling matrix $\boldsymbol{\mathcal{B}}\left(\boldsymbol{x}\right)$.
The $n\times p$ matrix $\boldsymbol{\mathcal{B}}\left(\boldsymbol{x}\right)$
can be written in terms of its $p$ column vectors as
\begin{align}
\boldsymbol{\mathcal{B}}\left(\boldsymbol{x}\right) & =\left(\begin{array}{ccc}
\boldsymbol{b}_{1}\left(\boldsymbol{x}\right), & \dots, & \boldsymbol{b}_{p}\left(\boldsymbol{x}\right)\end{array}\right),
\end{align}
or, written component wise, 
\begin{align}
\mathcal{B}_{ij}\left(\boldsymbol{x}\right) & =\left(\boldsymbol{b}_{j}\right)_{i}\left(\boldsymbol{x}\right)=b_{j,i}\left(\boldsymbol{x}\right).
\end{align}
From the relation
\begin{align}
\boldsymbol{\mathcal{Q}}\left(\boldsymbol{x}\right)\boldsymbol{\mathcal{B}}\left(\boldsymbol{x}\right) & =\boldsymbol{0},
\end{align}
or
\begin{align}
\sum_{j=1}^{n}\mathcal{Q}_{ij}\left(\boldsymbol{x}\right)\mathcal{B}_{jk}\left(\boldsymbol{x}\right) & =0,
\end{align}
follows that the vectors $\boldsymbol{b}_{i}\left(\boldsymbol{x}\right)$
are indeed eigenvectors to eigenvalue zero, 
\begin{align}
\sum_{j=1}^{n}\mathcal{Q}_{ij}\left(\boldsymbol{x}\right)b_{k,j}\left(\boldsymbol{x}\right) & =0,
\end{align}
or
\begin{align}
\boldsymbol{\mathcal{Q}}\left(\boldsymbol{x}\right)\boldsymbol{b}_{k}\left(\boldsymbol{x}\right) & =0.
\end{align}
Finally, the matrix $\boldsymbol{\mathcal{T}}\left(\boldsymbol{x}\right)$
becomes 
\begin{align}
\boldsymbol{\mathcal{T}}\left(\boldsymbol{x}\right) & =\left(\begin{array}{cccccc}
\boldsymbol{b}_{1}\left(\boldsymbol{x}\right), & \dots, & \boldsymbol{b}_{p}\left(\boldsymbol{x}\right), & \boldsymbol{q}_{1}\left(\boldsymbol{x}\right), & \dots, & \boldsymbol{q}_{n-p}\left(\boldsymbol{x}\right)\end{array}\right).
\end{align}
If the projectors $\boldsymbol{\mathcal{P}}\left(\boldsymbol{x}\right)$
and $\boldsymbol{\mathcal{Q}}\left(\boldsymbol{x}\right)$ are not
diagonal, splitting up the vector $\boldsymbol{x}\left(t\right)$
as
\begin{align}
\boldsymbol{x}\left(t\right) & =\boldsymbol{\mathcal{P}}\left(\boldsymbol{x}\left(t\right)\right)\boldsymbol{x}\left(t\right)+\boldsymbol{\mathcal{Q}}\left(\boldsymbol{x}\left(t\right)\right)\boldsymbol{x}\left(t\right)
\end{align}
results in the parts $\boldsymbol{\mathcal{P}}\left(\boldsymbol{x}\left(t\right)\right)\boldsymbol{x}\left(t\right)$
and $\boldsymbol{\mathcal{Q}}\left(\boldsymbol{x}\left(t\right)\right)\boldsymbol{x}\left(t\right)$
being nonlinear combinations of the components of $\boldsymbol{x}\left(t\right)$.
Due to this nonlinear mixing, it is not clear which state components
belong to which part. If the projectors are diagonal, 
\begin{align}
\boldsymbol{x}\left(t\right) & =\boldsymbol{\mathcal{P}}_{D}\boldsymbol{x}\left(t\right)+\boldsymbol{\mathcal{Q}}_{D}\boldsymbol{x}\left(t\right),
\end{align}
the parts $\boldsymbol{\mathcal{P}}_{D}\boldsymbol{x}\left(t\right)$
and $\boldsymbol{\mathcal{Q}}_{D}\boldsymbol{x}\left(t\right)$ are
linear combinations of the state components. Furthermore, only the
first $p$ components of $\boldsymbol{\mathcal{P}}_{D}\boldsymbol{x}\left(t\right)$
and the last $n-p$ components of $\boldsymbol{\mathcal{Q}}_{D}\boldsymbol{x}\left(t\right)$
are nonzero, and all other components vanish. Thus, diagonal projectors
allow a clear interpretation which components of $\boldsymbol{x}\left(t\right)$
belong to which part. If the projectors $\boldsymbol{\mathcal{P}}\left(\boldsymbol{x}\right)$
and $\boldsymbol{\mathcal{Q}}\left(\boldsymbol{x}\right)$ are not
diagonal, the matrix $\boldsymbol{\mathcal{T}}\left(\boldsymbol{x}\right)$
defines a coordinate transformation as follows. Let the vector $\boldsymbol{y}\left(t\right)$
be defined by
\begin{align}
\boldsymbol{y}\left(t\right) & =\boldsymbol{\mathcal{T}}^{-1}\left(\boldsymbol{x}\left(t\right)\right)\boldsymbol{x}\left(t\right).\label{eq:A107}
\end{align}
According to the construction above, the matrix $\boldsymbol{\mathcal{T}}\left(\boldsymbol{x}\left(t\right)\right)$
always exists and is invertible. The inverse relation of Eq. \eqref{eq:A107}
is 
\begin{align}
\boldsymbol{x}\left(t\right) & =\boldsymbol{\mathcal{T}}\left(\boldsymbol{x}\left(t\right)\right)\boldsymbol{y}\left(t\right).
\end{align}
Splitting up $\boldsymbol{x}\left(t\right)$ in Eq. \eqref{eq:A107}
with the help of the projectors $\boldsymbol{\mathcal{P}}\left(\boldsymbol{x}\left(t\right)\right)$
and $\boldsymbol{\mathcal{Q}}\left(\boldsymbol{x}\left(t\right)\right)$
gives 
\begin{align}
\boldsymbol{y}\left(t\right) & =\boldsymbol{\mathcal{T}}^{-1}\left(\boldsymbol{x}\left(t\right)\right)\boldsymbol{x}\left(t\right)\nonumber \\
 & =\boldsymbol{\mathcal{T}}^{-1}\left(\boldsymbol{x}\left(t\right)\right)\boldsymbol{\mathcal{P}}\left(\boldsymbol{x}\left(t\right)\right)\boldsymbol{x}\left(t\right)+\boldsymbol{\mathcal{T}}^{-1}\left(\boldsymbol{x}\left(t\right)\right)\boldsymbol{\mathcal{Q}}\left(\boldsymbol{x}\left(t\right)\right)\boldsymbol{x}\left(t\right)\nonumber \\
 & =\boldsymbol{\mathcal{T}}^{-1}\left(\boldsymbol{x}\left(t\right)\right)\boldsymbol{\mathcal{P}}\left(\boldsymbol{x}\left(t\right)\right)\boldsymbol{\mathcal{T}}\left(\boldsymbol{x}\left(t\right)\right)\boldsymbol{y}\left(t\right)+\boldsymbol{\mathcal{T}}^{-1}\left(\boldsymbol{x}\left(t\right)\right)\boldsymbol{\mathcal{Q}}\left(\boldsymbol{x}\left(t\right)\right)\boldsymbol{\mathcal{T}}\left(\boldsymbol{x}\left(t\right)\right)\boldsymbol{y}\left(t\right)\nonumber \\
 & =\boldsymbol{\mathcal{P}}_{D}\boldsymbol{y}\left(t\right)+\boldsymbol{\mathcal{Q}}_{D}\boldsymbol{y}\left(t\right).
\end{align}
Thus, in the new coordinates, the state can be separated in two parts
$\boldsymbol{\mathcal{P}}_{D}\boldsymbol{y}\left(t\right)$ and $\boldsymbol{\mathcal{Q}}_{D}\boldsymbol{y}\left(t\right)$
which are linear combinations of the state components of $\boldsymbol{y}$.
The first $p$ components of $\boldsymbol{y}$ belong to $\boldsymbol{\mathcal{P}}_{D}\boldsymbol{y}\left(t\right)$
and the last $n-p$ components of $\boldsymbol{y}$ belong to $\boldsymbol{\mathcal{Q}}_{D}\boldsymbol{y}\left(t\right)$.
Such a representation can be viewed as a normal form suitable for
computations with affine control systems. Note that Eq. \eqref{eq:A107}
yields an explicit expression for the new coordinates $\boldsymbol{y}$
in terms of the old coordinates $\boldsymbol{x}$. To obtain $\boldsymbol{x}$
in terms of $\boldsymbol{y}$, Eq. \eqref{eq:A107} must be solved
for $\boldsymbol{x}$.